\RequirePackage[l2tabu, orthodox]{nag}		
\documentclass[12pt]{book}

\usepackage[T1]{fontenc}										
\usepackage{lmodern}												
\usepackage{microtype}											
\usepackage[frenchmath,LGRgreek]{mathastext}
\MTgreekfont{lmtt} 													
\Mathastext
\let\varepsilon\epsilon	
\let\varphi\phi		
\let\vartheta\theta											
\usepackage[mathscr]{eucal}                 
\usepackage[sans]{dsfont}										
\usepackage{amsbsy}													
\usepackage{natbib}       									
\usepackage{amsmath}                        
\numberwithin{equation}{chapter}
\usepackage{graphicx}                       
\usepackage{tikz}                           
\usepackage{mdwlist}												
\usepackage[dvipsnames]{xcolor}							
\usepackage[plainpages=false, pdfpagelabels]{hyperref} 
	\hypersetup{
		colorlinks   = true,
		citecolor    = RoyalBlue,
		linkcolor    = RubineRed,
		urlcolor     = Turquoise
		}
\usepackage[paperwidth=8.5in,paperheight=11in,top=1.00in, bottom=1.00in, left=0.50in, right=0.50in]{geometry}
\usepackage{mathtools}                      
\mathtoolsset{showonlyrefs=true}            
\usepackage{amsthm}                         
\allowdisplaybreaks                         
\newtheoremstyle{plain}
  {}   
  {}   
  {\itshape}  
  {}       
  {\mdseries\scshape} 
  {.}         
  { } 
  {\thmname{#1}\thmnumber{ #2}\ifx#3\empty\else\ (#3)\fi}
\theoremstyle{plain}
\newtheorem{theorem}{\underline{Theorem}}
\numberwithin{theorem}{section}
\newtheorem{lemma}[theorem]{\underline{Lemma}}   		
\newtheorem{proposition}[theorem]{\underline{Proposition}}
\newtheorem{corollary}[theorem]{\underline{Corollary}}
\newtheoremstyle{definition}
  {}   
  {}   
  {}  
  {}       
  {\mdseries\scshape} 
  {.}         
  { } 
  {\thmname{#1}\thmnumber{ #2}\ifx#3\empty\else\ (#3)\fi}
\theoremstyle{definition}
\newtheorem{definition}[theorem]{\underline{Definition}}
\newtheorem{example}[theorem]{\underline{Example}}
\newtheorem{remark}[theorem]{\underline{Remark}}
\newtheorem{assumption}[theorem]{\underline{Assumption}}
\newtheorem{exercise}{\underline{Exercise}}
\numberwithin{exercise}{chapter}

\expandafter\let\expandafter\oldproof\csname\string\proof\endcsname
\let\oldendproof\endproof
\renewenvironment{proof}[1][\proofname]{%
  \oldproof[\scshape\underline{#1}]%
}{\oldendproof}

\usepackage{sectsty}
\allsectionsfont{\mdseries\scshape}

\newcommand{\<}{\langle}
\renewcommand{\>}{\rangle}
\renewcommand{\(}{\left(}
\renewcommand{\)}{\right)}
\renewcommand{\[}{\left[}
\renewcommand{\]}{\right]}

\newcommand\Cb{\mathds{C}}
\newcommand\Eb{\mathds{E}}
\newcommand\Fb{\mathds{F}}
\newcommand\Gb{\mathds{G}}
\newcommand\Ib{\mathds{1}} 	
\newcommand\Pb{\mathds{P}}

\newcommand\Rb{\mathds{R}}
\newcommand\Zb{\mathds{Z}}
\newcommand\Nb{\mathds{N}}
\newcommand\Vb{\mathds{V}}
\newcommand\Ub{\mathds{U}}

\newcommand\Ac{\mathscr{A}}
\newcommand\Bc{\mathscr{B}}
\newcommand\Cc{\mathscr{C}}
\newcommand\Dc{\mathscr{D}}
\newcommand\Ec{\mathscr{E}}
\newcommand\Fc{\mathscr{F}}
\newcommand\Gc{\mathscr{G}}
\newcommand\Hc{\mathscr{H}}
\newcommand\Lc{\mathscr{L}}

\newcommand\Nc{\mathscr{N}}
\newcommand\Oc{\mathscr{O}}
\newcommand\Pc{\mathscr{P}}
\newcommand\Rc{\mathscr{R}}

\newcommand\Tc{\mathscr{T}}
\newcommand\Uc{\mathscr{U}}

\newcommand\eps{\varepsilon}
\newcommand\om{\omega}
\newcommand\Om{\Omega}
\newcommand\sig{\sigma}
\newcommand\Sig{\Sigma}
\newcommand\Lam{\Lambda}
\newcommand\gam{\gamma}
\newcommand\Gam{\Gamma}
\newcommand\lam{\lambda}
\newcommand\del{\delta}
\newcommand\Del{\Delta}

\newcommand\Wb{\overline{W}}

\newcommand\rhob{\bar{\rho}}
\newcommand\taub{\bar{\tau}}

\newcommand\Fv{\pmb{F}}
\newcommand\ev{\pmb{e}}

\newcommand\Gv{\pmb{G}}
\newcommand\Iv{\pmb{I}}
\newcommand\Pv{\pmb{P}}
\newcommand\Qv{\pmb{Q}}

\newcommand\Hv{\pmb{H}}
\newcommand\Cv{\pmb{C}}
\newcommand\mv{\pmb{m}}
\newcommand\Uv{\pmb{U}}
\newcommand\muv{\pmb{\mu}}
\newcommand\nuv{\pmb{\nu}}
\newcommand\piv{\pmb{\pi}}
\newcommand\psiv{\pmb{\psi}}
\newcommand\lamv{\pmb{\lambda}}
\newcommand\Lamv{\pmb{\Lambda}}

\newcommand\Pbh{\widehat{\Pb}}
\newcommand\Ebh{\widehat{\Eb}}
\newcommand\varphih{\widehat{\varphi}}
\newcommand\uh{\widehat{u}}
\newcommand\Wh{\widehat{W}}

\newcommand\Nh{\widehat{N}}
\newcommand\Xh{\widehat{X}}
\newcommand\Yh{\widehat{Y}}
\newcommand\nuh{\widehat{\nu}}
\newcommand\muh{\widehat{\mu}}
\newcommand\ph{\widehat{p}}
\newcommand\qh{\widehat{q}}
\newcommand\sigh{\widehat{\sig}}
\newcommand\lamh{\widehat{\lam}}

\newcommand\Ebt{\widetilde{\Eb}}
\newcommand\Pbt{\widetilde{\Pb}}
\newcommand\Vbt{\widetilde{\Vb}}

\newcommand\Wt{\widetilde{W}}

\newcommand\Nt{\widetilde{N}}
\newcommand\Ft{\widetilde{F}}

\newcommand\pt{\widetilde{p}}

\newcommand\Xt{\widetilde{X}}
\newcommand\St{\widetilde{S}}

\newcommand\Wd{\dot{W}}

\renewcommand\d{\partial}
\newcommand{\ind}{\perp \! \! \! \perp}
\newcommand\ii{\mathtt{i}}
\newcommand\dd{\mathrm{d}}
\newcommand\ee{\mathrm{e}}
\newcommand\Id{\text{Id}}
\newcommand\es{\emptyset}
\newcommand\half{\frac{1}{2}}

\newcommand{\blu}[1]{\textcolor{blue}{#1}}

\newcommand{\norm}[1]{\left\|#1\right\|}
\renewcommand\limsup{\operatorname*{lim \, sup}}
\renewcommand\liminf{\operatorname*{lim \, inf}}

\begin{document}

\frontmatter
\title{An Introduction to Applied Probability and Stochastic Processes\\[0.5em] \large Open Lecture Notes}
\author{Matthew Lorig \thanks{Department of Applied Mathematics, University of Washington, Seattle, WA, USA.  \textbf{e-mail}: \url{mlorig@uw.edu}}}
\date{Version 1.70 (November 2026)}
\maketitle
\tableofcontents
\chapter{Preface} 
These notes provide a broad introduction to measure-theoretic probability and stochastic processes.
They were written for a two-quarter course for first-year PhD students in applied mathematics at the University of Washington, but they should also be useful to graduate students in statistics, engineering, economics and finance, and to anyone studying the subject on their own.
Readers are assumed to be familiar with undergraduate probability, real analysis and linear algebra; some exposure to ordinary and partial differential equations is helpful.
\\[0.5em]
The first half of the notes covers probability spaces, conditional expectation, generating and characteristic functions, Markov chains in discrete and continuous time, and convergence of random variables.
The second half covers Brownian motion, stochastic calculus, stochastic differential equations and their connection to partial differential equations, jump diffusions, stochastic control, optimal stopping and Monte Carlo methods.
Most chapters end with exercises, and solutions are collected at the end of the notes.
\\[0.5em]
The focus of these notes is on the \emph{applied} aspects of the subject.
Mathematical statements (definitions, theorems, propositions, lemmas and corollaries) are intended to be rigorously true.
Proofs, however, are sometimes only sketched, or make simplifying assumptions that avoid technical complications; where this is the case, we say so.
Each chapter begins with a list of references, which the interested reader can consult for complete proofs.
Some examples and exercises are adapted from these references.

\subsection*{An open and evolving set of notes}
These notes are released under the Creative Commons Attribution-ShareAlike 4.0 International License (CC BY-SA 4.0).
The \LaTeX\ source is available at\\
\url{https://github.com/mattlorig/intro-probability-stochastic-processes}.\\
Readers are encouraged to report errors, suggest improvements and contribute new material by opening an issue or a pull request; see the file \texttt{CONTRIBUTING.md} in the repository.
\blu{Text in blue} marks material that is still under development.
Information on how to cite these notes (including their arXiv identifier and DOI) can be found in the repository.

\subsection*{Acknowledgments}
The author of these notes wishes to express his sincere thanks to Weston Barger and Yu-Chen Cheng for checking and writing solutions to the exercises as well as making corrections and improvements to the text, and to the contributors listed on GitHub.

\subsection*{Donations}
If you find these notes useful and would like to make a donation to help me develop them further you can donate ETH or ERC20 tokens to the address below.\\[0.25em]
\underline{ETH address}: 0x3D6042D6a6f39735bbAd43ad1Aa676113D55BaAC

\mainmatter

%
%

\chapter{Review of probability}
\label{chapter:probability.review}
The notes from this chapter are taken primarily from \cite[Chapter 1]{shreve2004} and \cite[Chapters 1--5]{gs2001}.


\section{Events as sets}
\begin{definition}
The set of all possible outcomes of an experiment is called the \emph{sample space}.  We denote the sample space as $\Om$.
\end{definition}

\noindent
We will typically denote by $\om$ a generic element of $\Om$.

\begin{definition}
An \emph{event} is a subset of the sample space.  We usually denote events by capital roman letters $A, B, C, \ldots$.
\end{definition}

\begin{example}[Toss a coin once]
When one tosses a coin, there are two possible outcomes: heads ($H$) and tails ($T$).  Thus, we have
$\Om = \{ H, T \}$.  One element of $\Om$ is, e.g., $\om = H$.  Possible event: ``toss a head.''  $A = H$.
\end{example}

\begin{example}[Toss two distinguishable coins]
$\Om = \{ (HH), (HT), (TH), (TT) \}$.  One element of $\Om$ is, e.g.,  $\om = (HT)$.  Possible event: `'second toss a tail.'' $A = \{ (HT), (TT) \}$. 
\end{example}

\begin{example}[Toss two indistinguishable coins]
$\Om = \{ \{HH\}, \{HT\}, \{TT\} \}$.  One element of $\Om$ is, e.g.,  $\{HH\}$. Possible event: `'coins match.'' $A = \{ \{HH\},\{TT\} \}$.
\end{example}

\begin{example}[Roll a die]
$\Om = \{ 1, 2, 3, 4, 5, 6 \}$.  One element of $\Om$ is, e.g.,  $\om = 2$.  Possible event: ``roll an odd number.'' $A = \{ 1, 3, 5 \}$.
\end{example}

\noindent
Here, we use $( \cdot )$ to denote an ordered sequence and we use $\{ \cdot \}$ to denote an unordered set.  Thus, $(HT) \neq (TH)$ but $\{HT\}=\{TH\}$.
\\[0.5em]
If $A$ and $B$ are subsets of $\Om$, we can reasonably concern ourselves with events such as ``not $A$'' ($A^c$), ``$A$ or $B$'' ($A \cup B$), ``$A$ and $B$'' ($A \cap B$), etc.  A $\sig$-algebra is a mathematical way to describe all possible sets of interest for a given sample space $\Om$.

\begin{definition}
A collection $\Fc$ of subsets of $\Om$ is called a \emph{$\sig$-algebra} if it satisfies
\begin{enumerate*}
\item contains the empty set: $\emptyset \in \Fc$; \label{itm:empty}
\item is closed under countable unions: $A_1, A_2, A_3, \ldots \in \Fc$ $\Rightarrow$ $\cup_i A_i \in \Fc$; \label{itm:countable.union}
\item is closed under complements: $A \in \Fc$ $\Rightarrow$ $A^c \in \Fc$; \label{itm:complement}
\end{enumerate*}
\end{definition}

\noindent
Note that if $\Fc$ is a $\sig$-algebra then $\Om \in \Fc$ by items \ref{itm:empty} and \ref{itm:complement}.  Note also that $\Fc$ is closed under countable intersections since
\begin{align}
A_1, A_2, A_3, \ldots
	&\in \Fc &
	&\Rightarrow&
A_1^c, A_2^c, A_3^c, \ldots
	&\in \Fc \\ & & &\Rightarrow&
\cup_i A_i^c
	&\in \Fc \\ & & &\Rightarrow&
\Big( \bigcup_i A_i^c \Big)^c
	&\in \Fc \\ & & &\Rightarrow&
\bigcap_i A_i
	&\in \Fc .
\end{align}
Alternatively, one can define a $\sig$-algebra $\Fc$ as a set of subsets of $\Om$ that contains at least the empty set $\emptyset$ and is closed under countable set operations (though, \emph{not} necessarily closed under \emph{un}countable set operators).

\begin{example}[Trivial $\sig$-algebra]
The set of subsets $\Fc_0 := \{\emptyset, \Om \}$ of $\Om$ is commonly referred to as the \emph{trivial} $\sig$-algebra.
\end{example}

\begin{example}
If $A$ is a subset of $\Om$ then $\Fc_A := \{ \emptyset, \Om, A, A^c \}$ is a $\sig$-algebra.
\end{example}

\begin{example}
The \emph{power set} of $\Om$, written $2^\Om$ is the collection of all subsets of $\Om$.  The power set $\Fc = 2^\Om$ is a $\sig$-algebra.
\end{example}

\begin{definition}
Let $\Gc$ be a collection of subsets of $\Om$.  The $\sigma$-algebra \emph{generated} by $\Gc$, written $\sig(\Gc)$, is the smallest $\sig$-algebra that contains $\Gc$.
\end{definition}

\noindent
By ``smallest'' $\sig$-algebra we mean the $\sig$-algebra with the fewest sets.
One can show (although we will not do so in these notes) that $\sig(\Gc)$ is equal to the intersection of all $\sig$-algebras that contain $\Gc$.

\begin{example}
The collection of sets $\Gc = \{ \emptyset, A , \Om \}$ is not a $\sig$-algebra because it does not contain $A^c$.  However, we could create a $\sig$-algebra from $\Gc$ by simply adding the set $A^c$.  Thus, we have $\sig(\Gc) = \{ \emptyset, \Om, A, A^c  \}$.  
\end{example}

\begin{definition}
The pair $(\Om,\Fc)$ where $\Om$ is a sample space and $\Fc$ is a $\sig$-algebra of subsets of $\Om$ is called a \emph{measurable space}.
\end{definition}


\section{Probability}
So far, we have not yet talked about probabilities at all -- only outcomes of a random experiment (elements $\om \in \Om$) and events (subsets $A \subseteq \Om$).  A probability measure assigns probabilities to events.
\begin{definition}
\label{def:P}
A \emph{probability measure} defined on $(\Om,\Fc)$ is a function $\Pb : \Fc \to [0,1]$ that satisfies
\begin{enumerate*}
\item $\Pb( \Om ) = 1$; \label{itm:01}
\item if $A_i \cap A_j = \emptyset$ for $i \neq j$ then $\Pb( \cup_i A_i ) = \sum_i \Pb( A_i )$. (countable additivity) \label{itm:disjoint}
\end{enumerate*}
\end{definition}

\noindent
As we shall see in Section \ref{sec:coin-toss}, it is very important to recognize that Item \ref{itm:disjoint} holds only for \emph{countable} unions.  For an uncountable union it is \emph{not true} that $\Pb( \cup_\alpha A_\alpha ) = \sum_\alpha \Pb( A_\alpha )$.
\\[0,5em]
How can we see the well-known result $\Pb( A^c ) = 1 - \Pb( A )$ from the above definition?  Simply note that
\begin{align}
1
	&=	\Pb(\Om) = \Pb( A \cup A^c ) = \Pb( A ) + \Pb( A^c ) .
\end{align}

\noindent
A probability measure $\Pb$ does \emph{not} need to correspond to empirically observed probabilities!  For example, from experience, we know that if we toss a fair coin we have $\Pb( H ) = \Pb( T ) = 1/2$.  However, we can always define a measure $\Pbt$ that assigns different probabilities $\Pbt( H ) = p$ and $\Pbt( T ) = 1-p$.  As long as $p \in [0,1]$ the measure $\Pbt$ is a probability measure on $(\Om,\Fc)$ where $\Om = \{H,T\}$ and $\Fc = \{ \emptyset, \Om, H, T \}$.
\\[0.5em]
The triple $(\Om,\Fc,\Pb)$ is often referred to as a \emph{probability space} or \emph{probability triple}.  To review, the sample space $\Om$ is the collection of all possible outcomes of an experiment.  The $\sig$-algebra $\Fc$ is all sets of interest of an experiment.  And the probability measure $\Pb$ assigns probabilities to these sets.
\\[0.5em]
When a sample space is countable $\Om=\{\om_1, \om_2, \ldots \}$, we can always take the $\sig$-algebra as the power set $\Fc = 2^\Om$ and construct a probability measure $\Pb$ on $(\Om,\Fc)$ by specifying the probabilities of each individual outcome $\Pb(\om_i) = p_i$.  However, when the sample space $\Om$ is uncountable, choosing an appropriate $\sig$-algebra $\Fc$, and constructing a probability measure $\Pb$ on $(\Om,\Fc)$ is a more delicate procedure.


\section{Infinite probability spaces}

In this section, we construct a probability measure on two uncountable sample spaces.

\subsection{Uniform Lebesgue measure on $(0,1)$}
\label{sec:lebesgue.measure}
In our first example, we will construct a mathematical model for choosing a number at random on the open interval $(0,1)$.
Thus, we take $\Om = (0,1)$.  
Note that this sample space is uncountable.
A particular outcome is, e.g., $\om = 0.5$ or $\om = 1/\sqrt{2}$.
Clearly, the number of possible outcomes is infinite.
We would like to construct a probability $\Pb$ so that all numbers in the interval $(0,1)$ are equally likely to be chosen.
It makes sense, then, to choose
\begin{align}
\Pb( \{ \om : \om \in (a,b) \} )
	&=	\mu( (a,b) )
	:=	b	-a , &
0
	&< a \leq b < 1 , \label{eq:lebesgue.measure}
\end{align}
where we have defined the \emph{Lebesgue measure} $\mu$.  Equation \eqref{eq:lebesgue.measure} tells us how to determine the probability that $\om$ falls within an open interval.  But, in fact, \eqref{eq:lebesgue.measure} tells us more than that.  If $\Pb$ is to be a probability measure, then it must satisfy the countable additivity property given in Definition \ref{def:P}.  Thus, we also know, for example, that
\begin{align}
\Pb( \{ \om : \om \in (a,b) \cup (c,d) \} )
	&=	\Pb( \{ \om : \om \in (a,b) \} ) + \Pb( \{ \om : \om \in  (c,d) \} ) \\
	&=		(b-a) + (d-c) , &
0
	&< a < b < c < d < 1 .
\end{align}
It is natural to ask, \emph{what are all of the sub-sets of $(0,1)$ whose probabilities are uniquely determined by \eqref{eq:lebesgue.measure}}?
Surprisingly, the answer is \emph{not} the power set $2^{(0,1)}$.  It turns out that the power set $2^{(0,1)}$ has sets whose probabilities are \emph{not} determined by \eqref{eq:lebesgue.measure}.  The sets whose probabilities \emph{are} uniquely determined by \eqref{eq:lebesgue.measure} and Definition \ref{def:P} are the sets in $\sig$-algebra generated by the open intervals
\begin{align}
\Bc((0,1))
	&:=	\sig( \Oc ) , &
	&\text{where} &
\Oc
	&:=	\{ A \subseteq (0,1) : A = (a,b) , \,  0 \leq a < b \leq 1 \} .
\end{align}
We call $\Bc((0,1))$ the \emph{Borel $\sig$-algebra on $(0,1)$}.  Thus, the appropriate sample space for our experiment is $(\Om,\Fc,\Pb) = ((0,1),\Bc((0,1)),\mu)$.
\\[0.5em]
We can generalize the notion of Borel sets to all topological spaces.

\begin{definition}
Let $\Om$ be some topological space and let $\Oc(\Om)$ be the set of open sets in $\Om$.  We define we define the \emph{Borel $\sig$-algebra on $\Om$}, denoted $\Bc(\Om)$ by $\Bc(\Om) := \sig(\Oc(\Om))$.
\end{definition}

\begin{remark}
Do not worry too much about what exactly Borel $\sig$-algebras are.  Just think of them as ``reasonable'' sets.  In fact, you would have to think very hard to come up with a set that is not a Borel set.
\end{remark}


\subsection{Infinite sequence of coin tosses}
\label{sec:coin-toss}
In this section we consider and infinite sequence of coin tosses.  We define
\begin{align}
\Om
	&:=	\text{ the set of infinite sequences of  $H$s and $T$s } .
\end{align}
Note that this set is uncountable because there is a one-to-one correspondence between $\Om$ and the set of reals in $[0,1]$.  We will denote a generic element of $\Om$ as follows:
\begin{align}
\om
	&=	\om_1 \om_2 \om_3 \ldots
\end{align}
where $\om_i$ is the result of the $i$th coin toss.  We want to construct a $\sig$-algebra for this experiment.
\\[0.5em]
Let us define some $\sig$-algebras.  First, consider the trivial $\sig$-algebra
\begin{align}
\Fc_0
	&=	\{ \emptyset, \Om \} .
\end{align}
Given no information, I can tell if $\om$ is in the sets in $\Fc_0$ because we know $\om \in \Om$ and $\om \notin \emptyset$.
Next, define two sets
\begin{align}
A_H
	&= \{ \om \in \Om : \om_1 = H \} , &
A_T
	&= \{ \om \in \Om : \om_1 = T \} .
\end{align}
Noting that $A_H = A_T^c$ we see that
\begin{align}
\Fc_1
	&:=	\{ \emptyset, \Om, A_H , A_T \} ,
\end{align}
satisfies the conditions of $\sig$-algebra.  Given $\om_1$ it is possible to say whether or not $\om$ is in each of the sets in $\Fc_1$.  For example, if $\om_1 = H$ then $\om \in A_H$ and $\om \in \Om$, but $\om \notin A_T$ and $\om \notin \emptyset$.
Next define four sets
\begin{align}
A_{HH}
	&:=	\{ \om \in \Om : \om_1 = H , \om_2 = H \} , &
A_{HT}
	&:=	\{ \om \in \Om : \om_1 = H , \om_2 = T \} , \\
A_{TT}
	&:=	\{ \om \in \Om : \om_1 = T , \om_2 = T \} , &
A_{TH}
	&:=	\{ \om \in \Om : \om_1 = T , \om_2 = H \} .
\end{align}
We wish to construct a $\sig$-algebra that contains these sets and the sets in $\Fc_1$.  The smallest such $\sig$-algebra is
\begin{align}
\Fc_2
	&=	\left\{ \begin{aligned} 
			&	\emptyset, \Om, A_H , A_T, A_{HH}, A_{HT}, A_{TT}, A_{TH}, A_{HH}^c , A_{HT}^c, A_{TT}^c, A_{TH}^c \\ 
			&	A_{HH} \cup A_{TH}, A_{HH} \cup A_{TT} , A_{HT} \cup A_{TH} , A_{HT} \cup A_{TT}
			\end{aligned} \right\} .
\end{align}
Given $\om_1$ and $\om_2$, we can say if $\om$ belongs to each of the sets in $\Fc_2$.  Continuing in this way, we can define a $\sig$-algebra $\Fc_n$ for every $n \in \Nb$.  Finally, we take
\begin{align}
\Fc
	&:=	\sig( \Fc_\infty ) , &
\Fc_\infty
	&=	\cup_n \Fc_n	.
\end{align}
One might ask if we could have simply taken $\Fc = \Fc_\infty$?  Well, $\Fc_\infty$ contains every set that can be described in terms of \emph{finitely many} coin tosses.  However, we may be interested in sets such as ``sequences for which $x$ percent of coin tosses are heads,'' and these sets are not in $\Fc_\infty$.  It turns out such sets are in $\Fc$.
\\[0.5em]
Now, we want to construct a probability measure on $\Fc$.  Let us assume the coin tosses are independent (a term we will describe rigorously later on) and that the probability of a head is $p$.  Setting $q = 1-p$, it should be obvious that
\begin{align}
\Pb( \emptyset )
	&=	0 , &
\Pb( \Om ) 
	&=	1 , &
\Pb( A_H )
	&=	p , &
\Pb( A_T )
	&=	q , \\
\Pb( A_{HH} )
	&=	p^2 , &
\Pb( A_{HT} ) 
	&=	p q , &
\Pb( A_{TH} )
	&=	p q , &
\Pb( A_{TT} )
	&=	q^2 , \ldots
\end{align}
Continuing in this way, we can define $\Pb(A)$ for every $A \in \Fc_\infty$.  What about the sets that are in $\Fc$ but not in $\Fc_\infty$?  It turns out that once we have defined $\Pb$ for sets in $\Fc_\infty$ there is only one way to assign probabilities to those sets that are in $\Fc$ but not in $\Fc_\infty$.
We refer the interested reader to \emph{Carath\'eodory's Extension Theorem} for details.
\\[0.5em]
Now, let us define
\begin{align}
A
	&=	\left\{ \om : \lim_{n \to \infty} \frac{ \text{ $\#H$ in first $n$ coin tosses} }{ n } = \frac{1}{2} \right\} .
\end{align}
The strong law of large numbers (SLLN) tells us that $\Pb( A ) = 1$ if $p = 1/2$ and $\Pb( A ) = 0$ if $p \neq 1/2$ (if you have not yet seen the SLLN, you should be able to see this from intuition).  Now it should be clear why \emph{uncountable} additivity does \emph{not} hold for probability measures.  The probability of any given sequence of infinite coin tosses is zero: $\Pb(\om)=0$.  If we were to attempt to compute $\Pb(A)$ by adding up the probabilities $\Pb(\om)$ of all elements $\om \in A$ we would find
\begin{align}
\sum_{\om \in A} \Pb(\om) 
	&= \sum_{\om \in A} 0 = 0
	\neq 1 = \Pb(A), &
	&		\text{(when $p = 1/2$)} .
\end{align}
Thus, uncountable additivity clearly does \emph{not} hold.
\\[0.5em]
We finish this example (we will come back to it!) with the following definition
\begin{definition}
Let $(\Om,\Fc,\Pb)$ be a probability space. If a set  $A \in \Fc$ satisfies $\Pb(A)  = 1$, we say that the event $A$ occurs \emph{$\Pb$ almost surely} (written, $\Pb$-a.s.).
\end{definition}
\noindent
Note in the example above that, when $p = 1/2$ we have $\Pb(A)=1$ and thus $A$ occurs almost surely.  But it is important to recognize that $A \neq \Om$ and $A^c \neq \emptyset$.  The elements of $A^c$ are part of the sample space $\Om$, but they have zero probability of occurring.


\section{Random variables and distributions}
A random variable maps the outcome of an experiment to $\Rb$.  We capture this idea with the following definition.

%

\begin{definition}
\label{def:rv}
A \emph{random variable} defined on $(\Om,\Fc)$ is a function $X: \Om \to \Rb$ with the property that
\begin{align}
\{ X \in A \}
	&:=	\{ \om \in \Om : X(\om) \in A \} 
	\in \Fc , \label{eq:rv}
\end{align}
for all $A \in \Bc(\Rb)$.
\end{definition}

\noindent
Observe all random variables must be defined on a measurable space $(\Om,\Fc)$, as these appear in the definition.  Note, however, that the probability measure $\Pb$ does \emph{not} appear in the definition.  Random variables are defined \emph{independent} of a probability measure $\Pb$.
\\[0.5em]
What does Definition \ref{def:rv} mean?  Recall that a probability measure $\Pb$ defined on $(\Om,\Fc)$ maps $\Fc \to [0,1]$.  In order for us to answer the question: ``what is the probability that $X \in A$?'' we need for the set $\{ X \in A \} \in \Fc$.  And this is precisely what Definition \ref{def:rv} requires.  Why do we only consider sets $A \in \Bc(\Rb)$ rather than any set $A \subset \Rb$?  The answer is rather technical and, frankly, not worth exploring at the moment.
\\[0.5em]
\textbf{A word on notation}: the standard convention is to use capital Roman letters (typically, $X$, $Y$, $Z$) for random variables and lower case Roman letters ($x$, $y$, $z$) for real numbers.
\\[0.5em]
Let us look at some random variables.

\begin{example}[Discrete time model for stock prices]
\label{ex:stock}
Consider the infinite sequence of coin tosses in Section \ref{sec:coin-toss}.  We Define a sequence of random variables $(S_n)_{n \geq 0}$ via
\begin{align}
S_0(\om)
	&=	1 , &
S_{n+1}(\om)
	&=	\left\{ \begin{aligned}
			&	u S_n & \text{if $\om_{n+1} = H$}\\
			&	d S_n & \text{if $\om_{n+1} = T$}
			\end{aligned} \right. \label{eq:construction}
\end{align} 
Here, $S_n$ represents the value of a stock at time $n$.
Note that $\Pb( S_1 = u ) = \Pb( A_H ) = p$.  Likewise $\Pb( S_2 = u d ) = \Pb( A_{HT} \cup A_{TH} ) = 2 p q$.  More generally, one can show that
\begin{align}
\Pb( S_n = u^k d^{n-k} )
	&=	\binom{n}{k} p^k q^{n-k} . \label{eq:Sn.P}
\end{align}
Note if we had simply defined the random variables $(S_n)_{n \geq 1}$ as having probabilities given by \eqref{eq:Sn.P} we would have no information about how, e.g., $S_n$ relates to $S_{n-1}$.  From the above construction \eqref{eq:construction}, however, we know that if $S_n = u^n$ then $S_{n-1} = u^{n-1}$.  Thus, the structure of a given probability space, not just the probabilities of events, is very important.
\end{example}

\begin{example}
\label{ex:XY}
Let $(\Om, \Fc) = ((0,1),\Bc((0,1)))$.  Define random variables $X(\om) = \om$ and $Y(\om) = 1 - \om$.  Clearly, we have $X = 1-Y$.  Now, suppose we define $\Pb(\dd \om):= \Pb([\om,\om + \dd \om]) = \dd \om$ (i.e., the Lebesgue measure).  Then $X$ and $Y$ have the same distribution.  For $x \in [0,1]$ we have
\begin{align}
\Pb(X \leq x)
	&= \Pb(\om \leq x) = \int_0^x \Pb(\dd \om) = \int_0^x \dd \om = x , \\
\Pb(Y \leq x)
	&= \Pb(1-\om \leq x) = \int_{1-x}^1 \Pb(\dd \om) = \int_{1-x}^1 \dd \om = x .
\end{align}
However, if we define a new probability measure via $\Pbt( \dd \om ) = 2 \om \, \dd \om$ then $X$ and $Y$ have different distributions.  For $x \in [0,1]$ we have
\begin{align}
\Pbt(X \leq x)
	&= \Pbt(\om \leq x) = \int_0^x \Pbt(\dd \om) = \int_0^x 2 \om \dd \om = x^2 , \\
\Pbt(Y \leq x)
	&= \Pbt(1-\om \leq x) = \int_{1-x}^1 \Pbt(\dd \om) = \int_{1-x}^1 2 \om \dd \om = 1-(1-x)^2 .
\end{align}
Thus, the statistical properties of random variables depend not only on the structure of the probability space on which they live, but also on the probability measure.
\end{example}

\noindent
The distribution of a random variable $X$ is most easily described through its cumulative distribution function.

\begin{definition}
The \emph{distribution function} $F_X: \Rb \to [0,1]$ of a random variable $X$ defined on a probability space $(\Om,\Fc,\Pb)$ is given by
\begin{align}
F_X(x)
	&:=	\Pb( X \leq x ) .
\end{align}
\end{definition}

\noindent
Observe that, while a random variable $X$ is defined with respect to $(\Om,\Fc)$ (with no reference to $\Pb)$, the distribution $F_X$ is specific to a probability measure $\Pb$.
\\[0.5em]
Note that we put the random variable $X$ in the subscript of $F_X$ to remind us that $F_X$ is the distribution function corresponding to the random variable $X$ (and not, e.g., $Y$).  It is a good idea to do this.
\\[0.5em]
We give here some obvious properties of $F_X$.
\begin{enumerate*}
\item We have the following limits: $\lim_{x \to - \infty} F_X(x) = 0$ and $\lim_{x \to \infty} F_X(x) = 1$.
\item $F_X$ is non-decreasing: $x < y$ implies $F_X(x) \leq F_X(y)$.
\item $F_X$ is right-continuous and has left-limits: $F_X(x+) := \lim_{h \searrow 0} F_X(x+h) = F_X(x)$.
\end{enumerate*}

\noindent
Many (but not all) random variables fall in to  one of two categories: discrete and continuous.  We describe these two categories below.

\begin{definition}
A random variable $X$ is called \emph{discrete} if it takes values in some countable set $A := \{ x_1, x_2,\ldots \} \subset \Rb$.  We associate with  a discrete random variable a \emph{probability mass function} $f_X : A \to \Rb$, defined by $f_X(x_i) := \Pb( X = x_i )$. 
\end{definition}

\begin{definition}
A random variable $X$ is called \emph{continuous} if its distribution function $F_X$ can be written as
\begin{align}
F_X(x)
	&=	\int_{-\infty}^x \dd u \, f_X(u) , &
x
	&\in \Rb ,
\end{align}
for some $f_X : \Rb \to [0,\infty)$ called the \emph{probability density function}.
\end{definition}

\noindent
It may help to think of the density function $f_X$ as $f_X(x)\dd x = \Pb( X \in \dd x )$.
\\[0.5em]
Note that for a continuous random variable $X$ we have $f_X = F_X'$.
\\[0.5em]
If $X$ is either discrete or continuous, it is easy to compute $\Pb(X \in A)$ for any $A \in \Bc(\Rb)$.  We have
\begin{align}
\text{discrete}:&&
\Pb(X \in A)
	&=	\sum_{\{i : x_i \in A\}} f_X(x_i) , \\
\text{continuous}:&&
\Pb(X \in A)
	&=	\int_A \dd x \, f_X(x) .
\end{align}

\begin{remark}
Although we have defined $F_X : \Rb \to [0,1]$ by $F_X(x) := \Pb( X \leq x )$, it is common to also to utilize $F_X$ as a set function $F_X: \Bc(\Rb) \to [0,1]$, which means $F_X(B) := \Pb( X \in B)$.
In particular, for a continuous random variable $X$, we have
\begin{align}
F_X(\dd x)
	&:= F_X([x,x+\dd x]) = F_X(x + \dd x) - F_X(x) \\
	&=	F_X(x) + F_X'(x) \dd x + \Oc(\dd x^2) - F_X(x)
	=		f_X(x) \dd x , &
	&\text{as $\dd x \to 0$} .
\end{align}
It should always be clear from the argument of $F_X$, which of the multiple meanings we intend.
\end{remark}

\subsection*{Examples of discrete random variables}
The following discrete random variables frequently arise in applications in nature and social sciences.

\begin{example}
\label{ex:bernoulli}
If $X$ is distributed as a \emph{Bernoulli random variable} with parameter $p \in [0,1]$, written $X \sim \text{Ber}(p)$, then
\begin{align}
X
	&\in \{ 0 , 1 \} , &
f_X(k)
	&=	\left\{ \begin{aligned}
			& 1-p & k=0 , \\
			& p & k=1 .
			\end{aligned} \right. 
\end{align}
\end{example}

\begin{example}
\label{ex:binomial}
If $X$ is distributed as a \emph{Binomial random variable} with parameters $n \in \Nb$ and $p \in [0,1]$, written $X \sim \text{Bin}(n,p)$, then
\begin{align}
X
	&\in \{0,1,2,\ldots,n\} , &
f_X(k)
	&=	\binom{n}{k} p^k (1-p)^{n-k} .
\end{align}
Note that if $X_i \sim \text{Ber}(p)$ and independent of each other then $Y := \sum_{i=1}^n X_i \sim \text{Bin}(n,p)$.
\end{example}

\begin{example}
\label{ex:geometric}
If $X$ is distributed as a \emph{Geometric random variable} with parameter $p \in [0,1]$, written $X \sim \text{Geo}(p)$, then
\begin{align}
X
	&\in \Nb , &
f_X(k)
	&=	p (1-p)^{k-1} .
\end{align}
Note that $X_i \sim \text{Ber}(p)$ $i=1,2,\ldots$ are independent of each other then $Y := \inf\{ i : X_i = 1 \} \sim \text{Geo}(p)$.
\end{example}

\begin{example}
\label{ex:poisson}
If $X$ is distributed as a \emph{Poisson random variable} with parameter $\lam \in \Rb_+$, written $X \sim \text{Poi}(\lam)$, then
\begin{align}
X
	&\in \{ 0 \} \cup \Nb , &
f_X(k)
	&=	\frac{\lam^k}{k!}\ee^{-\lam} .
\end{align}
\end{example}

\subsection*{Examples of continuous random variables}

Before introducing some common continuous random variables, let us introduce a useful function.

\begin{definition}
Let $A$ be a set in some topological space $\Om$ (e.g., $\Om = \Rb^d$).  The \emph{indicator function} $\Ib_A : \Om \to \{0,1\}$ is defined as follows
\begin{align}
\Ib_A(x)
	&:= \left\{ \begin{aligned}
			&1 & \text{if $x \in A$}, \\
			&0 & \text{if $x \notin A$}.
			\end{aligned} \right.
\end{align}
\end{definition}

\noindent
\textbf{Notation}: We will sometimes write $\Ib_A(x) = \Ib_{\{ x \in A \}}$.
\\[0.5em]
We now introduce some continuous random variables that frequently arise in applications.

\begin{example}
\label{ex:uniform}
If $X$ is distributed as a \emph{Uniform random variable} on the interval $[a,b] \subset \Rb$, written $X \sim \Uc([a,b])$, then
\begin{align}
X
	&\in [a,b]  , &
f_X(x)
	&=	\Ib_{[a,b]}(x) \frac{1}{b-a}.
\end{align}
\end{example}

\begin{example}
\label{ex:exponential}
If $X$ is distributed as a \emph{Exponential random variable} with mean $\lam>0$, written $X \sim \Ec(\lam)$, then
\begin{align}
X
	&\in [0,\infty) , &
f_X(x)
	&=	\Ib_{[0,\infty)}(x) \lam \ee^{-\lam x} .
\end{align}
Note that $f_X(x) = 0$ if $x<0$ due to the presence of the indicator function.
\end{example}

\begin{example}
\label{ex:normal}
If $X$ is distributed as a \emph{Gaussian} or \emph{Normal random variable} with mean $\mu \in \Rb$ and variance $\sig^2>0$ (we will give a meaning for ``mean'' and ``variance'' below), written $X \sim \Nc(\mu,\sig^2)$, then
\begin{align}
X
	&\in \Rb , &
f_X(x)
	&=	\frac{1}{\sqrt{2 \pi \sig^2}}\exp\( - \frac{(x-\mu)^2}{2 \sig^2} \) .
\end{align}
A random variable $Z \sim \Nc(0,1)$ is referred to as \emph{standard normal}.
\end{example}


%
%


\section{Stochastic Processes}

Intuitively, we think of a stochastic process as a process that evolves randomly in time.  Now that we understand what a random variable is, we can define rigorously what we mean when we say \textit{stochastic process}.

\begin{definition}
A \textit{Stochastic process} is a collection of random variables $X = (X_t)_{t \in \mathds{T}}$ where $\mathds{T}$ is some index set.
If the index set $\mathds{T}$ is countable (e.g., $\mathds{T} = \Nb_0$) we say that $X$ is a \textit{discrete time} process.
If the index set $\mathds{T}$ is uncountable (e.g., $\mathds{T} = \Rb_+$) we say that $X$ is a \textit{continuous time} process.
The \textit{State Space} $S$ of a stochastic process $X$ is union of the state spaces of $(X_t)_{t \in \mathds{T}}$.
\end{definition}

\noindent
We can think of a stochastic process $X : \mathds{T} \times \Om \to \Rb$ in (at least) two ways.  First, for any $t \in \mathds{T}$ we have that $X_t: \Om \to \Rb$ is random variable.  Second, for any $\om \in \Om$, we have that $X_\cdot(\om) : \mathds{T} \to \Rb$ is a function of time.  Both interpretations can be useful.


\section{Expectation}

When we think of averaging we think of weighting outcomes by their probabilities.  The mathematical way to encode this is via the expectation.

\begin{definition}
\label{def:E}
Let $X$ be a random variable defined on $(\Om,\Fc,\Pb)$.  The \emph{expectation} of $X$, written $\Eb X$, is defined as
\begin{align}
\Eb X
	&:=	\int_\Om X(\om) \Pb(\dd \om) ,
\end{align}
where the integral is understood in the Lebesgue sense.
\end{definition}

\subsection{Integration in the Lebesgue sense}
For those who have not previously encountered Lebesgue integration, we now give a brief (\emph{very brief}!) overview of this concept.

\begin{definition}
Fix a probability space $(\Om, \Fc, \Pb)$.  Let $A \in \Fc$.
The \emph{indicator random variable},
denoted $\Ib_A$, is defined by
\begin{align}
\Ib_A(\om)
	&:=	\left\{ \begin{aligned}
			&	1	& \om \in A , \\
			& 0 &	\om \notin A .
			\end{aligned} \right.
\end{align} 
\end{definition}

\noindent
Observe that $\Ib_A \sim \text{Ber}(p)$ with $p = \Pb(A)$.
For disjoint sets $A$ and $B$ we have
\begin{align}
\Ib_{A \cup B}
	&=	\Ib_A + \Ib_B , &
A \cap B
	&= \emptyset .
\end{align}
And, for any two sets $A$ and $B$ we have
\begin{align}
\Ib_{A \cap B}
	&=	\Ib_A \Ib_B .
\end{align}

\begin{definition}
\label{def:partition}
A collection of non-empty sets $(A_i)$ is said to be a \textit{partition} of $\Om$ if $A_i \cap A_j = \emptyset$ for all $i \neq j$ and $\cup_i A_i = \Om$. 
\end{definition}

\begin{definition}
Let $(A_i)$ be a finite partition of $\Om$.
A non-negative random variable $X$, defined on a probability space $(\Om,\Fc,\Pb)$, which is of the form
\begin{align}
X(\om)
	&=	\sum_{i=1}^n x_i \Ib_{A_i}(\om) , &
x_i
	&\geq 0 , &
A_i
	&\in \Fc , \label{eq:X.simple}
\end{align}
is called \emph{simple}. 
\end{definition}

\noindent
Let $X$ be a simple random variable.  We define the expectation of $X$ as follows
\begin{align}
\Eb X
	&:=	\sum_{i=1}^n x_i \Pb(A_i) . &
	&\text{(if $X$ is simple)}
\end{align}
Note that, from this definition, we have
\begin{align}
\Eb \Ib_A
	&=	\Pb(A) .
\end{align}
Thus, \emph{we can always represent probabilities of sets as expectations of indicator random variables}.
\\[0.5em]
Now, consider a non-negative random variable $X$, which is not necessarily simple.  Let $(X_n)_{n \geq 0}$ be an increasing sequence of simple random variables that converges almost surely to $X$.  That is
\begin{align}
X_i
	&\leq X_{i+1} , &
\lim_{i \to \infty} X_i 
	&\to X , &
	&\text{$\Pb$-a.s.}.
\end{align}
We define the expectation of a non-negative random variable $X$ as the following limit
\begin{align}
\Eb X
	&:=	\lim_{i \to \infty} \Eb X_i , &
	&\text{(if $X$ is non-negative)} \label{eq:E.nonneg}
\end{align}
where each of the expectations on the right-hand side are well-defined because all of the $X_i$ are simple by construction.
Finally, consider a general random variable $X$ that could take either positive or negative values.  Define
\begin{align}
X^+
	&=	\max\{ X, 0 \} , &
X^-
	&=	\max\{ -X, 0 \} 	.
\end{align}
Note that $X^+$ and $X^-$ are non-negative and $X = X^+ - X^-$.  With this in mind, we define
\begin{align}
\Eb X
	&:=	\Eb X^+ - \Eb X^- ,
\end{align}
where the expectations of $X^+$ and $X^-$ are defined via \eqref{eq:E.nonneg}.
\\[0.5em]
Definition \ref{def:E} of $\Eb X$ makes sense if $\Eb|X|<\infty$ or if $\Eb X^\pm = \infty$ and $\Eb X^\mp < \infty$.  In the latter case, we have $\Eb X = \pm \infty$.
If both $\Eb X^+ = \infty$ and $\Eb X^- = \infty$, then we find ourselves in an $\infty - \infty$ situation and, in this case, $\Eb X$ is undefined.

\subsection{Computing expectations}
If $X$ is either discrete or continuous Definition \ref{def:E} reduces to the formulas one learns as an undergraduate.
\begin{align}
\text{discrete}:&&
\Eb X
	&=	\sum_i x_i \, f_X(x_i) , \\
\text{continuous}:&&
\Eb X
	&=	\int_\Rb \dd x \, x \,  f_X(x) .
\end{align}
In the discrete case, the sum runs over all possible values of $x$.  More generally, we can express the expected value of $X$ as
\begin{align}
\Eb X
	&=	\int_\Rb x \, F_X(\dd x) 
	:=	\lim_{\norm{\Pi} \to 0} \sum_i \Big( \frac{x_i + x_{i+1}}{2} \Big) \Big( F_X(x_{i+1}) - F_X(x_i) \Big) , \label{eq:stieltjes} 
\end{align}
where $\Pi$ is a partition of $\Rb$, meaning
\begin{align}
\Pi
	&=	\{ x_i, x_2, \ldots , x_n \} , &
x_i
	&< x_{i+1} \, \forall i , &
\norm{\Pi}
	&:= \sup_i ( x_{i+1} - x_i ) .
\end{align}
The expression on the right-hand side of \eqref{eq:stieltjes} is known as a \emph{Stieltjes integral}.
The advantage of using the Stieltjes integral $\int x F_X(\dd x)$ to compute an expectation is that every random variable $X$ -- whether it be discrete or continuous -- has a distribution $F_X$.  Thus, by using the Stieltjes integral, we avoid having to treat discrete and continuous cases separately.
\\[0.5em]
Note that $\Eb$ is a linear operator.  If $X$ and $Y$ are random variables and $a$ and $b$ are constants, then
\begin{align}
\Eb( a X + b Y )
	&=	a \Eb X + b \Eb Y .
\end{align}

\noindent
How does one compute $\Eb g(X)$ where $g:\Rb \to \Rb$?  Although we have not stated it explicitly, it should be obvious that if $X$ is a random variable, then $Y:=g(X)$ is also a random variable.
\footnote{Rigorously, $g$ should be a \emph{measurable function}, meaning $g^{-1}(A) \in \Bc(\Rb)$ for all $A \in \Bc(\Rb)$.  Do not concern yourself too much with this.}
Thus, we have
\begin{align}
\Eb Y
	&=	\Eb g(X)
	=		\int_\Om g(X(\om)) \Pb( \dd \om ) ,
\end{align}
which in the discrete and continuous cases become
\begin{align}
\text{discrete}:&&
\Eb g(X)
	&=	\sum_i g(x_i) f_X(x_i) , \\
\text{continuous}:&&
\Eb g(X)
	&=	\int_\Rb \dd x\, g(x) f_X(x) .
\end{align}




\section{Change of measure}
\label{sec:measure.change}
Consider two probability measures $\Pb$ and $\Pbt$ defined on a measurable space $(\Om,\Fc)$.  What is the relation between $\Pb$ and $\Pbt$?  The following theorem answers this question.

\begin{theorem}
\label{thm:measure.change}
Fix a probability space $(\Om,\Fc,\Pb)$ and let $Z$ be a random variable satisfying $Z \geq 0$ $\Pb$-a.s. and $\Eb Z = 1$.  
Define a $\Pbt : \Fc \to [0,1]$ by
\begin{align}
\Pbt(A)
	&:=	\Eb Z \Ib_A . \label{eq:PtA=EZIA}
\end{align}
Then $\Pbt$ is a probability measure on $(\Om,\Fc)$.  Denote by $\Ebt$ the expectation taken with respect to $\Pbt$.  Then
\begin{align}
\Ebt X
	&=	\Eb Z X , &
	&\text{and if}&
Z
	&>	0 , &
	&\text{then}&
\Eb X
	&=	\Ebt \frac{1}{Z} X . \label{eq:EtX=EZX}
\end{align}
where $X$ is a random variable defined on $(\Om,\Fc)$.
\end{theorem}
\begin{definition}
We call the random variable $Z$ in Theorem \ref{thm:measure.change} the \emph{Radon-Nikod\'ym derivative} of $\Pbt$ with respect to $\Pb$.
\end{definition}
\begin{proof}[Proof of Theorem \ref{thm:measure.change}]
First, we show that $\Pbt$ satisfies the properties of a probability measure given in Definition \ref{def:P}.  Let us check that $\Pbt(\Om) = 1$.  We have
\begin{align}
\Pbt(\Om)
	&=	\Eb Z \Ib_\Om
	=		\Eb Z
	=		1.
\end{align}
Next, we check that $\Pbt$ satisfies countable additivity.  Let $(A_i)_{i \geq 0}$ be a sequence of disjoint sets.  We have
\begin{align}
\Pbt(\cup_i A_i)
	&=	\Eb \Ib_{\cup_i A_i} Z
	=		\sum_i \Eb \Ib_{A_i} Z
	=		\sum_i \Pbt( A_i ) .
\end{align}
Note that interchanging the sum with the expectation $\Eb \sum_i \to \sum_i \Eb$ is allowed by Tonelli's Theorem.  Finally, to show equation \eqref{eq:EtX=EZX} holds, it is enough to check that it holds for simple random variables $X = \sum_i x_i \Ib_{A_i}$.  We have
\begin{align}
\Ebt X
	&=	\Ebt \sum_i x_i \Ib_{A_i}
	=		\sum_i x_i \Eb Z \Ib_{A_i}
	=		\sum_i x_i \Pbt(A_i) , 
\end{align}
which agrees with the definition of Expectation for simple random variables.  Finally, if $Z>0$, we have
\begin{align}
\Ebt \frac{1}{Z} X
	=		\Eb Z \frac{1}{Z} X
	=		\Eb X .
\end{align}
\end{proof}

\begin{definition}
A probability measure $\Pbt$ defined on $(\Om,\Fc)$ is \emph{absolutely continuous} with respect to another probability measure $\Pb$ (equivalently, $\Pb$ \emph{dominates} $\Pbt$), written $\Pb \gg \Pbt$, if
\begin{align}
\Pb(A) 
	&= 0 &
	&\Rightarrow& 
\Pbt(A) 
	&= 0.
\end{align}
\end{definition}

\noindent
In order for there to exist a Radon-Nikod\'ym derivative $Z = \dd \Pbt / \dd \Pb$ we must have $\Pb \gg \Pbt$.
The reason for this is that the relationship between $\Pb$ and $\Pbt$ is multiplicative $\Pbt(\dd \om) = Z(\om) \Pb(\dd \om)$.
If $\Pb(A) = 0$, then there is no random variable $Z$ that would result in $\Pbt(A) > 0$.
To see this more formally, note that
\begin{align}
\Pb(A)
	&=	\Eb \Ib_A , &
\Pbt(A)
	&=	\Eb Z \Ib_A . \label{eq:probabilities}
\end{align}
Now, suppose that $Z$ exists such that \eqref{eq:probabilities} holds and $\Pb \not\gg \Pbt$.
Then there exists some set $A$ such that $\Pb(A) = 0$ and $\Pbt(A) > 0$.
But, this is impossible because if $\Pb(A) = 0$ then $Z \Ib_A = 0$ $\Pb$-a.s., and thus $\Pbt(A) = \Eb Z \Ib_A = 0$.
Thus, if \eqref{eq:probabilities} holds, we must have $\Pb \gg \Pbt$.

\begin{definition}
Two probability measures $\Pb$ and $\Pbt$ on $(\Om,\Fc)$ are \emph{equivalent}, written $\Pb \sim \Pbt$, if 
\begin{align}
\Pb(A) 
	&= 0 &
	&\Leftrightarrow& 
\Pbt(A) 
	&= 0 .
\end{align}
\end{definition}

\noindent
Two probability measures are equivalent $\Pb \sim \Pbt$ if and only if the Radon-Nikod\'ym Derivative that relates them is strictly positive $Z > 0$.
Equivalent measures agree on which events will happen with probability zero (and thus, they agree on which events will happen with probability one).

\begin{example}
Let us return to Example \ref{ex:XY}. We set $(\Om, \Fc) = ((0,1),\Bc((0,1)))$.  On this measure space, we define two probability measures $\Pb(\dd \om) = \dd \om$ and $\Pbt(\dd \om) = 2 \om \dd \om$.  Note that we have
\begin{align}
\Pbt(A)
	&=	\Ebt \Ib_A
	=		\int_\Om \Ib_A(\om) \Pbt(\dd \om)
	=		\int_\Om \Ib_A(\om) 2 \om \dd \om 
	=		\int_\Om \Ib_A 2 \om \Pb(\dd \om)
	=		\Eb \Ib_A Z , &
Z(\om)
	&:=	2 \om .
\end{align}
One can easily check that $\Eb Z = 1$ and $Z>0$.  Defining $\Pbt$ by \eqref{eq:PtA=EZIA}, one can easily check that \eqref{eq:EtX=EZX} holds true.
\end{example}


\noindent
It is quite common to use the notation
\begin{align}
Z(\om)
	&=	\frac{\dd \Pbt}{\dd \Pb}(\om) , &
\Pbt(\dd \om)
	&=	\frac{\dd \Pbt}{\dd \Pb}(\om) \Pb(\dd \om) ,
\end{align}
as a reminder of how the Radon-Nikod\'ym Derivative $Z$ relates $\Pbt$ to $\Pb$.
For a finite probability space, it is true that $Z(\om) = \Pbt(\om)/\Pb(\om)$.  
However, for an infinite probability space, it makes no sense in general to define $Z(\om)=\Pbt(\om)/\Pb(\om)$ since it may be that $\Pb(\om) = 0$.  Nevertheless, the heuristic $Z(\om) = \Pbt(\om)/\Pb(\om)$ gives the correct intuition.  In particular, for the special case of an infinite probability space in which $\Pb(\dd \om) = p(\om) \dd \om$ and $\Pbt(\dd \om) = \pt(\om) \dd \om$ and $\Pb \gg \Pbt$, we have $Z(\om) = \pt(\om) / p(\om)$.  

\begin{example}[Change of measure Normal random variable]
\label{ex:normal.meas.change}
On $(\Om, \Fc, \Pb)$ let $X \sim \Nc(0,1)$ and define $Y = X + \theta$.  Clearly, we have $Y \sim \Nc(\theta,1)$.
Now, define a random variable $Z$ by
\begin{align}
Z
	&=	\ee^{ - \theta X - \tfrac{1}{2} \theta^2 } .
\end{align} 
Clearly $Z > 0$.  We also have $\Eb Z = 1$.  To see this, simply compute
\begin{align}
\Eb Z
	&=	\int_\Rb \dd x \, \ee^{ - \theta x - \tfrac{1}{2} \theta^2 } f_X(x) \\
	&=	\int_\Rb \dd x \, \ee^{ - \theta x - \tfrac{1}{2} \theta^2 } \frac{1}{\sqrt{2\pi}}\ee^{-x^2/2} \\
	&=	\int_\Rb \dd x \,  \frac{1}{\sqrt{2\pi}}\ee^{-(x+\theta)^2/2} = 1 .
\end{align}
Since $Z > 0$ and $\Eb Z = 1$, we can define a new probability measure $\Pbt$ with $Z = \dd \Pbt / \dd \Pb$ as the Radon-Nikod\'ym derivative.  Let us compute the distribution of $Y$ under $\Pbt$.  We have
\begin{align}
\Pbt( Y \leq b)
	&=	\Eb Z \Ib_{\{ Y \leq b \}} 
	=		\Eb \ee^{- \theta X - \tfrac{1}{2} \theta^2 } \Ib_{\{X \leq b - \theta \}} \\
	&=	\int_{-\infty}^{b - \theta} \dd x \,  \ee^{ - \theta x - \tfrac{1}{2} \theta^2 } f_X(x) \\
	&=	\int_{-\infty}^{b - \theta} \dd x \, \frac{1}{\sqrt{2\pi}}\ee^{-(x+\theta)^2/2} \\
	&=	\int_{-\infty}^{b} \dd z \, \frac{1}{\sqrt{2\pi}}\ee^{-z^2/2}
\end{align}
Thus, under $\Pbt$ we see that $Y \sim \Nc(0,1)$.  The Radon-Nikod\'ym derivative $Z$ changes the mean of $Y$ from $\theta$ to $0$, but it does not affect the variance of $Y$.
\end{example}

\begin{example}[Change of measure Exponential random variable]
On $(\Om, \Fc, \Pb)$, let $X \sim \Ec(\lam)$ and define
\begin{align}
Z
	&=	\frac{\mu}{\lam} \ee^{-(\mu - \lam)X} .
\end{align}
Clearly, we have $Z \geq 0$.  We also have
\begin{align}
\Eb Z
	&=	\Eb \frac{\mu}{\lam} \ee^{-(\mu - \lam)X}
	=		\int_0^\infty \dd x \, \frac{\mu}{\lam} \ee^{-(\mu - \lam)x} f_X(x) \\
	&=	\int_0^\infty \dd x \, \frac{\mu}{\lam} \ee^{-(\mu - \lam)x} \lam \ee^{- \lam x}
	=		\int_0^\infty \dd x \, \mu \ee^{-\mu x} = 1 .
\end{align}
Thus, we can define a new probability measure $\Pbt$ with $Z = \dd \Pbt / \dd \Pb$ as the Radon-Nikod\'ym derivative.
Let us compute the distribution of $X$ under $\Pbt$.  We have
\begin{align}
\Pbt(X \leq b)
	&=	\Eb Z \Ib_{\{X \leq b\}}
	=		\Eb \frac{\mu}{\lam} \ee^{-(\mu - \lam)X} \Ib_{\{X \leq b\}} \\
	&=		\int_0^b \frac{\mu}{\lam} \ee^{-(\mu - \lam)x}f_X(x)
	=		\int_0^b \dd x \, \frac{\mu}{\lam} \ee^{-(\mu - \lam)x} \lam \ee^{- \lam x}
	=		\int_0^b \dd x \, \mu \ee^{-\mu x}
\end{align}
Thus, under $\Pbt$, we have $X \sim \Ec(\mu)$.
\end{example}

%
%


\section{Exercises}
\label{sec:exercises-1}

\begin{exercise}
\label{exer:1.1}
Let $\Fc$ be a $\sig$-algebra of $\Om$.  Suppose $B \in \Fc$.  Show that $\Gc := \{ A \cap B : A \in \Fc \}$ is a $\sig$-algebra of $B$.
\end{exercise}

\begin{exercise}
\label{exer:1.2}
Let $\Fc$ and $\Gc$ be $\sig$-algebras of $\Om$.  
(a) Show that $\Fc \cap \Gc$ is a $\sig$-algebra of $\Om$.
(b) Show that $\Fc \cup \Gc$ is not necessarily a $\sig$-algebra of $\Om$.
\end{exercise}

\begin{exercise}
\label{exer:1.3}
Describe the probability space $(\Om,\Fc,\Pb)$ for the following three experiments:
(a) a biased coin is tossed three times;
(b) two balls are drawn without replacement from an urn which originally contained two blue and two red balls;
(c) a biased coin is tossed repeatedly until a head turns up.
\end{exercise}

\begin{exercise}
\label{exer:1.4}
Suppose $X$ is a continuous random variable with distribution $F_X$.  Let $g$ be a strictly increasing continuous function.  Define $Y = g(X)$. (a) What if $F_Y$, the distribution of $Y$? (b) What is $f_Y$, the density of $Y$?
\end{exercise}

\begin{exercise}
\label{exer:1.5}
Suppose $X$ is a continuous random variable with distribution $F_X$.  Find $F_Y$ where $Y$ is given by
(a) $X^2$
(b) $\sqrt{|X|}$
(c) $\sin X$
(d) $F_X(X)$.
\end{exercise}

\begin{exercise}
\label{exer:1.6}
Suppose $X$ is a continuous random variable defined on a probability space $(\Om,\Fc,\Pb)$.  Let $f$ be the density of $X$ under $\Pb$ and assume $f>0$.  Let $g$ be the density function of a random variable.  Define $Z := g(X)/f(X)$.  (a) Show that $Z \equiv \dd \Pbt/\dd \Pb$ defines a Radon-Nikod\'ym derivative.  (b) What is the density of $X$ under $\Pbt$?
\end{exercise}

\begin{exercise}
\label{exer:1.7}
Let $X$ be uniformly distributed on $[0,1]$.  For what function $g$ is the random variable $g(X)$ exponentially distributed with parameter 1 (i.e. $g(X) \sim \Ec(1)$)?
\end{exercise}

%
%

\chapter{Information and conditioning}
\label{chapter:conditioning}
The notes from this chapter are taken primarily from \cite[Chapter 2]{shreve2004}.


\section{Information and $\sig$-algebras}
Let us return to the coin-toss example of Section \ref{sec:coin-toss}.  If we are given no information about $\om$ what can we say about $\om$?  In other words, what are the subsets of $\Om$ for which we can say: ``$\om$ is in this set'' or ``$\om$ is not in this set''?  The answer is $\emptyset$ and $\Om$, which, together, form the trivial $\sig$-algebra $\Fc_0 = \{\emptyset,\Om\}$.
\\[0.5em]
Now suppose we are given the value of $\om_1$.  What are the subsets of $\Om$ for which we can say: ``$\om$ is in this set'' or ``$\om$ is not in this set''?  The answer is the sets in $\Fc_0$ as well as $A_H$ and $A_T$.  Together, these sets form the $\sig$-algebra $\Fc_1 = \{ \emptyset, \Om, A_H, A_T \}$.  We say the sets in $\Fc_1$ are \emph{resolved} by the first coin toss.
\\[0.5em]
Now suppose we are given the value of $\om_1$ and $\om_2$.  What are the subsets of $\Om$ for which we can say: ``$\om$ is in this set'' or ``$\om$ is not in this set''?  The answer is the sets in $\Fc_2$, given by
\begin{align}
\Fc_2
	&=	\left\{ \begin{aligned} 
			&	\emptyset, \Om, A_H , A_T, A_{HH}, A_{HT}, A_{TT}, A_{TH}, A_{HH}^c , A_{HT}^c, A_{TT}^c, A_{TH}^c \\ 
			&	A_{HH} \cup A_{TH}, A_{HH} \cup A_{TT} , A_{HT} \cup A_{TH} , A_{HT} \cup A_{TT}
			\end{aligned} \right\} .
\end{align}
The sets in $\Fc_2$ are resolved by the first two coin tosses.
\\[0.5em]
Continuing in this way, for each $n \in \Nb$ we can define $\Fc_n$ as the $\sig$-algebra containing the sets that are resolved by the first $n$ coin tosses.  Note that if a set $A \in \Fc_n$ then $A \in \Fc_{n+1}$.  Thus, $\Fc_n \subset \Fc_{n+1}$.  In other words, $\Fc_{n+1}$ contains more ``information'' than $\Fc_n$.  This kind of structure is encapsulated in the following definition.

\begin{definition}
\label{def:filtration}
Let $\Om$ be a nonempty set. Let $T $ be a fixed positive number (possibly infinity), and assume that for each $t \in [0,T]$ there is a $\sig$-algebra $\Fc_t$.  Assume further
that if $0 \leq s \leq  t \leq T$, then $\Fc_s \subseteq \Fc_t$. Then we call the sequence of $\sig$-algebras $\Fb = (\Fc_t)_{t \in [0,T]}$ a \emph{continuous time filtration}.
\end{definition}

\noindent
A \emph{discrete time filtration} is a sequence of $\sig$-algebras $\Fb=(\Fc_n)_{n \in \Nb_0}$ that satisfies $\Fc_n \subseteq \Fc_{n+1}$ for all $n$.

\begin{example}
Let $\Om = C_0[0,T]$, the set of continuous functions defined on $[0,T]$, starting from zero.  We denote by $\om = (\om_t)_{t \in [0,T]}$ and element of $\Om$.  Let $\Fc_t$ be the $\sig$-algebra generated by observing $\om$ over the interval $[0,t]$.  Mathematically, we write this as
\begin{align}
\Fc_t
	&:= \sig( \om_s , 0 \leq s \leq t ) .
\end{align}
It should be obvious that the sequence of $\sig$-algebras $\Fb=(\Fc_t)_{t \in [0,T]}$ forms a filtration.
Below, we define two sets, one of which is in $\Fc_t$, one of which is not.
\begin{align}
A := \{ \om : \sup_{0 \leq s \leq t} \om_s \leq 1 \}
	&\in \Fc_t , &
B:= \{ \om : \om_T \leq 1 \}
	&\notin \Fc_t .
\end{align}
The set $A$ is an element of $\Fc_t$ because, given the path of $\om$ over the interval $[0,t]$ one can answer the question: is the maximum of $\om$ over the interval $[0,t]$ less than $1$?  The set $B$ is not an element of $\Fc_t$ because one needs to know $\om_T$ in order to answer the question: is $\om_T \leq 1$?
\end{example}

\noindent
In the above example we generated a sequence of $\sig$-algebras by observing directly an element $\om \in \Om$.  Suppose that, instead of observing $\om$ we can observe only a random variable $X(\om)$.  We can use this information to generate a $\sig$-algebra as well.

\begin{definition}
\label{def:sigma.X}
Let  $X$  be a random variable defined on a nonempty sample space  $\Om$.  The \emph{$\sig$-algebra generated by $X$},  denoted  $\sig(X)$, is the collection of all
subsets of $\Om$ of the form  $\{X \in A\}$ where $A \in \Bc(\Rb)$.
\end{definition}

\begin{example}
Let us return to Example \ref{ex:stock}.  What is $\sig(S_2)$?  From the definition, we need to ask, which sets are of the form $\{S_2 \in A\}$?  Since $S_2$ can only take three values, $u^2$, $ud$ and $d^2$ we check the following sets
\begin{align}
\{ S_2  = u^2 \}
	&= A_{HH} , &
\{ S_2 = u d \} 
	&= A_{HT} \cup A_{TH} , &
\{ S_2 = d^2 \} 
	&= A_{TT} .
\end{align}
We add to these sets the sets that are necessary to form a $\sig$-algebra (i.e., $\emptyset$, $\Om$ and unions and complements of the above sets) to obtain
\begin{align}
\sig(S_2)
	&=	\sig( \{ A_{HH}, A_{TT}, A_{HT} \cup A_{TH} \} ) .
\end{align}
Note, that $\sig(S_2) \subset \Fc_2$ since $A_{HT}, A_{TH} \in \Fc_2$ but $A_{HT}, A_{TH} \notin \sig(S_2)$.  The reason is that, if $S_2 = ud$ we cannot say if $\om_1 = T$ or $\om_1 = H$.  
\end{example}

\begin{definition}
Let $X$ be a random variable defined on a nonempty sample space $\Om$.  Let $\Gc$ be a $\sig$-algebra of subsets of $\Om$.  If $\sig(X) \subset \Gc$ we say that \emph{$X$ is $\Gc$-measurable}, and we write $X \in \Gc$.
\end{definition}

\noindent
A random variable  $X$ is $\Gc$-measurable if and only if the information in $\Gc$ is sufficient to determine the value of $X$.  Obviously, if $X \in \Gc$ then $g(X) \in \Gc$ (assuming is $g$ is a measurable map from $(\Rb,\Bc(\Rb))$ to $(\Rb,\Bc(\Rb))$).
\\[0.5em]
Eventually, we will want to consider stochastic processes $X=(X_t)_{t \in [0,T]}$ and we will want to know at each time $t$ if $X_t$ is measurable with respect to $\sig$-algebra $\Fc_t$.

\begin{definition}
Let $\Om$ be a nonempty sample space equipped with a filtration $\Fb=(\Fc_t)_{t \in [0,T]}$.  Let $X=(X_t)_{t \in [0,T]}$ be a collection of random variables indexed by $t \in [0,T]$.  We say this collection of random variables is \emph{$\Fb$-adapted} if $X_t \in \Fc_t$ for all $t \in [0,T]$.
\end{definition}


\section{Independence}
When $X \in \Gc$ this means that the information in $\Gc$ is sufficient to determine the value of $X$.  On the other extreme, if $X$ is independent (a term we will define soon) of $\Gc$ this means that the information in $\Gc$ tells us nothing about the value of $X$.

\begin{definition}
\label{def:independent.sets}
Let  $(\Om,\Fc,\Pb)$ be a probability space. We say that two sets two sets $A$ and $B$ in $\Fc$ are \emph{independent}, written $A \ind B$, if $\Pb(A \cap B) = \Pb(A) \cdot \Pb(B)$.
\end{definition}

\begin{example}
Let us return to the coin-toss example of Section \ref{sec:coin-toss}.  Consider two sets and their intersection
\begin{align}
\{ \om_1 = H \}
	&=	A_H , &
\{ \om_2 = H \}
	&=	A_{HH} \cup A_{TH} &
\{ \om_1 = H \} \cap \{ \om_2 = H \}
	&=	A_{HH}
\end{align}
Since the coin tosses are independent, we should have $\{ \om_1 = H \} \ind \{ \om_2 = H \}$
Let us verify that these events are independent according to Definition \ref{def:independent.sets}.  We have
\begin{align}
\Pb( A_H )
	&=	p , &
\Pb( A_{HH} \cup A_{TH} )
	&=	p^2 + q p = p , &
\Pb( A_{HH} )
	&=	p^2 .
\end{align}
Thus, $\Pb( A_H ) \cdot \Pb( A_{HH} \cup A_{TH} ) = \Pb( A_{HH} )$, as expected.
\end{example}

\begin{example}
Can a set be independent of itself?  Surprisingly, the answer is ``yes.''  Suppose $A \ind A$.  Then, by the definition of independent sets, we have
$\Pb( A \cap A ) 	=	\Pb(A) \cdot \Pb(A)$.  We also have $\Pb( A \cap A ) = \Pb(A)$.  Combining these equations, we obtain $\Pb(A) \cdot \Pb(A) = \Pb(A)$.  This equation has two solutions $\Pb(A) = 1$ and $\Pb(A) = 0$.  Thus, a set is independent of itself if the probability of that set is zero or one.
\end{example}

\noindent
Having defined independent sets, we can now extend to independent $\sig$-algebras and random variables.

\begin{definition}
\label{def:independent.vars}
Let $(\Om,\Fc,\Pb)$  be a probability space, and let $\Gc$ and $\Hc$ be
sub-$\sig$-algebras of $\Fc$ (i.e., $\Gc,\Hc \subseteq \Fc$).  We say these \emph{two $\sig$-algebras are independent}, written $\Gc \ind \Hc$, if 
\begin{align}
\Pb(A \cap B) 
	&= \Pb(A) \cdot \Pb(B) , &
\forall \, A 
	&\in \Gc , &
\forall \, B
	&\in \Hc .
\end{align}
Let  $X$  and  $Y$  be random variables on $(\Om,\Fc,\Pb)$.  We say these \emph{two random variables are independent}, written $X \ind Y$, if $\sig(X) \ind \sig(Y)$.  Lastly, we say the random variable $X$ is \emph{independent} of the $\sig$-algebra $\Gc$, written $X \ind \Gc$, if $\sig(X) \ind \Gc$.
\end{definition}

\noindent
Recall from Definition \ref{def:sigma.X} that $\sig(X)$ contains all sets of the form $\{ X \in A \}$, where $A \in \Bc(\Rb)$.  Combining this with Definition \ref{def:independent.vars} we see that
\begin{align}
X
	&\ind Y  &
	&\Leftrightarrow&
\Pb( X \in A , Y \in B )
	&=	\Pb( X \in A) \cdot \Pb(Y \in B) , &
\forall \, A,B
	&\in \Bc(\Rb) . \label{eq:PAB=PA.PB}
\end{align}
It follows from \eqref{eq:PAB=PA.PB} that
\begin{align}
X
	&\ind Y  &
	&\Rightarrow&
\Eb X Y 
	&=	\Eb X \cdot \Eb Y . \label{eq:EAB=EA.EB}
\end{align}
Note that $\Eb XY = \Eb X \cdot \Eb Y$ does \emph{not} imply $X \ind Y$.
\\[0.5em] 
The above notion of independence is called \emph{pairwise} independence.  If $X \ind Y$ and $Y \ind Z$, this notion of independence \emph{not} imply $X \ind Z$ (for example, what if $Z=X$?).  Thus, at times, we may need a stronger notion of independence.

\begin{definition}
\label{def:ind.strong}
Let $(\Om,\Fc,\Pb)$  be a probability space, and let $\Gc_1, \Gc_2, \ldots , \Gc_n$ be
sub-$\sig$-algebras of $\Fc$.  We say the \emph{sequence of $\sig$-algebras are independent}, if 
\begin{align}
\Pb( \cap_{i=1}^n A_i) 
	&= \prod_{i=1}^n \Pb(A_i) , &
	&\forall \, A_1 \in \Gc_1 , \forall \, A_2 \in \Gc_2, \ldots ,  \forall A_n \in \Gc_n. 
\end{align}
Let  $X_1, X_2, \ldots, X_n$ be a sequence of random variables on $(\Om,\Fc,\Pb)$.  We say the \emph{sequence of random variables are independent} if the $\sig$-algebras $\sig(X_1), \sig(X_2), \ldots, \sig(X_n)$ are independent.
\end{definition}

\noindent
As with a pair of random variables, a sequence of random variables $(X_i)_{i \geq 1}$ is independent if and only if
\begin{align}
\Pb\( \cap_{i=1}^n \{X_i \in A_i\}  \)
	&=	\prod_{i=1}^n \Pb( X_i \in A_i) , &
\forall \, A_1 \in \Bc(\Rb), 
\forall \, A_2 \in \Bc(\Rb), \ldots ,
\forall \, A_n \in \Bc(\Rb). 
\end{align}

\noindent
We will often say that a sequence of random variables $(X_i)_{i \geq 0}$ is \emph{independent and identically distributed} (iid), by which me mean all $X_i$ have the same distribution and $(X_i)_{1 \leq i \leq n}$ are independent for every $n \in \Nb$.

\begin{example}
Let us return to the coin-toss example of Section \ref{sec:coin-toss}.  Let us define a sequence of random variables $(X_i)_{i \geq 1}$ via
\begin{align}
X_1(\om)
	&=	\left\{ \begin{aligned}
			& 1 & \om_1 = H , \\
			& 0 & \om_1 = T ,
			\end{aligned} \right. , &
X_2(\om)
	&=	\left\{ \begin{aligned}
			& 1 & \om_2 = H , \\
			& 0 & \om_2 = T ,
			\end{aligned} \right. , &
X_i(\om)
	&=	\left\{ \begin{aligned}
			& X_1 & \text{$i$ odd} , \\
			& X_2 & \text{$i$ even} .
			\end{aligned} \right. , &
\end{align}
Clearly, since the coin tosses are independent we have $X_1 \ind X_2$ and $X_i \ind X_j$ if $i$ is even and $j$ is odd.  But, the sequence $(X_i)_{1 \leq i \leq n}$ is not independent for any $n \geq 3$ since $X_i = X_{i+2n}$ for all $n \in \Nb$.
\end{example}

\noindent
It is not easy to verify if two random variables $X$ and $Y$ are independent using Expression \eqref{eq:PAB=PA.PB}, since the equation must be verified for \emph{all} Borel sets $A,B \in \Bc(\Rb)$.  In fact, there is an easier way to check independence.  

\begin{definition}
The \emph{joint distribution function} $F_{X,Y}: \Rb^2 \to [0,1]$ of two random variables $X$ and $Y$ defined on a probability space $(\Om,\Fc,\Pb)$ is given by
\begin{align}
F_{X,Y}(x,y)
	&:=	\Pb( X \leq x , Y \leq y) .
\end{align}
\end{definition}

\noindent
Again, we have two special cases for jointly discrete and jointly continuous random variables.

\begin{definition}
Two random variables $X$ and $Y$ are called \emph{jointly discrete} if the pair $(X,Y)$ takes values in some countable set $A = \{x_1,x_2,\ldots\} \times \{y_1,y_2,\ldots\} \subset \Rb^2$.  We associate is a discrete random variable a \emph{probability mass function} $f_{X,Y} : A \to \Rb$, defined by $f_{X,Y}(x_i,y_j) := \Pb( X = x_i, Y = y_j )$. 
\end{definition}

\begin{definition}
A pair of random variables $X$ and $Y$ is called \emph{jointly continuous} if its joint distribution function $F_{X,Y}$ can be written as
\begin{align}
F_{X,Y}(x,y)
	&=	\int_{-\infty}^x \int_{-\infty}^y \dd u \dd v \, f_{X,Y}(u,v) , &
(x,y)
	&\in \Rb^2 ,
\end{align}
for some $f_{X,Y} : \Rb^2 \to [0,\infty)$ called the \emph{joint probability density function}.
\end{definition}

\noindent
As in the one-dimensional case, it may help to think of the joint density function $f_{X,Y}$ as $f_{X,Y}(x,y) \dd x \dd y = \Pb( X \in \dd x , Y \in \dd y )$.
\\[0.5em]
Note that for jointly continuous random variables $X$ and $Y$ we have $f_{X,Y}(x,y) = \d_x\d_y F_X(x,y)$.
\\[0.5em]
If the pair $(X,Y)$ is either jointly discrete or jointly continuous, it is easy to compute $\Pb( (X,Y) \in A)$ for any $A \in \Bc(\Rb^2)$.  We have
\begin{align}
\text{discrete}:&&
\Pb( (X,Y) \in A)
	&=	\sum_{ \{i,j : (x_i,y_j) \in A \}} f_{X,Y}(x_i,y_j) , \\
\text{continuous}:&&
\Pb( (X,Y) \in A)
	&=	\int_A \dd x \dd y \, f_{X,Y}(x,y) .
\end{align}

\noindent
To recover the \emph{marginal} distribution $F_X$ from $F_{X,Y}$, simply note that
\begin{align}
F_X(x)
	&=	\Pb( X \leq x )
	=		\Pb( X \leq x, Y \leq \infty )
	=		F_{X,Y}(x,\infty) .
\end{align}
It follows that for the discrete and continuous cases, we have, respectively
\begin{align}
\text{discrete}:&&
f_X(x_i)
	&=	\sum_j f_{X,Y}(x_i,y_j) , \\
\text{continuous}:&&
f_X(x)
	&=	\int_\Rb \dd y \, f_{X,Y}(x,y) .
\end{align}

\noindent
The following theorem gives some easy-to-check conditions for independence.

\begin{theorem}
Let $X$ and $Y$ be random variables defined on a probability space $(\Om,\Fc,\Pb)$.  The following conditions are equivalent (that is, if one of them holds, all of them hold)
\begin{enumerate*}
\item $X \ind Y$.
\item $F_{X,Y}(x,y) = F_X(x) F_Y(y)$ for every $(x,y) \in \Rb^2$.
\item Discrete case: $f_{X,Y}(x,y) = f_X(x) f_Y(y)$ for every $(x,y) \in \Rb^2$.\\
	Continuous case: $f_{X,Y}(x,y) = f_X(x) f_Y(y)$ for `almost' every $(x,y) \in \Rb^2$.
\item $\Eb \ee^{ \ii u X + \ii v Y} = \Eb \ee^{u \ii X} \cdot \Eb \ee^{ \ii v Y}$ for all $(u,v) \in \Rb^2$.
\end{enumerate*}
\end{theorem}

\noindent
Together with expectation, the most important statistical properties of a random variable (or pair) are the variance and co-variance.

\begin{definition}
The \emph{variance} of a random variable $X$, written $\Vb X $ is defined by
\begin{align}
\Vb X 
	&=	\Eb (X - \Eb X)^2 
	=		\Eb  X^2  - (\Eb X )^2 , \label{eq:var}
\end{align}
whenever the expectation exists.
\end{definition}

\begin{definition}
The \emph{co-variance} of two random variables $X$ and $Y$, written $\text{Co}\Vb[ X , Y ]$ is defined by
\begin{align}
\text{Co}\Vb[ X, Y ]
	&=	\Eb (X - \Eb X)(Y - \Eb Y) 
	=		\Eb X Y  - \Eb X \cdot \Eb Y  , \label{eq:covar}
\end{align}
whenever the expectation exists.
\end{definition}

\noindent
Note that $\text{Co}\Vb[ X, X ] = \Vb X$.

\noindent
Note $\Vb \cdot $ is \emph{not} a linear operator, since
\begin{align}
\Vb[ a X + b Y ]
	&=	a^2 \Vb X  + b^2 \Vb Y  + 2 ab \, \text{Co}\Vb[X,Y] .
\end{align}
where $a$ and $b$ are constants.

\begin{definition}
We say two random variables $X$ and $Y$ are \emph{uncorrelated} if $\text{Co}\Vb[X,Y]=0$.
\end{definition}

\noindent
Note that $X \ind Y$ implies $X$ and $Y$ are uncorrelated.  However, the converse is \emph{not} true.  


\section{Conditional expectation}

Let $X$ be a random variable defined on $(\Om,\Fc,\Pb)$ and let $\Gc$ be a sub-$\sig$ algebra of $\Fc$.
When $X \in \Gc$ this means that the information in $\Gc$ is sufficient to determine the value of $X$.
When $X \ind \Gc$, this means that the information in $\Gc$ gives us no information at all about $X$.
Usually, however, the information in $\Gc$ gives us some information about $X$, but not enough to determine $X$ exactly.  And this brings us to the notion of conditioning.
\\[0.5em]
Presumably, you have run across the following formula for the conditional probability of a set $A$ given $B$
\begin{align}
\Pb( A | B )
	&=	\frac{ \Pb( A \cap B ) }{ \Pb( B ) } , &
\Pb(B)
	&> 	0 .
\end{align}
When $(X,Y)$ have are jointly discrete or jointly continuous, this readily leads to \emph{conditional probability mass function}
\begin{align}
\text{discrete}:&&
f_{X|Y}(x_i,y_j)
	&:=	\Pb( X = x_i | Y = y_j )
	=		\frac{ \Pb( X = x_i \cap Y = y_j ) }{ \Pb( Y = y_j ) }
	=		\frac{ f_{X,Y}(x_i,y_j)}{f_Y(y_j)} , \\
\text{continuous}:&&
f_{X|Y}(x,y) \dd x 
	&:=	\Pb( X \in \dd x | Y = y )
	=		\frac{ \Pb( X \in \dd x \cap Y \in \dd y ) }{ \Pb( Y \in \dd y ) }
	=		\frac{ f_{X,Y}(x,y)}{f_Y(y)} \dd x .
\end{align}
And from this, we can define $\Eb[ X | Y = y ]$, the \emph{conditional expectation of $X$ given $Y = y$}
\begin{align}
\text{discrete}:&&
\Eb[ X | Y = y_j ]
	&:=	\sum_i x_i \, f_{X|Y}(x_i,y_j) , \\
\text{continuous}:&&
\Eb[ X | Y = y ]
	&:=	\int_\Rb \dd x \, x \, f_{X|Y}(x,y)
\end{align}
Note that $\Eb[ X | Y = y ]$ is simply a function of $y$ -- there is nothing random about it.
\\[0.5em]
Unfortunately, there are cases for which the pair $(X,Y)$ are neither jointly discrete nor jointly continuous.  And, for these cases we need a more general notion of conditional expectation.  Here we will make two conceptual leaps:
\begin{enumerate*}
\item We will condition with respect to a $\sig$-algebra rather than conditioning on an event.
\item The conditional expectation will be a random variable.
\end{enumerate*}
We will just hop in with our new definition of conditional expectation and then we will see, through an example, that his new definition makes sense.

\begin{definition}
\label{def:cond.E}
Let $(\Om,\Fc,\Pb)$  be a probability space, let  $\Gc$ be a sub-$\sig$-algebra of $\Fc$, and let $X$ be a random variable that is either nonnegative or integrable.
The  conditional expectation of $X$ given $\Gc$, denoted $\Eb[X|\Gc]$, is any random
variable that satisfies
\begin{enumerate*}
\item Measurability: $\Eb[X|\Gc] \in \Gc$.
\item Partial averaging: $\Eb[ \Ib_A \Eb[X|\Gc] ] = \Eb[ \Ib_A X ]$ for all $A \in \Gc$.  \\
Alternatively, $\Eb[ Z \Eb[X|\Gc] ] = \Eb[ Z X ]$ for all bounded random variables $Z \in \Gc$.
\end{enumerate*}
When $\Gc = \sig(Y)$ we shall often use the short-hand notation $\Eb[X|Y]:=\Eb[X|\sig(Y)]$.
\end{definition}

\noindent
Conditional probabilities are defined from conditional expectations using
\begin{align}
\Pb( A | \Gc )
	&=	\Eb[ \Ib_A | \Gc ] .
\end{align}

\noindent
Admittedly, Definition \ref{def:cond.E} is rather abstract (and, for the purposes of computation, useless).  
In fact, it is not at all clear from Definition \ref{def:cond.E} that $\Eb[X|\Gc]$ even exists!  It does exist, though we will not prove this here.
\\[0.5em]
Conditional expectation has an interesting $L^2$ interpretation.  Consider a probability space $(\Om,\Fc,\Pb)$. Let $\Gc$ be a sub-$\sig$-algebra of $\Fc$ (i.e., $\Gc \subseteq \Fc$).  Define
\begin{align}
L^2(\Om,\Fc,\Pb)
	&:=	\{ X : \Om \to \Rb \text{ s.t. } X^{-1}(B) \in \Fc \,\, \forall B \in \Bc(\Rb) \text{ and } \Eb |X|^2 < \infty \} ,
\end{align}
and likewise for $L^2(\Om,\Gc,\Pb)$.  Clearly, since $\Gc \subseteq \Fc$ we have $L^2(\Om,\Gc,\Pb) \subseteq L^2(\Om,\Fc,\Pb)$.  Next, define an inner product
\begin{align}
\< X, Y \>
	&:=	\Eb X Y , &
X,Y
	&\in L^2(\Om,\Fc,\Pb).
\end{align}
From the definition of conditional expectation, we have
\begin{align}
0
	&=	\< X - \Eb[X|\Gc], Z \> 
	=		\Eb(X - \Eb[X|\Gc]) Z , &
	&\forall \, Z \in \Gc .
\end{align}
Thus, $\Eb[X|\Gc]$ is the projection of $X \in L^2(\Om,\Fc,\Pb)$ onto the subspace $L^2(\Om,\Gc,\Pb)$.
\\[0.5em]
When conditioning on the $\sig$-algebra generated by a random variable, it is easiest to use the following formula
\begin{align}
\Eb[ X | Y ]
	&=	\psi(Y) , &
\psi(y)
	&:=	\Eb[ X | Y = y ] . \label{eq:cond.E}
\end{align}
The following example should help to build some intuition for conditional expectation.

\begin{example}
Let $\Om = \{ a, b, c, d, e, f \}$, $\Fc = 2^\Om$ and $\Pb( \om ) = (1/6)$ for $\om=a,b,\ldots,f$.  Define two random variables $X$ and $Y$ on $(\Om,\Fc,\Pb)$ as follows
\begin{align}
\( \begin{array}{c} \om \\ X(\om) \\ Y(\om) 
\end{array} \)
	&=	\( \begin{array}{cccccc} 
			a & b & c & d & e & f \\ 
			1 & 3 & 3 & 3 & 5 & 7 \\ 
			2 & 2 & 1 & 1 & 7 & 7 
\end{array} \) .
\end{align}
Let $\Gc = \sig(Y)$.  Next, let us compute $\Eb[X|Y]$ using \eqref{eq:cond.E} and check if it agrees with Definition \ref{def:cond.E}.  We have
\begin{align}
\( \begin{array}{c} 
\om \\ Y(\om) \\ \Eb[X|Y](\om) 
\end{array} \)
	&=	\( \begin{array}{cccccc} 
			a & b & c & d & e & f \\ 
			2 & 2 & 1 & 1 & 7 & 7 \\
			2 & 2 & 3 & 3 & 6 & 6 
			\end{array} \) .
\end{align}
Le us check measurability: is $\Eb[X|Y] \in \sig(Y)$?  In other words, is $\sig(\Eb[X|Y]) \subseteq \sig(Y)$?  We have
\begin{align}
\sig(Y)
	&=	\sig( \{a,b\}, \{c,d\}, \{e,f\} ) , \\
\sig(\Eb[X|Y])
	&=	\sig( \{a,b\}, \{c,d\}, \{e,f\} ) .
\end{align}
So, yes, $\Eb[X|Y] \in \sig(Y)$.  Another way to think of measurability is to simply ask: given the value of $Y$, can one determine the value of $\Eb[X|Y]$?  Clearly, in this case the answer is ``yes.''
Next, let us check the partial averaging property: does $\Eb[ \Ib_A \Eb[X|Y] ] = \Eb[ \Ib_A X ]$ for all $A \in \sig(Y)$.  Rather than check this for \emph{every} $A \in \sig(Y)$, let us just check that this holds for the sets $\{a,b\}$, $\{c,d\}$ and $\{e,f\}$.  We have
\begin{align}
A
	&=	\{ a, b\}, &
\Eb[ \Ib_A \Eb[ X | Y ] ]
	&=	\Pb(a) 2 + \Pb(b) 2 = 4/6, &
\Eb[ \Ib_A X ]
	&=	\Pb(a) 1 + \Pb(b) 3 = 4/6, \\
A
	&=	\{ c, d\}, &
\Eb[ \Ib_A \Eb[ X | Y ] ]
	&=	\Pb(c) 3 + \Pb(d) 3 = 6/6, &
\Eb[ \Ib_A X ]
	&=	\Pb(c) 3 + \Pb(d) 3 = 6/6, \\
A
	&=	\{ e, f\}, &
\Eb[ \Ib_A \Eb[ X | Y ] ]
	&=	\Pb(e) 6 + \Pb(f) 6 = 12/6, &
\Eb[ \Ib_A X ]
	&=	\Pb(e) 5 + \Pb(f) 7 = 12/6 .
\end{align}
\end{example}

\noindent
The following properties are arguably more important to remember than the definition of conditional expectation.  Memorize them!

\begin{theorem}
\label{thm:cond.E}
Let $(\Om, \Fc, \Pb)$ be a probability space and let $\Gc$ be a sub-$\sig$-algebra of $\Fc$.  Conditional expectations satisfy the following properties.
\begin{enumerate*}
\item Linearity: $\Eb[ a X + b Y | \Gc ] = a \Eb[ X | \Gc] + b \Eb[ Y | \Gc ]$.
\item Taking out what is known: if $X \in \Gc$ and $\Eb |XY| < \infty$ then $\Eb[ X Y | \Gc ] = X \Eb[ Y | \Gc ]$.
\item Iterated conditioning: if $\Hc$ is a sub-$\sig$-algebra of $\Gc$ then $\Eb[\Eb[X|\Gc]|\Hc]=\Eb[X|\Hc]$.
\item Independence: if $X \ind \Gc$ then $\Eb[X|\Gc]=\Eb X$.
\end{enumerate*}
\end{theorem}

\noindent
Theorem \ref{thm:cond.E} can be proved directly from Definition \ref{def:cond.E}, though we will not do so here.  In addition to the above properties, the following Theorem, which we state without proof is often useful:

\begin{theorem}[Jensen's inequality]
\label{thm:jensen}
Let $X$ be a random variable defined on $(\Om,\Fc,\Pb)$ and let $\Gc$ be a sub-$\sig$-algebra of $\Fc$.  Suppose $\varphi : \Rb \to \Rb$ is a convex function and that $\Eb |X| < \infty$ and $\Eb |\varphi(X)| < \infty$.  Then we have
\begin{align}
\varphi ( \Eb[ X | \Gc ] )
	&\leq \Eb[ \varphi(X) | \Gc ] , &
	&\Pb-\text{a.s.} \label{eq:jensen}
\end{align}
\end{theorem}

\noindent
In order to keep straight which direction the inequality in \eqref{eq:jensen} goes, it is helpful to remember that $\phi(x) = x^2$ is a convex function and the conditional variance of a function satisfies
\begin{align}
\Vb[ X | \Gc]
	&=	\Eb[ X^2 | \Gc ] - \Eb[ X | \Gc ]^2
	\geq 0 .
\end{align}
Now that we have defined conditional expectation and established some of its key properties, we can define ``Markov process'' and ``martingale'' -- two seemingly similar, but very distinct concepts.

\begin{definition}
\label{def:martingale}
 Let $(\Om, \Fc, \Pb)$ be a probability space, let $T$ be a fixed positive
number, and let $\Fb=(\Fc_t)_{t \in [0,T]}$ be a filtration of sub-$\sig$-algebras of $\Fc$.
Consider an $\Fb$-adapted stochastic process $M=(M_t)_{t \in [0,T]}$.  
Suppose $\Eb |M_t| < \infty$ for all $t \in [0,T]$.  Then, we say that $M$ is
\begin{align}
\text{a \emph{martingale} if}&&
\Eb[ M_t | \Fc_s ]
	&=	 M_s , &
	&		\forall \, 0 \leq s \leq t \leq T , \\
\text{a \emph{sub-martingale} if}&&
\Eb[ M_t | \Fc_s ]
	&\geq	M_s , &
	&		\forall \, 0 \leq s \leq t \leq T , \\
\text{a \emph{super-martingale} if}&&
\Eb[ M_t | \Fc_s ]
	&\leq	M_s , &
	&		\forall \, 0 \leq s \leq t \leq T . 
\end{align}
\end{definition}

\noindent
We have given above the definition of a \emph{continuous time} martingale (resp. sub-, super-).  We can also define \emph{discrete-time} martingales by making the obvious modifications.  In fact, the condition for a process $M$ to be a discrete-time martingale is somewhat simpler
\begin{align}
\Eb[ M_{n+1} | \Fc_n ]
	&=	 M_n , &
	&		\forall \, n \in \Nb_0 . &
	&		\text{(discrete-time martingale condition)} \label{eq:disc-mart}
\end{align}
Note that if \eqref{eq:disc-mart} holds then we also have $\Eb[ M_{n+m} | \Fc_n ] = M_n$ for all $n,m \in \Nb_0$ because
\begin{align}
\Eb[M_{n+m}|\Fc_n]
	&=	\Eb[ \Eb[ M_{n+m}|\Fc_{n+m-1}]|\Fc_n]
	=		\Eb[M_{n+m-1}|\Fc_n] \\
	&=	\Eb[ \Eb[ M_{n+m-1}|\Fc_{n+m-2}]|\Fc_n]
	=		\Eb[M_{n+m-2}|\Fc_n] \\
	&=	\ldots \\
	&=	\Eb[ \Eb[ M_{n+1}|\Fc_{n}]|\Fc_n]
	=		\Eb[M_{n}|\Fc_n]
	=		M_n ,
\end{align}
where we have used the tower property of conditional expectation as well as taking out what is known.
\\[0.5em]
Admittedly, the definition of sub- and super-martingales seems backwards; sub-martingales tend to rise in expectation, whereas a super-martingales tend to fall in expectation.
\\[0.5em]
Note: when we say that a process $M$ is a martingale (or sub- or super-martingale) this is with respect to a fixed probability measure and filtration.  If $\Pb$ and $\Pbt$ are two probability measures and $\Fb$ and $\Gb$ are two filtrations, it is entirely possible that a process $M$ may be a martingale with respect to $(\Pb,\Fb)$ and may not be a martingale with respect to $(\Pbt,\Fb)$, $(\Pb,\Gb)$ or $(\Pbt,\Gb)$.

\begin{definition}
\label{def:markov}
Let  $(\Om,\Fc,\Pb)$  be a probability space, let $T$ be a fixed positive number, and let $\Fb=(\Fc_t)_{t \in [0,T]}$ be a filtration of sub-$\sig$-algebras of $\Fc$.
Consider an $\Fb$-adapted stochastic process $X=(X_t)_{t \in [0,T]}$. Assume that for all $0 \leq s \leq t \leq T$ and for every nonnegative, Borel-measurable function $f$, there is another Borel-measurable function $g$ (which depends on $s$, $t$, and $f$) such that
\begin{align}
\Eb[ f(X_t) | \Fc_s ]
	&=	g(X_s) .
\end{align}
Then we say $X$ is a \emph{Markov process} or simply ``$X$ is Markov.''
\end{definition}

\noindent
Identifying $g(X_s) \equiv \Eb[ f(X_t) | X_s ]$ we can write the Markov property as follows
\begin{align}
\Eb[ f(X_t) | \Fc_s ]
	&=	\Eb[ f(X_t) | X_s ]
\end{align}
A Markov process is a process for which the following holds: \emph{given the present (i.e., $X_s$), the future (i.e, $X_t$, $t \geq s$) is independent of the past (i.e, $\Fc_s$)}.  What this means in practice is that 
\begin{align}
\Pb( X_t \in A | \Fc_s )
	&=	\Pb(X_t \in A | X_s) , &
	&		\forall s \leq t ,
	&		\forall \, A \in \Bc(\Rb) . &
	&		\text{(if $X$ is Markov)} \label{eq:markov.prop}
\end{align}
If $X_t$ is a discrete or continuous random variable for every $t$ then we have a \emph{transition kernel}, written as $P$ in the discrete case and $\Gam$ in the continuous case.
\begin{align}
\text{discrete}:&&
P(s,x;t,y)
	&:= \Pb(X_t = y | X_s = x) , \\
\text{continuous}:&&
\Gam(s,x;t,y) \dd y
	&=	\Pb(X_t \in \dd y | X_s = x) .
\end{align}
If you can write the transition kernel of a process explicitly, then you have essentially proved that the process is Markov.
\\[0.5em]
Note that \emph{any process} that has independent increments is Markov since, if $X_t - X_s \ind X_s$ for $t \geq s$, then
\begin{align}
\Pb( X_t \in A | \Fc_s)
	&=	\Pb( X_t - X_s + X_s \in A | X_s) ,
\end{align}
Where $\Fc_s$ is the filtration generated by observing $X$ up to time $s$.
\\[0.5em]
Markov processes and Martingales are \emph{entirely separate} concepts.  A process $X$ can be both a martingale and a Markov process, it can be a martingale but not a Markov process, it can be a Markov process but not a martingale, and it can be neither a Markov process nor a martingale.  We illustrate the difference with an example.

\begin{example}
Let us return to the stock price Example \ref{ex:stock}.  Let us show that $S=(S_n)_{0 \leq n}$ is a Markov process.  Recall that $\Fc_m$ is the $\sig$-algebra generated by observing $\om_1, \om_2, \ldots, \om_m$.
Observe that $S_m \in \Fc_m$.  Next, note that
\begin{align}
\Pb( S_{n+m} = S_m u^k d^{n-k} | S_m )
	&=	\binom{n}{k} p^k q^{n-k} .
\end{align}
Since we have written the transition kernel explicitly, we have established that $S$ is Markov.  Let us also find the function $g$ in Definition \ref{def:markov}.
For any $f:\Rb \to \Rb$ we have
\begin{align}
\Eb[ f(S_{n+m}) | \Fc_m ]
	&=	\sum_{k=0}^n f(S_m u^k d^{n-k}) \cdot \binom{n}{k} p^k q^{n-q}
	=:	g(S_m) . \label{eq:Ef=g}
\end{align}
Thus, we have found $g$.
Now, to see if $S$ is a martingale note that
\begin{align}
\Eb[ S_{n+1} | \Fc_n ]
	&=	\Eb[ S_{n+1} | S_n ]
	=		p \cdot u S_n + q \cdot d S_n
	=		( p \cdot u + q \cdot d ) S_n .
\end{align} 
Thus, if $( p \cdot u + q \cdot d )=1$, then $\Eb[ S_{n+1} | \Fc_n ] = S_n$.  Let us assume that $( p \cdot u + q \cdot d )=1$.  Then we have
\begin{align}
\Eb[S_{n+m}|\Fc_n]
	&=	\Eb[ \Eb[ S_{n+m}|\Fc_{n+m-1}]|\Fc_n]
	=		\Eb[S_{n+m-1}|\Fc_n] \\
	&=	\Eb[ \Eb[ S_{n+m-1}|\Fc_{n+m-2}]|\Fc_n] 
	=		\Eb[S_{n+m-2}|\Fc_n] \\
	&=	\ldots \\
	&=	\Eb[ \Eb[ S_{n+1}|\Fc_{n}]|\Fc_n] 
	=		\Eb[S_{n}|\Fc_n]
	=		S_n ,
\end{align} 
Therefore, the process $S$ is a martingale.  Similarly, if $( p \cdot u + q \cdot d ) \geq 1$, then $S$ is a sub-martingale and if $( p \cdot u + q \cdot d ) \leq 1$, then $S$ is a super-martingale.
\end{example}



\section{Stopping times}
\label{sec:stopping.time}
Roughly speaking, a \emph{Stopping time} is a particular kind of random time $\tau$ marks the time at which a stochastic process $X=(X_t)_{t \geq 0}$ exhibits a particular kind of behavior.  A stopping time is often defined by a \emph{stopping rule}, e.g., when the process $X$ does \emph{this}, then \emph{stop}.  Mathematically, a stopping time is defined as follows:

\begin{definition}
\label{def:stopping.time}
Fix a probability space $(\Om,\Fc,\Pb)$ and a filtration $\Fb = (\Fc_t)_{t \geq 0}$.  A random time $\tau : \Om \to [0,\infty]$ is called a \emph{$\Fb$-stopping time} if it satisfies
\begin{align}
\{ \tau \leq t \}
	&\in \Fc_t , &
\forall \, t
	&\in [0,\infty) . \label{eq:stopping.time}
\end{align}
\end{definition}

\noindent
Above, we have focused on the continuous-time setting.  If we are working with a discrete-time filtration $\Fb = (\Fc_n)_{n \in \Nb_0}$, then a stopping time is a random time $\tau : \Om \to \Nb_0 \cup \{ \infty \}$ that satisfies
\begin{align}
\{ \tau \leq n \}
	&\in \Fc_n , &
\forall \, n
	&\in \Nb_0 . 
\end{align}
Observe that, stopping times, like martingales, are defined with respect to a specific filtration $\Fb$.  Also note that a stopping time may be infinite.  The meaning of a stopping time should be fairly clear from \eqref{eq:stopping.time}.  If $\tau$ is an $\Fb$-stopping time, then for any $t \geq 0$ we should be able to say whether or not $\tau$ has occurred given the information in $\Fc_t$.  Below, we give two examples of random times, one of which is a stopping time, one of which is not.

\begin{example}
Let $X = (X_t)_{t \geq 0}$ be a continuous time stochastic process on $(\Om,\Fc,\Pb)$ and let $\Fb = (\Fc_t)_{t \geq 0}$ be the filtration generated by observing the path of $X$, that is, $\Fc_t = \sig( X_s, 0 \leq s \leq t)$.  Define
\begin{align}
\underline{\tau}
	&:=	\inf\{ t \geq 0 : X_t = a \} , &
\overline{\tau}
	&=	\sup\{ t \geq 0 : X_t = a \} , &
a
	&\in \Rb .
\end{align}
The first random time $\underline{\tau}$ is clearly a $\Fb$-stopping time.  At any time $t \geq 0$, using the information in $\Fc_t$, we can clearly answer the question: \emph{has $X$ hit $a$}?  On the other hand, the second random time $\overline{\tau}$ is \emph{not} a stopping time.  At time $t$ we will not be able to answer the question: \emph{when is the last time that $X$ will hit $a$}?
\end{example}

\noindent
The following theorem, will prove to be incredibly useful.

\begin{theorem}
\label{thm:stopped.martingale}
Suppose a process $M = (M_t)_{t \geq 0}$ with right-continuous sample paths is a martingale with respect to a filtration $\Fb = (\Fc_t)_{t \geq 0}$.
Suppose further that $\tau$ is an $\Fb$-stopping time.  Then the \emph{stopped process} $M^\tau := (M_{t \wedge \tau})_{t \geq 0}$ is a martingale.  That is, for any $0 \leq t \leq T < \infty$, we have
\begin{align}
\Eb[ M_{T \wedge \tau} | \Fc_t ]
	&=	M_{t \wedge \tau} .
\end{align}
In particular, $M_0 = \Eb[ M_{t \wedge \tau} ]$ for all $t \geq 0$.
\end{theorem}

\noindent
We will not prove Theorem \ref{thm:stopped.martingale} rigorously.  However, we will provide a bit of intuition for why the theorem is true.  For any $t < \tau$ we have $M_{\tau \wedge t} = M_t$ and $M$ is a martingale.  Likewise, for any $t \geq \tau$ we have $M_t^\tau = M_\tau$, which is a constant (and thus trivially a martingale).
As the process $M^\tau$ is a martingale both prior to and after $\tau$, it is reasonable to expect that $M^\tau$ is in fact a martingale.
\\[0.5em]
Just as we can construct a $\sig$-algebra $\Fc_t$ of information up to a fixed time $t$, we can construct a $\sig$-algebra $\Fc_\tau$ of information up to a stopping time $\tau$.  

\begin{definition}
\label{def:Fc.tau}
Fix a probability space $(\Om,\Fc,\Pb)$ and a filtration $\Fb = (\Fc_t)_{t \geq 0}$.  Suppose $\tau$ is an $\Fb$-stopping time.  Then we define the $\sig$-algebra $\Fc_\tau$ at the stopping time $\tau$ as follows
\begin{align}
\Fc_\tau
	&:=	\{ A \in \Fc_\infty : A \cap \{ \tau \leq t \} \in \Fc_t \, \forall \, t \geq 0 \} . \label{eq:Fc.tau}
\end{align}
\end{definition}

\noindent
The idea underlying Definition \eqref{def:Fc.tau} is that if a set $A$ is observable at time $\tau$ then for any time $t$, its restriction to the set $\{\tau\leq t\}$ should be in $\Fc_t$.  The restriction to sets in $A \in \Fc_\infty$ is to take account of the possibility that the stopping time can be infinite and ensures that $A \cap \{ \tau \leq \infty \} \in \Fc_\infty$. 
If $A \notin \Fc_\infty$ then we could have $A \cap \{ \tau \leq \infty \} \notin \Fc_\infty$.
From the above definition, a random variable $X$ is $\Fc_\tau$-measurable if and only if 
$\Ib_{\{\tau \leq t\}} X \in \Fc_t$ for all $t \in [0,\infty]$.
\\[0.5em]
Now, suppose a process $M = (M_t)_{t \geq 0}$ is a martingale with respect to a filtration $\Fb = (\Fc_t)_{t \geq 0}$.
Suppose further that $M$ has a well-defined limit $M_\infty := \lim_{t \to \infty} M_t$.
As $M$ is a martingale we have (by definition) that $M_{t_1} = \Eb[ M_{t_2} | \Fc_{t_1}]$ for any $0 \leq t_1 \leq t_2 \leq \infty$.
Now, consider two $\Fb$-stopping times $\tau_1$ and $\tau_2$ satisfying $0 \leq \tau_1 \leq \tau_2 \leq \infty$.
One may ask: is it true that $M_{\tau_1} = \Eb[ M_{\tau_2} | \Fc_{\tau_1} ]$?
Unfortunately, the answer to this question is ``no, not in general.''  However, under certain conditions, provided in the following theorem, we will have $M_{\tau_1} = \Eb[ M_{\tau_2} | \Fc_{\tau_1} ]$.

\begin{theorem}[Doob's Optional Stopping]
\label{thm:sampling}
Let $M = (M_t)_{t \geq 0}$ be a martingale with right-continuous sample paths with respect to a filtration $\Fb = (\Fc_t)_{t \geq 0}$.
Suppose that for any $\eps > 0$ there exists a constant $K_\eps \in [0,\infty)$ such that
\footnote{We shall see in Definition \ref{def:UI} that this condition is called \textit{uniform integrability}.}
\begin{align}
\sup_{t \geq 0} \Eb \Ib_{\{|M_t|>K_\eps\}} |M_t|
	&< \eps . \label{eq:ui.condition}
\end{align}
Suppose that $\tau_1$ and $\tau_2$ are $\Fb$-stopping times and that $0 \leq \tau_1 \leq \tau_2$.
Then the limit $M_\infty := \lim_{t \to \infty} M_t$ exists almost surely (so that $M_{\tau_i}$ is well-defined on $\{ \tau_i = \infty \}$) and we have
\begin{align}
M_{\tau_1} 
	&=	\Eb[ M_{\tau_2} | \Fc_{\tau_1} ] . \label{eq:martingale.stopping}
\end{align}
In particular, we have $M_0 = \Eb M_{\tau_i}$ for $i = 1, 2$.
\end{theorem}

\noindent
We will not prove Theorem \ref{thm:sampling}.  Rather, we will demonstrate its usefulness through an example.

\begin{example}
Consider a symmetric random walk
\begin{align}
S_n
	&=	\sum_{i=1}^n X_i , &
S_0
	&=	x \in \Zb ,
\end{align}
where the $(X_i)$ are i.i.d. random variables with $\Pb( X_i = 1 ) = \Pb( X_i = - 1) = 1/2$.  
Let $\Fb = (\Fc_n)_{n \in \Nb_0}$ where $\Fc_n := \sig( S_i , 0 \leq i \leq n)$.
For any $a \in \Zb$ define the first hitting time to $a$ as follows
\begin{align}
\tau_a
	&:=	 \inf \{ n \in \Nb_0 : S_n = a \} .
\end{align}
Now, suppose $a < x < b$.  We wish to find $\Pb( \tau_a < \tau_b)$.  To this end, we define\
\begin{align}
\tau
	&:=	\tau_a \wedge \tau_b .
\end{align}
Observe that $\tau$ is an $\Fb$-stopping time and $S$ is an $\Fb$ martingale.  The stopped process $S^\tau := ( S_{n \wedge \tau} )_{n \in \Nb}$ is a martingale by Theorem \ref{thm:stopped.martingale}.  Moreover, because $a < S_n^\tau < b$ we see that $S^\tau$ satisfies \eqref{eq:ui.condition}.  Thus, we have by Theorem \ref{thm:sampling} that
\begin{align}
S_0^\tau = x
	&=	\Eb S_\tau^\tau
	=		a \Pb( \tau_a < \tau_b ) + b \Pb( \tau_b < \tau_a ) ,
\end{align}
where we have used the fact that $S_\tau^\tau = a$ when $\tau_a < \tau_b$ and $S_\tau^\tau = b$ when $\tau_b < \tau_a$.  We also have
\begin{align}
1
	&=		\Pb( \tau_a < \tau_b ) + \Pb( \tau_b < \tau_a ) .
\end{align}
Thus, using the above two equations, we find
\begin{align}
\Pb( \tau_a < \tau_b )
	&=	\frac{b - x}{b - a} , &
\Pb( \tau_b < \tau_a )
	&=	\frac{x - a}{b - a} .
\end{align}
\end{example}

\noindent
The condition given in \eqref{eq:ui.condition} is essential in order for \eqref{eq:martingale.stopping} to hold, as we shall see in the next example.

\begin{example}
Let $B=(B_i)_{i \in \Nb}$ be an iid sequence of Bernoulli random variables: $B \sim \text{Ber}(p)$ with $p = 1/2$.  Construct a sequence $X = (X_i)$ of random variables as follows
\begin{align}
X_0
	&=	1 , &
X_{n+1}
	&=	2 B_{n+1} X_n .
\end{align}
Define a filtration $\Fb = (\Fc_n)_{n \in \Nb_0}$ where $\Fc_n = \sig( B_i, 0 \leq i \leq n)$.  Observe that both $B$ and $X$ are $\Fb$-adapted.  It is easy to see that $X$ is a martingale with respect to $\Fb$ because
\begin{align}
\Eb[ X_{n+1} | \Fc_n ]
	&=	2 X_n \Eb[ B_{n+1} ]
	=		X_n .
\end{align}
Note, however, that condition \eqref{eq:ui.condition} does not hold.  Now, let us define the first hitting time of $X$ to zero: $\tau := \inf \{ n \geq 0 : X_n = 0  \}$.  We note that $\Pb(\tau < \infty ) = 1$.  We have
\begin{align}
\Eb X_\tau
	&=	\Eb 0
	=		0 \neq X_0 .
\end{align}
Thus, we see that, without condition \eqref{eq:ui.condition}, we cannot expect \eqref{eq:martingale.stopping} to hold.
\end{example}

%


\section{Exercises}
\label{sec:exercises-2}

\begin{exercise}
\label{exer:2.1}
Let $\Om = \{a,b,c,d\}$ and let $\Fc = 2^\Om$ (the set of all subsets of $\Om$).  We define a probability measure $\Pb$ as follows
\begin{align}
\Pb(a)
	&=	1/6, &
\Pb(b)
	&=	1/3, &
\Pb(c)
	&=	1/4, &
\Pb(d)
	&=	1/4, 
\end{align}
Next, define three random variables
\begin{align}
X(a)
	&=	1, &
X(b)
	&=	1, &
X(c)
	&=	-1, &
X(d)
	&=	-1,  \\
Y(a)
	&=	1, &
Y(b)
	&=	-1, &
Y(c)
	&=	1, &
Y(d)
	&=	-1 , 
\end{align}
and $Z = X+Y$. (a) List the sets in $\sig(X)$.  (b) What are the values of $\Eb[Y|X]$ for $\{a,b,c,d\}$?  Verify the partial averaging property: $\Eb[\Ib_A\Eb[Y|X]]=\Eb[\Ib_A Y]$ for all $A \in \sig(X)$.  (c) What are the values of $\Eb[Z|X]$ for $\{a,b,c,d\}$?  Verify the partial averaging property.
\end{exercise}

\begin{exercise}
\label{exer:2.2}
Fix a probability space $(\Om,\Fc,\Pb)$.  Let $Y$ be a square integrable random variable: $\Eb Y^2 < \infty$ and let $\Gc$ be a sub-$\sig$-algebra of $\Fc$.  Show that
\begin{align}
\Vb( Y - \Eb[Y|\Gc] )
	&\leq \Vb(Y - X) , &
\forall X
	&\in \Gc.
\end{align}
\end{exercise}

\begin{exercise}
\label{exer:2.3}
Give an example of a probability space $(\Om,\Fc,\Pb)$, a random variable $X$ and a function $f$ such that $\sig(f(X))$ is strictly smaller than $\sig(X)$ but $\sig(f(X)) \neq \{ \emptyset, \Om \}$.   Give a function $g$ such that $\sig(g(X)) = \{\emptyset,\Om\}$.
\end{exercise}

\begin{exercise}
\label{exer:2.4}
On a probability space $(\Om,\Fc,\Pb)$ define random variables $X$ and $Y_0, Y_1, Y_2, \ldots$ and suppose $\Eb |X|<\infty$.  Define $\Fc_n := \sig(Y_0, Y_1,\ldots,Y_n)$ and $X_n = \Eb[X|\Fc_n]$.  Show that the sequence $X_0,X_1,X_2,\ldots$ is a martingale under $\Pb$ with respect to the filtration $(\Fc_n)_{n \geq 0}$.
\end{exercise}

\begin{exercise}
\label{exer:2.5}
In this exercise, we will construct a discrete-time model for a stock $Z = (Z_n)_{n \geq 0}$.
Let $X_0, X_1, \ldots$ be a sequence of i.i.d Bernoulli random variables with parameter $p \in (0,1)$ (i.e., $\Pb(X_i = 1) = p$).  Define $S_n = \sum_{i=1}^n X_i$ where $S_0 = 0$.  Let $Z$ be given by
\begin{align}
Z_n
	&:=	r^{2 S_n - n} , &
r
	&> 0 , &
n
	&=0, 1, 2, \ldots .
\end{align}
(a) Let $\Fb = (\Fc_n)_{n \geq 0}$ where $\Fc_n := \sig(X_0,X_1,\ldots,X_n)$.  For what value(s) of $r$ is $Z$ a $(\Pb,\Fb)$ martingale, submartingale and supermartingale?  In particular, show that $Z$ is a $(\Pb,\Fb)$ martingale when $r = (1-p)/p$.\\
(b) Let $\Gb = (\Gc_n)_{n \geq 0}$ where $\Gc_n := \sig(S_0,S_1,\ldots,S_n)$. Suppose $r$ is such that $Z$ is a $(\Pb,\Fb)$ martingale.  Is $Z$ also a $(\Pb,\Gb)$ martingale?\\
(c) Suppose that $\Pbt( X_i = 1 ) = \pt \neq p$.  Suppose $r$ is such that $Z$ is a $(\Pb,\Fb)$ martingale.  Is $Z$ also a $(\Pbt,\Fb)$ martingale?
\end{exercise}

\begin{exercise}
\label{exer:2.6}
This exercise is similar to Exercise \ref{exer:2.2}.  Suppose $X$ and $Y$ are random variables defined on a probability space $(\Om,\Fc,\Pb)$.  Define the following \textit{functional}
\begin{align}
J(f)
	&:=	\Eb[ ( f(X) - h(X,Y) )^2 ] ,
\end{align}
where $\Eb h^2(X,Y) < \infty$.  We wish to find the function $f^*$ that minimizes $J$.  That is
\begin{align}
J(f^*)
	&=	\min_f J(f) .
\end{align}
Note that, for any $g$, the following must hold
\begin{align}
0
	&=	\frac{\dd}{\dd \eps} J(f^* + \eps g) \Big|_{\eps = 0} .
\end{align}
If you do not understand why, think about what happens to $J(f^* + \eps g)$ for values of $\eps$ near zero.
Use this information to show that $f^*(x) = \Eb[ h(X,Y) | X = x]$.
If you find it helpful, you can assume $(X,Y)$ are jointly continuous random variables, and therefore have a joint density $p(x,y)$.  Though, this is not necessary to prove the result.  Essentially, what this result shows, is that the $\sig(X)$-measurable random variable $Z$ that minimizes $\Eb[ (Z - h(X,Y))^2 ]$ is $Z = \Eb[h(X,Y)|X]$.
\end{exercise}


%
%

\chapter{Generating and Characteristic Functions}
\label{chapter:generating}
The notes from this chapter are taken primarily from \cite[Chapter 5]{gs2001}.


\section{Generating functions}

The generating function is a powerful tool for analyzing sums of random variables.

\begin{definition}
Suppose $X$ is a discrete random variable taking values in $\{0\} \cup \Nb$.  We define the \emph{probability generating function} of $X$, written $G_X(s)$, by
\begin{align}
G_X(s)
	&:=	\Eb s^X 
	=		\sum_k s^k f_X(k) .
\end{align}
where $f_X(k)$ is the probability mass function of $X$.
\end{definition}

\noindent
The generating function exists at least when $|s|\leq1$ (since $\sum_k f_X(k)=1$) and may exist in a larger interval.  As we will see, we usually only care about the value of $G_X$ and its derivatives at the point $s=1$.  You can always integrate and differentiate a convergent power series term-by-term within its radius of convergence.
\\[0.5em]
Why is $G_X$ call the ``probability generating function''?  The reason is that the coefficient $f_X(k)$ of the $\Oc(s^k)$ term in the series expansion of $G_X(s)$ is precisely $\Pb( X= k )$.  So, one can expand $G_X$ as a power series and obtain the probability mass function $f_X$.

\begin{example}
Suppose $X \sim \text{Ber}(p)$ as in Example \ref{ex:bernoulli}.  Then
\begin{align}
G_X(s) 
	&= 	q s^0 + p s^1 
	=		q + p s .
\end{align}
\end{example}

\begin{example}
\label{ex:G.exp}
Suppose $X \sim \text{Geo}(p)$ as in Example \ref{ex:geometric}.  Then
\begin{align}
G_X(s)
	&=	\sum_{k=1}^\infty s^k p q^{k-1}
	=		sp \sum_{k=1}^\infty (sq)^{k-1}
	=		\frac{sp}{1-sq} , &
|s|
	&<	\frac{1}{q} .
\end{align}
\end{example}

\begin{theorem}
\label{thm:G.moments}
Let $G_X^{(n)}$ denote the $n$th derivative of $G_X$ and let $G_X^{(n)}(1) := \lim_{s \nearrow 1} G_X^{(n)}(s)$, which may be $+\infty$.  Then
\begin{align}
G_X^{(n)}(1)
	&=	\Eb \frac{X!}{(X-n)!} , \label{eq:GX.1}
\end{align}
where $X!/(X-n)! = 0$ when $X < n$.
\end{theorem}
\begin{proof}
We have
\begin{align}
G_X^{(n)}(s)
	&=	\frac{\dd^n}{\dd s^n} G_X(s)
	=		\Eb  \frac{\dd^n}{\dd s^n} s^X 
	=		\Eb  \frac{X!}{(X-n)!} s^{X-n} , \label{eq:GX.s}
\end{align}
where passing the derivative through the expectation is allowed by dominated convergence.
Taking $s=1$ in \eqref{eq:GX.s} yields \eqref{eq:GX.1}
\end{proof}

\begin{example}
How would you find $\Vb X $ using $G_X$?  Note that
\begin{align}
G_X'(1)
	&=	\Eb X  , &
G_X''(1)
	&=	\Eb X (X-1) = \Eb X^2  - \Eb X  . \label{eq:V.G}
\end{align}
Using the above, we compute
\begin{align}
\Vb X 
	&=	\Eb X^2  - (\Eb X)^2
	=		G_X''(1) - (G_X'(1))^2 + G_X'(1) . 
\end{align}
\end{example}

\noindent
Given two independent random variables $X$ and $Y$ taking values in $\{0\} \cup \Nb$, one can compute $\Pb( X+Y = k)$ using
\begin{align}
\Pb( X+Y = n)
	&=	\sum_{k=0}^n \Pb( X = k ) \Pb( Y = n-k ) .
\end{align}
Alternatively, one can compute the probability generating function of $X+Y$ and then expand the generating function to obtain the probability $\Pb(X+Y = n)$.

\begin{theorem}
Suppose $X \ind Y$.  Then $G_{X+Y}(s) = G_X(s) G_Y(s)$.
\end{theorem}
\begin{proof}
We have
\begin{align}
G_{X+Y}(s) 
	&= \Eb s^{X+Y} 
	=		\Eb s^X  \cdot \Eb s^Y 
	= G_X(s) G_y(s) .
\end{align}
\end{proof}

\noindent
The real use of generating functions arises when one wants to compute $\Pb( \sum_{i=1}^n X_i = k)$ where $(X_i)_{i \geq 1}$ is an iid sequence of random variables.  In this case, computing the generating function $G_{S_n}(s)$ with $S_n = \sum_{i=1}^n X_i$ is relatively easy.  Upon computing the generating function $G_{S_n}$, one can compute probabilities of the form $\Pb( S_n = k )$ by expanding $G_{S_n}$.  By contrast, computing $\Pb( S_n = k )$ directly from the probability mass function of $X_i$ is difficult.

\begin{theorem}
Suppose $(X_i)_{i \geq 1}$ are iid with common distribution $X$.  Define $S_n := \sum_{i=1}^n X_i$.  Then
\begin{align}
G_{S_n}(s)
	&=	[ G_X(s) ]^n	.
\end{align}
\end{theorem}
\begin{proof}
We have
\begin{align}
G_{S_n}(s)
	&=	\Eb s^{X_1 + X_2 + \ldots + X_n}
	=		(\Eb s^X )^n
	=		[G_X(s)]^n .
\end{align}
\end{proof}

\begin{example}
Let $(X_i)_{i \geq 1}$ be iid with each $X_i \sim \text{Ber}(p)$.  Define $S_n := \sum_{i=1}^n X_i$.  Recall from Example \ref{ex:binomial} that $S_n \sim \text{Bin}(n,p)$.  We have
\begin{align}
G_{S_n}(s)
	&=	[G_X(s)]^n = (q + p s)^n .
\end{align}
Thus, we have obtained the probability generating function of the binomial random variable $S_n$.  Although we could have computed $G_{S_n}$ using the probability mass function $f_{S_n}$ directly, this approach would have been much more work.
\end{example}

\begin{theorem}
\label{thm:random.sum}
Suppose $(X_i)_{i \geq 1}$ are iid with common distribution $X$.  Define $S_n := \sum_{i=1}^n X_i$.  
Let $N$ be independent of $(X_i)_{i \geq 1}$.  Then
\begin{align}
G_{S_N}(s)
	&=	G_N( G_X(s) )	.
\end{align}
\end{theorem}
\begin{proof}
We have
\begin{align}
G_{S_N}(s)
	&=	\Eb s^{X_1 + X_2 + \ldots + X_N}
	=		\Eb \Eb[ s^{X_1 + X_2 + \ldots + X_N} | N] \\
	&=	\Eb (\Eb s^X)^N 
	=		\Eb (G_X(s))^N  
	=		G_N( G_X(s) ) ,
\end{align}
where, in the third equality, we have used the fact that the $X_i$ are independent of each other and independent of $N$.
\end{proof}

\noindent
One can extend the notion of a probability generating function to multiple random variables in the obvious way.

\begin{definition}
Suppose $X$ and $Y$ are discrete random variables taking values in $\{0\} \cup \Nb$.  We define the \emph{joint probability generating function} of $(X,Y)$, written $G_{X,Y}(s,t)$, as
\begin{align}
G_{X,Y}(s,t)
	&=	\Eb s^X t^Y 
	=		\sum_k \sum_m s^k t^m f_{X,Y}(k,m) .
\end{align}
where $f_{X,Y}(k,m)$ is the joint probability mass function of $(X,Y)$.
\end{definition}

\noindent
The coefficient $f_{X,Y}(k,m)$ of the $\Oc(s^k t^m)$ term in the power series expansion of $G_{X,Y}(s,t)$ about $(0,0)$ gives $\Pb( X=k, Y=m)$.

\begin{theorem}
Let $G_{X,Y}^{(n,m)}$ denote the $(n,m)$th partial derivative of $G_{X,Y}$ with respect to the first and second arguments.  Then
\begin{align}
G_{X,Y}^{(n,m)}(1,1)
	&=	\Eb \frac{X!}{(X-n)!} \frac{Y!}{(Y-m)!}  . \label{eq:GXY.1}
\end{align}
\end{theorem}
\begin{proof}
The proof follows the proof the Theorem \ref{thm:G.moments}.
\end{proof}


\section{Branching processes}
\label{sec:branching}

Suppose a population evolves in generations.  Let $Z_n$ be the size of the $n$th generation and assume $Z_0 = 1$.  We assume that the $i$th member of the $n$th generation gives birth to a random number $X_{n,i}$ of individuals.  Then, clearly, the number of individuals in the $(n+1)$th generation is given by
\begin{align}
Z_{n+1}
	&=	X_{n,1} + X_{n,2} + \ldots + X_{n,Z_n} . \label{eq:Zn}
\end{align}
We would like to know what the probability mass function, mean and variance of $Z_n$ are.  Our route for obtaining these will be through the generating function $G_{Z_n}$.

\begin{theorem}
Assume $Z_{n+1}$ is given by \eqref{eq:Zn} and that the collection of random variables $(X_{n,i})$ are iid.  Define $G_n(s) := \Eb s^{Z_n}$ and $G(s) := \Eb s^X$.  Then
\begin{align}
G_{n+m}(s)
	&=	G_n( G_m(s) ) = G_m( G_n(s) ) , &
	&\text{and thus}&
G_n(s)
	&=	\underbrace{ G( G( \ldots G(s) \ldots ) ) }_{ \text{$n$-fold iteration} } .
\end{align}
\end{theorem}

\begin{proof}
Let $Y_{m,i}$ denote the number of offspring in the $(m+n)$th generation that descend from the $i$th member of the $m$th generation.  Clearly, the $(Y_{m,i})$ are iid and $Y_{m,i} \sim Z_n$.  The number of members of the $(m+n)$th generation is given by
\begin{align}
Z_{m+n}
	&=	Y_{m,1} + Y_{m,2} + \ldots + Y_{m,Z_m} .
\end{align}
We have a random sum of iid random variables.
By Theorem \ref{thm:random.sum} it follows that $G_{m+n}(s) =	G_m(G_n(s))$.  Obviously, we can interchange $m \leftrightarrow n$ and obtain $G_{m+n}(s) =	G_n(G_m(s))$.
Finally, we have
\begin{align}
G_n(s)
	&=	G( G_{n-1}(s) )
	=		G( G( G_{n-2}(s) )
	= 	\ldots
	=		\underbrace{ G( G( \ldots G(s) \ldots ) ) }_{ \text{$n$-fold iteration} } ,
\end{align}
as claimed.
\end{proof}

\noindent
In principle, one can obtain the probability mass function $f_n$ of the $n$th generation $Z_n$ by expanding the generating function $G_n(s)$ as a power series about $s=0$.  In practice, this may be difficult to carry out.  However, moments of $Z_n$ can usually be computed with relative ease.

\begin{theorem}
Suppose $\Eb Z_1 = \mu$ and $\Vb Z_1 = \sig^2$.  Then
\begin{align}
\Eb Z_n
	&=	\mu^n , &
\Vb Z_n
	&=	\left\{ \begin{array}{ll}
			n \sig^2 & \mu = 1 , \\
			\frac{\sig^2 (\mu^n-1)\mu^{n-1}}{\mu-1} & \mu \neq 1 .
			\end{array} \right. \label{eq:mu.and.sig}
\end{align}
\end{theorem}

\begin{proof}
We have
\begin{align}
\Eb Z_n
	&=	\frac{\dd}{\dd s}G_n(s) \Big|_{s=1}
	=		G_{n-1}'(G(1)) \cdot G'(1)
	=		G_{n-1}'(1) \cdot \mu
	=		\Eb Z_{n-1} \cdot \mu ,
\end{align}
where we have used $G(1) = \Eb 1^{Z_1} = 1$.
Iterating, we obtain $\Eb Z_n = \mu^n$.  Next, from \eqref{eq:V.G}, we compute
\begin{align}
\Eb Z_n^2 - \Eb Z_n
	&=	\frac{\dd^2}{\dd s^2}G_n(s) \Big|_{s=1} 
	=		\frac{\dd}{\dd s} G_{n-1}'(G(s)) \cdot G'(s) \Big|_{s=1} \\
	&=	G_{n-1}''(G(s)) \cdot ( G'(s) )^2 +  G_{n-1}'(G(s)) \cdot G''(s) \Big|_{s=1} \\
	&=	G_{n-1}''(1) \cdot ( G'(1) )^2 +  G_{n-1}'(1) \cdot G''(1) \\
	&=	\( \Eb Z_{n-1}^2 - \Eb Z_{n-1} \) \cdot ( \Eb Z_1 )^2 + \Eb Z_{n-1} \cdot \( \Eb Z_1^2 - \Eb Z_1 \) . \label{eq:EZ2-EZ}
\end{align}
Thus, using \eqref{eq:EZ2-EZ}, as well as $\Eb Z_k = \mu^k$ and $\Eb Z_k^2 = \Vb Z_k + (\mu^k)^2$, we obtain
\begin{align}
\Vb Z_n
	&=	\Eb Z_n^2 - ( \Eb Z_n )^2 \\
	&=	\( \Eb Z_{n-1}^2 - \Eb Z_{n-1} \) \cdot ( \Eb Z_1 )^2 + \Eb Z_{n-1} \cdot \( \Eb Z_1^2 - \Eb Z_1 \) - ( \Eb Z_n )^2 + \Eb Z_n \\
	&=	\( \Vb Z_{n-1} + (\Eb Z_{n-1})^2 - \Eb Z_{n-1} \) \cdot ( \Eb Z_1 )^2 + \Eb Z_{n-1} \cdot \( \Vb Z_1 + (\Eb Z_1)^2 - \Eb Z_1 \) - ( \Eb Z_n )^2 + \Eb Z_n \\
	&=	\( \Vb Z_{n-1} + \mu^{2(n-1)} - \mu^{n-1} \) \cdot \mu^2 + \mu^{n-1} \cdot \( \sig^2 + \mu^2 - \mu \) - \mu^{2n} + \mu^n .
\end{align}
Thus, we have obtained an expression for $\Vb Z_n$ in terms of $\Vb Z_{n-1}$, $\mu$ and $\sig^2$.  Solving this explicitly yields \eqref{eq:mu.and.sig}.
\end{proof}


\section{Characteristic functions}

Suppose a random variable $X$ has a state space that is not contained in the non-negative integers.  When we cannot use the probability generating function $G_X(s)$ to analyze the properties of $X$.  A natural candidate to replace the probability generating function is the \emph{moment generating function} $M_X(t):=\Eb \ee^{t X}$.  Moment generating functions, however, are limited in their use because, for many random variables, $\Eb \ee^{t X}= \infty$ for all $t \neq 0$.  However, there is a simple fix to this.  We consider instead the characteristic function.

\begin{definition}
The \emph{characteristic function} of a random variable $X$ is the function $\phi_X : \Rb \to \Cb$ defined by
\begin{align}
\phi_X(t)
	&:=	\Eb \ee^{\ii t X} , &
\ii
	&=	\sqrt{-1} .
\end{align}
\end{definition}  

\noindent
The characteristic function (unlike the moment generating function) always exists since $\Eb |\ee^{\ii t X } | = 1$.  
Clearly, if $G_X$ and $\phi_X$ both exist, then we have
\begin{align}
\phi_X(t) 
	&=	G_X(\ee^{\ii t}) &
	&\text{since}&
\Eb \ee^{\ii t X}
	&=	\Eb s^X , &
s
	&=	\ee^{\ii t} .
\end{align}
We have the following obvious properties.

\begin{theorem}
\label{thm:char.props}
Characteristic functions have the following properties.
\begin{enumerate*}
\item Let $a,b$ be constants.  Then $\phi_{a X + b}(t) = \ee^{\ii b t}\phi(at)$.
\item If $X$ and $Y$ are independent then $\phi_{X+Y}(t) = \phi_X(t)\phi_Y(t)$.
\item If $(X_i)_{1 \leq i \leq n}$ are iid and $S_n = \sum_{i=1}^n X_i$, then $\phi_{S_n}(t) = (\phi_X(t))^n$.
\end{enumerate*}
\end{theorem}

\begin{example}
\label{ex:char.exp}
Let $X \sim \Ec(\lam)$ as in Example \ref{ex:exponential}.  Then
\begin{align}
\phi_X(t)
	&=	\int_0^\infty \dd x \, \ee^{\ii t x} \lam \ee^{-\lam x}
	=		\frac{\lam}{\lam-\ii t} .
\end{align}
\end{example}

\begin{example}
\label{ex:char.normal}
Let $X \sim \Nc(\mu,\sig^2)$ as in Example \ref{ex:normal}.  Then
\begin{align}
\phi_X(t)
	&=	\int_{-\infty}^\infty \dd x \, \ee^{\ii t x} \frac{1}{\sqrt{2 \pi \sig^2}} \exp\( \frac{-(x-\mu)^2}{2 \sig^2} \)
	=		\exp\( \ii \mu t - \frac{1}{2} \sig^2 t^2 \) .
\end{align}
\end{example}

\noindent
Characteristic functions have several uses.  First, they can be used to capture the moments of a random variable when they exist.

\begin{theorem}
Let $\phi_X(t)$ be the characteristic function of a random variable $X$.  If $\Eb |X|^n < \infty$, then
\begin{align}
\phi^{(n)}(0)
	&=	\ii^n \Eb X^n .
\end{align}
\end{theorem}
\begin{proof}
We have
\begin{align}
\frac{\dd^n}{\dd t^n} \phi_X(t)
	&=	\Eb \frac{\dd^n}{\dd t^n} \ee^{\ii t X}
	=		\ii^n \Eb X^n \ee^{\ii t X} ,
\end{align}
Now set $t=0$ to complete the proof.
\end{proof}

\noindent
The characteristic function uniquely determines the distribution of a random variable.  In other words, there is a one-to-one correspondence between $F_X$ and $\phi_X$.  We show this for a continuous random variable.

\begin{theorem}[Inversion]
\label{thm:inversion}
Suppose the characteristic function $\phi_X$ of a random variable $X$ is absolutely integrable: $\int_\Rb \dd t \, |\phi_X(t)| < \infty$.  Then $X$ is a continuous random variable with a continuous density $f_X$, which is given by
\begin{align}
f_X(x)
	&=	\frac{1}{2\pi} \int_\Rb \dd t \, \ee^{- \ii t x } \phi_X(t) , &
x
	&\in \Rb .
\end{align}
\end{theorem}
\noindent
To obtain $F_X$ from $f_X$ simply use $F_X(x) = \int_{-\infty}^x \dd y \, f_X(y)$.
\begin{proof}[Sketch of proof]
The proof of Theorem \ref{thm:inversion} follows from standard Fourier results.  Assuming $X$ has a density $f_X$, Fourier inversion gives
\begin{align}
f_X(x)
	&=	\frac{1}{2\pi} \int_\Rb \dd t \, \ee^{- \ii t x } \int_\Rb \dd y \, \ee^{\ii t y} f_X(y) 
	=	\frac{1}{2\pi} \int_\Rb \dd t \, \ee^{- \ii t x } \Eb \ee^{\ii t X}
	=		\frac{1}{2\pi} \int_\Rb \dd t \, \ee^{- \ii t x } \phi_X(t) .
\end{align}
\end{proof}

\noindent
An inversion theorem for random variables that are not continuous also exists.  Though, it is not particularly useful for the purposes of computation.
\\[0.5em]
Perhaps the most important property of characteristic functions is they can be used to prove the convergence of a sequence of random variables to a limiting distribution.

\begin{definition}
\label{def:Fn.convergence}
We say that a sequence of distribution functions $(F_n)_{n \geq 1}$ \emph{converges to a distribution $F$}, written $F_n \to F$, if $\lim_{n \to \infty} F_n(x) = F(x)$ at all points $x$ where $F$ is continuous.    
\end{definition}

\begin{theorem}[Continuity Theorem]
\label{thm:continuity}
Let $(F_n)_{n \geq 1}$ be a sequence of distributions functions with corresponding characteristic functions $(\phi_n)_{n \geq 1}$.  
\begin{enumerate*}
\item If $F_n \to F$ where $F$ is a distribution function with corresponding characteristic function $\phi$, then $\phi_n \to \phi$ (pointwise).
\item Conversely, if $\phi(t):=\lim_{n \to \infty} \phi_n(t)$ exists and is continuous at $t=0$, then $\phi$ is the characteristic function of a distribution $F$ and $F_n \to F$. \label{itm:phi.to.F}
\end{enumerate*}
\end{theorem}

\noindent
Item \ref{itm:phi.to.F} in Theorem \ref{thm:continuity} if particularly powerful.  If $F_n$ and $\phi_n$ are, respectively, the distribution and characteristic function of a sum of $n$ independent random variables, it is often easier to compute $\phi_n$ than $F_n$.  If we can compute $\phi_n$ and find its limit, then we can obtain $F$.

\begin{example}
In Example \ref{ex:G.exp} we found that 
\begin{align}
X 
	&\sim \text{Geo}(p)&
	&\Rightarrow&
G_X(s)
	&= \frac{sp}{1-sq} . 
\end{align}
From this, one can easily show that $\Eb X = 1/p$.  Now consider a sequence of random variables $(Y_n)_{n \geq 0}$ defined by
\begin{align}
Y_n
	&= 	\frac{1}{n} X_n , &
X_n
	&\sim \text{Geo}(\lam/n) ,
\end{align}
where $\lam>0$ is a fixed constant.  For $n$ large enough, $\lam/n < 1$.  
Note that $\Eb Y_n = \frac{1}{n} \Eb X_n = \frac{1}{n} \frac{n}{\lam}= \frac{1}{\lam}$ for all $n$.
We would like to know what the limiting distribution of the sequence $(Y_n)_{n \geq 0}$ is; we will use the Continuity Theorem \ref{thm:continuity} to do this.
We have
\begin{align}
\phi_{Y_n}(t)
	&=	\phi_{X_n}(t/n)
	=		\frac{ \ee^{\ii t/n}(\lam/n)}{1-\ee^{\ii t/n}(1-\lam/n)}
	=		\frac{ \ee^{\ii t/n}\lam }{n -\ee^{\ii t/n}(n-\lam)}
	=		\frac{ \ee^{\ii t/n}\lam }{n(1 -\ee^{\ii t/n}) + \ee^{\ii t/n} \lam} \\
	&=		\frac{ (1 + \ii t/n)\lam }{n(1 - 1 - \ii t/n - \ldots ) + (1 + \ii t/n + \ldots ) \lam}
	\to	\frac{\lam}{\lam - \ii t} = \phi_Z(t) ,
\end{align}
where $Z \sim \Ec(\lam)$ (see Example \ref{ex:char.exp}).  By the Continuity Theorem \ref{thm:continuity}, since $\phi_{Y_n} \to \phi_Z$ it follows that $F_{Y_n} \to F_Z$.
\end{example}


\section{LLN and the CLT}
\label{sec:slln.clt}
In this section we will use Theorem \ref{thm:continuity} to prove the Law of Large Numbers (LLN) and the Central Limit Theorem (CLT).  First, we must define the notion of convergence in distribution.

\begin{definition}
\label{def:D.convergence}
We say that a sequence of random variables $(X_n)_{n \geq 1}$ \emph{converges in distribution} to a random variable $X$, written $X_n \stackrel{\Dc}{\to} X$, if $F_{X_n} \to F_X$.
\end{definition}

\begin{theorem}[Law of Large Numbers]
\label{thm:lln}
Let $(X_n)_{n \geq 1}$ be a sequence of iid random variables with $\Eb X_n = \mu$.  
Define a sequence of random variables $(S_n)_{n \geq 1}$ by $S_n := \frac{1}{n}\sum_{i=1}^n X_i$.  Then $S_n \stackrel{\Dc}{\to} \mu$.
\end{theorem}
\begin{proof}
From Theorem \ref{thm:char.props} we have
\begin{align}
\phi_{S_n}(t)
	&=	\( \phi_X(t/n) \)^n
	=		\( 1 + \frac{\ii \mu t}{n} + \Oc((t/n)^2) \)^n
	\to		\ee^{\ii \mu t} , &
	&\text{as $n\to  \infty$} .
\end{align}
Note that $\ee^{\ii \mu t} = \phi_\mu(t)$, where $\phi_\mu$ is the characteristic function of a constant random variable $\mu$.
Thus, we have shown that $\phi_{S_n} \to \phi_\mu$.  By Theorem \ref{thm:continuity} we have $F_{S_n} \to F_\mu$ and thus, from Definition \ref{def:D.convergence}, we have $S_n \stackrel{\Dc}{\to} \mu$.
\end{proof}

\noindent
In the proof of Theorem \ref{thm:lln} we actually assumed that $\Eb X_n^2 < \infty$ when we wrote the $\Oc((t/n)^2)$.  In fact, the Law of Large Numbers (LLN) holds even when $\Eb X_n^2 = \infty$.

\begin{theorem}[Central Limit Theorem]
\label{thm:clt}
Let $(X_n)_{n \geq 1}$ be a sequence of iid random variables with $\Eb X_n = \mu$ and $\Vb X_n = \sig^2 \in (0,\infty)$.  
Define two sequences of random variables $(S_n)_{n \geq 1}$ and $(U_n)_{n \geq 1}$ by
\begin{align}
S_n 
	&= \sum_{i=1}^n X_i , &
U_n 
	&=	\frac{S_n - n \mu}{\sqrt{n \sig^2}} .
\end{align}
Then $U_n \stackrel{\Dc}{\to} Z$ where $Z \sim \Nc(0,1).$
\end{theorem}
\begin{proof}
First, we re-write $U_n$ as follows
\begin{align}
U_n
	&=	\frac{1}{\sqrt{n}} \sum_{i=1}^n Y_i , &
Y_i
	&:=	\frac{X_i - \mu}{ \sig } .
\end{align}
Note that $\Eb Y_n = 0$ and $\Vb Y_n = 1$.  Next, we compute the characteristic function of $U_n$.  From Theorem \ref{thm:char.props} we have
\begin{align}
\phi_{U_n}(t)
	&=	\( \phi_Y(t/\sqrt{n}) \)^n
	=		\( 1 - \frac{ t^2}{2 n} + \Oc((t/n)^3) \)^n
	\to		\exp \( - \frac{1}{2} t^2 \) , &
n
	&\to  \infty .
\end{align}
From Example \ref{ex:char.normal}, we see that $\exp \( - \frac{1}{2} t^2 \) = \phi_Z(t)$ where $Z \sim \Nc(0,1)$.  Thus, we have shown that $\phi_{U_n} \to \phi_Z$.  By Theorem \ref{thm:continuity} we have $F_{U_n} \to F_Z$ and thus, from Definition \ref{def:D.convergence}, we have $U_n \stackrel{\Dc}{\to} Z$, as claimed.
\end{proof}

\noindent
In the proof of Theorem \ref{thm:clt} we assume $\Eb |X_i|^3 < \infty$ when we wrote the $\Oc((t/n)^3)$ term.  The Central Limit Theorem (CLT) holds even when $\Eb |X_i|^3 = \infty$.


\section{Large Deviations Principle}
\label{sec:ldp}
Let $S_n$ be the sum of $n$ iid random variables $(X_i)$ with common distribution $X$ whose mean a variance are $\mu$ and $\sig$, respectively
\begin{align}
S_n
	&=	\sum_{i=1}^n X_i , &
\Eb X
	&=	\mu , &
\Vb X
	&=	\sig^2 .
\end{align}
The SLLN asserts that, as $n \to \infty$, we roughly have $S_n \sim n \mu$.
Similarly, the CLT asserts that, as $n \to \infty$, we roughly have $| S_n - n \mu | \sim \sqrt{n}$.
But, of course, $S_n$ is a random variable, not a constant.  And it is quite possible that $| S_n - n \mu|$ is $\Oc( n^\alpha )$ where $\alpha > 1/2$.  
The \emph{theory of large deviations} studies the asymptotic behavior of $\Pb( |S_n - n \mu| > n^\alpha )$ as $n \to \infty$ for $\alpha > 1/2$.
\\[0.5em]
Let us introduce $\Lam_X$, the \emph{cumulant generating function of $X$}, which is given by
\begin{align}
\Lam_X(t)
	&=	\log M_X(t) , &
M_X(t)
	&=	\Eb \ee^{t X} , \label{eq:cumulant}
\end{align}
assuming $M_X(t)$ exists in a neighborhood of zero.
We note that $\Lam_X'(0) = \mu$ because
\begin{align}
\Lam_X'(0)
	&=	\frac{M_X'(0)}{M_X(0)} =  \frac{ \Eb X \ee^{t X} }{\Eb \ee^{t X} } \Big|_{t=0} = \Eb X = \mu .
\end{align}
We also note that $\Lam_X$ is convex because
\begin{align}
\Lam_X''(t)
	&=	\frac{ M_X(t) M_X''(t) - ( M_X'(t))^2 }{ M_X^2(t) }
	=		\frac{ \Eb \ee^{t X} \Eb X^2 \ee^{t X} - ( \Eb X \ee^{t X} )^2 }{ M_X^2(t) } \geq 0 , \label{eq:note0}
\end{align}
where we have used the Cauchy-Schwarz inequality
\begin{align}
( \Eb Y Z )^2
	&\leq \Eb Y^2 \Eb Z^2, &
	&\text{with}&
Y
	&= X \ee^{t X/2}, &
Z
	&=	\ee^{t X/2} .
\end{align}
Next, we define $\Lam_X^*$, the \emph{Fenchel-Legendre transform} for $\Lam_X$ as follows
\begin{align}
\Lam_X^*(a)
	&=	\sup_{t \in \Rb} \{ a t - \Lam_X(t) \} , &
a
	&\in \Rb . \label{eq:legendre.transform}
\end{align}
We are now in a position to state the main theorem of this section

\begin{theorem}[Large Deviations Principle]
\label{thm:ldp}
Let $(X_i)$ be a sequence of iid random variables with common distribution $X$.  Define $\mu := \Eb X$ and suppose the moment generating function $M_X(t) := \Eb \ee^{t X}$ is finite in some neighborhood of $t = 0$.  Let $\Lam_X$ and $\Lam_X^*$ be given by \eqref{eq:cumulant} and \eqref{eq:legendre.transform}, respectively.  Suppose $a > \mu$ and $\Pb(X > a ) > 0$.  Then $\Lam_X^*(a) > 0$ and
\begin{align}
\lim_{n \to \infty} \frac{1}{n} \log \Pb ( S_n > n a ) 
	&= - \Lam_X^*(a) , &
S_n
	&=	\sum_{i=1}^n X_i . \label{eq:ldp}
\end{align}
\end{theorem}

\noindent
Theorem \ref{thm:ldp} asserts that, under appropriate conditions, we have $\Pb(S_n > n a) \sim \exp ( - n \Lam_X^*(a) )$.
Although the theorem appears to deal only with deviations of $S_n$ in excess of its mean, the corresponding result for deviations of $S_n$ below the mean can be obtained by considering the sequence of iid random variables $(-X_i)$.

\begin{proof}[Proof of Theorem \ref{thm:ldp}]
Without loss of generality, we may assume that $\mu = 0$ (if $\mu \neq 0$, we can define $Y_i = X_i - \mu$ and translate the result accordingly).
We begin by proving that $\Lam_X^*(a) > 0$.  We have
\begin{align}
a t - \Lam_X(t)
	&=	\log \Big( \frac{\ee^{at}}{M_X(t)} \Big)
	=		\log \Big( \frac{1 + a t + \Oc(t^2)}{1 + \sig^2 t^2 / 2 + \Oc(t^3)} \Big) , &
	&\text{where}&
\sig^2
	&:= \Vb X .
\end{align}
For $t>0$ sufficiently small we have
\begin{align}
\frac{1 + a t + \Oc(t^2)}{1 + \sig^2 t^2 / 2 + \Oc(t^3)}
	> 1 .
\end{align}
Hence, we have
\begin{align}
\Lam_X^*(a)
	&=	\sup_{t \in \Rb} \{ a t - \Lam_X(t) \}
	\geq \sup_{t \in \Rb_+} \log \Big( \frac{1 + a t + \Oc(t^2)}{1 + \sig^2 t^2 / 2 + \Oc(t^3)} \Big)
	> 0 .
\end{align}
We make two notes for future use.  First, by assumption we have $a > \mu = 0$.  As $\Lam_X$ is convex with $\Lam_X'(0) = \mu = 0$ and $\Lam_X(0) = 0$ it follows that
\begin{align}
\Lam_X^*(a) 
	&= \sup_{t > 0 } \{ a t - \Lam_X(t) \} , &
a
	&> 0 . \label{eq:note1}
\end{align}
To see \eqref{eq:note1}, it may help to draw a picture.  Second, 
\begin{align}
	\text{$\Lam_X$ is strictly convex at points where $\Lam_X''$ exists}. \label{eq:note2}
\end{align}
To see this, note that $\Vb X > 0$ under the hypotheses of the theorem, implying by \eqref{eq:note0} and the Cauchy-Schwarz inequality that $\Lam_X''(t) > 0$.
\\[0.5em]
We now proceed to derive an upper bound for $\Pb( S_n > n a )$.  Using the fact that $\ee^{t S_n} > \ee^{n a t} \Ib_{\{S_n > n a\}}$ for all $t>0$ we obtain
\begin{align}
\Pb( S_n > n a )
	&=	\Eb \Ib_{\{S_n > n a\}}
	\leq \ee^{- n a t} \Eb \ee^{t S_n}
	=	\Big( \ee^{-a t} M_X(t) \Big)^n 
	=	\ee^{-n (a t - \Lam_X(t))} , &
	&\forall \, t > 0 . \label{eq:upper.bound}
\end{align}
Note that \eqref{eq:upper.bound} is valid for all $t > 0$ and thus
\begin{align}
\frac{1}{n} \log \Pb(S_n > na) 
	&\leq - \sup_{t > 0} \{ at - \Lam_X(t) \} 
	=	- \Lam_X^*(a) , \label{eq:ldp.ub}
\end{align}
where the equality follows from \eqref{eq:note1}.
\\[0.5em]
Before obtaining lower bound for $\Pb( S_n > na )$ let us define
\begin{align}
T
	:= \sup \{ t : M_X(t) < \infty \} .
\end{align}
Throughout this proof we will assume that 
\begin{align}
t^*
	&:=	\operatorname*{arg max}_{t} \{ a t - \Lam_X(t) \}
	\in (0,T) . \label{eq:additional.assumption}
\end{align}
This assumption is not necessary, but simplifies the proof considerably.
As $a t - \Lam_X(t)$ has a maximum at $t^*$, and as $\Lam_X$ is $C^\infty$ on $(0,T)$, the derivative of $a t - \Lam_X(t)$ equals zero at $t = t^*$ and hence
\begin{align}
\Lam_X'(t^*)
	&=	a .
\end{align}
Let us denote by $F_X$ the distribution of $X$.  We can define an \emph{exponentially tilted} distribution $F_{\Xt}$ as follows
\begin{align}
F_{\Xt}(x)
	&=	\frac{1}{M_X(t^*)} \int_{-\infty}^x \ee^{t^* y} \dd F_X(y) . \label{eq:F.Ftilde}
\end{align}
Let $(\Xt_i)$ be a sequence of iid random variables with common distribution $F_{\Xt}$ and define $\St_n = \sum_{i=1}^n \Xt_i$.  Observe that
\begin{align}
M_{\Xt}(t)
	&:= \int_\Rb  \ee^{t x} \dd F_{\Xt}(x) 
	= 	\int_\Rb \frac{\ee^{(t+t^*) x}}{M_X(t^*)} \dd F_X(x)
	=		\frac{M_X(t + t^*)}{M_X(t^*)} . 
\end{align}
The first two moments of $\Xt$ can be obtained by direct computation
\begin{align}
\Eb \Xt 
	&=	M_{\Xt} (0) 
	=		\frac{M_X'(t^*)}{M(t^*)} 
	=		\Lam_X'(t^*) = a  , \\
\Vb \Xt
	&=	\Eb \Xt^2 - (\Eb \Xt)^2 
	=		M_{\Xt}''(0) - ( M_{\Xt}'(0) )^2
	=		\Lam_X''(t^*) \in (0,\infty) .
\end{align}
where, the fact that $0 < \Lam_X''(t^*) < \infty$ follow from \eqref{eq:note2} and the assumption \eqref{eq:additional.assumption}.
Noting that $S_n$ and $\St_n$ are sums of iid random variables, we have
\begin{align}
M_{S_n}(t)
	&=	M_X^n(t) , &
M_{\St_n}(t)
	&=	M_{\Xt}^n(t) .
\end{align}
Using the above and \eqref{eq:F.Ftilde}, it is easy to show that $F_{S_n}$ and $F_{\St_n}$ are related as follows
\begin{align}
F_{\St_n}(x)
	&=	\frac{1}{M_X^n(t^*)} \int_{-\infty}^x \ee^{t^* y} \dd F_{S_n}(y) .
\end{align}
Now, let $b > a$.  We have
\begin{align}
\Pb(S_n > na )
	&=	\int_{na}^\infty \dd F_{S_n}(x) \\
	&=	\int_{na}^\infty M_X^n(t^*) \ee^{- t^* x} \dd F_{\St_n}(x) \\
	&\geq M_X^n(t^*) \ee^{- n t^* b}  \int_{na}^{nb} \dd F_{\St_n}(x) \\
	&\geq	\ee^{-n(t^* b - \Lam_X(t^*))} \Pb( na < \St_n < n b ) .
\end{align}
As $\Eb \Xt = a$ and $\Vb \Xt > 0$ we have by the CLT that $\Pb( \St_n > n a ) \to 1/2$ as $n \to \infty$,
and by the SLLN we have that $\Pb( \St_n < nb) \to 1$.  Therefore, we have
\begin{align}
\frac{1}{n} \log \Pb( S_n > na )
	&\geq - (t^* b - \Lam_X(t^*) ) + \frac{1}{n} \log \Pb( na < \St_n < nb ) \\
	&\to - (t^* b - \Lam_X(t^*) ) &
	&\text{as $n \to \infty$} \\
	&\to - (t^* a - \Lam_X(t^*) ) = - \Lam_X^*(a) &
	&\text{as $b \to a$}. \label{eq:ldp.lb}
\end{align}
Equation \eqref{eq:ldp} follows from \eqref{eq:ldp.ub} and \eqref{eq:ldp.lb}.
\end{proof}

\begin{example}
Let $(X_i)$ be a sequence of iid random variables with $\Pb(X_i = 1) = \Pb(X_i = -1) = 1/2$.  We claim that 
\begin{align}
\Pb( S_n > a n)^{1/n} 
	&\to \frac{1}{\sqrt{(1+a)^{1+a}(1-a)^{1-a}}} , &
0
	<a < 1 ,
\end{align}
as $n \to \infty$.  To see this, we have from \eqref{eq:ldp} that
\begin{align}
\Pb(S_n > na)^{1/n}
	&\to \ee^{- \Lam_X^*(a)}
\end{align}
We must compute $\Lam_X^*(a)$.  First, we have
\begin{align}
\Lam_X(t)
	&=	\log \Eb \ee^{t X}
	=		\log \Big( \tfrac{1}{2} \ee^{-t} + \tfrac{1}{2} \ee^t \Big) 
	=		\log \cosh t .
\end{align}
Next, to find $t^* := \operatorname*{arg max} \{ a t - \Lam_X(t) \}$ we compute
\begin{align}
0
	&=	\d_t \Big[ a t - \Lam_X(t) \Big|_{t = t^*}
	=		a - \tanh t^* , &
	&\Rightarrow&
t^*
	&=	\tanh^{-1} a .
\end{align}
Therefore, we have
\begin{align}
\Lam_X^*(a)
	&=	\sup_t \{ a t - \Lam_X(t) \}
	=		a t^* - \Lam_X(t^*)
	=		a \tanh^{-1} a - \log \cosh \tanh^{-1} a .
\end{align}
Lastly, using $2 \tanh^{-1} a = \log (1+a) - \log (1-a)$, we have
\begin{align}
\ee^{- \Lam_X^*(a)}
	&=	\ee^{ - a \tanh^{-1} a + \log \cosh \tanh^{-1} a}
	=		\frac{\ee^{-a \tanh^{-1} a}}{\sqrt{1-a^2}} \\
	&=	\frac{\sqrt{(1+a)^{-a} (1-a)^a}}{\sqrt{(1+a)(1-a)}}
	=		\frac{1}{\sqrt{(1+a)^{1+a}(1-a)^{1-a}}} ,
\end{align}
as claimed.  Note that, if $a \notin (0,1)$, then we have $\Pb( S_n > n a) = 0$ for all $n$.
\end{example}


\section{Exercises}
\label{sec:exercises-3}

\begin{exercise}
\label{exer:3.1}
Let $X \sim \text{Bin}(n,U)$ where $U \sim \Uc((0,1))$.  What is the probability Generating function $G_X(s)$ of $X$?  What is $\Pb( X = k)$ where $k \in \{0, 1, 2, \ldots , n \}$?
\end{exercise}

\begin{exercise}
\label{exer:3.2}
Let $Z_n$ be the size of the $n$th generation in an ordinary branching process with $Z_0 = 1$, $\Eb Z_1 = \mu$ and $\Vb Z_1 > 0$.  Show that $\Eb Z_n Z_m = \mu^{n-m} \Eb Z_m^2$ for $m \leq n$.  Use this to find the correlation coefficient $\rho(Z_m,Z_n)$ in terms of $\mu$, $n$, and $m$.  Consider the case $\mu = 1$ and the case $\mu \neq 1$.
\end{exercise}

\begin{exercise}
\label{exer:3.3}
Consider a branching process with generation sizes $Z_n$ satisfying $Z_0 = 1$ and $\Pb(Z_1 = 0) = 0$.  Pick two individuals at random with replacement from the $n$th generation and let $L$ be the index of the generation which contains their most recent common ancestor.  Show that $\Pb(L=r) = \Eb Z_r^{-1} - \Eb Z_{r+1}^{-1}$ for $0 \leq r < n$.  
\end{exercise}

\begin{exercise}
\label{exer:3.4}
Consider a branching process with immigration
\begin{align}
Z_0
	&=	1 , &
Z_{n+1}
	&=	\sum_{i=1}^{Z_n} X_{n,i} + Y_n 
\end{align}
where the $(X_{n,i})$ are iid with common distribution $X$, the $(Y_n)$ are iid with common distribution $Y$ and the $(X_{n,i})$ and $(Y_n)$ are independent.  What is $G_{Z_{n+1}}(s)$ in terms of $G_{Z_n}(s)$, $G_X(s)$ and $G_Y(s)$?  Write $G_{Z_2}(s)$ explicitly in terms of $G_X(s)$ and $G_Y(s)$.
\end{exercise}

\begin{exercise}
\label{exer:3.5}
Find $\phi_{X^2}(t) := \Eb \ee^{\ii t X^2}$ where $X \sim \Nc(\mu,\sig^2)$.
\end{exercise}

\begin{exercise}
\label{exer:3.6}
Let $X_n$ have cumulative distribution function
\begin{align}
F_{X_n}(x)
	&=	\Big( x - \frac{\sin (2n\pi x)}{2 n\pi} \Big) \Ib_{\{0 \leq x \leq 1 \}} + \Ib_{\{x > 1 \}}
\end{align}
\textbf{(a)} Show that $F_{X_n}$ is a distribution function and find the corresponding density function $f_{X_n}$.
\textbf{(b)}  Show that$F_{X_n}$ converges to the uniform distribution function $F_U$ as $n \to \infty$, but that the density functions $f_{X_n}$ do not converge to $f_U$.  Here, $U \sim \Uc((0,1))$.
\end{exercise}

\begin{exercise}
\label{exer:3.7}
A coin is tossed repeatedly, with heads turning up with probability $p$ on each toss.  Let $N$ be the minimum number of tosses required to obtain $k$ heads.  Show that, as $p \to 0$, the distribution function of $2Np$ converges to that of a gamma distribution.  Note that, if $X \sim \Gam(\lam,r)$ then
\begin{align}
f_X(x)
	&=	\frac{1}{\Gam(r)} \lam^r x^{r-1} \ee^{-\lam x} \Ib_{\{x \geq 0\}} .
\end{align}
\end{exercise}

\begin{exercise}
\label{exer:3.8}
Show that, as $n  \to \infty$, we have
\begin{align}
2^{-n} \sum_{k \in A_n} \binom{n}{k} \to \int_{-x}^x \dd y \frac{1}{\sqrt{2\pi}} \ee^{-y^2/2}
\end{align}
where $A_n := \{ k : |k - n/2| \leq x \sqrt{n} /2 \}$.
\end{exercise}

\begin{exercise}
\label{exer:3.9}
Let $(X_j)_{j \in \Nb}$ be a sequence of i.i.d standard normal $N(0,1)$ random variables.  Define
\begin{align}
Y_n
	&:=	\frac{1}{n} \sum_{j=1}^n \sqrt{j} X_j , &
n
	&\in \Nb .
\end{align}
Show that $Y_n \stackrel{\Dc}{\to} Y$ and identify the distribution of $Y$.  \textit{Hint}: recall that if
\begin{align}
Z 
	&\sim \Nc(\mu,\sig^2) , &
	&\text{then}&
\Eb \ee^{\ii u Z}
	&=	\ee^{\ii u \mu - \tfrac{1}{2} u^2 \sig^2} .
\end{align}
\end{exercise}

\begin{exercise}
\label{exer:3.10}
Fix a probability space $(\Om,\Fc,\Pb)$.  Recall that, if $Z$ is a Poisson random variable with parameter $\lam>0$ and $X$ is a Bernoulli random variable with parameter $p \in (0,1)$ then their probability mass functions are given by 
\begin{align}
f_Z(k)
	&=	\frac{\lam^k \ee^{-\lam}}{k!} , &
k
	&\in \{0, 1, 2, \ldots \}, &
f_X(k)
	&=	\left\{ \begin{aligned}
			&p , & & k = 1 , \\
			&1 - p , & & k = 0 .
			\end{aligned} \right.
\end{align}
(a) Compute the characteristic functions $\phi_Z(t) := \Eb \ee^{\ii t Z}$ and $\phi_X(t) := \Eb \ee^{\ii t X}$.
\\[0.5em]
(b) For every $n \in \Nb$, let $(X_j^{(n)})_{1 \leq j \leq n}$ be a sequence of i.i.d. Bernoulli random variables with parameter $\lam/n$.
Define 
\begin{align}
S_n 
	&:= \sum_{j=1}^n X_j^{(n)} .
\end{align}
Show that $S_n \stackrel{\Dc}{\to} Z$.  The following may be useful
\begin{align}
\ee^x
	&=	\sum_{n=0}^\infty \frac{x^n}{n!} , &
\lim_{n \to \infty} \Big( 1 + \frac{x}{n} \Big)^n
	&=	\ee^x .
\end{align}
\end{exercise}

%
%

\chapter{Discrete time Markov chains}
\label{chapter:dtmc}
The notes from this chapter are primarily taken from \cite[Chapter 6]{gs2001}.


\section{Overview of discrete time Markov chains}
We have already defined what we mean by Markov process in Definition \ref{def:markov}.  In this chapter, we will derive some results that are specific for Markov discrete time processes that have a countable state space. 

\begin{definition}
A \emph{discrete time Markov chain} $X = (X_n)_{n \in \Nb_0}$ is a discrete-time Markov process with a countable state space $S$.
\end{definition}

\noindent
Recall from \eqref{eq:markov.prop} that a Markov process defined on a probability space $(\Om,\Fc,\Pb)$, which is equipped with a filtration $\Fb = (\Fc_n)_{n \in \Nc_0}$, is a process $X = (X_n)_{n \in \Nb_0}$ that satisfies $\Pb( X_{n+k} \in A | \Fc_n) = \Pb( X_{n+k} \in A | X_n)$.  
What the Markov property means in practice is that the evolution of a Markov chain is described by its one-step transition probabilities $\Pb( X_{n+1} = j | X_n = i)$.  In general, such probabilities may depend on the time $n$.  Throughout this Chapter, we will restrict our attention to Markov chains whose transition probabilities do not depend on $n$.

\begin{definition}
A discrete time Markov chain $X$ is called \emph{homogeneous} if 
\begin{align}
\Pb( X_{n+1} = j | X_n = i)
	&=	\Pb( X_{1} = j | X_0 = i) =: p(i,j) &
	&		\forall \, n \in \Nb_0 , \forall \, i,j \in S .
\end{align}
We call the $|S| \times |S|$ matrix $\Pv = ( p(i,j) )$ the \emph{one-step transition matrix}.
\end{definition}

\noindent
The one-step transition matrix $\Pv$ satisfies some obvious properties:
\begin{enumerate*}
\item $\sum_j p(i,j) = 1$.
\item $p(i,j) \geq 0$ for all $i,j \in S$.
\end{enumerate*}

\noindent
Since $\Pb( X_{n+1} = j | X_n = i) = \Pb( X_1 = j | X_0 = i)$, it follows that $\Pb( X_{n+m} = j | X_m = i) = \Pb( X_{n} = j | X_0 = i)$.

\begin{definition}
Let $X$ be a discrete time Markov chain and let 
\begin{align}
p_n(i,j)
	&:=	\Pb( X_{n} = j | X_0 = i)
	&		\forall \, n \in \Nb_0 , \forall \, i,j \in S .
\end{align}
We call the $|S| \times |S|$ matrix $\Pv_n = ( p_n(i,j) )$ the \emph{$n$-step transition matrix}.
\end{definition}

\noindent
From the one-step transition matrix $\Pv$, one can easily derive the $n$-step transition matrix $\Pv_n$.

\begin{theorem}[Chapman-Kolmogorov Equation]
\label{thm:chapman-kolmogorov}
Let $\Pv$ and $\Pv_n$ be the one-step and $n$-step transition matrices of a homogeneous discrete-time Markov chain.  Then
\begin{align}
\Pv_{m+n} 
	&= \Pv_m\Pv_n , &
\Pv_n
	&=	\Pv^n ,
\end{align}
where $\Pv^n$ is the $n$th power of $\Pv$.
\end{theorem}
\begin{proof}
For any $i,j \in S$ we have
\begin{align}
p_{m+n}(i,j)
	&=	\Pb( X_{m+n} = j | X_0 = i) \\
	&=	\sum_{k} \Pb( X_{n+m} = j | X_m = k) \Pb( X_{m} = k | X_0 = i) \\
	&=	\sum_{k} \Pb( X_{n} = j | X_0 = k) \Pb( X_{m} = k | X_0 = i)	
	=		\sum_{k} p_m(i,k) p_n(k,j) 
\end{align}
Thus, we have established that $\Pv_{m+n} = \Pv_m\Pv_n$.  Next, note that
\begin{align}
\Pv_n
	&=	\Pv_{n-1}\Pv
	=		\Pv_{n-2} \Pv^2
	=		\ldots
	=		\Pv^n ,
\end{align}
as claimed.
\end{proof}

\begin{lemma}
Let $X$ be a homogeneous discrete time Markov chain.
Denote the probability mass function of $X_n$ by
\begin{align}
\mu_n(i)
	&:=	\Pb( X_n = i ) , &
\muv_n
	&=	( \mu_n(1), \mu_n(2), \ldots , \mu_n(|S|) ) .
\end{align}
Then $\muv_{n+m} = \muv_{n} \Pv_m$ and thus $\muv_m = \muv_0 \Pv^m$.
\end{lemma}
\begin{proof}
For any $j \in S$ we have
\begin{align}
\mu_{n+m}(j)
	&=	\Pb( X_{n+m} = j ) 
	=		\sum_i \Pb(X_{n+m} = j | X_n = i) \Pb(X_n = i)
	=		\sum_i \mu_n(i) p_m(i,j) .
\end{align}
Thus, we have established that $\muv_{n+m} = \muv_{n} \Pv_m$.  To finish the proof, take $n=0$ and use $\Pv_m = \Pv^m$.
\end{proof}

\begin{example}
Consider a random walk on a circle with $n$ nodes.  At each step a particle moves clockwise with probability $p$ and counter-clockwise with probability $q$.  Let us write the one-step transition matrix.  We have
\begin{align}
\Pv
	&=	\( \begin{array}{cccccccc}
			0	& p & 0 & 0 & \ldots & 0 & 0 & q \\
			q	& 0 & p & 0 & \ldots & 0 & 0 & 0 \\
			0	& q & 0 & p & \ldots & 0 & 0 & 0 \\
			\vdots & \vdots & \vdots & \vdots & \ddots & \vdots & \vdots & \vdots \\
			0	& 0 & 0 & 0 & \ldots & q & 0 & p \\
			p	& 0 & 0 & 0 & \ldots & 0 & q & 0 
			\end{array} \) .
\end{align}
\end{example}

\begin{example}
Consider a random walk on $\Nb_0$.  At each step a particle moves to the up one unit with probability $p$ and returns to the origin with probability $q$.  Let us write the one-step transition matrix.  We have
\begin{align}
\Pv
	&=	\( \begin{array}{cccccc}
			q & p & 0 & 0 & 0 & \ldots \\
			q & 0 & p & 0 & 0 & \ldots \\
			q & 0 & 0 & p & 0 & \ldots \\
			\vdots  & \vdots & \vdots & \vdots & \vdots & \ddots 
			\end{array} \) .
\end{align}
\end{example}


\section{Classification of states}
\label{sec:state.classification}

\begin{definition}
\label{def:persistent}
Let $X$ be a Markov chain.  We say that state $i$ is \emph{persistent} or \emph{recurrent} if
\begin{align}
\Pb( X_n = i \text{ for some } n | X_0 = i )
	&=	1 .
\end{align}
Otherwise, we say that state $i$ is \emph{transient}.
\end{definition}

\noindent
If a Markov chain visits a persistent state once, then it is guaranteed to return to that state.  In fact, it is guaranteed to return to that state \emph{infinitely often}.
\\[0.5em]
We would like to find conditions on the transition matrix $\Pv$ that enable us to classify a state as either persistent or transient.

\begin{definition}
\label{def:hitting.time}
Let $X$ be a Markov chain.  We define the \emph{first passage to state $j$} by
\begin{align}
\tau_j
	&:=	\inf\{ n \geq 1 : X_n = j \}
\end{align}
and we denote by $f_{ij}$ the probability mass function of $\tau_j$ given $X_0 = i$.  That is
\begin{align}
f_{ij}(n)
	&=	\Pb( \tau_j = n | X_0 = i ).
\end{align}
\end{definition}

\begin{theorem}
Define the following generating functions
\begin{align}
P_{ij}(s) 
	&:= \sum_{n=0}^\infty s^n p_n(i,j) , & 
p_0(i,j)
	&:=	\del_{i,j} , &
F_{ij}(s)
	&:=	\sum_{n=0}^\infty s^n f_{ij}(n) , &
f_{ij}(0)
	&:=	0 .
\end{align}
Then, for any $|s|<1$, we have
\begin{align}
P_{ii}(s)
	&=	1 + F_{ii}(s)P_{ii}(s) , &
P_{ij}(s)
	&=	F_{ij}(s)P_{jj}(s) , &
i
	&\neq j . \label{eq:Pij=FijPjj}
\end{align}
\end{theorem}

\noindent
Note that we must restrict $|s|<1$ because $P$ is not a \emph{probability} generating function.  When we need $s=1$ we can take a limit as $s \nearrow 1$ and use \emph{Abel's Theorem}. which states that
\begin{align}
\lim_{s \nearrow 1} P_{ij}(s)
	&=	\sum_{n=0}^\infty p_n(i,j) . \label{eq:Abel}
\end{align}

\begin{proof}
Fix $i,j \in S$ and define
\begin{align}
A_m
	&=	\{ X_m = j \} , &
B_m
	&=	\{ \tau_j = m \} .
\end{align}
Since the sets $(B_m)_{m \geq 1}$ are disjoint we have
\begin{align}
p_m(i,j)
	&=	\Pb( A_m | X_0 = i )
	 =	\sum_{k=1}^m \Pb( A_m \cap B_k | X_0 = i) . \label{eq:P.Am}
\end{align}
Note that
\begin{align}
\Pb( A_m \cap B_k | X_0 = i)
	&=	\Pb( A_m | B_k , X_0 = i ) \cdot \Pb(B_k | X_0 = i )  \\
	&=	\Pb( A_m | X_k = j ) \cdot \Pb( B_k | X_0 = i ) 
	=		p_{m-k}(j,j) \cdot f_{ij}(k)  . \label{eq:P.Am.Bk}
\end{align}
Inserting  \eqref{eq:P.Am.Bk} into \eqref{eq:P.Am} we obtain
\begin{align}
p_m(i,j)
	&=	\sum_{k=1}^m f_{ij}(k) p_{m-k}(j,j) .
\end{align}
Thus
\begin{align}
P_{ij}(s)
	&=	\sum_{m=0}^\infty s^m p_m(i,j) \\
	&=	\del_{ij} + \sum_{m=1}^\infty s^m \sum_{k=1}^m f_{ij}(k) p_{m-k}(j,j) \\
	&=	\del_{ij} + \sum_{m=1}^\infty \sum_{k=1}^m s^k f_{ij}(k) s^{m-k} p_{m-k}(j,j) \\
	&=	\del_{ij} + \sum_{k=1}^\infty s^k f_{ij}(k)  \cdot \sum_{m=k}^\infty  s^{m-k} p_{m-k}(j,j) \\
	&=	\del_{ij} + F_{ij}(s)P_{jj}(s) ,
\end{align}
which is precisely the claimed result.
\end{proof}

\begin{corollary}
\label{cor:trans.pers.cond}
The following holds:
\begin{enumerate*}
\item State $j$ is persistent if $\sum_n p_n(j,j) = \infty$, and if this holds then $\sum_n p_n(i,j) = \infty$ for all $i$ such that $F_{ij}(1)>0$.
\item State $j$ is transient if $\sum_n p_n(j,j) < \infty$, and if this holds then $\sum_n p_n(i,j) < \infty$ for all $i$.
\end{enumerate*}
\end{corollary}

\begin{proof}
We will show that $j$ is persistent if and only if $\sum_n p_n(j,j) = \infty$.
First, by Definition \ref{def:persistent}, we have
\begin{align}
	&\text{$j$ is persistent}&
	&\Leftrightarrow&
F_{jj}(1)
	&=	\sum_{n=1}^\infty f_{jj}(n) = 1.
\end{align}
Next, from \eqref{eq:Pij=FijPjj} and Abel's Theorem \eqref{eq:Abel}, we have
\begin{align}
\sum_n p_n(j,j)
	&=	\lim_{s \nearrow 1} P_{jj}(s)
	=		\lim_{s \nearrow 1}\frac{1}{1-F_{jj}(s)} 
	=		\frac{1}{1-F_{jj}(1)} .
\end{align}
Thus, $j$ is persistent if and only if $\sum_n p_n(j,j)=\infty$.
We also have from \eqref{eq:Pij=FijPjj} and Abel's Theorem \eqref{eq:Abel} that
\begin{align}
\sum_n p_n(j,j)
	&=	\lim_{s \nearrow 1} P_{jj}(s)
	=		\frac{ \lim_{s \nearrow 1}  P_{ij}(s)}{F_{ij}(1)} 
	=		\frac{ \sum_n p_n(i,j)}{F_{ij}(1)} .
\end{align} 
Thus, if $j$ is persistent then $\sum_n p_n(i,j) = \infty$ whenever $F_{ij}(1) > 0$.
\end{proof}

\begin{corollary}
\label{cor:pnij.to.0}
If $j$ is transient then $p_n(i,j) \to 0$ as $n \to \infty$ for all $i$.
\end{corollary}

\begin{proof}
From Corollary \ref{cor:trans.pers.cond}, if $j$ is transient, we have $\sum_n p_n(i,j) < \infty$ for all $i$.  It follows that $p_n(i,j) \to 0$ as $n \to \infty$ for all $i$.
\end{proof}

\noindent
There are a few other state classifications that are important

\begin{definition}
\label{def:tau.bar}
The \emph{mean recurrence time} of a state $j$ is defined as
\begin{align}
\taub_j
	&=	\Eb[ \tau_j | X_0 = j ] ,
\end{align}
where $\tau_j$ is defined in Definition \ref{def:hitting.time}.
\end{definition}

\noindent
Clearly, $\taub_j = \infty$ if a state is transient.  However, $\taub_j$ may be infinite even if $j$ is recurrent.

\begin{definition}
\label{def:null}
A recurrent state $j$ is said to be \emph{null} if $\taub_j = \infty$ and \emph{non-null} or \emph{positive} if $\taub_j < \infty$.
\end{definition}

\noindent
This is a simple condition which differentiates null recurrent from positive recurrent states.

\begin{theorem}
\label{thm:null}
A persistent state $j$ is null if and only if $p_n(j,j) \to 0$ as $n \to \infty$; if this holds then $p_n(i,j) \to 0$ as $n \to \infty$ for all $i$.
\end{theorem}

\noindent
We will not prove of Theorem \ref{thm:null}.

\begin{definition}
The \emph{period} of a state $i$, denoted $d(i)$, is defined by
$d(i) = \text{gcd}\{ n : p_n(i,i)> 0 \}$.
If $d(i) = 1$ we say that state $i$ is \emph{aperiodic}.
Here $\text{gcd}$ means ``greatest common divisor.'' 
\end{definition}

\begin{example}
Let $X$ be a Markov Chain with one-step transition matrix
\begin{align}
\Pv
	&=	\( \begin{array}{ccc}
			0 & 1 & 0 \\
			0 & 0 & 1 \\
			1 & 0 & 0
			\end{array} \) .
\end{align}
One can easily see that
\begin{align}
p_n(i,i)
	&=	\left\{ \begin{array}{ll}
			1 & n=3, 6, 9, \ldots \\
			0 & n \neq 3, 6, 9, \ldots
			\end{array} \right.
\end{align}
Since $\text{gcd}\{3,6,9,12,\ldots\} = 3$, we have $d(i) = 3$ for all $i$.
\end{example}

\begin{definition}
A state is called \emph{ergodic} if it is persistent, non-null, and aperiodic.
\end{definition}


\section{Classification of chains}
\label{sec:chain.classification}

\begin{definition}
We say state $i$ \emph{communicates} with state $j$, written $i \to j$, if $p_n(i,j)>0$ for some $n \geq 1$ (in other words, it is possible to reach state $j$ from state $i$).  We say that states $i$ and $j$ \emph{intercommunicate}, written $i \leftrightarrow j$, if $i \to j$ and $j \to i$.
\end{definition}

\begin{theorem}
\label{thm:equivalence}
Suppose $i \leftrightarrow j$.  Then 
\begin{enumerate*}
\item $i$ is transient if an only if $j$ is transient; \label{itm:transient}
\item $i$ is null persistent if and only if $j$ is null persistent. \label{itm:null.persistent}
\item $i$ and $j$ have the same period: $d(i) = d(j)$.
\end{enumerate*}
\end{theorem}

\begin{proof}
We will prove Item \ref{itm:transient}.
If $i \leftrightarrow j$ then there exists $k,n$ such that $\alpha := p_k(i,j) p_n(j,i) >0$.  Therefore, for any $m$, we have
\begin{align}
p_{k+m+n}(i,i)
	&\geq p_k(i,j) p_m(j,j) p_n(j,i)
	=	\alpha p_m(j,j) .
\end{align}
Therefore, if $i$ is transient. by Corollary \ref{cor:trans.pers.cond}, we have
\begin{align}
\infty
	&>	\sum_m p_m(i,i)
	\geq \sum_m p_{k+m+n}(i,i)
	\geq	\alpha \sum_m p_m(j,j) ,
\end{align}
and thus $j$ is transient is well.  To show the converse simply switch $i$ and $j$ in the above argument.
\end{proof}

\begin{definition}
Let $C \subset S$.  We say that $C$ is \emph{closed} if $p(i,j) = 0$ for all $i \in C$ and $j \notin C$.
We say that $C$ is \emph{irreducible} if $i \leftrightarrow j$ for all $i,j \in C$.
\end{definition}

\noindent
If a closed set $C$ consists of a single state (e.g., $C = \{j\}$), then we call this state \emph{absorbing}.

\begin{theorem}[Markov chain decomposition]
\label{thm:mc.decomp}
The state space $S$ of a Markov chain can be uniquely partitioned as follows
\begin{align}
S
	&=	T \cup C_1 \cup C_2 \cup \ldots ,
\end{align}
where $T$ is the set of transient states and $C_1, C_2, \ldots$ are closed, irreducible sets of persistent states.
\end{theorem}

\noindent
The importance of the decomposition Theorem is as follows.  Suppose $X_0 \in C_k$, where $C_k$ is a set of closed persistent states.  Then $X$ remains in $C_k$ forever; as such, we may as well analyze only $C_k$ rather than the entire states space $S$.  On the other hand, if $X_0 \in T$, then the chain $X$ either remains in $T$ forever or it eventually visits one of the closed persistent states $C_k$, where it remains forever.  If $|S|<\infty$, it turns out that it is impossible for $X$ to remain in $T$ forever.

\begin{lemma}
\label{lem:persistent.non-null}
If $S$ is finite, then at least one state is persistent and all persistent states are non-null.
\end{lemma}

\begin{proof}
First, suppose all state are transient.  Then, by Corollary \ref{cor:pnij.to.0}, we have
\begin{align}
0
	&=	\sum_j \lim_{n \to \infty} p_n(i,j)
	=		\lim_{n \to \infty} \sum_j p_n(i,j)
	=		1 .
\end{align}
Thus, we have a contradiction and we conclude that at least one state is persistent.  Now, suppose $j \in C_k$, where $C_k$ is a closed set of persistent states.  Suppose $j$
is null.  Then by Theorem \ref{thm:equivalence}, all states $i \in C_k$ are null.  Then, by Theorem \ref{thm:null} we have
\begin{align}
0
	&=	\sum_j \lim_{n \to \infty} p_n(i,j)
	=		\lim_{n \to \infty}  \sum_j p_n(i,j)
	=		1 .
\end{align}
Thus, $j$ cannot be null.
\end{proof}

\begin{example}
Let $S = \{ 1, 2, 3, 4, 5, 6 \}$ and
\begin{align}
\Pv
	&=	\( \begin{array}{cccccc}
			p & q & 0 & 0 & 0 & 0 \\
			q & p & 0 & 0 & 0 & 0 \\
			r & s & t & u & 0 & 0 \\
			0 & 0 & u & t & s & r \\
			0 & 0 & 0 & 0 & p & q \\
			0 & 0 & 0 & 0 & q & p \\
			\end{array} \) , &
q+p
	&= 1 , &
r+s+t+u
	&=	1 .
\end{align}
The states subsets $C_1=\{1,2\}$ and $C_2 = \{5,6\}$ are closed, persistent, non-null sets, since, once the chain visits these sets, it cannot escape.  The subset $T=\{3,4\}$ is transient.  If $X_0 \in T$, it will eventually move to either $C_1$ or $C_2$, where it will remain forever.
\end{example}


\section{Stationary distributions and the limit theorem}
\label{sec:stationary}
In this section, we will be interested in studying the long-run behavior of a Markov chain $X$.

\begin{definition}
Let $X$ be a Markov chain with one-step transition matrix $\Pv$.  We say that a row vector
\begin{align}
\piv
	&=	( \pi(1), \pi(2), \ldots , \pi(|S|) ) ,
\end{align}
is a \emph{stationary} or \emph{invariant distribution} of $X$ if it satisfies
\begin{align}
\piv 
	&=	\piv \Pv , &
\sum_i \pi(i)
	&=	1 , 
\end{align}
and $\pi(i) \geq 0$ for all $i$
\end{definition}

\noindent
The vector $\piv$ is called a stationary distribution for the following reason:
suppose the Markov chain $X$ has an initial distribution $\muv_0 = \piv$.  Then we have
\begin{align}
\muv_n
	&=	\piv \Pv^n 
	=		\piv \Pv \Pv^{n-1} \\
	&=	\piv \Pv^{n-1} 
	=		\piv \Pv \Pv^{n-2} \\
	&=	\ldots \\
	&=	\piv \Pv
	=		\piv .
\end{align}
Thus, if $\muv_0 = \piv$ then $\muv_n = \piv$ for all $n$.

\begin{example}
Consider a two-state Markov chain with state space $S = \{ 1, 2 \}$ and with one-step transition matrix
\begin{align}
\Pv
	&=	\( \begin{array}{cc}
			1-p & p \\
			q & 1-q
			\end{array} \) , &
q, p
	&\in (0,1) .
\end{align}
Let us find $\piv$.  From $\piv = \piv \Pv$ we derive
\begin{align}
\pi(1)
	&=	\pi(1) (1-p) + \pi(2) q , \\
\pi(2)
	&=	\pi(1) p + \pi(2) (1-q) .
\end{align}
We also have $\pi(1) + \pi(2) = 1$.  Solving for $\pi(1)$ and $\pi(2)$, we obtain
\begin{align}
\pi(1)
	&=	\frac{q}{p+q} , &
\pi(2)
	&=	\frac{p}{p+q} .
\end{align}
\end{example}

\noindent
For a given Markov chain $X$, a stationary distribution may not exist.  And, if it does exist, it may not be unique.

\begin{theorem}
\label{thm:pi.tau.bar}
An irreducible chain $X$ has a stationary distribution $\piv$ if and only if all the states are non-null persistent; in this case $\piv$ is unique and is given by
\begin{align}
\pi(i)
	&=	1/ \taub_i , &
	&		\forall \, i \in S ,
\end{align}
where $\taub_i$ is given in Definition \ref{def:tau.bar}.
\end{theorem}

\noindent
It is typically \emph{much} easier to compute $\piv$ rather than the recurrence times $\taub_i$.  Theorem \ref{thm:pi.tau.bar} gives us a method of computing $\taub_i$ from $\pi(i)$.

\begin{example}
Consider a Markov chain with state space $S = \{ 0, 1, 2, \ldots \}$ and with one-step transition matrix
\begin{align}
\Pv
	&=	\( \begin{array}{cccccc}
			q & p & 0 & 0 & 0 & \ldots \\
			q & 0 & p & 0 & 0 & \ldots \\
			0 & q & 0 & p & 0 & \ldots \\
			0 & 0 & q & 0 & p & \ldots \\
			\vdots & \vdots & \vdots & \vdots & \vdots & \ddots \\
			\end{array} \) , &
p + q
	&=	1 .
\end{align}
Note that the chain is irreducible.  Thus, if an invariant distribution $\piv$ exists, it is unique.
Let us find $\piv$ (if it exists).  From $\piv = \piv \Pv$ we have
\begin{align}
\pi(0)
	&=	q \pi(0) + q \pi(1) , &
	&\Rightarrow&
\pi(1)
	&=	\frac{1-q}{q} \pi(0) = \frac{p}{q} \pi(0) , \\
\pi(1)
	&=	p \pi(0) + q \pi(2) , &
	&\Rightarrow&
\pi(2)
	&=	\frac{1}{q} \Big( \pi(1) - p \pi(0) \Big)
	=		\( \frac{p}{q} \)^2 \pi(0) , \\
\pi(2)
	&=	p \pi(1) + q \pi(3) , &
	&\Rightarrow&
\pi(3)
	&=	\frac{1}{q} \Big( \pi(2) - p \pi(1) \Big)
	=		\( \frac{p}{q} \)^3 \pi(0) . 
\end{align}
This suggests that $\pi(n) = (p/q)^n \pi(0)$.  Assume this is true for every $i \leq n$.  We show it holds for $n+1$ (and thus, for every $n$).  We have
\begin{align}
\pi(n)
	&=	p \pi(n-1) + q \pi(n+1) , &
	&\Rightarrow&
\pi(n+1)
	&=	\frac{1}{q} \Big( \pi(n) - p \pi(n-1) \Big) \\ & & & &
	&=	\frac{1}{q} \( \( \frac{p}{q} \)^n \pi(0) - p \( \frac{p}{q} \)^{n-1} \pi(0) \) \\ & & & &
	&=	\( \frac{p}{q} \)^{n+1} \pi(0) .
\end{align}
Next, we use $\sum_n \pi(n) = 1$ to find $\pi(0)$.  We have
\begin{align}
1
	&=	\sum_{n=0}^\infty \pi(n)
	=		\pi(0) \sum_{n=0}^\infty (p/q)^n
	=		\pi(0) \frac{1}{1-(p/q)} , &
	&		\text{if $(p/q) < 1$} .
\end{align}
Thus, if $(p/q) < 1$ then $\pi(0) = 1-(p/q)$.  Finally, let us find the mean recurrence time of the $n$th state.  We have
\begin{align}
\taub_n
	&=	\frac{1 }{ \pi(n) } = \frac{1}{\pi(0)} \( \frac{q}{p} \)^n
	=		\frac{1}{1-(p/q)} \( \frac{q}{p} \)^n .
\end{align}
Note that if $(p/q) > 1$, then there is no invariant distribution $\piv$.  In this case, the chain is transient and the mean recurrence time for every state is infinite.
What happens if $p=q$?  It turns out that in this case, all states are null recurrent \blu{(need to check)}.
\end{example}


\section{Reversibility}
\label{sec:reversibility}

\begin{definition}
\label{def:Y.markov}
Let $X=(X_n)_{0 \leq n \leq N}$ be a Markov chain.  The \emph{time reversal of $X$} is the process $Y=(Y_n)_{0 \leq n \leq N}$, where
\begin{align}
Y_n
	&=	X_{N-n} .
\end{align}
\end{definition}

\begin{theorem}
\label{thm:Y.markov}
Let $X=(X_n)_{0 \leq n \leq N}$ be an irreducible Markov chain with a one-step transition matrix $\Pv$ and invariant distribution $\piv$.
Suppose $X_0 \sim \piv$ (so that $\muv_n = \piv$ for every $n$).
Then $Y$, the time reversal of $X$, is a Markov chain and its one-step transition matrix, denoted $\Qv = (q(i,j))_{ij}$, is given by
\begin{align}
q(i,j)
	&=	\frac{\pi(j)}{\pi(i)}p(j,i) . \label{eq:q.ij}
\end{align}
\end{theorem}

\begin{proof}
Let $\Fc_n = \sig( Y_k , 0 \leq k \leq n)$.
We cannot assume, a priori that $Y$ is Markov.  Thus, in computing $\Pb( Y_{n+1} = i_{n+1} | \Fc_n )$ we must condition on the entire history $\Fc_n$ of $Y$ rather than conditioning on $Y_n$ only.  We have
\begin{align}
&\Pb( Y_{n+1} = i_{n+1} | Y_k = i_k, 0 \leq k \leq n ) \\
	&=	\frac{ \Pb(  Y_k = i_k, 0 \leq k \leq n+1 ) }{ \Pb(  Y_k = i_k, 0 \leq k \leq n ) } \\
	&=	\frac{ \Pb(  X_{N-k} = i_k, 0 \leq k \leq n+1 ) }{ \Pb(  X_{N-k} = i_k, 0 \leq k \leq n ) } \\
	&=	\frac	{ \Pb(  X_{N-k} = i_k, 0 \leq k \leq n-1 | X_{N-n} = i_n) \Pb(X_{N-n} = i_n | X_{N-(n+1)} = i_{n+1}) \Pb(X_{N-(n+1)} = i_{n+1} )}
						{ \Pb(  X_{N-k} = i_k, 0 \leq k \leq n-1 | X_{N-n} = i_n) \Pb(X_{N-n} = i_n ) } \\
	&=	\frac	{ \Pb(X_{N-n} = i_n | X_{N-(n+1)} = i_{n+1}) \Pb(X_{N-(n+1)} = i_{n+1} )}
						{ \Pb(X_{N-n} = i_n ) } \\
	&=	\frac	{ p(i_{n+1},i_n) \pi(i_{n+1})}
						{ \pi(i_n) } .
\end{align}
Thus, we have show that $Y$ is Markov and its one-step transition matrix is defined by \eqref{eq:q.ij}.
\end{proof}

\begin{definition}
\label{def:reversible}
Let $X=(X_n)_{n \geq 0}$ be an irreducible Markov chain with a one-step transition matrix $\Pv$ and invariant distribution $\piv$.
Suppose $X_0 \sim \piv$ (so that $\muv_n = \piv$ for every $n$).
Let $Y$ be the time reversal of $X$ and denote by $\Qv$ the one-step transition matrix of $Y$.
We say that $X$ is \emph{reversible} if $\Pv = \Qv$, or equivalently (by Theorem \ref{thm:Y.markov}), if
\begin{align}
\pi(i) p(i,j)
	&=	\pi(j) p(j,i) , &
	&		\forall \, i,j \in S . 
\end{align}
\end{definition}

\begin{definition}
\label{def:detailed.balance}
Let $\Pv$ be an $|S|\times|S|$ one step transition matrix.
We say that a distribution $\lamv = (\lam(1), \lam(2), \ldots, \lam(|S|))$ is in \emph{detailed balance} with $\Pv$ if
\begin{align}
\lam(i) p(i,j)
	&=	\lam(j) p(j,i) , &
	&		\forall \, i,j \in S . \label{eq:detailed.balance}
\end{align}
\end{definition}


\begin{theorem}
Let $\Pv$ be the one step transition matrix of an irreducible Markov chain $X$.
Suppose there exists a distribution $\lamv$ such that \eqref{eq:detailed.balance} holds.
Then $\lamv$ is a stationary distribution for $X$.  That is, $\lamv = \lamv \Pv$.
\end{theorem}

\begin{proof}
Assuming \eqref{eq:detailed.balance} holds, we have
\begin{align}
\sum_i \lam(i) p(i,j)
	&=	\sum_i \lam(j) p(j,i)
	=		\lam(j) \sum_i p(j,i)
	=		\lam(j) .
\end{align}
Thus, we have proved $\lamv = \lamv \Pv$, and thus $\lamv$ is a stationary distribution for $X$.
\end{proof}


\section{Chains with finitely many states}
\label{sec:finite.chains}

In this section, we will give some important results that apply to Markov chains with finitely many states.
By Theorem \ref{thm:equivalence} and Lemma \ref{lem:persistent.non-null}, if the state space $S$ is finite and irreducible then all states are necessarily non-null persistent and have the same period.

\begin{theorem}[Perron-Frobenius]
If $\Pv$ is a one-step transition matrix of a finite irreducible Markov chain $X$ with period $d$, then
\begin{enumerate*}
\item $\lam_1 = 1$ is an eigenvalue of $\Pv$.
\item The $d$ complex roots of unity
\begin{align}
\lam_n
	&=	\om^{n-1} , &
\om 
	&= \ee^{2 \pi \ii / d} , &
n
	&=	1, 2, \ldots, d, 
\end{align}
are eigenvalues of $\Pv$.
\item The remaining eigenvalues $\lam_{d+1}, \lam_{d+2}, \ldots , \lam_{|S|}$, satisfy $|\lam_j| < 1.$
\end{enumerate*} 
\end{theorem}

\noindent
Suppose the eigenvalues of $\Pv$ are distinct.  Then it is well-known that $\Pv$ has the decomposition
\begin{align}
\Pv
	&=	\Uv^{-1} \Lamv \Uv , &
\Uv
	&=	\( \begin{array}{ccc}
			{} & \ev_1 & {} \\
			{} & \ev_2 & {} \\
			{} & \vdots & {} \\
			{} & \ev_{|S|} 
			\end{array} \) , &
\Lamv
	&=	\( \begin{array}{cccc}
			\lam_1 & 0 & \ldots & 0 \\
			0 & \lam_2 & \ldots & 0 \\
			\vdots & \vdots & \ddots & \ldots \\
			0 & 0 & \ldots & \lam_{|S|}
			\end{array} \) ,
\end{align}
where $(\ev_i)_{1 \leq i \leq |S|}$ are the \emph{left} eigenvectors of $\Pv$.
\begin{align}
\Pv^n
	&=	\underbrace{\Uv^{-1} \Lamv \Uv \cdot \Uv^{-1} \Lamv \Uv \cdot \ldots \cdot \Uv^{-1} \Lamv \Uv }_{\text{$n$ times}}
	=	\Uv^{-1} \Lamv^n \Uv , \label{eq:Pn}
\end{align}
where we have used the fact that $\Uv^{-1} \Uv = \Iv$.
\\[0.5em]
Expression \eqref{eq:Pn} allows us to study the long-run behavior $n \to \infty$ of a Markov chain $X$.  
In what follows, assume the period of $X$ is one so that the eigenvalues of $\Pv$ satisfy
$1 = \lam_1 > |\lam_2| > |\lam_3| > \ldots > |\lam_{|S|}|$.
Let $X_0 = \muv_0$.  As the eigenvectors $(\ev_i)_{1 \leq i \leq |S|}$ form a basis for $\Rb^{|S|}$ we can express $\muv_0$ as
\begin{align}
\muv_0
	&=	\sum_{i=1}^{|S|} c_i \ev_i , \label{eq:mu.0}
\end{align}
for some constants $(c_i)$.  Next, using \eqref{eq:Pn} and \eqref{eq:mu.0}, we obtain
\begin{align}
\muv_n
	&=	\muv_0 \Pv^n 
	=	\sum_{i=1}^{|S|} c_i \ev_i \Uv^{-1} \Lamv^n \Uv
	=	\sum_{i=1}^{|S|} c_i \lam_i^n \ev_i .
	=	 c_1 \alpha \piv +  \sum_{i=2}^{|S|} c_i \lam_i^n \ev_i ,
\end{align}
where, in the last step, we used $\lam_1 = 1$ and $\ev_1 = \alpha \piv$ for some constant $\alpha$.  As $n \to \infty$, the terms in the sum go to zero at least as fast as $|\lam_2|^n = \exp( n \log |\lam_2| )$.  Thus, $\muv_n \to \piv$ exponentially fast with a rate determined by the second eigenvalue.


\section{Exercises}
\label{sec:exercises-4}

\begin{exercise}
\label{exer:4.1}
A six-sided die is rolled repeatedly.  Which of the following a Markov chains?  For those that are, find the one-step transition matrix.
(a) $X_n$ is the largest number rolled up to the $n$th roll.
(b) $X-n$ is the number of sixes rolled in the first $n$ rolls.
(c) At time $n$, $X_n$ is the time since the last six was rolled.
(d) At time $n$, $X_n$ is the time until the next six is rolled.
\end{exercise}

\begin{exercise}
\label{exer:4.2}
Let $Y_n = X_{2n}$.  Compute the transition matrix for $Y$ when (a) $X$ is a simple random walk (i.e., $X$ increases by one with probability $p$ and decreases by $1$ with probability $q$ and (b) $X$ is a branching process where $G$ is the generating function of the number of offspring from each individual.
\end{exercise}

\begin{exercise}
\label{exer:4.3}
Let $X$ be a Markov chain with state space $S$ and absorbing state $k$ (i.e., $p(k,j)=0$ for all $j \in S$).  Suppose $j \to k$ for all $j \in S$.  Show that all states other than $k$ are transient.
\end{exercise}

\begin{exercise}
\label{exer:4.4}
Suppose two distinct states $i,j$ satisfy
\begin{align}
\Pb(\tau_j < \tau_i | X_0 = i)
	&=	\Pb(\tau_i < \tau_j | X_0 = j) ,
\end{align}
where $\tau_k := \inf\{ n \geq 1 : X_n = j\}$.  Show that, if $X_0 = i$, the expected number of visits to $j$ prior to re-visiting $i$ is one.
\end{exercise}

\begin{exercise}
\label{exer:4.5}
Let $X$ be a Markov chain with transition matrix 
\begin{align}
\Pv
	&=	\( \begin{array}{ccc}
			1-2p & 2p & 0 \\
			p & 1-2p & p \\
			0 & 2p & 1-2p
			\end{array} \) , &
p
	&\in (0,1) .
\end{align}
Find $\Pv^n$, the invariant distribution $\piv$ and the mean-recurrence times $\taub_j$ for $j=1,2,3$.
\end{exercise}

\begin{exercise}
\label{exer:4.6}
Let $X_n$ be the number of mistakes in the $n$th addition of a book.  Between the $n$th and the $(n+1)$th addition and editor corrects each mistake independently with probability $p$ and introduces $Y_n$ new mistakes where the $(Y_n)$ are iid and Poisson distributed with parameter $\lam$.  Find the invariant distribution $\piv$ of the number of mistakes in the book.
\end{exercise}

\begin{exercise}
\label{exer:4.7}
Give an example of a transition matrix $\Pv$ that admits multiple stationary distributions $\piv$.
\end{exercise}

\begin{exercise}
\label{exer:4.8}
A Markov chain on $S = \{0,1,2,\ldots,n\}$ has transition probabilities
$p(0,0)=1-\lam_0$,
$p(i,i+1)=\lam_i$ and $p(i+1,i)=\mu_{i+1}$ for $i = 0, 1, \ldots n-1$,
and
$p(n,n)=1-\mu_n$.
Show that the process is reversible in equilibrium.
\end{exercise}

%
%

\chapter{Continuous time Markov chains}
\label{chapter:ctmc}
The notes from this chapter are primarily taken from \cite[Chapter 6]{gs2001}.


\section{The Poisson process}
\label{sec:poisson}

Many stochastic process change values at random times rather than at fixed intervals.  For example, the value of a stock does not change once per second every second.  Rather, a stock price changes when there is a buy or sell order.  And these orders arrive at random times.  Thus, it makes no sense to model a stock price as a discrete time Markov chain; rather, we must move to continuous time.  The simplest continuous time process is a \emph{counting process}.  

\begin{definition}
A \emph{counting process} is a stochastic process $N = (N_t)_{t \geq 0}$ taking values in $S = \{ 0, 1, 2, \ldots \}$ such that the following hold:
\begin{enumerate*}
\item $N_0 = 0$.
\item If $s<t$ then $N_s \leq N_t$.
\end{enumerate*}
\end{definition}

\noindent
As its name suggests, a counting process simply counts the number of times an event occurs prior to a given time.  For example, we could define $N_t := \{ \text{number of buy orders prior to time $t$} \}$.  Of course, there are many ways to model a counting process.  And, depending on how we model $N$, the counting process may or may not be Markov.  However, because Markov processes allow for many analytically tractable computations, we will focus on Markov counting processes; the most basic of these is the Poisson process.

\begin{definition}
\label{def:poisson.process}
A \emph{Poisson process with intensity $\lam$} is a stochastic process $N = (N_t)_{t \geq 0}$ taking values in $S = \{ 0, 1, 2, \ldots \}$ such that the following hold.
\begin{enumerate*}
\item $N_0 = 0$.
\item If $s<t$ then $N_s \leq N_t$.
\item If $s<t$ then $(N_t - N_s) \ind \Fc_s$, where $\Fc_s=\sigma(N_r, 0\leq r\leq s)$. \label{itm:ind.inc}
\item Lastly, 
\begin{align}
\Pb( N_{t+s} = n+m | N_t = n) 
	&=	\left\{ \begin{array}{ll}
			\lam s + \Oc(s^2) & m=1 ,\\
			\Oc(s^2)					& m\geq 2 ,\\
			1-\lam s + \Oc(s^2) & m = 0 ,
			\end{array} \right.
\end{align}
as $s \to 0^+$.
\end{enumerate*}
\end{definition}

\noindent
Note that Item \ref{itm:ind.inc} implies that $N$ is Markov.  To see this, define $\Fc_t := \sig(N_s, 0 \leq s \leq t)$ and observe that
\begin{align}
\Eb[ g(N_t) | \Fc_s ]
	&=	\Eb[ g(N_t - N_s + N_s) | N_s ]
	=		\sum_{n=0}^\infty g( n + N_s) f_{N_t-N_s}(n) .
\end{align}
Let us find the probability mass function of $N_t$.

\begin{theorem}
Let $N = (N_t)_{t \geq 0}$ be a Poisson process with parameter $\lam$.
Then, for all $t \geq 0$ we have $N_t \sim \text{\emph{Poi}}(\lam t)$.  That is
\begin{align}
p_t(j) := \Pb( N_t = j )
	&=	\frac{(\lam t)^j}{j!} \ee^{- \lam t} . \label{eq:poisson.pmf}
\end{align}
\end{theorem}

\begin{proof}
We have
\begin{align}
\Pb( N_{t+s} = j )
	&=	\sum_{i=0}^j \Pb(N_{t+s} = j | N_t = i ) \Pb( N_t = i ) \\
	&=	\sum_{i=0}^j \Pb(N_{t+s} - N_t = j - i ) \Pb( N_t = i ) \\
	&=	\Pb(N_{t+s} - N_t = 0  ) \Pb( N_t = j ) + \Pb(N_{t+s} - N_t = 1  ) \Pb( N_t = j-1 ) \\ &\quad 
			+ \sum_{i=0}^{j-2} \Pb(N_{t+s} - N_t = j - i ) \Pb( N_t = i ) \\
	&=	( 1 - \lam s ) \Pb( N_t = j ) + \lam s \Pb( N_t = j-1 ) + \Oc(s^2) .
\end{align}
Subtracting $\Pb( N_t = j )=:p_t(j)$ from both sides, dividing by $s$ and taking a limit as $s \to 0$ we obtain
\begin{align}
\d_t p_t(0)
	&=	- \lam p_t(0) , &
p_0(0)
	&=	1 , \label{eq:poisson.pde.0} \\
\d_t p_t(j)
	&=	- \lam p_t(j)  + \lam p_t(j-1) , &
p_0(j)
	&=	0 , &
j
	&\geq 1 , \label{eq:poisson.pde.j}
\end{align}
where we have included the obvious boundary conditions.
Thus, we have obtained a sequence of nested ODEs for $(p_t(j))_{j \geq 0}$.  
One can easily verify that \eqref{eq:poisson.pmf} satisfies \eqref{eq:poisson.pde.0}-\eqref{eq:poisson.pde.j}.
\end{proof}

For any counting process (whether or not it is Poisson), we may be interested to know the arrival time of the $n$th event and the time between the $n$th and $(n+1)$th events.

\begin{definition}
\label{def:inter-arrival}
Let $N = (N_t)_{t \geq 0}$ be a counting process.  We define $S_n$, the \emph{the $n$th arrival time}, by
\begin{align}
S_0
	&=	0 , &
S_n
	&=	\inf\{ t \geq 0 : N_t = n \} , &
n
	&\geq 1 , \label{eq:S.def}
\end{align}
with $\inf\{ \emptyset \} = \infty$.
We define $\tau_n$, the \emph{$n$th inter-arrival time}, by
\begin{align}
\tau_n
	&:=	S_n - S_{n-1} , &
n
	&\geq 1 . \label{eq:tau.def}
\end{align}
\end{definition}

\noindent
Given the full history of $N$, we can construct the sequence $(\tau_n)_{n \geq 1}$ using \eqref{eq:S.def}--\eqref{eq:tau.def}.  Alternatively, given the sequence $(\tau_n)_{n \geq 1}$, we can construct the path of $N$ using
\begin{align}
S_n
	&=	\sum_{i=1}^n \tau_i , &
n
	&\geq 1 , &
N_t
	&=	\sum_{n=1}^\infty n \Ib_{\{S_n \leq t < S_{n+1}\}} \label{eq:S=sum.tau}
\end{align}
We now undertake this alternative construction.

\begin{theorem}
\label{thm:construction.poisson}
Suppose the inter-arrival times $\tau_i$ (see Definition \ref{def:inter-arrival}) of a counting process $N = (N_t)_{t \geq 0}$ are iid and exponentially distributed with parameter $\lam$.  Then $N$ is a Poisson process with parameter $\lam$.
\end{theorem}

\begin{proof}
We will show that if the inter-arrival times $(\tau_i)_{i \geq 1}$ are iid with $\tau_i \sim \Ec(\lam)$ then $N_t$, given by \eqref{eq:S=sum.tau} has a probability mass function given by \eqref{eq:poisson.pmf}.
First, we show that $S_n$, as defined in \eqref{eq:S=sum.tau}, has a \emph{Gamma} density
\begin{align}
\Pb(S_n \in \dd x)
	&=	f_{S_n}(x) \dd x , &
f_{S_n}(x)
	&=	\Ib_{\{x \geq 0\}} \frac{(\lam x)^{n-1}}{(n-1)!}\lam \ee^{- \lam x} . \label{eq:gamma.density}
\end{align}
Note that $S_1 = \tau_1$ and $f_\tau(x) = f_{S_1}(x) = \Ib_{\{x \geq 0\}} \lam \ee^{- \lam x}$.  Thus, \eqref{eq:gamma.density} is correct for $n=1$.  We now assume \eqref{eq:gamma.density} holds for $n$ and show that it holds for $n+1$.  As $S_{n+1} = S_n + \tau_{n+1}$ we have
\begin{align}
f_{S_{n+1}}(x)
	&=	\int_{\Rb^2} \dd y \dd z \, f_{S_n}(y) \cdot f_{\tau}(z) \cdot \del(y+z-x) \\
	&=	\int_\Rb \dd y \, f_{S_n}(y) \cdot f_{\tau}(x-y) \\
	&=	\int_\Rb \dd y \, \Ib_{\{y \geq 0\}} \frac{(\lam y)^{n-1}}{(n-1)!}\lam \ee^{- \lam y} \cdot \Ib_{\{x-y \geq 0\}} \lam \ee^{- \lam (x-y)} \\
	&=	\Ib_{\{x \geq 0\}} \frac{\lam^{n+1}\ee^{- \lam x}}{(n-1)!} \int_0^x \dd y \, y^{n-1}
	=		\Ib_{\{x \geq 0\}}\frac{(\lam x)^{n}}{(n)!}\lam \ee^{- \lam x} ,
\end{align}
which agrees with \eqref{eq:gamma.density}.  Thus, by induction, \eqref{eq:gamma.density} holds for every $n$.
Using the above result, we now show that $N_t$, as defined in \eqref{eq:S=sum.tau}, has probability mass function \eqref{eq:poisson.pmf}.  For $n \geq 1$ we clearly have
\begin{align}
\Pb( N_t \geq  n )
	&=	\Pb( S_n \leq t )
	=		\int_0^t \dd x \, f_{S_n}(x) 
	= 	\int_0^t \dd x \, \frac{(\lam x)^{n-1}}{(n-1)!}\lam \ee^{- \lam x} .
\end{align}
Using this result, we compute
\begin{align}
\Pb( N_t \geq  n + 1 )
	&=	\int_0^t \dd x \, \frac{(\lam x)^{n}}{(n)!}\lam \ee^{- \lam x}
	=		- \frac{(\lam t)^n}{n!}\ee^{- \lam t} + \Pb( N_t \geq n ) ,
\end{align}
where we have used integration by parts to obtain the second equality.  Thus, we have
\begin{align}
\Pb( N_t = n )
	&=	\Pb( N_t \geq n) - \Pb(N_t \geq n+1) 
	=		\frac{(\lam t)^n}{n!}\ee^{- \lam t} ,
\end{align}
which agrees with \eqref{eq:poisson.pmf}.  To complete the proof, we observe that
\begin{align}
\Pb(N_t = 0)
	&=	\Pb( S_1 > t )
	=		\Pb( \tau_1 > t )
	=		\ee^{- \lam t} ,
\end{align}
which shows that \eqref{eq:poisson.pmf} holds for $n=0$.
\end{proof}

\noindent
The inter-arrival times $(\tau_i)_{i \geq 1}$ of a Poisson process $N$ are \emph{memoryless} in the sense that
\begin{align}
\Pb(\tau_i > s+t | \tau_i > s)
	&=	\frac{\Pb(\tau_i > s+t , \tau_i > s)}{\Pb(\tau_i > s)}
	=		\frac{\Pb(\tau_i > s+t)}{\Pb(\tau_i > s)}
	=		\frac{\ee^{-\lam(t+s)}}{\ee^{-\lam s}}
	=		\ee^{-\lam t}
	=		\Pb( \tau_i > t ) . \label{eq:memoryless}
\end{align}
The memoryless property of the inter-arrival times is what gives the Poisson process its Markov property.  If we were to construct a counting process whose inter-arrival times were not exponentially distributed, we would find that this process is not Markov.  The reason is that if inter-arrival times are not exponentially distributed, then the amount of time we have to wait for the next event arrival will depend on how much time has passed since the last event arrival.
\\[0.5em]
It should be clear from the construction of $N$ that the Poisson process has stationary increments
\begin{align}
N_t - N_s
	&\sim N_{t-s} ,
\end{align}
since $N_t - N_s \ind N_s$.


\section{Overview of continuous time Markov chains}
\label{sec:ctmc}
We have already defined what we mean by Markov process in Definition \ref{def:markov}.  In this section, we will derive some results that are specific for Markov continuous time processes that have a countable state space.

\begin{definition}
A \emph{continuous time Markov chain} $X = (X_t)_{t  \geq 0}$ is a Markov process with a countable state space $S$.
\end{definition}

\noindent
Recall from \eqref{eq:markov.prop} that a Markov process satisfies $\Pb( X_{t+s} \in A | \Fc_s) = \Pb( X_{t+s} \in A | X_s)$.  Here, we take the filtration $\Fc_s$ to be the \emph{natural filtration} for $X$.  That is, $\Fc_s = \sig(X_u,0\leq u \leq s)$.   What the Markov property means in practice is that the evolution of a Markov chain is described by its one-step transition probabilities $\Pb( X_{t+s} = j | X_s = i)$.  In general, such probabilities may depend on the times $t,s$.  Throughout this Chapter, we will restrict our attention to Markov chains whose transition probabilities do not depend on $s$.

\begin{definition}
A continuous time Markov chain $X$ is called \emph{homogeneous} if 
\begin{align}
\Pb( X_{t+s} = j | X_s = i)
	&=	\Pb( X_{t} = j | X_0 = i) =: p_t(i,j) &
	&		\forall \, t,s \in \Rb_+ , \forall \, i,j \in S .
\end{align}
Let $\Pv_t$ be the $|S| \times |S|$ matrix $\Pv_t = ( p_t(i,j) )$.  We call the collection of matrices $(\Pv_t)_{t \geq 0}$ the \emph{transition semigroup}.
\end{definition}

\begin{theorem}
The transition semigroup $(\Pv_t)_{t \geq 0}$ satisfies the following properties.
\begin{enumerate*}
\item $\Pv_0 = \Iv$.\label{itm:I}
\item $\sum_j p_t(i,j) = 1$.\label{itm:sum}
\item $\Pv_t \Pv_s = \Pv_{t+s}$.\label{itm:semigroup}
\end{enumerate*}
\end{theorem}

\begin{proof}
Item \ref{itm:I} follows from $\Pb( X_0 = j | X_0 = i) = \del_{i,j}$.
Item \ref{itm:sum} follows from $\sum_j \Pb( X_t = j | X_0 = i) = 1$.
Lastly, to show Item \ref{itm:semigroup} note that
\begin{align}
p_{t+s}(i,j)
	&=	\Pb( X_{t+s} = j | X_0 = i) 
	=	\sum_k \Pb( X_{t+s} = j | X_t = k) \cdot \Pb( X_t = k | X_0 = i) \\
	&=	\sum_k \Pb( X_{s} = j | X_0 = k) \cdot \Pb( X_t = k | X_0 = i) 
	=	\sum_k p_{t}(i,k) p_{s}(k,j) . \label{eq:chapman-kolmogorov}
\end{align}
\end{proof}

\noindent
We now construct a continuous-time Markov chain in very much the same manner we constructed the Poisson process in Definition \ref{def:poisson.process}.
The following derivation is \emph{purely formal} (i.e., \emph{not} rigorous).
Some of what we write is not true in general, but holds in many practical applications.

\begin{definition}
\label{def:ctmc}
A continuous time Markov chain with \emph{generator} $\Gv = (g(i,j))_{i,j}$ is a stochastic process satisfying
\begin{enumerate*}
\item Markov property: if $s<t$ then $\Pb( X_t = j | \Fc_s^X ) = \Pb( X_t = j | X_s )$ for all $j \in S$, where $\Fc_s^X := \sig(X_u, 0 \leq u \leq s)$.
\item As $s \to 0$ we have $\Pb( X_{t+s} = j | X_t = i ) =	\del_{ij} + g(i,j) s + \Oc(s^2)$. \label{itm:del+g}
\end{enumerate*}
\end{definition}

\noindent
The essence of Definition \ref{def:ctmc} is as follows.  In a small time interval $[t,t+s)$ the process $X$ either remains at its current state $X_t=i$ or it jumps to a new state $j \neq i$.  The probability of remaining at $i$ is $1+g(i,i)s + \Oc(s^2)$.  The probability of jumping to state $j$ is $g(i,j)s + \Oc(s^2)$.  Since probabilities should always fall in $[0,1]$ we must have
\begin{align}
g(i,i)
	&\leq 0 , &
g(i,j)
	&\geq 0 , &
i
	&\neq j .
\end{align}
Moreover, since $X$ should be found \emph{somewhere} with probability one we must also have
\begin{align}
1
	&=	\sum_j \Pb( X_{t+s} = j | X_t = i ) \\
	&=	\sum_j ( \del_{ij} + g(i,j) s ) + \Oc(s^2)
	=		1 + s \sum_j g(i,j) + \Oc(s^2) &
	&		\Rightarrow &
0
	&=	\sum_j g(i,j) .
\end{align}
Note that item \ref{itm:del+g} in Definition \ref{def:ctmc} can be written more compactly as
\begin{align}
\Gv
	&= \lim_{s \to 0^+} \frac{1}{s} \Big( \Pv_s - \Iv \Big) .
\end{align}
Thus, it is clear that we can obtain $\Gv$ from knowledge of $\Pv_t$.  We can also obtain $\Pv_t$ from $\Gv$.

\begin{theorem}
\label{thm:kfe.kbe}
Let $X = (X_t)_{t \geq 0}$ be a continuous time Markov chain with generator $\Gv$ and a finite state space $S$.  The transition semigroup of $X$ satisfies the following ODEs
\begin{align}
\text{Kolmogorov Forward equation}:&&
\frac{\dd}{\dd t} \Pv_t
	&=	\Pv_t \Gv , \label{eq:kfe} \\
\text{Kolmogorov Backward equation}:&&
\frac{\dd}{\dd t} \Pv_t
	&=	\Gv \Pv_t . \label{eq:kbe}
\end{align}
\end{theorem}

\begin{proof}
To derive the forward equation \eqref{eq:kfe} we computed $\Pb(X_{t+s} = j | X_0 = i )$ by conditioning on the value of $X_t$.  We have 
\begin{align}
p_{t+s}(i,j)
	&=	\sum_k \Pb(X_{t+s}=j | X_t = k) \Pb(X_t = k | X_0 = i)
	=	\sum_k p_{t}(i,k) p_{s}(k,j) \\
	&=	p_t(i,j) (1 + g(j,j)s) + \sum_{k: k\neq j} p_{t}(i,k)  g(k,j)s + \Oc(s^2) .
\end{align}
Subtracting $p_t(i,j)$ from both sides, dividing by $s$ and taking a limit as $s$ goes to zero, we obtain
\begin{align}
\lim_{s \searrow 0 } \frac{1}{s} \Big( p_{t+s}(i,j) - p_t(i,j) \Big) = \frac{\dd}{\dd t} p_t(i,j)
	&=	\sum_{k} p_{t}(i,k)  g(k,j) ,
\end{align}
which is the component-wise representation of \eqref{eq:kfe}.
To derive the backward equation \eqref{eq:kbe} we computed $\Pb(X_{t+s} = j | X_0 = i )$ by conditioning on the value of $X_s$.  We have 
\begin{align}
p_{t+s}(i,j)
	&=	\sum_k \Pb(X_{t+s}=j | X_s = k) \Pb(X_s = k | X_0 = i)
	=	\sum_k p_{s}(i,k) p_{t}(k,j) \\
	&=	(1 + g(i,i)s)p_t(i,j)  + \sum_{k: k\neq i} g(i,k)s \, p_{t}(k,j)   + \Oc(s^2) .
\end{align}
Subtracting $p_t(i,j)$ from both sides, dividing by $s$ and taking a limit as $s$ goes to zero, we obtain
\begin{align}
\lim_{s \searrow 0 } \frac{1}{s} \Big( p_{t+s}(i,j) - p_t(i,j) \Big) = \frac{\dd}{\dd t} p_t(i,j)
	&=	\sum_{k} g(i,k) p_{t}(k,j) ,
\end{align}
which is the component-wise representation of \eqref{eq:kbe}.
\end{proof}

\noindent
If $\Pv_t$ solves the forward equation \eqref{eq:kfe} then it also solves the backward equation \eqref{eq:kbe} and vice versa.  The solution to \eqref{eq:kfe} and \eqref{eq:kbe}, subject to the initial condition $\Pv_0 = \Iv$ is
\begin{align}
\Pv_t
	&=	\ee^{t \Gv} :=	\sum_{n=0}^\infty \frac{1}{n!} t^n \Gv^n .
\end{align}
If $\Gv$ has an eigenvector decomposition $\Gv = \Uv^{-1} \Lamv \Uv$ as in Section \ref{sec:finite.chains} then
\begin{align}
\Pv_t
	&=	\Uv^{-1} \Big( \sum_{n=0}^\infty \frac{1}{n!} (t \Lamv)^n \Big) \Uv 
	=		\Uv^{-1} \ee^{t \Lamv} \Uv . \label{eq:Pt.spectral}
\end{align}

\begin{example}
\label{ex:ctmc}
Let the generator $\Gv$ of a continuous time Markov chain $X$ be given by
\begin{align}
\Gv
	&=	\( \begin{array}{cc}
			-\alpha & \alpha \\
			\beta & -\beta
			\end{array} \) , &
\alpha, \beta > 0 . \label{eq:G.2state}
\end{align}
We wish to find the semigroup $\Pv_t$ generated by $\Gv$.  We note that $\Gv$ has left eigenvectors and eigenvalues given by
\begin{align}
\lam_0
	&=	0 , &
\ev_0
	&= \frac{1}{\sqrt{\alpha^2 + \beta^2}} (\beta , \alpha), &
\lam_1
	&= - \alpha - \beta	, &
\ev_1
	&=	\frac{1}{\sqrt{2}} ( 1 ,- 1 ) .
\end{align} 
Thus, from \eqref{eq:Pt.spectral} we have
\begin{align}
\Pv_t
	&=		\Uv^{-1} \ee^{t \Lamv} \Uv 
	=			\( \begin{array}{c}
				\ev_0 \\ \ev_1
				\end{array} \)^{-1}
				\( \begin{array}{cc}
				\ee^{t \lam_0} & 0 \\
				0 & \ee^{t \lam_1}
				\end{array} \)
				\( \begin{array}{c}
				\ev_0 \\ \ev_1
				\end{array} \) \\
	&=		\frac{1}{\alpha + \beta} 
				\( \begin{array}{cc}
				\beta + \alpha \ee^{-(\alpha+\beta)t} & \alpha - \alpha \ee^{-(\alpha+\beta)t} \\
				\beta - \beta \ee^{-(\alpha+\beta)t} & \alpha + \beta \ee^{-(\alpha+\beta)t}
				\end{array} \) . \label{eq:Pt.2state}
\end{align}
\end{example}

\noindent
For a Poisson process with intensity $\lam$, we found that we could construct $N$ from a sequence of iid inter-arrival times $(\tau_i)_{i \geq 1}$ that were exponentially distributed: $\tau_i \sim \Ec(\lam)$.  A continuous time Markov chain can be constructed in a similar manner.  
\\[0.5em]
Let $\Gv = (g(i,j))_{i,j \in S}$ be the generator of continuous time Markov chain.
Consider a discrete time Markov chain $Y = (Y_n)_{n \in \Nb_0}$ with one-step transition matrix $\Pv = (p(i,j))_{i,j \in S}$.  Suppose the entries of $\Pv$ are given by
\begin{align}
p(i,j)
	&=	- g(i,j)/g(i,i) , &
i
	&\neq j , &
p(i,i)
	&=	0 .
\end{align}
Observe that $p(i,j) \geq 0$ $\forall i,j \in S$ and $\sum_j p(i,j) = 1$ $\forall i$ as required.
Now, let $(\tau_i)_{i \in \Nb}$ be a sequence of independent random variables satisfying
\begin{align}
\tau_n | Y_n
	&\sim \Ec( -g(Y_n,Y_n) ) .
\end{align}
Let $(S_n)_{n \in \Nb_0}$ be given by
\begin{align}
S_n
	&:=	\sum_{i=1}^n \tau_i , &
S_0
	&:=	0 .
\end{align}
We can construct a continuous time process $X = (X_t)_{t \geq 0}$ as follows
\begin{align}
X_t
	&=	\sum_{n=0}^\infty Y_n \Ib_{\{S_n \leq t < S_{n+1}\}} .
\end{align}
We claim that the process $X$ is a continuous time Markov chain with generator $\Gv$.  We call $(S_n)_{n \in \Nb_0}$ the \emph{jump times} and $(\tau_n)_{n \in \Nb}$ the \emph{inter-jump times}.  Compare this construction of a Markov chain to the second construction we gave for a Poisson process (see Theorem \ref{thm:construction.poisson} and the text preceding it).
\\[0.5em]
The above construction makes clear what the dynamics of a CTMC $X$ with generator $\Gv$ are.  Specifically, once $X$ jumps into state $i$, the amount of time it remains in that state is exponentially distributed with parameter $-g(i,i)$.  When a jump occurs, the probability that $X$ jumps from $i$ to $j$ is $-g(i,j)/g(i,i)$.
\\[0.5em]  
To see how this construction is equivalent to Definition \ref{def:ctmc}, note that, for small $s$, the probability of two jumps in the time interval $[0,s)$ is $\Oc(s^2)$ because
\begin{align}
\Pb( S_2 \leq s )
	&\leq \Pb(\tau_1 \leq s) \Pb( \tau_2 \leq s )
	\leq (1 - \ee^{\underline{g} s}) (1 - \ee^{\underline{g} s} )
	=	(1 - (1 + \underline{g} s) + \Oc(s^2) )^2
	=	\Oc(s^2) ,
\end{align} 
where we have defined $\underline{g} := \inf_{ i \in S} \{ -g(i,i) \}$.
As the probability of two jumps is $\Oc(s^2)$, for small $s$ we have
\begin{align}
p_s(i,i)
	&=	\Pb( X_{s} = i |  X_0 = i)
	=		\Pb( \tau_1 > s | X_0 = i ) + \Oc(s^2) \\
	&=	\exp( g(i,i) s) + \Oc(s^2)
	=		1 + g(i,i) s + \Oc(s^2) , \\
p_s(i,j)
	&=	\Pb( X_{s} = j |  X_0 = i)
	=		\Pb( \tau_1 \leq s | X_0 = i ) \cdot \Pb( X_{\tau_1} = j |  X_0 = i)  + \Oc(s^2) \\
	&=	( 1 - \exp( g(i,i) s) ) \cdot \frac{g(i,j)}{-g(i,i)} + \Oc(s^2) \\
	&=		- g(i,i) s \cdot \frac{g(i,j)}{-g(i,i)} + \Oc(s^2)
	=		g(i,j) s + \Oc(s^2) ,
\end{align}
which agrees with Item \ref{itm:del+g} in Definition \ref{def:ctmc}.  The independent increments property $X_t - X_s \ind X_s$ follows from the memoryless property of the exponentially distributed waiting times.
\\[0.5em]
For discrete time Markov chains, we found that $\muv_{n+m} = \muv_n \Pv_m$.  We defined an invariant distribution $\piv$ as any probability mass function that satisfied $\piv = \piv \Pv$, which implied that $\piv = \piv \Pv_n$ for all $n$.  In continuous time, we have the following analogs
\begin{align}
\muv_{t+s}
	&=	\muv_t \Pv_s , &
	& 	\forall \, t,s \geq 0 , \\
\piv
	&=	\piv \Pv_t , & 
	&		\forall \, t \geq 0 .
\end{align}
It is not always easy to find an invariant distribution $\piv$ from $\Pv_t$.  However, $\piv$ can also be obtained from the generator $\Gv$.

\begin{theorem}
Let $X$ be a continuous time Markov chain with generator $\Gv$ and a finite state space $S$, and let $\Pv_t$ be the semigroup generated by $\Gv$.  Then
\begin{align}
\piv
	&=	\piv \Pv_t \quad \forall \, t \geq 0 , &
	&\Leftrightarrow &
\piv \Gv
	&=	0 .
\end{align}
\label{ctmc-stat}
\end{theorem}

\begin{proof}
Recalling that $\Gv^0 \equiv \Iv$ we have
\begin{align}
\piv \Gv
	&=	0 &
	&\Leftrightarrow&
\piv \Gv^n
	&=	0 , &
	&		\forall \, n \geq 1 , \\ & & &\Leftrightarrow&
\sum_{n=1}^\infty \frac{1}{n!} t^n \piv \Gv^n
	&=	0 \\ & & &\Leftrightarrow&
\piv  \sum_{n=0}^\infty \frac{1}{n!} t^n \Gv^n
	&=	\piv \\ & & &\Leftrightarrow&
\piv \Pv_t
	&=	\piv .
\end{align}
\end{proof}

\begin{example}
Let us return to Example \ref{ex:ctmc}.  Let us find the stationary distribution for the Markov process $X$ with generator $\Gv$ given by \eqref{eq:G.2state}.  We have
\begin{align}
0
	&=	\( \pi(0) , \pi(1) \) \( \begin{array}{cc}
			-\alpha & \alpha \\
			\beta & -\beta
			\end{array} \) , &
1
	&= \pi(0) + \pi(1) , &
\pi(i)
	&\geq 0 , &
	&\Rightarrow&
\piv
	&=  \(  \frac{\beta}{\alpha + \beta} ,  \frac{\alpha}{\alpha + \beta} \) .
\end{align}
Note that from \eqref{eq:Pt.2state} we have for any $\muv_0$ that
\begin{align}
\lim_{t \to \infty} \muv_t
	&=		\muv_0 \lim_{t \to \infty} \Pv_t
	=		\( \mu_0(0) , \mu_0(1) \) 
				\frac{1}{\alpha + \beta} 
				\( \begin{array}{cc}
				\beta  & \alpha  \\
				\beta  & \alpha
				\end{array} \)
	=		\piv ,
\end{align}
where we have used $\mu_0(0) + \mu_1(1) = 1 $.
\end{example}

\noindent
As in the discrete time case, for a given Markov chain $X$ a stationary distribution $\piv$ may not exist.  If a stationary distribution does exist, then it may not be unique.  Given the existence of a stationary distribution, the condition for uniqueness is that
\begin{align}
\text{for all $i,j \in S$, there exists $t>0$ such that  $p_t(i,j)>0$.} \label{irreducibility}
\end{align} 
If \eqref{irreducibility} is satisfied, then we say the chain $X$ is \emph{irreducible}.

\begin{example}
\label{ex:birth-death}
Consider a continuous time Markov chain $X$ with generator
\begin{align}
\Gv
	&=	\( \begin{array}{cccccc}
			-\lam_0 & \lam_0 & 0 & 0 & 0 & \ldots \\
			\mu_1 & - (\mu_1 + \lam_1) & \lam_1 & 0 & 0 & \ldots \\
			0 & \mu_2 & - (\mu_2 + \lam_2) & \lam_2 & 0 & \ldots \\
			0 & 0 & \mu_3 & - (\mu_3 + \lam_3) & \lam_3 & \ldots \\
			\vdots & \vdots & \vdots & \vdots & \vdots & \ddots 
			\end{array} \) .
\end{align}
The chains $X$ is called a a \emph{birth-death} process since, when $X_t = i$ it either jumps to $i+1$ (i.e., a birth) or jumps down to $i-1$ (i.e., a death).  The expected waiting time $\tau$ in state $i$ is exponentially distributed $\tau \sim \Ec(\lam_i + \mu_i)$.  When the process does jump out of state $i$, the probability of an up-jump is $\lam_i / (\lam_i+\mu_i)$ and the probability of a down jump is $\mu_i / (\lam_i+\mu_i)$.  If $\lam_0 >0$ then the chain is irreducible.  Let us see if we can find the stationary distribution $\piv$.  We see $\piv$ such that $\piv \Gv = 0$, which leads to
\begin{align}
0
	&=	-\lam_0 \pi(0) + \mu_1 \pi(1) , \\
0
	&=	\lam_{n-1} \pi(n-1) - (\mu_n + \lam_n) \pi(n) + \mu_{n+1} \pi(n+1) , &
n
	&\geq 1 .
\end{align}
One can check by direct substitution that
\begin{align}
\pi(n)
	&=	\frac{\lam_0 \lam_1 \ldots \lam_{n-1}}{\mu_1 \mu_2 \ldots \mu_n}\pi(0) .
\end{align}
In order for $\piv$ to be a distribution, the following normalization condition must be satisfied
\begin{align}
1
	&=	\sum_{n=0}^\infty \pi(n) 
	=		\pi(0) + \pi(0) \sum_{n=1}^\infty \frac{\lam_0 \lam_1 \ldots \lam_{n-1}}{\mu_1 \mu_2 \ldots \mu_n} , &
	&		\Rightarrow &
\pi(0)
	&=	\( 1 + \sum_{n=1}^\infty \frac{\lam_0 \lam_1 \ldots \lam_{n-1}}{\mu_1 \mu_2 \ldots \mu_n}  \)^{-1} .
\end{align}
Clearly, we can only normalize if
\begin{align}
\sum_{n=1}^\infty \frac{\lam_0 \lam_1 \ldots \lam_{n-1}}{\mu_1 \mu_2 \ldots \mu_n} < \infty .
\end{align}
\end{example}

\noindent
While a stationary distribution $\piv$ can often be computed with relative ease, it is almost always difficult (or impossible) to compute $\muv_t$ explicitly.  Nevertheless, we can still obtain moments $\Eb X_t^n$ if we can compute the generating function $G_{X_t}(s):=\Eb s^{X_t}$.

\begin{example}
Consider the birth-death process given in Example \ref{ex:birth-death} and assume 
\begin{align}
\lam_n
	&=	n \lam , &
\mu_n
	&=	n \mu , &
\lam , \mu
	&> 0 .
\end{align}
Note that the state $X_t = 0$ is an absorbing state since, $\lam_0 = 0$.  Assume $X_0 = i \geq 1$.  The Kolmogorov forward equations \eqref{eq:kfe} become
\begin{align}
\frac{\dd}{\dd t}p_t(i,0)
	&=	\mu p_t(i,1) , \\
\frac{\dd}{\dd t}p_t(i,j)
	&=	\lam (j-1) p_t(i,j-1) - (\lam + \mu) j p_t(i,j) + \mu ( j + 1) p_t(i,j+1) , &
j
	& \geq 1 .
\end{align}
Multiplying the $j$th equation by $s^j$ and summing all of the equations, we obtain
\begin{align}
\frac{\d }{\d t}G_{X_t}
	&=	\lam s^2 \frac{\d G_{X_t}}{\d s} - (\lam + \mu) s \frac{\d G_{X_t}}{\d s} + \mu \frac{\d G_{X_t}}{\d s} , \label{eq:GX.pde}
\end{align}
where we have used
\begin{align}
G_{X_t}(s)
	&:=	\Eb[ s^{X_t} | X_0 = i ] 
	=		\sum_{j=0}^\infty s^j p_t(i,j) , 
\end{align}
from which it follows that
\begin{align}
\frac{\d G_{X_t}(s)}{\d t}
	&=	\sum_{j=0}^\infty s^j \frac{\d p_t(i,j) }{\d t} , &
\frac{\d G_{X_t}(s)}{\d s}
	&=	\sum_{j=1}^\infty j s^{j-1} p_t(i,j) .
\end{align}
One can check by direct substitution that the solution to PDE \eqref{eq:GX.pde} with initial condition $G_{X_0}(s) = s^i$ is given by
\begin{align}
G_{X_t}(s)
	&=	\( \frac{\mu(1-s) - (\mu - \lam s)\ee^{-(\lam-\mu)t}}{\lam(1-s) - (\mu - \lam s)\ee^{-(\lam-\mu)t}} \)^i , &
\mu 
	&\neq \lam .
\end{align}
From the generating function we can compute the first to moments of $X_t$ easily.  We have
\begin{align}
\Eb X_t
	&=	i \ee^{(\lam-\mu)t} , &
\Vb X_t
	&=	i \frac{\lam + \mu}{\lam-\mu}\ee^{(\lam-\mu)t}\(\ee^{(\lam-\mu)t} - 1 \) , &
\mu
	&\neq \lam .
\end{align}
\end{example}

\section{Reversibility}
\label{sec:ctmc_RV}

In Section \ref{sec:reversibility}, we introduced a notion of \textit{reversibility} for discrete time Markov chains.  In this section, we will discuss reversibility of continuous time Markov chains.  As we shall see, most of the definitions and theorems discussed below are direct analogs of the definitions and theorems given in Section \ref{sec:reversibility}.

\begin{definition}
\label{def:reversal.cont}
Let $X = (X_t)_{0 \leq t \leq T}$ be a continuous time Markov chain.  The \textit{time reversal of $X$} is the process $Y = (Y_t)_{0 \leq t \leq T}$ defined as follows
\begin{align}
Y_t
	&:=	\lim_{\del \searrow 0} X_{T-t+\del} = X_{(T-t)^+}.
\end{align}
\end{definition}

\noindent
Compare Definition \ref{def:reversal.cont} for the time reversal of a continuous time Markov chain with Definition \ref{def:Y.markov} for the time reversal of a discrete time Markov Chain.  Notice the similarities.

\begin{theorem}
\label{thm:generator.reversal}
Let $X = (X_t)_{0 \leq t \leq T}$ be an irreducible continuous time Markov chain with 
generator $\Gv$, transition semigroup $\Pv_t = \ee^{t \Gv}$ and invariant distribution $\piv$ and suppose $X_0 \sim \piv$.  
Let $Y = (Y_t)_{0 \leq t \leq T}$ be the time reversal of $X$.  
Then $Y$ is a continuous time Markov chain with generator $\Hv$ and transition semigroup $\Qv_t = \ee^{t \Hv}$, which satisfy
\begin{align}
h(i,j)
	&=	\frac{\pi(j)}{\pi(i)}g(j,i) ,  &
q_t(i,j)
	&=	\frac{\pi(j)}{\pi(i)}p_t(j,i) . \label{eq:h.g}
\end{align}
\end{theorem}

\noindent
Rather than prove Theorem \ref{thm:generator.reversal}, we will explain intuitively why it is true.  Fix a large $N \in \Nb$, define $\del := T/N \ll 1$ and consider the discrete-time processes $\Xh = (\Xh_n)_{0 \leq n \leq N}$ where $\Xh_n := X_{n \del}$.  The process $\Xh$ is a discrete time Markov chain with transition probabilities  $\ph(i,j) := \Pb( \Xh_{n+1} = j | \Xh_n = i)$ given by
\begin{align}
\ph(i,j)
	&=	\del_{i,j} + g(i,j) \del + \Oc(\del^2 ) . \label{eq:p-hat-1}
\end{align}
Now, let $\Yh = (\Yh_n)_{0 \leq n \leq N}$ be the time-reversal of $\Xh$.  That is $\Yh_n := \Xh_{N-n} = X_{T - n \del} = Y_{n \del}$.  We know from Theorem \ref{thm:Y.markov} that $\Yh$ is a Markov process with transition probabilities $\qh(i,j) := \Pb( \Yh_{n+1} = j | \Yh_n = i)$ given by
\begin{align}
\qh(i,j)
	&=	\frac{\pi(j)}{\pi(i)} \ph(j,i) . \label{eq:q-hat-2}
\end{align}
If we accept that $Y$, like its discrete-time counterpart $\Yh$, is a Markov process, then we must have
\begin{align}
\qh(i,j)
	&=	\del_{i,j} + h(i,j) \del  + \Oc(\del^2 ) . \label{eq:q-hat-1}
\end{align}
for some generating matrix $H = (h(i,j))_{i,j \in S}$.  Inserting \eqref{eq:p-hat-1} and \eqref{eq:q-hat-1} in \eqref{eq:q-hat-2} we obtain
\begin{align}
\del_{i,j} + h(i,j) \del 
	&=	\frac{\pi(j)}{\pi(i)} (\del_{j,i} + g(j,i) \del ) + \Oc(\del^2 ) 
\end{align}
The terms of order $\Oc( 1 )$ cancel (you can check they are equal for both for $i = j$ and for $i \neq j$).  In order for the $\Oc(\del)$ terms to match, we must have
\begin{align}
h(i,j)
	&=	\frac{\pi(j)}{\pi(i)} g(j,i) , 
\end{align}
which is the left-hand side of equation \eqref{eq:h.g}.  Now, note that
\begin{align}
\pi(i)h(i,j) 
	&= \pi(j)g(j,i) &
	&\Rightarrow &
\pi(i) (\Hv^n)_{ij} 
	&= \pi(j) (\Gv^n)_{ji} , &
	& \forall n \in \Nb_0 ,
\end{align}
which implies
\begin{align}
\pi(i) (\ee^{t \Hv}  )_{ij} 
	&= \pi(j) (\ee^{t \Gv}  )_{ji},  &
	&\forall \, t \in [0 ,T] .
\end{align}
The right-hand side of \eqref{eq:h.g} follows from $(\ee^{t \Gv}  )_{ij} = ( \Pv_t)_{i,j} = p_t(i,j)$ and $(\ee^{t \Hv}  )_{ij} = ( \Qv_t)_{ij} = q_t(i,j)$.

\begin{definition}
\label{def:reversible.2}
Let $X = (X_t)_{0 \leq t \leq T}$ be an irreducible continuous-time Markov Chain with transition semigroup $P_t = \ee^{t \Gv}$ and invariant distribution $\piv$.  Suppose $X_0 \sim \piv$.  Let $Y$ be the time reversal of $X$ and denote by $Q_t = \ee^{t \Hv}$ its time-reversal.  We say that $X$ is \textit{reversible} if $\Pv_t = \Qv_t$ for all $t \in [0,T]$, or equivalently (by Theorem \ref{thm:generator.reversal}), if either of the following holds
\begin{align}
\pi(i) g(i,j) 
	&=	\pi(j) g(j,i) ,  &
	&\text{or} &
\pi(i) p_t(i,j)
	&=	\pi(j) p_t(j,i) , &
	&\forall \, t \in [0,T] .
\end{align}
\end{definition}

\noindent
Once again, it is instructive to compare Definition \ref{def:reversible.2} with its discrete time analog Definition \ref{def:reversible}.

\begin{definition}
\label{def:detailed.balance.2}
Let $\Pv_t = \ee^{t \Gv}$ the transition semigroup of a continuous time Markov process $X$ with state space $S$.
We say that a distribution $\lamv = (\lam(1), \lam(2), \ldots, \lam(|S|))$ is in \emph{detailed balance} with $\Pv_t$ or $\Gv$ if
\begin{align}
\lam(i) g(i,j) 
	&=	\lam(j) g(j,i) ,  &
	&\text{or} &
\lam(i) p_t(i,j)
	&=	\lam(j) p_t(j,i) , &
	&\forall \, t \in [0,T] . \label{eq:detailed.balance.2}
\end{align}
\end{definition}

\begin{theorem}
\label{theorem: rev.stat.}
Let $\Pv_t = \ee^{t \Gv}$ be transition semigroup of an irreducible Markov chain $X$.
Suppose there exists a distribution $\lamv$ such that \eqref{eq:detailed.balance.2} holds.
Then $\lamv$ is a stationary distribution for $X$.  That is, $\lamv \Gv = 0$ and $\lamv = \lamv \Pv_t$ for all $t$.
\end{theorem}

\begin{proof}
Assuming \eqref{eq:detailed.balance.2} holds, we have
\begin{align}
\sum_i \lam(i) p_t(i,j)
	&=	\sum_i \lam(j) p_t(j,i)
	=		\lam(j) \sum_i p_t(j,i)
	=		\lam(j) .
\end{align}
Thus, we have proved $\lamv = \lamv \Pv_t$, and thus $\lamv$ is n stationary distribution for $X$.
\end{proof}

\subsection{Potential and entropy production}

\begin{definition}
\label{def:cycle}
The generator $\Gv$ of a continuous time Markov chain $X$ is said to satisfy the \textit{Kolmogorov cycle condition} if, for any finite and feasible cycle of states $(s_0, s_1, \ldots, s_n, s_0)$ we have
\begin{align}
\frac{g(s_0, s_1)g(s_1, s_2) \ldots g(s_{n-1}, s_n)g(s_n, s_0)}	{g(s_1, s_0)g(s_2, s_1) \ldots g(s_n, s_{n-1})g(s_0, s_n)}
	&= 1. \label{eq:kcc}
\end{align}
Here, ``feasible'' means that $X$ can jump directly from $s_{k}$ to $s_{k+1}$ for all $k$ and also from $s_n$ to $s_0$.
\end{definition}

\begin{theorem}
\label{thm:kolmogorov.cycle.cond}
Let $X=(X_t)_{0 \leq t \leq T}$ be an irreducible continuous-time Markov Chain with 
state space $S$, transition semigroup $\Pv_t = \ee^{t\Gv}$ and stationary distribution $\piv$. Then
\begin{align}
	&\text{$\piv$ is in detailed balance with $\Gv$}&
	&\Leftrightarrow&
	&\text{$\Gv$ satisfies the Kolmogorov cycle condition}.
\end{align}
\end{theorem}

\begin{proof}
$(\Rightarrow)$ 
Suppose $\piv$ is in detailed balance with $\Gv$, then
$\pi(i) g(i ,j) = \pi(j)g(j, i)$.  Thus
\begin{align}
1
	&= \frac{\pi(s_0)g(s_0, s_1) \pi(s_1)g(s_1, s_2)\cdots \pi(s_{n-1})g(s_{n-1}, s_n) \pi(s_n) g(s_n, s_0)}
	{\pi(s_1)g(s_1, s_0)\pi(s_2)g(s_2, s_1)\cdots \pi(s_n)g(s_n, s_{n-1})\pi(s_0)g(s_0, s_n)} \\ 
	&=	\frac{g(s_0, s_1)g(s_1, s_2)\cdots g(s_{n-1}, s_n)g(s_n, s_0)}
	{g(s_1, s_0)g(s_2, s_1)\cdots g(s_n, s_{n-1})g(s_0, s_n)} . 
\end{align}
$(\Leftarrow)$ 
The Kolmogorov cycle condition \eqref{eq:kcc} is equivalent to a conservative vector field:
\begin{align}
0
	&=	\sum_{cycle \{s_k\}} \log \( \frac{g(s_{k}, s_{k+1})}{g(s_{k+1}, s_{k})}\) ,
\end{align}
Therefore, we can define a \textit{potential function} $U : S \to \Rb$ relative to a reference state $s_0 \in S$ via
\begin{align}
U(s_j) 
	&:=	\sum_{k=0}^{j-1} \log \(\frac{g(s_k, s_{k+1})}{g(s_{k+1}, s_k)} \). \label{proof:kolmogorov.cond}
\end{align}
From the potential function $U(s)$, we can construct a distribution $\lamv = (\lam(1), \lam(2), \ldots, \lam(|S|))$ for any state $s \in S$ as follows
\begin{align}
\lam(s) 
	&:= 	A \ee^{U(s)} , &
A
	&=		\Big( \sum_{s_k \in S} \ee^{U(s_k)} \Big)^{-1} .
\end{align}
Note that $A$ is just a normalization factor.  It is easy to check that $\lambda(i)$ and $\lambda(j)$ satisfy 
\begin{align}
\lambda(i)g(i,j)
	&=	\lambda(j)g(j,i) ,
\end{align}
for any states $i$ and $j$.  Therefore, by Definition (\ref{def:detailed.balance.2}), $\lamv$ is in detailed balance with $\Gv$. Additionally, by Theorem (\ref{theorem: rev.stat.}), $\lamv$ is a stationary distribution for $X$.  As $X$ is irreducible, the stationary distribution is unique, which implies $\lamv = \piv$.  Therefore, $\piv$ is in detailed balance with $\Gv$.
\end{proof}

\noindent
For any distributions $\muv$ and $\nuv$ on some countable state space $S$ we can define
\begin{align}
	&\textit{Entropy}:&
H[\muv]
	&:=	- \sum_{i \in S} \mu(i) \log \mu(i) , \\
	&\textit{Relative Entropy}:&
K[\muv||\nuv]
	&:=	\sum_{i \in S} \mu(i) \log \frac{\mu(i)}{\nu(i)} .
\end{align}
Relative Entropy is sometimes called \textit{Kullback–Leibler divergence}.

\noindent
Consider a physical system with $N \gg 1$ particles, each of which must be in some state $s \in S$.
Let us define 
\begin{align}
\mu_t(i) 
	&:= \frac{\# \{ \text{particles in state $i$ at time $t$} \}}{N}   .
\end{align}
Then $\muv_t$ is a distribution on $S$.
Suppose this system has an invariant distribution $\piv$.
We can compute the relative entropy of $\muv_t$ with respect to $\piv$.  We have
\begin{align}
\underbrace{ K[\muv_t||\piv] }_{\text{free energy of $\muv_t$}}
	&=	- \underbrace{ \Big( - \sum_{i \in S} \mu_t(i) \log \mu_t(i) \Big) }_{\text{entropy of $\muv_t$}} + \underbrace{ \Big( - \sum_{i \in S} \mu_t(i) \log \pi(i) \Big) }_{\text{mean energy of $\mu_t$}}.
\end{align}
where we have indicated the names that the above quantities are given in physics.

\begin{theorem}
Let $X$ be an irreducible Markov chain with semigroup $\Pv_t = \ee^{t \Gv}$ and invariant distribution $\piv$.  Suppose $\piv$ is in detailed balance with $\Gv$.  Then we have
\begin{align}
\frac{\dd}{\dd t} K[\muv_t||\piv] &\leq 0 .
\end{align}
That is, the relative entropy $K[\muv_t||\piv]$ is non-increasing.
\end{theorem}

\begin{proof}
Suppose $X_0 \sim \del_i$ (i.e., all mass at $i$).  Then $\mu_t(j) = p_t(i,j)$ and we have
\begin{align}
\frac{\dd}{\dd t} K[\muv_t||\piv]
	&=	\frac{\dd}{\dd t} \sum_{j \in S} \mu_t(j) \log \frac{\mu_t(j)}{\pi(j)} \\
	&=	\sum_{j \in S} \Big( 1 + \log \frac{\mu_t(j)}{\pi(j)} \Big) \frac{\dd}{\dd t} \mu_t(j) \\
	&=	\sum_{j \in S} \Big( 1 + \log \frac{p_t(i,j)}{\pi(j)} \Big) \frac{\dd}{\dd t} p_t(i,j) \\
	&=	\frac{\dd}{\dd t} \sum_{j \in S}  p_t(i,j) + \sum_{j \in S} \Big( \log \frac{p_t(i,j)}{\pi(j)} \Big) \frac{\dd}{\dd t} p_t(i,j) \\
	&=	\sum_{j \in S} \Big(  \log \frac{p_t(i,j)}{\pi(j)}\Big)   \frac{\dd}{\dd t} p_t(i,j) &
	&\text{(as $\sum_{j \in S}  p_t(i,j) = 1$)} \\
	&=	\sum_{j \in S} \Big( \log \frac{p_t(i,j)}{\pi(j)} \Big)\sum_{k\in S} \Big( p_t{(i, k)}g(k, j) \Big)  &
	&\text{(by the KFE)} \\ 
	&= \sum_{j \in S}\sum_{k\in S}  \Big( \log \frac{p_t(i,j)}{\pi(j)} \Big) \Big( p_t{(i, k)}g(k, j)
		\Big)   \\ 
	&= \sum_{j \in S} \sum_{k\in S}  \Big( \log \frac{p_t(i,j)}{\pi(j)} \Big) \Big( p_t{(i, k)}g(k, j)
		\Big)  \\ & \quad +  \underbrace{ \sum_{k\in S}  \Big( \log \frac{\pi(k)}{p_t(i.k)} \Big)p_t{(i, k)} \sum_{j \in S}g(k, j) }_{\text{this term is zero since} \  \sum_{j \in S}g(k, j) = 0 } \\	
		&= \sum_{j \in S}\sum_{k\in S}  \Big( \log \frac{p_t(i,j) \pi(k) }{ \pi(j) p_t(i, k)} \Big) \Big( p_t{(i, k)}g(k, j) \Big)  &
	& \text{(algebra)} \\ 
	 &= \sum_{j \in S} \sum_{k \in S} \Big(  \log \frac{p_t(i,j)g(j,k) }{ p_t(i,k) g(k,j)} \Big) \Big( p_t{(i, k)}g(k, j)\Big) &
	&\text{($\pi$ is in detailed balance with $\Gv$)} \\
	  &= \frac{1}{2}\sum_{j \in S} \sum_{k \in S} \Big(  \log \frac{p_t(i,j)g(j,k) }{ p_t(i,k) g(k,j)} \Big) \\ & \quad
				\times \Big( p_t{(i, k)}g(k, j)+ p_t{(i, k)}g(k, j)\Big) \\
	  &= \frac{1}{2}\sum_{j \in S} \sum_{k \in S} \Big(  \log \frac{p_t(i,j)g(j,k) }{ p_t(i,k) g(k,j)} \Big)  \\ & \quad
				\times \Big( -p_t{(i, j)}g(j, k) + p_t{(i, k)}g(k, j)\Big) &
				&\text{(algebra)} \\ 
	 &= -\frac{1}{2}\sum_{j \in S} \sum_{k \in S} \underbrace{\Big(  \log \frac{p_t(i,j)g(j,k) }{p_t(i,k) g(k,j)} \Big)}_A  \\ & \quad
				\times \underbrace{ \Big( p_t{(i, j)}g(j, k) - p_t{(i, k)}g(k, j) \Big)}_B\\
	 &=: f(i) \leq 0, &
	&\text{(as $A$ and $B$ have the same sign)} 
\end{align}
Thus, for a general starting condition $X_0 \sim \muv_0$ we have
\begin{align}
\frac{\dd}{\dd t} K[\muv_t||\piv]
	&=\sum_{i \in S} \mu_0(i) f(i) \leq 0 .
\end{align}
\end{proof}

\noindent
In statistical physics, the term $A$ is called the \textit{thermodynamic force}, the term $B$ is called {\em thermodynamic flux}, the produce $A \cdot B$ is called \textit{entropy production} or \textit{free energy dissipation} and the result $\frac{\dd}{\dd t} K[\muv_t||\piv] \leq 0$ is the Second Law of Thermodynamics.

\section{Exercises}
\label{sec:exercises-5}

\begin{exercise}
\label{exer:5.1}
Patients arrive at an emergency room as a Poisson process with intensity $\lam$.  The time to treat each patient is an independent exponential random variable with parameter $\mu$.  Let $X = (X_t)_{t \geq 0}$ be the number of patients in the system (either being treated or waiting).  Write down the generator of $X$.  Show that $X$ has an invariant distribution $\piv$ if and only if $\lam < \mu$.   Find $\piv$.  What is the total expected time (waiting + treatment) a patient waits when the system is in its invariant distribution?
\end{exercise}
%

\begin{exercise}
\label{exer:5.2}
Let $X = (X_t)_{t \geq 0}$ be a Markov chain with stationary distribution $\piv$.  Let $N$ be an independent Poisson process with intensity $\lam$ and denote by $\tau_n$ the time of the $n$th arrival of $N$.  Define $Y_n := X_{\tau_n +}$ (i.e., $Y_n$ is the value of $X$ immediately after the $n$th jump).  Show that $Y$ is a discrete time Markov chain with the same stationary distribution as $X$.
\end{exercise}
%

\begin{exercise}
\label{exer:5.3}
Let $X = (X_t)_{t \geq t}$ be a Markov chain with state space $S=\{0,1,2,\ldots\}$ and with a generator $\Gv$ whose $i$th row has entries
\begin{align}
g_{i,i-1}
	&=	i \mu , &
g_{i,i}
	&=	- i \mu - \lam , &
g_{i,i+1} 
	&=	\lam ,
\end{align}
with all other entries being zero (the zeroth row has only two entries: $g_{0,0}$ and $g_{0,1}$).  Assume $X_0 = j$.  Find $G_{X_T}(s) := \Eb s^{X_t}$.  What is the distribution of $X_t$ as $t \to \infty$?
\end{exercise}
%

\begin{exercise}
\label{exer:5.4}
Let $N$ be a time-inhomogeneous Poisson process with intensity function $\lam(t)$.  That is, the probability of a jumps of size one in the time interval $(t,t+\dd t)$ is $\lam(t) \dd t$ and the probability of two jumps in that interval of time is $\Oc(\dd t^2)$.  Write down the Kolmogorov forward and backward equations of $N$ and solve them.  Let $N_0 = 0$ and let $\tau_1$ be the time of the first jump of $N$.  If $\lam(t) = c/(1+t)$ show that $\Eb \tau_1 < \infty$ if and only if $c>1$.
\end{exercise}
%

\begin{exercise}
\label{exer:5.5}
Let $N$ be a Poisson process with a random intensity $\Lam$ with is equal to $\lam_1$ with probability $p$ and $\lam_2$ with probability $1-p$.  Find $G_{N_t}(s) = \Eb s^{N_t}$.  What is the mean and variance of $N_t$?
\end{exercise}

\begin{exercise}
\label{exer:5.6}
Let $X=(X_t)_{0 \leq t \leq T}$ be an irreducible continuous-time Markov Chain with 
state space $S$, transition semigroup $\Pv_t = \ee^{t\Gv}$ and stationary distribution $\piv$. Prove the following six statements are equivalent\\
(i) Its stationary distribution satisfies detailed balance: $\pi(i)g(i,j) = \pi(j)g(j,i), \ \forall i,j \in S$\\
(ii) Any path connecting states $i$ and $j$: $i \equiv s_0, s_1, s_2, ..., s_n \equiv j$, has a path independent
$$\log\left( \frac{g(s_0, s_1)}{g(s_1, s_0)} \right) + \log\left( \frac{g(s_1, s_2)}{g(s_2, s_1)} \right) + \cdots + \log\left( \frac{g(s_{n-1}, s_n)}{g(s_n, s_{n-1})} \right) = \log \pi(s_n) - \log \pi(s_0) $$
(iii) It defines a time reversible stationary Markov process.\\
(iv) Its $\Gv$ matrix satisfies Kolmogorov cycle condition for every sequence of states.\\
(v) There exists a positive diagonal matrix $\Pi$ such that matrix $\Gv \Pi $ is symmetric.\\
(vi) Its stationary process has zero entropy production rate. 
Entropy production rate, $e_p (\mu_t)$, is defined by
$$e_p(\mu_t) = \frac{1}{2} \sum_{i, j} \left(\mu_t{(i)} g(i,j)-\mu_t{(j)}g(j,i) \right) 
\log\left(\frac{\mu_t(i) g(i,j)}{\mu_t{(j)}g(j,i)}\right)    $$
\end{exercise}


%
%

\chapter{Convergence of random variables}
\label{chapter:convergence}
The notes from this chapter are primarily taken from \cite[Chapter 7]{gs2001}.  The goals of this chapter are (i) to introduce various modes of convergence of random variables and (ii) to establish when convergence of one mode implies convergence of another mode.


\section{A motivating example}
\label{sec:motivation}
We begin this chapter with a motivating example.  

\begin{example}
\label{ex:double-or-nothing}
Consider the infinite sequence of coin tosses described in Section \ref{sec:coin-toss}.  Suppose you have an initial wealth $X_0$.  Just before the $n^\text{th}$ coin toss, you bet your entire wealth on the outcome of the $n^\text{th}$ toss being a heads.  It is easy to see that, at time $n$ your wealth is given by
\begin{align}
X_n(\om)
	&=	\left\{ \begin{aligned}
			& 2^n X_0 && \text{if $\om_i = H$ for $i = 1, 2, \ldots n$}, \\
			& 0				&&	\text{else} .
			\end{aligned}  \right.
\end{align}
Suppose the probability of a heads is $p \in (0,1)$.  As each of the coin tosses are independent, we deduce that
\begin{align}
\Pb( X_n = 2^n X_0)
	&=	p^n , &
\Pb( X_n = 0 )
	&=	1 - p^n .
\end{align}
As $p^n \to 0$ as $n \to \infty$ we see that
\begin{align}
\lim_{n \to \infty} 
	\Pb( X_n = 0 ) = \lim_{n \to \infty} 1 - p^n = 1  &
\end{align}
Thus, it seems reasonable to say $X_n \to X_\infty$ where $X_\infty(\om) := 0$.  
Note that
\begin{align}
\lim_{n \to \infty} \Pb( |X_n - X_\infty | > \eps )
	&=	0 , &
\forall \, \eps
	&> 0 . \label{eq:conv.in.P}
\end{align}
From this computation one might expect that
$\Eb | X_n - X_\infty | \to 0$ as $n \to \infty$.  However, this is not correct, as
\begin{align}
\Eb | X_n - X_\infty |
	&=	\Eb X_n
	=		2^n X_0 \cdot p^n
	\to \left\{ \begin{aligned}
	&\infty && \text{if $p \in (1/2,1)$} , \\
	&X_0		&& \text{if $p = 1/2$}, \\
	&0			&& \text{if $p \in (0,1/2)$} .
	\end{aligned} \right.
\end{align}
Thus, we see that \eqref{eq:conv.in.P} does \textit{not} imply $\Eb | X_n - X_\infty | \to 0$.  This simple example forces us to consider more carefully what we mean by ``convergence'' of a random variable.
\end{example}


\section{Modes of convergence defined}
\label{sec:convergence.defs}
There are many (many!) ways to define convergence of a sequence of random variables.  The most important of these modes are the four we define here
\begin{definition}
Consider a sequence of random variables $(X_n)_{n \geq 0}$ defined on a probability space $(\Om,\Fc,\Pb)$.
We say that $X_n$ converges to $X_\infty$ \emph{almost surely}, written $X_n \stackrel{\text{a.s.}}{\to} X_\infty$, if
\begin{align}
 \Pb( \lim_{n \to \infty} |X_n - X_\infty| = 0) 
	&= 1 . \label{eq:as.convergence}
\end{align}
We say that $X_n$ converges to $X_\infty$ \emph{in $L^p$}, written $X_n \stackrel{L^p}{\to} X_\infty$, if $\Eb | X_n |^p < \infty$ for all $n$ and
\begin{align}
\lim_{n \to \infty} \Eb | X_n - X_\infty |^p 
	&= 0 . \label{eq:lp.convergence}
\end{align}
We say that $X_n$ converges to $X_\infty$ \emph{in probability}, written $X_n \stackrel{\Pb}{\to} X_\infty$, if
\begin{align}
\lim_{n \to \infty} \Pb( |X_n - X_\infty | > \eps ) 
	&=	0 , &
\forall \, \eps
	&> 0 . \label{eq:p.convergence}
\end{align}
And we say that $X_n$ converges to $X_\infty$ \emph{in distribution}, written $X_n \stackrel{\Dc}{\to} X_\infty$, if
\begin{align}
\lim_{n \to \infty} | F_{X_n}(x) - F_{X_\infty}(x) | 
	&= 0 , \label{eq:d.convergence}
\end{align}
for all points $x \in \Rb$ at which $F_{X_\infty}$ is continuous.
\end{definition}

\noindent
It will be helpful to comment a bit more on the above definitions.
Observe that almost sure convergence does \emph{not} require that $X_n(\om) \to X_\infty(\om)$ for \emph{every} $\om \in \Om$.  It may be helpful to define
\begin{align}
A
	&:=	\{ \om : \lim_{n \to \infty} | X_n(\om) - X_\infty(\om) | = 0 \} .
\end{align}
Almost sure convergence simply asks that $\Pb( A ) = 1$; it does not ask that $A = \Om$.  There may be infinitely many $\om \in A^c$ (uncountably infinite, in fact).  But, if $X_n \stackrel{\text{a.s.}}{\to} X_\infty$ then $\Pb(A^c) = 0$.  At times, it is helpful to write that $X_n \to X_\infty$ $\Pb$-a.s. to indicate that $X_n$ converges to $X_\infty$ almost surely under a specific probability measure $\Pb$.  Note that, that $X_n \to X_\infty$ $\Pb$-a.s. does \emph{not} imply that $X_n \to X_\infty$ $\Pbt$-a.s..  However, if two measure are equivalent $\Pb \sim \Pbt$ then almost sure convergence under $\Pb$ implies almost sure convergence under $\Pbt$ and vice-versa because $\Pb(A) = 1 \Leftrightarrow \Pbt(A)=1$.  Finally, we mention that the phrases \emph{convergence with probability one} and \emph{convergence almost everywhere} (a.e.) are sometimes used to indicate almost sure convergence.
\\[0.5em]
For any $p \geq 1$ we can define $L^p(\Om,\Fc,\Pb)$, the space of random variables $X : \Om \to \Rb$ with a finite $p^\text{th}$ moment.
\begin{align}
\Eb |X|^p
	&< \infty .
\end{align}
We can define the $p$-norm $\norm{ \cdot }_p$ of a random variable $X \in L^p(\Om,\Fc,\Pb)$ as
\begin{align}
\norm{X}_p
	&:=	\big( \Eb |X|^p \big)^{1/p} .
\end{align}
Then, from \eqref{eq:lp.convergence} we see that convergence in $L^p$ is equivalent to convergence in the $p$-norm
\begin{align}
\lim_{n \to \infty} \norm{X_n - X_\infty}_p
	&=	0 .
\end{align}
As with almost sure convergence, it is sometimes necessary to specify a particular probability measure under which $L^p$ convergence occurs by writing $X_n \to X_\infty$ in $L^p(\Pb)$ or, even more specifically, $L^p(\Om,\Fc,\Pb)$.  It is common to say $X_n \to X_\infty$ \emph{in mean} and \emph{in mean square} when $X_n \to X_\infty$ in $L^1$ and $L^2$, respectively.
\\[0.5em]
Convergence in probability is based on the following intuition: two random variables $X$ and $Y$ are ``close to each other'' if there is a low probability that their difference is very large: $\Pb( | X - Y | > \eps ) \approx 0$ for $\eps \ll 1$.  In some sense, the ``distance'' between $X$ and $Y$ is measured by $d_\eps(X,Y) := \Pb( |X-Y| > \eps)$.  Convergence in probability simply asks that this measure of distance goes to zero for any $\eps > 0$ as $n \to \infty$.  That is, $\lim_{n \to \infty} d_\eps(X_n,X_\infty) = 0$.
\\[0.5em]
Note that we have already defined convergence in distribution (see Definition \ref{def:D.convergence}).  We have simply repeated the definition here for convenience.  
In a sense, convergence in distribution is the weakest of the four modes of convergence, as the definition makes no reference to the underlying probability space $(\Om,\Fc,\Pb)$.  Indeed, suppose $Z \sim \Nc(0,1)$ and $X_n = - Z$ for all $n$.  Then it is clear that $F_{X_n} = F_Z$ for all $n$.  But, it is easy to see that $X_n$ does not converge to $Z$ in probability, in $L^p$ or almost surely.  Convergence in distribution is sometimes alternatively called \emph{weak convergence}, written $X_n \stackrel{w}{\to} X_\infty$, or \emph{convergence in law}, written $X_n \stackrel{\Lc}{\to} X_\infty$.
\\[0.5em]
Sometimes, we may wish to see if a sequence of random variables $(X_n)$ converges (in some sense) to a random variable $X_\infty$, but we do not know what the random variable $X_\infty$ is.  In such cases, it can be helpful to check for \emph{Cauchy convergence} of some sort.  

\begin{definition}
\label{def:convergence.RVs.cauchy}
Consider a sequence of random variables $(X_n)_{n \geq 0}$ defined on a probability space $(\Om,\Fc,\Pb)$.
We say that $(X_n)$ is \emph{Cauchy convergent almost surely} if
\begin{align}
 \Pb( \lim_{n,m \to \infty} |X_n - X_m| = 0) 
	&= 1 . \label{eq:as.convergence.cauchy}
\end{align}
We say that $(X_n)$ is \emph{Cauchy convergent in $L^p$} if $\Eb | X_n |^p < \infty$ for all $n$ and
\begin{align}
\lim_{n,m \to \infty} \Eb | X_n - X_m |^p 
	&= 0 . \label{eq:lp.convergence.cauchy}
\end{align}
And we say that $X_n$ is \emph{Cauchy convergent in probability} if
\begin{align}
\lim_{n,m \to \infty} \Pb( |X_n - X_m | > \eps ) 
	&=	0 , &
\forall \, \eps
	&> 0 . \label{eq:p.convergence.cauchy}
\end{align}
\end{definition}

\noindent
We will not prove the following Theorem.

\begin{theorem}
\label{thm:cauchy.convergent}
Consider a sequence of random variables $(X_n)$ defined on a probability space $(\Om,\Fc,\Pb)$ and let $X_\infty$ be a random variable defined on the same probability space.
\begin{align}
&\exists \, X_\infty \text{ s.t. } X_n \stackrel{\textup{a.s.}}{\to} X_\infty&
&\Leftrightarrow&
&\text{$(X_n)$ is Cauchy convergent almost surely.} \\
&\exists \, X_\infty \text{ s.t. } X_n \stackrel{\textup{$L^p$}}{\to} X_\infty&
&\Leftrightarrow&
&\text{$(X_n)$ is Cauchy convergent in $L^p$.} \\
&\exists \, X_\infty \text{ s.t. } X_n \stackrel{\Pb}{\to} X_\infty&
&\Leftrightarrow&
&\text{$(X_n)$ is Cauchy convergent in probability.} 
\end{align}
\end{theorem}

\noindent
The usefulness of Theorem \ref{thm:cauchy.convergent} is that it allows us to establish the existence of a limit $X_\infty$ of a sequence of random variables $(X_n)$ (in some sense) without having knowledge of what $X_\infty$ is.
\\[0.5em]
In Example \ref{ex:double-or-nothing} we saw that the wealth process $(X_n)$ converged in probability to zero: $X_n \stackrel{\Pb}{\to} 0$.  Suppose $\Pb$ was such that $\Pb(\om_i=H) = p \in (0,1/2)$.  Then the wealth process $(X_n)$ also converged in mean to zero: $X_n \to 0$ in $L^1(\Pb)$.  However, if $\Pbt(\om_i=H) = \widetilde{p} \in [1/2,1)$, then $(X_n)$ did not converge in $L^1(\Pbt)$.  It would be very helpful to know when one mode of convergence implies convergence in another mode.  This is the subject of the next section.


\section{Convergence in one mode implies convergence in another?}
Let us get right to it and state the main result.
\begin{theorem}
\label{thm:convergence.of.RVs}
Let $(X_n)$ be a sequence of random variables defined on $(\Om,\Fc,\Pb)$ and let $X_\infty$ be a random variable on the same probability space.  
Suppose $1 \leq q \leq p$.
Then the following implications hold
\begin{align}
\left. \begin{array}{r}
X_n \stackrel{a.s.}{\to} X_\infty \\
X_n \stackrel{L^p}{\to} X_\infty 
\Rightarrow
X_n \stackrel{L^q}{\to} X_\infty 
\end{array} \right\}
	&\Rightarrow
	X_n \stackrel{\Pb}{\to} X_\infty
	\Rightarrow
	X_n \stackrel{\Dc}{\to} X_\infty . \label{eq:modes}
\end{align}
No other implications hold in general.
\end{theorem}

\noindent
Note: the right bracket ``$\}$'' indicates that convergence in probability holds if $X_n$ converges to $X_\infty$ either in $L^q$ \textit{or} $\Pb$-a.s..
One does \textit{not} need for $X_n$ to converge to $X_\infty$ in both $L^q$ and $\Pb$-a.s. in order for convergence in probability to hold.
We will prove Theorem \ref{thm:convergence.of.RVs} in a series of steps.  Along the way, we will introduce some related (and useful!) inequalities.

\begin{proof}[Proof that $X_n \stackrel{\Pb}{\to} X_\infty \Rightarrow X_n \stackrel{\Dc}{\to} X_\infty$]
First we note that, if $\eps > 0$, then
\begin{align}
F_{X_n}(x)
	&=	\Pb( X_n \leq x) \\
	&=	\Pb( X_n \leq x, X_\infty \leq x + \eps) + \Pb( X_n \leq x, X_\infty > x + \eps)  \\
	&\leq	\Pb( X_\infty \leq x + \eps ) + \Pb(|X_n - X_\infty| > \eps ) \\
	&=	F_{X_\infty}(x + \eps) + \Pb(|X_n - X_\infty| > \eps ) . \label{ineq:1}
\end{align}
If you find the above inequality difficult to derive, it may help to draw a pictures of regions in the $X_n$-$X_\infty$ plane to which the above probabilities correspond.  Likewise, we have
\begin{align}
F_{X_\infty}(x-\eps)
	&=	\Pb( X_\infty \leq x - \eps ) \\
	&=	\Pb( X_\infty \leq x - \eps, X_n \leq x ) + \Pb( X_\infty \leq x - \eps, X_n > x ) \\
	&\leq \Pb( X_n \leq x) + \Pb(|X_n - X_\infty| > \eps ) \\
	&=	F_{X_n}(x) + \Pb(|X_n - X_\infty| > \eps ) . \label{ineq:2}
\end{align}
Combining \eqref{ineq:1} with \eqref{ineq:2} we obtain
\begin{align}
F_{X_\infty}(x-\eps) - \Pb(|X_n - X_\infty| > \eps )
	&\leq F_{X_n}(x)
	 \leq F_{X_\infty}(x + \eps) + \Pb(|X_n - X_\infty| > \eps ) .
\end{align}
As $X_n \stackrel{\Pb}{\to} X_\infty$ by assumption, we have that $\lim_{n \to \infty} \Pb(|X_n - X_\infty| > \eps ) = 0$.  Thus
\begin{align}
F_{X_\infty}(x-\eps)
	&\leq \operatorname*{lim inf}_{n \to \infty} F_{X_n}(x)
	 \leq \operatorname*{lim sup}_{n \to \infty} F_{X_n}(x)
	 \leq F_{X_\infty}(x + \eps) .
\end{align}
Finally, if $F_{X_\infty}(x)$ is continuous at $x$, then $\lim_{\eps \to 0} F_{X_\infty}(x \pm \eps) = F_{X_\infty}(x)$.  Therefore, letting $\eps \to 0$ we obtain
\begin{align}
F_{X_\infty}(x)
	&\leq \operatorname*{lim inf}_{n \to \infty} F_{X_n}(x)
	 \leq \operatorname*{lim sup}_{n \to \infty} F_{X_n}(x)
	 \leq F_{X_\infty}(x) ,
\end{align}
which implies that $\lim_{n \to \infty} F_{X_n}(x) = F_{X_\infty}(x)$ at any point $x$ where $F_{X_\infty}$ is continuous.  And this establishes that $X_n \stackrel{\Dc}{\to} X_\infty$.
\end{proof}

\noindent
We shall need the following result, which we state without proof.

\begin{theorem}[H\"older's inequality]
Let $X$ and $Y$ be random variables defined on $(\Om,\Fc,\Pb)$ and let $p,q \geq 1$ satisfy $1/p + 1/q = 1$.  Then we have
\begin{align}
\Eb | XY |
	&\leq ( \Eb |X|^p )^{1/p} ( \Eb |Y|^q )^{1/q} . \label{eq:holder}
\end{align}
\end{theorem}

\begin{corollary}[Lyapunov's inequality]
Let $Z$ be a random variable defined on $(\Om,\Fc,\Pb)$ and suppose $1 \leq q \leq p$.  Then we have
\begin{align}
( \Eb |Z|^q )^{1/q}
	&\leq ( \Eb |Z|^p )^{1/p} . \label{eq:lyapunov}
\end{align}
\end{corollary}

\begin{proof}
Let $r,p \geq 1$.  Applying H\"older's inequality \eqref{eq:holder} with $Y=1$ and $X = |Z|^r$ we obtain
\begin{align}
\Eb |Z|^r
	&\leq ( \Eb |Z|^{rp} )^{1/p} .
\end{align}
Defining $q := rp$, we have
\begin{align}
\Eb |Z|^r
	&\leq ( \Eb |Z|^q )^{r/q} &
	&\Rightarrow&
( \Eb |Z|^r )^{1/r}
	&\leq ( \Eb |Z|^q )^{1/q} .
\end{align}
Noting the $q \geq r$, we have proved \eqref{eq:lyapunov}.
\end{proof}

\begin{proof}[Proof that $X_n \stackrel{L^p}{\to} X_\infty$ $\Rightarrow$ $X_n \stackrel{L^q}{\to} X_\infty$ with $p \geq q \geq 1$]
Apply \eqref{eq:lyapunov} with $Z = X_n - X_\infty$, take a limit as $n \to \infty$ and use the fact that $\Eb | X_n - X_\infty |^p \to 0$ by assumption.
\end{proof}

\begin{lemma}[Markov's inequality]
Let $Z$ be a random variable defined on $(\Om,\Fc,\Pb)$.  For any $a > 0$ we have
\begin{align}
\Pb( |Z| > a )
	&\leq \tfrac{1}{a} \Eb |Z| . \label{ineq:markov}
\end{align}
\end{lemma}

\begin{proof}
By direct computation, we have
\begin{align}
\Pb( |Z| > a )
	&=	\Eb \Ib_{\{|Z| > a \}} 
	 =	\tfrac{1}{a} \Eb a \Ib_{\{|Z| > a \}} 
	\leq \tfrac{1}{a} \Eb |Z| , 
\end{align}
where we have used $|Z| \geq a \Ib_{\{|Z| > a \}}$.
\end{proof}

\begin{proof}[Proof that $X_n \stackrel{L^p}{\to} X_\infty$ $\Rightarrow$ $X_n \stackrel{\Pb}{\to} X_\infty$ for $p \geq 1$]
Fix $\eps > 0$.  Applying Markov's inequality \eqref{ineq:markov} with $Z = X_n - X_\infty$ and $a = \eps$ and taking a limit as $n \to \infty$ we obtain
\begin{align}
\lim_{n \to \infty} \Pb( | X_n - X_\infty | > \eps )
	&\leq \tfrac{1}{\eps} \lim_{n \to \infty} \Eb | X_n - X_\infty | = 0 ,
\end{align}
where, to deduce the last equality, we have used the fact that convergence in $X_n \stackrel{L^p}{\to} X_\infty$ $\Rightarrow$ $X_n \stackrel{L^1}{\to} X_\infty$.
\end{proof}

\noindent
We will not prove the following result.

\begin{theorem}[Lebesgue's dominated convergence theorem]
\label{thm:ldc}
Let $(X_n)$ be a sequence of random variables defined on $(\Om,\Fc,\Pb)$ and suppose $X_n \stackrel{\text{a.s.}}{\to} X_\infty$.  If there exists non-negative random variable $Y$ such that (i) $|X_n| \leq Y$ for every $n$ and (ii) $\Eb Y < \infty$, then
\begin{align}
\lim_{n \to \infty} \Eb X_n
	&=	\Eb \lim_{n \to \infty} X_n 
	=		\Eb X_\infty  . \label{eq:ldc}
\end{align}
\end{theorem}

\noindent
We can now prove the remaining part of Theorem \ref{thm:convergence.of.RVs}.

\begin{proof}[Proof that $X_n \stackrel{\text{a.s.}}{\to} X_\infty$ $\Rightarrow$ $X_n \stackrel{\Pb}{\to} X_\infty$]
Fix $\eps > 0$ and define $Y_n := \Ib_{\{ | X_n - X_\infty| > \eps \} }$.
As $X_n \stackrel{\text{a.s.}}{\to} X_\infty$ it follows that $Y_n \stackrel{\text{a.s.}}{\to} 0$.  Moreover, $Y_n \leq 1$ for all $n$ and $\Eb 1 < \infty$ (clearly!).  Thus, using Lebesgue's dominated convergence Theorem \eqref{eq:ldc} we obtain
\begin{align}
\lim_{n \to \infty} \Pb( | X_n - X_\infty| > \eps )
	&=	\lim_{n \to \infty} \Eb Y_n
	=		\Eb \lim_{n \to \infty} Y_n
	=		 0 ,
\end{align}
which implies $X_n \stackrel{\Pb}{\to} X_\infty$.
\end{proof}


\section{Continuity of probability measures}
\label{sec:continuity.P}
Much like we can define the limit (if it exists) of a sequence of real numbers $(x_n)$, we can define the limit (if it exists) of a sequence of sets.

\begin{definition}
We say a sequence of sets $(A_n)$ \emph{increases} to $A$, written $A_n \uparrow A$, if
\begin{align}
A_n 
	&\subseteq A_{n+1} &
	&\forall \, n \in \Nb , &
	&\text{and}&
A
	&=	\bigcup_{n=1}^\infty A_n .
\end{align}
We say a sequence of sets $(A_n)$ \emph{decreases} to $A$, written $A_n \downarrow A$, if
\begin{align}
A_n 
	&\supseteq A_{n+1} &
	&\forall \, n \in \Nb , &
	&\text{and}&
A
	&=	\bigcap_{n=1}^\infty A_n .
\end{align}
\end{definition}

\noindent
A very important property of probability measures is that they are continuous in the sense described by the following Theorem.

\begin{theorem}[Continuity of Probability Measures]
Fix a probability space $(\Om,\Fc,\Pb)$.
Suppose $(A_n)$ is a sequence of subsets of $\Fc$.  Then
\begin{align}
A_n
	&\uparrow A , &
	&\Rightarrow&
\lim_{n \to \infty} \Pb(A_n)
	&=	\Pb \big( \lim_{n \to \infty} \bigcup_{k=1}^n A_k \big) = \Pb(A) , \\
A_n
	&\downarrow A , &
	&\Rightarrow&
\lim_{n \to \infty} \Pb(A_n)
	&=	\Pb \big( \lim_{n \to \infty} \bigcap_{k=1}^n A_k \big) = \Pb(A) .
\end{align}
\end{theorem}

\begin{proof}
Suppose $A_n \uparrow A$.  Define
\begin{align}
B_n
	&:=	A_n \setminus A_{n-1} , &
	&\forall \, n \in \Nb , &
A_0
	&= \emptyset ,
\end{align}
and observe that the sets $(B_n)$ are disjoint and
\begin{align}
A_n
	&=	\bigcup_{k=1}^n B_k , &
A
	&:=	\bigcup_{n=1}^\infty A_n = \bigcup_{n=1}^\infty B_n .
\end{align}
Hence, using the countably additive property of $\Pb$ we have
\begin{align}
\Pb(A)
	&=	\sum_{k=1}^\infty \Pb(B_k)
	=		\lim_{n \to \infty} \sum_{k=1}^n \Pb(B_k)
	=	\lim_{n \to \infty} \Pb \big( \bigcup_{k=1}^n B_k \big)
	=	\lim_{n \to \infty} \Pb ( A_n ) .
\end{align}
Next, suppose $A_n \downarrow A$.  Then $A_n^c \uparrow A^c$.  Thus, we have
\begin{align}
\Pb(A)
	&=	1-\Pb(A^c) 
	=		1- \lim_{n \to \infty} \Pb ( A_n^c ) 
	=		\lim_{n \to \infty} \big( 1 - \Pb ( A_n^c ) \big)
	=		\lim_{n \to \infty} \Pb ( A_n ) .
\end{align}
\end{proof}


\section{$\lim \sup$ and $\lim \inf$ of sets}
\label{sec:liminf.limsup}
Much like we can define the $\lim \sup$ and $\lim \inf$ of a sequence of real numbers $(x_n)$, we can define the $\lim \sup$ and $\lim \inf$ of a sequence of sets.

\begin{definition}
\label{def:liminf-limsup}
Let $(A_n)$ be a sequence of sets.  We define
\begin{align}
\limsup_{n \to \infty} A_n
	&:=	\bigcap_{m=1}^\infty \bigcup_{n=m}^\infty A_n , &
\liminf_{n \to \infty} A_n
	&:=	\bigcup_{m=1}^\infty \bigcap_{n=m}^\infty A_n .
\end{align}
\end{definition}

\noindent
To get a more intuitive understanding of the above definitions, it may be helpful to note the following
\begin{align}
\limsup_{n \to \infty} A_n
	&=	\{ \om \in \Om : \forall \, m \, \exists \, n(\om) \geq m \text{ s.t. } \om \in A_{n(\om)} \} \\
	&=	\{ \om \in \Om : \om \in A_n \text{ for infinitely many $n$} \} , \\
\liminf_{n \to \infty} A_n
	&=	\{ \om \in \Om : \exists \, m(\om) \text{ s.t. } \om \in A_n \, \forall \, n \geq m(\om) \} \\
	&=	\{ \om \in \Om : \om \in A_n \text{ for all but finitely many $n$} \} .
\end{align}
For this reason, we sometimes write
\begin{align}
\limsup_{n \to \infty} A_n
	&=	\{ A_n, \text{i.o.} \} , &
\text{i.o.}
	&=	\text{``infinitely often,''} \\
\liminf_{n \to \infty} A_n
	&=	\{ A_n, \text{a.b.f.m.} \} , &
\text{a.b.f.m}
	&=	\text{``all but finitely many.''}
\end{align}
Clearly, we have
\begin{align}
\liminf_{n \to \infty} A_n
	&\subseteq \limsup_{n \to \infty} A_n .
\end{align}
It is also helpful to note that
\begin{align}
( A_n, \text{i.o})^c
	&=	(A_n^c, \text{a.b.f.m}) , &
(A_n, \text{a.b.f.m.})^c
	&=	(A_n^c, \text{i.o.}) .
\end{align}
Note that the $\lim \inf$ and $\lim \sup$ of a sequence of sets always exists (just like the $\lim \inf$ and $\lim \sup$ of a sequence of real numbers always exists).
To see this, observe that, if we define
\begin{align}
B_m
	&:= \bigcup_{n=m}^\infty A_n , &
C_m
	&:=	\bigcap_{n=m}^\infty A_n ,
\end{align}
then, for all $n \in \Nb$ we have $B_n \supseteq B_{n+1}$ and $C_n \subseteq C_{n+1}$.  We further have
\begin{align}
\limsup_{n \to \infty} A_n
	&=	\bigcap_{m=1}^\infty B_m , &
\liminf_{n \to \infty} A_n
	&=	\bigcup_{m=1}^\infty C_m .
\end{align}
Stated another way, we have
\begin{align}
B_m
	&\downarrow \limsup_{n \to \infty} A_n , &
C_n
	&\uparrow \liminf_{n \to \infty} A_n .
\end{align}
It follows by the continuity of $\Pb$  that
\begin{align}
\Pb \Big( \limsup_{n \to \infty} A_n \Big)
	&=	\lim_{n \to \infty} \Pb( B_n ) , &
\Pb \Big( \liminf_{n \to \infty} A_n \Big)
	&=	\lim_{n \to \infty} \Pb( C_n ) .
\end{align}

\begin{example}
Consider the following sequence of sets
\begin{align}
A_{2n}
	&=	B , &
A_{2n+1}
	&=	C , &
n
	&=	0, 1, 2, \ldots ,
\end{align}
where $B$ and $C$ are some arbitrary sets in a sample space $\Om$.
To compute the $\liminf_n A_n$ consider the following question: for which $\om \in \Om$ does the event $\om \in A_n$ occur all but finitely many times?  The answer is all of the $\om \in B \cap C$.  
To compute the $\limsup_n A_n$ consider the following question: for which $\om \in \Om$ does the event $\om \in A_n$ infinitely often?  The answer is all of the $\om \in B \cup C$.  
Thus, we have
\begin{align}
\liminf_{n \to \infty} A_n
	&=	B \cap C , &
\limsup_{n \to \infty} A_n
	&=	B \cup C .
\end{align}
\end{example}


\section{Borel-Cantelli Lemma}
\label{sec:borel-cantelli}

\begin{lemma}[Borel-Cantelli]
\label{thm:borel-cantelli}
Let $(A_n)$ be a sequence of sets.  Then
\begin{align}
\sum_{n=1}^\infty \Pb(A_n)
	&<\infty &
	&\Rightarrow&
\Pb(A_n, \text{i.o})
	&=	0 . \label{eq:bc1}
\end{align}
If the events $(A_n)$ are independent, then
\begin{align}
\sum_{n=1}^\infty \Pb(A_n)
	&=	\infty &
	&\Rightarrow&
\Pb(A_n, \text{i.o.})
	&= 1 . \label{eq:bc2}
\end{align}
\end{lemma}

\begin{proof}
To prove \eqref{eq:bc1} we note that
\begin{align}
\{ A_n, \text{i.o.} \}
	&=	\bigcap_{m=1}^\infty \bigcup_{n=m}^\infty A_n 
	\subseteq \bigcup_{n=m}^\infty A_n .
\end{align}
Therefore, we have
\begin{align}
\Pb(A_n, \text{i.o.})
	&\leq \Pb\Big( \bigcup_{n=m}^\infty A_n \Big)
	\leq	\sum_{n=m}^\infty \Pb( A_n ) .
\end{align}
The right-hand side above goes to zero as $m \to \infty$ if $\sum_{n=1}^\infty \Pb(A_n) < \infty$, which proves \eqref{eq:bc1}.
To prove \eqref{eq:bc2} we note that
\begin{align}
\Pb \Big( \bigcap_{n=m}^\infty A_n^c \Big)
	&=	\lim_{k \to \infty} \Pb \Big( \bigcap_{n=m}^k A_n^c \Big) \\
	&=	\prod_{n=m}^\infty \big( 1 - \Pb(A_n) \big) &
	&		\text{(by independence)}\\
	&\leq \prod_{n=m}^\infty \exp \big( - \Pb(A_n) \big) &
	&		\text{(as $1-x \leq \ee^{-x}$ if $x \geq 0$)} \\
	&=	\exp \Big( - \sum_{n=m}^\infty \Pb(A_n) \Big) .
\end{align}
The last term equals zero if $\sum_{n=1}^\infty \Pb(A_n) = \infty$, and in this case we have
\begin{align}
\Pb( \{ A_n, \text{i.o.} \}^c )
	&=	\Pb \Big( \bigcup_{m=1}^\infty \bigcap_{n=m}^\infty A_n^c \Big) \
	=		\lim_{m \to \infty} \Pb \Big( \bigcap_{n=m}^\infty A_n^c \Big) = 0 ,
\end{align}
which establishes that $\Pb( A_n, \text{i.o.}) = 1$, as claimed.
\end{proof}

\noindent
We can use the Borel-Cantelli Lemma to establish that \emph{fast} convergence in probability implies almost sure convergence.

\begin{theorem}
Suppose $X_n \stackrel{\Pb}{\to} X_\infty$ and
\begin{align}
\sum_{n=1}^\infty \Pb( |X_n - X_\infty| > \eps )
	&< \infty , &
\forall \, \eps
	&> 0 .
\end{align}
Then $X_n \stackrel{\textup{a.s.}}{\to} X_\infty$.
\end{theorem}

\begin{proof}
Fix $\eps > 0$ and take $A_n = \{ |X_n - X_\infty| > \eps \}$.  
We have by assumption that $\sum_n \Pb( A_n ) < \infty$.  
It follows from the Borel-Cantelli Lemma \ref{thm:borel-cantelli}, that
\begin{align}
\Pb( A_n, \text{i.o.} ) =
\Pb( |X_n - X_\infty| > \eps , \text{i.o.} )
	&=	0 .
\end{align}
Thus, $X_n \stackrel{\text{a.s.}}{\to} X_\infty$.
\end{proof}


\section{Martingale convergence theorems}
\label{sec:martingale.convergence}
It is often the case that a sequence of random variables $(X_n)$ represents the value of a process at discrete periods of time.  This was the case in Example \ref{ex:double-or-nothing}, in which $X_n$ represented the wealth of a gambler after the $n^\text{th}$ flip of a coin.  A very special class of processes that we have previously encountered is the class of \emph{martingale} processes (see Definition \ref{def:martingale}).  For the reader's convenience, we restate the definition of a martingale for a discrete-time process here.

\begin{definition}
\label{def:martingale.2}
Let $X=(X_n)_{n \in \Nb_0}$ be a sequence of random variables defined on a probability space $(\Om,\Fc,\Pb)$.  The sequence $X$ is said to be a \emph{martingale} with respect to a filtration $(\Fc_n)_{n \in \Nb_0}$ if the following hold:\\
(i) $\Eb |X_n| < \infty$ for all $n \in \Nb_0$,\\
(ii) $\Eb[ X_{n+m} | \Fc_n ] = X_n$ for all $n,m \in \Nb_0$.  Alternatively, $\Eb[ X_{n+1} | \Fc_n ] = X_n$ for all $n \in \Nb_0$.
\end{definition}

\noindent
The special structure of martingales will allow us to prove some very powerful theorems about the behavior of a martingale $X$ as $n \to \infty$.  Before stating and proving these theorems, however, let us take a look at some example of discrete-time martingales.

\begin{example}[Martingales from branching processes]
Suppose $Z_n$ is the size of the $n^\text{th}$ generation of a branching process with $Z_0 = 1$ (see Section \ref{sec:branching}).  We have
\begin{align}
Z_{n+1} 
	&=	\sum_{i=1}^{Z_n} X_{n,i} ,
\end{align}
where the $(X_{n,i})$ are i.i.d with common distribution $X$.  Let $\Fc_n := \sig(Z_0,Z_1,\ldots,Z_n)$.  We have
\begin{align}
\Eb [ Z_{n+1} | \Fc_n ]
	&=	\Eb [ \sum_{i=1}^{Z_n} X_{n,i} | Z_n ] 
	=		Z_n \Eb X
	=	Z_n \mu , &
\mu
	&:=	\Eb X .
\end{align}
From this, one obtains by induction that $\Eb Z_n = \mu^n$.  Clearly the process $Z$ is not a martingale when $\mu \neq 1$.  However, consider $W_n := \mu^{-n}Z_n$.  We have
\begin{align}
\Eb [ W_{n+1} | \Fc_n ]
	&=	\mu^{-(n+1)} \Eb [ Z_{n+1} | \Fc_n ] 
	=	\mu^{-(n+1)} Z_n \mu 
	=	\mu^{-n} Z_n 
	=		W_n .
\end{align}
Thus, the process $W = (W_n)$ is a martingale with respect to $(\Fc_n)$.  Next, suppose $\eta$ is the smallest solution of
\begin{align}
G(\eta)
	&=	\eta , &
	&\text{where}&
G(s)
	&:=	\Eb s^X .
\end{align}
Defining $V_n := \eta^{Z_n} $, we find
\begin{align}
\Eb[ V_{n+1} | \Fc_n ]
	&=	\Eb[ \eta^{Z_{n+1}} | \Fc_n ]
	=		\Eb[ \eta^{\sum_{i=1}^{Z_n}X_{n,i}} | \Fc_n ]
	=		(\Eb \eta^X)^{Z_n} \\
	&=	( G(\eta) )^{Z_n}
	=		\eta^{Z_n}
	=		V_n .
\end{align}
Thus, the process $V = (V_n)$ is a martingale with respect to $(\Fc_n)$.
\end{example}

\begin{example}[Martingales from Markov chains]
Let $X = (X_n)$ be a discrete time Markov chain with state space $S$ and one-step transition matrix $\Pv = (p(i,j))$.  Suppose $\psiv$ is an a right eigenvector of $\Pv$ with corresponding eigenvalue $\lam = 1$.  That is
\begin{align}
\Pv \psiv
	&=	\psiv , &
	&\text{or component-wise}&
\sum_{j \in S} p(i,j) \psi(j)
	&=	\psi(i) , &
\forall \, i
	&\in S .
\end{align}
Define $Y_n := \psi(X_n)$.  Letting $\Fc_n := \sig(X_0,X_1,\ldots,X_n)$, we have
\begin{align}
\Eb[ Y_{n+1} | \Fc_n ]
	&=	 \Eb[ \psi(X_{n+1}) | X_n ]
	=		\sum_{j \in S} p(X_n,j) \psi(j)
	=		\psi(X_n)
	=		Y_n .
\end{align}
This, the process $Y = (Y_n)$ is a martingale with respect to $(\Fc_n)$.
\end{example}

\noindent
To prove the main result of this section (Theorem \ref{thm:martingale.convergence}) we require the following result:

\begin{theorem}[Doob-Kolmogorov inequality]
\label{thm:doob-kolmogorov}
Let $M = (M_n)$ be a martingale with respect to a filtration $(\Fc_n)$.  Then
\begin{align}
\Pb( \sup_{0 \leq i \leq n} |M_i| \geq \eps )
	&\leq	\tfrac{1}{\eps^2} \Eb M_n^2 , &
\forall \, \eps
	&>	0 . \label{ineq:doob-kolmogorov}
\end{align}
\end{theorem} 

\begin{proof}
We begin by defining the following sets
\begin{align}
A_0
	&:=	\Om , &
A_k
	&:=	\{ \sup_{0 \leq i \leq k} |M_i| < \eps \} , &
B_k
	&:=	A_{k-1} \cap \{ |M_k| \geq \eps \} , &
k
	&\geq 1 .
\end{align}
Observe that 
$A_k$ is the event that $|M|$ does not reach $\eps$ prior to time $k$ and
$B_k$ is the event that $|M|$ first reaches or exceeds $\eps$ for the first time at time $k$.  Thus, for any $k$, the events $A_k, B_1, B_2, \ldots, B_k$ are a partition of $\Om$.  That is, for all $k$, we have
\begin{align}
\Om
	&=	A_k \cup \Big( \bigcup_{i=1}^k B_i \Big) , &
A_k \cap B_i
	&=	\emptyset , &
i 
	&\leq k , &
B_i \cap B_j
	&=	\emptyset , &
i
	&\neq j .
\end{align}
It follows that
\begin{align}
\Eb M_n^2
	&=	\Eb M_n^2 \Ib_{A_n} + \sum_{i=1}^n \Eb M_n^2 \Ib_{B_i} 
	\leq \sum_{i=1}^n \Eb M_n^2 \Ib_{B_i} .
\end{align}
Moreover, we have
\begin{align}
\Eb M_n^2 \Ib_{B_i}
	&=	\Eb (M_n - M_i + M_i)^2 \Ib_{B_i} \\
	&=	\underbrace{ \Eb (M_n - M_i)^2 \Ib_{B_i} }_{\alpha}
			+ \underbrace{ 2 \Eb (M_n - M_i) M_i \Ib_{B_i} }_{2 \beta}
			+ \underbrace{ \Eb M_i^2 \Ib_{B_i} }_{\gam} .
\end{align}
Let us look at these three terms one-by-one.
Clearly, the first term satisfies $\alpha \geq 0$.  Next, using the fact that $M$ is a martingale, the second term satisfies
\begin{align}
\beta
	&=	\Eb (M_n - M_i) M_i \Ib_{B_i} 
	=		\Eb M_i \Ib_{B_i} \Eb[ (M_n - M_i) | \Fc_i ]
	=		0 .
\end{align}
And the last term satisfies
\begin{align}
\gam
	&=	\Eb M_i^2 \Ib_{B_i}
	\geq \Eb \eps^2 \Ib_{B_i}
	=		\eps^2 \Pb( B_i ) .
\end{align}
Using the above results, we have
\begin{align}
\Eb M_n^2 \Ib_{B_i}
	&=	\alpha + 2 \beta + \gam
	\geq \eps^2 \Pb( B_i ) 
\end{align}
and hence
\begin{align}
\Eb M_n^2
	&\leq \sum_{i=1}^n \Eb M_n^2 \Ib_{B_i} 
	\leq \sum_{i=1}^n  \eps^2 \Pb( B_i ) 
	=	\eps^2 \Pb( \sup_{0 \leq i \leq n} |M_i| \geq \eps ) ,
\end{align}
which proves \eqref{ineq:doob-kolmogorov}.
\end{proof}

\noindent
We are now ready to prove the main result of this section:

\begin{theorem}
\label{thm:martingale.convergence}
Let $M = (M_n)$ be a martingale with respect to a filtration $(\Fc_n)$.  Suppose that $\sup_n \Eb M_n^2 < \infty$.  Then there exists a random variable $M_\infty$ such that 
\begin{align}
M_n
	&\stackrel{\textup{a.s.}}{\to} M_\infty , &
	&\text{and}&
M_n
	&\stackrel{L^2}{\to} M_\infty .
\end{align}  
\end{theorem}

\begin{proof}
Observe that, for any $n,m \in \Nb_0$, we have
\begin{align}
\Eb M_{n+m}^2
	&=	\Eb ( M_{n+m} - M_n + M_n )^2 \\
	&=		\Eb ( M_{n+m} - M_n)^2  + \Eb M_n^2 + 2 \Eb M_n ( M_{n+m} - M_n) \\
	&=	\Eb ( M_{n+m} - M_n)^2  + \Eb M_n^2 + 2 \Eb M_n \Eb[ ( M_{n+m} - M_n) | \Fc_n ] \\
	&=		\Eb ( M_{n+m} - M_n)^2  + \Eb M_n^2
	\geq  \Eb M_n^2 .
\end{align}
Thus, the sequence $(\Eb M_n^2)$ is strictly increasing and (by assumption) bounded.  The sequence therefore has a limit $\ell := \lim_{n \to \infty} \Eb M_n^2$.
We will now show that the sequence $(M_n)$ is Cauchy convergent almost surely (see Definition \ref{def:convergence.RVs.cauchy}).  By Theorem \ref{thm:cauchy.convergent}, this will establish that $M_n \stackrel{\text{a.s.}}{\to} M_\infty$.  Define $C$ as the set of $\om \in \Om$ for which the sequence $(M_n)$ is Cauchy convergent almost surely
\begin{align}
C
	&:= \{ \forall \eps > 0 \, \exists \, m \, \text{s.t.} |M_{m+i} - M_{m+j} | < \eps \, \forall \, i,j \in \Nb \} \\
	&=	\bigcap_{\eps > 0} \bigcup_{m=1}^\infty \{ |M_{m+i} - M_m | < \eps \, \forall \, i \in \Nb \}
\end{align}
The complement of $C$ is given by
\begin{align}
C^c
	&=	\bigcup_{\eps > 0} \bigcap_{m=1}^\infty \{ \exists \, i \in \Nb \, \text{s.t.} \, |M_{m+i} - M_m | \geq \eps \} 
	= 	\bigcup_{\eps > 0} \bigcap_{m=1}^\infty A_m(\eps) , \\
A_m(\eps)
	&:=	\{ \exists \, i \in \Nb \, \text{s.t.} \, |M_{m+i} - M_m | \geq \eps \} .
\end{align}
Noting that $A_m(\eps) \subseteq A_m(\eps')$ if $\eps \geq \eps'$, it follows  that
\begin{align}
\Pb(C^c)
	&=	\lim_{\eps \searrow 0} \Pb \Big( \bigcap_{m=1}^\infty A_m(\eps) \Big)
	\leq	\lim_{\eps \searrow 0} \lim_{m \to \infty} \Pb( A_m(\eps) ) .
\end{align}
Thus, to prove that $\Pb(C^c)=0$, it is sufficient to show that $\Pb(A_n(\eps)) \to 0$ as $n \to \infty$ for all $\eps > 0$.  
To this end, for a fixed $m \in \Nb$ define the sequence $Y = (Y_n)$ by $Y_n := M_{n+m} - M_m$.  The process $(Y_n)$ is a martingale with respect to $(\Gc_n)$  where $\Gc_n = \sig(Y_1, Y_2, \ldots, Y_n)$ since
\begin{align}
\Eb[ Y_{n+1} | \Gc_n ]
	&=	\Eb[  M_{n+1+m}- M_m  | \Gc_n ]
	=		\Eb[ \Eb[ M_{n+1+m}- M_m| \Fc_{n+m}] | \Gc_n ] \\
	&=	\Eb[ M_{n++m} - M_m | \Gc_n ]
	=		\Eb[ Y_n | \Gc_m ]
	=		Y_n .
\end{align}
As the process $Y$ is a martingale, we have by the Doob-Kolmogorov inequality (Theorem \ref{thm:doob-kolmogorov}) that
\begin{align}
&&
\Pb( \sup_{0 \leq i \leq n} |Y_i| \geq \eps )
	&\leq	\tfrac{1}{\eps^2} \Eb Y_n^2 , &
\forall \, \eps
	&>	0 , \\
&\Rightarrow&
\Pb( \sup_{0 \leq i \leq n} |M_{m+i}-M_m| \geq \eps )
	&\leq	\tfrac{1}{\eps^2} \Eb (M_{m+n}-M_n)^2 \\ & &
	&=		\tfrac{1}{\eps^2} ( \Eb M_{m+n}^2 - \Eb M_m^2), &
\forall \, \eps
	&>	0 .
\end{align} 
Taking a limit as $n \to \infty$ we obtain
\begin{align}
&&
\Pb( A_m(\eps))
	&\leq \tfrac{1}{\eps^2}( \ell - \Eb M_m^2 ) , \\
&\Rightarrow&
\lim_{m \to \infty} \Pb( A_m(\eps))
	&\leq \tfrac{1}{\eps^2}( \ell - \ell ) = 0 ,
\end{align}
which establishes that $(M_n)$ is Cauchy convergent almost surely.  Hence, there exists $M_\infty$ such that $M_n \stackrel{\text{a.s.}}{\to} M_\infty$.
Next using Fatou's Lemma
\footnote{
\emph{Fatou's Lemma} states that, for a sequence of non-negative random variables $(X_n)$ the following holds:
$\Eb \operatorname*{lim\,inf}_{n \to \infty} X_n \leq \operatorname*{lim\,inf}_{n \to \infty} \Eb X_m$.
}
we have
\begin{align}
\Eb (M_n - M_\infty)^2
	&=	\Eb \operatorname*{lim\,inf}_{m \to \infty} (M_n - M_m)^2
	\leq \operatorname*{lim\,inf}_{m \to \infty} \Eb  (M_n - M_m)^2 \\
	&=		\operatorname*{lim\,inf}_{m \to \infty} \Eb  (M_n - M_m)^2 
	=	\lim_{m \to \infty} \Eb M_m^2 - \Eb M_n^2 \\
	&=	\ell - \Eb M_n^2 \to 0 \quad \text{as $n \to \infty$,}
\end{align}
which establishes that $M_n \stackrel{L^2}{\to} M_\infty$.
\end{proof}

\noindent
Let us see an example of how Theorem \ref{thm:martingale.convergence} can be applied.

\begin{example}
Suppose $(X_n)$ is a Markov chain with one-step transition probabilities given by
\begin{align}
p(i,j)
	&=	\Pb( X_{n+1} = j | X_n = i)
	=		\binom{N}{j} \( \frac{i}{N} \)^j \( 1 - \frac{i}{N} \)^{N-j} , &
i,j
	&\in \{0,1,\ldots,N\} .
\end{align}
Let $\Fc_n = \sig(X_1,X_2,\ldots,X_n)$.  The process $(X_n)$ is a martingale with respect to the filtration $(\Fc_n)$ because
\begin{align}
\Eb[ X_{n+1} | \Fc_n ]
	&=	\sum_{j=0}^N j \,  p(X_n,j) \\
	&=	\sum_{j=0}^N j \binom{N}{j} \( \frac{X_n}{N} \)^j \( 1 - \frac{X_n}{N} \)^{N-j} \\
	&=	\( 1 - \frac{X_n}{N} \)^{-N} \( 1 - \frac{X_n}{N} \)^N X_n
	=		X_n .
\end{align}
Clearly, as $(X_n)$ is bounded, we have $\sup_n \Eb X_n^2 < \infty$.  As such, by Theorem \ref{thm:martingale.convergence} there exists a random variable $X_\infty$ such that $X_n \stackrel{\text{a.s.}}{\to} X_\infty$ and $X_n \stackrel{L^2}{\to} X_\infty$.
To find the distribution of $X_\infty$, we note that both $\{0\}$ and $\{N\}$ are absorbing states.  As there are no other absorbing states, we conclude that $(X_n)$ will eventually end up in one of these two states.  Thus, using the martingale property $X_0 = \Eb X_n$ we find
\begin{align}
X_0
	&=	\lim_{n \to \infty} \Eb X_n
	=		\Eb \lim_{n \to \infty} X_n 
	=		\Eb X_\infty \\
	&=	N \cdot \Pb( X_\infty = N ) + 0 \cdot \Pb( X_\infty = 0 ) ,
\end{align}
where passing the limit through the expectation is allowed by Lebesgue's dominated convergence (Theorem \ref{thm:ldc}).  From the above computation, we conclude that
\begin{align}
X_0/N 
	&=	\Pb( X_\infty = N ) =  1 - \Pb( X_\infty = 0 ) .
\end{align}
\end{example}


\section{Uniform Integrability}
\label{sec:ui}
We begin this section with a motivating example.

\begin{example}
\label{ex:not-UI}
Suppose $(\Om,\Fc,\Pb) = ([0,1],\Bc([0,1]),\text{Leb})$ so that $\Pb(\dd \om) = \dd \om$.  Define a sequence of random variables $(X_n)$ by
\begin{align}
X_n(\om)
	&:=	n \Ib_{\{[0,1/n]\}}(\om) .
\end{align}
It is easy to see that $X_n \stackrel{\Pb}{\to} 0$ as, for any $\eps > 0$, we have
\begin{align}
\lim_{n \to \infty} \Pb( |X_n - 0| > \eps )
	&=	\lim_{n \to \infty} (1/n) = 0 .
\end{align}
In fact, in this case, we have $X_n \stackrel{\text{a.s.}}{\to} 0$ as well.
However, we clearly do not have $L^1$ convergence of $X_n \to 0$ because
\begin{align}
\lim_{n \to \infty} \Eb |X_n - 0|
	&= \lim_{n \to \infty} 1 = 1 .
\end{align}
It is natural to ask: are there conditions under which $X_n \stackrel{\Pb}{\to} X_\infty$ implies $X_n \stackrel{L^1}{\to} X_\infty$?  It turns out that the answer is ``yes.''  We introduce these conditions in this section.
\end{example}

\begin{definition}
\label{def:UI}
A collection of random variables $(X_n)$ (possibly uncountable) is said to be \emph{uniformly integrable} (UI) if, for every $\eps > 0$ there exists $K_\eps \in [0,\infty)$ such that
\begin{align}
\sup_{n} \Eb |X_n| \Ib_{\{|X_n|>K_\eps\}}
	< \eps .
\end{align}
\end{definition}

\begin{example}
Let us return to Example \ref{ex:not-UI}, discussed above.  For every $n$ we have $\Eb | X_n | = 1$.  Thus, the collection of random variables $(X_n)$ is integrable.  However, we have
\begin{align}
\Eb |X_n| \Ib_{\{|X_n|>K\}}
	&=	1 , &
\forall \, n
	&> K .
\end{align}
Thus, for a given $\eps$ there is no $K_\eps \in [0,\infty)$ for which $\sup_n \Eb |X_n| \Ib_{\{|X_n|>K_\eps \}} < \eps$.  As such, we conclude that the collection $(X_n)$ is \emph{not} UI.
\end{example}

\noindent
Let us provide some simple conditions which, if satisfied by a collection of random variables, ensure uniform integrability.

\begin{lemma}
Suppose a collection of random variables $(X_n)$ satisfies
\begin{align}
\sup_{n} \Eb |X_n|^p
	< \infty
\end{align}
for some $p>1$.  Then the collection $(X_n)$ is UI.
\end{lemma}

\begin{proof}
Assuming $K>0$ and $p>1$ we have $|X_n|^{1-p} \Ib_{\{|X_n| > K \}}  \leq K^{1-p} \Ib_{\{|X_n| > K \}}$ and thus $|X_n| \Ib_{\{|X_n| > K \}}  \leq K^{p-1} |X_n|^p \Ib_{\{|X_n| > K \}} $.  Hence
\begin{align}
\Eb |X_n| \Ib_{\{|X_n|>K\}}
	&\leq K^{1-p} \Eb |X_n|^p \Ib_{\{|X_n|>K\}}
	\leq K^{1-p} \Eb |X_n|^p .
\end{align}
Taking $\sup$'s on both sides, we have
\begin{align}
\sup_{n} \Eb |X_n| \Ib_{\{|X_n|>K\}}
	&\leq K^{1-p}  \sup_{n} \Eb |X_n|^p  =: K^{1-p} c , &
c
	&:=	\sup_{n} \Eb |X_n|^p < \infty
\end{align}
The right-hand side can be made smaller than $\eps$ by choosing $K  > K_\eps := (c/\eps)^{1/(p-1)}$.  Thus, we have that $(X_n)$ is UI.
\end{proof}

\begin{lemma}
Suppose a collection of random variables $(X_n)$ satisfies
\begin{align}
\sup_{n} |X_n| 
	&\leq Y , &
\Eb Y
	&< \infty ,
\end{align}
for some random variable $Y$.  Then the collection $(X_n)$ is UI.
\end{lemma}

\begin{proof}
By assumption $|X_n| \leq Y$, which implies $ \Ib_{\{|X_n|>K\}} \leq  \Ib_{\{Y>K\}}$.  Hence, 
we have
\begin{align}
\sup_{n} \Eb |X_n| \Ib_{\{|X_n|>K\}}
	&\leq \Eb Y \Ib_{\{Y>K\}}  .
\end{align}
Because $\Eb Y < \infty$ by assumption, for any $\eps > 0$, there exists $K_\eps$ such that $\Eb Y \Ib_{\{Y>K_\eps\}} < \eps$.
Thus, by choosing $K > K_\eps$ above, we can ensure $\sup_n \Eb |X_n| \Ib_{\{|X_n|>K\}} < \eps$.  It follows that the collection of random variables $(X_n)$ is UI.
\end{proof}

\begin{theorem}[Bounded convergence]
\label{thm:bdd}
Let $(X_n)$ be a sequence of random variables and let $X_\infty$ be a random variable.  Suppose that $X_n \stackrel{\Pb}{\to} X_\infty$ and that, for some $K<\infty$ we have
\begin{align}
\sup_n |X_n| \leq K .
\end{align}
Then $X_n \stackrel{L^1}{\to} X_\infty$.
\end{theorem}

\begin{proof}
First, we check that $\Pb( |X_\infty| \leq K ) = 1$.  For any $k \in \Nb$ we have
\begin{align}
\Pb( |X_\infty| > K + 1/k) 
	&\leq \Pb( |X_\infty-X_n| > 1/k) .
\end{align}
The right-hand side goes to zero as $n \to \infty$ because (by assumption) $X_n \stackrel{\Pb}{\to} X_\infty$.  Thus, we have $\Pb( |X| > K + 1/k) = 0$.  Now, observe that
\begin{align}
\Pb( |X_\infty| > K )
	&=	\Pb \Big( \bigcup_k \{ |X_\infty| > K + 1/k \} \Big)
	=		\lim_{k \to \infty} \Pb ( |X_\infty| > K + 1/k )
	= 	0 .
\end{align}
Hence, $\Pb( |X_\infty| \leq K ) = 1$.  Next, because $X_n \stackrel{\Pb}{\to} X_\infty$, for any $\eps>0$ we can choose $n_\eps$ such that
\begin{align}
\Pb( |X_n - X_\infty | > \eps/3) 
	&< \eps/(3K) , &
\forall \, n
	&\geq n_\eps .
\end{align}
Then, for $n \geq n_\eps$, we have
\begin{align}
\Eb |X_n - X_\infty|
	&=	\Eb |X_n - X_\infty| \Ib_{\{|X_n - X_\infty|>\eps/3\}} + \Eb |X_n - X_\infty| \Ib_{\{|X_n - X_\infty| \leq \eps/3\}} \\
	&\leq 2 K \Pb( |X_n - X_\infty| > \eps/3) + \eps/3 \leq \eps ,
\end{align}
where we have used $|X_n - X_\infty| \leq 2 K$.  As $\eps$ was arbitrary, we can make $\Eb |X_n - X_\infty|$ as small as we like.  Hence, $X_n \stackrel{L^1}{\to} X_\infty$.
\end{proof}

\noindent
We can now state and prove the main result of this section.

\begin{theorem}
\label{thm:UI}
Let $(X_n)$ be a sequence of random variables and let $X_\infty$ be a random variable.  Suppose that $\Eb |X_n| < \infty$ for every $n$ and $\Eb |X_\infty | < \infty$.  Then
\begin{align}
&X_n \stackrel{L^1}{\to} X_\infty&
&\Leftrightarrow&
&\text{$X_n \stackrel{\Pb}{\to} X_\infty$ and $(X_n)$ is UI}.
\end{align}
\end{theorem}

\begin{proof}
We have already proved the $\Rightarrow$ part
\footnote{More specifically, we have proved that convergence in $L^1$ implies convergence in probability, but it also implies uniform integrability.}
(see Theorem \ref{thm:convergence.of.RVs}).  So, we focus now on the $\Leftarrow$ part.  To this end, suppose $X_n \stackrel{\Pb}{\to} X_\infty$ and $(X_n)$ is UI.  For any $K \in [0,\infty)$, define $\phi_K : \Rb \to [-K,K]$ as follows
\begin{align}
\phi_K(x)
	&:= \left\{ \begin{aligned}
			& -K && x < -K , \\
			& x && |x| \leq K ,\\
			& K && x > K .
			\end{aligned} \right.
\end{align}
Observe that
\begin{align}
|\phi_K(x) - x| = (|x|-K) \Ib_{\{|x|>K\}} \leq |x| \Ib_{\{|x|>K\}} . \label{ineq:phi}
\end{align}
Fix $\eps > 0$.  From \eqref{ineq:phi} and the fact that $(X_n)$ is UI and $X_\infty$ integrable, there exists a $K_\eps$ such that
\begin{align}
\sup_n \Eb | \phi_{K_\eps}(X_n) - X_n | &< \eps/3 , &
\Eb | \phi_{K_\eps}(X_\infty) - X_\infty | &< \eps/3 .
\end{align}
Now, $X_n \stackrel{\Pb}{\to} X_\infty$ implies that $\phi_{K_\eps}(X_n) \stackrel{\Pb}{\to} \phi_{K_\eps}(X_\infty)$.  And furthermore, $\sup_n|\phi_{K_\eps}(X_n)| \leq {K_\eps}$.  Thus, it follows from Theorem \ref{thm:bdd} that $\phi_{K_\eps}(X_n) \stackrel{L^1}{\to} \phi_{K_\eps}(X_\infty)$.  Thus, there exists an $n_\eps$ such that
\begin{align}
\Eb | \phi_{K_\eps}(X_n) - \phi_{K_\eps}(X_\infty) | 
	&< \eps/3 , &
\forall \, n
	&\geq n_\eps .
\end{align}
Combining the above results, we obtain
\begin{align}
\Eb | X_n - X_\infty |
	&\leq \Eb | X_n - \phi_{K_\eps}(X_n) | + \Eb | \phi_{K_\eps}(X_n) - \phi_{K_\eps}(X_\infty) | + \Eb | \phi_{K_\eps}(X_\infty) - X_\infty | 
	 < \eps .
\end{align}
As $\eps$ was arbitrary, it follows that $X_n \stackrel{L^1}{\to} X_\infty$, as claimed.
\end{proof}

\begin{remark}
From Theorem \ref{thm:convergence.of.RVs}, we know that $X_n \stackrel{\text{a.s.}}{\to} X_\infty$ implies $X_n \stackrel{\Pb}{\to} X_\infty$.  From this and Theorem \ref{thm:UI}, it follows that
\begin{align}
&\text{$X_n \stackrel{\text{a.s.}}{\to} X_\infty$ and $(X_n)$ is UI}
&\Rightarrow&
&X_n \stackrel{L^1}{\to} X_\infty.
\end{align}
\end{remark}

\noindent
We will show below that the set of conditional expectations of an integrable random variable form a class of uniformly integrable random variables.  To prove this result, we require the following lemma.

\begin{lemma}
\label{thm:L1.eps.del}
Suppose $\Eb |X| < \infty$.  Then, for all $\eps > 0$ there exists $\del_\eps > 0$ such that $\Pb(A) < \del_\eps$ implies $\Eb |X|\Ib_A < \eps$.
\end{lemma}

\begin{proof}
Suppose, by contradiction, that $\Eb |X| < \infty$ and, for a given $\eps>0$, there exists a sequence of sets $(A_n)$ such that
\begin{align}
\Pb(A_n) 
	&< 1/2^n , &
	&\text{and}&
\Eb |X|\Ib_{A_n}
	&\geq \eps ,&
	&\forall \, n .
\end{align}
Define $A := \{ A_n, \text{i.o.} \}$.  As $\sum_{n=1}^\infty \Pb(A_n) < \infty$, we have by the Borel-Cantelli Lemma (see \eqref{eq:bc1}) that $\Pb(A) = 0$ and thus $\Eb |X| \Ib_A = 0$.  On the other hand we have
\begin{align}
\eps
	&\leq \limsup_{n \to \infty} \Eb |X| \Ib_{A_n}
	\leq \Eb \limsup_{n \to \infty} |X| \Ib_{A_n}
	=	\Eb |X| \Ib_A ,
\end{align}
which is a contradiction.
\end{proof}

\begin{theorem}
\label{thm:ui.cond.exp}
Let $X$ be a random variable on $(\Om,\Fc,\Pb)$ and suppose $\Eb |X|<\infty$.  Then the collection $\Cc$ of random variables defined by
\begin{align}
\Cc 
	&=	\{ \Eb(X|\Gc) : \text{$\Gc$ is a sub-$\sig$-algebra of $\Fc$} \} ,
\end{align}
is uniformly integrable.
\end{theorem}

\begin{proof}
Fix $\eps > 0$.  From Lemma \ref{thm:L1.eps.del} there exists $\del_\eps > 0$ such that $\Pb(A) < \del_\eps$ implies $\Eb|X| \Ib_A < \eps$.
As $\Eb |X|< \infty$ we can choose $K_\eps > 0$ so that $\Eb |X| / K_\eps < \del_\eps$.  Now, for any $\Gc \subseteq \Fc$, we can define $Y = \Eb(X|\Gc)$.  By Jensen's inequality (see Theorem \ref{thm:jensen}) we have
\begin{align}
|Y|
	&=	| \Eb(X|\Gc) |
	\leq \Eb(|X|\,|\Gc) , &
	& \Pb-\text{a.s.}. \label{ineq:Abs.Y}
\end{align}
Hence, we have $\Eb |Y| \leq \Eb \Eb( |X| | \Gc) = \Eb |X|$ and, by the Markov inequality
\begin{align}
K_\eps \Pb( |Y| > K_\eps)
	\leq \Eb |Y| \leq \Eb |X| ,
\end{align}
so that
\begin{align}
\Pb( |Y| > K_\eps ) 
	&\leq \Eb |X| / K_\eps < \del_\eps .
\end{align}
Note that $\{ |Y| > K_\eps \} \in \Gc$.  Thus, from \eqref{ineq:Abs.Y} we have
\begin{align}
\Eb |Y| \Ib_{\{|Y|> K_\eps\}}
	&\leq \Eb \Eb(|X|\,|\Gc) \Ib_{\{|Y|> K_\eps \}} \\
	&=	\Eb \Eb(|X| \Ib_{\{|Y|> K_\eps\}}|\Gc) 
 	 = 	\Eb |X| \Ib_{\{|Y|> K_\eps\}}
	< \eps . 
\end{align}
As the above inequality holds for any $\Gc \subseteq \Fc$, the class $\Cc$ is UI, as claimed.
\end{proof}

\begin{example}[Uniformly integrable martingale]
Suppose $X_\infty$ is a random variable defined on $(\Om,\Fc,\Pb)$ and $\Eb |X_\infty| < \infty$.  Suppose $(\Fc_n)$ is a filtration and define a sequence of random variables $(X_n)$ by $X_n := \Eb[ X_\infty | \Fc_n ]$.  Then the process $X = (X_n)$ is martingale because
\begin{align}
\Eb[ X_{n+1} | \Fc_n ]
	&=	\Eb[ \Eb[ X_\infty | \Fc_{n+1} ] | \Fc_n ]
	 =	\Eb[ X_\infty | \Fc_n ]
	=		X_n ,
\end{align}
and, moreover, $(X_n)$ is uniformly integrable by Theorem \ref{thm:ui.cond.exp}.
\end{example}


\section{Exercises}
\label{sec:exercises-6}

\begin{exercise}
\label{exer:6.1}
Let $X_1, X_2, X_3, ...$ be a sequence of random variables such that 
$$ X_n \sim \text{Geo}(\lambda/n),  \quad for  \ n = 1,2,3,...$$ 
where $\lambda > 0$ is a constant. Define a new sequence as 
$$ Y_n  := \frac{1}{n} X_n, \quad for \ n=1,2,3,...$$
Show that $Y_n$ converges in distribution to $\Ec(\lambda)$.
\end{exercise}

\begin{exercise}
\label{exer:6.2}
Consider the sample space $\Om = [0,1]$ with uniform probability distribution, i.e., 
$$ P( \om \in [a,b]) = b - a, \quad \forall \ 0 \leq a \leq b \leq 1 $$
Define the sequence $\{X_n, n=1,2,...\}$ as $X_n(\om) = \frac{n}{n+1}\om + (1-\om)^n$. Also, define the random variable $X$ on this sample space as $X(\om) = \om$. Show that $X_n \stackrel{\text{a.s.}}{\rightarrow} X$.
\end{exercise}

\begin{exercise}
\label{exer:6.3}
Let $\{ X_n, n=1,2,...\}$ and $\{  Y_n, n=1,2,...  \}$ be two sequences of random variables, defined on some probability space $(\Om,\Fc,\Pb)$. Suppose that we know 
\begin{align}
X_n &\stackrel{\text{a.s.}}{\rightarrow} X , &
Y_n &\stackrel{\text{a.s.}}{\rightarrow} Y .
\end{align}
Prove that $X_n + Y_n \stackrel{\text{a.s.}}{\rightarrow} X+Y$.
\end{exercise}

\begin{exercise}
\label{exer:6.4}
Let $\{ X_n, n=1,2,...\}$ and $\{  Y_n, n=1,2,...  \}$ be two sequences of random variables, defined on some probability space $(\Om,\Fc,\Pb)$. Suppose that we know 
\begin{align}
X_n &\stackrel{\Pb}{\rightarrow} X , &
Y_n &\stackrel{\Pb}{\rightarrow} Y .
\end{align}
Prove that $X_n + Y_n \stackrel{\Pb}{\rightarrow} X+Y$.
\end{exercise}

\begin{exercise}
\label{exer:6.5}
Show that if $X_n$ is any sequence of random variables, there are constants $c_n \rightarrow \infty$ so that $X_n/c_n \stackrel{\text{a.s.}}{\rightarrow} 0$.  
\end{exercise}

\begin{exercise}
\label{exer:6.6}
Let $X_1, X_2, ...,$ be independent with $P(X_n = 1) = p_n$ and $P(X_n=0) = 1-p_n$. Show that\\
(a) $X_n \stackrel{\Pb}{\rightarrow} 0$ if and only if $p_n \rightarrow 0$.\\
(b) $X_n \stackrel{\text{a.s.}}{\rightarrow} 0$ if and only if $\sum_n p_n < \infty $.
\end{exercise}

\begin{exercise}
\label{exer:6.7}
Suppose that $X_1, X_2, ..., $ are independent with $P(X_n > x) = x^{-5}$ for all $x\geq 1$ and $n = 1,2, ...,$ Show that  limsup$_{n\rightarrow \infty} (\log X_n)/ \log n = c$ almost surely for some number c, and find c. 
\end{exercise}

%
%

\chapter{Brownian motion}
\label{chapter:bm}
The notes from this chapter are primarily taken from \cite[Chapter 3]{shreve2004}.  The goals of this chapter are (i) to define what we mean by ``Brownian motion'' and (ii) to develop important properties of Brownian motion.


\section{Scaled random walks}
\label{sec:random.walks}
In order to construct Brownian motion, we will begin with a random walk.  We return to the setting of Section \ref{sec:coin-toss}.  Our sample space $\Om$ and a generic element $\om$ in this sample space are given by
\begin{align}
\Om
	&=	\text{ the set of infinite sequences of $H$s and $T$s }, &
\om
	&=	\om_1 \om_2 \om_3 \ldots ,
\end{align}
where $\om_i$ is the result of the $i$th coin toss.  We take $\Fc_n$ to be the $\sig$-algebra generated by observing the first $n$ coin tosses and we set $\Fc = \sig( \cup_n \Fc_n)$.  The coin tosses are assumed to be independent and we take $\Pb(\om_i = H) = \Pb(\om_i = T) = (1/2)$.

\subsection{Symmetric random walk}
We construct a \emph{symmetric random walk} $M = (M_i)_{i \in \Nb_0}$ by
\begin{align}
M_0
	&=	0 , &
M_k
	&=	\sum_{j=1}^k X_j , &
X_i
	&=	\left\{ \begin{array}{ll}
			+1 & \om_i = H , \\
			-1 & \om_i = T .
			\end{array} \right.
\end{align}
Let $k_{i} < k_{i+1}$ for all $i \in \Nb_0$.  As the $(X_i)_{i \in \Nb}$ are independent, we clearly have
\begin{align}
\text{independent increments}:&&
(M_{k_2} - M_{k_1})
	&\ind (M_{k_4} - M_{k_3}) .
\end{align}
Moreover, as $\Eb X_j = 0$ and $\Vb X_j = 1$ for all $j$ we also have
\begin{align}
\Eb (M_{k_1} - M_{k_2})
	&=	\sum_{j=k_1+1}^{k_2} \Eb X_j
	=	0 , &
\Vb (M_{k_2} - M_{k_1})
	&=	\sum_{j=k_1+1}^{k_2} \Vb X_j
	=	k_2 - k_1 . \label{eq:dM.mean.var}
\end{align}
Thus, for the discrete time process $M = (M_i)_{i \geq 0}$ we see that variance accumulates at a rate of one per unit time.
\\[0.5em]
Next, we note that $M$ is a martingale with respect to the filtration $\Fb = (\Fc_n)_{n \in \Nb_0}$.  To see this, let $k \leq l$ and note that
\begin{align}
\Eb[ M_l | \Fc_k ]
	&= 	\Eb[ M_l - M_k + M_k | \Fc_k ]	
	= 	\Eb[ M_l - M_k | \Fc_k ] + \Eb[ M_k | \Fc_k ] \\
	&= 	\Eb ( M_l - M_k ) +  M_k .
	=		M_k ,
\end{align}
where we have used the independent increments property $(M_l - M_k) \ind \Fc_k$, equation \eqref{eq:dM.mean.var} and the fact that $M_k \in \Fc_k$.
\\[0.5em]
The \emph{quadratic variation} of 
the symmetric random walk $M$ up to time $k$, denoted $[M,M]_k$, is defined as
\begin{align}
[M,M]_k
	&:=	\sum_{j=1}^k (M_j - M_{j-1} )^2
	=		\sum_{j=1}^k X_j^2
	=		k ,
\end{align}
where we have used $X_j^2 = 1$.  The astute reader will notice that $[M,M]_k = \Vb M_k$.  However, it is important to note that \emph{the computation of variance and the computation of quadratic variation are different}!  To see this, note that if $\Pb(\om_i = H) = p \neq (1/2)$ then $\Vb X_i \neq 1$ and thus $\Vb M_k \neq k$.  However, since $X_j^2 = 1$ is unaffected by the value of $p = \Pb(\om_i=H)$, the computation of $[M,M]_k$ is also unaffected by $p$.  Another way to see that $\Vb M_k$ and $[M,M]_k$ are different is to note that $\Vb M_k$ is a \emph{statistical} quantity (i.e., it is an average over all $\om$) whereas $[M,M]_k$ is computed $\om$-by-$\om$ (it just turns out that, for each $\om$ we have $[M,M]_k(\om)=k$.  

\subsection{Scaled symmetric random walk}
We now fix a positive integer $n$.  We construct a \emph{scaled symmetric random walk} $W^{(n)} = (W_t^{(n)})_{t \geq 0}$ by
\begin{align}
W_t^{(n)}
	&=	\frac{1}{\sqrt{n}} M_{nt} , &
	&		\text{if $nt \in \Nb_0$}.
\end{align}
If $nt$ is not an integer, we define $W_t^{(n)}$ to be the linear interpolation of $M_{\lfloor nt \rfloor}$ and $M_{\lceil nt \rceil}$.  Note that $M$ was a discrete time process whereas $W^{(n)}$ is a continuous time process.
\\[0.5em]
Let $0 = t_0 < t_1 < t_2 < \ldots $ and suppose $n t_j \in \Nb$ for every $j$.  Then we have
\begin{align}
\text{independent increments}:&&
(W_{t_2}^{(n)} - W_{t_1}^{(n)})
	&\ind (W_{t_4}^{(n)} - W_{t_3}^{(n)}) ,
\end{align}
because non-overlapping increments depend on different coin tosses.  For example
\begin{align}
W_{0.2}^{(100)} - W_{0.0}^{(100)}&&
&\text{depends on coin tosses 1 through 20}, \\
W_{0.7}^{(100)} - W_{0.2}^{(100)}&&
&\text{depends on coin tosses 21 through 70} .
\end{align}
Next, suppose $0 \leq s \leq t$ are such that $ns \in \Nb$ and $nt \in \Nb$.  Then
\begin{align}
\Eb (W_t^{(n)} - W_s^{(n)} )
	&=	\frac{1}{\sqrt{n}} \Eb (M_{nt} - M_{ns})
	=		0 , &
\Vb (W_t^{(n)} - W_s^{(n)} )
	&=	\frac{1}{n} \Vb (M_{nt} - M_{ns})
	=		t-s. \label{eq:dW.mean.var}
\end{align}
Now, let us define the scaled filtration $\Fb^{(n)} = (\Fc_t^{(n)})_{t \geq 0}$ where
\begin{align}
\Fc_t^{(n)}
	&:= \Fc_{nt} , &
	& 	\text{if $nt \in \Nb_0$}.
\end{align}
The scaled random walk $W^{(n)}$, restricted to the set of $t$ for which $n t \in \Nb_0$, is a martingale with respect to the filtration $\Fb^{(n)}$.  
To see this, assume $0 \leq s \leq t$ are such that $ns \in \Nb$ and $nt \in \Nb$.  Then we have
\begin{align}
\Eb[ W_t^{(n)} | \Fc_s^{(n)} ]
	&= 	\Eb[ W_t^{(n)} - W_s^{(n)} + W_s^{(n)} | \Fc_s^{(n)} ]
	= 	\Eb[ W_t^{(n)} - W_s^{(n)} | \Fc_s^{(n)} ] + \Eb[ W_s^{(n)} | \Fc_s^{(n)} ] \\
	&= 	\Eb ( W_t^{(n)} - W_s^{(n)} ) +  W_s^{(n)} 
	=		W_s^{(n)} ,
\end{align}
where we have used $(W_t^{(n)} - W_s^{(n)}) \ind \Fc_s^{(n)}$, equation \eqref{eq:dW.mean.var} and $W_s^{(n)} \in \Fc_s^{(n)}$.
\\[0.5em]
The \emph{quadratic variation} of the scaled symmetric random walk $W^{(n)}$ up to time $t$, denoted $[W^{(n)},W^{(n)}]_t$, is defined as
\begin{align}
[W^{(n)},W^{(n)}]_t
	&:=	\sum_{j=1}^{nt} (W_{j/n}^{(n)} - W_{(j-1)/n}^{(n)} )^2
	=		\sum_{j=1}^{nt} \( \tfrac{1}{\sqrt{n}}X_j \)^2
	=		\sum_{j=1}^{nt} \tfrac{1}{n}
	=		t ,
\end{align}
where we again assume $nt \in \Nb_0$.
Thus, for the scaled symmetric random walk, we see that $\Vb W_t^{(n)} =  [W^{(n)},W^{(n)}]_t$.  However, we emphasize one more time that the \textit{computation of variance is a statistical average over all $\om$ }and the \textit{computation of quadratic variation is done $\om$-by-$\om$}.

\begin{theorem}
Fix $t \geq 0$.  Then $W_t^{(n)} \stackrel{\Dc}{\to} Z$ as $n \to \infty$, where $Z \sim \Nc(0,t)$.
\end{theorem}

\begin{proof}
We will use the Continuity Theorem \ref{thm:continuity}.  We have
\begin{align}
\phi_{W_t^{(n)}}(u)
	&=	\Eb \exp \( \ii u W_t^{(n)} \)
	=		\Eb \exp \Big( \frac{ \ii u}{\sqrt{n}} \sum_{j=1}^{nt} X_j \Big) \\
	&=	\prod_{j=1}^{nt} \Eb \exp \Big( \frac{\ii u}{\sqrt{n}}  X_j \Big) \\
	&=	\Big( \tfrac{1}{2} \ee^{\ii u/\sqrt{n}} + \tfrac{1}{2} \ee^{-\ii u/\sqrt{n}} \Big)^{nt}  \\
	&=	\Big( \frac{1}{2} + \frac{\ii u}{2\sqrt{n}} - \frac{u^2}{4n} + \Oc(\tfrac{1}{n^{3/2}})
				+ \frac{1}{2} - \frac{\ii u}{2\sqrt{n}} - \frac{u^2}{4n} + \Oc(\tfrac{1}{n^{3/2}}) \Big)^{nt} \\
	&=	\Big( 1 - \frac{u^2}{2n } + \Oc(\tfrac{1}{n^{3/2}}) \Big)^{nt}
	\to \ee^{-\tfrac{1}{2} t u^2 } , &
	&		\text{as $n \to \infty$}.
\end{align}
From Example \ref{ex:char.normal}, we know that if a random variable $Z$ is normally distributed $Z \sim \Nc(0,t)$ then its characteristic function $\phi_Z$ is given by
\begin{align}
\phi_{Z}(u)
	&=	\ee^{-\tfrac{1}{2} t u^2 } .
\end{align}
As $\lim_{n \to \infty} \phi_{W_t^{(n)}} \to \phi_{Z}$, we have from the Continuity Theorem \ref{thm:continuity} that $W_t^{(n)} \stackrel{\Dc}{\to} \Nc(0,t)$, as claimed.
\end{proof}


\section{Brownian motion}

\begin{definition}
\label{def:bm}
Let $(\Om, \Fc, \Pb)$ be a probability space.
A \emph{Brownian motion} is a stochastic process $W = (W_t)_{t \geq 0}$ that satisfies:
\begin{enumerate*}
\item $W_0 = 0$.
\item Independent increments: if $0 \leq t_0 < t_1 < \ldots < t_n < \infty$ then the increments $W_{t_1} - W_{t_0}, W_{t_2} - W_{t_1}, \ldots, W_{t_n} - W_{t_{n-1}}$ are independent.
\item If $0 \leq r < s$ then $W_s - W_r \sim \Nc(0,s-r)$.
\item The map $t \to W_t$ is continuous for every $\om$.
\end{enumerate*}
\end{definition}

\noindent
It is clear from the previous sections that we can construct a Brownian motion as a limit of a scaled symmetric random walk.  Had we simply given Definition \ref{def:bm} at the beginning of this chapter with no further introduction, one might have legitimately asked if there exists a process that satisfies the properties of a Brownian motion.  There are other methods to prove the existence of Brownian motion.  But, the scaled random walk construction is perhaps the most intuitive.
\\[0.5em]
What is $\Om$ in Definition \ref{def:bm}?  It could be an infinite series of $H$s and $T$s, representing movements up and down of a scaled random walk.  Or, it could be $\Om = C_0(\Rb_+)$, the set of continuous functions on $\Rb_+$, starting from zero.  In this case, an element of $\Om$ is a continuous function $t \to \om(t)$ and one could simply take $W_t(\om) = \om(t)$.  Whatever the sample space, the probability of any single element $\om$ is zero: $\Pb(\om)=0$, but probabilities such as $\Pb(W_t \leq 0)$ are well-defined.  
\\[0.5em]
It is useful to note that
\begin{align}
\text{Co}\Vb[W_T,W_t]
	&=	\Eb W_T W_t
	=		\Eb (W_T - W_t + W_t) W_t \\
	&=	\Eb (W_T - W_t) W_t + \Eb W_t^2 
	=	\Eb W_t \Eb[ W_T - W_t | W_t ] + t
	=		t . \label{eq:BM-Cov}
\end{align}  
Thus, $\text{Co}\Vb[W_s,W_t] = s \wedge t$.

\begin{definition}
\label{def:bm.filtration}
Let $(\Om, \Fc, \Pb)$ be on probability space on which a Brownian motion $W=(W_t)_{t \geq 0}$ is defined.  A \emph{filtration for the Brownian motion $W$} is a collection of $\sig$-algebras $\Fb=(\Fc_t)_{t \geq 0}$ satisfying:
\begin{enumerate*}
\item Information accumulates: if $0 \leq s < t$ then $\Fc_s \subset \Fc_t$.
\item Adaptivity: for all $t \geq 0$, we have $W_t \in \Fc_t$.
\item Independence of future increments: if $u>t \geq 0$ then $(W_u - W_t) \ind \Fc_t$.
\end{enumerate*}
\end{definition}

\noindent
The most natural choice for this filtration $\Fb$ is the \emph{natural filtration} for $W$.  That is $\Fc_t = \sig(W_u,0\leq u \leq t)$.  In principle the filtration $(\Fc_t)_{t \geq 0}$ could contain more than the information obtained by observing $W$.  However, the information in the filtration is not allowed to destroy the independence of future increments of Brownian motion.
\\[0.5em]
Not surprisingly, if $\Fb = (\Fc_t)_{t\geq 0}$ is a filtration for a Brownian motion $W$ then $W$ is a martingale with respect to this filtration.  We see this, let $0 \leq s < t$ and observe that
\begin{align}
\Eb[ W_t | \Fc_s ]
	&=	\Eb[ W_t - W_s + W_s | \Fc_s ]
	=		\Eb[ W_t - W_s | \Fc_s ] + \Eb[ W_s | \Fc_s ]
	=		\Eb[ W_t - W_s ] + W_s
	=		W_s .
\end{align}


\section{Quadratic variation}
\label{sec:QV}
In Section \ref{sec:random.walks}, we computed the quadratic variation of a symmetric random walk $M$ and a scaled symmetric random walk $W^{(n)}$.  However, we did not rigorously define the notion of quadratic variation.  In this Section we will define what we mean by ``quadratic variation'' and we will compute this quantity for a Brownian motion $W$.  Before defining quadratic variation we first introduce the ``first variation.''
\\[0.5em]
Let $\Pi$ be a partition of $[0,T]$ and define $\norm{\Pi}$ as
\begin{align}
\Pi
	&=	\{ t_0 , t_1 , \ldots , t_n \} , &
0
	&= t_0 < t_1 < \ldots < t_n = T , &
\norm{\Pi}
	&=	\max_i ( t_{i+1} - t_i ) .  \label{eq:Pi.norm}
\end{align}
Let $f :[0,T] \to \Rb$.  We define the \emph{first variation of $f$ up to time $T$}, denoted $\text{FV}_T(f)$, by
\begin{align}
\text{FV}_T(f)
	&:=	\lim_{\norm{\Pi} \to 0} \sum_{j=0}^{n-1} | f(t_{j+1}) - f(t_j) | .
\end{align}
Suppose that $f \in C([0,T])$ and $f'(t)$ exists and is finite for all $t \in (0,T)$.  Then, by the \emph{Mean Value Theorem}, there exits $t_j^* \in [t_j,t_{j+1}]$ such that
\begin{align}
f(t_{j+1}) - f(t_j)
	&=	f'(t_j^*) ( t_{j+1} - t_j ) . \label{eq:t.star}
\end{align}
Thus, if $f \in C([0,T])$ and $f'(t)$ exists and is finite for all $t \in (0,T)$, we have
\begin{align}
\text{FV}_T(f)
	&:=	\lim_{\norm{\Pi} \to 0} \sum_{j=0}^{n-1} | f'(t_j^*) | ( t_{j+1} - t_j )
	=		\int_0^T \dd t \, |f'(t)| .
\end{align}

\begin{definition}
\label{def:QV}
Let $f :[0,T] \to \Rb$. We define the \emph{quadratic variation of $f$ up to time $T$}, denoted $[f,f]_T$ as
\begin{align}
[f,f]_T
	&:=	\lim_{\norm{\Pi} \to 0} \sum_{j=0}^{n-1} [ f(t_{j+1}) - f(t_j) ]^2 ,
\end{align}
where $\Pi$ and $\norm{\Pi}$ are as defined in \eqref{eq:Pi.norm}.  
\end{definition}

\begin{proposition}
\label{prp:zero.QV}
Suppose $f:[0,T] \to \Rb$ has a continuous first derivative on the closed interval $[0,T]$: $f \in C^1([0,T])$.  Then $[f,f]_T = 0$.
\end{proposition}

\begin{proof}
For any partition $\Pi$ of $[0,T]$ we have
\begin{align}
\sum_{j=0}^{n-1} \[ f(t_{j+1}) - f(t_j) \]^2
	&=	\sum_{j=0}^{n-1} [ f'(t_j^*) ]^2 ( t_{j+1} - t_j )^2
	\leq \norm{\Pi} \sum_{j=0}^{n-1} [ f'(t_j^*) ]^2 ( t_{j+1} - t_j ) .
\end{align}
Inserting the above inequality into the definition of $[f,f]_T$, we obtain
\begin{align}
[f,f]_T
	&\leq \lim_{\norm{\Pi} \to 0} \norm{\Pi} \cdot \lim_{\norm{\Pi}\to 0} \sum_{j=0}^{n-1} [ f'(t_j^*) ]^2 ( t_{j+1} - t_j ) 
	=	0 \cdot \int_0^T \dd t \, |f'(t)|^2 = 0 ,
\end{align}
where we have used $\int_0^T \dd t \, |f'(t)|^2 < \infty$ as, by assumption, $f \in C^1((0,T))$.
\end{proof}

\noindent
In ordinary calculus, we typically deal with functions $f \in C^1$, and hence $[f,f]_T = 0$.  For this reason, quadratic variation never arises in usual calculus.  However, it turns out that for almost every $\om \in \Om$, we have that $t \mapsto W_t(\om)$ is \emph{not} differentiable. We can see this from the scaled random walk construction of Brownian motion.  The slope of the scaled random walk $W_t^{(n)}$ at any $t$ for which $\frac{\dd}{\dd t} W_t^{(n)}(\om)$ is defined is
\begin{align}
\lim_{\eps \to 0} \frac{ W_{t+\eps}^{(n)} - W_{t}^{(n)} }{\eps}
	&=	\pm \sqrt{n} \to \pm \infty &
	&\text{as $n \to \infty$}.
\end{align}
Thus, Brownian motion, which we constructed as a limit of a scaled random walk $W(\om) := \lim_{n \to \infty} W^{(n)}(\om)$ is not differentiable at any $t$, $\Pb$-a.s..
When a function $f \notin C^1((0,T))$, then Proposition \ref{prp:zero.QV} can fail.  
Indeed, as we will show, paths of BM have strictly positive quadratic variation.  
It is for this reason that stochastic calculus is different from ordinary calculus.

\begin{theorem}
\label{thm:bm.QV}
Let $W$ be a Brownian motion.  Then, for all $T \geq 0$ we have $[W,W]_T = T$, where the limit in Definition \ref{def:QV} is taken in $L^2(\Om,\Fc,\Pb)$.  Moreover, for any sequence of partitions $(\Pi_k)_{k \geq 1}$ of $[0,T]$ satisfying $\norm{\Pi_k} \to 0$, there exists a subsequence $(\Pi_{k_j})_{j \geq 1}$ along which the convergence holds $\Pb$-almost surely.
\end{theorem}

\begin{proof}
For a fixed partition $\Pi = \{ 0 = t_0, t_1, \ldots, t_n = T \}$ we define the \emph{sampled quadratic variation of $W$}, denoted $Q_\Pi$, by
\begin{align}
Q_\Pi
	&:=	\sum_{j=0}^{n-1} (W_{t_{j+1}} - W_{t_j})^2 .
\end{align}
\underline{Step 1: $L^2$ convergence.}
We will show that $\Eb Q_\Pi = T$ and $\Vb Q_\Pi \to 0$ as $\norm{\Pi} \to 0$, which implies that $Q_\Pi \to T$ in $L^2(\Om,\Fc,\Pb)$.
To begin, we use the fact that $W_{t_{j+1}} - W_{t_j} \sim \Nc(0,t_{j+1}-t_j)$ to compute
\begin{align}
\Eb (W_{t_{j+1}} - W_{t_j})^2
	&=	t_{j+1} - t_j , \label{eq:W2.mean} \\
\Vb (W_{t_{j+1}} - W_{t_j})^2
	&=	\Eb (W_{t_{j+1}} - W_{t_j})^4 - \( \Eb (W_{t_{j+1}} - W_{t_j})^2 \)^2 \\
	&=	3 ( t_{j+1} - t_j )^2 - (t_{j+1} - t_j )^2 
	=		2 ( t_{j+1} - t_j )^2 . \label{eq:W2.var}
\end{align}
Next, using \eqref{eq:W2.mean}, \eqref{eq:W2.var} and the independence of the increments of $W$, we obtain
\begin{align}
\Eb Q_\Pi
	&=	\sum_{j=0}^{n-1} \Eb (W_{t_{j+1}} - W_{t_j})^2 
	=		\sum_{j=0}^{n-1}  ( t_{j+1} - t_j ) 
	=		T , \\
\Vb Q_\Pi
	&=	\sum_{j=0}^{n-1} \Vb (W_{t_{j+1}} - W_{t_j})^2 
	=		\sum_{j=0}^{n-1}2 ( t_{j+1} - t_j )^2
	\leq \norm{\Pi} \sum_{j=0}^{n-1}2 ( t_{j+1} - t_j )
	=		2 \norm{\Pi} T .
\end{align}
As $\Eb Q_\Pi = T$, we have
\begin{align}
\Eb ( Q_\Pi - T )^2
	&=	\Vb Q_\Pi 
	\leq	2 \norm{\Pi} T 
	\to		0 , &
	&\text{as $\norm{\Pi} \to 0$} , \label{eq:QV.L2.bound}
\end{align}
which proves that $[W,W]_T = \lim_{\norm{\Pi}\to 0} Q_{\Pi} = T$ in $L^2(\Om,\Fc,\Pb)$.
\\[0.5em]
\underline{Step 2: $\Pb$-a.s. convergence along a subsequence.}
Let $(\Pi_k)_{k \geq 1}$ be a sequence of partitions of $[0,T]$ with $\norm{\Pi_k} \to 0$.  Choose a subsequence $(k_j)_{j \geq 1}$ such that $\norm{\Pi_{k_j}} \leq 2^{-j}$ for all $j$.  Fix $\eps > 0$.  Using Chebyshev's inequality (i.e., Markov's inequality \eqref{ineq:markov} applied to $(Q_\Pi - T)^2$) and \eqref{eq:QV.L2.bound}, we have
\begin{align}
\Pb ( | Q_{\Pi_{k_j}} - T | > \eps )
	&\leq	\frac{1}{\eps^2} \Eb ( Q_{\Pi_{k_j}} - T )^2
	\leq	\frac{2 T \norm{\Pi_{k_j}}}{\eps^2}
	\leq	\frac{2 T}{\eps^2} 2^{-j} ,
\end{align}
which is summable in $j$.  Thus, by the Borel-Cantelli Lemma \ref{thm:borel-cantelli}, we have $\Pb ( | Q_{\Pi_{k_j}} - T | > \eps \text{ i.o.} ) = 0$.  Taking $\eps = 1/m$ and using the fact that a countable union of null sets is a null set, we obtain
\begin{align}
\Pb \Big( \bigcup_{m=1}^\infty \{ | Q_{\Pi_{k_j}} - T | > 1/m \text{ i.o.} \} \Big)
	&=	0 ,
\end{align}
which is equivalent to $Q_{\Pi_{k_j}} \to T$ $\Pb$-a.s. as $j \to \infty$.  (The same argument shows that no subsequence is needed if $\sum_k \norm{\Pi_k} < \infty$.)
\end{proof}

\noindent
Suppose $\dd t \ll 1$ and define $\dd W_t := W_{t + \dd t}- W_t$.  The above computations show that $\Eb( (\dd W_t)^2 ) = \dd t$ and $\Vb( (\dd W_t)^2 ) = 2 \dd t^2$.  Since $\dd t^2$ is practically zero for $\dd t \ll 1$, one can imagine that $(\dd W_t)^2$ is \emph{almost} equal to a constant $\dd t$.  \emph{Informally}, we write this as
\begin{align}
\dd W_t \dd W_t
	&=	\dd t . \label{eq:dW2.dt}
\end{align}
This informal statement, while not rigorously correct, captures the spirit of the quadratic variation computation for $W$.

\begin{definition}
\label{def:CoV}
Let $f,g :[0,T] \to \Rb$. We define the \emph{covariation of $f$ and $g$ up to time $T$}, denoted $[f,g]_T$ as
\begin{align}
[f,g]_T
	&:=	\lim_{\norm{\Pi} \to 0} \sum_{j=0}^{n-1} \[ f(t_{j+1}) - f(t_j) \] \[ g(t_{j+1}) - g(t_j) \] ,
\end{align}
where $\Pi$ and $\norm{\Pi}$ are as defined in \eqref{eq:Pi.norm}.  
\end{definition}

\begin{theorem}
\label{thm:CoV}
Let $W$ be a Brownian motion and let \emph{Id} be the identity function: $\text{\emph{Id}}(t)=t$.  Then, for all $T \geq 0$ we have $[W,\text{\emph{Id}}]_T = 0$ almost surely and $[\text{\emph{Id}},\text{\emph{Id}}]_T = 0$.
\end{theorem}

\begin{proof}
For a fixed partition $\Pi = \{ 0 = t_0, t_1, \ldots, t_n = T \}$ we have
\begin{align}
\big| [W,\Id]_T \big|
	&=	\Big| \lim_{\norm{\Pi} \to 0 } \sum_{j=0}^{n-1}  (W_{t_{j+1}} - W_{t_j}) \cdot ( t_{j+1} - t_j ) \Big| \\
	&=	\lim_{\norm{\Pi} \to 0 } \Big| \sum_{j=0}^{n-1}  (W_{t_{j+1}} - W_{t_j}) \cdot ( t_{j+1} - t_j ) \Big| \\
	&\leq \lim_{\norm{\Pi} \to 0 } \sum_{j=0}^{n-1} \Big|  (W_{t_{j+1}} - W_{t_j}) \cdot ( t_{j+1} - t_j ) \Big| \\
	&\leq \lim_{\norm{\Pi} \to 0 } \max_{0 \leq j \leq n-1} |(W_{t_{j+1}} - W_{t_j})| \sum_{j=0}^{n-1}  ( t_{j+1} - t_j ) \\
	&=	 T  \cdot \lim_{\norm{\Pi} \to 0 } \max_{0 \leq j \leq n-1} | W_{t_{j+1}} - W_{t_j} |
	= T \cdot 0 = 0.
\end{align}
where, in the last step we are using the fact that $W$ is continuous.  To see that $[\Id,\Id]_T = 0$, simply note that $g \in C^1([0,T])$ and use Proposition \ref{prp:zero.QV}.
\end{proof}

\noindent
Just as \eqref{eq:dW2.dt} captures the spirit of the computation of $[W,W]_T$, the following equations
\begin{align}
\dd W_t \dd t
	&=	0 , &
\dd t \dd t
	&=	0 ,
\end{align}
\emph{informally} capture the spirit of the $[W,\Id]_T$ and $[\Id,\Id]_T$ computations.

\section{Markov property of Brownian motion}

\begin{theorem}
Let $W = (W_t)_{t \geq 0}$ be a Brownian motion and let $(\Fc_t)_{t \geq 0}$ be a filtration for this Brownian motion. Then $W$ is a Markov process.
\end{theorem}

\begin{proof}
According to Definition \ref{def:markov}, we must show that there exists a function $g$ such that
\begin{align}
\Eb[ f(W_T) | \Fc_t ]
	&=	g(W_t) ,
\end{align}
where $T \geq t \geq 0$.  Noting that $W_t \in \Fc_t$ and $W_T-W_t \ind \Fc_t$, we have
\begin{align}
\Eb[ f(W_T) | \Fc_t ]
	&=	\Eb[ f(W_T-W_t + W_t) | \Fc_t ]
	=		\int_\Rb \dd y \, f(y + W_t) \Gam(t,0;T,y) =: g(W_t) ,
\end{align}
where $\Gam(t,x;T,\cdot)$ is the density of a normal random variable with mean $x$ and variance $T-t$.
\end{proof}


\section{First hitting time of Brownian motion}
\label{sec:hitting.bm}

\begin{theorem}
Let $W=(W_t)_{t \geq 0}$ be a Brownian motion with a filtration $(\Fc_t)_{t \geq 0}$.  Define a process $Z=(Z_t)_{t \geq 0}$ by
\begin{align}
Z_t
	&=	\ee^{ -\tfrac{1}{2} \sig^2 t + \sig W_t } . \label{eq:exp.martingale}
\end{align}
Then $Z$ is a martingale with respect to $(\Fc_t)_{t \geq 0}$.
\end{theorem}

\begin{proof}
The proof is a simple computation.  Fixing $T \geq t$ we have
\begin{align}
\Eb[ Z_T | \Fc_t ]
	&=	\Eb[ \ee^{ -\tfrac{1}{2} \sig^2 T + \sig W_T }| \Fc_t ] 
	=		\ee^{ -\tfrac{1}{2} \sig^2 T + \sig W_t }\Eb[ \ee^{ \sig ( W_T - W_t ) }| \Fc_t ] \\
	&=	\ee^{ -\tfrac{1}{2} \sig^2 T + \sig W_t } \ee^{ \tfrac{1}{2} \sig^2 ( T - t ) } 
	=		\ee^{ -\tfrac{1}{2} \sig^2 t + \sig W_t }
	=		Z_t .
\end{align}
Thus, $Z$ is a martingale.
\end{proof}

\noindent
The process $Z$ is sometimes referred to as an \emph{exponential martingale}.  The exponential martingale will be used to compute the distribution of 
\begin{align}
\tau_m
	&:=	\inf\{ t \geq 0 : W_t = m \} . \label{eq:tau.m}
\end{align}
We call $\tau_m$ the \emph{first hitting time} or \emph{first passage time} of a Brownian motion $W$ to level $m$.  

\begin{theorem}
For any $m \in \Rb$ and $\alpha > 0$ we have
\begin{align}
\Eb \ee^{- \alpha \tau_m}
	&=	\ee^{- |m|\sqrt{2 \alpha} } . \label{eq:laplace.tau.m}
\end{align}
where $\tau_m$ is defined in \eqref{eq:tau.m}.
\end{theorem}

\begin{proof}
First, assume $m \geq 0$.
Let $Z^{(m)} = (Z_t^{(m)})_{t \geq 0}$ be given by
\begin{align}
Z_t^{(m)}
	&:=	Z_{t \wedge \tau_m} ,
\end{align}
where $Z$ is given by \eqref{eq:exp.martingale}.
We call $Z^{(m)}$ a \emph{stopped} process as it remains at $Z_{\tau_m}$ forever after $W$ hits $m$.
As $Z$ is a martingale it follows that the stopped process $Z^{(m)}$ is also a martingale.
Thus, we have
\begin{align}
1
	&=	Z_0^{(m)}
	=		\Eb Z_t^{(m)}
	=		\Eb \ee^{ -\tfrac{1}{2} \sig^2 t \wedge \tau_m + \sig W_{t \wedge \tau_m} } .
\end{align}
Now, it turns out that $\Pb( \tau_m < \infty) = 1$ (a fact that can be proved with relative ease).  
As a result, we have $\lim_{t \to \infty} t \wedge \tau_m = \tau_m$ and $W_{\tau_m}=m$.  Thus, we obtain
\begin{align}
1
	&=	\lim_{t \to \infty} \Eb \ee^{ -\tfrac{1}{2} \sig^2 t \wedge \tau_m + \sig W_{t \wedge \tau_m} }
	=		\Eb \lim_{t \to \infty} \ee^{ -\tfrac{1}{2} \sig^2 t \wedge \tau_m + \sig W_{t \wedge \tau_m} }
	=		\Eb \ee^{ -\tfrac{1}{2} \sig^2 \tau_m + \sig m } ,
\end{align}
where passing the limit through the expectation is allowed because $Z_t^{(m)} \leq \ee^{\sig m}$ (dominated convergence).
Setting $\sig = \sqrt{2 \alpha}$ we obtain $\Eb \ee^{- \alpha \tau_m} =	\ee^{- m\sqrt{2 \alpha} }$, which agrees with \eqref{eq:laplace.tau.m} for $m \geq 0$.  To obtain \eqref{eq:laplace.tau.m} for $m<0$, simply note that, since Brownian motion is symmetric about zero, the distribution of $\tau_m$ is the same as the distribution of $\tau_{-m}$.
\end{proof}

\noindent
We can compute $\Eb \tau_m$ as follows
\begin{align}
\Eb \tau_m
	&=	\lim_{\alpha \searrow 0} -\frac{\dd}{\dd \alpha} \Eb \ee^{- \alpha \tau_m}
	=		\lim_{\alpha \searrow 0} -\frac{\dd}{\dd \alpha} \ee^{- |m|\sqrt{2 \alpha} } 
	=		\lim_{\alpha \searrow 0} \frac{|m|}{\sqrt{2 \alpha}} \ee^{- |m|\sqrt{2 \alpha} }
	=		\infty , &
m
	&\neq 0 .
\end{align}
Thus, while $\Pb( \tau_m < \infty) = 1$ we have $\Eb \tau_m = \infty$.


\section{Reflection principle}
\label{sec:reflection}

\begin{figure}
\centering
\begin{tikzpicture}[x=0.75cm,y=0.75cm,line cap=round,line join=round]
	\draw[->] (0,0) -- (14.8,0) node[right] {$s$};
	\draw[->] (0,0) -- (0,6.8) node[above] {$W_s$};
	\foreach \y/\lab in {2/w, 4/m, 6/2m-w} {
		\draw[dotted, thick] (0,\y) -- (14.4,\y);
		\node[left] at (0,\y) {$\lab$};
	}
	\draw (6,0.1) -- (6,-0.1) node[below] {$\tau_m$};
	\draw (14,0.1) -- (14,-0.1) node[below] {$t$};
	\draw[thick] (0,0) -- (2,2) -- (3,1) -- (6,4) -- (7,3) -- (10,6) -- (14,2);
	\draw[thick, dashed] (6,4) -- (7,5) -- (10,2) -- (14,6);
\end{tikzpicture}
\caption{An illustration of the reflection principle for Brownian motion.}
\label{fig:ReflectionPinciple}
\end{figure}
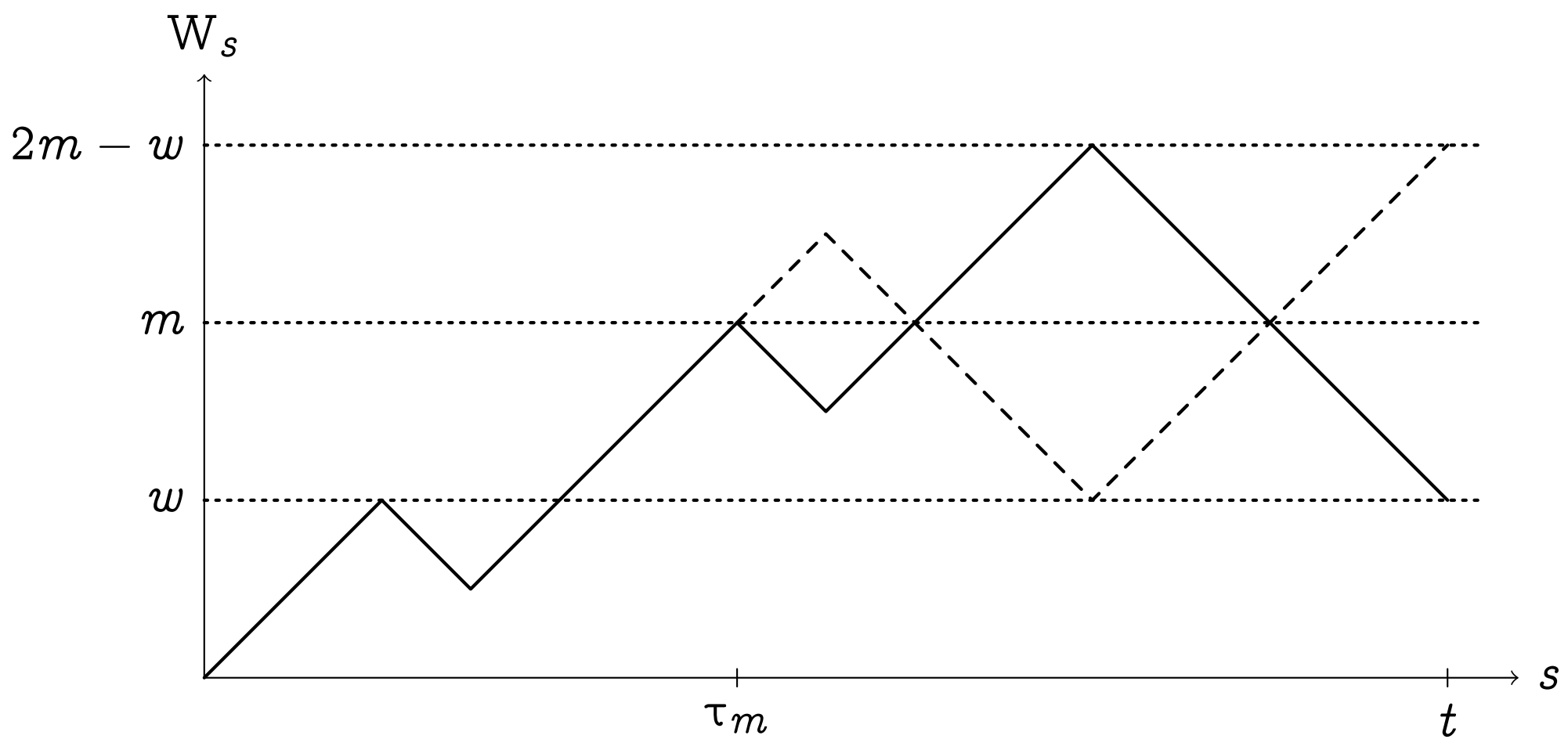

For every sample path of Brownian motion that hits a level $m>0$ prior to time $t$ and finishes at level $W_t = w \leq m$ there is an equally likely path that finishes at a level $W_t = 2 m - w$.  In Figure \ref{fig:ReflectionPinciple} these two paths are represented by the solid and dashed lines, respectively.  The dotted lines from top to bottom represent the levels $2 m - w$, $m$ and $w$, and $\tau_m$ is the first time the path hits level $m$.
It follows that
\begin{align}
\Pb( \tau_m \leq t , W_t \leq w )
	&=	\Pb( W_t \geq 2m-w) , &
m
	&>	0 , &
w
	&\leq m . \label{eq:ReflectionPrinciple}
\end{align}
We will use this fact to find the density of $\tau_m$.

\begin{theorem}
For all $m \neq 0$, the first hitting time $\tau_m$ of Brownian motion to level $m$ has a density $f_{\tau_m}$, which is given by
\begin{align}
f_{\tau_m}(t)
	&=	\Ib_{\{ t \geq 0 \}}\frac{|m|}{t\sqrt{2\pi t}}\ee^{-m^2/(2t)} . \label{eq:tau.m.density}
\end{align}
\end{theorem}

\begin{proof}
Consider the case $m>0$.  Substitute $w=m$ into \eqref{eq:ReflectionPrinciple} to obtain
\begin{align}
\Pb( \tau_m \leq t, W_t \leq m )
	&= \Pb( W_t \geq m ) .
\end{align}
Now, note that $W_t \geq m$ implies that $\tau_m \leq t$.  And thus
\begin{align}
\Pb( \tau_m \leq t, W_t \geq m )
	&= \Pb( W_t \geq m ) .
\end{align}
Adding the above two equations, we obtain
\begin{align}
2 \Pb( W_t \geq m )
	&=	\Pb( \tau_m \leq t, W_t \leq m ) + \Pb( \tau_m \leq t, W_t \geq m )
	=		\Pb( \tau_m \leq t ) .
\end{align}
The density \eqref{eq:tau.m.density} is obtained using
\begin{align}
f_{\tau_m}(t)
	&=	\frac{\dd}{\dd t} \Pb( \tau_m \leq t )
	=		\frac{\dd}{\dd t} 2 \Pb( W_t \geq m )
	=		\frac{\dd}{\dd t} 2 \int_m^\infty \dd x \, \frac{1}{\sqrt{2 \pi t}} \exp \Big( - \frac{x^2}{2t} \big) \\
	&=	\frac{\dd}{\dd t} 2 \int_{m/\sqrt{t}}^\infty \dd y \, \frac{1}{\sqrt{2 \pi }} \exp \Big( - \frac{y^2}{2} \big) 
	=		\frac{m}{t\sqrt{2\pi t}}\ee^{-m^2/(2t)}
\end{align}
The case $m \leq 0$ is proved in a similar manner.
\end{proof}

\noindent
Now, let us define the \emph{running maximum}, denoted $\Wb=(\Wb_t)_{t \geq 0}$, of Brownian motion
\begin{align}
\Wb_t
	&=	\max_{0 \leq s \leq t} W_s .
\end{align}
Note that $\Wb_t \geq m$ if an only if $\tau_m \leq t$.  Thus, we have from \eqref{eq:ReflectionPrinciple} that
\begin{align}
\Pb( \Wb_t \geq m , W_t \leq w)
	&=	\Pb(W_t \geq 2m-w) , &
m
	&> 0 , &
w
	&\leq m . \label{eq:ReflectionPrinciple2}
\end{align}
We can use \eqref{eq:ReflectionPrinciple2} to obtain the joint density of $(W_t,\Wb_t)$.

\begin{theorem}
For any $t>0$, the joint density of Brownian motion $W_t$ and its running maximum $\Wb_t$ is 
\begin{align}
f_{W_t,\Wb_t}(w,m)
	&=	\frac{2(2m-w)}{t \sqrt{2\pi t}} \ee^{-(2m-w)^2/(2t)} , &
m
	&>0 , &
w
	&\leq m . \label{eq:joint.density}
\end{align}
\end{theorem}

\begin{proof}
Note that
\begin{align}
\Pb(W_t \leq w , \Wb_t \leq m)
	&=	\Pb( W_t \leq w) - \Pb(W_t \leq w , \Wb_t \geq m) \\
	&=	\Pb( W_t \leq w) - \Pb(W_t \geq 2m-w)
	=		\Pb( W_t \leq w) - \int_{2m-w}^\infty \dd x \, f_{W_t}(x) ,
\end{align}
where $f_{W_t}$ is the density of $W_t$, which is a $\Nc(0,t)$ random variable.  
To obtain the density \eqref{eq:joint.density} use
\begin{align}
f_{W_t,\Wb_t}(w,m)
	&=	\frac{\d^2}{\d m \d w} \Pb(W_t \leq w , \Wb_t \leq m)
	=		- \frac{\d^2}{\d m \d w} \int_{2m-w}^\infty \dd x \, f_{W_t}(x) .
\end{align}
The rest of the computation is algebra.
\end{proof}

\begin{corollary}
For any $t>0$, the conditional density of $\Wb_t$ given $W_t$ is given by
\begin{align}
f_{\Wb_t|W_t}(m,w)
	&=	\frac{2(2m-w)}{t}\ee^{-2m(m-w)/t} , &
m
	&> 0 , &
w
	&\leq m . \label{eq:conditional.max}
\end{align}
\end{corollary}

\begin{proof}
Expression \eqref{eq:conditional.max} follows directly from
\begin{align}
f_{\Wb_t|W_t}(m,w)
	&=	\frac{f_{W_t,\Wb_t}(w,m)}{f_{W_t}(w)} .
\end{align}
where $f_{W_t}$ is the density of $W_t$, which is $\Nc(0,t)$.
\end{proof}


\section{Exercises}
\label{sec:exercises-7}

\begin{exercise}
\label{exer:7.1}
Let $W$ be a Brownian motion and let $\Fb$ be a filtration for $W$.  Show that $W_t^2 - t$ is a martingale with respect to the filtration $\Fb$.
\end{exercise}

\begin{exercise}
\label{exer:7.2}
Compute the characteristic function of $W_{N_t}$ where $N$ is a Poisson process with intensity $\lambda$ and the Brownian motion $W$ is independent of the Poisson process $N$.
\end{exercise}

\begin{exercise}
\label{exer:7.3}
The $m$th variation of a function $f$, over the interval $[0,T]$ is defined as
\begin{align}
V_T(m,f)
	&:=\lim_{\left\| \Pi \right\| \to 0} \sum_{j=0}^{n-1} | f(t_{j+1}) - f(t_j) |^m , &
\Pi
	&= \{ 0 = t_0, t_1, \ldots t_n = T \} , &
\left\| \Pi \right\|
	&= \max_j ( t_{j+1} - t_j ) .
\end{align}
Show that $V_T(1,W) = \infty$ and $V_T(3,W) = 0$, where $W$ is a Brownian motion.
\end{exercise}

\begin{exercise}
\label{exer:7.4}
Define
\begin{align}
X_t
	&=	\mu t + W_t , &
\tau_m
	&:=	\inf\{ t \geq 0 : X_t = m \} ,
\end{align}
where $W=(W_t)_{t \geq 0}$ is a Brownian motion.
Let $\Fb = (\Fc_t)_{t \geq 0}$ be a filtration for $W$.
Show that $Z$ is a martingale with respect to $\Fb$ where
\begin{align}
Z_t
	&=	\exp \( \sig X_t- (\sig \mu + \sig^2/2) t \) .
\end{align}
Assume $\mu > 0$ and $m \geq 0$.  Assume further that $\tau_m < \infty$ with probability one and the stopped process $Z_{t \wedge \tau_m}$ is a martingale.  Find the Laplace transform $\Eb \ee^{-\alpha \tau_m}$.
\end{exercise}

\begin{exercise}
\label{exer:7.5}
Prove that
\begin{align}
[f,W]_T
	&=	0 ,
\end{align}
where $f \in C^1([0,T])$ and $W$ is a Brownian motion.
\end{exercise}

\begin{exercise}
\label{exer:7.6}
Fix a probability space $(\Om, \Fc, \Pb)$ and let $\Fb = (\Fc_t)_{0 \leq t \leq T}$ be a filtration for a Brownian motion $W = (W_t)_{0 \leq t \leq T}$.  Define a process $Z = (Z_t)_{0 \leq t \leq T}$ by
\begin{align}
Z_t
	&:=	\ee^{u W_t - u^2 t / 2} , &
u
	&\in \Rb .
\end{align}
(a) Prove that the process $Z$ is a martingale under $\Pb$ with respect to the filtration $\Fb$.
\\[0.5em]
(b) Using the fact that $Z$ is a martingale, show that the following processes $M = (M_t)_{0 \leq t \leq T}$ and $N = (N_t)_{0 \leq t \leq T}$ are martingales, where
\begin{align}
M_t
	&:= W_t^2 - t , &
N_t
	&=	W_t^3 - 3 t W_t .
\end{align}
Note: if $X \sim \Nc(\mu,\sig^2)$ then $\Eb \ee^{u X} = \ee^{\mu u + u^2 \sig^2 / 2}$.
\end{exercise}

\begin{exercise}
\label{exer:7.7}
Fix a probability space $(\Om, \Fc, \Pb)$ and a filtration $\Fb = (\Fc_t)_{0 \leq t \leq T}$.  We say that a $\Cb$-valued process is a \textit{$(\Pb,\Fb)$-martingale} if, for any $0 \leq s \leq t \leq T$ the following holds $\Eb( Z_t | \Fc_s ) = Z_s$.  Note that if we decompose $Z$ into its Real and Imaginary parts $Z = X + \ii Y$, then we must have $\Eb( X_t | \Fc_s ) = X_s$ and $\Eb( Y_t | \Fc_s ) = Y_s$. Use this information to prove that the processes $M=(M_t)_{0 \leq t \leq T}$ and $N=(N_t)_{0 \leq t \leq T}$, defined as follows
\begin{align}
M_t
	&:= \ee^{ u^2 t /2} \cos (u W_t) , &
N_t
	&:=	\ee^{ u^2 t /2} \sin (u W_t) , &
u
	&\in \Rb ,
\end{align}
are $(\Pb,\Fb)$-martingales, where $W = (W_t)_{0 \leq t \leq T}$ is a Brownian motion and $\Fb$ is a filtration for $W$.
\\[0.5em]
Note: if $V \sim \Nc(\mu,\sig^2)$ then $\Eb \ee^{\ii u V} = \ee^{\ii \mu u - u^2 \sig^2 / 2}$.
\end{exercise}

\begin{exercise}
\label{exer:7.8}
Fix a probability space $(\Om, \Fc, \Pb)$ and let $W = (W_t)_{t \geq 0}$ be a Brownian motion.  Define the running maximum $\Wb = (\Wb_t)_{t \geq 0}$ of $W$ by
\begin{align}
\Wb_t
	&:=	\max_{0 \leq s \leq t} W_s .
\end{align}
Recall that, as a consequence of the reflection principle, the joint density of $(W_t,\Wb_t)$ is given by
\begin{align}
f_{W_t,\Wb_t}(w,m)
	&=	\frac{2(2m-w)}{t \sqrt{2\pi t}} \ee^{-(2m-w)^2/(2t)} , &
m
	&>0 , &
w
	&\leq m .
\end{align}
(a) Compute the probability density function $f_{\Wb_t}$ of $\Wb_t$.  Show that $\Wb_t$ has the same probability density function as $|W_t|$.
\\[0.5em]
(b) Compute $\Eb \Wb_t$.
\\[0.5em]
(c) Fix $m > 0$ and $w \leq m$.  Compute $\Pb( \Wb_t \geq m | W_t = w )$.
\\[0.5em]
Note: the density of $W_t$ is $f_{W_t}(w) = \frac{1}{\sqrt{2 \pi t}} \ee^{-w^2/(2t)}$.
\end{exercise}

%
%

\chapter{Stochastic calculus}
\label{chapter:stoch.calc}

The notes from this chapter are taken primarily from \cite[Chapters 4 and 5]{shreve2004}.


\section{It\^o Integrals}
\label{sec:ito.integral}

If a function $g : [0,T] \to \Rb$ satisfies $g \in C^1([0,T])$, then we can define
\begin{align}
\int_0^T f(t) \dd g(t)
	&:=	\int_0^T f(t) g'(t) \dd t .
\end{align}
where the right-hand side is a Riemann integral with respect to $t$.  However, as Brownian motion is not differentiable, this method of integration will not work when we consider integrals of the form $\int_0^T \Del_t \dd W_t$.  Our goal is to make sense of $\int_0^T \Del_t \dd W_t$.

\begin{assumption}
In what follows $W=(W_t)_{t \geq 0}$ will always represent a Brownian motion and $\Fb=(\Fc_t)_{t \geq 0}$ will always be a filtration for this Brownian motion.
We shall assume the integrand $\Del = (\Del_t)_{t \geq 0}$ is adapted to $\Fb$, meaning $\Del_t \in \Fc_t$ for all $t$.
\end{assumption}

\noindent
Note that the process $\Del$ \emph{can} and, in many cases, \emph{will} be random.  However, the information available in $\Fc_t$ will always be sufficient to determine the value of $\Del_t$ at time $t$.  Also note, since $(W_T - W_t) \ind \Fc_t$ for $T>t$, it follows that $(W_T - W_t) \ind \Del_t$.  In other words, future increments of Brownian motion are independent of the $\Del$ process.

\subsection*{It\^o integrals for simple integrands}
To begin let us assume that $\Del$ is a \emph{simple process}, meaning $\Del$ is of the form
\begin{align}
\Del_t
	&=	\sum_{j=0}^{n-1} \Del_{t_j} \Ib_{\{t_j \leq t < t_{j+1}\}} , &
0
	&=	t_0 < t_1 < \ldots t_n = T , &
\Del_{t_j}
	&\in \Fc_{t_j} . \label{eq:Delta.simple}
\end{align}
Since the process $\Del$ is constant over intervals of the form $[t_j,t_{j+1})$, it makes sense to define
\begin{align}
I_T = \int_0^T \Del_t \dd W_t
	&:=	\sum_{j=0}^{n-1} \Del_{t_j} ( W_{t_{j+1}} - W_{t_j} ) . &
	&		\text{(for $\Del$ a simple process)} \label{eq:Ito.simple}
\end{align}
Let us establish some properties of the process $I = (I_t)_{t \geq 0}$.

\begin{theorem}
The process $I = (I_t)_{t \geq 0}$ defined in \eqref{eq:Ito.simple} is a martingale with respect to the filtration $(\Fc_t)_{t \geq 0}$.
\end{theorem}

\begin{proof}
Without loss of generality assume $T = t_{n}$ and $t = t_i$ for some $0 \leq i \leq n-1$ (we can always re-define our time grid so that this is true).  Then we have
\begin{align}
\Eb[ I_T | \Fc_t ]
	&=	\sum_{j=0}^{n-1} \Eb[ \Del_{t_j} ( W_{t_{j+1}} - W_{t_j} ) | \Fc_{t_i} ] \\
	&=	\sum_{j=0}^{i-1} \Del_{t_j} ( W_{t_{j+1}} - W_{t_j} ) + \sum_{j=i}^{n-1} \Eb[ \Del_{t_j} ( W_{t_{j+1}} - W_{t_j} ) | \Fc_{t_i} ] \\
	&=	\sum_{j=0}^{i-1} \Del_{t_j} ( W_{t_{j+1}} - W_{t_j} ) + \sum_{j=i}^{n-1} \Eb[ \Del_{t_j} \Eb[ W_{t_{j+1}} - W_{t_j} | \Fc_{t_j}] | \Fc_{t_i} ] \\
	&=	\sum_{j=0}^{i-1} \Del_{t_j} ( W_{t_{j+1}} - W_{t_j} ) = I_t .
\end{align}
Thus, the process $I$ is a martingale, as claimed.
\end{proof}

\noindent
Note that $I_0 = 0$ and therefore $\Eb I_t = 0$ for all $t \geq 0$.

\begin{theorem}[It\^o Isometry]
The process $I = (I_t)_{t \geq 0}$ defined in \eqref{eq:Ito.simple} satisfies
\begin{align}
\Vb I_T = \Eb I_T^2 
	&= \Eb \int_0^T \Del_t^2 \dd t. \label{eq:Ito.var}
\end{align}
\end{theorem}

\begin{proof}
We have
\begin{align}
\Eb I_T^2
	&=	\sum_{j=0}^{n-1} \sum_{i=0}^{n-1} \Eb \Del_{t_i} \Del_{t_j} ( W_{t_{i+1}} - W_{t_i} )  ( W_{t_{j+1}} - W_{t_j} ) \\
	&=	\sum_{j=0}^{n-1} \Eb \Del_{t_j}^2 ( W_{t_{j+1}} - W_{t_j} )^2
			+ 2 \sum_{j=0}^{n-1} \sum_{i=0}^{j-1} \Eb \Del_{t_i} \Del_{t_j} ( W_{t_{i+1}} - W_{t_i} )  ( W_{t_{j+1}} - W_{t_j} ) \\
	&=	\sum_{j=0}^{n-1} \Eb \Del_{t_j}^2 \Eb[  ( W_{t_{j+1}} - W_{t_j} )^2 | \Fc_{t_j}]
			+ 2 \sum_{j=0}^{n-1} \sum_{i=0}^{j-1} \Eb \Del_{t_i} \Del_{t_j} ( W_{t_{i+1}} - W_{t_i} )  \Eb[ ( W_{t_{j+1}} - W_{t_j} ) | \Fc_{t_j} ] \\
	&=	\sum_{j=0}^{n-1} \Eb \Del_{t_j}^2 ( {t_{j+1}} - {t_j} )
	=		\Eb \int_0^T \Del_t^2 \dd t .
\end{align}
\end{proof}

\begin{theorem}
\label{thm:ito.integral.QV}
The process $I = (I_t)_{t \geq 0}$ defined in \eqref{eq:Ito.simple} satisfies
\begin{align}
[I,I]_T 
	&=  \int_0^T \Del_t^2 \dd t. \label{eq:Ito.QV}
\end{align}
\end{theorem}

\begin{proof}
Fix a partition $\Pi$ of $[t_j,t_{j+1}]$
\begin{align}
 t_j 
	&= s_0 < s_1 < \ldots < s_m = t_{j+1} .
\end{align}
From the definition of Quadratic variation, we compute
\begin{align}
[I,I]_{t_{j+1}}-[I,I]_{t_j}
	&=	\lim_{\norm{\Pi} \to 0 } \sum_{i=0}^{m-1} \( I_{s_{i+1}} - I_{s_i} \)^2
	=		\Del_{t_j}^2 \lim_{\norm{\Pi} \to 0 }  \sum_{i=0}^{m-1} \( W_{s_{i+1}} - W_{s_i} \)^2 \\
	&=	\Del_{t_j}^2 \( [W,W]_{t_{j+1}} - [W,W]_{t_j} \)
	=		\Del_{t_j}^2 \( {t_{j+1}} - {t_j} \) .
\end{align}
Thus, we obtain
\begin{align}
[I,I]_T
	&=	\sum_{j=0}^{n-1} \( [I,I]_{t_{j+1}}-[I,I]_{t_j} \)
	=		\sum_{j=0}^{n-1} \Del_{t_j}^2 \( {t_{j+1}} - {t_j} \)
	=		\int_0^T \Del_t^2 \dd t ,
\end{align}
as claimed.
\end{proof}

\noindent
When we computed $\Vb W_T$ and $[W,W]_T$ we found that these two quantities were equal, even though the computations for these quantities were completely different.
From \eqref{eq:Ito.var} and \eqref{eq:Ito.QV} we now see how the variance and quadratic variation of a stochastic process can be different.  Note that $\Vb I_T$ is a non-random constant, whereas $[I,I]_T$ is random.

\subsection*{It\^o integrals for general integrands}
Clearly, it is rather restrictive to limit ourselves to integrands $\Del$ that are simple processes.  We now allow the process $\Del$ to be any process that is adapted to $\Fb=(\Fc_t)_{t \geq 0}$ and which satisfies the following integrability condition
\begin{align}
\Eb \int_0^T \Del_t^2 \dd t 
	< \infty . \label{eq:L2.condition}
\end{align}
To construct an It\^o integral with $\Del$ as the integrand, we first approximate $\Del$ by a simple process
\begin{align}
\Del_t \approx \Del_t^{(n)}
	&:=	\sum_{j=0}^{n-1} \Del_{t_j} \Ib_{\{ t_j \leq t < t_{j+1} \}} , &
0
	&\leq t_0 < t_1 < \ldots < t_n = T .
\end{align}
As $n \to \infty$ the process $\Del^{(n)}$ converges to $\Del$ in the sense that
\begin{align}
\lim_{n \to \infty} \Eb \int_0^T \( \Del_t - \Del_t^{(n)} \)^2 \dd t
	&=	0 . \label{eq:L2.condition}
\end{align}
We now \emph{define} the It\^o integral for a general integrand $\Del$ by
\begin{align}
I_T
	\equiv \int_0^T \Del_t \dd W_t
	&:=			\lim_{n \to \infty} \int_0^T \Del_t^{(n)} \dd W_t . \label{eq:Ito.def}
\end{align}
Note that the integrals $\int_0^T \Del_t^{(n)} \dd t$ are well-defined for every $n$, since $\Del^{(n)}$ is a simple process.  
Furthermore, the condition \eqref{eq:L2.condition} ensures that the limit exists in $L^2(\Om,\Fc,\Pb)$.
The It\^o integral for general integrands inherits the properties we established for simple integrands.

\begin{theorem}
\label{thm:ito.props}
Let $W$ be a Brownian motion and let $\Fb=(\Fc_t)_{t \geq 0}$ be a filtration for this Brownian motion.
Let $\Del = (\Del_t)_{0 \leq t \leq T}$ be adapted to the filtration $\Fb$ and satisfy \eqref{eq:L2.condition}.  
Let $I = (I_t)_{0 \leq t \leq T}$  be given by $I_t = \int_0^t \Del_s \dd W_s$, where the integral is defined as in \eqref{eq:Ito.def}.  Then the process $I$ has the following properties.
\begin{enumerate*}
\item The sample paths of $I$ are continuous.
\item The process $I$ is adapted to the filtration $\Fb$.  That is, $I_t \in \Fc_t$ for all $t$.
\item If $\Gam = (\Gam_t)_{0 \leq t \leq T}$ satisfies the same conditions as $\Del$, then
\begin{align}
\int_0^T \( a \Del_t + b \Gam_t \) \dd W_t
	&=	a \int_0^T \Del_t \dd W_t + b \int_0^T \Gam_t \dd W_t ,
\end{align}
where $a$ and $b$ are constants.
\item The process $I$ is a martingale with respect to the filtration $\Fb$.  
\item We have the It\^o isometry $\Eb I_T^2 = \Eb \int_0^T \Del_t^2 \dd t$.
\item The quadratic variation of $I$ is given by $[I,I]_T = \int_0^T \Del_t^2 \dd t$.
\end{enumerate*}
\end{theorem}


\section{It\^o formula}
\label{sec:ito.integral}
If $f$ and $g$ are differentiable functions then we can compute
\begin{align}
\frac{\dd}{\dd t}f(g(t))
	&=	f'(g(t)) \cdot g'(t) .
\end{align}
Thus, multiplying through by $\dd t$ we obtain
\begin{align}
\dd f(g(t))
	&=	f'(g(t)) \cdot g'(t) \dd t = f(g(t)) \dd g(t) , &
\dd g(t)
	&=	g'(t) \dd t .
\end{align}
This formula is correct in the sense that
\begin{align}
f(g(T)) - f(g(0))
	&=	\int_0^T \dd f(g(t)) = \int_0^T f'(g(t)) \cdot g'(t) \dd t 
	=		\int_0^T f'(g(t)) \dd g(t) .
\end{align}
Although paths of a Brownian motion $W$ are not differentiable, we might suspect that 
\begin{align}
f(W_T) - f(W_0)
	&=	\int_0^T \dd f(W_t) \stackrel{?}{=} \int_0^T f'(W_t) \dd W_t . &
	&   \text{(Not correct!)} \label{eq:not.correct}
\end{align}
Unfortunately, as indicated above, equation \eqref{eq:not.correct} is not correct.

\begin{theorem}
\label{thm:ito.bm}
Let $W = (W_t)_{t \geq 0}$ be a Brownian motion and suppose $f : \Rb \to \Rb$ satisfies $f \in C^2(\Rb)$ and $\Eb \int_0^T ( f'(W_t) )^2 \dd t < \infty$ for all $T \geq 0$.  Then, for any $T \geq 0$ we have
\begin{align}
f(W_T) - f(W_0)
	&=	\int_0^T f'(W_t) \dd W_t + \tfrac{1}{2} \int_0^T f''(W_t) \dd t . \label{eq:ito.bm}
\end{align}
\end{theorem}

\begin{proof}
We shall simply sketch the proof of Theorem \ref{thm:ito.bm}.
Suppose for simplicity that $f$ is analytic (i.e., that $f$ is equal to its power series expansion at every point).  
Let $0 = t_0 < t_1 < \ldots < t_n = T$ be a partition $\Pi$ of $[0,T]$.
Then
\begin{align}
\int_0^T \dd f(W_t)
	&=	f(W_T) - f(W_0) 
	=		\sum_{j=0}^{n-1} \( f(W_{t_{j+1}}) - f(W_{t_j}) \) 
	=		\sum_{j=0}^{n-1} \( A_j + B_j + C_j \) ,
\end{align}
where we have defined
\begin{align}
A_j
	&:=	f'(W_{t_j}) \( W_{t_{j+1}} - W_{t_j} \) , \\
B_j
	&:=  \tfrac{1}{2} f''(W_{t_j}) \( W_{t_{j+1}} - W_{t_j} \)^2 , \\
C_j
	&:=	\tfrac{1}{3!} f'''(W_{t_j}) \( W_{t_{j+1}} - W_{t_j} \)^3 + \ldots .
\end{align}
In the limit as $\norm{\Pi} \to 0$ we have
\begin{align}
\sum_{j=0}^{n-1} A_j
	&\to \int_0^T f'(W_t) \dd W_t , &
\sum_{j=0}^{n-1} B_j
	&\to \int_0^T \tfrac{1}{2} f''(W_t) \dd t , &
\sum_{j=0}^{n-1} C_j
	&\to 0 .
\end{align}
\end{proof}

\begin{example}
\label{ex:WdW}
What is $\int_0^T W_t \dd W_t$?  To answer this question, consider $f(W_t)$ with $f(x) = x^2$.  According to equation \eqref{eq:ito.bm} we have
\begin{align}
W_T^2 - W_0^2
	&=	\int_0^T 2 W_t \dd W_t + \int_0^T \dd t
\end{align}
where we have used $f'(x) = 2 x$ and $f''(x) = 2$.  Noting that $W_0 = 0$ and solving for $\int_0^T W_t \dd W_t$ we obtain
\begin{align}
\int_0^T W_t \dd W_t
	&=	\tfrac{1}{2} W_T^2 - \tfrac{1}{2} T .
\end{align}
\end{example}

\noindent
Not surprisingly, there are stochastic processes that are not adequately described by Brownian motion alone.  However, a large class of stochastic processes can be constructed from Brownian motion.

\begin{definition}
\label{def:ito.process}
Let $W = (W_t)_{t \geq 0}$ be a Brownian motion and let $\Fb=(\Fc_t)_{t \geq 0}$ be a filtration for this Brownian motion.  An \emph{It\^o process} is any process $X = (X_t)_{t \geq 0}$ of the form
\begin{align}
X_t
	&= X_0 + \int_0^t \Theta_s \dd s + \int_0^t \Del_s \dd W_s  , \label{eq:ito.process}
\end{align}
where $\Theta = (\Theta_t)_{t \geq 0}$ and $\Del = (\Del_t)_{t \geq 0}$ are adapted to the filtration $\Fb$ and satisfy
\begin{align}
\int_0^T | \Theta_t | \dd t 
	&< \infty , &
\Eb \int_0^T \Del_t^2 \dd t
	&< \infty , &
	&\forall \, T \geq 0 ,
\end{align}
and $X_0 \in \Fc_0$.
\end{definition}

\noindent
We sometimes write an It\^o process in \emph{differential form}
\begin{align}
\dd X_t
	&=	\Theta_t \dd t + \Del_t \dd W_t . \label{eq:differential.form}
\end{align}
Expression \eqref{eq:differential.form} literally means that $X$ satisfies \eqref{eq:ito.process}.  Informally, the differential form can be understood as follows: in a small interval of time $\del t$ , the process $X$ changes according to
\begin{align}
X_{t + \del t} - X_t
	&\approx	\Theta_t \del t + \Del_t \( W_{t+\del t} - W_t \) . \label{eq:euler}
\end{align}
In fact, noting that $W_{t+\del t} - W_t \sim \Nc(0,\del t)$ and $W_{t+\del t} - W_t \ind \Fc_t$, one can use expression \eqref{eq:euler} to simulate the increment $X_{t + \del t} - X_t$.  This way of simulating $X$ is called the \emph{Euler scheme} and is the workhorse of many Monte Carlo methods.

\begin{lemma}
\label{lem:ito.process.QV}
The quadratic variation $[X,X]_T$ of an It\^o process \eqref{eq:ito.process} is given by
\begin{align}
[X,X]_T
	&=	\int_0^T \Del_t^2 \dd t .
\end{align}
\end{lemma}

\begin{proof}
We sketch the proof of Lemma \ref{lem:ito.process.QV}.
Let $0 = t_0 < t_1 < \ldots < t_n = T$ be a partition $\Pi$ of $[0,T]$.  By definition we have
\begin{align}
[X,X]_T
	&=	\lim_{\norm{\Pi}\to 0} \sum_{j=0}^{n-1} \( X_{t_{j+1}} - X_{t_j}\)^2 
	=		\lim_{\norm{\Pi}\to 0} \sum_{j=0}^{n-1} \( A_j + B_j + C_j \)
\end{align}
where we have defined
\begin{align}
A_j
	&:=	 \( I_{t_{j+1}} - I_{t_j} \)^2 , \\
B_j
	&:= \( J_{t_{j+1}} - J_{t_j} \)^2 , \\
C_j
	&:=	2 \( I_{t_{j+1}} - I_{t_j} \) \( J_{t_{j+1}} - J_{t_j} \)   .
\end{align}
with
\begin{align}
I_t
	&=	\int_0^t \Del_s \dd W_s , &
J_t
	&=	\int_0^t \Theta_s \dd s .
\end{align}
In the limit as $\norm{ \Pi } \to 0$ we obtain
\begin{align}
\sum_{j=0}^{n-1} A_j
	& \to [I,I]_T , &
\sum_{j=0}^{n-1} B_j
	& \to 0 , &
	\sum_{j=0}^{n-1} C_j
	& \to 0 . &
\end{align}
The proof of these limits is similar to the proof of Theorem \ref{thm:ito.integral.QV}.
\end{proof}

\begin{definition}
Let $X = (X_t)_{t \geq 0}$ be an It\^o process, as described in Definition \ref{def:ito.process}.  Let $\Gam = (\Gam_t)_{t \geq 0}$ be adapted to the filtration of the Brownian motion $\Fb=(\Fc_t)_{t \geq 0}$.  We define 
\begin{align}
\int_0^T \Gam_t \dd X_t
	&:=	\int_0^T \Gam_t \Theta_t \dd t + \int_0^T \Gam_t \Del_t \dd W_t ,
\end{align}
where we assume
\begin{align}
\int_0^T | \Gam_t \Theta_t | \dd t 
	&< \infty , &
\Eb \int_0^T \( \Gam_t \Del_t \)^2 \dd t
	&< \infty , &
	& \forall \, T \geq 0 .
\end{align}
\end{definition}

\begin{theorem}[It\^o formula in one dimension]
\label{thm:ito.X}
Let $X = (X_t)_{t \geq 0}$ be an It\^o process of the form \eqref{eq:ito.process} and suppose $f : \Rb \to \Rb$ satisfies $f \in C^2(\Rb)$ and $\Eb \int_0^T ( f'(X_t) \Del_t )^2 \dd t < \infty$ for all $T \geq 0$.  Then, for any $T \geq 0$ we have
\begin{align}
f(X_T) - f(X_0)
	&=	\int_0^T f'(X_t) \dd X_t + \tfrac{1}{2} \int_0^T f''(X_t) \dd [X,X]_t . \label{eq:ito.X}
\end{align}
\end{theorem}

\begin{proof}
The proof of Theorem \ref{thm:ito.X} is very similar to the proof of Theorem \ref{thm:ito.bm}.  We outline the proof here.
Suppose for simplicity that $f$ is analytic (i.e., that $f$ is equal to its power series expansion at every point).  
Let $0 = t_0 < t_1 < \ldots < t_n = T$ be a partition $\Pi$ of $[0,T]$.
Then
\begin{align}
\int_0^T \dd f(X_t)
	&=	f(X_T) - f(X_0) 
	=		\sum_{j=0}^{n-1} \( f(X_{t_{j+1}}) - f(X_{t_j}) \) 
	=		\sum_{j=0}^{n-1} \( A_j + B_j + C_j \) ,
\end{align}
where we have defined
\begin{align}
A_j
	&:=	f'(X_{t_j}) \( X_{t_{j+1}} - X_{t_j} \) , \\
B_j
	&:=  \tfrac{1}{2} f''(X_{t_j}) \( X_{t_{j+1}} - X_{t_j} \)^2 , \\
C_j
	&:=	\tfrac{1}{3!} f'''(X_{t_j}) \( X_{t_{j+1}} - X_{t_j} \)^3 + \ldots .
\end{align}
In the limit as $\norm{\Pi} \to 0$ we have
\begin{align}
\sum_{j=0}^{n-1} A_j
	&\to \int_0^T f'(X_t) \dd X_t , &
\sum_{j=0}^{n-1} B_j
	&\to \int_0^T \tfrac{1}{2} f''(X_t) \dd [X,X]_t , &
\sum_{j=0}^{n-1} C_j
	&\to 0 .
\end{align}
\end{proof}

\noindent
In differential form, with $X$ given by \eqref{eq:differential.form}, It\^o's formula becomes
\begin{align}
\dd f(X_t)
	&=	f'(X_t) \dd X_t + \tfrac{1}{2} f''(X_t) \dd [X,X]_t \\
	&=	f'(X_t) \( \Theta_t \dd t + \Delta_t \dd W_t \) + \tfrac{1}{2} f''(X_t) \Del_t^2 \dd t , \label{eq:ito.differential}
\end{align}
where we have used $\dd [X,X]_t = \Del_t^2 \dd t$.  Perhaps the easiest way to remember \eqref{eq:ito.differential} is to use the following two-step procedure:
\begin{enumerate*}
\item Expand $f(X_t + \dd X_t) - f(X_t)$ to second order about the point $X_t$
\begin{align}
\dd f(X_t) = f(X_t + \dd X_t) - f(X_t)
	&=	f'(X_t) \dd X_t + \tfrac{1}{2} f''(X_t) (\dd X_t)^2 , \label{eq:df.X}
\end{align}
\item Insert the differential $\dd X_t = \Theta_t \dd t + \Del_t \dd W_t$ into \eqref{eq:df.X}, expand $(\dd X_t)^2$ and use the rules
\begin{align}
\dd W_t \dd W_t
	&=	\dd t , &
\dd W_t \dd t
	&= 0 , &
\dd t \dd t
	&=	0 .
\end{align}
\end{enumerate*}
The resulting formula gives the correct expression for $\dd f(X_t)$.

\begin{example}
Let $X$ be an It\^o process with the following dynamics
\begin{align}
\dd X_t
	&=	\mu_t X_t \dd t + \sig_t X_t \dd W_t . \label{eq:GBM.general}
\end{align}
Assuming $\mu = (\mu_t)_{t \geq 0}$ and $\sig = (\sig_t)_{t \geq 0}$ are bounded above and below and $X_0>0$, the process $X$ remains strictly positive.
We call $X$ a \emph{generalized geometric Brownian motion}.  The ``geometric'' part refers to the fact that the relative step size $\dd X_t / X_t$ has dynamics $\mu_t \dd t + \sig_t \dd W_t$.  The ``generalized'' part refers to the fact that the processes $\sig$ and $\mu$ are stochastic rather than constant.
Define $Y_t = X_t^p$.  What is $\dd Y_t$?  Let $f(x) = x^p$.  Then $f'(x) = p x^{p-1}$ and $f''(x) = p(p-1)x^{p-2}$.  Thus, we have
\begin{align}
\dd Y_t
	&=	\dd f(X_t)
	=		p X_t^{p-1} \dd X_t + \tfrac{1}{2} p (p-1) X_t^{p-2} (\dd X_t)^2 \\
	&=	p X_t^{p-1} \( \mu_t X_t \dd t + \sig_t X_t \dd W_t \) + \tfrac{1}{2} p (p-1) X_t^{p-2} \( \mu_t X_t \dd t + \sig_t X_t \dd W_t \)^2 \\
	&=	p X_t^{p-1} \( \mu_t X_t \dd t + \sig_t X_t \dd W_t \) + \tfrac{1}{2} p (p-1) X_t^{p-2} \sig_t^2 X_t^2 \dd t  \\
	&=	\( p \mu_t + \tfrac{1}{2} p (p-1) \sig_t^2 \) X_t^p \dd t + p \sig_t X_t^p \dd W_t \\
	&=	\( p \mu_t + \tfrac{1}{2} p (p-1) \sig_t^2 \) Y_t \dd t + p \sig_t Y_t \dd W_t .
\end{align}
We see from the last line that $Y = (Y_t)_{t \geq 0}$ is also a generalized geometric Brownian motion.
\end{example}

\begin{example}
Let $X$ have generalized geometric Brownian motion dynamics as in \eqref{eq:GBM.general}.  
We would like to find an explicit expression for $X_t$ (i.e., an expression of the form $X_t = \ldots$ where $\ldots$ does not contain $X$).  To this end,we let $Y_t = \log X_t$.  With $f(x) = \log x$ we have $f'(x) = 1/x$ and $f''(x) = - 1/x^2$.  Thus, we have
\begin{align}
\dd Y_t
	&=	\frac{1}{X_t} \dd X_t + \frac{1}{2} \frac{-1}{X_t^2} (\dd X_t)^2 
	 =	\( \mu_t - \tfrac{1}{2} \sig_t^2 \) \dd t + \sig_t \dd W_t .
\end{align}
Thus, we have
\begin{align}
X_T
	&=	\exp( Y_T ) = \exp \( Y_0 + \int_0^T \( \mu_t - \tfrac{1}{2} \sig_t^2 \) \dd t + \int_0^T \sig_t \dd W_t \) \\
	&=	X_0 \exp \( \int_0^T \( \mu_t - \tfrac{1}{2} \sig_t^2 \) \dd t + \int_0^T \sig_t \dd W_t \) , \label{eq:X.exp}
\end{align}
where we have used $Y_0 = \log X_0$.
\end{example}

\begin{proposition}
\label{prop:normal.ito}
Let $W = (W_t)_{t \geq 0}$ be a Brownian motion.
Suppose $g: \Rb_+ \to \Rb$ is a deterministic function satisfying $\int_0^T g^2(t) \dd t < \infty$.  Then
\begin{align}
I_T
	&:= \int_0^T g(t) \dd W_t \sim \Nc( 0, v(T)) , &
v(T)
	&=	\int_0^T g^2(t) \dd t .
\end{align}
\end{proposition}

\begin{proof}
Set $\mu_t = 0$ and $\sig_t = u g(t)$ in \eqref{eq:GBM.general} where $u$ is a constant.  And suppose $X_0 = 1$.  Then we have
\begin{align}
X_T
	= 1 + \int_0^T u g(t) X_t \dd W_t ,
\end{align}
which is a martingale since It\^o integrals are martingales.  From \eqref{eq:X.exp}, we know that $X_T$ can be written explicitly as
\begin{align}
X_T
	&=	\exp \( - \tfrac{1}{2} \int_0^T u^2 g^2(t) \dd t + \int_0^T u g(t) \dd W_t \) 
	=		\exp \( - u^2 \tfrac{1}{2} \int_0^T g^2(t) \dd t + u I_T \) .
\end{align}
Since $X$ is a martingale, we have
\begin{align}
1
	&= \Eb X_T 
	=		\Eb \exp \( - u^2 \tfrac{1}{2} \int_0^T g^2(t) \dd t + u I_T \) &
	&\Rightarrow&
\Eb \exp \( u I_T \)
	&=	\exp \( u^2 \tfrac{1}{2} \int_0^T g^2(t) \dd t \) .
\end{align}
Note that $\Eb \exp \( u I_T \) = M_{I_T}(u)$ is the moment generating function of $I_T$ and $\exp \( u^2 \tfrac{1}{2} \int_0^T g^2(t) \dd t \)$ is the moment generating function of a normal random variable with mean zero and variance $v(T) = \int_0^T g^2(t) \dd t$.  Thus, $I_T \sim \Nc(0,v(T))$, as claimed.
\end{proof}


\section{Multivariate stochastic calculus}
\label{sec:multivariate}

\begin{definition}
\label{def:d.bm}
A  \emph{$d$-dimensional Brownian motion} is a process
\begin{align}
W 
	&=	( W_t^1 , W_t^2, \ldots, W_t^d)_{t \geq 0}
\end{align}
with the following properties.
\begin{enumerate*}
\item Each $W^i = (W_t^i)_{t \geq 0}$, $i = 1, 2, \ldots d$, is a one-dimensional Brownian motion.
\item The processes $(W^i)_{1 \leq i \leq d}$ are independent.
\end{enumerate*}
A \emph{filtration for $W$} is a collection of $\sig$-algebras $\Fb=(\Fc_t)_{t \geq 0}$ such that
\begin{enumerate*}
\item Information accumulates: $\Fc_s \subset \Fc_t$ for all $0 \leq s < t$.
\item Adaptivity: $W_t \in \Fc_t$ for all $t \geq 0$.
\item Independent increments: for $0 \leq s < t$ we have $W_t - W_s \ind \Fc_s$.
\end{enumerate*}
\end{definition}

\begin{theorem}
\label{thm:bm.CoV}
Let $W=( W_t^1 , W_t^2, \ldots, W_t^d)_{t \geq 0}$ be a $d$-dimensional Brownian motion.  The covariation of independent components of $W$ is zero: $[W^i,W^j]_T = 0$ for all $i \neq j$ and $T \geq 0$.
\end{theorem}

\begin{proof}
Let $0 = t_0 < t_1 < \ldots < t_n = T$ be a partition $\Pi$ of $[0,T]$.  The sampled covariation $C_\Pi$ of $W^i$ and $W^j$ is given by
\begin{align}
C_\Pi
	&=	\sum_{k=0}^{n-1} \( W_{t_{k+1}}^i - W_{t_k}^i \)\( W_{t_{k+1}}^j - W_{t_k}^j \)
\end{align}
Since $\Eb \( W_{t_{k+1}}^i - W_{t_k}^i \)\( W_{t_{k+1}}^j - W_{t_k}^j \) = 0$, we clearly have $\Eb C_\Pi = 0$.  
Next, we compute the variance of $C_\Pi$.  We have
\begin{align}
\Vb C_\Pi
	&=	\Eb C_\Pi^2
	=		\sum_{k=0}^{n-1} \sum_{l=0}^{n-1} 
			\Eb \( W_{t_{k+1}}^i - W_{t_k}^i \) \( W_{t_{k+1}}^j - W_{t_k}^j \) \( W_{t_{l+1}}^i - W_{t_l}^i \) \( W_{t_{l+1}}^j - W_{t_l}^j \) \\
	&=	\sum_{k=0}^{n-1} \Eb \( W_{t_{k+1}}^i - W_{t_k}^i \)^2 \( W_{t_{k+1}}^j - W_{t_k}^j \)^2 \\ &\quad
			+ 2 \sum_{k=0}^{n-1} \sum_{l=0}^{k-1} 
			\Eb \( W_{t_{k+1}}^i - W_{t_k}^i \) \( W_{t_{k+1}}^j - W_{t_k}^j \) \( W_{t_{l+1}}^i - W_{t_l}^i \) \( W_{t_{l+1}}^j - W_{t_l}^j \) \\
	&=	\sum_{k=0}^{n-1}\({t_{k+1}} - {t_k} \) \( {t_{k+1}} -{t_k}\) \\
	&\leq \norm{\Pi} \sum_{k=0}^{n-1} \({t_{k+1}} - {t_k} \)
	=		\norm{\Pi} T .	
\end{align}
Thus, $\Eb C_\Pi \to 0$ and $\Vb C_\Pi \to 0$ as $\norm{\Pi} \to 0$, which proves that $[W^i,W^j]_T := \lim_{\norm{\Pi} \to 0} C_\Pi = 0$.
\end{proof}

\noindent
Theorem \ref{thm:bm.CoV} can be used to derive the covariation of two It\^o processes $X^i$ and $X^j$.


\begin{theorem}
\label{thm:CoVar.XiXj}
Let $X^i=(X_t^i)_{t \geq 0}$ $i = 1, 2, \ldots, n$ be the It\^o processes given by
\begin{align}
\dd X_t^i
	&=	\Theta_t^i \dd t + \sum_{j=1}^d \sig_t^{ij} \dd W_t^j , &
i
	&=	1, 2, \ldots, n , \label{eq:dX.d} 
\end{align} 
where $W = (W_t^1, W_t^2, \ldots, W_t^{d})_{t \geq 0}$ is a $d$-dimensional Brownian motion.  Then
\begin{align}
\dd [ X^i, X^j ]_t
	&=	\sum_{k=1}^d \sig_t^{ik} \sig_t^{jk} \dd t . 
\end{align}
\end{theorem}

\noindent
We will not prove Theorem \ref{thm:CoVar.XiXj}.  Rather, we simply remark that it can be obtained \emph{informally} by writing
\begin{align}
\dd [X^i,X^j]_t
	&=	\dd X_t^i \dd X_t^j , \label{eq:CoVar.XY}
\end{align}
inserting expression \eqref{eq:dX.d} into \eqref{eq:CoVar.XY} and using the multiplication rules
\begin{align}
\dd W_t^i \dd W_t^j
	&=	\del_{ij} \dd t , &
\del_{ij}
	&=	\left\{ \begin{array}{ll}
			1 , & i=j , \\
			0 , & i\neq j ,
			\end{array} \right. &
\dd W_t^j \dd t
	&= 0 , &
\dd t \dd t
	&=	0 . \label{eq:rules.2}
\end{align}
Note that $\dd [X^i,X^j]_t = 0$ unless $X^j$ and $X^j$ are driven by at least one common one-dimensional Brownian motion.
\\[0.5em]
We can now give a $n$-dimensional version of It\^o's Lemma.  We present the formula in differential form, as it is written more compactly in this way.


\begin{theorem}[It\^o formula in $n$ dimensions]
\label{thm:ito.XY}
Let $X = (X_t^1, X_t^2, \ldots, X_t^n)_{t \geq 0}$ be an $n$-dimensional It\^o process whose components are of the form \eqref{eq:dX.d}, and suppose $f : \Rb^n \to \Rb$ satisfies $f \in C^2(\Rb^n)$ and
\begin{align}
\Eb \int_0^T \Big( \sum_{i=1}^n \frac{\d f(X_t)}{\d x_i} \sig_t^{ij} \Big)^2 \dd t
	&<	\infty , &
j
	&=	1, 2, \ldots, d , &
T
	&\geq 0 .
\end{align}
Then we have
\begin{align}
\dd f(X_t)
	&=	\sum_{i=1}^n \frac{\d f(X_t)}{\d x_i} \dd X_t^i + \frac{1}{2} \sum_{i=1}^n \sum_{j=1}^n \frac{\d^2 f(X_t)}{\d x_i \d x_j}  \dd[X^i,X^j]_t .
\end{align}
\end{theorem}

\noindent
The proof of Theorem \ref{thm:ito.XY} is a straightforward extension of Theorem \ref{thm:ito.X} to the $n$-dimensional case and will not be presented here.  
\\[0.5em]
To obtain an explicit expression for $\dd f(X_t)$ in terms of $\dd W_t^1, \dd W_t^2, \ldots, \dd W_t^d$ and $\dd t$ we can repeat the same informal procedure we used in the one-dimensional case.
\begin{enumerate*}
\item Expand $\dd f(X_t) = f(X_t + \dd X_t) - f(X_t)$ about the point $X_t$ to second order
\begin{align}
\dd f(X_t)
	&=	\sum_{i=1}^n \frac{\d f(X_t)}{\d x_i} \dd X_t^i + \frac{1}{2} \sum_{i=1}^n \sum_{j=1}^n \frac{\d^2 f(X_t)}{\d x_i \d x_j}  \dd X_t^i \dd X_t^j . \label{eq:df.dX.dY}
\end{align}
\item Insert expression for $\dd X_t^i$ into \eqref{eq:df.dX.dY} and use the multiplication rules given in \eqref{eq:rules.2}.
\end{enumerate*}

\begin{example}[Product rule]
To compute $\dd (X_t Y_t)$ where $X$ and $Y$ are one-dimensional It\^o processes, we define $f(x,y) = xy$ and use $f_x = y$, $f_y = x$, $f_{xy}=1$ and $f_{xx} = f_{yy} = 0$ to compute
\begin{align}
\dd (X_t Y_t)
	&=	Y_t \dd X_t + X_t \dd Y_t + \dd [X,Y]_t .
\end{align}
\end{example}

\begin{example}[OU process]
\label{ex:OU}
An \emph{Ornstein-Uhlenbeck process} (OU process, for short) is an It\^o process $X=(X_t)_{t \geq 0}$ that satisfies
\begin{align}
\dd X_t
	&=	\kappa ( \theta - X_t) \dd t + \sig \dd W_t , \label{eq:dX.OU}
\end{align}
where $W = (W_t)_{t \geq 0}$ is a one-dimensional Brownian motion and $\kappa, \theta > 0$.
The OU process is \emph{mean-reverting} in the following sense.
If $X_t > \theta$ then $\kappa ( \theta - X_t) <0$ and the deterministic part of \eqref{eq:dX.OU} (i.e., the $\dd t$-term) pushes the process down towards $\theta$.
If $X_t < \theta$ then $\kappa ( \theta - X_t) > 0$ and the deterministic part of \eqref{eq:dX.OU} pushes the process up towards $\theta$.
The OU process mean-reverts to  the \emph{long-run mean} $\theta$.
We often call $\kappa$ the \emph{rate of mean reversion}, though this nomenclature is somewhat misleading since the instantaneous rate of mean reversion is actually $\kappa (\theta - X_t)$.
\\[0.5em]
We will find an explicit expression for $X_t$ and also compute $\Eb X_t$ and $\Vb X_t$.  To this end, let us define $Y_t = X_t - \theta$ so that
\begin{align}
\dd Y_t
	&=	- \kappa Y_t \dd t + \sig \dd W_t
\end{align}
Note that $Y$ is an OU process that mean-reverts to zero.  Next, we define $Z_t = f(t,Y_t) = \ee^{\kappa t} Y_t$.  We can use the two-dimensional It\^o formula to compute $\dd Z_t$.  Using $f_{yy}=0$ and the heuristic rules $\dd t \dd W_t = 0$ and $\dd t \dd t = 0$ we have
\begin{align}
\dd Z_t
	&=	f_t \dd t + f_y \dd Y_t + \tfrac{1}{2} f_{yy} \dd [Y,Y]_t \\
	&=	\kappa \ee^{\kappa t} Y_t \dd t + \ee^{\kappa t} \dd Y_t  \\
	&=	\kappa \ee^{\kappa t} Y_t \dd t + \ee^{\kappa t} \( - \kappa Y_t \dd t + \sig \dd W_t \) \\
	&=	\ee^{\kappa t} \sig \dd W_t . 
\end{align}
Thus, we have obtained an expression for $Z_t$:
\begin{align}
Z_t
	&=	Z_0 + \int_0^t \ee^{\kappa s} \sig \dd W_s .
\end{align}
Next, we use $Y_t = \ee^{- \kappa t} Z_t$ and $X_t = Y_t + \theta$ to obtain
\begin{align}
Y_t
	&=	\ee^{- \kappa t} Y_0 + \int_0^t \ee^{-\kappa(t-s)} \sig \dd W_s , \\
X_t
	&=	\theta + \ee^{- \kappa t} ( X_0 - \theta ) + \int_0^t \ee^{-\kappa(t-s)} \sig \dd W_s .
\end{align}
Note that $X_t$ has a normal distribution at every time $t>0$ because It\^o integrals with deterministic integrands are normally distributed random variables; see Proposition \ref{prop:normal.ito}.  Thus, the distribution of $X_t$ is completely determined by its mean and variance.  We have
\begin{align}
\Eb X_t
	&=	\theta + \ee^{- \kappa t} ( X_0 - \theta ) + \Eb \int_0^t \ee^{-\kappa(t-s)} \sig \dd W_s \\
	&=	\theta + \ee^{- \kappa t} ( X_0 - \theta ) , \\
\Vb X_t
	&=	\Vb \int_0^t \ee^{-\kappa(t-s)} \sig \dd W_s \\
	&=	\int_0^t \( \ee^{-\kappa(t-s)} \sig \)^2 \dd s
	=		\frac{\sig^2}{2 \kappa} \( 1 - \ee^{-2 \kappa t} \) .
\end{align}
where we have used Proposition \ref{prop:normal.ito}.
\end{example}

\noindent
In some cases, it will be useful to see if a process is in fact a Brownian motion.

\begin{theorem}[L\'evy characterization of Brownian motion]
\label{thm:Levy}
Let $M = (M_t)_{t \geq 0}$ be a martingale with respect to a filtration $\Fb=(\Fc_t)_{t \geq 0}$.  Suppose $M$ has continuous sample paths and satisfies $M_0 = 0$ and $[M,M]_t = t$ for all $t \geq 0$.  Then $M$ is a Brownian motion.
\end{theorem}

\begin{proof}
Let $f: \Rb_+ \times \Rb \to \Rb$ be $C^{1,2}$.  Note that
\begin{align}
\dd f(t,M_t)
	&=	\d_t f(t,M_t) \dd t + \d_m f(t,M_t) \dd M_t + \tfrac{1}{2} \d_m^2 f(t,M_t) \dd [M,M]_t \\
	&=	\( \d_t + \tfrac{1}{2} \d_m^2 \) f(t,M_t) \dd t + \d_m f(t,M_t) \dd M_t .
\end{align}
In integral form, we have
\begin{align}
f(T,M_T)
	&=	f(t,M_t) + \int_t^T \( \d_s + \tfrac{1}{2} \d_m^2 \) f(s,M_s) \dd s + \int_t^T \d_m f(s,M_s) \dd M_s .
\end{align}
Although we have not proved it, since $M$ is a martingale (by assumption), it follows that any integral of the form $I_t := \int_0^t \Del_s \dd M_s$, where $\Del = (\Del_t)_{t \geq 0}$ is adapted to $(\Fc_t)_{t \geq 0}$, is a martingale.  Thus, using the short-hand $\Eb_t \, \cdot \,:= \Eb[ \, \cdot \, | \Fc_t]$, we have
\begin{align}
\Eb_t f(T,M_T)
	&=	f(t,M_t) + \Eb_t \int_t^T \( \d_s + \tfrac{1}{2} \d_m^2 \) f(s,M_s) \dd s + \Eb \int_t^T \d_m f(s,M_s) \dd M_s \\
	&=	f(t,M_t) + \Eb_t \int_t^T \( \d_s + \tfrac{1}{2} \d_m^2 \) f(s,M_s) \dd s .
\end{align}
Now, let us fix a constant $u$ and define
\begin{align}
f(t,m)
	&=	\ee^{ u m - \tfrac{1}{2} u^2 t} .
\end{align}
It is easy to verify that
\begin{align}
\( \d_t + \tfrac{1}{2} \d_m^2 \) f
	&=	0 .
\end{align}
It follows that
\begin{align}
\Eb_t \ee^{ u M_T - \tfrac{1}{2} u^2 T}
	&=	f(t,M_t) = \ee^{ u M_t - \tfrac{1}{2} u^2 t} &
	&		\Rightarrow &
\Eb_t \ee^{ u (M_T - M_t ) }
	&=	\ee^{ \tfrac{1}{2} u^2 (T -t ) }. \label{eq:shows.ind}
\end{align}
By definition, $\Eb_t \ee^{ u (M_T - M_t) }$ is the $\Fc_t$-conditional moment generating function of $(M_T - M_t)$.  And $\ee^{ \tfrac{1}{2} u^2 (T-t) }$ is the moment generating function of a $\Nc(0,T-t)$ random variable.  It follows that $M_T - M_t \sim \Nc(0,T-t)$ for all $0 \leq t \leq T < \infty$, just like a Brownian motion.  Furthermore, we see that $M_T - M_t \ind \Fc_t$, as the right-hand-side of \eqref{eq:shows.ind} does not depend on $\Fc_t$.  As $M$ satisfies all properties of a Brownian motion, it must be a Brownian motion.
\end{proof}

\begin{example}[Correlated Brownian motions]
Let $B=(B_t)_{t \geq 0}$ be given by
\begin{align}
B_t
	&=	\rho W_t^1 + \rhob W_t^2 , &
\rhob
	&=	\sqrt{ 1 - \rho^2} , &
\rho
	&=	[-1,1] ,
\end{align}
where $W = (W_t^1,W_t^2)$ is a two-dimensional Brownian motion.  We will use Theorem \ref{thm:Levy} to show that $B$ is a Brownian motion.  It is clear that $B_0 = 0$ and $B$ has sample paths that are continuous.  It is also clear that $B$ is a martingale since $W^1$ and $W^2$ are martingales.  What remains is to show that $[B,B]_t = t$.  We have
\begin{align}
[B,B]_t
	&=	[ \rho W^1 + \rhob W^2 , \rho W^1 + \rhob W^2]_t \\
	&=	\rho^2 [W^1,W^1]_t + 2 \rho \rhob [W^1,W^2] + \rhob^2 [W^2,W^2]_t \\
	&=	\rho^2 t + 0 + \rhob^2 t = t .
\end{align}
Alternatively, one can simply compute $(\dd B_t)^2 = ( \rho \dd W_t^1 + \rhob^2 \dd W_t^2)^2$ and use $\dd W_t^i \dd W_t^j = \del_{ij} \dd t$ in order to show that $\dd B_t = \dd t$.  Note that
\begin{align}
[B,W^1]_t
	&=	[\rho W^1 + \rhob W^2,W^1]_t \\
	&=	\rho [W^1,W^1] + \rhob [W^2,W^1] = \rho t .
\end{align}
where we have used $[W^i,W^j]_t=\del_{ij}t$.
\end{example}


\section{Brownian bridge}
\label{sec:bridge}

\begin{definition}
A \emph{Gaussian process} $X = (X_t)_{t \geq 0}$ is any process that, for arbitrary times $t_1 < t_2 < \ldots t_n$, the joint distribution of $(X_{t_1},X_{t_2},\ldots,X_{t_n})$ is jointly normal.
\end{definition}

\noindent
As previously mentioned, the distribution of a normally distributed random vector is uniquely determined by its mean vector $\mv$ and covariance matrix.  Thus, for a Gaussian process, we are interested in
\begin{align}
m(t)
	&=	\Eb X_t , &
c(s,t)
	&=	\Eb \(X_t - m(t) \) \(X_s - m(s) \) .
\end{align}

\begin{example}
\label{eq:bm.gaussian}
A Brownian motion $W$ is a Gaussian process.  To see this, fix an arbitrary sequence of times $0=t_0 < t_1 < t_2 < \ldots t_n$.  Note that
\begin{align}
W_{t_k}
	&=	\sum_{j=0}^{k-1} ( W_{t_{j+1}} - W_{t_j} ) .
\end{align}
The increments are independent and normally distributed $W_{t_{j+1}} - W_{t_j} \sim \Nc(0,t_{j+1}-t_j)$.  It follows that the vector $(W_{t_1},W_{t_2},\ldots,W_{t_n})$ is jointly normal.  The mean vector is given by $\mv = (m(t_1),m(t_2),\ldots,m(t_n)) = (0,0,\ldots,0)$ and the covariance matrix $\Cv = (c(t_i,t_j))_{1 \leq i,j \leq n}$ has entries $c(t_i,t_j) = t_i \wedge t_j$; see equation \eqref{eq:BM-Cov}.
\end{example}

\begin{example}
\label{eq:ito.gaussian}
A stochastic integral with a deterministic integrand
\begin{align}
I_t
	&=	\int_0^t g(s) \dd W_s ,
\end{align}
is a Gaussian process.  To see this, fix an arbitrary sequence of times $0=t_0 < t_1 < t_2 < \ldots t_n$.  Note that
\begin{align}
I_{t_k}
	&=	\sum_{j=0}^{k-1} \int_{t_j}^{t_{j+1}} g(s) \dd W_s .
\end{align}
The increments are independent and normally distributed with mean zero and variance:
\begin{align}
\Vb \int_{t_j}^{t_{j+1}} g(s) \dd W_s
	&=	\int_{t_j}^{t_{j+1}} g^2(s) \dd s
\end{align}
It follows that the vector $(I_{t_1},I_{t_2},\ldots,I_{t_n})$ is jointly normal.  The mean vector is given by $\mv = (m(t_1),m(t_2),\ldots,m(t_n)) = (0,0,\ldots,0)$.  To compute the covariance matrix $\Cv = (c(t_i,t_j))_{1 \leq i,j \leq n}$, assume without loss of generality that $t_i < t_j$.  Then we have 
\begin{align}
c(t_i,t_j) 
	&=  \Eb \int_0^{t_i}g(s) \dd W_s \, \int_0^{t_j} g(u)  \dd W_u \\
	&=	\Eb \Big( \int_0^{t_i} g(s) \dd W_s \Big)^2 + \Eb \int_0^{t_i} g(s) \dd W_s \Eb \Big[ \int_{t_i}^{t_j} g(u) \dd W_u \Big| \Fc_{t_i} \Big] \\
	&=	\int_0^{t_i} g^2(s) \dd s
	=		\int_0^{t_i \wedge t_j} g^2(s) \dd s ,
\end{align}
where we have used Proposition \ref{prop:normal.ito}.
\end{example}

\begin{definition}[Brownian bridge, version I]
\label{def:bridge.bm}
Let $W=(W_t)_{t \geq 0}$ be a Brownian motion and fix $T>0$.  We define $X^{a \to b} = (X_t^{a \to b})_{0 \leq t \leq T}$, a \emph{Brownian bridge from $a$ to $b$ on $[0,T]$}, by
\begin{align}
X_t^{a \to b}
	&=	a + \frac{t}{T}(b-a) + W_t - \frac{t}{T} W_T , &
0
	&\leq t \leq T . \label{eq:X.ab}
\end{align}
\end{definition}

\noindent
Clearly, we have $X_0^{a \to b} = a$ and $X_T^{a \to b} = b$.

\begin{theorem}
The Brownian Bridge from $a$ to $b$, defined in \eqref{eq:X.ab}, is a Gaussian process.  The mean vector and covariance matrix have entries given by
\begin{align}
m^{a \to b}(t)
	&=	a + \frac{t}{T}(b-a) , \\
c^{a \to b}(t,s)
	&=	t \wedge s - \frac{s t}{T} .
\end{align}
\end{theorem}

\begin{proof}
Since the sum of two (possibly correlated) normal random variables is again normal, and
since for every $t \in [0,T]$, the value of $X_t^{a \to b}$ is a linear combination of $W_t$ and $W_T$, both of which are normal,
it follows that $(X_{t_1}^{a \to b},X_{t_2}^{a \to b},\ldots,X_{t_n}^{a \to b})$ is jointly normal.
To compute the mean and variance, we note that
\begin{align}
m^{a \to b}(t)
	&=	\Eb X_t^{a \to b}
	=		a + \frac{t}{T}(b-a) + \Eb W_t - \frac{t}{T} \Eb W_T
	=		a + \frac{t}{T}(b-a) , \\
c^{a \to b}(t,s)
	&=	\Eb \( X_t^{a \to b} - m^{a \to b}(t) \) \( X_s^{a \to b} - m^{a \to b}(s) \) \\
	&=	\Eb \( W_t - \frac{t}{T} W_T \) \( W_s - \frac{s}{T} W_T \) \\
	&=	\Eb W_t W_s - \frac{t}{T} \Eb W_s W_T - \frac{s}{T} \Eb W_t W_T + \frac{ts}{T^2} \Eb W_T^2 \\
	&=	t \wedge s - \frac{st}{T} - \frac{st}{T} + \frac{st}{T} \\
	&=	t \wedge s - \frac{st}{T} .
\end{align}
\end{proof}

\noindent
Suppose $\Fb=(\Fc_t)_{t \geq 0}$ is a filtration for a Brownian motion $W$.  Then $X^{a \to b}$ is clearly \emph{not} adapted to $\Fb$.  Indeed, since $X_t^{a \to b}$ is expressed in terms of $W_T$, we require the information in $\Fc_T$ in order to write the value of $X_t^{a \to b}$.  There is an alternative definition of a Brownian bridge, which \emph{is} adapted to $\Fb$.

\begin{definition}[Brownian bridge, version II]
\label{def:bridge.bm.2}
Let $W=(W_t)_{t \geq 0}$ be a Brownian motion and fix $T>0$.  We define $Y^{a \to b} = (Y_t^{a \to b})_{0 \leq t \leq T}$, a \emph{Brownian bridge from $a$ to $b$ on $[0,T]$}, by
\begin{align}
Y_t^{a \to b}
	&=	a + \frac{t}{T}(b-a) + (T-t) \int_0^t \frac{1}{T-s} \dd W_s , &
0
	&\leq t \leq T . \label{eq:dY.ab}
\end{align}
\end{definition}

\begin{theorem}
The process $Y^{a \to b} = (Y_t^{a \to b})_{0 \leq t \leq T}$ given in Definition \ref{def:bridge.bm.2} is a Gaussian process.  The mean vector and covariance matrix have entries given by
\begin{align}
m^{a \to b}(t)
	&=	a + \frac{t}{T}(b-a) , \\
c^{a \to b}(t,s)
	&=	t \wedge s - \frac{s t}{T} ,
\end{align}
just as $X^{a \to b} = (X_t^{a \to b})_{0 \leq t \leq T}$.
\end{theorem}

\begin{proof}
The stochastic integral $I_t:=\int_0^t \frac{1}{T-s} \dd W_s$ in \eqref{eq:dY.ab} has a deterministic integrand.  Therefore, $I_t$ is Gaussian, as is $(T-t)I_t$.  It follows that the process $Y$ is a Gaussian process.  Next, we compute
\begin{align}
m^{a \to b}(t)
	&=	\Eb Y_t^{a \to b} 
	=		a + \frac{t}{T}(b-a) + (T-t)  \Eb \int_0^t \frac{1}{T-s} \dd W_s
	=		a + \frac{t}{T}(b-a) .
\end{align}
Finally, to compute $c(t,s)$, assume without loss of generality that $s \leq t$.  Then we have
\begin{align}
c^{a \to b}(t,s)
	&=	\Eb \( Y_t^{a \to b}- m(t) \) \( Y_s^{a \to b}- m(s) \) \\
	&=	(T-t)(T-s) \Eb \int_0^t \frac{1}{T-u} \dd W_u \cdot \int_0^s \frac{1}{T-r} \dd W_r \\
	&=	(T-t)(T-s) \Eb \Big( \int_0^s \frac{1}{T-u} \dd W_u \Big)^2
			+ (T-t)(T-s) \Eb \int_0^s \frac{1}{T-u} \dd W_u \Eb \Big[ \int_s^t \frac{1}{T-r} \dd W_r \Big | \Fc_s \Big] \\
	&=	(T-t)(T-s) \int_0^s \Big( \frac{1}{T-u} \Big)^2 \dd u  \\
	&=	(T-t)(T-s) \( \frac{1}{T-s} - \frac{1}{T} \) 
	=		s - \frac{st}{T} .
	=		s \wedge t - \frac{st}{T} .
\end{align}
This completes the proof.
\end{proof}

\noindent
It is interesting to compute the differential of $Y_t^{a \to b}$.  For simplicity, consider $Y^{0 \to 0}$.  We have
\begin{align}
\dd Y_t^{0 \to 0}
	&=	 \Big( \int_0^t \frac{1}{T-s} \dd W_s \Big) \dd(T- t) + (T-t) \dd \int_0^t \frac{1}{T-s} \dd W_s \\
	&=	- \Big( \frac{1}{T-t} Y_t^{0 \to 0} \Big) \dd t + (T-t) \frac{1}{T-t} \dd W_t \\
	&=	- \frac{Y_t^{0 \to 0}}{T-t} \dd t + \dd W_t .
\end{align}
Note that $Y^{0 \to 0}$ mean reverts to zero with an instantaneous rate $\frac{|Y_t^{0 \to 0}|}{T-t}$ that blows up as $t \to T$.  This blow up guarantees that $Y_T^{0 \to 0} = 0$.

\begin{definition}[Brownian bridge, version III]
\label{def:bridge.bm.3}
Let $W=(W_t)_{t \geq 0}$ be a Brownian motion and fix $T>0$.  We define $Z^{a \to b} = (Z_t^{a \to b})_{0 \leq t \leq T}$, a \emph{Brownian bridge from $a$ to $b$ on $[0,T]$}, by
\begin{align}
Z_t^{a \to b}
	&=	a + W_t | W_T = b -a , &
0
	&\leq t \leq T . \label{eq:Z.ab}
\end{align}
\end{definition}

\begin{theorem}
\label{thm:Z.bridge}
The process $Z^{a \to b} = (Z_t^{a \to b})_{0 \leq t \leq T}$ given in Definition \ref{def:bridge.bm.3} is a Gaussian process.  The mean vector and covariance matrix have entries given by
\begin{align}
m^{a \to b}(t)
	&=	a + \frac{t}{T}(b-a) , \\
c^{a \to b}(t,s)
	&=	t \wedge s - \frac{s t}{T} ,
\end{align}
just as $X^{a \to b} = (X_t^{a \to b})_{0 \leq t \leq T}$.
\end{theorem}

\begin{proof}
We will compute the density of $Z_t^{a \to b}$ directly.  For simplicity, let $a=0$.  We have
\begin{align}
\Pb( Z_t^{0 \to b} \in \dd x )
	&=	\Pb( W_t \in \dd x | W_T = b)
	=		\frac{\Pb( W_t \in \dd x , W_T \in \dd b)}{\Pb(W_T \in \dd b)} \\
	&=		\frac{\Pb( W_t \in \dd x ) \Pb(W_T - W_t + x \in \dd b)}{\Pb(W_T \in \dd b)} .		
\end{align}
Letting $f_{\mu,\sig^2}$ denote the density of a $\Nc(\mu,\sig^2)$ random variable we have
\begin{align}
\Pb( Z_t^{0 \to b} \in \dd x )
	&=	\frac{f_{0,t}(x)\dd x \, f_{x,T-t}(b) \dd b}{f_{0,T}(b) \dd b} 
	 =	f_{m^{0 \to b}(t),c^{0 \to b}(t,T)}(x) \dd x
\end{align}
where $m^{a \to b}(t)$ and $c^{a \to b}(t,s)$ are as given in Theorem \ref{thm:Z.bridge}.  We can clearly see then, that $Z_t^{0 \to b}$ is normally distributed at every $t$ with mean $m^{0 \to b}(t)$.  We leave computation of the entries of the covariance matrix $c^{0 \to b}(t,s)$ as an exercise for the reader.
\end{proof}

%
%

%
%

%


\section{Girsanov's Theorem for a single Brownian motion}
We briefly recall some results from Section \ref{sec:measure.change}.
Suppose that, on a probability space $(\Om,\Fc,\Pb)$ we have a random variable $Z \geq 0$ that has expectation $\Eb Z = 1$.
Then we can define a new probability measure $\Pbt$ via
\begin{align}
\Pbt(A)
	&=	\Eb Z \Ib_A , &
A
	&\in \Fc ,
\end{align}
and we call $Z=\frac{\dd \Pbt}{\dd \Pb}$ the Radon-Nikod\'ym derivative of $\Pbt$ with respect to $\Pb$.  If $Z$ is strictly positive $Z>0$, then we also have
\begin{align}
\Pb(A)
	&=	\Ebt \frac{1}{Z} \Ib_A , &
A
	&\in \Fc ,
\end{align}
and we call $\frac{1}{Z}=\frac{\dd \Pb}{\dd \Pbt}$ the Radon-Nikod\'ym derivative of $\Pb$ with respect to $\Pbt$.
\\[0.5em]
In Example \ref{ex:normal.meas.change}, on a probability space $(\Om,\Fc,\Pb)$, we defined $X \sim \Nc(0,1)$ and a Radon-Nikod\'ym derivative $Z = \ee^{- \theta X - \tfrac{1}{2} \theta^2}$.  We showed that $Y := X + \theta$ was $\Nc(\theta,1)$ under $\Pb$ and $\Nc(0,1)$ under $\Pbt$.  Thus, $Z$ had the effect of changing the mean of $Y$.
\\[0.5em]
We would like to extend this idea from a static to a dynamics setting.  Specifically, we would like to find a measure change that modifies the dynamics of a stochastic process $X = (X_t)_{t \geq 0}$.

\begin{definition}
\label{def:rndp}
Let $(\Om,\Fc,\Pb)$ be a probability space and let $\Fb=(\Fc_t)_{0 \leq t \leq T}$ be a filtration on this space.  A \emph{Radon-Nikod\'ym derivative process} $(Z_t)_{0 \leq t \leq T}$ is any process of the form
\begin{align}
Z_t
	&:=	\Eb[ Z | \Fc_t ]
\end{align}
where $Z$ is an $\Fc_T$-measurable random variable satisfying $\Eb Z = 1$ and $Z > 0$.
\end{definition}

\noindent
Note that $Z$ in Definition \ref{def:rndp} satisfies the conditions of a Radon-Nikod\'ym derivative.  As such, one can define a measure change $\frac{\dd \Pbt}{\dd \Pb}
$ from $Z$.

\begin{lemma}
A \emph{Radon-Nikod\'ym derivative process} $Z = (Z_t)_{0 \leq t \leq T}$ is a martingale.
\end{lemma}

\begin{proof}
For $0 \leq s \leq t \leq T$ we have
\begin{align}
\Eb[ Z_t | \Fc_s ]
	&=	\Eb[ \Eb[ Z | \Fc_t ] | \Fc_s ]
	=		\Eb[  Z | \Fc_s ]
	=		Z_s .
\end{align}
\end{proof}

\begin{lemma}
\label{lem:1}
Let $(Z_t)_{0 \leq t \leq T}$ be a Radon-Nikod\'ym derivative process and define $\frac{\dd \Pbt}{\dd \Pb} = Z$.
Suppose $Y \in \Fc_s$, where $s \in [0,T]$.  Then
\begin{align}
\Ebt Y
	&=	\Eb Z_s Y . 
\end{align}
\end{lemma}

\begin{proof}
The proof is a simple exercise in iterated conditioning.  We have
\begin{align}
\Ebt Y
	&=	\Eb Z Y
	=		\Eb Y \Eb[ Z | \Fc_s ]
	=		\Eb Y Z_s .
\end{align}
\end{proof}

\begin{lemma}
\label{lem:2}
Let $(Z_t)_{0 \leq t \leq T}$ be a Radon-Nikod\'ym derivative process and define $\frac{\dd \Pbt}{\dd \Pb} = Z$.
Suppose $Y \in \Fc_t$ where $0 \leq s \leq t \leq T$.  Then
\begin{align}
\Ebt [ Y | \Fc_s ]
	&=	\frac{1}{Z_s} \Eb [ Z_t Y | \Fc_s ] .
\end{align}
\end{lemma}

\begin{proof}
From Definition \ref{def:cond.E}, we recall that a conditional expectation $\Ebt [ Y | \Fc_s ]$ must satisfy two properties:\\
(i) $\Ebt[ Y | \Fc_s] \in \Fc_s $.\\
(ii) $\Ebt \Ib_A \Ebt[Y | \Fc_s ] = \Ebt \Ib_A Y$ for all $A \in \Fc_s$.\\
It is clear that property (i) is satisfied since
\begin{align}
\Ebt [ Y | \Fc_s ]
	&=	\frac{1}{Z_s} \Eb [ Z_t Y | \Fc_s ]	
	=		\frac{1}{\Eb[Z|\Fc_s]} \Eb [ Z_t Y | \Fc_s ]	
	\in \Fc_s .
\end{align} 
We must check property (ii).  Let $A \in \Fc_s$.  Then
\begin{align}
\Ebt \Ib_A \Ebt[Y|\Fc_s] 
	&=	\Ebt \Ib_A \frac{1}{Z_s} \Eb [ Z_t Y | \Fc_s ] \\
	&=	\Eb Z_s \Ib_A \frac{1}{Z_s} \Eb [ Z_t Y | \Fc_s ]  &
	&		\text{(because $\Ib_A \frac{1}{Z_s} \Eb [ Z_t Y | \Fc_s ] \in \Fc_s$ and Lemma \ref{lem:1})} \\
	&=	\Eb \Ib_A \Eb [ Z_t Y | \Fc_s ] \\
	&=	\Eb \Eb [ \Ib_A Z_t Y | \Fc_s ] \\
	&=	\Eb  \Ib_A Z_t Y \\
	&=	\Ebt  \Ib_A Y . &
	&		\text{(because $\Ib_A Y \in \Fc_t$ and Lemma \ref{lem:1})}
\end{align}
\end{proof}

\begin{theorem}[Girsanov]
\label{thm:girsanov}
Let $W = (W_t)_{0 \leq t \leq T}$ be a Brownian motion on a probability space $(\Om,\Fc,\Pb)$ and let $\Fb=(\Fc_t)_{0 \leq t \leq T}$ be a filtration for $W$.  Suppose $\Theta = (\Theta_t)_{0 \leq t \leq T}$ is adapted to the filtration $\Fb$.  
Define $(Z_t)_{0 \leq t \leq T}$ and $\Wt = (\Wt_t)_{0 \leq t \leq T}$ by
\begin{align}
Z_t
	&=	\exp \Big( - \int_0^t \tfrac{1}{2} \Theta_s^2 \dd s - \int_0^t \Theta_s \dd W_s \Big) , &
\dd \Wt_t
	&=	\Theta_t \dd t + \dd W_t , &
\Wt_0
	&=	0 .
\end{align}
Assume that
\begin{align}
\Eb \int_0^T \Theta_t^2 Z_t^2 \dd t
	&< \infty .
\end{align}
Define a Radon-Nikod\'ym derivative $Z \equiv \frac{\dd \Pbt}{\dd \Pb} := Z_T$.  Then the process $\Wt$ is a Brownian motion under $\Pbt$.
\end{theorem}

\begin{proof}
We will use L\'evy's Theorem \ref{thm:Levy} to show that $\Wt$ is a Brownian motion under $\Pbt$.  By definition, $\Wt_0 = 0$.  Also, we see that 
\begin{align}
\dd[ \Wt, \Wt ]_t 
	&= (\dd \Wt_t)^2
	=		(\dd W_t + \Theta_t \dd t)^2
	=		(\dd W_t)^2 + 2 \Theta_t \dd W_t \dd t + \Theta_t^2 (\dd t)^2 
	=		\dd t .  
\end{align}
Lastly, we must verify that $\Wt$ is a martingale under $\Pbt$.  To do this, we first show that $(Z_t)_{t \geq 0}$ is a martingale under $\Pb$.  Setting
\begin{align}
Z_t
	&=	\ee^{X_t} , &
X_t
	&=	- \int_0^t \tfrac{1}{2} \Theta_s^2 \dd s - \int_0^t \Theta_s \dd W_s ,
\end{align}
we find (using the It\^o formula) that
\begin{align}
\dd Z_t
	&=	\ee^{X_t} \dd X_t + \tfrac{1}{2} \ee^{X_t} \dd [X,X]_t \\
	&=	Z_t ( - \tfrac{1}{2} \Theta_t^2 \dd t - \Theta_t \dd W_t ) + \tfrac{1}{2} Z_t \Theta_t^2 \dd t
	=		- Z_t \Theta_t \dd W_t .
\end{align}
Since It\^o integrals are martingales, it follows that $(Z_t)_{0 \leq t \leq T}$ is a martingale under $\Pb$.  In particular we have
$\Eb Z = \Eb Z_T = Z_0 = 1$.  We also have
\begin{align}
Z_t
	&=	\Eb[ Z_T | \Fc_t] = \Eb[ Z | \Fc_t ] ,
\end{align}
for all $0 \leq t \leq T$, which shows that $(Z_t)_{0 \leq t \leq T}$ is a Radon-Nikod\'ym derivative process.
Next, we show that $(\Wt_t Z_t)_{0 \leq t \leq T}$ is a martingale under $\Pb$.  We have
\begin{align}
\dd (\Wt_t Z_t)
	&=	\Wt_t \dd Z_t + Z_t \dd \Wt_t + \dd[ \Wt, Z]_t \\
	&=	\Wt_t ( - Z_t \Theta_t \dd W_t )  + Z_t ( \Theta_t \dd t + \dd W_t )  +  (\Theta_t \dd t + \dd W_t )(- Z_t \Theta_t \dd W_t) \\
	&=	\Wt_t ( - Z_t \Theta_t \dd W_t )  + Z_t ( \Theta_t \dd t + \dd W_t )  - Z_t \Theta_t \dd t \\
	&=	( - \Wt_t \Theta_t + 1) Z_t \dd W_t .
\end{align}
Again, since It\^o integrals are martingales, it follows that $(\Wt_t Z_t)_{0 \leq t \leq T}$ is a martingale under $\Pb$.  Finally, we can show that $\Wt$ is a martingale under $\Pbt$.  Assuming $0 \leq s \leq t \leq T$ we have
\begin{align}
\Ebt[ \Wt_t | \Fc_s ]
	&=	\frac{1}{Z_s} \Eb[ Z_t \Wt_t | \Fc_s ] &
	&		\text{(by Lemma \ref{lem:2})} \\
	&=	\frac{1}{Z_s} Z_s \Wt_s
	=		\Wt_s .
\end{align}
Thus, $\Wt$ is a martingale, and therefore, a Brownian motion under $\Pbt$.
\end{proof}

\begin{theorem}[Martingale representation]
\label{thm:martingale.representation}
Let $W = (W_t)_{0 \leq t \leq T}$ be a Brownian motion on a probability space $(\Om,\Fc,\Pb)$ and let $\Fb=(\Fc_t)_{0 \leq t \leq T}$ be a filtration generated by $W$ (that is, $\Fc_t = \sig( W_s, 0 \leq s \leq t)$).  Let $M = (M_t)_{0 \leq t \leq T}$ be a martingale with respect to the filtration $\Fb$ satisfying $\Eb M_T^2 < \infty$.  Then there exists a process $\Gam = (\Gam_t)_{0 \leq t \leq T}$ that is adapted to the filtration $\Fb$ and satisfies $\Eb \int_0^T \Gam_t^2 \dd t < \infty$ such that
\begin{align}
M_t
	&=	M_0 + \int_0^t \Gam_s \dd W_s , &
t
	&\in [0,T] .
\end{align}
\end{theorem}

\noindent
The proof of Theorem \ref{thm:martingale.representation} is rather involved, and will not be given here.  The theorem plays a fundamental role in mathematical finance, where a martingale $M$ represents the value of an asset and $\Gam$ represents a trading strategy that replicates the asset's value.  For more precise details, take a course in mathematical finance!

\section{Girsanov's Theorem for $d$-dimensional Brownian motion}

\noindent
We conclude this chapter by stating (without proof) Girsanov's Theorem and the martingale representation theorem for multi-dimensional Brownian motions.

\begin{theorem}[Girsanov]
\label{thm:girsanov-d}
Let $W = (W_t^1,W_t^2,\ldots,W_t^d)_{0 \leq t \leq T}$ be a $d$-dimensional Brownian motion on a probability space $(\Om,\Fc,\Pb)$ and let $\Fb=(\Fc_t)_{0 \leq t \leq T}$ be a filtration for $W$.  Suppose $\Theta = (\Theta_t^1,\Theta_t^2,\ldots,\Theta_t^d)_{0 \leq t \leq T}$ is adapted to the filtration $\Fb$.  
Define $(Z_t)_{0 \leq t \leq T}$ and $\Wt = (\Wt_t^1,\Wt_t^2,\ldots,\Wt_t^d)_{0 \leq t \leq T}$ by
\begin{align}
Z_t
	&=	\exp \Big( - \int_0^t \tfrac{1}{2} \< \Theta_s, \Theta_s \> \dd s - \int_0^t \< \Theta_s,  \dd W_s \> \Big) , &
\dd \Wt_t
	&=	\Theta_t \dd t + \dd W_t , &
\Wt_0
	&=	0 ,
\end{align}
where $\< \cdot , \cdot \>$ denotes a $d$-dimensional Euclidean inner product.  Assume that
\begin{align}
\Eb \int_0^T \< \Theta_t, \Theta_t \> Z_t^2 \dd t
	&< \infty .
\end{align}
Define a Radon-Nikod\'ym derivative $Z \equiv \frac{\dd \Pbt}{\dd \Pb} := Z_T$.  Then the process $\Wt$ is a $d$-dimensional Brownian motion under $\Pbt$.
\end{theorem}

\noindent
It is interesting to note that the components of $\Wt$ in Theorem \ref{thm:girsanov-d} could be co-dependent under $\Pb$ (as $\Theta_j$ could depend on any of $W^1, W^2, \ldots, W^d$).  Nevertheless, under $\Pbt$, the components of $\Wt$ are independent of each other.
 
\begin{theorem}[Martingale representation]
\label{thm:martingale.representation-d}
Let $W = (W_t^1,W_t^2,\ldots,W_t^d)_{0 \leq t \leq T}$ be a $d$-dimensional Brownian motion on a probability space $(\Om,\Fc,\Pb)$ and let $\Fb=(\Fc_t)_{0 \leq t \leq T}$ be a filtration generated by $W$ (that is, $\Fc_t = \sig( W_s, 0 \leq s \leq t)$).  Let $M = (M_t)_{0 \leq t \leq T}$ be a martingale with respect to the filtration $\Fb$ satisfying $\Eb M_T^2 < \infty$.  Then there exists a process $\Gam = (\Gam_t^1, \Gam_t^2,\ldots,\Gam_t^d)_{0 \leq t \leq T}$ that is adapted to the filtration $\Fb$ and satisfies $\Eb \int_0^T \< \Gam_t, \Gam_t \> \dd t < \infty$ such that
\begin{align}
M_t
	&=	M_0 + \int_0^t \< \Gam_s , \dd W_s \> , &
t
	&\in [0,T] .
\end{align}
\end{theorem}


\section{Exercises}
\label{sec:exercises-8}

\begin{exercise}
\label{exer:8.1}
Compute $\dd (W_t^4)$.  Write $W_T^4$ as an integral with respect to $W$ plus an integral with respect to $t$.  Use this representation of $W_T^4$ to show that $\Eb W_T^4 = 3 T^2$.  Compute $\Eb W_T^6$ using the same technique.
\end{exercise}

\begin{exercise}
\label{exer:8.2}
Find an explicit expression for $Y_T$ where
\begin{align}
\dd Y_t
	&=	r \dd t + \alpha Y_t \dd W_t .
\end{align}
Hint: compute $\dd( Y_t Z_t)$ where $Z_t := \exp ( - \alpha W_t + \tfrac{1}{2} \alpha^2 t )$.
\end{exercise}

\begin{exercise}
\label{exer:8.3}
Suppose $X$, $\Del$ and $\Pi$ are given by 
\begin{align}
\dd X_t
	&=	\sig X_t \dd W_t , &
\Del_t
	&=	\frac{ \d f}{\d x}(t,X_t) , &
\Pi_t
	&=	X_t \Del_t
\end{align}
where $f$ is some smooth function.  Show that if $f$ satisfies
\begin{align}
\( \frac{\d}{\d t} + \frac{1}{2}\sig^2 x^2 \frac{\d^2}{\d x^2} \) f(t,x) = 0 ,
\end{align}
for all $(t,x)$, then $\Pi$ is a martingale with respect to a filtration $\Fc_t$ for $W$.
\end{exercise}

\begin{exercise}
\label{exer:8.4}
Suppose $X$ is given by
\begin{align}
\dd X_t
	&=	\mu(t,X_t) \dd t + \sig(t,X_t) \dd W_t .
\end{align}
For any smooth function $f$ define
\begin{align}
M_t^f
	&:=	f(t,X_t) - f(0,X_0) - \int_0^t \Big( \frac{\d}{\d s} + \mu(s,X_s)\frac{\d}{\d x} + \frac{1}{2} \sig^2(s,X_s) \frac{\d^2}{\d x^2} \Big) f(s,X_s) \dd s .
\end{align}
Show that $M^f$ is a martingale with respect to a filtration $\Fc_t$ for $W$.
\end{exercise}

\begin{exercise}
\label{exer:8.5}
Let $X = (X_t)_{0 \leq t \leq T}$ be an OU process on a probability space $(\Omega, \Fc, \Pb)$ 
\begin{align}
\dd X_t
	&=	K(\theta - X_t)\dd t + \sigma \dd W_t .
\end{align}
Where $W=(W_t)_{0 \leq t \leq T}$ is a \em Brownian motion \em  under probability measure $\Pb$. Then we can define a new probability measure $\widetilde{\Pb}$ such that the process $\widetilde{W} =(\widetilde{W}_t )_{0\leq t \leq T}$ is a \em Brownian motion \em under $\widetilde{\Pb}$. Then the OU process $X = (X_t)_{0 \leq t \leq T}$ on the new probability space $(\Omega, \Fc, \widetilde{\Pb})$  will be
\begin{align}
\dd X_t
	&=	K(\theta^* - X_t)\dd t + \sigma \dd \widetilde{W}_t  .
\end{align}
Find the \em Radon-Nikod\'ym \em derivative $\frac{\dd \widetilde{\Pb} }{\dd \Pb }$.
\end{exercise}

\begin{exercise}
\label{exer:8.6}
Note, to prove that a process is a Brownian Bridge, you simply need to show that it is a Gaussian process and has the correct mean $m(t)$ and covariance $c(t,s)$.
\\
(a) Let $X$ be a Brownian bridge from zero to zero on the interval $t=0$ to $t=1$.  
Prove that the process $X_{1-t}, 0 \leq t \leq 1$, is also a Brownian bridge.
\\
(b) Prove that if $W$ is a Brownian motion, then the processes $ (1-t)W_{t/(1-t)}$ and $t W_{(1/t)-1},  0 \leq t \leq 1$,
 are both Brownian bridges.
\\
(c) If we write $\bar{t}$ for $t$ modulo 1, prove that the process $Y_t = X_{\overline{t+s}} - X_s,  0 \leq t \leq 1$ is a Brownian bridge for every fixed $s \in (0,1)$. This is called the cyclic invariance of Brownian bridge.
\end{exercise}

\begin{exercise}
\label{exer:8.7}
Suppose that the dynamics of $X$ are given by
\begin{align}
\dd X_t
	&=	\mu_t \dd t + \sig_t \dd W_t ,
\end{align}
where $W$ is a Brownian motion under $\Pb$.  Find a change of measure $\dd \Pbt / \dd \Pb$ so that
\begin{align}
\dd X_t
	&=	\widetilde{\mu}_t \dd t + \sig_t \dd \Wt_t ,
\end{align}
where $\Wt$ is a Brownian motion under $\Pbt$.  You may assume $\sig_t \neq 0$ for all $t$ and that any process you introduce satisfies the conditions needed for Girsanov's Theorem to apply.
\end{exercise}

\begin{exercise}
\label{exer:8.8}
Consider a market with stochastic interest rates $R = (R_t)_{0 \leq t \leq T}$.  The \textit{Cox-Ingersoll-Ross} (CIR) model for $R$ is
\begin{align}
\dd R_t
	&=	\kappa ( \theta - R_t ) \dd t + \del \sqrt{R_t} \dd W_t ,
\end{align}
where $W$ is a Brownian motion under the real-world probability measure $\Pb$.  Note that, in this model, $R_t \geq 0$ for all $t \in [0,T]$.  Show that, under the following change of measure
\begin{align}
\frac{\dd \Pbt}{ \dd \Pb}
	&=	\exp \Big( -\frac{1}{2} \int_0^T \lam^2 R_t \dd t - \int_0^T \lam \sqrt{R_t } \dd W_t \Big) ,
\end{align}
the dynamics of $R$ have the form
\begin{align}
\dd R_t
	&=	\widetilde{\kappa} ( \widetilde{\theta} - R_t ) \dd t + \del \sqrt{R_t} \dd \Wt_t ,
\end{align}
and find $\widetilde{\kappa}$ and $\widetilde{\theta}$.
\end{exercise}

\begin{exercise}
\label{exer:8.9}
Suppose that $X = (X_t)_{t \geq 0}$ and $Z = (Z_t)_{t \geq 0}$ satisfy
\begin{align}
\dd X_t
	&=	\frac{1}{X_t} \dd t + \alpha X_t \dd W_t , &
X_0
	&=	x > 0 , &
\dd Z_t
	&=	- \alpha Z_t \dd W_t , &
Z_0
	&=	1 ,
\end{align}
where $W = (W_t)_{t \geq 0}$ is a Brownian motion.
\\[0.5em]
(a) Compute $\dd (X_t Z_t)$.
\\[0.5em]
(b) Compute $\dd ( (X_t Z_t)^2 )$.
\\[0.5em]
(c) Find an explicit expression for $Z_t$.
\\[0.5em]
(d) Find an explicit solution for $X_t$.  \textit{Hint}: the solution to 
\begin{align}
\dd Y_t
	&=	\Big(G_t - D Y_t \Big) \dd t, &
Y_0
	&=	y , &
	&\text{is}&
Y_t
	&=	\ee^{- D t} y + \int_0^t \ee^{- D(t-s)} G_s \dd s . \label{eq:exer.8.9.result}
\end{align}
\end{exercise}

\begin{exercise}
\label{exer:8.10}
Fix a probability space $(\Om, \Fc, \Pb)$ and a filtration $\Fb = (\Fc_t)_{t \in [0,T]}$ for a Brownian motion $W = (W_t)_{t \in [0,T]}$.  Define the process $Z=(Z_t)_{t\in[0,T]}$ by
\begin{align}
Z_t=\exp \Big( - \tfrac{1}{2} \int_0^t \sig_s^2 \dd s - \int_0^t \sig_s \dd W_s \Big),
\end{align}
where $\sig = (\sig_t)_{t \in [0,T]}$ is $\Fb$-adapted and bounded and introduce the Radon-Nikodým derivative
\begin{align}
\frac{\dd \Pbt}{\dd \Pb}
	&:=	Z_T
 \end{align}
Recall that the process $\Wt = (\Wt_t)_{t \in [0,T]}$ defined as follows
\begin{align}
\Wt_t
	&:=	\int_0^t \sig_s \dd s + W_t ,\quad \widetilde{W_0}= 0
\end{align}
is a Brownian motion under $\Pbt$.  Show that
\begin{align}
\Eb( Z_t \log Z_t )
	&=	\Ebt \Big( \tfrac{1}{2} \int_0^t \sig_s^2 \dd s \Big) ,
\end{align}
where $\Ebt$ denotes expectation under $\Pbt$.
\end{exercise}

\begin{exercise}
\label{exer:8.11}
Fix a probability space $(\Om, \Fc, \Pb)$ and let $\Fb = (\Fc_t)_{0 \leq t \leq T}$ be a filtration for a Brownian motion $W = (W_t)_{0 \leq t \leq T}$.  Define
\begin{align}
Z \equiv \frac{\dd \Pbt}{\dd \Pb}
	&:=	\frac{\ee^{a W_T}}{\Eb \ee^{a W_T}} , &
a
	&\in \Rb .
\end{align}
(a) Show that $Z$ satisfies the properties of a Radon-Nikod\'ym derivative (i.e., a valid change of measure).
\\[0.5em]
(b) Compute $\Ebt W_T$, where $\Ebt$ denotes expectation under $\Pbt$.  
\\[0.5em]
(c) Compute $\Ebt W_t$, where $0 \leq t \leq T$.
\\[0.5em]
Note: if $X \sim \Nc(\mu,\sig^2)$ then $\Eb \ee^{b X} = \ee^{\mu b + \sig^2 b^2 / 2}$.
\end{exercise}

\begin{exercise}
\label{exer:8.12}
Fix a probability space $(\Om, \Fc, \Pb)$ and a filtration $\Fb = (\Fc_t)_{0 \leq t \leq T}$.  Suppose $W$ and $B$ are independent Brownian motions.  Define processes $X = (X_t)_{t \geq 0}$, $Y=(Y_t)_{t \geq 0}$, $R=(R_t)_{t \geq 0}$ and $\Theta = (\Theta_t)_{t \geq 0}$ as follows
\begin{align}
\dd X_t
	&=	\dd W_t , &
\dd Y_t
	&=	\dd B_t , &
R_t
	&:=	\sqrt{X_t^2 + Y_t^2} , &
\Theta_t
	&=	\arctan ( Y_t / X_t ) .
\end{align}
(a) Derive the dynamics of $R$, i.e., $\dd R_t = \ldots \dd W_t + \ldots \dd B_t + \ldots \dd t$ where $\ldots$ is written in terms of $(R_t,\Theta_t)$ and not $(X_t,Y_t)$.
\\[0.5em]
(b) Derive the dynamics of $\Theta$ (i.e., $\dd \Theta_t = \ldots$).  Your expression should \textit{not} involve $X_t$ or $Y_t$.
\\[0.5em]
In case you forgot, $\d_z \arctan(z) = 1/(1+z^2)$.
\end{exercise}

\begin{exercise}
\label{exer:8.13}
Fix a probability space $(\Om, \Fc, \Pb)$ and let $\Fb = (\Fc_t)_{0 \leq t \leq T}$ be a filtration for a Brownian motion $W = (W_t)_{0 \leq t \leq T}$.
Consider two processes $S = (S_t)_{0 \leq t \leq T}$ and $M = (M_t)_{0 \leq t \leq T}$ whose dynamics are
\begin{align}
\dd S_t
	&=	\mu S_t \dd t + \sig S_t \dd W_t , &
\dd M_t
	&=	r M_t \dd t ,
\end{align}
where $\mu$, $r$ and $\sig > 0$ are constants.
(In finance, this is the \textit{Black-Scholes} model, in which $S$ is the price of a stock and $M$ is the value of a money market account.)
Let $\Del = (\Del_t)_{0 \leq t \leq T}$ be an $\Fb$-adapted process and let $X = (X_t)_{0 \leq t \leq T}$ satisfy
\begin{align}
\dd X_t
	&=	\Del_t \dd S_t + \frac{X_t - \Del_t S_t}{M_t} \dd M_t .
\end{align}
(In finance, $X$ is the value of a \textit{self-financing portfolio} that holds $\Del_t$ shares of $S$ at time $t$ and invests the rest of its value in $M$.)\\
(a) Find a change of measure $\dd \Pbh/ \dd \Pb$ under which $X/S$ is a $(\Pbh,\Fb)$ martingale for every process $\Del$.\\
(b) What are the dynamics of $M/S$ under $\Pbh$?  That is, find $\Theta_t$ and $\Sig_t$ such that $\dd (M_t/S_t) = \Theta_t \dd t + \Sig_t \dd \Wh_t$ where $\Wh$ is a $(\Pbh,\Fb)$ Brownian motion.
\end{exercise}

%
%

\chapter{SDEs and PDEs}
\label{chapter:SDEPDE}

The notes from this chapter are taken primarily from \cite[Chapter 6]{shreve2004}, \cite[Chapters 5, 7, 8 and 9]{oksendal} and \cite{linetskybook}.
Another good reference is \cite[Chapter 15]{karlin}.


\section{Stochastic differential equations}
\label{sec:sde}

\noindent
We will begin this chapter by focusing on scalar stochastic differential equations driven by a single Brownian motion.  The extension to $d$-dimensional stochastic differential equations driven by $m$ Brownian motions will be given in Section \ref{sec:sde-d}.

\begin{definition}
A \emph{stochastic differential equation} (SDE) is an equation of the form
\begin{align}
\dd X_s
	&=	\mu(s,X_s) \dd s + \sig(s,X_s) \dd W_s , &
X_t
	&=	x , \label{eq:SDE}
\end{align}
We call functions $\mu$ and $\sig$ the \emph{drift} and \emph{diffusion}, respectively, and we call $X_t = x$ the \emph{initial condition}.  A \emph{(strong) solution} of an SDE is a stochastic process $X = (X_s)_{s \geq t}$ such that
\begin{align}
X_T
	&=	x + \int_t^T \mu(s,X_s) \dd s + \int_t^T \sig(s,X_s) \dd W_s , \label{eq:SDE.soln}
\end{align}
for all $T \geq t$.
\end{definition}

\noindent
One way to envision a strong solution of an SDE is as follows: think of a sample path $W_\cdot(\om) : [t,\infty) \to \Rb$ as input.  From this input, we can construct a unique sample path $X_\cdot(\om) : [t,\infty) \to \Rb$.   
\\[0.5em]
Ideally, we would like to write $X_T$ as an \textit{explicit} functional of the Brownian path $(W_s)_{s \geq t}$.  Unfortunately, this is typically not possible.  Still, it will help to build intuition if we see some explicitly solvable examples.

\begin{example}[Geometric Brownian motion]
\label{ex:gbm}
A \emph{geometric Brownian motion} is a process $Z = (Z)_{t \geq 0}$ that satisfies
\begin{align}
\dd Z_t
	&=	\mu(t) Z_t \dd t + \sig(t) Z_t \dd W_t , &
Z_0
	&=	z ,
\end{align}
where $\mu$ and $\sig$ are \emph{deterministic} functions of $t$.
To solve this SDE, we consider $X_t = \log Z_t$.  Using the It\^o formula, we obtain
\begin{align}
\dd X_t
	&=	\dd \log Z_t
	=		\frac{1}{Z_t} \dd Z_t + \frac{1}{2} \( \frac{-1}{Z_t^2} \) \dd[Z,Z]_t \\
	&=	\( \mu(t) - \tfrac{1}{2} \sig^2(t) \) \dd t + \sig(t) \dd W_t .
\end{align}
Integrating from $0$ to $T$, we obtain
\begin{align}
X_T
	&=	x + \int_0^T \big( \mu(t) - \tfrac{1}{2} \sig^2(t) \big) \dd t + \int_0^T \sig(t)  \dd W_t , &
x
	&= \log z .
\end{align}
Finally, we obtain our expression for $Z_T$.
\begin{align}
Z_T
	&=	\exp \( X_T \)
	=		z \exp \Big( \int_0^T \big( \mu(t) - \tfrac{1}{2} \sig^2(t) \big) \dd t + \int_0^T \sig(t)  \dd W_t \Big) .
\end{align}
\end{example}

\begin{example}[Linear SDE]
Consider the following SDE with linear coefficients
\begin{align}
\dd X_t
	&=	\( b(t) + \mu(t) X_t \) \dd t + \( a(t) + \sig(t) X_t \) \dd W_t , &
X_0
	&=	x . \label{eq:linear.SDE}
\end{align}
Define two processes $Z = (Z_t)_{t \geq 0}$ and $Y=(Y_t)_{t \geq 0}$ by
\begin{align}
\dd Z_t
	&=	\mu(t) Z_t \dd t + \sig(t) Z_t \dd W_t , &
Z_0
	&=	1 , \\
\dd Y_t
	&=	\frac{b(t) - \sig(t) a(t) }{Z_t} \dd t + \frac{a(t)}{Z_t} \dd W_t , &
Y_0
	&=	x .
\end{align}
We will show that $X_t = Y_t Z_t$ solves \eqref{eq:linear.SDE}.
Using the product rule, we have
\begin{align}
\dd X_t
	&=	\dd ( Y_t Z_t ) \\
	&=	Y_t \dd Z_t + Z_t \dd Y_t + \dd [Y,Z]_t \\
	&=	\Big( \mu(t) Y_t Z_t \dd t + \sig(t) Y_t Z_t \dd W_t \Big) 
			+ \Big( \( b(t) - \sig(t) a(t) \) \dd t + a(t) \dd W_t \Big) 
			+ a(t) \sig(t) \dd t \\
	&=	\( b(t) + \mu(t) Y_t Z_t \) \dd t + \( a(t) + \sig(t) Y_t Z_t \) \dd W_t \\
	&=	\( b(t) + \mu(t) X_t \) \dd t + \( a(t) + \sig(t) X_t \) \dd W_t .
\end{align}
We also have $X_0 = Y_0 Z_0 = x \cdot 1 = x$.  Thus, we have shown that $YZ$ solves \eqref{eq:linear.SDE}.
Now, we note that
\begin{align}
Y_T
	&=	x + \int_0^T \frac{b(t) - \sig(t) a(t) }{Z_t} \dd t + \int_0^T \frac{a(t)}{Z_t} \dd W_t ,
\end{align}
and thus
\begin{align}
X_T
	&=	Y_T Z_T 
	=	x Z_T + Z_T \int_0^T \frac{b(t) - \sig(t) a(t) }{Z_t} \dd t + Z_T \int_0^T \frac{a(t)}{Z_t} \dd W_t .
\end{align}
where, from Example \ref{ex:gbm}, we have
\begin{align}
Z_t
	&=	\exp \Big( \int_0^t \big( \mu(s) - \tfrac{1}{2} \sig^2(s) \big) \dd s + \int_0^t \sig(s)  \dd W_s \Big) .
\end{align}
\end{example}

\noindent
Although we have mostly thrown mathematical rigor out the window in these notes, in an effort to be responsible mathematicians, we should at least state (even if we do not prove) an existence and uniqueness result.

\begin{theorem}[Existence and Uniqueness of SDEs]
\label{thm:exist.unique}
Consider the following SDE
\begin{align}
\dd X_t
	&=	\mu(t,X_t) \dd t + \sig(t,X_t) \dd W_t . \label{eq:dX.exist.unique}
\end{align}
Suppose the functions $\mu: \Rb_+ \times \Rb \to \Rb$ and $\sig : \Rb_+ \times \Rb \to \Rb$ satisfy
\begin{align}
\text{At most linear growth}:&&
|\mu(t,x)| + |\sig(t,x)|
	&\leq C_1 ( 1 + |x| ) , &
	&\forall \, t, x, \label{eq:linear.growth} \\
\text{Lipschitz continuity}:&&
| \mu(t,x) - \mu(t,y) | + | \sig(t,x) - \sig(t,y) |
	&\leq C_2 | x - y | , &
	&\forall \, t, x, y , \label{eq:lipschitz}
\end{align}
for some constants $C_1,C_2 < \infty$.  Then, for any initial condition $X_0 = x \in \Rb$, SDE \eqref{eq:dX.exist.unique} has a unique solution, which is adapted to the filtration $\Fb=(\Fc_t)_{t \geq 0}$ generated by $W = (W_t)_{t \geq 0}$ and satisfies $\Eb \int_0^T X_t^2 \dd t < \infty$ for all $T< \infty$.
\end{theorem}

\begin{remark}
Theorem \ref{thm:exist.unique} actually refers to a \emph{strong} solution of an SDE.  There is another notion of a solution of an SDE called a \emph{weak solution}.  We will not discuss weak solutions here.
\end{remark}

\noindent
We will not prove Theorem \ref{thm:exist.unique}.  We refer the interested reader, instead, to \cite[Theorem 5.21]{oksendal}.  However, we will illustrate with two examples what can go wrong if equations \eqref{eq:linear.growth} and \eqref{eq:lipschitz} are not satisfied.

\begin{example}
Let $X$ satisfy
\begin{align}
\dd X_t
	&=	X_t^2 \dd t , &
X_0
	&=	1 . \label{eq:X.blow.up}
\end{align}
We identify $\mu(t,x) = x^2$, which does \emph{not} satisfy the linear growth condition \eqref{eq:linear.growth}.  The unique solution to \eqref{eq:X.blow.up} is
\begin{align}
X_t
	&=	\frac{1}{1-t} , &
0
	&\leq t < 1 .
\end{align}
Note that $X_t$ blows up as $t \to 1$ and that $X_t$ is not defined for $t \geq 1$.  Thus, it is impossible to find a solution that is defined for all $t \geq 0$.
\end{example}

\begin{example}
Let $X$ satisfy
\begin{align}
\dd X_t
	&=	3 X_t^{2/3} \dd t , &
X_0
	&=	0 . \label{eq:X.not.unique}
\end{align}
We identify $\mu(t,x) = 3 x^{2/3}$, which does \emph{not} satisfy the Lipschitz condition \eqref{eq:lipschitz} at $x=0$.  One can check directly that any $X^{(a)}$ of the form
\begin{align}
X_t^{(a)}
	&=	\left\{ \begin{array}{ll}
			0 & t \leq a \\
			(t-a)^3 & t > a ,
			\end{array} \right.
\end{align}
satisfies \eqref{eq:X.not.unique}.  Thus, solutions of \eqref{eq:X.not.unique} are not unique.
\end{example}

\begin{theorem}[Markov property of solutions of an SDE]
\label{thm:Markov.prop}
Let $X = (X_t)_{t \geq 0}$ be the unique solution of an SDE of the form \eqref{eq:SDE} whose coefficients satisfy the conditions of Theorem \ref{thm:exist.unique}.  Then $X$ is a Markov process.  That is, for $t \leq T$ and for any bounded measurable function $\varphi$, there exists a function $g$ (which depends on $t$, $T$ and $\varphi$) such that
\begin{align}
\Eb[ \varphi( X_T ) | \Fc_t ]
	&=	g(X_t) ,
\end{align}
where $\Fb=(\Fc_t)_{t \geq 0}$ is a filtration for the Brownian motion $W$ (see Definition \ref{def:bm.filtration}).
\end{theorem}

\noindent
The proof of Theorem \ref{thm:Markov.prop} is somewhat technical and will not be given here.  But, the intuitive idea for why the theorem is true is rather simple.
From \eqref{eq:SDE.soln}, we see that the value of $X_T$  depends only on the path of the Brownian motion over the interval $[t,T]$ and the initial value $X_t = x$.
The path that $X$ took to arrive at $X_t = x$ plays no role.  In other words, \emph{given the present $X_t = x$, the future $(X_T)_{T>t}$ is independent of the past $\Fc_t$}.  With this in mind, the process $X$ should admit a transition density
\begin{align}
\Pb( X_T \in \dd y | X_t = x)
	&=	\Gam(t,x;T,y) \dd y ,
\end{align}
and thus, the function $g$ should be given by
\begin{align}
g(X_t)
	&=	\Eb [ \varphi( X_T ) | \Fc_t ]
	=		\Eb [ \varphi( X_T ) | X_t ]
	=		\int \dd y \, \Gam(t,X_t;T,y) \varphi( y ) .
\end{align}
Of course, finding an explicit representation of the transition density $\Gam$ may not be possible.


\section{Relation between SDEs and Langevin Equations}
\label{sec:langevin}
Suppose we \textit{formally} divide both sides of SDE \eqref{eq:dX.exist.unique} by $\dd t$.  Then we obtain
\begin{align}
\frac{\dd X_t}{\dd t}
	&=	\mu(t,X_t) + \sig(t,X_t) \Wd_t , &
\Wd_t
	&:= \frac{\dd W_t}{\dd t} , \label{eq:langevin}
\end{align}
where we have introduced \textit{white noise} $\Wd$, the time derivative for Brownian motion.  Ignore, for the moment, that paths of Brownian motion are not differentiable and, thus, that $\dd W_t / \dd t$ is therefore not well-defined.  Equations of the form \eqref{eq:langevin} are known as \textit{Langevin equations}, in honor of the French physicist Paul Langevin.  You can view a Langevin equation as an ODE with a random source term.
\\[0.5em]
As physicists often prefer to describe random systems using Langevin equations rather than SDEs, it is worth taking a moment to learn the mathematical formalism of Langevin equations.
Introducing the physicists' notation $\< \, \cdot \, \> := \Eb \, \cdot \,$, we note that
\begin{align}
\< \Wd_t \>
	&=	\frac{\dd}{\dd t} \Eb W_t 
	= 	0 . \label{eq:noise-mean}
\end{align}
Next, denoting by $\Theta$ the \textit{Heaviside Step Function}, we have
\begin{align}
\< \Wd_t \Wd_T \>
	&=	\frac{\d^2}{\d t \d T} \Eb W_t W_T
	=		\frac{\d^2}{\d t \d T} ( t \wedge T ) \\
	&=	\frac{\d^2}{\d t \d T} ( t \Theta(T-t) + T \Theta(t-T) ) \\
	&=	\frac{\d}{\d t} ( t \del(T-t) + \Theta(t-T) - T \del(t-T) ) \\
	&=	\frac{\d}{\d t} \Theta(t-T) 
	=		\del(t-T) , \label{eq:noise-covariance}
\end{align}
where we have used the facts that $t \del(T-t) = T \del(t-T)$ and $\Theta' = \del$.  We can use \eqref{eq:noise-mean} and \eqref{eq:noise-covariance} to recapture properties we have already established about It\^o integrals.
\\[0.5em]
For example, suppose that
\begin{align}
\frac{\dd I_t }{\dd t}
	&=	\sig(t) \Wd_t , &
I_0
	&=	0 , &
	&\text{so that}&
I_T
	&=	\int_0^T \sig(t) \Wd_t \dd t ,
\end{align}
where $\sig$ is a deterministic (i.e., non-random) function of time.  Then we have
\begin{align}
\< I_T \>
	=		\int_0^T \sig(t)\< \Wd_t \> \dd t
	=		0 .
\end{align}
Thus, we recover the fact that It\^o integrals have expectation zero.  We can also recover the It\^o isometry
\begin{align}
\< I_T^2 \>
	&=	\int_0^T \int_0^T \sig(t) \sig(s) \< \Wd_t \Wd_s \> \dd s \dd t \\
	&=	\int_0^T \int_0^T \sig(t) \sig(s) \del(t-s) \dd s \dd t \\
	&=	\int_0^T \sig^2(t) \dd t .
\end{align}
The above results can be extenuated to the case where the diffusion coefficient $\sig$ is a function of $t$ and $X_t$, but this is a bit more involved.


\section{Connection to partial differential equations}
\label{sec:pde}

In Chapter \ref{chapter:ctmc}, we showed that, when $X = (X_t)_{t \geq 0}$ is a continuous time Markov chain, its semigroup $\Pv = (\Pv_t)_{t \geq 0}$ satisfies coupled ODEs (see Theorem \ref{thm:kfe.kbe}).  Similarly, when $X$ is the solution of an SDE, conditional expectations of the path of $X$ can be related to certain PDEs.  We discuss these relations below.

\begin{theorem}[Kolmogorov Backward equation]
\label{thm:kbe}
Let $X$ be the solution of SDE \eqref{eq:SDE}.  For some suitable function $\varphi$, define
\begin{align}
u(t,X_t)
	&:=	\Eb[ \varphi(X_T) | \Fc_t ] . \label{eq:u.def}
\end{align}
If the function $u \in C^{1,2}$, then it satisfies the \emph{Kolmogorov Backward Equation} (KBE), a linear PDE of the form
\begin{align}
( \d_t + \Ac(t) ) u(t,\cdot)
	&=	0 , &
u(T,\cdot)
	&=	\varphi , &
\Ac(t)
	&=	\mu(t,x) \d_x + \tfrac{1}{2} \sig^2(t,x) \d_x^2 . \label{eq:u.pde}
\end{align}
\end{theorem}

\begin{proof}
First, we note that the process $(u(t,X_t))_{0 \leq t \leq T}$ is a martingale since, for any $0 \leq s \leq t \leq T$, we have
\begin{align}
\Eb[ u(t,X_t) | \Fc_s ]
	&=	\Eb[\Eb[ \varphi(X_T) | \Fc_t ] | \Fc_s ] 
	=		\Eb[ \varphi(X_T) | \Fc_s ]
	=		u(s,X_s) .
\end{align}
Next we take the differential of $u(t,X_t)$ and find using the It\^o formula that
\begin{align}
\dd u(t,X_t)
	&=	\d_t u(t,X_t) \dd t + \d_x u(t,X_t) \dd X_t + \tfrac{1}{2} \d_x^2 u(t,X_t) \dd[X,X]_t \\
	&=	\( \d_t + \Ac(t) \) u(t,X_t) \dd t + \sig(t,X_t) \d_x u(t,X_t) \dd W_t .
\end{align}
The only way for $(u(t,X_t))_{0 \leq t \leq T}$ to be a martingale is for the $\dd t$ term to be zero for every possible path of $X$.  Thus, the function $u$ must satisfy PDE \eqref{eq:u.pde}.  To see why $u(T,x) = \varphi(x)$ simply use the fact that $\varphi(X_T) \in \Fc_T$ to write
\begin{align}
u(T,X_T)
	&=	\Eb[ \varphi(X_T) | \Fc_T ]
	=		\varphi( X_T ) ,
\end{align}
which must hold for all possible $X_T$.
\end{proof}

\noindent
Theorem \ref{thm:kbe} tells us that the function $u$ defined in \eqref{eq:u.def} satisfies a PDE \eqref{eq:u.pde}.  Alternatively, the \emph{Feynman-Kac formula} says that the solution $u$ of the PDE \eqref{eq:u.pde} has the stochastic representation $u(t,x) = \Eb[ \varphi(x) | X_t = x ]$.
\\[0.5em]
The methods outlined in the proof above can be applied more generally to find PDE representations for more complicated functionals of the path of $X$.  The basic steps are as follows
\begin{enumerate*}
\item Find a martingale process.
\item Take the differential of it.
\item Set the $\dd t$ term to zero.
\end{enumerate*}
We illustrate this method with the following Theorem

\begin{theorem}
Define
\begin{align}
u(t,X_t)
	&:=	\Eb \Big[ \ee^{-A(t,T)} \varphi(X_T) + \int_t^T \dd s \, \ee^{-A(t,s)} g(s,X_s) \Big| \Fc_t \Big] , &
A(t,s)
	&:=	\int_t^s \dd u \, \gam(u,X_u) .
\end{align}
Then $u$ solves the following PDE
\begin{align}
(\d_t + \Ac ) u + g
	&=	0 , &
u(T,\cdot)
	&=	\varphi , &
\Ac
	&=	\tfrac{1}{2} \sig^2(t,x) \d_x^2 + \mu(t,x) \d_x - \gam(t,x) . \label{eq:u.pde.2}
\end{align}
\end{theorem}

\begin{proof}
Note that the function $u$ is \emph{not} a martingale as, for $0 \leq s \leq t \leq T$ we have
\begin{align}
\Eb[ u(t,X_t) | \Fc_s]
	&=	\Eb[ \Eb[ \ee^{-A(t,T)} \varphi(X_T) + \int_t^T \dd r \, \ee^{-A(t,r)} g(r,X_r) | \Fc_t ] | \Fc_s ] \\
	&=	\Eb[ \ee^{-A(t,T)} \varphi(X_T) + \int_t^T \dd r \, \ee^{-A(t,r)} g(r,X_r) | \Fc_s ] \\
	&\neq  \Eb[ \ee^{-A(s,T)} \varphi(X_T) + \int_s^T \dd r \, \ee^{-A(s,r)} g(r,X_r) | \Fc_s ]
	=		u(s,X_s) .
\end{align}
However, the  process $M=(M_t)_{0 \leq t \leq T}$, defined by
\begin{align}
M_t
	&:=	\ee^{-A(0,t)} u(t,X_t) + \int_0^t \dd s \, \ee^{-A(0,s)} g(s,X_s) , \\
	&=	\Eb[ \ee^{-A(0,T)} \varphi(X_T) + \int_0^T \dd s \, \ee^{-A(0,s)} g(s,X_s) | \Fc_t ]
\end{align}
\emph{is} a martingale as, for $0 \leq s \leq t \leq T$ we have
\begin{align}
\Eb[ M_t | \Fc_s]
	&=	\Eb[ \Eb[ \ee^{-A(0,T)} \varphi(X_T) + \int_0^T \dd r \, \ee^{-A(0,r)} g(r,X_r) | \Fc_t ] | \Fc_s ] \\
	&=	\Eb[ \ee^{-A(0,T)} \varphi(X_T) + \int_0^T \dd r \, \ee^{-A(0,r)} g(r,X_r) | \Fc_s ] 
	=	M_s .
\end{align}
Now, we take the differential of $M$.  We compute
\begin{align}
\dd M_t
	&=	u(t,X_t) \dd \ee^{-A(0,t)} +  \ee^{-A(0,t)} \dd u(t,X_t) + \dd \int_0^t \dd r \, \ee^{-A(0,r)} g(r,X_r) \\
	&=	- \ee^{-A(0,t)} \gam(t,X_t) u(t,X_t) \dd t 
			+	\ee^{-A(0,t)} \Big( \d_t + \mu(t,X_t) \d_x + \tfrac{1}{2}\sig^2(t,X_t) \d_x^2 \Big) u(t,X_t) \dd t \\ &\quad
			+	\ee^{-A(0,t)} \sig(t,X) \d_x u(t,X_t) \dd W_t 
			+ \ee^{-A(0,t)} g(t, X_t) \dd t \\
	&=	\ee^{-A(0,t)} \Big( ( \d_t + \Ac ) u(t,X_t) + g(t,X_t) \Big) \dd t + \ee^{-A(0,t)} \sig(t,X) \d_x u(t,X_t) \dd W_t  .
\end{align}
Setting the $\dd t$ term equal to zero, we obtain \eqref{eq:u.pde.2}.  The terminal condition is obtained using
\begin{align}
u(T,X_T)
	&=	\Eb[ \ee^{-A(T,T)} \varphi(X_T) + \int_T^T \dd s \, \ee^{-A(t,s)} g(s,X_s) | \Fc_T ] \\
	&= 	\Eb[ \varphi(X_T) | \Fc_T ] 
	= 	\varphi(X_T) ,
\end{align}
completing the proof.
\end{proof}

\subsection*{Killing a diffusion}

On a probability space $(\Om, \Fc, \Pb)$, consider the following model for a diffusion $X$ and a random time $\tau$
\begin{align}
\dd X_t
	&=	\mu(t,X_t) \dd t + \sig(t,X_t) \dd W_t , &
X_0
	&=	x , \label{eq:X.again} \\
\tau
	&=	\inf\{ t \geq 0 : \int_0^t \gam(s,X_s) \dd s \geq E \} &
E
	&\sim \Ec(1) , &
E
	&\ind X . \label{eq:killing.tau}
\end{align}
where $E$ exponentially distributed and independent of the $X$.  
Here, we assume $\gam : \Rb_+ \times \Rb \to \Rb_+$ so that the integral $\int_0^t \gam(s,X_s) \dd s$ is non-decreasing.
Note that $\tau$ is the first time the integral $\int_0^t \gam(s,X_s) \dd s$ equals or exceeds $E$.  The integral depends on the path of $X$, and therefore is stochastic.  The exponentially distributed random variable $E$ is also random (obviously).  For these reasons, we say that $\tau$ is \emph{doubly stochastic}.  The random time $\tau$ called the \emph{killing time} as it sometimes is used to model the lifetime of the process $X$.
\\[0.5em]
Let $\Fb^X=(\Fc_t^X)_{t \geq 0}$ be the filtration generated by observing the $X$ process.  Note that 
\begin{align}
\Ib_{\{ \tau > t \} } 
	& \notin \Fc_t^X .
\end{align}
In order to keep track of the information obtained by observing $\tau$, we introduce an auxiliary process $D$ as follows
\begin{align}
D_t
	&=	\Ib_{ \{ \tau \leq t \} } .
\end{align}
Note that $D_t = 0$ if $\tau > t$ and $D_t = 1$ if $\tau \leq t$.  As a result
\begin{align}
\Ib_{\{ \tau > t \} } 
	&\in	\Fc_t^D ,
\end{align}
where $\Fb^D=(\Fc_t^D)_{t \geq 0}$ is the filtration generated by observing $D$.  We denote by $\Fb=(\Fc_t)_{t \geq 0}$ the filtration generated by observing both $D$ and $X$; this is sometimes written $\Fb = \Fb^X \vee \Fb^D$, which means $\Fc_t = \Fc_t^X \vee \Fc_t^D := \sig( \Fc_t^X \cup \Fc_t^D)$.

\begin{theorem}
Let $X = (X_t)_{t \geq 0}$ and $\tau$ be as given in \eqref{eq:X.again} and \eqref{eq:killing.tau}.  Then
\begin{align}
\Eb[ \Ib_{\{ \tau > T \}} \varphi(X_T) | \Fc_t ]
	&=	\Ib_{\{\tau > t\}} \Eb[ \ee^{- \int_t^T \gam(s,X_s) \dd s} \varphi(X_T)  | \Fc_t ] , \label{eq:killing.result}
\end{align}
where $T \geq t$.
\end{theorem}

\begin{proof}
Noting that $1 = \Ib_{\{\tau \leq t\}} + \Ib_{\{\tau > t\}}$ and $\Ib_{\{\tau \leq t\}} \Ib_{\{\tau > T\}} = 0$, we have
\begin{align}
\Eb[ \Ib_{\{ \tau > T \}} \varphi(X_T) | \Fc_t ]
	&=	\Eb[ \Ib_{\{\tau > t\}} \Ib_{\{ \tau > T \}} \varphi(X_T) | \Fc_t ] +  \Eb[ \Ib_{\{\tau \leq t\}} \Ib_{\{ \tau > T \}} \varphi(X_T) | \Fc_t ] \\
	&=	\Eb[ \varphi(X_T) \Ib_{\{\tau > t\}}  \Eb[ \Ib_{\{ \tau > T \}}| \Fc_T^X, \Fc_t ]  | \Fc_t ] .  \label{eq:killed.E}
\end{align}
Now, Noting that $\Pb( E > x ) = \ee^{-x}$, we compute
\begin{align}
\Pb( \tau > t | \Fc_t^X )
	&=	\Pb\Big( E > \int_0^t \gam(s,X_s) \dd s \Big| \Fc_t^X \Big) 
	=	\ee^{- \int_0^t \gam(s,X_s) \dd s} .
\end{align}
From this, we find
\begin{align}
\Ib_{\{\tau > t\}} \Eb[ \Ib_{\{ \tau > T \}}| \Fc_T^X, \Fc_t ]
	&=	\Ib_{\{\tau > t\}} \Pb( \tau > T | \Fc_T^X, \tau > t ) \\
	&=	\Ib_{\{\tau > t\}} \Pb \Big( E > \int_0^T \gam(s,X_s) \dd s \Big| E > \int_0^t \gam(s,X_s) \dd s , \Fc_T^X \Big) \\
	&=	\Ib_{\{\tau > t\}} \ee^{- \int_t^T \gam(s,X_s) \dd s} , \label{eq:killed.2.exp}
\end{align}
where, in the last equality, we have used the memorylessness property of the exponential distribution.
Inserting \eqref{eq:killed.2.exp} into \eqref{eq:killed.E} we obtain \eqref{eq:killing.result}.
\end{proof}


\section{Kolmogorov forward and backward Equations}
\label{eq:kbe.kfe}
Let $X = (X_t)_{t \geq 0}$ be the solution of the following SDE
\begin{align}
\dd X_t
	&=	\mu(t,X_t) \dd t + \sig(t,X_t) \dd W_t .
\end{align}
Let $\Gam$ denote the \emph{transition density} of the process $X$.
\begin{align}
\Gam(t,x;T,y) \dd y
	&:= \Pb( X_T \in \dd y | X_t = x ) .
\end{align}
In this Section, we will derive two PDEs which are satisfied by the transition density $\Gam$: the \emph{Kolmogorov Backward Equation} (KBE), which is a PDE in the \emph{backward variables} $(t,x)$, and the \emph{Kolmogorov Forward Equation} (KFE), which is a PDE satisfied by the \emph{forward variables} $(T,y)$.  Physicists and biologists sometimes call the KFE the \emph{Fokker-Planck Equation}.
\\[0.5em]
We have already seen the KFE and KBE in the continuous time Markov chain setting, discussed in Section \ref{sec:ctmc}.  The development of the KFE and KBE for diffusion processes is remarkably similar.  

\begin{definition}
The \emph{two-parameter semigroup} $(\Pc(t,T))_{0 \leq t \leq T < \infty}$, of a Markov diffusion $X$, is defined as
\begin{align}
\Pc(t,T) \varphi(x)
	&=	\Eb[  \varphi(X_T) | X_t = x ]
	=		\int \dd y \, \Gam(t,x;T,y) \varphi(y) .
\end{align}
where $\varphi$ is integrable with respect to the transition density $\Gam$.
\end{definition}

\begin{theorem}
The semigroup satisfies
\begin{align}
\Pc(t,t) 
	&= \Iv , &
\Pc(t,s) \Pc(s,T)
	&=	\Pc(t,T) , &
	&0 \leq t \leq s \leq T .
\end{align}
\end{theorem}

\begin{proof}
To see that $\Pc(t,t) = \Iv$ simply observe that
\begin{align}
\Pc(t,t) \varphi(x)
	&=	\Eb[  \varphi(X_t) | X_t = x ] = \varphi(x) = \Iv \varphi(x) .
\end{align}
To see the \emph{semigroup property} $\Pc(t,s) \Pc(s,T) = \Pc(t,T)$, note that
\begin{align}
\Pc(t,s) \Pc(s,T) \varphi(x)
	&=	\Pc(t,s) \Eb[ \varphi(X_T) | X_s = x ]
	=		\Eb[  \Eb[ \varphi(X_T) | X_s ] | X_t = x ] \\
	&=		\Eb [ \varphi(X_T) | X_t = x ]
	=		\Pc(t,T) \varphi(x) .
\end{align}
The semigroup property can, alternatively, be derived from the Chapman-Kolmogorov equation
\begin{align}
\Gam(t,x;T,y)
	&=	\int \dd z \, \Gam(t,x;s,z) \Gam(s,z;T,y) , &
	& 0 \leq t \leq s \leq T < \infty .
\end{align}
We have
\begin{align}
\Pc(t,s) \Pc(s,T) \varphi(x)
	&=	\Pc(t,s) \int \dd y \, \Gam(s,x;T,y) \varphi(y)  \\
	&=	\int \dd z \, \Gam(t,x;s,z) \int \dd y \, \Gam(s,z;T,y) \varphi(y) \\
	&=	\int \dd y \, \Gam(t,x;T,y) \varphi(y) 
	=		\Pc(t,T) \varphi(x) .
\end{align}
\end{proof}

\noindent
We associate with a semigroup, its generator.

\begin{definition}
The \emph{infinitesimal generator} or simply the \emph{generator}, of a semigroup of operators $(\Pc(t,s))_{0 \leq t \leq T < \infty}$ is defined as
\begin{align}
\Ac(t) \varphi(x)
	&:=	\lim_{s \searrow t} \frac{1}{s-t} \Big( \Pc(t,s) \varphi(x) - \varphi(x) \Big) \\
	&=		\lim_{s \searrow t} \frac{1}{s-t} \Big( \Eb[ \varphi(X_s) | X_t = x ] - \varphi(x) \Big)
\end{align}
where $\varphi$ is any function for which the limit exists.
\end{definition}

\begin{theorem}
If $\varphi \in C_0^2$ (bounded and twice differentiable), then and $X$ is the solution of
\begin{align}
\dd X_t
	&=	\mu(t,X_t) \dd t + \sig(t,X_t) \dd W_t
\end{align}
then the generator $\Ac(t)$ of the semigroup $(\Pc(t,s))_{0 \leq t \leq s < \infty}$ of $X$ is given by
\begin{align}
\Ac(t)
	&=	\mu(t,x) \d_x + \tfrac{1}{2} \sig^2(t,x) \d_x^2 .
\end{align}
\end{theorem}

\begin{proof}
From the It\^o formula we have
\begin{align}
\dd \varphi(X_t)
	&=	\d_x \varphi(X_t) \dd X_t + \tfrac{1}{2} \d_x^2 \varphi(X_t) \dd [X,X]_t \\
	&=	\Big( \mu(t,X_t) \d_x \varphi(X_t) + \tfrac{1}{2} \sig^2(t,X_t) \d_x^2 \Big) \dd t + \sig(t,X_t) \d_x \varphi(X_t) \dd W_t \\
	&=	\Ac(t) \varphi(X_t) \dd t + \sig(t,X_t) \d_x \varphi(X_t) \dd W_t .
\end{align}
Therefore, we have
\begin{align}
\varphi(X_s)
	&=	\varphi(X_t) + \int_t^s \Ac(r) \varphi(X_r) \dd r +  \int_t^s \sig(r,X_r) \varphi(X_r) \dd W_r , \\
\Eb [ \varphi(X_s) | X_t = x ]
	&=	\varphi(x) + \int_t^s \Eb[ \Ac(r) \varphi(X_r) | X_t = x ] \dd r .
\end{align}
Finally,
\begin{align}
\lim_{s \searrow t} \frac{1}{s-t}\Big( \Eb[ \varphi(X_s) | X_t = x ] - \varphi(x) \Big) 
	&=	\lim_{s \searrow t} \frac{1}{s-t} \int_t^s \Eb[ \Ac(r) \varphi(X_r) | X_t = x ] \dd r 
	=	\Ac(t) \varphi(x) .
\end{align}
\end{proof}

\begin{remark}
For all intents and purposes, when $X$ is the solution of an SDE the generator $\Ac(t)$ is simply the operator that acts on $\varphi$ in the $\dd t$ term of $\dd \varphi(X_t)$.  Thus, we can write the It\^o formula more compactly as
\begin{align}
\dd \varphi(X_t)
	&=	\Ac(t) \varphi(X_t) \dd t + \sig(t,X_t) \d_x \varphi(X_t) \dd W_t .
\end{align}
Of course, if $X$ is not the solution of an SDE, we cannot write the differential $\dd \varphi(X_t)$ in this more compact form.
\end{remark}

\noindent
We can now state and prove the KBE and KFE

\begin{theorem}
Let $X = (X_t)_{t \geq 0}$ be the unique solution (assumed to exist) of
\begin{align}
\dd X_t
	&=	\mu(t,X_t) \dd t + \sig(t,X_t) \dd W_t .
\end{align}
Let $\Gam$ denote the transition density of $X$
\begin{align}
\Gam(t,x;T,y) \dd y
	&= \Pb(X_T \in \dd y | X_t = x ) .
\end{align}
Then, seen as a function of the backwards variables $(t,x)$, the transition density $\Gam(\,\cdot\,,\,\cdot\,;T,y)$ satisfies
\begin{align}
\text{KBE}:&&
 (\d_t + \Ac(t) ) \Gam(t,\,\cdot \,;T,y)
	&=	0 , &
\Gam(T,\,\cdot \,,T,y)
	&=	\del_y ,
\end{align}
where $\Ac(t)$ is the infinitesimal generator of $X$
\begin{align}
\Ac(t)
	&=	\mu(t,x) \d_x + \tfrac{1}{2} \sig^2(t,x) \d_x^2 .
\end{align}
Seen as a function of the forward variables $(T,y)$, the transition density $\Gam(t,x;\,\cdot\,,\,\cdot\,)$ satisfies
\begin{align}
\text{KFE}:&&
 (-\d_T + \Ac^*(T) ) \Gam(t,x;T,\,\cdot \,)
	&=	0 , &
\Gam(t,x,t,\,\cdot \,)
	&=	\del_x ,
\end{align}
where $\Ac^*(T)$ is the $L^2(\dd x)$ adjoint of $\Ac(T)$
\begin{align}
\Ac^*(T)
	&=	- \d_y \mu(T,y) + \tfrac{1}{2}  \d_y^2 \sig^2(T,y)  .
\end{align}
\end{theorem}

\begin{proof}
We have already established in Theorem \ref{thm:kbe} that
\begin{align}
u(t,x)
	&:=	\Eb[ \varphi(X_T) | X_t = x ]
	=		\int \dd y \, \Gam(t,x;T,y) \varphi(y) ,
\end{align}
satisfies \eqref{eq:u.pde}.  Thus, we must have
\begin{align}
0
	&=	(\d_t + \Ac(t)) u(t,\,\cdot \,)
	=		\int \dd y \, \varphi(y) (\d_t + \Ac(t)) \Gam(t,\, \cdot \,;T,y)  , \label{eq:must.have.pde} \\
\varphi
	&=	u(T,\,\cdot \,) = \int \dd y \, \varphi(y) \Gam(T,\, \cdot \,;T,y) . \label{eq:must.have.bc}
\end{align}
The above equations must hold for \emph{all} functions $\varphi$.  It follows that
\begin{align}
(\d_t + \Ac(t) ) \Gam(\, \cdot \,,t;T,y)
	&=	0 , &
\Gam(T,\, \cdot \,,T,y)
	&=	\del_y ,
\end{align}
which is precisely the KBE for $\Gam$.
\\[0.5em]
Proving the KFE requires a little more effort.  First, we note that, by definition, the generator $\Ac(t)$ and its $L^2(\dd x)$ adjoint satisfy
\begin{align}
\int \dd y \, f(y) \Ac(t) g(y)
	&=	\int \dd y \, g(y) \Ac^*(t) f(y) ,
\end{align}
where, both $\Ac(t)$ and $\Ac^*(t)$ act on the $y$ variable and $f,g \to 0$ as $y \to \pm \infty$.  Now, observe that
\begin{align}
\int \dd y \, \varphi(y) \Gam(t,x;T,y)
	&=	\Eb[ \varphi(X_T) | X_t = x ] \\
	&=	\varphi(x) + \int_t^T \dd s \, \Eb[ \Ac(s) \varphi(X_s) | X_t = x ] + \Eb[ \int_t^T \sig(s,X_s) \d_x \varphi(X_s) \dd W_s | X_t = x ] \\
	&=	\varphi(x) + \int_t^T \dd s \int \dd y \, \Gam(t,x;s,y) \Ac(s) \varphi(y) \\
	&=	\varphi(x) + \int_t^T \dd s \int \dd y \, \varphi(y) \Ac^*(s) \Gam(t,x;s,y) . 
\end{align}
Now, we differentiate the above expression with respect to $T$ and find
\begin{align}
\int \dd y \, \varphi(y) \d_T \Gam(t,x;T,y)
	&=	\int \dd y \, \varphi(y) \Ac^*(T) \Gam(t,x;T,y) .
\end{align}
We also have
\begin{align}
\int \dd y \, \varphi(y) \Gam(t,x;t,y)
	&=	\Eb[ \varphi(X_t) | X_t = x ] = \varphi(x) .
\end{align}
Again, the above expressions must hold for \emph{all} $\varphi$.  Thus, we see that
\begin{align}
0
	&= (-\d_T + \Ac^*(T) ) \Gam(t,x;T,\,\cdot \,) , &
\Gam(t,x;t,\,\cdot \,)
	&=	\del_x ,
\end{align}
which is exactly the KFE.
\end{proof}

\begin{example}
Let us check that the KBE and KFE are satisfied in the following simple case:
\begin{align}
\dd X_t
	&=	\sig \dd W_t .
\end{align}
Since $(X_T | X_t = x) \sim \Nc(x,\sig^2(T-t))$ we have
\begin{align}
\Gam(t,x;T,y)
	&=	\frac{1}{\sqrt{2\pi \sig^2(T-t)}} \exp \Big( \frac{-(y-x)^2}{2 \sig^2(T-t)} \Big) .
\end{align}
Once can check by direct computation that
\begin{align}
(\d_t + \Ac) \Gam(t,\cdot;T,y)
	&=	0 , &
\lim_{t \nearrow T} \Gam(t,\cdot;T,y)
	&=	\del_y , \\
(-\d_T + \Ac^*) \Gam(t,x;T,\cdot)
	&=	0 , &
\lim_{T \searrow t} \Gam(t,x;T,\cdot)
	&=	\del_x ,
\end{align}
where $\Ac = \tfrac{1}{2} \sig^2 \d_x^2$ and $\Ac^* = \tfrac{1}{2} \sig^2 \d_y^2$.
\end{example}


\section{Poisson equations}
The KFE and KBE are parabolic PDEs (think heat-type equations).  It turns out we can also relate the expected value of a functional of $X$ to an elliptic PDE (think, potential equations).  Consider a one-dimensional time-homogeneous diffusion $X = (X_t)_{t \geq 0}$ that lives on some interval $I = (l,r)$, where $-\infty \leq l < r \leq \infty$.  We suppose $X$ satisfies
\begin{align}
\dd X_t
	&=	\mu(X_t) \dd t + \sig(X_t) \dd W_t . \label{eq:dX.th}
\end{align}
Let us define the hitting time
\begin{align}
\tau
	&=	\inf\{ t \geq 0 : X_t \notin I \} , \label{eq:tau.dI}
\end{align}
where $\d I = \{l\} \cup \{r\}$ denotes the boundary of $I$.  

\begin{theorem}
\label{thm:Poisson}
Let $X = (X_t)_{t \geq 0}$ and $\tau$ be given by \eqref{eq:dX.th} and \eqref{eq:tau.dI}. Assume $\Pb(\tau < \infty) = 1$.  Let $\lam \geq 0$, and let $\varphi : \d I \to \Rb$ and $g: I \to \Rb$ be bounded measurable functions.  For any $t \geq 0$ and $x \in I$, define
\begin{align}
u(x)
	&:=	\Eb[ \ee^{-\lam(\tau-t)}\varphi(X_\tau) + \int_t^\tau \ee^{-\lam(s-t)}g(X_s) \dd s | X_t = x ] , &
t &< \tau . \label{eq:u-of-x}
\end{align}
Then $u(x)$ does not depend on $t$ (by time-homogeneity of $X$).  Moreover, if $u \in C^2(I)$ and $u$ is continuous on $I \cup \d I$, then the function $u$ satisfies
\begin{align}
(\Ac - \lam ) u + g
	&=	0 &
	&\text{in $I$} , \label{eq:Poisson.PDE} \\
u
	&=	\varphi , &
	&\text{on $\d I$} , \label{eq:Poisson.BC}
\end{align}
where $\Ac = \mu(x) \d_x + \tfrac{1}{2} \sig^2(x) \d_x^2$ is the infinitesimal generator of $X$.
\end{theorem}

\begin{proof}
We first verify that $u(x)$ is independent of $t$.  By the strong Markov property and time-homogeneity of $X$, conditioning on $X_t = x$ is equivalent to restarting the process at time $t$ from position $x$, and the future evolution is identical in distribution to starting at time $0$ from position $x$.  More precisely, for any $0 \leq s < t < \tau$, we have
\begin{align}
&\Eb[ \ee^{-\lam(\tau-t)}\varphi(X_\tau) + \int_t^\tau \ee^{-\lam(r-t)}g(X_r) \dd r | X_t = x ] \\
	&=	\Eb[ \ee^{-\lam(\tau-s)}\varphi(X_\tau) + \int_s^\tau \ee^{-\lam(r-s)}g(X_r) \dd r | X_s = x ] .
\end{align}
Thus, $u(x)$ is well-defined independently of $t$.
\\[0.5em]
Now, let $\Fb=(\Fc_t)_{0 \leq t < \tau}$ be the filtration generated by $X$.
The process $M=(M_t)_{0 \leq t < \tau}$, defined by
\begin{align}
M_t
	&:=	\ee^{-\lam t} u(X_t) + \int_0^t \ee^{-\lam s} g(X_s) \dd s 
\end{align}
is a martingale.  To see this, first observe that
\begin{align}
M_t
	&=	\ee^{-\lam t} \Eb[ \ee^{-\lam(\tau-t)}\varphi(X_\tau) + \int_t^\tau \ee^{-\lam(s-t)}g(X_s) \dd s | X_t ] + \int_0^t \ee^{-\lam s} g(X_s) \dd s \\
	&=	\ee^{-\lam t} \Eb[ \ee^{-\lam(\tau-t)}\varphi(X_\tau) + \int_t^\tau \ee^{-\lam(s-t)}g(X_s) \dd s | \Fc_t ] + \int_0^t \ee^{-\lam s} g(X_s) \dd s \\
	&=	\Eb[ \ee^{- \lam \tau} \varphi(X_\tau) + \int_t^\tau \ee^{-\lam s} g(X_s) \dd s | \Fc_t ] + \int_0^t \ee^{-\lam s} g(X_s) \dd s \\
	&=	\Eb[ \ee^{- \lam \tau} \varphi(X_\tau) + \int_0^\tau \ee^{-\lam s} g(X_s) \dd s | \Fc_t ] .
\end{align}
Thus, for $0 \leq s \leq t < \tau$ we have
\begin{align}
\Eb[ M_t | \Fc_s ]
	&=	\Eb[ \Eb[ \ee^{- \lam \tau} \varphi(X_\tau) + \int_0^\tau \ee^{-\lam r} g(X_r) \dd r | \Fc_t] | \Fc_s ] \\
	&=	\Eb[ \ee^{- \lam \tau} \varphi(X_\tau) + \int_0^\tau \ee^{-\lam r} g(X_r) \dd r | \Fc_s ]
	=		M_s ,
\end{align}
which proves that $M$ is a martingale.
\\[0.5em]
Assuming $u \in C^2(I)$, we take the differential of $M$ and obtain
\begin{align}
\dd M_t
	&=	\dd  \( \ee^{-\lam t} u(X_t) \) + \dd \int_0^t \ee^{-\lam s} g(X_s) \dd s \\
	&=	\ee^{-\lam t} \dd u(X_t) + u(X_t) \dd \ee^{-\lam t} + \ee^{-\lam t} g(X_t) \dd t \\
	&=	\ee^{- \lam t} \Big[ \d_x u(X_t) \dd X_t + \tfrac{1}{2} \d_x^2 u(X_t) \dd [X,X]_t \Big] - \lam \ee^{-\lam t} u(X_t) \dd t + \ee^{-\lam t} g(X_t) \dd t \\
	&=	\ee^{- \lam t} \Big[ \mu(X_t) \d_x u(X_t) \dd t + \sig(X_t) \d_x u(X_t) \dd W_t + \tfrac{1}{2} \sig^2(X_t) \d_x^2 u(X_t) \dd t \Big] \notag \\
	&\quad - \lam \ee^{-\lam t} u(X_t) \dd t + \ee^{-\lam t} g(X_t) \dd t \\
	&=	\ee^{- \lam t} \Big( \Ac u(X_t) - \lam u(X_t) + g(X_t) \Big) \dd t + \ee^{- \lam t} \sig(X_t) \d_x u(X_t) \dd W_t .
\end{align}
Since $M$ is a martingale, the $\dd t$ term must be zero for all $X_t \in I$.  Thus, we have
\begin{align}
(\Ac - \lam ) u + g
	&=	0 , &
	&\text{in $I$} , 
\end{align}
in agreement with \eqref{eq:Poisson.PDE}.  
\\[0.5em]
The boundary condition \eqref{eq:Poisson.BC} follows from observing that if $X_t \in \d I$ (i.e., $X_t = l$ or $X_t = r$), then $\tau = t$ since the process has already exited $I$.  Thus,
\begin{align}
u(x)
	&=	\Eb[ \ee^{-\lam(\tau-t)}\varphi(X_\tau) + \int_t^\tau \ee^{-\lam(s-t)}g(X_s) \dd s | X_t = x \in \d I] \\
	&=	\Eb[ \ee^{\lam(t-t)}\varphi(x) + \int_t^t \ee^{-\lam(s-t)}g(X_s) \dd s | X_t = x ]
	=		\varphi(x) .
\end{align}
This completes the proof.
\end{proof}

\begin{remark}
The key insight is that for a \emph{time-homogeneous} diffusion (where $\mu$ and $\sig$ do not depend on $t$), the function $u(x)$ defined in \eqref{eq:u-of-x} is independent of $t$. This is because the future evolution of $X$ starting from position $x$ at time $t$ is identical in distribution to starting from position $x$ at time $0$.  This allows us to work with $u$ as a function of $x$ alone, making the PDE elliptic rather than parabolic.
\end{remark}

\begin{remark}
We can generalize Theorem \ref{thm:Poisson} as follows: if we replace $\lam(\tau-t)$ with $\int_t^\tau \lam(X_s) \dd s$ (and $\lam(s-t)$ with $\int_t^s \lam(X_r) \dd r$) then $u$ solves $(\Ac - \lam(x) ) u(x) + g(x) = 0$.  We leave the proof of this as an exercise for the reader.
\end{remark}

\begin{corollary}
\label{cor:laplace}
The Laplace transform of $\tau$, given by
\begin{align}
u(x)
	&:=	\Eb[ \ee^{- \lam \tau } | X_0 = x] , 
\end{align}
satisfies
\begin{align}
(\Ac - \lam) u
	&=	0 , &
	&\text{in $I$} , &
u
	&=	1 , &
	&\text{on $\d I$}.
\end{align}
\end{corollary}

\begin{proof}
Take $\varphi(x) = 1$, $g(x) = 0$ and $t=0$ in Theorem \ref{thm:Poisson}.
\end{proof}

\begin{example}[First hitting time of Brownian motion]
Let us define a process $X$ and a hitting time $\tau$ by
\begin{align}
\dd X_t
	&=	\dd W_t , &
X_0
	&=	x < r < \infty , &
\tau
	&=	\inf\{ t \geq 0 : X_t \notin (-\infty,r) \} .
\end{align}
According to Corollary \ref{cor:laplace}, the function
\begin{align}
u(x)
	&:= \Eb[ \ee^{- \lam \tau} | X_0 = x ] ,
\end{align}
should satisfy
\begin{align}
(\Ac - \lam) u
	&=	0 , &
	&\text{in $(-\infty,r)$} , &
\Ac
	&=	\tfrac{1}{2} \d_x^2 , &
u(r)
	&=	1 . \label{eq:u.laplace.pde}
\end{align}
Let us check that this is the case.  We showed in Section \ref{sec:hitting.bm} that
\begin{align}
\Eb \ee^{- \lam \tau_m}
	&=	\ee^{- |m| \sqrt{2 \lam}} , &
\tau_m
	&:=	\inf\{ t \geq 0 : W_t = m \} ,
\end{align}
where $W_0 = 0$ by definition.  Since Brownian motion is invariant to translations, it follows that
\begin{align}
u(x)
	&=	\ee^{-(r - x)\sqrt{2 \lam}} . \label{eq:u.laplace}
\end{align}
One can easily see that \eqref{eq:u.laplace} satisfies \eqref{eq:u.laplace.pde}, since
\begin{align}
u(r)
	&=	1 , &
\tfrac{1}{2} \d_x^2 u(x) 
	&= \lam u(x) .
\end{align}
\end{example}


\section{In depth look: scalar time-homogeneous diffusions}
In this section we will specialize to \emph{time-homogeneous} scalar diffusions.  Specifically, we consider a diffusion $X = (X_t)_{t \geq 0}$ that lives on some interval $I$ with endpoints $l$ and $r$, where $-\infty \leq l < r \leq \infty$.  We suppose $X$ satisfies
\begin{align}
\dd X_t
	&=	\mu(X_t) \dd t + \sig(X_t) \dd W_t .
\end{align}
The end-points $l$ and $r$ may or may not be part of the interval $I$.  The generator $\Ac$ of $X$ is given by
\begin{align}
\Ac
	&=	\mu(x) \d_x + \tfrac{1}{2} \sig^2(x) \d_x^2 .
\end{align}
We can alternatively write $\Ac$ as
\begin{align}
\Ac
	&=	\frac{1}{m(x)} \d_x \( \frac{1}{s(x)} \d_x \cdot \) , \label{eq:self.adjoint}
\end{align}
where $s$ and $m$, called the \emph{scale} and \emph{speed} densities, respectively, are given by
\begin{align}
s(x)
	&=	\exp \( - \int \dd x \frac{2 \mu(x) }{\sig^2(x)} \) , &
m(x)
	&=	\frac{2}{\sig^2(x) s(x)} 
	=		\frac{2}{\sig^2(x)} \exp \( \int \dd x \frac{2 \mu(x) }{\sig^2(x)} \) . \label{eq:s.and.m}
\end{align}
The constant of integration is arbitrary.  
\\[0.5em]
Note that, if the constant of integration in the indefinite integral in $m$ can be chose such that
\begin{align}
\int_I \dd y \, m(y)
	&=	1 ,
\end{align}
then $m$ is a time-homogeneous solution of the KFE.  To see this simply observe that
\begin{align}
(-\d_T + \Ac^* )m(y)
	&=	\Ac^* m(y)
	=	\d_y \( - \mu(y) m(y) + \d_y \frac{\sig^2(y)}{2} m(y) \) \\
	&=	\d_y \( - \frac{2 \mu(y)}{\sig^2(y)} \exp\( \int \dd y \, \frac{2 \mu(y)}{\sig^2(y)} \) + \d_y \exp\( \int \dd y \, \frac{2 \mu(y)}{\sig^2(y)}\) \)
	=		0 .
\end{align}
Thus, $m$ is a stationary density for $X$.


\subsection*{Useful results from functional analysis}
Before continuing further, it will be convenient to recall a few definitions from functional analysis.

\begin{definition}
Let $L^2(E,\rho)$ denote the set of square-integrable functions on $E$ weighted by $\rho$
\begin{align}
L^2(E,\rho)
	&:=	\{ f : \< f, f \>_\rho < \infty \} , &
\< f, g \>_\rho
	&:=	\int_E \dd x \, \rho(x) f(x) g(x) .
\end{align}
The \emph{adjoint of $\Lc$}, is the operator $\Lc^*$ that satisfies
\begin{align}
\< f , \Lc g \>_\rho
	&=	\< \Lc^* f , g \>_\rho , &
	&\forall \, g \in \text{dom}(\Lc) , &
	&\forall \, f \in \text{dom}(\Lc^*) .
\end{align}
We say that $\Lc$ is \emph{self-adjoint} in $L^2(E,\rho)$ if $\Lc = \Lc^*$ and $\text{dom}(\Lc) = \text{dom}(\Lc^*)$.
\end{definition}

\noindent
Expression \eqref{eq:self.adjoint} is called the \emph{self-adjoint form of $\Ac$} as, for functions $f$ and $g$ that satisfy appropriate boundary conditions (to be determined below) we have
\begin{align}
\< f , \Ac g \>_m
	&=	\< \Ac f, g \>_m , &
\< f, g \>_m
	&:=	\int_l^r \dd x \, m(x) f(x) g(x) .
\end{align}
To see this, simply note that
\begin{align}
&\< f , \Ac g \>_m
	=	\int_I \dd x \, m(x) f(x) \cdot \frac{1}{m(x)} \d_x \( \frac{1}{s(x)} \d_x g(x) \) \\
	&=	\int_I \dd x \, f(x) \cdot \d_x \( \frac{1}{s(x)} \d_x g(x) \) \\
	&=	- \int_I \dd x \, \( \frac{1}{s(x)} \d_x g(x) \) \cdot \d_x f(x) 
				+ \[ f(x) \cdot \( \frac{1}{s(x)} \d_x g(x) \) \]_{\d I} \\
	&=	- \int_I \dd x \, \( \d_x g(x) \) \cdot  \frac{1}{s(x)} \d_x f(x) 
				+ \[ f(x) \cdot \( \frac{1}{s(x)} \d_x g(x) \) \]_{\d I} \\
	&=	\int_I \dd x \, g(x) \cdot \d_x \( \frac{1}{s(x)} \d_x f(x) \)  
				+ \[ f(x) \cdot \( \frac{1}{s(x)} \d_x g(x) \) - g(x) \cdot \( \frac{1}{s(x)} \d_x f(x) \) \]_{\d I} \\ 
	&=	\int_I \dd x \, m(x) g(x) \cdot \frac{1}{m(x)} \d_x \( \frac{1}{s(x)} \d_x f(x) \) 
				+ \[ f(x) \cdot \( \frac{1}{s(x)} \d_x g(x) \) - g(x) \cdot \( \frac{1}{s(x)} \d_x f(x) \) \]_{\d I} \\
	&=	\< \Ac f, g \>_m 
				+ \[ f(x) \cdot \( \frac{1}{s(x)} \d_x g(x) \) - g(x) \cdot \( \frac{1}{s(x)} \d_x f(x) \) \]_{\d I} .
\end{align}
Thus, we have
\begin{align}
\< f , \Ac g \>_m
	&=	\< \Ac f , g \>_m , &
	&\text{assuming}&
\[ f(x) \cdot \( \frac{1}{s(x)} \d_x g(x) \) - g(x) \cdot \( \frac{1}{s(x)} \d_x f(x) \) \]_{\d I}
	&=	0 .
\end{align}
If we require that $\Ac$ acts only on functions $h$ that satisfy
\begin{align}
h(l) + \alpha \frac{h'(l)}{s(l)} 
	&=	0 , &
h(r) + \gam \frac{h'(r)}{s(r)} 
	&=	0 , &
	&\text{for some constants $\alpha, \gam, \in \Rb \cup \{ \infty \} $}, \label{eq:bcs-h}
\end{align}
then the boundary terms will vanish because we have
\begin{align}
f(l)\frac{g'(l)}{s(l)} - g(l) \frac{f'(l)}{s(l)}
	&=	f(l)g(l) \Big( \frac{g'(l)}{s(l)g(l)} - \frac{f'(l)}{s(l)f(l)} \Big)
	=		f(l)g(l)  \Big( \frac{-1}{\alpha} - \frac{-1}{\alpha} \Big) 
	= 	0 ,
\end{align}
and likewise at $x=r$. Thus, the operator $\Ac$ with domain
\begin{align}
\text{dom}(\Ac)
	&=	\{ h \in L^2(I,m) : \Ac h \in L^2(I,m) \text{ and \eqref{eq:bcs-h} holds} \} .
\end{align}
is self-adjoint in $L^2(I,m)$.  
Note that $\Ac$ is \emph{not} self-adjoint in $L^2(I,\dd x)$.
Also note, when we talk about the operator $\Ac$ we are really talking about the pair $(\Ac, \text{dom}(\Ac))$.  The domain of $\Ac$ includes boundary conditions that must be satisfied in order of $\Ac$ to be self-adjoint in $L^2(I,m)$.

\subsection*{Boundary classification and boundary conditions}
The endpoints $l$ and $r$ of the interval $I$ fall in to one of four categories, which we describe below.
In order to determine which of these four categories the end-points belong to we introduce $S$, the \emph{scale measure}
\begin{align}
S[x,y]
	&:=	\int_x^y \dd z \, s(z) , &
x,y
	&\in (l,r) , &
S(l,y]
	&=	\lim_{x \searrow l} S[x,y] , &
S[x,r)
	&=	\lim_{y \nearrow r} S[x,y] . &
\end{align}
Next we introduce four useful integrals related to the scale measure and the speed density $m$ given in \eqref{eq:s.and.m}
\begin{align}
I_l
	&=  \int_l^x \dd z \, m(z) S(l,z] , &
I_r
	&=	\int_x^r \dd z \, m(z) S[z,r) , \label{eq:I} \\
J_l
	&=	\int_l^x \dd z \, m(z) S[z,x] , &
J_r
	&=	\int_x^r \dd z \, m(z) S[x,z] . \label{eq:J}
\end{align}
Here, $x$ is any point in the open interval $(l,r)$.

\begin{definition}
Let $X$ be a time-homogeneous scalar diffusion, which lives on an interval $I$ with endpoints $l$ and $r$.
An endpoint $l$ or $r$ is said to be
\begin{itemize*}
\item \emph{regular} if $I<\infty$ and $J<\infty$.
\item \emph{exit} if $I<\infty$ and $J=\infty$.
\item \emph{entrance} if $I=\infty$ and $J<\infty$.
\item \emph{natural} if $I=\infty$ and $J=\infty$.
\end{itemize*}
where $I$ and $J$ are given in \eqref{eq:I} and \eqref{eq:J}.
\end{definition}

\noindent
The definition above, while precise, gives us very little intuition as to what the four different boundary classifications mean.  Thus, we elaborate a bit.
\begin{itemize*}
\item \emph{Regular boundary} The process $X$ can be started from a regular boundary and $X$ can reach a regular boundary in finite time.
\item \emph{Exit boundary}: The process $X$ cannot be started from an exit boundary but $X$ can reach an exit boundary in finite time.  If $X$ reaches an exit boundary, it does not return.
\item \emph{Entrance boundary}: The process $X$ can be started from an entrance boundary but $X$ cannot reach an entrance boundary in finite time.
\item \emph{Natural boundary}: The process $X$ cannot be started from a natural boundary nor can $X$ reach a natural boundary in finite time.
\end{itemize*}
We can further classify the boundaries as follows:
\begin{itemize*}
\item \emph{Unattainable}: $I=\infty$.  Entrance and natural boundaries are unattainable.  A diffusion starting in $(l,r)$ cannot reach an unattainable boundary in finite time.
\item \emph{Attainable}: $I<\infty$.  Regular and exit boundaries are attainable. A a diffusion starting in $(l.r)$ can reach an attainable boundary in finite time.
\end{itemize*}
We must specify the behavior of $X$ at a regular boundary.
Different boundary behaviors correspond to different boundary conditions for $\text{dom}(\Ac)$.
The two most common behaviors are killing and reflecting.
\\[0.5em]
If a regular boundary is specified as \emph{killing} then the process $X$ is killed as soon as it hits this boundary and it cannot return to the state space $I$.  Thus, if $l$ and $r$ are regular killing boundaries we clearly have
\begin{align}
\Pb( X_T \in \dd y | X_t = l ) \equiv \Gam(t,l;T,y) \dd y
	&=	0 , &
\Pb( X_T \in \dd y | X_t = r ) \equiv \Gam(t,r;T,y) \dd y
	&=	0 .
\end{align}
Since $\Gam(t,\cdot;T,y) \in \text{dom}(\Ac)$, we clearly want $\Ac$ to act on functions $f$ that satisfy $f(l) = f(r) = 0$.
\\[0.5em]
If a regular boundary is specified as \emph{reflecting}, then the process $X$ is instantaneously reflected back into $I$ if it hits this boundary.
Thus, if $l$ and $r$ are regular reflecting boundaries we clearly have
\begin{align}
1
	&=	\int_I \dd y \, \Gam(t,x;T,y) , \\
0
	&=	\int_I \dd y \, \d_T \Gam(t,x;T,y) \\
	&=	\int_I \dd y \, \Ac^* \Gam(t,x;T,y) \\
	&=	\int_I \dd y \, \d_y \( - \mu(y) \Gam(t,x;T,y) + \tfrac{1}{2} \d_y \sig^2(y) \Gam(t,x;T,y) \) \\
	&=	\[ - \mu(y) \Gam(t,x;T,y) + \tfrac{1}{2} \d_y \sig^2(y) \Gam(t,x;T,y) \]_{\d I} \\
	&=	\[ - \mu(y) m(y) \frac{\Gam(t,x;T,y)}{m(y)} + \tfrac{1}{2} \d_y \sig^2(y) m(y) \frac{\Gam(t,x;T,y)}{m(y)} \]_{\d I} \\
	&=	\[ - \frac{2 \mu(y)}{\sig^2(y)} \exp\Big( \int \dd y \frac{2 \mu(y)}{\sig^2(y)} \Big) \frac{\Gam(t,x;T,y)}{m(y)} 
				+ \d_y \exp \Big( \int \dd y \frac{2 \mu(y)}{\sig^2(y)} \Big) \frac{\Gam(t,x;T,y)}{m(y)} \]_{\d I} \\
	&=	\[ \exp \Big( \int \dd y \frac{2 \mu(y)}{\sig^2(y)} \Big) \d_y \frac{\Gam(t,x;T,y)}{m(y)} \]_{\d I} \\
	&=	\[ \frac{1}{s(y)} \d_y \frac{\Gam(t,x;T,y)}{m(y)} \]_{\d I}.
\end{align}
Since $\Gam(t,x;T,\cdot) \in \text{dom}(\Ac^*)$, if follows that $\Ac^*$ acts on functions $f$ that satisfy $\d_y(f(y)/m(y))=0$ at $y = l$ and $y=r$.  But we are interested in the domain of $\Ac$.  Recall that $\Ac^*$ is the $L^2(I,\dd x)$ adjoint of $\Ac$ -- \emph{not} the $L^2(I,m)$ adjoint of $\Ac$.  Thus, we seek boundary conditions for $\Ac$ such that $\<f, \Ac g \> = \< \Ac^* f, g\>$.  Note that
\begin{align}
\< f , \Ac g\>
	&=	\int_I \dd x \, f(x) \frac{1}{m(x)} \d_x \( \frac{1}{s(x)}\d_x g(x) \) \\
	&=	- \int_I \dd x \, \(\d_x \frac{f(x)}{m(x)} \) \cdot \( \frac{1}{s(x)}\d_x g(x) \)
			+ \[ \frac{f(x)}{m(x)} \frac{1}{s(x)} \d_x g(x) \]_{\d I} \\
	&=	\int_I \dd x \, \d_x \( \frac{1}{s(x)} \d_x \frac{f(x)}{m(x)} \) \cdot g(x)
			+ \[ \frac{f(x)}{m(x)} \frac{1}{s(x)} \d_x g(x) - \frac{g(x)}{s(x)} \d_x \frac{f(x)}{m(x)}\]_{\d I} \\
	&=	\< \Ac^* f , g \>
			+ \[ \frac{f(x)}{m(x)} \frac{1}{s(x)} \d_x g(x) - \frac{g(x)}{s(x)} \d_x \frac{f(x)}{m(x)}\]_{\d I} .
\end{align}
Since
\begin{align}
\[ \d_x \frac{f(x)}{m(x)} \]_{\d I}
	&=	0 , &
	&\forall \, f \in \text{dom}(\Ac^*) ,
\end{align}
in order for $\<f, \Ac g\> = \< \Ac^* f , g \>$, we must have
\begin{align}
\[ \frac{1}{s(x)} \d_x g(x) \]_{\d I}
	&=	0 , &
&\forall \, g \in \text{dom}(\Ac) .
\end{align}
\textbf{To summarize}: If the endpoints of $I$ are regular, we must specify if the boundary is killing or reflecting.  The correct boundary conditions to impose for the generator $\Ac$ are
\begin{align}
\text{killing}:&&
\Big[ f(x) \Big]_{\d I}
	&=	0 , &
\text{reflecting}:&&
\Big[ \frac{1}{s(x)} \d_x f(x) \Big]_{\d I}
	&=	0 .
\end{align}
Although we will not derive it, the $\text{dom}(\Ac)$ should also satisfy certain boundary conditions at exit and entrance boundaries
\begin{align}
\text{exit}:&&
\Big[ f(x) \Big]_{\d I}
	&=	0 , &
\text{entrance}:&&
\Big[ \frac{1}{s(x)} \d_x f(x) \Big]_{\d I}
	&=	0 .
\end{align}
No boundary conditions are required for natural boundaries.

\subsection*{Eigenfunction expansions}
Consider the following eigenvalue equation
\begin{align}
\Ac \psi_n
	&=	\lam_n \psi_n , &
\psi_n
	\in \text{dom}(\Ac) .
\end{align}
Here, and throughout this subsection, we have implicitly assumed that the spectrum of $\Ac$ is discrete and simple (i.e., no repeated eigenvalues).

\begin{theorem}
Eigenfunctions corresponding to different eigenvalues of $\Ac$ are orthogonal in $L^2(I,m)$.
\end{theorem}

\begin{proof}
Let $n \neq k$.  Then
\begin{align}
\lam_n \< \psi_n , \psi_k \>_m
	&=	\< \Ac \psi_n , \psi_k \>_m
	=		\< \psi_n , \Ac \psi_k \>_m
	=		\lam_k \< \psi_n , \psi_k \>_m .
\end{align}
It follows that $\< \psi_n , \psi_k \>_m = 0$.
\end{proof}

\begin{theorem}
\label{thm:eigenfunction-kbe}
Suppose the eigenfunctions of $\Ac$ are complete in $L^2(I,m)$ and $\varphi \in L^2(I,m)$.  Then the function $u(t,x):= \Eb[ \varphi(X_T) | X_t = x]$ is given by 
\begin{align}
u(t,x)
	&=	\sum_n \ee^{(T-t)\lam_n} \< \psi_n, \varphi \>_m \psi_n(x) , \label{eq:u.eigen}
\end{align}
where $\Ac \psi_n = \lam_n \psi_n$, the eigenfunctions are normalized (i.e., $\<\psi_n,\psi_n\>_m = 1$ for all $n$) and $m$ is the speed measure of $\Ac$.
\end{theorem}

\begin{proof}
We must show that the function $u$ satisfies the KBE.  We have
\begin{align}
(\d_t + \Ac) u(t,x)
	&=	\sum_n \< \psi_n, \varphi \>_m \( \d_t \ee^{(T-t)\lam_n} \psi_n(x)  + \ee^{(T-t)\lam_n} \Ac \psi_n(x) \) \\
	&=	\sum_n \< \psi_n, \varphi \>_m \( -\lam_n \ee^{(T-t)\lam_n} \psi_n(x)  + \ee^{(T-t)\lam_n} \lam_n \psi_n(x) \)
	=		0 . \\
u(T,x)
	&=	\sum_n \< \psi_n, \varphi \>_m \psi_n(x)
	=		\varphi(x) .
\end{align}
where, to establish the terminal condition $u(T,x) = \varphi(x)$ we have used the fact that the eigenfunctions $(\psi_n)_{n \geq 0}$ are complete in $L^2(I,m)$.
\end{proof}

\begin{corollary}
Under the conditions of Theorem \ref{thm:eigenfunction-kbe}, the transition density $\Gam(t,x;T,y)$ of $X$ has the following eigenfunction expansion
\begin{align}
\Gam(t,x;T,y)
	&=	m(y) \sum_n \ee^{(T-t)\lam_n} \psi_n(y) \psi_n(x) . \label{eq:Gamma.eigen}
\end{align}
\end{corollary}

\begin{proof}
The transition density $\Gam(t,x;T,y)$ satisfies the KBE with a terminal condition $\Gam(T,x;T,y) = \del_y(x)$.  Setting $\varphi = \del_y$ in \eqref{eq:u.eigen},we obtain
\begin{align}
\Gam(t,x;T,y)
	&=	\sum_n \ee^{(T-t)\lam_n} \< \psi_n, \del_y \>_m \psi_n(x) \\
	&=	m(y) \sum_n \ee^{(T-t)\lam_n} \psi_n(y) \psi_n(x) ,
\end{align}
where we have used
\begin{align}
\< \psi_n, \del_y \>_m
	&= \int_I \dd x \, \psi_n(x) \del_y(x) m(x)
	=		\psi_n(y) m(y) .
\end{align}
\end{proof}

\begin{example}[Brownian motion in a finite interval]
Consider a diffusion $X = (X_t)_{t \geq 0}$ that lives on a finite interval $(l,r)$ and satisfies the SDE
\begin{align}
\dd X_t
	&=	\dd W_t .
\end{align}
We identify the drift $\mu(x) = 0$ and the diffusion coefficient $\sig(x) = 1$.  The generator $\Ac$ of $X$ and the speed density $m$ are given by
\begin{align}
\Ac
	&=	\tfrac{1}{2} \d_x^2 , &
m(x)
	&=	\frac{1}{r-l} .
\end{align}
One can easily check that the endpoints $l$ and $r$ are regular.  As such, we must specify the behavior at the endpoints as either killing or reflecting.  In the case of killing boundaries we have
\begin{align}
\text{killing}:&&
\Ac \psi_n
	&=	\lam_n \psi_n , &
\psi_n(l)
	&=	\psi_n(r) = 0 , \\ &&
\psi_n(x)	
	&=	\sqrt{2} \sin \( \frac{n \pi}{r-l}(x-l) \) , &
\lam_n
	&=	- \frac{1}{2} \( \frac{n \pi}{r-l} \)^2 , &
n
	&=	1, 2, \ldots ,
\end{align}
where we have normalized the eigenfunctions so that $\< \psi_n, \psi_k \>_m = \del_{n,k}$.
In the case of reflecting boundaries we have
\begin{align}
\text{reflecting}:&&
\Ac \psi_n
	&=	\lam_n \psi_n , &
\psi_n'(l)
	&=	\psi_n'(r) = 0 , \\ &&
\psi_n(x)	
	&=	\sqrt{2} \cos \( \frac{n \pi}{r-l}(x-l) \) , &
\lam_n
	&=	- \frac{1}{2} \( \frac{n \pi}{r-l} \)^2 , &
n
	&=	0, 1, 2, \ldots ,
\end{align}
where, again, we have normalized the eigenfunctions.
One can immediately write down the transition density $\Gam(t,x;T,y)$ using \eqref{eq:Gamma.eigen}.  Of course, the $\sin$ and $\cos$ expansions should be familiar to anybody that has solved the heat equation on a finite domain.
\end{example}

\begin{example}[Ornstein-Uhlenbeck process]
Consider now an OU process
\begin{align}
\dd X_t
	&=	- X_t \dd t + \dd W_t ,
\end{align}
which has a state space $I = \Rb$.  We identify the drift $\mu(x) = - x$ and volatility coefficient $\sig(x) = 1$.  The generator and speed density of $X$ are
\begin{align}
\Ac
	&=	- x \d_x + \tfrac{1}{2} \d_x^2 , &
m(x)
	&=	\ee^{-x^2} .
\end{align}
One can easily verify that the endpoints $l= -\infty$ and $r = + \infty$ are natural.  As such, we do not need to specify any boundary conditions.  The eigenfunctions and eigenvalues of $\Ac$ are given by
\begin{align}
\psi_n(x)
	&=	H_n(x) , &
\lam_n
	&= - n , &
n
	&=	0, 1, 2, \ldots ,
\end{align}
where $(H_n)_{n \geq 0}$ are the Hermite polynomials, properly normalized so that $\< \psi_n, \psi_k \>_m = \del_{n,k}$.
\end{example}

\subsubsection{Heuristic computations}
We have used the fact that $\Ac$ is self-adjoint in $L^2(I,m)$ in order to derive to write an eigenfunction expansion for $u(t,x) = \Eb[ \varphi(X_T) | X_t = x]$, the solution of the KBE.  In fact, this is a special case of a more general method of solving PDEs and ODEs involving a self-adjoint operator $\Ac$.  In what follows, we continue to assume that the generator $\Ac$ of a diffusion $X$ has a discrete spectrum, and that the eigenfunctions of $\Ac$ are a complete basis in $L^2(I,m)$.
\\[0.5em]
For any Borel-measurable function $g : \Rb \to \Rb$ we define the operator
\begin{align}
g(\Ac) \cdot
	&=	\sum_n g(\lam_n) \< \psi_n , \cdot \>_m \psi_n , &
\text{dom}(g(\Ac))
	&:=	\{ f \in L^2(I,m) : \sum_n g^2(\lam) \< \psi_n , f \>_m^2 < \infty \} . \label{eq:g.A}
\end{align}
We call \eqref{eq:g.A} the \emph{eigenfunction} or \emph{spectral representation} of $g(\Ac)$.
For example, taking $g(\lam)=1$ gives the spectral representation of the identity operator
\begin{align}
\Iv \, \cdot
	&=	\sum_n \< \psi_n , \cdot \>_m \psi_n , &
\text{dom}(\Iv)
	&:=	\{ f \in L^2(I,m) : \sum_n \< \psi_n , f \>_m^2 < \infty \} .
\end{align}
Similarly, taking $g(\lam) = (\lam-\mu)^{-1}$ and defining $ \Rc_\mu := g(\Ac)$ gives
\begin{align}
\Rc_\mu \cdot
	&=	\sum_n \frac{1}{\lam_n - \mu} \< \psi_n , \cdot \>_m \psi_n , &
\text{dom}(\Rc_\mu)
	&:=	\{ f \in L^2(I,m) : \sum_n \(\frac{1}{\lam_n-\mu}\)^2 \< \psi_n , f \>_m^2 < \infty \} ,
\end{align}
and taking $g(\lam) = \ee^{(T-t)\lam}$ and defining $\Pc(t,T) := g(\Ac)$ gives
\begin{align}
\Pc(t,T) \cdot
	&=	\sum_n \ee^{(T-t)\lam_n} \< \psi_n , \cdot \>_m \psi_n , &
\text{dom}(\Pc(t,T))
	&:=	\{ f \in L^2(I,m) : \sum_n \ee^{(T-t)\lam_n} \< \psi_n , f \>_m^2 < \infty \} ,
\end{align}
Note that $\Pc(t,T)$ as defined above is in fact the semigroup operator we introduced previously.  To see this note that
\begin{align}
\Pc(t,T) \varphi(x)
	&=	\sum_n \ee^{(T-t)\lam_n} \< \psi_n , \varphi \>_m \psi_n(x) \\
	&=	\sum_n \ee^{(T-t)\lam_n} \int_I \dd y \, \psi_n(y) \varphi(y) m(y) \psi_n(x) \\
	&=	\int_I \dd y \Big( m(y) \sum_n \ee^{(T-t)\lam_n} \psi_n(y) \psi_n(x) \Big) \varphi(y) \\
	&=	\int_I \dd y \, \Gam(t,x;T,y) \varphi(y) ,
\end{align}
where we have used \eqref{eq:Gamma.eigen}.
\\[0.5em]
The usefulness in defining operators of the form $g(\Ac)$ is as follows.  Consider a linear ODE or PDE for a function $u$ involving the operator $\Ac$.  One can obtain a solution to this ODE or PDE as follows:
\begin{enumerate*}
\item Solve the ODE or PDE assuming $\Ac$ is a constant; the solution will involve terms of the form $g(\Ac)f$.
\item Replace $g(\Ac)f$ with its spectral representation.
\end{enumerate*}
Let us see how this works with a few examples.

\begin{example}
Let us check that the above two-step method works for the KBE.  We want to solve
\begin{align}
(\d_t + \Ac)u(t,\, \cdot \,)
	&=	0 , &
u(T,\, \cdot \,)
	&= \varphi , \label{eq:kbe.again}
\end{align}
where $\varphi \in L^2(I,m)$.
Treating $\Ac$ like a constant, we have an ODE in $t$.  Solving this ODE and then using the spectral representation any expression involving $\Ac$ we obtain
\begin{align}
u(t,\, \cdot \,)
	&=	\ee^{(T-t)\Ac} \varphi
	=		\sum_n \ee^{(T-t)\lam_n} \< \psi_n, \varphi\>_m \psi_n ,
\end{align}
which agrees with the expression given in \eqref{eq:u.eigen}.  Noting that $\ee^{(T-t)\Ac} = \Pc(t,T)$ we can write the solution to \eqref{eq:kbe.again} compactly as
$u(t,\, \cdot \,) = \Pc(t,T)\varphi$.
\end{example}

\begin{example}
Consider the following PDE involving $\Ac$
\begin{align}
(\d_t + \Ac) u(t,\, \cdot \,) + g(t,\, \cdot \,)
	&=	0 , &
u(T,\, \cdot \,)
	&=	\varphi , \label{eq:fk.again}
\end{align}
where $\varphi \in L^2(I,m)$ and $g(s,\, \cdot \,)\in L^2(I,m)$ for all $s \in [0,T]$.
Treating $\Ac$ as a constant, we have an ODE in $t$.  The solution is given by
\begin{align}
u(t,\, \cdot \,)
	&=	\ee^{(T-t)\Ac} \varphi + \int_t^T \dd s \, \ee^{(s-t)\Ac} g(s,\, \cdot \,) .
\end{align}
Using the spectral representation of $\ee^{(s-t)\Ac}$ we obtain
\begin{align}
u(t,\, \cdot \,)
	&=	\sum_n \ee^{(T-t)\lam_n} \< \psi_n, \varphi\>_m \psi_n + \int_t^T \dd s \, \sum_n \ee^{(s-t)\lam_n} \< \psi_n, g(s,\, \cdot \,) \>_m \psi_n .
\end{align}
Again, recognizing that $\ee^{(s-t)\Ac} = \Pc(t,s)$ we can write the solution of \eqref{eq:fk.again} as
\begin{align}
u(t,\, \cdot \,)
	&=	\Pc(t,T) \varphi + \int_t^T \dd s \, \Pc(t,s) g(s,\, \cdot \,) .
\end{align}
\end{example}

\begin{example}
Consider one last example.  We wish to solve
\begin{align}
(\Ac - \mu) u
	&=	g . \label{eq:u.poisson}
\end{align}
Treating $\Ac$ as a constant, we have an algebraic equation for $u$.  Solving this, and using the spectral representation for any expression involving $\Ac$ we obtain
\begin{align}
u
	&=	\frac{g}{\Ac - \mu} = \sum_n \frac{1}{\lam_n - \mu} \< \psi_n, g \>_m \psi_n .
\end{align}
Noting that $\frac{1}{\Ac - \mu} = \Rc_\mu$ we can write the solution $u$ of \eqref{eq:u.poisson} compactly as $u = \Rc_\mu g$.  Note that his solution makes sense only if $g \in \text{dom}(\Rc_\mu)$.  In particular, this means that we must have $\lam_n \neq \mu$ for all $n$.
\end{example}

\subsection*{The meaning of ``scale'' and ``speed''}
From the scale density $s$ we can define a \textit{scale function} $S$, which is given by
\begin{align}
S(x)
	&:=	\int \dd x \, s(x) .
\end{align}
By construction, the scale function $S$ is one-to-one.  It is interesting to note that the process $S(X)=(S(X_t))_{t \geq 0}$ is a martingale, as the following computation shows
\begin{align}
\dd S(X_t)
	&=	\Ac S(X_t) \dd t + \sig(X_t) \d_x S(X_t) \dd W_t \\
	&=	\Big( \mu(X_t) s(X_t) + \tfrac{1}{2} \sig^2 (X_t) \d_x s(X_t) \Big) \dd t + \sig(X_t) s(X_t) \dd W_t \\
	&=	\Big( \mu(X_t) s(X_t) + \tfrac{1}{2} \sig^2 (X_t) \Big( \frac{- 2 \mu(X_t)}{\sig^2(X_t)} \Big) s(X_t) \Big) \dd t + \sig(X_t) s(X_t) \dd W_t \\
	&=	\sig(X_t) s(X_t) \dd W_t . \label{eq:dSX}
\end{align}
Now,  consider two regular points $\{l\}$ and $\{r\}$ of a scalar diffusion $X$.  Let $\tau_l$ and $\tau_r$ be the first hitting times of $X$ to $l$ and $r$, respectively
\begin{align}
\tau_l
	&:=	\inf\{ t \geq 0 : X_t = l \} , &
\tau_r
	&:=	\inf\{ t \geq 0 : X_t = r \} .
\end{align}
We wish to compute $\Pb( \tau_l < \tau_r | X_0 = x)$, the probability that $X$ hits $l$ prior to hitting $r$.  Let us define
\begin{align}
\tau
	&:=	\tau_l \wedge \tau_r .
\end{align}
As the process $S(X)$ is a martingale it follows from Theorems \ref{thm:stopped.martingale} and \ref{thm:sampling} that
\begin{align}
S(x)
	&=	\Eb S( X_\tau)
	 =	\Pb( \tau_l < \tau_r | X_0 = x ) S( l ) + \Pb( \tau_r < \tau_l | X_0 = x ) S( r ) .
\end{align}
Also, because $\{ \tau_l < \tau_r \} \cup \{ \tau_l < \tau_r \} = \Om$ we have trivially that
\begin{align}
1
	&=	\Pb( \tau_l < \tau_r | X_0 = x ) + \Pb( \tau_r < \tau_l | X_0 = x ) .
\end{align}
Using the above two equations, we find
\begin{align}
\Pb( \tau_l < \tau_r | X_0 = x )
	&=	\frac{S(r) - S(x)}{S(r) - S(l)} , &
\Pb( \tau_r < \tau_l | X_0 = x )
	&=	\frac{S(x) - S(l)}{S(r) - S(l)} . \label{eq:tl.tr}
\end{align}
If we set $Y = S(X)$, let $L = S(l)$ and $R = S(r)$, and define
\begin{align}
\tau_L
	&:=	\inf\{ t \geq 0 : Y_t = L ) , &
\tau_R
	&=	\inf\{ t \geq 0 : Y_t = R ) , \label{eq:tL.tR}
\end{align}
then equation \eqref{eq:tl.tr} becomes
\begin{align}
\Pb( \tau_L < \tau_R | Y_0 = y )
	&=	\frac{R - y}{R - L} , &
\Pb( \tau_R < \tau_L | Y_0 = y )
	&=	\frac{y - L}{R - L} . \label{def:natural.scale}
\end{align}
Any process that satisfies \eqref{eq:tL.tR}-\eqref{def:natural.scale} is said to be \textit{on a natural scale}.  Thus, given a scalar diffusion $X$ we can always associate with it a diffusion $S(X)$ that is on a natural scale. 
\\[0.5em]
\blu{To do: meaning of \textit{speed}.}



\section{Extensions to higher dimensions}
\label{sec:sde-d}
Throughout this Chapter, we have mostly dealt with Scalar SDEs and PDEs with one spacial dimension.  To conclude this chapter, we simply restate some of the main results in higher dimensions, without providing the conditions under which they hold rigorously.
\\[0.5em]
Throughout this section, we consider $d$-dimensional diffusion process $X = (X_t)_{t \geq 0}$ that satisfies the following SDE
\begin{align}
\dd X_t
	&=	\mu(t,X_t) \dd t + \sig(t,X_t) \dd W_t , \label{eq:dX.multi}
\end{align}
where $W$ is an $m$-dimensional Brownian motion and
\begin{align}
\mu 
	&: \Rb_+ \times \Rb^d \to \Rb^d , &
\sig
	&: \Rb_+ \times \Rb^d \to \Rb_+^{d \times m} .
\end{align}
SDE \eqref{eq:dX.multi} has a unique strong solution, which is square integrable for all $t$ (i.e., $\Eb X_t^2 < \infty$) and adapted to the filtration $\Fb=(\Fc_t)_{t \geq 0}$ generated by $W$ (i.e., $X_t \in \Fc_t$) if the following are satisfies
\begin{align}
&\text{Linear Growth}:&
| \mu(t,x) | + | \sig(t,x) |
	&\leq C_1 (1 + |x| ) , \label{eq:lin.growth} \\
&\text{Lipschitz Continuous}:&
| \mu(t,x) - \mu(t,y) | + | \sig(t,x)-\sig(t,y) |
	&\leq C_2  |x - y| , \label{eq:lipschitz.continuous}
\end{align}
for some constants $C_1$ and $C_2$, independent of $t,x,y$, where
\begin{align}
|\sig|^2
	&:= \sum_{i=1}^d \sum_{j=1}^m \sig_{ij}^2 , &
|\mu|^2
	&:=	\sum_{i=1}^d \mu_i^2 , &
|x|^2
	&:=	\sum_{i=1}^d x_i^2 .
\end{align}
For any continuous-time Markov process $X$ that lives in $\Rb^d$, we can define an \textit{infinitesimal generator} $\Ac(t)$ as follows
\begin{align}
\Ac(t) \varphi(x)
	&:= \lim_{T \searrow t} \frac{1}{T-t} \Big( \Eb ( \varphi(X_T) | X_t = x ) - \varphi(x) \Big) ,
\end{align}
where $\varphi: \Rb^d \to \Rb$ may be any function for which the limit exists.
When $X$ is given by \eqref{eq:dX.multi}, the generator $\Ac(t)$ and its $L^2(\Rb^d,\dd x)$ adjoint are given by
\begin{align}
\Ac(t)
	&=	\sum_{i=1}^d \mu_i(t,x)\d_{x_i} + \tfrac{1}{2} \sum_{i=1}^d \sum_{j=1}^d (\sig \sig^\text{T})_{i,j}(t,x) \d_{x_i} \d_{x_j} , \\
\Ac^*(t)
	&=	- \sum_{i=1}^d \d_{y_i} \mu_i(t,y) + \tfrac{1}{2} \sum_{i=1}^d \sum_{j=1}^d \d_{y_i} \d_{y_j} (\sig \sig^\text{T})_{i,j}(t,y) . 
\end{align}
We can write the It\^o Formula most compactly as
\begin{align}
\dd \varphi(X_t)
	&=	\Ac(t) \varphi(X_t) \dd t + (\nabla \varphi)(X_t) \sig(t,X_t) \dd W_t .
\end{align}
Thus, we can think of $\Ac(t)$ as the operator that appears in the $\dd t$-term when we take the differential $\dd \varphi(X_t)$.
\\[0.5em]
The transition density $\Gam$, defined by $\Gam(t,x;T,y) \dd y = \Pb(X_T \in \dd y | X_t = x)$ satisfies the KBE in the backwards variables $(t,x)$ and the KFE in the forward variables $(T,y)$
\begin{align}
\text{KBE}:&&
\(\d_t + \Ac(t) \) \Gam(t,\cdot;T,y)
	&=	0 , &
\Gam(T,\cdot;T,y)
	&=	\del_y , \\
\text{KFE}:&&
\(-\d_T + \Ac^*(T) \) \Gam(t,x;T,\cdot)
	&=	0 , &
\Gam(t,x;t,\cdot)
	&=	\del_x .
\end{align}
If we define a function $u : \Rb_+ \times \Rb \to \Rb$ by
\begin{align}
u(t,x)
	&:=	\Eb \Big[ \ee^{-A(t,T)}\varphi(X_T) + \int_t^T \ee^{-A(t,s)}g(s,X_s) \dd s \Big| X_t = x \Big] , 
\end{align}
where $A(t,T)$ is given by
\begin{align}
A(t,T)
	&=	\int_t^T \gam(s,X_s) \dd s, &
\gam
	&:	\Rb_+ \times \Rb^d \to \Rb_+ .
\end{align}
Then $u$ satisfies
\begin{align}
\Big( \d_t - \gam(t,\cdot) + \Ac(t) \Big) u(t,\cdot) + g(t,\cdot)
	&=	0 , &
u(T,\cdot)
	&= \varphi .
\end{align}
Lastly, consider the time-homogeneous diffusion: $\mu(t,x) = \mu(x)$ and $\sig(t,x) = \sig(x)$.  Let $D \subset \Rb^d$ be an open, connect set.  Denote by $\d D$ the boundary of $D$.  Assume $X_0 \in D$ and define the hitting time
\begin{align}
\tau_D
	&:=	\inf \{ t \geq 0 : X_t \notin D \} .
\end{align}
Assume $\Pb( \tau_D < \infty ) = 1$ and let 
\begin{align}
u(x)
	&=	\Eb \Big[ \ee^{-\lam(\tau_D-t)}\varphi( X_{\tau_D} ) + \int_t^{\tau_D} \ee^{- \lam(s-t)} g(X_s) \dd s \Big| X_t = x \Big] , &
t
	< \tau_D .
\end{align}
Then $u$ satisfies
\begin{align}
(\Ac - \lam) u + g
	&= 0 , &
	&\text{in $D$} , \\
u
	&=	\varphi , &
	&\text{on $\d D$} .
\end{align}



\section{Exercises}
\label{sec:exercises-9}

\begin{exercise}
\label{exer:9.1}
A (one-dimensional) \emph{backward stochastic differential equation} (BSDE), defined on a probability space filtered probability space $(\Om,\Fc,\Fb=(\Fc_t)_{t \geq 0},\Pb)$, is an equation of the form
\begin{align}
\dd Y_t
	&=	- f(\om,t,Y_t,Z_t) \dd t + Z_t \dd W_t , &
Y_T
	&=	\xi , \label{eq:bsde}
\end{align}
where $W=(W_t)_{0 \leq t \leq T}$ is an $(\Pb,\Fb)$-Brownian motion, the process $Y=(Y_t)_{0 \leq t \leq T}$ lives in $\Rb$, the process $Z=(Z_t)_{0 \leq t \leq T}$ lives in $\Rb$, the \emph{driver} $f: \Om \times [0,T] \times \Rb \times \Rb \to \Rb$, satisfies $f(\om,t,Y_t,Z_t) \in \Fc_t$ for all $t \in [0,T]$ and the random variable $\xi$ is $\Fc_T$-measurable.  A \emph{solution} of BSDE \eqref{eq:bsde} is any pair of $\Fb$-adapted processes $(Y,Z)$ such that the terminal condition $Y_T = \xi$ is satisfied.  A \emph{forward-backward stochastic differential equation} (FBSDE), is a BSDE of the form
\begin{align}
\begin{aligned}
\dd X_t
	&=	\mu(t,X_t,Y_t,Z_t) \dd t + \sig(t,X_t,Y_t) \dd W_t , &
X_0
	&=	x , \\
\dd Y_t
	&=	- f(t,X_t,Y_t,Z_t) \dd t + Z_t \dd W_t , &
Y_T
	&=	\varphi(X_T) \in \Rb, 
\end{aligned} \label{eq:fbsde}
\end{align}
where the processes $W$, $Y$, and $Z$ are as described above, the process $X = (X_t)_{0 \leq t \leq T}$ lives in $\Rb$, and the functions $\mu$, $\sig$, $f$ and $\varphi$ are maps
\begin{align}
\mu
	&: [0,T] \times \Rb \times \Rb \times \Rb \to \Rb , &
\sig
	&: [0,T] \times \Rb \times \Rb \to \Rb , \\
f
	&: 	[0,T] \times \Rb \times \Rb \times \Rb \to \Rb , &
\varphi
	&:	\Rb \to \Rb .
\end{align}
In general, the coefficient $\sig$ in \eqref{eq:fbsde} could depend on $Z$ as well.  However, for simplicity, we do not consider this case here.
FBSDEs naturally arise in mathematical finance where components of $X$ are assets in a market, components of $Y$ are values of hedging portfolios, and components $Z$ are the associated hedges.
\\[0.5em]
We wish to solve FBSDE \eqref{eq:fbsde} -- meaning we wish to find a pair of $\Fb$-adapted processes $(Y,Z)$ such that $Y_T = \varphi(X_T)$.
Let us supposed that $Y_t = u(t,X_t)$ for some function $u: [0,T] \times \Rb \to \Rb$, which is to be determined.  Let us further suppose that $u \in C^{1,2}$.  
Compute $\dd u(t,X_t)$ and compare your result to the expression given for $\dd Y_t$ in \eqref{eq:fbsde}.
Conclude that, if the function $u$ satisfies the semilinear PDE
\begin{align}
0
	&=	(\d_t + \Ac^u(t)) u(t) + f^u(t) , &
u(T)
	&=	\varphi , \label{eq:nonlin-fk}
\end{align}
where $\Ac^u$ and $f^u$ are given by
\begin{align}
\Ac^u(t)
	&=	\mu(t,x,u(t,x),\sig(t,x) \d_x u(t,x) ) \d_x + \frac{1}{2} \sig^2(t,x,u(t,x)) \d_x^2, \\
f^u(t)
	&=	f(t,x,u(t,x),\sig(t,x)\d_x u(t,x)) ,
\end{align}
then the pair $(Y,Z)$ defined by
\begin{align}
Y_t 
	&:=  u(t,X_t) , &
Z_t
	&:= \sig(t,X_t) \d_x u(t,X_t)   , \label{eq:YZ}
\end{align}
is (at least formally) a solution of FBSDE \eqref{eq:fbsde}.
\end{exercise}

\begin{exercise}
\label{exer:9.2}
Let $W$ be a Brownian motion on a probability space $(\Om, \Fc, \Pb)$ and let $\Fb = (\Fc_t)_{0 \leq t \leq T}$ be a filtration for $W$.
Let $X$ be the solution of the following SDE
\begin{align}
\dd X_t
	&=	\kappa( \theta - X_t) \dd t + \del \sqrt{X_t} \dd W_t .
\end{align}
The process $X$ is called a \textit{Cox-Ingersoll-Ross} (CIR) process.  Define
\begin{align}
u(t,x)
	&:=	\Eb\[ \exp \( - \int_t^T X_s \dd s \) \Big| X_t = x \] .
\end{align}
(In quantitative finance, the CIR process is often used to model a stochastic short-term interest rate, in which case $u(t,X_t)$ is the price at time $t$ of a zero-coupon bond that pays one unit of currency at time $T$.)
\\[0.5em]
(a) Derive a PDE and terminal condition for the function $u$.\\
(b) To solve the PDE for $u$, try a solution of the form
\begin{align}
u(t,x)
	&=	\exp( -x A(t) - B(t) ) , \label{eq:exer.9.2.guess}
\end{align}
where $A$ and $B$ are deterministic functions of $t$.  Show that $A$ and $B$ must satisfy a pair of coupled ODEs.  Write these ODEs as well as their terminal conditions at time $T$.
\textit{Hint}: insert the expression \eqref{eq:exer.9.2.guess} into the PDE for $u$ and collect terms of like powers in $x$.  As the PDE for $u$ must hold for all values of $x$, what does this tell you about the terms that are of order $x^0$, $x^1$, etc.?\\
(c) \textit{Bonus question}: solve the ODEs (it may be helpful to note that one of the ODEs is a Riccati equation).
\end{exercise}

\begin{exercise}
\label{exer:9.3}
For $i=1,2,\ldots,d$, let $X^{(i)}$ satisfy
\begin{align}
\dd X_t^{(i)}
	&=	-\frac{b}{2} X_t^{(i)} \dd t + \frac{1}{2} \sig \dd W_t^{(i)} ,
\end{align}
where the $(W^{(i)})_{i=1}^d$ are independent Brownian motions.  Define
\begin{align}
R_t
	&:=	\sum_{i=1}^d (X_t^{(i)})^2 , &
B_t
	&:=	\sum_{i=1}^d \int_0^t \frac{1}{\sqrt{R_s}} X_s^{(i)} \dd W_s^{(i)}.
\end{align}
Show that $B$ is a Brownian motion.  Derive an SDE for $R$ that involves only $\dd t$ and $\dd B_t$ terms (i.e., no $\dd W_t^{(i)}$ terms should appear).
\end{exercise}

\begin{exercise}
\label{exer:9.4}
Consider the following 2-dimensional SDE
\begin{align}
\dd X_t
	&=	-\frac{1}{2} Z_t \dd t + \sqrt{Z_t} \dd W_t , \\
\dd Z_t
	&=	\kappa (\theta - Z_t) \dd t + \sig \sqrt{Z_t} \( \rho \dd W_t + \sqrt{1-\rho^2} \dd B_t \) ,
\end{align}
where $W$ and $B$ are independent Brownian motions on a probability space $(\Om, \Fc, \Pb)$.  Define
\begin{align}
u(t,x,z)
	&:=	\Eb[ \varphi(X_T) | X_t = x, Z_t = z ] .
\end{align}
(In quantitative finance, the process $(X,Z)$ is known as the \textit{Heston model}.  In this model, $\ee^{X}$ represents the price of a stock and $Z$ represents its stochastic variance.  If interest rates are zero and $\Pb$ is the market's pricing measure, then $u(t,X_t,Z_t)$ is the price at time $t$ of an option that pays $\varphi(X_T)$ at time $T$.)
\\[0.5em]
(a) Derive a PDE and terminal condition for the function $u$.\\
(b) Let $\uh$ be the Fourier transform of $u$ in the $x$ variable
\begin{align}
\uh(t,\xi,z)
	&=	\int_\Rb \dd x \, \ee^{- \ii \xi x} u(t,x,z) .
\end{align}
Show that $\uh$ satisfies a PDE in $(t,z)$ with a terminal condition $\uh(T,\xi,z) = \varphih(\xi)$ where $\varphih$ is the Fourier transform of $\varphi$.
\textit{Hint}: using integration by parts, we have
\begin{align}
\int_\Rb \dd x \, \ee^{- \ii \xi x} \d_x u(t,x,z)
	&=	- \int_\Rb \dd x \, u(t,x,z) \d_x \ee^{- \ii \xi x}
	=		\ii \xi \int_\Rb \dd x \, u(t,x,z) \ee^{- \ii \xi x}
	=		\ii \xi \uh(t,\xi,z) .
\end{align}
(c) To solve the PDE from part (b) assume that $\uh$ is of the form
\begin{align}
\uh(t,\xi,z)
	&=	\varphih(\xi) \ee^{A(t,\xi) + z B(t,\xi)} ,
\end{align}
where $A$ and $B$ are deterministic functions of $t$ and $\xi$.  Show that $A$ and $B$ satisfy a pair of coupled ODEs in $t$; write these ODEs as well as their terminal conditions at time $T$ (it may be a good idea to look at the \textit{hint} given in Exercise \ref{exer:9.2}).  Note, if one solves the ODEs for $A$ and $B$, then the function $u$ can be obtained by taking the inverse Fourier transform of $\uh$, that is
\begin{align}
u(t,x,z)
	&=	\frac{1}{2\pi} \int_\Rb \dd \xi \, \ee^{ \ii \xi x} \uh(t,\xi,z) .
\end{align}
Typically, this integral must be computed numerically.\\
(d) \textit{Bonus question}: solve the ODEs (it may be helpful to note that one of the ODEs is a Riccati equation).
\end{exercise}

\begin{exercise}
\label{exer:9.5}
Consider a diffusion $X = (X_t)_{t\geq 0}$ that lives on a finite interval $(l,r), \quad  0< l < r < \infty $ and satisfies the SDE
$$\dd X_t =  \mu X_t \dd t + \sigma X_t \dd W_t $$
One can easily check that the endpoints $l$ and $r$ are regular (you do not have to prove it here). Assume both endpoints are killing. Find the transition density $\Gamma(t,x;T, y)$ of $X$.
\end{exercise}

\begin{exercise}
\label{exer:9.6}
Consider a two-dimensional diffusion processes $X=(X_t)_{t \geq 0}$ and $Y=(Y_t)_{t \geq 0}$ that satisfy the SDEs
\begin{align}
\dd X_t &= \dd W_t^1\\
\dd Y_t &= \dd W_t^2
\end{align}
where $W_t^1$ and $W_t^2$ are two independent Brownian motions.  Define a function $u$ as follows
\begin{align}
u(x,y) &= \Eb [\phi (X_\tau) | X_t = x, Y_t = y] ,\\
\tau &= \text{inf} \{s \geq t : Y_s = a\}
\end{align}
(1) State a PDE and boundary conditions satisfied by the function $u$.\\
(2) Let us define the Fourier transform and inverse Fourier transform, respectively, as follows
\begin{align}
\text{Fourier Transform}:&&
\widehat{f}(\omega) 
	&:= \int e^{- i \omega x} f(x) \dd x , \\
\text{Inverse Transform}:&&
f(x) 
	&:= \frac{1}{2 \pi} \int e^{i \omega x} \widehat{f}(\omega) \dd \omega. 
\end{align}
Use Fourier transforms and a conditioning argument to derive an expression for $u(x,y)$ as an inverse Fourier transform.  Use this result to derive an explicit form for $\Pb(X_\tau \in \dd z | X_t=x, Y_t= y)$ (i.e., an expression involving no integrals).\\
(3) Show the expression you derived in part (2) for $u(x,y)$ satisfies the PDE and BCs you stated in part (1).
\end{exercise}

\begin{exercise}
\label{exer:9.7}
Suppose $W = (W_t^1,W_t^2,\ldots,W_t^d)_{t \geq 0}$ is a $d$-dimensional Brownian motion.  Define
\begin{align}
R_t
	&=	\Big( \sum_{i=1}^d (W_t^i)^2 \Big)^{1/2} .
\end{align}
Clearly, $R$ lives in an interval $I$ with endpoints $\{0\}$ and $\{\infty\}$.  Show that, when $d=1$ the origin is a regular endpoint.  Show that when $d \geq 2$, the origin is an entrance point.
\end{exercise}

\begin{exercise}
\label{exer:9.8}
Fix a probability space $(\Om, \Fc, \Pbt)$ and let $\Fb = (\Fc_t)_{0 \leq t \leq T}$ be a filtration for a Brownian motion $\Wt = (\Wt_t)_{0 \leq t \leq T}$.  Let $S = (S_t)_{0 \leq t \leq T}$ satisfy
\begin{align}
\dd S_t
	&=	\sig S_t \dd \Wt_t , &
\sig
	&> 0 .
\end{align}
(In finance, this is the Black-Scholes model for a stock $S$ under the risk-neutral measure $\Pbt$ when the risk-free rate of interest is zero.)
Fix constants $p, q \in \Rb$ and define
\begin{align}
u(t,S_t)
	&:=	\Ebt( p S_T^q | \Fc_t ) .
\end{align}
(In finance, $u(t,S_t)$ is the value at time $t$ of a \textit{power option}, which pays $p S_T^q$ at time $T$.)\\
(a) Find the function $u$ by solving the SDE for $S$ and computing the expectation above directly.  \textit{Hint}: once you solve the SDE for $S$, use that fact that $\Wt_T - \Wt_t \sim \Nc( 0 ,T-t)$ under $\Pbt$. It may be helpful to use
\begin{align}
\Ebt \ee^{a Z}
	&=	\ee^{\mu a + \tfrac{1}{2} \sig^2 a^2} , &
Z
	&\sim \Nc( \mu , \sig^2 ) .
\end{align}
(b) Derive a PDE and terminal condition for the function $u$.  Show that the expression you obtained for $u$ in part (a) satisfies this PDE and terminal condition.
\end{exercise}

\begin{exercise}
\label{exer:9.9}
Fix a probability space $(\Om,\Fc,\Pb)$ and let $\Fb = (\Fc_t)_{t \geq 0}$ be a filtration for a Brownian motion $W$.
Let $X = (X_t)_{t \geq 0}$ be the unique strong solution of the following SDE
\begin{align}
\dd X_t
	&=	\mu(t,X_t) \dd t + \sig(t,X_t) \dd W_t .
\end{align}
Define
\begin{align}
u(t,X_t)
	&:=	\Eb \Big( \ee^{- \lam (T-t)} h(X_T) + \int_t^T \ee^{- \lam (s-t)} g(X_s) \dd s \Big| \Fc_t \Big) ,
\end{align}
where $\lam \geq 0$ is a constant.
\\[0.5em]
(a) Define $M = (M_t)_{t \geq 0}$ by
\begin{align}
M_t
	&=	\ee^{- \lam t} u(t,X_t) + \int_0^t \ee^{- \lam s} g(X_s) \dd s .
\end{align}
Show that $M$ is a martingale.
\\[0.5em]
(b) Derive a partial differential equation (PDE) and terminal condition (TC) for the function $u$.
\end{exercise}

\begin{exercise}
\label{exer:9.10}
(a) Suppose $g,h : \Rb_+ \to \Rb$ are deterministic functions.  Define
\begin{align}
X
	&:=	\int_0^T g(t) \dd W_t , &
Y
	&:=	\int_0^T h(t) \dd W_t .
\end{align}
Compute the mean vector $\mu$ and the covariance matrix $\Sig$ for $(X, Y)$.
\\[0.5em]
(b) Consider the process $(X,Y) = (X_t,Y_t)_{t \geq 0}$, defined by the following SDE
\begin{align}
\dd X_t
	&=	\dd W_t , &
\dd Y_t
	&=	X_t \dd g(t) ,
\end{align}
where $g: \Rb_+ \to \Rb$ is a deterministic function with a continuous first derivative.  Solve this SDE.  That is, give explicit expressions for $X_T$ and $Y_T$ in terms of $X_t$, $Y_t$ and some functional of the Brownian motion $W = (W_t)_{t \geq 0}$.  Your expressions for $X_T$ and $Y_T$ should \textit{not} involve $(X_s,Y_s)_{t < s \leq T}$.
\\[0.5em]
(c) Define the marginal transition density of the $(X,Y)$ process
\begin{align}
\Gam(t,x,y;T,\eta)\dd \eta
	&:= \Pb( Y_T \in \dd \eta | X_t = x, Y_t = y) .
\end{align}
Write the Kolmogorov backward equation for $\Gam$ and give an explicit expression for $\Gam$.
\end{exercise}

\begin{exercise}
\label{exer:9.11}
Fix a probability space $(\Om, \Fc, \Pb)$ and let $\Fb = (\Fc_t)_{t \geq 0}$ be a filtration for a scalar Brownian motion $W = (W_t)_{t \geq 0}$.
Consider a scalar diffusion $X = (X_t)_{t \geq 0}$, which is the unique strong solution of the following SDE
\begin{align}
\dd X_t
	&=	\mu(X_t) \dd t + \sig(X_t) \dd W_t .
\end{align}
Assume that the state space of $X$ is $\Rb$. 
Let $Y = (Y_t)_{t \geq 0}$ be given by $\dd Y_t = \dd t$.
Define
\begin{align}
p(x,y;t) \dd t 
	&:= \Pb( Y_\tau \in \dd t | X_0 = x, Y_0 = y ) , &
\tau
	&:= \inf \{ t \geq 0 : X_t = m \} .
\end{align}
(a) Derive a PDE and boundary condition for $p(x,y;t)$.
\\[0.5em]
(b) Suppose $\mu(x) = 0$ and $\sig(x) = 1$.  Find $p(x,y;t)$.
\end{exercise}

\begin{exercise}
\label{exer:9.12}
Let $X$ and $Y$ be given by
\begin{align}
\dd X_t
	&=	\dd W_t , &
\dd Y_t
	&=	\dd B_t , &
0
	&<	a <  \sqrt{ X_0^2 + Y_0^2 } < b < \infty ,
\end{align}
where $W$ and $B$ are independent Brownian motions.  Define $R$ and $\Phi$ as follows
\begin{align}
R_t
	&:=	\sqrt{ X_t^2 + Y_t^2 } , &
\Phi_t
	&:=	\tan^{-1} (Y_t/X_t) .
\end{align}
(a) Find an expression for $\dd R_t$ and $\dd \Phi_t$ that does not involve $X_t$ or $Y_t$.
\\[0.5em]
(b) Define the first exit time $\tau$ as follows
\begin{align}
\tau
	&:=	\inf\{ t \geq 0 : R_t \notin (a,b) \} ,
\end{align}
Note that $\Pb( \tau < \infty) = 1$.  Find $\Pb( R_\tau = a | R_0 = r, \Phi_0 = \phi)$.
\\[0.5em]
(c) Prove the following
\begin{align}
\Eb_{r,\phi} f(R_\tau,\Phi_\tau)
	&=	f(r,\phi) + \Eb_{r,\phi} \int_0^\tau \Big( \frac{1}{2 r} \d_r + \frac{1}{2} \d_r^2 + \frac{1}{2r^2} \d_\phi^2 \Big) f(R_t,\Phi_t) \dd t . \label{eq:exer.9.12.dynkin}
\end{align}
where we have introduced the short-hand notation $\Eb_{r,\phi} \, \cdot \, = \Eb( \, \cdot \, | R_0 = r, \Phi_0 = \phi )$.
\\[0.5em]
(d) Use \eqref{eq:exer.9.12.dynkin} to find $\Eb_{r,\phi} \, \tau$.
\end{exercise}

%
%

\chapter{Jump diffusions}
Notes from this chapter are taken primarily from \cite[Chapter 1]{oksendal2}.
Notes for Section \ref{sec:hawkes} on Hawkes processes follow \cite{hawkes}.


\section{Basic definitions and results on L\'evy processes}

\begin{definition}
\label{def:levy}
A $d$-dimensional stochastic process $\eta = (\eta_t)_{t \geq 0}$, defined on a probability space $(\Om,\Fc,\Pb)$, is called a \emph{L\'evy process} if it satisfies
\begin{enumerate*}
\item $\eta_0 = 0$. \label{itm:0}
\item Independent increments: for any $0 \leq t_0 < t_1 < \ldots < t_n < \infty$, the increments $\eta_{t_1} - \eta_{t_0}, \eta_{t_2} - \eta_{t_1}, \ldots, \eta_{t_n} - \eta_{t_{n-1}}$ are independent. \label{itm:ind}
\item Stationary increments: for any $0 \leq t_1 < t_2 < \infty$, we have $\eta_{t_2} - \eta_{t_1} \sim \eta_{t_2-t_1}$. \label{itm:station}
\item Continuity in probability: for any $\eps > 0$ and $t \geq 0$, we have $\lim_{s \searrow 0} \Pb( |\eta_{t+s}-\eta_t|>\eps ) = 0$. \label{itm:cont}
\end{enumerate*}
\end{definition}

\noindent
Note, Item \ref{itm:cont} in Definition \ref{def:levy} does \emph{not} mean that a L\'evy process cannot jump.  For example, consider a Poisson process $N = (N_t)_{t \geq 0}$ with intensity $\lam$.  A Poisson process is a jump process and yet it is easy to see that it is continuous in probability as, for any $\eps \in (0,1)$, we have
\begin{align}
\Pb( |N_{t+s}-N_t|>\eps ) = \Pb( N_{t+s}-N_t \geq 1 ) = 1 - \Pb(  N_{t+s}-N_t = 0 ) = 1 - \ee^{- \lam s}
	&\to 0&
	&\text{as $s \searrow 0$}.
\end{align}
Item \ref{itm:cont} simply means that, at a fixed $t$, the probability that a L\'evy process has a discontinuity at $t$ is zero. 
\\[0.5em]
We can and do assume that any L\'evy process is right-continuous with left limits (RCLL).  That is
\begin{align}
\lim_{s \searrow 0} \eta_{t+s}
	&=	\eta_t , &
	&\forall \, t \geq 0 .
\end{align}
A process that is RCLL is sometimes called \emph{c\`adl\`ag} (for those who speak French: \emph{continue \`a droite, limite \`a gauche}).
We have already encountered two examples of L\'evy processes.

\begin{example}
A Brownian motion $W = (W_t)_{t \geq 0}$ is a continuous L\'evy process with $W_{t+s}-W_t \sim \Nc(0,s)$.
\end{example}

\begin{example}
A Poisson process $N = (N_t)_{t \geq 0}$ with intensity $\lam$ is a pure jump L\'evy process with $N_{t+s}-N_t \sim \text{Poi}(\lam s)$.
\end{example}

\noindent
A filtration for a L\'evy process is defined just as a filtration for a Brownian motion.

\begin{definition}
\label{def:levy.filtration}
Let $(\Om, \Fc, \Pb)$ be a probability space on which a L\'evy process $\eta=(\eta_t)_{t \geq 0}$ is defined.  A \emph{filtration for the L\'evy process $\eta$} is a collection of $\sig$-algebras $\Fb=(\Fc_t)_{t \geq 0}$ satisfying:
\begin{enumerate*}
\item Information accumulates: if $0 \leq s < t$ then $\Fc_s \subset \Fc_t$.
\item Adaptivity: for all $t \geq 0$, we have $\eta_t \in \Fc_t$.
\item Independence of future increments: if $u>t \geq 0$ then $(\eta_u - \eta_t) \ind \Fc_t$.
\end{enumerate*}
\end{definition}

\noindent
The most natural choice for this filtration $\Fb$ is (not surprisingly) the \emph{natural filtration} for $\eta$, that is, $\Fc_t = \sig(\eta_s,0\leq s \leq t)$.  In principle the filtration $\Fb$ could contain more than the information obtained by observing $\eta$.  However, the information in the filtration is not allowed to destroy the independence of future increments of the L\'evy process.

\begin{definition}
The \emph{jump} of a L\'evy process $\eta$ at time $t \geq 0$ is defined as
\begin{align}
\Del \eta_t
	&:=	\eta_t - \eta_{t-} , &
\eta_{t-}
	&:=	\lim_{s \nearrow t} \eta_s .
\end{align}
\end{definition}

\noindent
Obviously, if $\eta$ does not experience a jump at time $t$, then $\Del \eta_t = 0$.

\begin{definition}
We define $N: \Rb_+ \times \Bc_0^d \times \Om \to \Nb_0$ the \emph{Poisson random measure of $\eta$} by
\begin{align}
N(t,U,\om)
	&=	\sum_{s:0 < s \leq t} \Ib_{\{ \Del \eta_s(\om) \in U \}} , &
U
	&\in \Bc_0^d := \{ B \in \Bc(\Rb^d) : 0 \notin \overline{B} \} ,
\end{align}
where $\overline{B}$ denotes the closure of $B$.  That is, the sets $U \in \Bc_0^d$ are bounded away from the origin.
\end{definition}

\noindent
As with most random variables and stochastic processes, we will typically omit the dependence of $N$ on $\om$, writing simply $N(t,U)$ as opposed to $N(t,U,\om)$.
The Poisson random measure $N(t,U)$ counts the number of jumps of size $\Del \eta_s \in U$ prior to time $t$.  It will be convenient to introduce the following differential form
\begin{align}
N(\dd t, \dd z)
	&=	\sum_{s>0} \del_{s}(\dd t) \del_{\Del \eta_s}(\dd z) , &
	&\text{so that}&
N(t,U)
	&=	\int_0^t \int_U N(\dd s, \dd z) .
\end{align}
Roughly speaking, $N(\dd t, \dd z)$ counts the number of jumps of size $\dd z$ over the time interval $\dd t$.

\begin{definition}
Let $N$ be the Poisson random measure of a L\'evy process $\eta$.  We define $\nu : \Bc_0^d \to \Rb_+$, the \emph{L\'evy measure of $\eta$}, as follows
\begin{align}
\nu(U)
	&:=	\Eb N(1,U) , &
U
	&\in \Bc_0^d .
\end{align}
\end{definition}

\begin{theorem}
\label{thm:poisson}
Let $U \in \Bc_0^d$.  Then the process $(N(t,U))_{t \geq 0}$ is a Poisson process with intensity $\nu(U)$.
\end{theorem}

\noindent
Because $N(t,U)$ is a Poisson process with intensity $\nu(U)$, we have
\begin{align}
\Eb N(t,U)
	&= \nu(U) t , &
	&\Rightarrow&
\frac{1}{t} \Eb N(t,U)
	&=	\nu(U) .
\end{align}
Thus, $\nu(U)$ is the expected number of times $\eta$ has a jump of size $z \in U$ per unit time.
Note that
\begin{align}
\nu(U) \dd t = \Eb N(\dd t, U)
	&=	\sum_{k=0}^\infty k \Pb( N(t,U) = k )
	=		\Pb( N( \dd t, U) = 1 ) + \Oc( \dd t^2 ) , &
	&\text{as $\dd t \to 0$}.
\end{align}
Hence, for small $\dd t$, one can think of $\nu( U) \dd t$ as the probability that $\eta$ experiences a jump of size $z \in U$ in the time interval $\dd t$.

\begin{example}[Compound Poisson process]
Let $(X_n)_{n \in \Nb}$ be a sequence of iid random vectors in $\Rb^d$ with distribution $F_X$.  Let $(P_t)_{t \geq 0}$ be a Poisson process with intensity $\lam$. Assume $(P_t)_{t \geq 0} \ind (X_n)_{n \in \Nb}$.  We define a \emph{compound Poisson process} $\eta = (\eta_t)_{t \geq 0}$ by
\begin{align}
\eta_t
	&=	\sum_{n=1}^{P_t} X_n .
\end{align}
The increments of $\eta$ are given by
\begin{align}
\eta_t - \eta_s
	&=	\sum_{n = P_s + 1}^{P_t} X_n .
\end{align}
The distribution $F_{\eta_t - \eta_s}$ depends only on $(t-s)$ and $F_X$.  As such we see that $\eta$ has stationary increments.  Also, non-overlapping increments of $\eta$ are clearly independent, as they depend on different $(X_i)$.  Thus, $\eta$ is a L\'evy process in $\Rb^d$.  Let us find the L\'evy measure $\nu$ corresponding to $\eta$.  For any $U \in \Bc_0^d$ we have
\begin{align}
\nu(U)
	&=	\Eb N(1,U)
	=		\Eb \sum_{s:0 < s \leq 1} \Ib_{\{ \Del \eta_s \in U \}} \\
	&=	\Eb \sum_{n=1}^{P_1} \Ib_{ \{X_n \in U \} } 
	=		\Eb \sum_{n=1}^{P_1} \Eb[ \Ib_{ \{X_n \in U \} } | P_1 ] \\
	&=	\Eb \sum_{n=1}^{P_1} \Eb[ \Ib_{ \{X_n \in U \} } ]
	=		\Eb \sum_{n=1}^{P_1} \Pb( X_n \in U ) \\
	&=	F_X(U) \Eb P_1
	=		\lam F_X(U) .
\end{align}
Thus, the L\'evy measure of $\eta$ is $\nu = \lam F_X$.
\end{example}

\noindent
A pure-jump L\'evy process $\eta$ has a finite L\'evy measure $\nu(\Rb^d) < \infty$ if and only if it can be represented by a compound Poisson process.  In this case, we can express $\eta$ in one of two ways
\begin{align}
\eta_t
	&=	\int_{\Rb^d} z N(t ,\dd z) , &
	&\text{or}&
\eta_t	
	&=	\sum_{n=1}^{P_t} X_n . \label{eq:two.ways}
\end{align}
However, there exist L\'evy processes  for which $\nu(\Rb^d) = \infty$.  We call a L\'evy process for which $\nu(\Rb^d) = \infty$ an \emph{infinite activity} L\'evy process.  For an infinite activity L\'evy process, neither of the representations in \eqref{eq:two.ways} make sense.  
To write the most general form of a L\'evy process, we introduce the compensated Poisson random measure.

\begin{definition}
Let $N$ be a Poisson random measure with associated L\'evy measure $\nu$.  The \emph{compensated Poisson random measure}, denoted $\Nt$ is defined as
\begin{align}
\Nt( t , A) 
	&:=	N( t , A) - \nu(A) t .
\end{align} 
In differential form, we have
\begin{align}
\Nt( \dd t , \dd z) 
	&:=	N( \dd t , \dd z) - \nu( \dd z) \dd t .
\end{align}
\end{definition}

\noindent
For any $A \in \Bc_0^d$ the process $(\Nt(t,A))_{t \geq 0}$ is a martingale (with respect to its own filtration) as
\begin{align}
\Eb [ \Nt(t,A) | \Fc_s ]
	&=	\Nt(s,A) + \Eb [ \Nt(t,A) - \Nt(s,A) | \Fc_s ] \\
	&=	\Nt(s,A) + \Eb [ N(t,A) - N(s,A) | \Fc_s] - \nu(A) (t - s) \\
	&=	\Nt(s,A) + \nu(A) (t-s) - \nu(A) (t-s)
	=		\Nt(s,A) ,
\end{align}
where we have used $\Eb N(t,A) = t \nu(A)$.
It follows that, for any fixed $R>0$, the process $M^{(k)}=(M_t^{(k)})_{t \geq 0}$, defined as
\begin{align}
M_t^{(k)}
	&:= \int_{1/k \leq |z| < R} z \Nt(t,\dd z) , &
k
	&=	1, 2, \ldots ,
\end{align}
is a ($d$-dimensional) martingale satisfying $\Eb |M_t^{(k)}|^2 < \infty$.  One can show that, as $k \to \infty$ the sequence of processes $(M^{(k)})_{k \in \Nb}$ converges in $L^2(\Om,\Fc,\Pb)$ to a process $M = (M_t)_{t \geq 0}$ defined as
\begin{align}
M_t
	&\equiv	\int_{ |z| < R} z \Nt(t,\dd z)
	:=	\lim_{k \to \infty} \int_{1/k \leq |z| < R} z \Nt(t,\dd z) ,
\end{align}
which is also a martingale and satisfies $\Eb |M_t|^2 < \infty$.

\begin{remark}
In general we cannot separate an integral with respect to $\Nt(t,\dd z)$ into an integral with respect to $N(t,\dd z)$ and an integral with respect to $t \nu(\dd z)$.  The reason is that we may have $\int_{|z|<R} |z| \nu(\dd z) = \infty$, and in this case, we have
\begin{align}
\int_{|z|<R} z \Nt(t,\dd z)
	&\neq \int_{|z|<R} z N(t,\dd z) - \int_{|z|<R} z \nu(\dd z) t .
\end{align}
\end{remark}

\noindent
The following theorem gives the most general form of a L\'evy process.

\begin{theorem}[It\^o-L\'evy decomposition]
Let $\eta = (\eta_t)_{t \geq 0}$ be a L\'evy process in $\Rb^d$.  Then $\eta$ has the decomposition
\begin{align}
\eta_t
	&=	\mu_R t + \sig W_t + \int_{|z| < R} z \Nt( t , \dd z ) + \int_{|z| \geq R} z N( t, \dd z) , \label{eq:eta.general}
\end{align}
for some vector $\mu_R \in \Rb^d$, matrix $\sig \in \Rb^{d \times d}$, constant $R \in [0,\infty]$,
Poisson random measure $N$ on $\Rb^d$, and $d$-dimensional Brownian motion $W$ that is independent of $N$.
The choice $R = 0$ is possible if and only if $\int_{|z| < 1} |z| \nu(\dd z) < \infty$, and the choice $R = \infty$ is possible if and only if $\int_{|z| \geq 1} |z| \nu(\dd z) < \infty$, where $\nu$ is the L\'evy measure associated with $N$.
\end{theorem}

\noindent
We can always choose $R = 1$ in \eqref{eq:eta.general}.  It is useful to recognize when we can choose $R=0$ and $R=\infty$ as, in these cases, the 
right-hand side of \eqref{eq:eta.general} reduces from four to three terms.
\begin{example}
If the L\'evy measure $\nu$ satisfies
\begin{align}
\int_{|z| \geq 1}  |z| \nu( \dd z)
	&< \infty , \label{eq:large.levy}
\end{align}
then we can write \eqref{eq:eta.general} as
\begin{align}
\eta_t
	&=	\Big( \mu_1 + \int_{|z| \geq 1} z \nu(\dd z) \Big)t + \sig W_t + \int_{|z| < 1} z \Nt( t , \dd z ) + \int_{|z| \geq 1} z \( N( t, \dd z)  - \nu(\dd z) t \) \\
	&=	\mu_\infty t + \sig W_t +  \int_{\Rb^d}  z \Nt( t , \dd z ) , \label{eq:levy.compensated}
\end{align} 
where $\mu_\infty = \mu_1 + \int_{|z| \geq 1} z \nu(\dd z)$.
We call $\mu_\infty$ the \emph{center} of $\eta$ since $\Eb \eta_t = \mu_\infty t$.
It is straightforward to show that if $\mu_\infty = 0$ then $\eta$ is a martingale.
\end{example}

\begin{example}
If the L\'evy measure $\nu$ satisfies
\begin{align}
\int_{|z|<1} |z| \nu( \dd z)
	&< \infty , \label{eq:small.levy}
\end{align}
then we can write \eqref{eq:eta.general} as
\begin{align}
\eta_t
	&=	\mu_1 t + \sig W_t + \int_{|z| < 1} z N( t , \dd z ) - \int_{|z| < 1} z \nu( \dd z ) t + \int_{|z| \geq 1} z N( t, \dd z) \\
	&=	\mu_0 t + \sig W_t + \int_{\Rb^d} z N(t, \dd z) , \label{eq:levy.uncompensated}
\end{align}
where $\mu_0 = \mu_1 - \int_{|z| < 1} z \nu( \dd z )$.  We call $\mu_0$ the \emph{drift} of $\eta$.
\end{example}

\begin{example}
If the L\'evy measure $\nu$ satisfies
\begin{align}
\nu(\Rb^d) 
	&< \infty ,
\end{align}
then \eqref{eq:small.levy} holds, and we can write \eqref{eq:levy.uncompensated}  as a compound Poisson process
\begin{align}
\eta_t
	&=	\mu_0 t + \sig W_t + \sum_{n=1}^{P_t} X_n ,
\end{align}
where $P$ is a Poisson process with intensity $\nu(\Rb^d)$ and the $(X_n)_{n \in \Nb}$ are iid random vectors in $\Rb^d$ with common distribution $F_X = \nu / \nu(\Rb^d)$.
\end{example}
\noindent
We will typically drop the subscript from $\mu_R$, writing it instead simply as $\mu$.


\begin{theorem}
L\'evy processes are Markov processes.
\end{theorem}

\noindent
This should come as no surprise.  We have already seen a Brownian motion and a Poisson process are Markov processes.  As these processes serve as building blocks for more general L\'evy processes, it follows that L\'evy processes are Markov as well.

\begin{theorem}[L\'evy-Khintchine formula]
Let $\eta = (\eta_t)_{t \geq 0}$ be the following $d$-dimensional L\'evy process
\begin{align}
\eta_t
	&=	\mu t + \sig W_t + \int_{|z| < R} z \Nt( t , \dd z ) + \int_{|z| \geq R} z N( t, \dd z) ,
\end{align}
where $\mu \in \Rb^d$, $\sig \in \Rb^{d \times d}$, $W$ is a $d$-dimensional Brownian motion and $N$ is a Poisson random measure on $\Rb^d$.
Then the L\'evy measure $\nu$ associated with $N$ satisfies
\begin{align}
\int_{\Rb^d} ( 1 \wedge |z|^2 ) \nu(\dd z)
	&<	\infty , \label{eq:levy.meas.cond}
\end{align}
and the characteristic function $\phi_{\eta_t} : \Rb^d \to \Cb$ is given by
\begin{align}
\phi_{\eta_t}(u) := \Eb \ee^{\ii \< u , \eta_t \> }
	&=	\ee^{t \psi(u)} , &
\< a , b \>
	&:=	\sum_{i=1}^d a_i b_i . \label{eq:levy.char}
\end{align}
where $\psi$, called the \emph{characteristic exponent of $\eta$}, is given by
\begin{align}
\psi(u)
	&=	\ii \< \mu , u \>  - \tfrac{1}{2} \< u , \sig \sig^\textup{T} u \> 
			+ \int_{|z| < R} (\ee^{\ii \< u , z\> } - 1 - \ii \< u , z\> ) \nu( \dd z ) + \int_{|z| \geq R} (\ee^{\ii \< u , z \>} - 1 ) \nu( \dd z ) . 
	\label{eq:levy.exp} 
\end{align}
Conversely, given triplet $(\mu, a, \nu)$ with $a = \sig \sig^\textup{T}$ and $\nu$ satisfying \eqref{eq:levy.meas.cond},
there exists a L\'evy process $\eta$ satisfying \eqref{eq:levy.char}-\eqref{eq:levy.exp}.
\end{theorem}

\begin{proof}
We will show why \eqref{eq:levy.char}-\eqref{eq:levy.exp} holds for a scalar L\'evy process of the compound Poisson type.  In this case, $\eta$ has the decomposition
\begin{align}
\eta_t
	&=	\mu t + \sig W_t + \sum_{n = 1}^{P_t} X_n ,
\end{align}
where $W$ is a scalar Brownian motion, $P$ is a Poisson process with intensity $\lam$ and the $(X_n)_{n \in \Nb}$ are iid random variables with common distribution $F_X$.  We compute
\begin{align}
\Eb \ee^{\ii u \eta_t}
	&=	\Eb \ee^{ \ii u \mu t + \ii u \sig W_t + \ii u \sum_{n = 1}^{P_t} X_n} \\
	&=	\ee^{\ii u \mu t} \Eb \ee^{ \ii u \sig W_t } \Eb \ee^{ \ii u \sum_{n = 1}^{P_t} X_n} \\
	&=	\ee^{\ii u \mu t - \tfrac{1}{2} \sig^2 u^2 t} \Eb \Eb[ \ee^{\ii u \sum_{n = 1}^{P_t} X_n} | P_t ] \\
	&=	\ee^{\ii u \mu t - \tfrac{1}{2} \sig^2 u^2 t} \Eb ( \Eb[ \ee^{\ii u X } ] )^{P_t} \\
	&=	\ee^{\ii u \mu t - \tfrac{1}{2} \sig^2 u^2 t} \Eb \Big( \int \ee^{\ii u z} F_X(\dd z) \Big)^{P_t} \\
	&=	\ee^{\ii u \mu t - \tfrac{1}{2} \sig^2 u^2 t} \exp \Big( \lam t \big( \int \ee^{\ii u z} F_X(\dd z) - 1 \big) \Big) \\
	&=	\ee^{\ii u \mu t - \tfrac{1}{2} \sig^2 u^2 t} \exp \Big( \lam t \int ( \ee^{\ii u z} - 1 ) F_X(\dd z) \Big) \\
	&=	\ee^{t \psi(u)} ,
\end{align}
where we have used $\Eb s^{P_t} = \exp( \lam t ( s - 1 ) )$ and $\nu = \lam F_X$.
The computation for a more general L\'evy process in multiple dimensions is similar.
\end{proof}

\noindent
Let us provide a few examples of L\'evy measures $\nu$ and compute the integrals that appear in expression \eqref{eq:levy.exp} for the characteristic exponent $\psi$.

\begin{example}[Dirac comb]
A Dirac comb L\'evy measure on $\Rb$ is a measure of the form
\begin{align}
\nu( \dd z)
	&=	\sum_{i=1}^n \lam_i \del_{z_i}(\dd z) ,
\end{align}
where $\lam_i$ is the intensity of a jump of size $z_i$.  In this case, both \eqref{eq:large.levy} and \eqref{eq:small.levy} are satisfied, so we can choose either $R=0$ or $R = \infty$.  Suppose we choose $R = 0$.  Then the third term in \eqref{eq:levy.exp} disappears and the last term becomes
\begin{align}
\int_\Rb (\ee^{\ii u z} - 1 ) \nu(\dd z)
	&=	\sum_{i=1}^n \lam_i (\ee^{\ii u z_i} - 1 ) .
\end{align}
\end{example}

\begin{example}[Gaussian Jumps]
A Gaussian L\'evy measure on $\Rb$ is a measure of the form
\begin{align}
\nu( \dd z)
	&=	\lam \frac{1}{\sqrt{2 \pi s^2}} \exp \Big( \frac{-(z-m)^2}{2 s^2} \Big) \dd z ,
\end{align}
where $m$ is the mean jump size, $s^2$ is the variance of the jumps, and the intensity of all jumps is $\nu(\Rb) = \lam$.  In this case, both \eqref{eq:large.levy} and \eqref{eq:small.levy} are satisfied, so we can choose either $R=0$ or $R = \infty$.  Suppose we choose $R = 0$.  Then the third term in \eqref{eq:levy.exp} disappears and the last term becomes
\begin{align}
\int_\Rb (\ee^{\ii u z} - 1 ) \nu(\dd z)
	&=	\lam \Big( \ee^{\ii u m - \tfrac{1}{2} s^2 u^2 } - 1 \Big) .
\end{align}
\end{example}

\begin{example}[Generalized tempered stable jumps]
A Generalized tempered stable L\'evy measure on $\Rb$ is a measure of the  form
\begin{align}
\nu( \dd z)
	&=	\Ib_{\{ z < 0 \}} \frac{A_- \ee^{-\beta_- |z|}}{|z|^{1+\alpha_-}} + \Ib_{\{ z > 0 \}} \frac{A_+ \ee^{-\beta_+ z}}{|z|^{1+\alpha_+}} ,
\end{align}
where the parameters $(A_\pm,\alpha_\pm,\beta_\pm)$ are all positive.  In this case \eqref{eq:large.levy} is satisfied by \eqref{eq:small.levy} is not.  Thus, we must choose $R = \infty$.  The fourth term in \eqref{eq:levy.exp} then disappears and the third term becomes
\begin{align}
\int_\Rb (\ee^{\ii u z} - 1 - \ii u z) \nu(\dd z)
	&=	A_-  \Gam(-\alpha_-) \beta_-^{\alpha_-} \Big[ \Big( 1 + \frac{\ii u}{\beta_-} \Big)^{\alpha_-} - 1 - \frac{\ii u \alpha_-}{\beta_-} \Big] \\ &\quad
			+ A_+  \Gam(-\alpha_+) \beta_+^{\alpha_+} \Big[ \Big( 1 - \frac{\ii u}{\beta_+} \Big)^{\alpha_+} - 1 + \frac{\ii u \alpha_+}{\beta_+} \Big] .
\end{align}
assuming $\alpha_\pm \neq 0$ and $\alpha_\pm \neq 1$.  The above integral can also be computed explicitly for the cases $\alpha_\pm = 0$ and $\alpha_\pm = 1$ (see, e.g., \cite[Proposition 4.2]{contbook}).
\end{example}

\begin{assumption}
From this point onward we assume $\Eb |\eta_t|^2 < \infty$, and thus $\Eb |\eta_t| < \infty$.  Thus, we can take $R = \infty$ and express $\eta$ as follows
\begin{align}
\eta_t
	&=	\mu t + \sig W_t + \int_{\Rb^d} z \Nt(t,\dd z) . \label{eq:eta.form}
\end{align}
\end{assumption}

\noindent
L\'evy process are \emph{semimartingales}, which (roughly speaking) form the class of stochastic processes that can be used to construct an It\^o integral.  In Section \ref{sec:ito.integral} we defined It\^o integral of an $(\Fc_t)_{t \geq 0}$-adapted process $\Gam = (\Gam_t)_{t \geq 0}$ with respect to a Brownian motion $W = (W_t)_{t \geq 0}$ by introducing a sequence of simple processes $\Gam^{(n)} = (\Gam_t^{(n)})_{t \geq 0}$ that converged to $\Gam$ in the following sense
\begin{align}
\Gam_t^{(n)}
	&:=	\sum_{j=0}^{n-1} \Gam_{t_j} \Ib_{\{ t_j \leq t < t_{j+1} \}} , &
0
	&\leq t_0 < t_1 < \ldots < t_n = T , &
\lim_{n \to \infty} \Eb \int_0^T | \Gam_t - \Gam_t^{(n)} |^2 \dd t
	&=	0 .
\end{align}
We then by defined the It\^o integral as the following limit
\begin{align}
\int_0^T \Gam_t \dd W_t
	&:=	\lim_{n \to \infty} \int_0^T \Gam_t^{(n)} \dd W_t
	=		\lim_{n \to \infty} \sum_{j=0}^{n-1} \Gam_{t_j} (W_{t_{j+1}} - W_{t_j}) .
\end{align}
We can define the It\^o integral of $\Gam$ with respect to a L\'evy process $\eta$ in the same manner
\begin{align}
\int_0^T \Gam_t \dd \eta_t
	&:=	\lim_{n \to \infty} \int_0^T \Gam_t^{(n)} \dd \eta_t
	=		\lim_{n \to \infty} \sum_{j=0}^{n-1} \Gam_{t_j} (\eta_{t_{j+1}} - \eta_{t_j}) .
\end{align}
Here, we must add a technical condition that the integrand $\Gam$ be a \emph{left-continuous} process.  
If a process $\Gam$ is not left-continuous, then we can still define a stochastic integral as follows
\begin{align}
I_T
	&= \int_0^T \Gam_{t-} \dd \eta_t , &
\Gam_{t-}
	&=	\lim_{s \nearrow t} \Gam_s , 
\end{align}
where $(\Gam_{t-})_{t \geq 0}$ is the left-continuous version of $(\Gam_t)_{t \geq 0}$.  
Note that, if $\Gam$ is left-continuous then $\Gam_t = \Gam_{t-}$ for all $t \geq 0$.
The process $(I_t)_{t \geq 0}$ will still be a right-continuous process as $\eta$ is right-continuous.
In view of \eqref{eq:eta.form},
for a left-continuous process $\Gam$, we can separate an It\^o integral into three terms
\begin{align}
\int_0^T \Gam_{t} \dd \eta_t
	&= 	\int_0^T \Gam_t \mu \dd t + \int_0^T \Gam_t \sig \dd W_t + \int_0^T \int \Gam_{t} z \Nt( \dd t , \dd z ) . \label{eq:I.again}
\end{align}
Expression \eqref{eq:I.again} suggests that we consider general processes of the form
\begin{align}
\dd X_t
	&=	\mu_t \dd t + \sig_t \dd W_t + \int \gam_{t}(z) \Nt( \dd t , \dd z ) , \label{eq:levy.ito}
\end{align}
where the processes $(\mu_t)_{t \geq 0}$, $(\sig_t)_{t \geq 0}$ and $(\gam_{t-}(z))_{t \geq 0}$, must be adapted to a filtration $\Fb=(\Fc_t)_{t \geq 0}$ obtained by observing the processes $W$ and $\int z N( \cdot ,\dd z)$.  Note that we have written things in \emph{differential form}; to obtain the integral form, simply integrate and add an initial condition.  We have been somewhat vague about the dimension of the various objects appearing in \eqref{eq:levy.ito}.  In general, for $X \in \Rb^d$, we could have
\begin{align}
\mu
	&: \Rb_+ \times \Om \to  \Rb^d , &
\sig
	&: \Rb_+ \times \Om \to \Rb^{d \times n} , &
\gam
	&: \Rb_+ \times \Rb^m \times \Om \to \Rb^{d \times k} .
\end{align}
where $W$ is an $n$-dimensional Brownian motion (with $W^{(i)} \ind W^{(j)}$ for $i \neq j$) and $\Nt$ is a $k$-dimensional Poisson random measure (with $N^{(i)} \ind N^{(j)}$ for $i \neq j$) on $\Rb^m$.  Component-wise, we have
\begin{align}
\dd X_t^{(i)}
	&=	\mu_t^{(i)} \dd t + \sum_{j=1}^n \sig_t^{(i,j)} \dd W_t^{(j)} + \sum_{j=1}^k \int_{\Rb^m} \gam_{t-}^{(i,j)}(z) \Nt^{(j)}( \dd t , \dd z ) , &
i
	&= 1, 2, \ldots, d .
\end{align}
We call any process of the form \eqref{eq:levy.ito} a \emph{L\'evy-It\^o process}.  
\\[0.5em]
We showed (for the scalar case) in Theorem \ref{thm:ito.props} that if $\sig$ satisfies $\Eb \int_0^T |\sig_t|^2 \dd t < \infty$ then the It\^o integral $I_t := \int_0^t \sig_s \dd W_s$ is a martingale.  Likewise, if
\begin{align}
\Eb \int_0^T \int_\Rb |\gam_t(z)|^2 \nu( \dd z) \dd t
	&< \infty  ,
\end{align}
then the process $M = (M_t)_{0 \leq t \leq T}$, defined by
\begin{align}
M_t
	&:=	\int_0^t \int_\Rb \gam_s(z) \Nt(\dd s,\dd z) ,
\end{align}
is a martingale.  It follows that if $\mu = 0$ in \eqref{eq:levy.ito}, and $\sig$ and $\gam(z)$ satisfy the above square integrability conditions, then the process $X = (X_t)_{t \geq 0}$ is a martingale.


\section{The It\^o formula for L\'evy-It\^o processes.}
In this section, we suppose $X$ is given by \eqref{eq:levy.ito} and we show how to compute $\dd f(X_t)$.  We simply state the result and then explain heuristically why it is true.
To keep the statements of the results in this chapter readable, we make the following standing assumption.

\begin{assumption}[Integrability]
\label{ass:integrability}
Whenever we write $\int_0^T H_t \dd W_t$ or $\int_0^T \int H_t(z) \Nt(\dd t,\dd z)$, we assume that
\begin{align}
\Eb \int_0^T H_t^2 \dd t
	&<	\infty , &
	&\text{respectively}, &
\Eb \int_0^T \int H_t^2(z) \nu(\dd z) \dd t
	&<	\infty ,
\end{align}
for every $T > 0$.  Whenever we write $\int_0^T H_t \dd t$ or $\int_0^T \int H_t(z) \nu(\dd z) \dd t$, we assume the integral is finite almost surely.  In particular, under this assumption, all stochastic integrals with respect to $W$ and $\Nt$ are martingales.
\end{assumption}

\noindent
For example, in the It\^o formula below, Assumption \ref{ass:integrability} requires that each of the integrands appearing in \eqref{eq:df} (e.g., $\sig_t f'(X_t)$ and $f(X_{t-} + \gam_t(z)) - f(X_{t-})$) satisfy the relevant integrability condition.

\begin{theorem}[One-dimensional It\^o formula for a L\'evy-It\^o processes]
\label{thm:ito.formula.levy}
Suppose Assumption \ref{ass:integrability} holds and $X = (X_t)_{t \geq 0}$ is a one-dimensional L\'evy-It\^o process of the form
\begin{align}
\dd X_t
	&=	\mu_t \dd t + \sig_t \dd W_t + \int_\Rb \gam_{t-}(z) \Nt( \dd t, \dd z) . \label{eq:dX.thm}
\end{align}
where the drift $\mu$ volatility $\sig$ and $\gam$ are maps
\begin{align}
\mu
	&: \Rb_+ \times \Om \to \Rb , &
\sig
	&: \Rb_+ \times \Om \to \Rb_+ , &
\gam
	&: \Rb_+ \times \Om \times \Rb \to \Rb ,
\end{align}
and $N$ is a one-dimensional Poisson random measure on $\Rb$, which is independent of a one-dimensional Brownian motion $W$.
Let $f : \Rb \to \Rb$ satisfy $f \in C^2(\Rb)$.  Then
\begin{align}
\dd f(X_t)
	&=	\Big( \mu_t  f'(X_t) + \tfrac{1}{2} \sig_t^2 f''(X_t)  \Big) \dd t 
			+ \sig_t f'(X_t) \dd W_t \\ &\quad
			+	\int_\Rb \Big( f\(X_{t-} + \gam_{t}(z) \) - f\(X_{t-}\) \Big) \Nt( \dd t, \dd z) \\ &\quad
			+	\int_\Rb \Big( f\(X_{t-} + \gam_{t}(z) \) - f\(X_{t-}\) - \gam_{t}(z) f'(X_{t-}) \Big) \nu(\dd z) \dd t . \label{eq:df}
\end{align}
\end{theorem}

\noindent
We are not going to prove Theorem \ref{thm:ito.formula.levy}.  However, we will attempt to understand why the formula is correct.  Suppose for simplicity that $\nu(\Rb)< \infty$.  In this case we have a finite activity L\'evy process and we can separate the compensated Poisson random measure $\Nt$ into two parts  $N(\dd t, \dd z)$ and $\nu(\dd z)\dd t$.  In this case, we can write $\dd X_t$ as
\begin{align}
\dd X_t
	&=	\mu_t \dd t + \sig_t \dd W_t + \int_\Rb \gam_t(z) N( \dd t, \dd z) - \int_\Rb \gam_t(z) \nu( \dd z) \dd t \\
	&=	\Big( \mu_t - \int_\Rb \gam_t(z) \nu( \dd z) \Big)  \dd t + \sig_t \dd W_t + \int_\Rb \gam_t(z) N( \dd t, \dd z) \\
	&=	\mu_t' \dd t + \sig_t \dd W_t + \int_\Rb \gam_t(z) N( \dd t, \dd z) , \\
\mu_t'
	&:=	\mu_t - \int_\Rb \gam_t(z) \nu( \dd z) ,
\end{align}
where we can now identify the drift of $X$ as $\mu_t'$.
Similarly, we can write $\dd f(X_t)$ as
\begin{align}
\dd f(X_t)
	&=	\Big( \mu_t f'(X_t) + \tfrac{1}{2} \sig_t^2 f''(X_t) \Big) \dd t 
			+ \sig_t  f'(X_t) \dd W_t \\ &\quad
			+	\int_\Rb \Big( f\(X_{t-} + \gam_{t}(z) \) - f\(X_{t-}\) \Big) N( \dd t, \dd z) \\ &\quad
			- \int_\Rb \Big( f\(X_{t-} + \gam_{t}(z) \) - f\(X_{t-}\) \Big) \nu( \dd z) \dd t \\ &\quad
			+	\int_\Rb \Big( f\(X_{t-} + \gam_{t}(z) \) - f\(X_{t-}\) - \gam_{t}(z) f'(X_t) \Big) \nu(\dd z) \dd t \\
	&=	\bigg( \Big( \mu_t - \int_\Rb \gam_t(z) \nu(\dd z) \Big) f'(X_t)  + \tfrac{1}{2} \sig_t^2 f''(X_t) \bigg) \dd t 
			+ \sig_t f'(X_t) \dd W_t \\ &\quad
			+	\int_\Rb \Big( f\(X_{t-} + \gam_{t}(z) \) - f\(X_{t-}\) \Big) N( \dd t, \dd z) \\
	&=	\Big( \mu_t' f'(X_t) + \tfrac{1}{2} \sig_t^2 f''(X_t) \Big) \dd t 
			+ \sig_t f'(X_t) \dd W_t \\ &\quad
			+	\int_\Rb \Big( f\(X_{t-} + \gam_{t}(z) \) - f\(X_{t-}\) \Big) N( \dd t, \dd z) . \label{eq:df.2}
\end{align}
Things are looking a bit more familiar now.  The non-integral terms in \eqref{eq:df.2} arise from the $\mu_t' \dd t + \sig_t \dd W_t$ part of $\dd X_t$.  To understand the integral term in \eqref{eq:df.2}, define
\begin{align}
\eta_t
	&:=	\int_\Rb z N(t, \dd z) .
\end{align}
Observe that
\begin{align}
\Del \eta_t
	&=	\int_\Rb z \Del N(t, \dd z) ,  \label{eq:eta-jump}
\end{align}
Hence, we have
\begin{align}
\Del N( t, \dd z )
	&=	\del_{\Del \eta_t}(z) \dd z .
\end{align}
It follows that
\begin{align}
\Del X_t
	&:=	X_t - X_{t-}
	=	\int_\Rb \gam_t(z) \Del N( t, \dd z) 
	=		\int_\Rb \gam_t(z) \del_{\Del \eta_t}(z) \dd z
	=		\gam_t(\Del \eta_t) .
\end{align}
Finally, we have
\begin{align}
\Del f(X_t)
	&=	f(X_t) - f(X_{t-}) \\
	&=	f(X_{t-} + \Del X_t) - f(X_{t-}) \\
	&= 	f(X_{t-} + \gam_t({\Del \eta_t}) ) - f(X_{t-}) \\
	&=	\int_\Rb \Big( f(X_{t-} + \gam_t(z) ) - f(X_{t-}) \Big) \del_{\Del \eta_t}(z) \dd z \\
	&=	\int_\Rb \Big( f(X_{t-} + \gam_t(z) ) - f(X_{t-}) \Big) \Del N(t, \dd z) .
\end{align}
This last expression agrees with the integral term in \eqref{eq:df.2}.  

\begin{example}
Consider the following \emph{Geometric L\'evy process}
\begin{align}
\dd X_t
	&=	b X_t \dd t + a X_t \dd W_t + \int_\Rb \( \ee^{c(z)} - 1 \)  X_{t-}  \Nt( \dd t, \dd z) . \label{eq:dX.ex}
\end{align}
Note, if $X_0 > 0$ then this process remains strictly positive since
\begin{align}
X_{t-} +  \( \ee^{c(z)} - 1 \)  X_{t-}
	&=	\ee^{c(z)} X_{t-} > 0 .
\end{align}
Comparing \eqref{eq:dX.ex} with \eqref{eq:dX.thm} we identify
\begin{align}
\mu_t
	&=	b X_t , &
\sig_t
	&=	a X_t , &
\gam_{t}(z)
	&=	\( \ee^{c(z)} - 1 \) X_{t-} .
\end{align}
We see an explicit expression for $X_t$.  To this end, we define
\begin{align}
f(X_t)
	&=	\log X_t , &
f'(X_t)
	&=	\frac{1}{X_t} , &
f''(X_t)
	&=	- \frac{1}{X_t^2} ,
\end{align}
Plugging all of the above into \eqref{eq:df} we obtain
\begin{align}
\dd f(X_t)
	&=	\dd \log X_t \\
	&=	\Big( b X_t  \frac{1}{X_t} - \tfrac{1}{2} a^2 X_t^2 \frac{1}{X_t^2}  \Big) \dd t 
			+ a X_t \frac{1}{X_t} \dd W_t \\ &\quad
			+	\int_\Rb \Big( \log \(X_{t-} +\( \ee^{c(z)} - 1 \) X_{t-} \) - \log \( X_{t-} \) \Big) \Nt( \dd t, \dd z) \\ &\quad
			+	\int_\Rb \Big( \log \(X_{t-} + \( \ee^{c(z)} - 1 \) X_{t-} \) - \log \( X_{t-} \) - \( \ee^{c(z)} - 1 \) X_t \frac{1}{X_t} \Big) \nu(\dd z) \dd t \\
	&=	\Big( b - \tfrac{1}{2} a^2  \Big) \dd t 
			+ a \dd W_t 
			+	\int_\Rb  c(z) \Nt( \dd t, \dd z) 
			+	\int_\Rb \Big( c(z) - \ee^{c(z)} + 1 \Big) \nu(\dd z) \dd t \\
	&=	b' \dd t 
			+ a \dd W_t 
			+	\int_\Rb c(z) \Nt( \dd t, \dd z) , \\
b'
	&=	b - \tfrac{1}{2} a^2 - \int_\Rb \big( \ee^{c(z)} - 1 - c(z) \big) \nu( \dd z) .
\end{align}
Integrating from $0$ to $t$, we obtain
\begin{align}
\log X_t - \log X_0
	&=	b' t + a W_t +	\int_\Rb c(z) \Nt( t, \dd z) .
\end{align}
Finally, adding $\log X_0$ to both sides and exponentiating, we have
\begin{align}
X_t
	&=	X_0 \exp \( b' t + a W_t +	\int_\Rb c(z) \Nt( t, \dd z) \) .
\end{align}
\end{example}

\begin{theorem}[L\'evy-It\^o Isometry]
\label{thm:isometry}
Suppose Assumption \ref{ass:integrability} holds and let $X = (X_t)_{t \geq 0}$ satisfy
\begin{align}
\dd X_t
	&=	\sig_t \dd W_t + \int_\Rb \gam_t(z) \Nt( \dd t , \dd z ) , &
X_0
	&=	0 ,
\end{align}
and assume that $\sig = (\sig_t)_{t \geq 0}$ and $\gam(z) = (\gam_t(z))_{t \geq 0}$ are $\Fb$-adapted.  Then
\begin{align}
\Eb X_T^2
	&=	\Eb \int_0^T \Big( \sig_t^2 + \int_\Rb \gam_t^2(z) \nu(\dd z) \Big) \dd t ,
\end{align}
assuming the right-hand side is finite.
\end{theorem}

\begin{proof}
One could prove Theorem \ref{thm:isometry} directly from the definition of the It\^o integral.  However, we will use the It\^o formula instead.  Setting $f(x) = x^2$ and noting that
\begin{align}
f'(X_t)
	&=	2 X_t , &
f''(X_t) 
	&=	2 , &
f(X_{t-} + \gam_t(z)) - f(X_{t-})
	&=	2 X_{t-} \gam_t(z) + \gam_t^2(z) ,
\end{align}
we obtain from \eqref{eq:df} that
\begin{align}
\dd f(X_t)
	&=	\tfrac{1}{2} \sig_t^2 \cdot 2 \dd t 
			+ \sig_t 2 X_t \dd W_t \\ &\quad
			+	\int_\Rb \Big( 2 X_{t-} \gam_t(z) + \gam_t^2(z) \Big) \Nt( \dd t, \dd z) \\ &\quad
			+	\int_\Rb \Big( 2 X_{t-} \gam_t(z) + \gam_t^2(z) - \gam_{t}(z) \cdot 2 X_{t-} \Big) \nu(\dd z) \dd t \\
	&=	\Big( \sig_t^2 + 	\int_\Rb \gam_t^2(z) \nu(\dd z) \Big) \dd t  + (\cdots) \dd W_t + \int_\Rb (\cdots ) \Nt( \dd t, \dd z) ,
\end{align}
Integrating from $0$ to $T$ and taking an expectation, we obtain
\begin{align}
\Eb X_T^2
	&=	\Eb \int_0^T \Big( \sig_t^2 + 	\int_\Rb \gam_t^2(z) \nu(\dd z) \Big) \dd t 
			+ \Eb \int_0^T (\cdots) \dd W_t + \Eb \int_0^T \int_\Rb (\cdots ) \Nt( \dd t, \dd z) \\
	&=	\Eb \int_0^T \Big( \sig_t^2 + 	\int_\Rb \gam_t^2(z) \nu(\dd z) \Big) \dd t ,
\end{align}
where we have used $\Eb \int_0^T (\cdots) \dd W_t=0$ and $\Eb \int_0^T \int_\Rb (\cdots ) \Nt( \dd t, \dd z)=0$.
\end{proof}

\noindent
We now provide the multi-dimensional version of the It\^o formula.

\begin{theorem}
Suppose Assumption \ref{ass:integrability} holds.  Let $X = (X_t)_{t \geq 0}$ be a $d$-dimensional L\'evy-It\^o process with components satisfying
\begin{align}
\dd X_t
	&=	\mu_t \dd t + \sig_t \dd W_t + \int_{\Rb^m} \gam_t(z) \Nt(\dd t, \dd z) ,
\end{align}
where $W$ is an $n$-dimensional Brownian motion,
$\Nt(\dd t , \dd z) = N(\dd t , \dd z) - \nu(\dd z) \dd t$ is $k$-dimensional compensated Poisson random measure on $\Rb^m$ 
and
\begin{align}
\mu
	&: \Rb_+ \times \Om \to \Rb^d , &
\sig
	&: \Rb_+ \times \Om \to \Rb^{d \times n} , &
\gam
	&: \Rb_+ \times \Om \times \Rb^m \to \Rb^{d \times k} .
\end{align}
Component-wise, we can express the dynamics of $X$ as 
\begin{align}
\dd X_t^{(i)}
	&=	\mu_t^{(i)} \dd t + \sum_{j=1}^n \sig_t^{(i,j)} \dd W_t^{(j)} + \sum_{j=1}^k \int_{\Rb^m} \gam_t^{(i,j)}(z) \Nt^{(j)}( \dd t, \dd z) , &
i
	&=	1, 2, \ldots, d ,
\end{align}
Let $f : \Rb^d \to \Rb$ satisfy $f \in C^2(\Rb^d)$.  Then
\begin{align}
\dd f(X_t)
	&=	\Big( \sum_{i=1}^d \mu_t^{(i)} \d_{x_i} f(X_t) 
			+ \tfrac{1}{2} \sum_{i=1}^d \sum_{j=1}^d (\sig_t\sig_t^\text{T})_{i,j} \d_{x_i} \d_{x_j} f(X_t) \Big) \dd t
			+ \sum_{i=1}^d \sum_{j=1}^n \sig_t^{(i,j)} \d_{x_i} f(X_t)  \dd W_t^{(j)} \\ &\qquad
			+  \sum_{j=1}^k \int_{\Rb^m} \Big( f(X_{t-}+\gam_t^{(\,\cdot\, ,j)}(z)) - f(X_{t-}) \Big) \Nt^{(j)}(\dd t , \dd z) + \\ &\qquad
			+  \sum_{j=1}^k \int_{\Rb^m} \Big( f(X_{t-}+\gam_t^{(\,\cdot\, ,j)}(z)) - f(X_{t-}) - \gam_t^{(\,\cdot\, ,j)}(z) \cdot \nabla f(X_{t-}) \Big) \nu^{(j)}(\dd z) \dd t . 
\end{align}
\end{theorem}

\noindent
I would not recommend memorizing the above formula.  Simply book-mark this page.  
If you work with L\'evy-It\^o processes for a sufficient amount of time (about 5 years for the author of these notes), you will eventually have the above formula fixed in your mind.

\begin{theorem}[Quadratic covariation]
\label{thm:covar}
Suppose Assumption \ref{ass:integrability} holds and let $X = (X_t)_{t \geq 0}$ and $Y = (Y_t)_{t \geq 0}$ be given by
\begin{align}
\dd X_t
	&=	\mu_t \dd t + \sum_{j=1}^n \sig_t^{(j)} \dd W_t^{(j)} + \sum_{j=1}^k \int_{\Rb^m} \gam_t^{(j)}(z) \Nt^{(j)}( \dd t , \dd z ) , \label{eq:dX.QV} \\
\dd Y_t
	&=	b_t \dd t + \sum_{j=1}^n a_t^{(j)} \dd W_t^{(j)} + \sum_{j=1}^k \int_{\Rb^m} g_t^{(j)}(z) \Nt^{(j)}( \dd t , \dd z ) .
\end{align}
Here, $X$ and $Y$ are both one-dimensional processes, $W$ is a $n$-dimensional Brownian motion, $\Nt$ is an $k$-dimensional compensated Poisson random measure on $\Rb^m$ and
\begin{align}
\mu_t, b_t
	& \in \Rb , &
\sig_t, a_t
	&\in \Rb^n , &
\gam_t(\cdot), g_t(\cdot)
	&: \Rb^m \to \Rb^k .
\end{align}
Recall the Definition \ref{def:CoV} of quadratic covariation.
The {quadratic covariation} of $X$ and $Y$ up to time $T$, denoted $[X,Y]_T$, is given by
\begin{align}
[X,Y]_T
	&=	\int_0^T \sum_{j=1}^n \sig_t^{(j)} a_t^{(j)} \dd t + \int_0^T \sum_{j=1}^k \int_{\Rb^m} \gam_t^{(j)}(z) \cdot g_t^{(j)}(z) N^{(j)}(\dd t, \dd z) \\
	&=	\int_0^T \Big( \sum_{j=1}^n \sig_t^{(j)} a_t^{(j)} + \sum_{j=1}^k \int_{\Rb^m} \gam_t^{(j)}(z) \cdot g_t^{(j)}(z) \nu^{(j)}(\dd z) \Big) \dd t \\ &\quad
			+ \int_0^T \sum_{j=1}^k \int_{\Rb^m} \gam_t^{(j)}(z) \cdot g_t^{(j)}(z) \Nt^{(j)}(\dd t, \dd z) .
\end{align}
\end{theorem}

\noindent
We will not prove Theorem \ref{thm:covar}. Rather, we simply note that the above result relies on the following facts
\begin{align}
[W^{(i)}, N^{(j)}( \cdot, A ) ]_t
	&=	0 , &
[N^{(i)}(\cdot, A), N^{(j)}(\cdot, B)]_t
	&=	\del_{i,j} N^{(i)}(t, A \cap B) , &
[ N^{(i)}(\cdot, B),\Id]_t
	&=	0 ,
\end{align}
in addition to the previously derived facts 
\begin{align}
[W^{(i)},W^{(j)}]_t 
	&= \del_{i,j} t, &
[W^{(i)}, \Id]_t 
	&= 0 , &
[\Id,\Id]]_t 
	&= 0 .
\end{align}
Heuristically, one can derive Theorem \ref{thm:covar} by writing $\dd[X,Y]_t = \dd X_t \dd Y_t$ and using the rules
\begin{align}
\dd W_t^{(i)} \dd W_t^{(j)}
	&=	\del_{i,j} \dd t , &
\Nt^{(i)}(\dd t, \dd z) \Nt^{(j)}(\dd t, \dd y)
	&=	\del_{i,j} \del(z-y) N^{(i)}(\dd t, \dd z) \dd y, \\
\dd W_t^{(i)} \dd t
	&=	0 , &
\dd W_t^{(i)} \Nt^{(j)}(\dd t, \dd z)
	&=	0  , \\
\dd t \, \dd t
	&=	0 , &
\Nt^{(i)}(\dd t, \dd z) \, \dd t
	&=	0.
\end{align}

\begin{example}
Suppose $X$ is given by \eqref{eq:dX.QV}.  Then, from Theorem \ref{thm:covar}, the quadratic variation of $X$, denoted $[X,X]_T$, is given by
\begin{align}
[X,X]_T
	&=	\int_0^T \sum_{j=1}^n (\sig_t^{(j)} )^2 \dd t + \int_0^T \sum_{j=1}^k \int_{\Rb^m} (\gam_t^{(j)}(z) )^2 N^{(j)}(\dd t, \dd z) \\
	&=	\int_0^T \Big( \sum_{j=1}^n (\sig_t^{(j)} )^2 + \sum_{j=1}^k \int_{\Rb^m} (\gam_t^{(j)}(z) )^2 \nu^{(j)}(\dd z) \Big) \dd t \\ &\quad
			+ \int_0^T \sum_{j=1}^k \int_{\Rb^m} (\gam_t^{(j)}(z) )^2 \Nt^{(j)}(\dd t, \dd z) .
\end{align}
\end{example}

\begin{example}
Some texts \emph{define} quadratic covariation of two semimartingales $X$ and $Y$ as the unique process $[X,Y]_T$ that satisfies
\begin{align}
X_T Y_T
	&=	X_0 Y_0 + \int_0^T X_{t-} \dd Y_t + \int_0^T Y_{t-} \dd X_t + [X,Y]_T .
\end{align}
Let us check that this definition yields the same result as Theorem \ref{thm:covar}.  For simplicity, let us take
\begin{align}
\dd X_t
	&=	\mu_t \dd t + \sig_t \dd W_t + \int_{\Rb^m} \gam_t(z) \Nt( \dd t , \dd z ) , \\
\dd Y_t
	&=	b_t \dd t + a_t \dd W_t + \int_{\Rb^m} g_t(z) \Nt(\dd t ,\dd z ) ,
\end{align}
where $W$ is a one-dimensional Brownian motion and $\Nt$ is a one-dimensional compensated Poisson random measure on $\Rb^m$.
Let $f(X_t,Y_t) = X_t Y_t$ and observe that
\begin{align}
\d_x f(X_t,Y_t)
	&=	Y_t , &
\d_y f(X_t,Y_t)
	&=	X_t , &
\d_x \d_y f(X_t,Y_t)
	&=	1 ,  &
\d_x^2 f(X_t,Y_t)
	&=	0 , &
\d_y^2 f(X_t,Y_t)
	&=	0 , 
\end{align}
and also that
\begin{align}
f(X_{t-} + \gam_t(z), Y_{t-} + g_t(z) ) - f(X_{t-}, Y_{t-})
	&=	X_{t-} g_t(z)  + Y_{t-} \gam_t(z) + \gam_t(z) g_t(z) .
\end{align}
Plugging all of this into the multidimensional It\^o formula yields
\begin{align}
\dd ( X_t Y_t )
	&=	Y_t \( \mu_t \dd t + \sig_t \dd W_t \) + X_t \( b_t \dd t + a_t \dd W_t \) \\ &\quad
			+ \int_{\Rb^m} \Big( X_{t-} g_t(z)  + Y_{t-} \gam_t(z) + \gam_t(z) g_t(z) \Big) \Nt( \dd t, \dd z) \\ &\quad
			+	\int_{\Rb^m} \Big( X_{t-} g_t(z)  + Y_{t-} \gam_t(z) + \gam_t(z) g_t(z) - X_{t-} g_t(z) - Y_{t-} \gam_t(z) \Big) \nu( \dd z) \dd t \\
	&=	Y_{t-} \Big( \mu_t \dd t + \sig_t \dd W_t + \int_{\Rb^m} \gam_t(z) \Nt(\dd t, \dd z) \Big) 
			+ X_{t-} \Big( b_t \dd t + a_t \dd W_t + \int_{\Rb^m} g_t(z) \Nt(\dd t, \dd z) \Big) \\ &\quad
			+ \int_{\Rb^m} \gam_t(z) g_t(z) \Nt( \dd t, \dd z)
			+ \int_{\Rb^m} \gam_t(z) g_t(z) \nu( \dd z) \dd t \\
	&=	Y_{t-} \dd X_t + X_{t-} \dd Y_t + \dd[ X , Y ]_t .
\end{align}
Thus, we have verified the new definition of quadratic covariation agrees (at least, for the processes in this example) with Theorem \ref{thm:covar}.
\end{example}


\section{L\'evy-It\^o SDE}
We now introduce a \emph{L\'evy-It\^o Stochastic Differential equation (SDE)}, which is an equation of the form
\begin{align}
\dd X_t
	&=	\mu(t,X_t) \dd t + \sig(t,X_t) \dd W_t + \int_{\Rb^m} \gam(t,X_{t-},z) \Nt( \dd t , \dd z) , &
X_0
	&=	x \in \Rb^d , \label{eq:levy.SDE}
\end{align}
where $W$ is a $n$-dimensional Brownian motion, $\Nt(\dd t , \dd z)$ is a $k$-dimensional compensated Poisson random measure on $\Rb^m$, and
\begin{align}
\mu
	&: \Rb_+ \times \Rb^d \to \Rb^d , &
\sig
	&: \Rb_+ \times \Rb^d \to \Rb^{d \times n} , &
\gam
	&:	\Rb_+ \times \Rb^d \times \Rb^m \to \Rb^{d \times k} .
\end{align}
For $i = 1, 2, \ldots, d$, the dynamics of the $i$th component of $X$ are then given by
\begin{align}
\dd X_t^{(i)}
	&=	\mu^{(i)}(t,X_t) \dd t + \sum_{j=1}^n \sig^{(i,j)}(t,X_t) \dd W_t^{(j)} + \sum_{j=1}^k \int_{\Rb^m} \gam^{(i,j)}(t,X_{t-},z) \Nt^{(j)}( \dd t , \dd z) , 
\end{align}
with $X_0^{(i)} =	x_i \in \Rb^d$.

\begin{definition}
A \emph{solution} (specifically, a \emph{strong solution}) of \eqref{eq:levy.SDE} is a functional $F$ of the Brownian motion $W$ and compensated Poisson random measure $\Nt$
\begin{align}
X_t
	&=	F( W_s, \Nt(s,\dd z); 0 \leq s \leq t ) ,
\end{align}
such that
\begin{align}
X_T
	&=	x + \int_0^T \mu(t,X_t) \dd t + \int_0^T \sig(t,X_t) \dd W_t + \int_0^T \int_{\Rb^m} \gam(t,X_{t-},z) \Nt( \dd t , \dd z) ,
	&\forall \, T \geq 0 .
\end{align}
\end{definition}

\begin{theorem}
Consider the one-dimensional process $X$ driven by a single Brownian motion and a single Poisson random measure on $\Rb$.  If there exist constants $C_1$ and $C_2$ such that the \emph{linear growth} condition
\begin{align}
\mu^2(t,x) + \sig^2(t,x) + \int_\Rb \gam^2(t,x,z) \nu(\dd z)
	&\leq C_1 (1 + x^2) , &
	&\forall \, t \geq 0 , 
\end{align}
and the \emph{Lipschitz continuity} condition
\begin{align}
&| \mu(t,x) - \mu(t,y) |^2 + | \sig(t,x) - \sig(t,y) |^2 \\
&\quad + \int_\Rb |\gam(t,x,z) - \gam(t,y,z) |^2 \nu(\dd z)
	\leq C_2 |x-y|^2 , &
	&\forall \, t \geq 0, \, x,y \in \Rb ,
\end{align}
then, for any initial condition $X_0 = x \in \Rb$, SDE \eqref{eq:levy.SDE} has a unique strong solution, which is adapted to the filtration generated by $W$ and $N$.  Moreover, the solution $X$ is a Markov process.  The conditions for higher dimensions are analogous.
\end{theorem}

\noindent
Since $X$, the solution \eqref{eq:levy.SDE} is a Markov we can define a corresponding semigroup $\Pc(t,s)$ and infinitesimal generator $\Ac(t)$.  For a sufficiently nice function $\varphi: \Rb^d \to \Rb$, we have
\begin{align}
\Pc(t,s) \varphi(x)
	&:=		\Eb[ \varphi(X_s) | X_t = x ] , &
0
	&\leq t \leq s < \infty , \\
\Ac(t) \varphi(x)
	&:= \lim_{s \searrow t} \frac{1}{s-t} \Big( \Pc(t,s)\varphi(x) - \varphi(x) \Big) .
\end{align}

\begin{theorem}
Suppose Assumption \ref{ass:integrability} holds and $X$ is given by SDE \eqref{eq:levy.SDE}.  Then
\begin{align}
\Ac(t) 
	&=	\sum_{i=1}^d \mu_i(t,x) \d_{x_i} + \tfrac{1}{2} \sum_{i=1}^d \sum_{j=1}^d (\sig \sig^\text{T} )_{i,j}(t,x) \d_{x_i} \d_{x_j} \\ & \quad
			+ \sum_{j=1}^k \int_{\Rb^m} \nu^{(j)}( \dd z )  \Big( \theta_{\gam^{(\,\cdot\,,j)}(t,x,z)} - 1 - \gam^{(\,\cdot\, ,j)}(t,x,z) \cdot \nabla \Big) ,
\end{align}
where $\theta_\gam$ is the \emph{shift operator}: $\theta_\gam f(x) := f(x + \gam)$ and $C_0^2(\Rb^d) \subseteq \text{dom}(\Ac(t))$.
\end{theorem}

\begin{proof}
The proof is the same as in the no-jump case.  First, write $\varphi(X_s) = \varphi(x) + \int_t^s \dd \varphi(X_u)$.  Next, use the It\^o formula to write $\dd \varphi(X_u)$ explicitly as terms involving $\dd u$, $\dd W_u$ and $\Nt(\dd u,\dd z)$.  Finally, take an expectation, note that integrals with respect to $W$ and $\Nt$ are martingales, and send $s \searrow t$.  We omit the details.
\end{proof}


\noindent
We can now write the It\^o formula for the solution $X$ of L\'evy-It\^o SDE \eqref{eq:levy.SDE} in a more compact form
\begin{align}
\dd \varphi(X_t)
	&=	\Ac(t) \varphi(X_t) \dd t + \nabla \varphi \cdot \sig(t,X_t) \dd W_t 
			+ \int_{\Rb^d} \Big( \varphi( X_{t-} + \gam(t,X_{t-},z) ) - \varphi( X_{t-}) \Big) \Nt( \dd t, \dd z ) .
\end{align}

\begin{theorem}
Let $X = (X_t)_{t \geq 0}$ satisfy
\begin{align}
\dd X_t
	&=	\mu(t,X_t) \dd t + \sig(t,X_t) \dd W_t + \int_{\Rb^d} \gam(t,X_{t-},z) \Nt( \dd t, \dd z) .
\end{align}
Define $u(t,x) := \Eb[ \varphi(X_T) | X_t = x ]$.  If $u \in C^{1,2}$ and Assumption \ref{ass:integrability} holds, then $u$ satisfies the \emph{Kolmogorov Backward Equation}
\begin{align}
(\d_t + \Ac(t)) u(t,\cdot)
	&=	0 , &
u(T,\cdot)
	&=	\varphi .
\end{align}
\end{theorem}

\begin{proof}
The proof is exactly analogous to the diffusion case.  First, show that $\Eb[ \varphi(X_T) | \Fc_t] = u(t,X_t)$ is a martingale.  Take the differential of $u(t,X_t)$ and set the $\dd t$-term equal to zero to obtain the partial integro-differential equation that must be satisfied by $u$.  The terminal condition is obtained from $\Eb[\varphi(X_T)|\Fc_T] = \varphi(X_T) = u(T,X_T)$.
\end{proof}

\begin{example}
Suppose $X = (X_t)_{t \geq 0}$ is the following one-dimensional L\'evy process
\begin{align}
\dd X_t
	&=	\mu \dd t + \sig \dd W_t + \int_\Rb z \Nt( \dd t , \dd z ) ,
\end{align}
where $W$ is a one-dimensional Brownian motion and $\Nt$ is a one-dimensional compensated Poisson random measure on $\Rb$.
We wish to find an expression for $u(t,x) := \Eb[ \varphi(X_T) | X_t = x ]$, which satisfies the KBE; in this case
\begin{align}
(\d_t + \Ac ) u(t, \cdot)
	&=	0 , &
u(T,\cdot)
	&=	\varphi , &
\Ac
	&=	\mu \d_x + \tfrac{1}{2} \sig^2 \d_x^2 + \int_\Rb \nu(dz) ( \theta_z - 1 - z \d_x) .
\end{align}
It will be useful to define the \emph{Fourier} and \emph{inverse Fourier transforms}
\begin{align}
\text{Fourier Transform}:&&
\Fv [\varphi](\xi) = \varphih(\xi)
	&:=	\int_\Rb \ee^{- \ii \xi x} \varphi(x) \dd x , \\
\text{Inverse Transform}:&&
\Fv^{-1}[ \varphih](x)
	&:=	\frac{1}{2 \pi} \int_\Rb \ee^{\ii \xi x} \varphih(\xi) \dd \xi .
\end{align}
Now, note that
\begin{align}
\Fv[ \d_t u(t,\cdot) ](\xi)
	&=	\int_\Rb \ee^{- \ii \xi x} \d_t u(t,x) \dd x \\
	&=		\d_t \int_\Rb \ee^{- \ii \xi x} u(t,x) \dd x 
	=		\d_t \uh(t,\xi) \\
\Fv[ \Ac u(t,\cdot) ](\xi)
	&=	\int  \ee^{- \ii \xi x} \Ac u(t,x) \dd x 
	=		\int u(t,x) {\Ac}^* \ee^{- \ii \xi x} \dd x \\
	&=	\int_\Rb u(t,x) \psi(\xi) \ee^{- \ii \xi x} \dd x 
	=		\psi(\xi) \uh(t,\xi) .
\end{align}
where $\Ac^*$ is the $L^2(\Rb,\dd x)$ adjoint of $\Ac$ and $\psi(\xi)$ is the characteristic exponent of $X$.  Specifically
\begin{align}
\Ac^*
	&=	- \mu \d_x + \tfrac{1}{2} \sig^2 \d_x^2 + \int_\Rb \nu(\dd z) ( \theta_{-z} - 1 + z \d_x) , \\
\psi(\xi)
	&=	\mu \ii \xi - \tfrac{1}{2} \sig^2 \xi^2 + \int_\Rb \nu(\dd z) ( \ee^{\ii \xi z } - 1 - \ii \xi z ) ,
\end{align}
(check that you can derive $\Ac^*$ and $\psi(\xi)$ for yourself!).  Taking the Fourier transform of the KBE and using the above results, we obtain an ODE in $t$ for $\uh(t,\xi)$
\begin{align}
(\d_t + \psi(\xi)) \uh(t,\xi)
	&=	0 , &
\uh(T,\xi)
	&=	\varphih(\xi) .
\end{align} 
The solution of this ODE is
\begin{align}
\uh(t,\xi)
	&=	\ee^{(T-t)\psi(\xi)} \varphih(\xi) .
\end{align}
Finally, we obtain $u$ as the inverse Fourier transform of $\uh$
\begin{align}
u(t,x)
	&=	\Fv^{-1}[ \uh(t,\cdot) ](x)
	=		\frac{1}{2 \pi} \int_\Rb \ee^{\ii \xi x + (T-t)\psi(\xi)} \varphih(\xi) .
\end{align}
where, we have assumed, for simplicity, that $\varphi = \Fv^{-1}[\varphih]$.  This is usually satisfied if $\varphi$ is continuous and in $L^1(\Rb,\dd x)$.
\end{example}

\section{Change of measure}
\label{sec:change.of.measure}

\begin{assumption}
Throughout Section \ref{sec:change.of.measure}, $W$ is a one-dimensional Brownian motion and $\Nt$ is a one-dimensional compensated Poisson random measure on $\Rb$.  All processes are scalar, unless specifically stated otherwise.
\end{assumption}

\subsection*{Review of the Girsanov change of measure for an It\^o process}
Fix a probability space $(\Om,\Fc,\Pb)$ and let $\Fb=(\Fc_t)_{t \geq 0}$ be a filtration on this space.  In Section \ref{sec:measure.change} we showed that if $Z$ is a random variable that satisfies
\begin{align}
Z
	&\geq 0 , &
\Eb Z
	&=	1 ,
\end{align}
Then we can define a new probability measure $\Pbt$ using
\begin{align}
\Pbt(A)
	&=	\Eb Z \Ib_A , &
	&\forall \, A \in \Fc ,
\end{align}
and we call $Z = \frac{\dd \Pbt}{\dd \Pb}$ the \emph{Radon-Nikod\'ym derivative of $\Pbt$ with respect to $\Pb$}.  Moreover, if $Z>0$, then probabilities under $\Pb$ can be obtained from probabilities under $\Pbt$ using
\begin{align}
\Pb(A)
	&=	\Ebt \frac{1}{Z} \Ib_A , &
	&\forall \, A \in \Fc ,
\end{align}
and thus, we identify $\frac{1}{Z} = \frac{\dd \Pb}{\dd \Pbt}$ as the Radon-Nikod\'ym derivative of $\Pbt$ with respect to $\Pb$.
\\[0.5em]
In Chapter \ref{chapter:stoch.calc}, we defined a \emph{Radon-Nikod\'ym derivative process} (Definition \ref{def:rndp}) by fixing a time horizon $T$ and setting
\begin{align}
Z_t
	&:=	\Eb[ Z_T | \Fc_t ] , &
0
	&\leq t \leq T ,
\end{align}
We showed (Lemma \ref{lem:2}) that 
\begin{align}
\Ebt[ Y | \Fc_s ]
	&=	\frac{1}{Z_s} \Eb[ Z_t Y | \Fc_s ] , &
Y
	&\in \Fc_t , &
0
	&\leq s \leq t \leq T .
\end{align}
In particular if a process $Y = (Y_t)_{0 \leq t \leq T}$ is $\Fb$-adapted and $Y Z := (Y_t Z_t)_{0 \leq t \leq T}$ is a martingale under $\Pb$ then the process $Y$ is a martingale under $\Pb$ since
\begin{align}
\Ebt[ Y_t | \Fc_s ]
	&=	\frac{1}{Z_s} \Eb[ Z_t Y_t | \Fc_s ] 
	=	\frac{1}{Z_s} Z_s Y_s 
	=	Y_s . \label{eq:Etilde.martingale}
\end{align}
As a specific example of this machinery, we showed (Girsanov's Theorem \ref{thm:girsanov}) that if we defined
\begin{align}
\frac{\dd \Pbt}{\dd \Pb} 
	&= Z_T , &
Z_t
	&=	\Eb[Z_T| \Fc_t ] = \exp \Big( - \int_0^t \tfrac{1}{2} \Theta_s^2 \dd s - \int_0^t \Theta_s \dd W_s \Big) ,
\end{align}
then, under $\Pbt$, the process $\Wt$ defined by
\begin{align}
\Wt_t
	&=	W_t + \int_0^t \Theta_s \dd s , &
0
	&\leq t \leq T ,
\end{align}
is a martingale (in fact, a Brownian motion).  We are now going to present an analogous result for Poisson random measures.

\subsection*{Girsanov change of measure for Poisson random measures}

\begin{theorem}[Girsanov Theorem for Poisson random measures]
\label{thm:girsanov.jump}
Fix a probability space $(\Om,\Fc,\Pb)$ and a filtration $\Fb=(\Fc_t)_{0 \leq t \leq T}$.  Let $\eta(z) = (\eta_t(z))_{0 \leq t \leq T}$ be $\Fb$-adapted and define a process $Z = (Z_t)_{0 \leq t \leq T}$ by
\begin{align}
Z_t
	&:=	\exp \Big( - \int_0^t \int_\Rb \( \ee^{\eta_s(z)} - 1 - \eta_s(z) \) \nu(\dd z) \dd s 
			+ \int_0^t \int_\Rb \eta_s(z) \Nt(\dd s , \dd z) \Big) , \label{eq:Z.jump} .
\end{align}
Assume that $\Eb Z_T = 1$ and set $\frac{\dd \Pbh}{\dd \Pb} = Z_T$.  Define
\begin{align}
\Nh( \dd t, \dd z)
	&:=		\Nt( \dd t , \dd z) - \( \ee^{\eta_t(z)} - 1\)\nu(\dd z) \dd t
	=			N( \dd t , \dd z ) - \ee^{\eta_t(z)} \nu(\dd z) \dd t . \label{eq:N.hat}
\end{align}
Then $\Nh$ is a compensated Poisson random measure under $\Pbh$ in the sense that the process $X = (X_t)_{0 \leq t \leq T}$, defined by
\begin{align}
X_t
	:=	\int_0^t \int_\Rb \gam_s(z) \Nh( \dd s, \dd z) , 
\end{align}
is a martingale for all $\Fb$-adapted processes $(\gam_t(z))_{0 \leq t \leq T}$ satisfying
\begin{align}
\Ebh \int_0^T \int_\Rb \gam_t^2(z) \ee^{\eta_t(z)} \nu(\dd z ) \dd t
	&< \infty .
\end{align}
\end{theorem}


\begin{proof}
First, we observe that $(Z_t)_{0 \leq t \leq T}$ is a martingale under $\Pb$.  To see this, simply note that
\begin{align}
\dd Z_t
	&=	\int_\Rb Z_{t-} \( \ee^{\eta_t(z)} - 1 \) \Nt( \dd t , \dd z ) .
\end{align}
Thus, we have
\begin{align}
1
	&=	Z_0
	=		\Eb Z_T
	=		\Eb Z , &
Z_t
	&=	\Eb[ Z_T | \Fc_t ]
	=		\Eb[ Z | \Fc_t ] .
\end{align}
As $\Eb Z_T = 1$ and, by construction, $Z > 0$, the random variable $Z_T$ defines a Radon-Nikod\'ym derivative $\frac{\dd \Pbh}{\dd \Pb}$.
Also, as $Z_t = \Eb[ Z | \Fc ]$ we see that $Z = (Z_t)_{0 \leq t \leq T}$ is a Radon-Nikod\'ym derivative process.
Now, note that $X=(X_t)_{0 \leq t \leq T}$ is adapted to the filtration $\Fb$.
In light of \eqref{eq:Etilde.martingale}, to show that $X$ is a martingale under $\Pbh$ we need only to show that $XZ=(X_t Z_t)_{0 \leq t \leq T}$ is a martingale under $\Pb$.  To this end we note
\begin{align}
\dd X_t
	&=	\int_\Rb \gam_t(z) \Nh( \dd t, \dd z) 
	 =	- \int_\Rb \gam_t(z) \( \ee^{\eta_t(z)} - 1\)\nu(\dd z) \dd t + \int_\Rb \gam_t(z) \Nt( \dd t, \dd z) .
\end{align}
Next, we compute
\begin{align}
\dd ( X_t Z_t)
	&=	X_{t-} \dd Z_t + Z_{t-} \dd X_t + \dd [ X, Z]_t \\
	&=	X_{t-} \dd Z_t + Z_{t-} \Big( - \int_\Rb \gam_t(z) \( \ee^{\eta_t(z)} - 1\)\nu(\dd z) \dd t  + \int_\Rb \gam_t(z) \Nt(\dd t, \dd z) \Big) \\ & \quad
			+ Z_{t-} \int_\Rb \gam_t(z) \( \ee^{\eta_t(z)} - 1 \) N( \dd t , \dd z ) \\
	&=	X_{t-}\dd Z_t + Z_{t-} \int_\Rb \gam_t(z) \Nt(\dd t, \dd z) \\ & \quad
			+ Z_{t-} \int_\Rb \gam_t(z) \( \ee^{\eta_t(z)} - 1 \) \Nt( \dd t , \dd z ) .
\end{align}
As integrals with respect to $\dd Z_t$  and $\Nt(\dd t, \dd z)$ are martingales under $\Pb$, we see that $XZ$ is a martingale under $\Pb$.  And thus, $(X_t)_{0 \leq t \leq T}$ is a martingale under $\Pbh$.
\end{proof}

\begin{remark}
\label{rmk:intensity}
The takeaway message from Theorem \ref{thm:girsanov.jump} is that if $\Eb N(\dd t,\dd z) = \nu(\dd z) \dd t$ then $\Ebh[ N(\dd t, \dd z)| \Fc_t ] = \ee^{\eta_t(z)} \nu(\dd z) \dd t$. In other words, the theorem tells us how we can change the intensity of jumps through an appropriate change of measure
\begin{align}
	&\text{$\nu(\dd z)$ under $\Pb$}&
	&\to &
	&\text{$\nuh_t(\dd z) := \ee^{\eta_t(z)} \nu(\dd z)$ under $\Pbh$}.
\end{align}
\end{remark}

\noindent
Let us see an example.

\begin{example}
Consider a pure-jump L\'evy process $X=(X_t)_{0 \leq t \leq T}$
\begin{align}
\dd X_t
	&=	\mu \dd t + \int_\Rb z \Nt( \dd t, \dd z ) , &
\Nt( \dd t, \dd z)
	&=	N( \dd t, \dd z) - \nu( \dd z ) \dd t .
\end{align}
We wish to find a change of measure $Z_T = \frac{\dd \Pbt}{\dd \Pb}$ such that, under $\Pbh$ jumps of size $\dd z$ arrive with intensity $\nuh(\dd z)$.
From Theorem \ref{thm:girsanov.jump} we see that the process $Z$ must be of the form
\begin{align}
Z_t
	&:=	\exp \Big( - \int_0^t \int_\Rb \( \ee^{\eta_s(z)} - 1 - \eta_s(z) \) \nu(\dd z) \dd s + \int_0^t \int_\Rb \eta_s(z) \Nt(\dd s , \dd z) \Big) .
	\label{eq:Z.girsanov.ex}
\end{align}
With $Z$ given by \eqref{eq:Z.girsanov.ex} we have that
\begin{align}
\Nh( \dd t, \dd z)
	&=	\Nt( \dd t , \dd z) -\( \ee^{\eta_t(z)} - 1\)\nu(\dd z) \dd t
	=	N( \dd t , \dd z ) - \ee^{\eta_t(z)} \nu(\dd z) \dd t . \label{eq:N.hat.ex}
\end{align}
is a compensated Poisson random measure under $\Pbh$.  We want
\begin{align}
\ee^{\eta_t(z)} \nu(\dd z)
	&=	\nuh(\dd z) .
\end{align}
Thus, we identify
\begin{align}
\eta_t(z)
	&=	\log \frac{\nuh(\dd z)}{\nu(\dd z)} .
\end{align}
Note that, in order for $\eta$ to be well-defined we must have $0 < \frac{\nuh(\dd z)}{\nu(\dd z)} < \infty$, which holds if an only if $\nuh \sim \nu$ (i.e., the measures are equivalent).
Re-writing the dynamics of $X$ as follows
\begin{align}
\dd X_t
	&=	\mu \dd t + \int_\Rb z \( \Nh( \dd t, \dd z ) + \( \ee^{\eta_t(z)} - 1\) \nu(\dd z) \dd t \) \\
	&=	\mu \dd t + \int_\Rb z \( \Nh( \dd t, \dd z ) + \( \frac{\nuh(\dd z)}{\nu(\dd z)} - 1\) \nu(\dd z) \dd t \) \\
	&=	\mu \dd t + \int_\Rb z \( \Nh( \dd t, \dd z ) + \( \nuh(\dd z) - \nu(\dd z) \) \dd t \) \\
	&=	\muh \dd t + \int_\Rb z \Nh( \dd t, \dd z ) \\
\muh
	&=	\mu + \int_\Rb z \( \nuh(\dd z) - \nu(\dd z) \) ,
\end{align}
we see that $X$ will be a martingale under $\Pbh$ if and only if $\muh = 0$.  Suppose, for example, that $\nu = \lam F$, where $\lam = \nu(\Rb)$ is the intensity of all jumps (assumed finite) and $F(\dd z) = \nu(\dd z)/\nu(\Rb)$ is the jump distribution.  If  we wish to change only the intensity of jumps (say, from $\lam$ to $\lamh$) and not the distribution of jumps $F$, then, if we want $X$ to be a martingale under $\Pbh$ we must have
\begin{align}
\lamh - \lam
	&=	\mu \cdot \Big( \int_\Rb z F(\dd z) \Big)^{-1} .
\end{align}
Note that, had we included a Brownian component in the $\Pb$ dynamics of $X$, we would have complete freedom to change the jump measure since any non-zero drift could be absorbed into the drift of the Brownian motion under the change of measure.
\end{example}

\section{Hawkes processes}
\label{sec:hawkes}
Recall that if $P = (P_t)_{t \geq 0}$ is a Poisson process with intensity $\lam$, then we have
\begin{align}
\Eb[ \dd P_t | \Fc_t ]
	&=	\Eb \dd P_t  
	=		\lam \dd t .
\end{align}
Note that the probability of a jump $\Pb( \dd P_t = 1 ) = \Eb \dd P_t$ in the instant $\dd t$ is entirely \emph{independent} of this history $\Fc_t$ of $P$.  This memorylessness property is convenient from a computational point of view and is necessary in order for $P$ to be Markov.  But, from a modeling perspective, limiting our analysis to Poisson processes (or more generally, Poisson random measures) is very restrictive.  For example, if we want to model the arrival of earthquakes in Seattle, it is well-known that the occurrence of an earthquake today will increase the probability of an earthquake tomorrow.  Thus, it would be useful to have some way of modeling events whose arrival intensities are history-\emph{dependent}.  This is precisely what a Hawkes process will allow us to do.

\begin{definition}
A \emph{Hawkes process} $N=(N_t)_{t \geq 0}$ is a simple counting process satisfying
\begin{align}
\Eb[ \dd N_t | \Fc_{t-} ]
	&=	\lam_t \dd t , &
\lambda_t
	&:=	\Lam \Big( \int_0^{t-} h(t-s) \dd N_s \Big), \label{eq:lambda-0}
\end{align}
where $h(\cdot):\Rb_+ \rightarrow \Rb_+$ is the \emph{kernel} (sometimes also called an \emph{exciting function}) and $\Lam(\cdot):\Rb_+ \to \Rb_+$ is locally integrable.  The intensity $\lambda_t$ is an $\Fc_t$-measurable random variable where $\Fc_t$ is the natural filtration for $N$.  When the function $\Lam(\cdot)$ is linear, the process $N_t$ is known as the \emph{linear Hawkes process}, otherwise, it is called a \emph{nonlinear Hawkes process}.
\end{definition}

\noindent
If we denote by $\tau_1, \tau_2, \ldots $ the (random) times that $N$ jumps, then the intensity at time $t$ is given by
\begin{align}
\lam_t
	&=	\Lam \Big( \sum_{\tau_i < t} h(t - \tau_i) \Big) .
\end{align}
Typically, $h$ is a strictly decreasing function and $\Lam$ is increasing.  Under these assumptions, the intensity $\lam_t$ is decreasing over the intervals of the form $[\tau_i,\tau_{i+1})$.
\\[0.5em]
For reason of analytic tractability, in what follows, we concentrate on the linear Hawkes processes.
Specifically, let us assume 
\begin{align}
\Lam(z)
	&=	\nu + z ,
\end{align}
where $\nu > 0$ is a strictly positive constant.
Thus, from \eqref{eq:lambda-0}, we have
\begin{align}
\lambda_t
	&=	\nu + \int_0^{t-}h(t-s) \dd N_s \label{eq:lam.1} \\
	&=	\nu + \int_0^{t-}h(t-s) \dd \Nt_s + \int_0^{t-}h(t-s) \lam_s \dd s , \label{eq:lam.2}
\end{align}
where we have introduced the \emph{compensated Hawkes process}, whose differential is given by
\begin{align}
\dd \Nt_s
	&=	\dd N_s - \lam_s \dd s .
\end{align}
As, by construction, the expected infinitesimal change of $\Nt$ is zero
\begin{align}
\Eb[ \dd \Nt_t | \Fc_t] 
	&= \Eb[ \dd N_t | \Fc_t] - \lam_t \dd t = 0 ,
\end{align}
 it follows that the compensated Hawkes process is a martingale
\begin{align}
\Eb [ \Nt_t | \Fc_s ]
	&=	\int_s^t \Eb[ \Eb( \dd \Nt_u | \Fc_u ) | \Fc_s ] + \Nt_s
	 =	\Nt_s , &
0
	&\leq s \leq t < \infty .
\end{align}
Thus, integrals with respect to the compensated Hawkes process are martingales as well.
\\[0.5em]
Suppose the intensity $\lam_t$ of a Hawkes process approaches a stationary distribution as $t \to \infty$.  This distribution must have a mean, which we shall call $\bar{\lam}$.  If we start the process from time $-\infty$, which is equivalent to looking at the behavior of the process as $t \to \infty$, then we have $\Eb \lam_t = \bar{\lam}$ for all $t>-\infty$ and thus, taking an expectation of \eqref{eq:lam.2}, we find
\begin{align}
\Eb \lambda_t
	&=	\nu + \int_{-\infty}^{t-}h(t-s) \Eb \dd \Nt_s + \int_{-\infty}^{t-}h(t-s) \Eb \lam_s \dd s \\
	&=	\nu + \int_{-\infty}^{t-}h(t-s) \Eb \lam_s \dd s , \\
\bar{\lam}
	&=	\nu + \bar{\lam} \int_0^\infty h(s) \dd s .
\end{align}
Hence, we find that the mean $\bar{\lam}$ of the stationary distribution of the intensity is given by
\begin{align}
\bar{\lam}
	&=	\nu \cdot \Big( 1 - \int_0^\infty h(s) \dd s \Big)^{-1} .
\end{align}
Clearly, a stationary distribution only exists if $h$ satisfies
\begin{align}
\int_0^\infty h(s) \dd s
	< 1 . \label{eq:h.condition}
\end{align}
Thus, some texts \emph{require} that $h$ satisfy \eqref{eq:h.condition}
\\[0.5em]
In general, a given Hawkes process $N = (N_t)_{t \geq 0}$ and the associated intensity process $\lam = (\lam_t)_{t \geq 0}$ are non-Markovian.  However, when the kernel $h(\cdot)$ has an exponential form
\begin{align}
h(t)
	&=	a \ee^{-bt}, &
a,b
	&>	0, \label{eq:exp.form}
\end{align}
then the intensity process $\lam$ by itself and the pair $(\lam,N)$ are Markov.  To see this, observe from \eqref{eq:lam.1} and \eqref{eq:exp.form} that
\begin{align}
\lam_t
	&=	\nu + \int_0^{t-} a \ee^{-b(t-s)} \dd N_s .
\end{align}
Taking the differential of $\lam_t$ we find
\begin{align}
\dd \lam_t
	&=	- b \Big( \int_0^{t-} a \ee^{-b(t-s)} \dd N_s \Big) \dd t + a \dd N_t \\
	&=	- b ( \lam_t - \nu ) \dd t + a \dd N_t , \\
	&=	- b ( \lam_t - \nu ) \dd t + a \dd \Nt_t + a \lam_t \dd t \\
	&=	\Big( (a - b) \lam_t + b \nu \Big) \dd t + a \dd \Nt_t .
\end{align}
Since the intensity of $N_t$ is $\lam_t$, we see that the future dynamics of $\lam$ depend only on the value of $\lam_t$ and not the entire history $\Fc_t$.  Thus, the process $\lam$ is Markov and its generator is
\begin{align}
\Ac
	&=	\Big( (a -b) \lam + b \nu \Big) \d_\lam  + \lam \( \theta_a - 1 - a \d_\lam  \) \\
	&=	 b ( \nu - \lam ) \d_\lam + \lam \( \theta_a - 1 \) ,
\end{align}
where, we remind the reader that $\theta_a$ is the shift operator: $\theta_a f(\lam) = f(\lam+a)$.  Since $\lam$ alone is Markov, it follows that the pair $(\lam,N)$ is Markov and has generator
\begin{align}
\Ac
	&=	b ( \nu - \lam ) \d_\lam + \lam \( \theta_{(a,1)} - 1 \) , \label{eq:A.lam.n}
\end{align}
where $\theta_{(a,1)}f(\lam,n) = f(\lam+a,n+1)$.
\\[0.5em]
We cannot find the probability mass function of $N_T$.  However, we can find its Laplace transform $\Eb \ee^{ - \eta N_T}$.
Let us define
\begin{align}
u(t,\lam,n)
	&=	\Eb[ \ee^{- \eta N_T} | \lam_t = \lam , N_t = n ] .
\end{align}
The function $u$ satisfies the Kolmogorov backward equation
\begin{align}
0
	&=	( \d_t + \Ac ) u , &
u(T,\lam,n)
	&=	\ee^{- \eta n} . \label{eq:hawkes.pde}
\end{align}
where $\Ac$ is given by \eqref{eq:A.lam.n}.  We seek a solution to \eqref{eq:hawkes.pde} in exponential affine form
\begin{align}
u(t, \lam, n)
	&=	\ee^{- \eta n + A(t) + B(t) \lam } . \label{eq:ansatz}
\end{align}
Inserting \eqref{eq:ansatz} into \eqref{eq:hawkes.pde}, we obtain
\begin{align}
0
	&=	u(t,\lam,n) \Big( (A'(t) + \lam B'(t)) + b ( \nu - \lam ) B(t) + \lam (\ee^{-\eta + a B(t) }-1) \Big) , &
\ee^{- \eta n + A(T) + B(T) \lam }
	&=	\ee^{- \eta n  } .
\end{align}
This must hold for all $\lam$.  Thus we collect terms of like order in $\lam$ and obtain a pair of coupled ODEs for $A$ and $B$
\begin{align}
\Oc(\lam^0):&&
0
	&=	A'(t) + b \nu B(t) , &
A(T)
	&=	0 , \label{eq:ODE.A} \\
\Oc(\lam^1):&&
0
	&=	B'(t) - b B(t) + (\ee^{-\eta + a B(t) }-1) , &
B(T)
	&=	0 . \label{eq:ODE.B}
\end{align}
These ODEs must be solved numerically.
\\[0.5em]
We finish this Section by defining a $d$-dimensional Hawkes process.

\begin{definition}
A \emph{$d$-dimensional Hawkes process} $N=(N_t^{(1)}, N_t^{(2)}, \ldots, N_t^{(d)})_{t \geq 0}$ is a simple counting process satisfying
\begin{align}
\Eb[ \dd N_t^{(i)} | \Fc_{t-} ]
	&=	\lam_t^{(i)} \dd t , &
\lambda_t^{(i)}
	&:=	\Lam^{(i)} \Big( \sum_{j=1}^d \int_0^{t-} h^{(i,j)}(t-s) \dd N_s^{(j)} \Big), \label{eq:lambda}
\end{align}
where $h^{(i,j)}(\cdot):\Rb_+ \rightarrow \Rb_+$ and $\Lam^{(i)}(\cdot):\Rb_+ \to \Rb_+$ is locally integrable.  The intensity $\lambda_t^{(i)}$ is an $\Fc_t$-measurable random variable where $\Fc_t$ is the natural filtration for the $d$-dimensional process $(N^{(i)})_{i=1}^d$.  When the function $\Lam^{(i)}(\cdot)$ is linear for all $i$, the process $N_t$ is known as the \emph{linear $d$-dimensional Hawkes process}, otherwise, it is called a \emph{nonlinear $d$-dimensional Hawkes process}.
\end{definition}

\noindent
As in the one-dimensional case, linear $d$-dimensional Hawkes processes with exponential kernels are Markov and admit many analytically tractable results.

\section{Exercises}
\label{sec:exercises-10}

\begin{exercise}
\label{exer:10.1}
Let $P=(P_t)_{t \geq 0}$ be a Poisson process with intensity $\lam$.\\
(a) What is the L\'evy measure $\nu$ of $P$?\\
(b) Let $\dd X_t = \dd P_t$.  Define $u(t,x) := \Eb[ \varphi(X_T) | X_t = x ]$.
Find $u(t,x)$ and verify that it solves the Kolmogorov Backward equation.
\end{exercise}

\begin{exercise}
\label{exer:10.2}
Let $X = (X_t)_{t \geq 0}$ and $Y=(Y_t)_{t \geq 0}$ be L\'evy processes defined by
\begin{align}
X_t
	&=	b t + \sig W_t + \int_{\Rb} z \Nt(t, \dd z) , &
Y_t
	&=	\gam t + \int_{\Rb_+} z M(t, \dd z) ,
\end{align} 
where $b \in \Rb$, $\sig, \gam \in \Rb_+$, $\Nt$ is a compensated Poisson random measure with L\'evy measure $\nu$ and $M$ is a Poisson random measure with L\'evy measure $\mu$.  Assume that $N$ and $M$ are independent and $\nu(\Rb) < \infty$ and $\mu(\Rb_+) < \infty$.  As the process $Y$ is strictly increasing, we can constrict a process $Z = (Z_t)_{t \geq 0}$ as follows
\begin{align}
Z_t
	&:= X_{Y_t} .
\end{align}
(a) Show that $Z$ is a L\'evy process.  For this problem, you may assume that $Z$ is continuous in probability.\\
(b) Find $\psi_Z$ where $\Eb \ee^{\ii \xi Z_t} = \ee^{t \psi_Z(\xi)}$.  You may leave your answer as a composition of functions if you find it helpful.\\
(c) Suppose $\nu \equiv 0$ (i.e., the $X$ process experiences no jumps).  Find the drift $A$, volatility $\Sig$ and L\'evy measure $\Pi(\dd z)$ corresponding to $Z$.
Write the L\'evy-It\^o decomposition of $Z_T$ in terms of a Brownian motion $B$ and a Poisson random measure $P$ whose L\'evy measure is $\Pi$.
\end{exercise}

\begin{exercise}
\label{exer:10.3}
Let $X = (X_t)_{t \geq 0}$ and $Y=(Y_t)_{t \geq 0}$ be processes defined by
\begin{align}
\dd X_t
	&=	\mu_t X_t \dd t + \sig_t X_t \dd W_t + \int_\Rb (\ee^{\gam_t(z)} - 1 ) X_{t-} \Nt( \dd t, \dd z) , \\
\dd Y_t
	&=	b_t Y_t \dd t + a_t Y_t \dd W_t + \int_\Rb (\ee^{g_t(z)} - 1 ) Y_{t-} \Nt( \dd t, \dd z) ,
\end{align} 
where $W$ is a one-dimensional Brownian motion, $\Nt$ is a one-dimensional compensated Poisson random measure on $\Rb$, and $\mu, b, \sig, a, \gam, g$ and $\Fb$-adapted stochastic processes.\\
(a) Define $Z_t := X_t / Y_t$.  Compute the differential $\dd Z_t$.  Your answer should not involve $X_t$ or $Y_t$.\\
(b) Find $\mu_t$ so that $Z$ is a martingale.
\end{exercise}

\begin{exercise}
\label{exer:10.4}
Let $\eta = (\eta_t)_{t \geq 0}$ be a one-dimensional L\'evy process and define $X = (X_t)_{t \geq 0}$ by
\begin{align}
\dd X_t
	&=	\kappa (\theta - X_t) \dd t + \dd \eta_t
\end{align} 
(a) Find $X_t$ explicitly as a functional of $\eta$.\\
(b) Assume $\eta_t = \sig W_t + \int_\Rb z \Nt(t, \dd z)$.  Compute $m(t) := \Eb X_t$ and $c(t,s) := \Eb (X_t - m(t))(X_s-m(s))$.
\end{exercise}

\begin{exercise}
\label{exer:10.5}
Let $X$ be the following one-dimensional jump-diffusion
\begin{align}
\dd X_t
	&=	\mu(t,X_t) \dd t + \sig(t,X_t) \dd W_t + \int_\Rb \gam(t,X_{t-},z) \Nt(\dd t,\dd z) ,
\end{align}
where $W$ is a one-dimensional Brownian motion and $\Nt$ is a one-dimensional compensated Poisson random measure on $\Rb$.  Derive using the L\'evy-It\^o formula the infinitesimal generator $\Ac(t)$ of the $X$ process
\begin{align}
\Ac(t)\varphi(x)
	&:=	\lim_{s \searrow t} \frac{\Eb[ \varphi(X_s)| X_t = x] - \varphi(x)}{s-t} .
\end{align}
\end{exercise}

\begin{exercise}
\label{exer:10.6}
Consider a financial market defined on a probability space $(\Om,\Fc,\Pbt)$ consisting of two assets: a stock $S$ and a money market account $M$ with
\begin{align}
M_t
	&= M_0 \ee^{r t} , &
S_t
	&=	S_0 \ee^{X_t} , &
X_t
	&=	\mu t + \sig \Wt_t + \int_\Rb z \Nt(t, \dd z) ,
\end{align}
where $\Wt$ is a Brownian motion under $\Pbt$ and $\Nt(\dd t, \dd z) = N(\dd t, \dd z) - \nu( \dd z) \dd t$ with $\Ebt N(\dd t, \dd z) = \nu(\dd z) \dd t$.\\
(a) Suppose that $\Pbt$ is the market's pricing measure with the money market account $M$ as num\'eraire.  Then $S/M$ must be a martingale under $\Pbt$.  Find $\mu$. \\
(b) Let $V_t$ be the price of a derivative asset that pays $S_T^p$ at time $T$.  Assuming that $V/M$ is also a martingale under $\Pbt$, show that $V_t = v(t,S_t)$ and find the function $v$.
\end{exercise}

\begin{exercise}
\label{exer:10.7}
Let $P = (P_t)_{0 \leq t \leq T}$ be a Poisson process with intensity $\lam$ under $\Pb$.  Under what change of measure $\dd \Pbh / \dd \Pb$ is $P$ a Poisson process with intensity $\widetilde{\lam}$ under $\Pbh$?
\end{exercise}

\begin{exercise}
\label{exer:10.8}
Define processes $X=(X_t)_{t \geq 0}$, $Y=(Y_t)_{t \geq 0}$, and $Z=(Z_t)_{t \geq 0}$ as follows
\begin{align}
X_t 
	&:= P_t^{(1)} , &
Y_t 
	&:= \mu t + P_t^{(2)} , &
Z_t 
	&:= X_{Y_t} 
\end{align}
where $P^{(1)}$ and $P^{(2)}$ are independent Poisson processes with intensities $\lambda_1$ and $\lambda_2$, respectively, and $\mu > 0$.
\\[0.5em]
(a) Find the characteristic functions of $\phi_{X_t}(u) := \Eb \ee^{\ii u X_t}$ and $\phi_{Y_t}(u) := \Eb \ee^{\ii u Y_t}$, where $u \in \Rb$.
\\[0.5em]
(b) Find the characteristic function of $\phi_{Z_t}(u) := \Eb \ee^{\ii u Z_t}$, where $u \in \Rb$. 
\\[0.5em]
(c) Define the probability generating function $G_{Z_t}(s) := \Eb s^{Z_t}$ where $|s|< 1$.  Use $G_{Z_t}$ to find $\Pb( Z_t = 0)$.
\\[0.5em]
Hint: write $\Eb s^{Z_t}$ as an infinite sum and compare to the Taylor series of $G_{Z_t}(s)$ about the point $s=0$.
\end{exercise}

\begin{exercise}
\label{exer:10.9}
Let $X = (X_t)_{t \geq 0}$ satisfy
\begin{align}
\dd X_t
	&=	\mu X_t \dd t + \sig X_t \dd W_t + \int X_{t-} (\ee^{\gam(z)}-1) \Nt(\dd t, \dd z) .
\end{align}
Define $Y_t := X_t^p$ where $p \in \Rb$.  Compute $\dd Y_t$.  Your answer should \textit{not} involve $X$.
\end{exercise}

\begin{exercise}
\label{exer:10.10}
Consider the following L\'evy-It\^o process
\begin{align}
\dd S_t
	&=	\sig_t S_t \dd W_t + S_{t-} \int_\Rb (\ee^z - 1 ) \Nt(\dd t, \dd z) ,
\end{align}
where $W$ is a Brownian motion, $\Nt$ is a compensated Poisson random measure with associated L\'evy measure $\nu$, and the process $\sig$ is continuous and independent of $W$ and $N$.  You may assume $\nu(\Rb)<\infty$.
Let $\Fb = (\Fc_t)_{t \geq 0}$ be a filtration for $W$ and $\int_\Rb z N( \cdot , \dd z)$ and assume $\sig$ is $\Fb$-adapted.  Define
\begin{align}
X_t
	&:=	\log S_t , &
Y_t
	&:=	[X,X]_t .
\end{align}
(a) Find expressions for $\dd X_t$ and $\dd Y_t$.
\\[0.5em]
(b) Define a process $Q = (Q_t)_{t \geq 0}$ as follows
\begin{align}
Q_t
	&:=	\Eb( S_T^q | \Fc_t ) , &
0
	&\leq t \leq T ,
\end{align}
For a general process $Z = (Z_t)_{t \geq 0}$ define the jump of $Z$ at time $t$ as $\Del Z_t := Z_t - Z_{t-}$.
Compute $\Del Q_t$.  It may help to note that
\begin{align}
\Eb( \ee^{q (X_T^j -  X_{t}^j)} | \Fc_{t} )
	&=	\Eb( \ee^{q (X_T^j -  X_{t-}^j)} | \Fc_{t-} ) ,
\end{align}
where $X^j$ is the jump component of $X$.
Please give your answer in terms of $Q_{t-}$ and an integral of the form $\int (\cdots ) \Del N(t, \dd z)$.
\\[0.5em]
(c) Consider the case $\nu \equiv 0$.  Is there a function $u : \Rb_+ \times \Rb \to \Rb$ such that $u(t,X_t) = Q_t$?  If so, find this function.  If not, explain why such a function does not exist in general.
\\[0.5em]
(d) Consider the case $\sig \equiv 0$.  Is there a function $u : \Rb_+ \times \Rb \to \Rb$ such that $u(t,X_t) = Q_t$?  If so, find this function.  If not, explain why such a function does not exist in general.
\end{exercise}

\begin{exercise}
\label{exer:10.11}
Suppose that $X = (X_t)_{t \geq 0}$ is a non-decreasing L\'evy process with a L\'evy measure $\nu$ satisfying 
\begin{align}
\int_\Rb |z| \nu(\dd z)<\infty .
\end{align}
(a) Compute $\phi(\lam) := -(1/t) \log \Eb \ee^{- \lam X_t}$ where $\lam > 0$.
\\[0.5em]
(b) Define $\tau_m := \inf \{ t \geq 0 : W_t = m \}$ where $m > 0$ and $W = (W_s)_{s \geq 0}$ is a Brownian motion.  Show that the process $\tau = (\tau_t)_{t \geq 0}$ is a non-decreasing L\'evy process and give an expression for its L\'evy measure $\nu(\dd z)$.  You may find the following useful
\begin{align}
v^r
	&=	\frac{r}{\Gam(1-r)} \int_0^\infty \frac{1}{z^{1+r}}(1 - \ee^{-zv}) \dd z , &
r
	&\in (0,1) ,
\end{align}
where $\Gam$ is an Euler-Gamma function.  (You may leave $\Gam$ in your answer; you do not need to evaluate it.)
\end{exercise}

\begin{exercise}
\label{exer:10.12}
Fix a probability space $(\Om,\Fc,\Pb)$.
Let $X = (X_t)_{t \geq 0}$ be the unique strong solution of the following L\'evy-It\^o SDE
\begin{align}
\dd X_t
	&=	\mu(t,X_t) \dd t + \sig(t,X_t) \dd W_t + \int_\Rb \gam(t,X_{t-},z) \Nt(\dd t, \dd z) ,
\end{align}
where $W$ is a Brownian motion and $\Nt( \dd t, \dd z) := N(\dd t, \dd z) - \nu(\dd z) \dd t$ is a compensated Poisson random measure with associated L\'evy measure $\nu$.  Let $\Fb = (\Fc_t)_{t \geq 0}$ be the filtration generated by observing $X$ (so that $X_t \in \Fc_t$ for all $t \geq 0$).
Define
\begin{align}
u(t,X_t)
	&:=	\Eb \Big( \ee^{- \lam (T-t)} h(X_T) + \int_t^T \ee^{- \lam (s-t)} g(X_s) \dd s \Big| \Fc_t \Big) ,
\end{align}
where $\lam \geq 0$ is a constant.
\\[0.5em]
(a) Define $M = (M_t)_{t \geq 0}$ by
\begin{align}
M_t
	&=	\ee^{- \lam t} u(t,X_t) + \int_0^t \ee^{- \lam s} g(X_s) \dd s .
\end{align}
Show that $M$ is a martingale.
\\[0.5em]
(b) Derive a partial integro-differential equation (PIDE) and terminal condition (TC) for the function $u$.
\end{exercise}

%
%

\chapter{Stochastic Control}
Notes from this chapter are taken primarily from \cite[Chapter 19]{bjork2004arbitrage} and \cite[Chapter 6]{vanhandel}.

\section{Problem formulation}
\label{sec:formulation}
Fix a probability space $(\Om,\Fc,\Pb)$ and a filtration $\Fb=(\Fc_t)_{t \geq 0}$ and consider the following \emph{controlled diffusion} $X^u = (X_t^u)_{t \geq 0}$, whose dynamics are given by
\begin{align}
\dd X_t^u
	&=	\mu(t,X_t^u,u_t) \dd t + \sig(t,X_t^u,u_t) \dd W_t , &
X_0^u
	&=	x \in \Rb^d , \label{eq:dXu}
\end{align}
where $W$ is an $n$-dimensional Brownian motion with respect to $(\Pb,\Fb)$ and
\begin{align}
\mu
	&:	\Rb_+ \times \Rb^d \times \Ub \to \Rb^d , &
\sig
	&:	\Rb_+ \times \Rb^d \times \Ub \to \Rb^{d \times n}
\end{align}
Here, $X^u$ is called the \emph{state process} and $u=(u_t)_{t \geq 0}$ is the \emph{control}.  The control must live in some \emph{control set} $\Ub$.  Frequently, we take the control set $\Ub = \Rb^k$, but this is not required.  The superscript $u$ on $X^u$ indicates the dependence of state process $X^u$ on the control $u$ (clearly, if you change $u$ you change $X^u$).

\begin{definition}
We call a control $u$ is \emph{admissible} if the following holds:
\begin{enumerate*}
\item the control is $\Fb$-adapted, that is, $u_t \in \Fc_t$ for all $t \geq 0$,
\item the control remains in the control set: $u_t(\om) \in \Ub$ for all $(t,\om)$,
\item the process $X^u$ is well-defined.
\end{enumerate*}
We denote by $\Uc$ the set of admissible controls.
\end{definition}

\noindent
Within the class of admissible controls is a subset of controls that we will find very useful:

\begin{definition}
We call an admissible control $u$ a \emph{feedback} or \emph{Markov} control if it is of the form: $u_t = \alpha(t,X_t^u)$ for some function $\alpha : \Rb_+ \times \Rb^d \to \Ub$. 
\end{definition}

\noindent
Note Markov controls are a strict subset of admissible controls, since a Markov control at time $t$ is only allowed to depend on $(t,X_t^u)$ whereas, in general, an admissible control at time $t$ may depend on the entire history $\Fc_t$.  Nevertheless, Markov controls are important because they are easy to implement, and (as we shall see) they can actually be computed!  Moreover, it often turns out that the optimal Markov control coincides with the optimal admissible control.
\\[0.5em]
If we consider a Markov control $u_t = \alpha(t,X_t^u)$ then $X^u$ satisfies an SDE
\begin{align}
\dd X_t^u
	&=	\mu(t,X_t^u,\alpha(t,X_t^u)) \dd t + \sig(t,X_t^u,\alpha(t,X_t^u)) \dd W_t , &
X_0^u
	&=	x \in \Rb^d .
\end{align}
Since the solution of an SDE is a Markov process, a Markov control yields a Markov state process $X^u$.
\\[0.5em]
Below, we introduce some \emph{cost} or \emph{gain functionals} $J(\cdot)$ that assign to each admissible control strategy $u$ a cost or gain $J(u)$. The goal of \emph{optimal control} is to find an \emph{optimal strategy} $u^*$ (obviously!), which minimizes or maximizes this cost functional.  
Three common types of cost functionals are:
\begin{enumerate*}
\item \emph{Finite time}.  Here, the cost functional is given by
\begin{align}
J(u)
	&=	\Eb \Big[ \int_0^T F(t,X_t^u,u_t) \dd t + \varphi(X_T^u) \Big] , \label{eq:finite.time}
\end{align}
where $F: [0,T] \times \Rb^d \times \Ub \to \Rb$ and $\varphi : \Rb^d \to \Rb$.
\item \emph{Indefinite time}.  Here, the cost functional is given by
\begin{align}
J(u)
	&=	\Eb \Big[ \int_0^{\tau^u} F(X_t^u,u_t) \dd t +  \varphi(X_{\tau^u}^u) \Big], &
\tau^u
	&=	\inf\{ t \geq 0 : X_t^u \notin A \} ,
\end{align}
where $A \subset \Rb^d$ is some open set, $F: A \times \Ub \to \Rb$ and $\varphi : \d A \to \Rb$.
\item \emph{Infinite time}.  In this case, the cost functional is given by
\begin{align}
J(u)
	&=	\Eb \int_0^\infty \ee^{- \lam t} F(X_t^u,u_t) \dd t ,
\end{align}
where $F : \Rb^d \times \Ub \to \Rb$ and $\lam > 0$.
\end{enumerate*}
Note, for the Indefinite time and Infinite time functionals, the dynamics of $X^u$ should be time-homogeneous (i.e., the coefficient functions $\mu$ and $\sig$ should not depend on $t$).
\\[0.5em]
Of course, one might construct other cost functionals, but the above three are by far the most common.  As stated above, the optimal control (if it exists) is the strategy $u^*$ that minimizes or maximizes a given cost or gain functional $J(u)$.  Thus, we make the following definition:

\begin{definition}
For a given cost functional $J : \Uc \to \Rb$, we define the \emph{optimal control} $u^*$, if it exists, by 
\begin{align}
u^*
	&:=	\operatorname*{argmax}_{u \in \Uc} J(u) , &
	&\text{if our goal is to maximize $J(u)$}, \\
u^*
	&:=	\operatorname*{argmin}_{u \in \Uc} J(u) , &
	&\text{if our goal is to minimize $J(u)$}.
\end{align}
\end{definition}

\section{The Dynamic Programming Principle and the HJB PDE}
In this Section we focus on the finite time horizon control problem and we shall attempt to find the control $u^*$ that \emph{maximizes} $J(u)$ in \eqref{eq:finite.time}.  The minimization procedure is completely analogous.  We shall restrict our controls to Markov controls: $u_t = \alpha(t,X_t^u)$.  
We shall make the following assumption.

\begin{assumption}
\label{ass:markov}
We assume that there exists an admissible Markov strategy $u^*$, which is optimal.
\end{assumption}

\noindent
Clearly, this assumption is \emph{not} always justified.  Nevertheless, we will go with it for the time being.  
We define the \emph{reward-to-go} function $J^u$ (or \emph{cost-to-go} in the case of minimization) as
\begin{align}
J^u(t,X_t^u)
	&:=	\Eb \Big[ \int_t^T F(s,X_s^u,u_s) \dd s + \varphi(X_T^u) \Big| \Fc_t \Big] 
	=		\Eb \Big[ \int_t^T F(s,X_s^u,u_s) \dd s + \varphi(X_T^u) \Big| X_t^u \Big] ,
\end{align}
where we have used the Markov property of $X^u$ to replace $\Fc_t$ with the $\sig$-algebra generated by $X_t^u$.  
We also define the \emph{value function} $V$ as
\begin{align}
V(t,x)
	&:=	J^{u^*}(t,x) = \max_{u \in \Uc} J^u(t,x) . \label{eq:Value-Function}
\end{align}
The idea of the \emph{Dynamic Programming Principle} (DPP) is to split the optimization problem into two intervals $[t,t+\del)$ and $[t+\del,T]$, where $\del$ is small and positive.  Note that, for any $\del \geq 0$ we have
\begin{align}
J^u(t,x)
	&=	\Eb \Big[ \int_t^{t + \del} F(s,X_s^u,u_s) \dd s 
			+ \Eb \Big[ \int_{t + \del}^T F(s,X_s^u,u_s) \dd s + \varphi(X_T^u) \Big| \Fc_{t+\del} \Big] 
			\Big| X_t^u = x \Big] \\
	&=	\Eb \Big[ \int_t^{t + \del} F(s,X_s^u,u_s) \dd s 
			+ J^u(t+\del,X_{t+\del}^u) 
			\Big| X_t^u = x \Big] 
\end{align}
Now, fix an arbitrary control $\uh$ and set
\begin{align}
u_s
	&=	\left\{ \begin{array}{ll}
			\uh_s 	& s \in [t,t+\del) , \\
			u_s^* & s \in [t+\del, T] .
			\end{array} \right.
\end{align}
In words, the control $u$ may be sub-optimal over the interval $[t,t+\del)$ and it is optimal over the interval $[t+\del,T]$.  Clearly, we have
\begin{align}
V(t,x)
	&\geq	\Eb \Big[ \int_t^{t + \del} F(s,X_s^u,u_s) \dd s 
			+ V(t+\del,X_{t+\del}^u) 
			\Big| X_t^u = x \Big] . \label{eq:V.1}
\end{align}
The inequality arises from the fact that the strategy $u$ is not necessarily optimal over in interval $[t,t+\del)$.  If we had $\uh = u^*$, then we would have obtained an equality since, in this case, we have $u = u^*$.  Now observe that
\begin{align}
V(t+\del,X_{t+\del}^u)
	&=	V(t,X_t^u) + \int_t^{t+\del} ( \d_s + \Ac^u(s) ) V(s,X_s^u) \dd s + \text{martingale part} , 
\end{align}
where $\Ac^u(s)$ is the infinitesimal generator of $X^u$.  
Taking an expectation, we have
\begin{align}
\Eb[ V(t+\del,X_{t+\del}^u) | X_t^u = x ]
	&=	V(t,x) + \Eb \Big[  \int_t^{t+\del} ( \d_s + \Ac^u(s) )V(s,X_s^u) \dd s \Big| X_t^u = x \Big] . \label{eq:V.2}
\end{align}
Inserting \eqref{eq:V.2} into \eqref{eq:V.1} and obtain
\begin{align}
0
	&\geq \Eb \Big[ 
			\int_t^{t + \del} \Big( F(s,X_s^u,u_s) + ( \d_s + \Ac^u(s) ) V(s,X_s^u)  \Big) \dd s
			\Big| X_t^u = x \Big] .
\end{align}
Finally, we divide by $\del$ and take a limit as $\del \to 0$.  Assuming there is no problem with passing a limit through an expectation, we have
\begin{align}
0
	&\geq \lim_{\del \to 0} \frac{1}{\del} \Eb \Big[ 
			\int_t^{t + \del} \Big( F(s,X_s^u,u_s) + ( \d_s + \Ac^u(s) ) V(s,X_s)  \Big) \dd s
			\Big| X_t^u = x \Big] \\
	&=	\Eb \Big[ \lim_{\del \to 0} \frac{1}{\del} 
			\int_t^{t + \del} \Big( F(s,X_s^u,u_s) + ( \d_s + \Ac^u(s) ) V(s,X_s)  \Big) \dd s
			\Big| X_t^u = x \Big] \\
	&=	F(t,x,u) + (\d_t + \Ac^u(t)) V(t,x) .
\end{align}
Once again, if $u$ is optimal, we obtain an equality.  Thus, we have arrived at the \emph{Hamilton-Jacobi-Bellman (HJB) PDE}
\begin{align}
0
	&=	\d_t V + \max_{u \in \Ub} \( F^u(t) + \Ac^u(t) V \) , &
V(T,x)
	&=	\varphi(x) , \label{eq:HJB.PDE}
\end{align}
where we have defined
\begin{align}
F^u(t)
	&=	F(t,x,u) , &
\Ac^u(t)
	&=	\sum_{i=1}^d \mu_i(t,x,u) \d_{x_i} + \tfrac{1}{2} \sum_{i=1}^d \sum_{j=1}^d \( \sig \sig^\text{T} \)_{i,j}(t,x,u) \d_{x_i} \d_{x_j} .
\end{align}
We summarize this result with the following Theorem

\begin{theorem}
\label{thm:hjb.I}
If (i) an optimal control $u^*$ exists and is Markov, (ii) the value function $V$, defined in \eqref{eq:Value-Function}, satisfies $V \in C^{1,2}$ and (iii) the limiting procedures performed above are allowed, then $V$ satisfies the HJB PDE \eqref{eq:HJB.PDE} and the optimal control $u^*$ is given by
\begin{align}
u_t^*
	&=	\alpha(t,X_t^{u^*}) , &
	&\text{where}&
\alpha(t,x)
	&:=	\operatorname*{argmax}_{u \in \Ub} \( F(t,x,u) + \Ac^u V(t,x) \) . \label{eq:u.star}
\end{align}
\end{theorem}

\noindent
Theorem \ref{thm:hjb.I} is so incredibly unsatisfactory!  It only tells us that, \textit{if} an optimal controls exists and is Markov, then (modulo some technical conditions) the value function $V$ satisfies the HJB PDE and the optimal control is given by \eqref{eq:u.star}.  But, this is not what we want.  What we would \textit{like}, is solve the HJB PDE \eqref{eq:HJB.PDE} and \textit{conclude} that the solution actually \emph{is} the value function and $u^*$, given by \eqref{eq:u.star}, is the optimal control.  The following Theorem gives us such a result.

\begin{theorem}[Verification Theorem]
Suppose $H \in C^{1,2}([0,T] \times \Rb^d)$ solves the HJB PDE \eqref{eq:HJB.PDE}, including the terminal condition $H(T,\cdot) = \varphi$.
Suppose the maximizer
\begin{align}
g(t,x)
	=	\operatorname*{argmax}_{u \in \Ub} \( F(t,x,u) + \Ac^u H(t,x) \)
\end{align}
exists and that the feedback control $u_t^* := g(t,X_t^{u^*})$ is admissible.  Suppose further that, for every admissible control $u \in \Uc$, the process
\begin{align}
M_t^u
	&:=	H(t,X_t^u) - \int_0^t ( \d_s + \Ac^u(s) )H(s,X_s^u) \dd s
\end{align}
is a true martingale (not just a local martingale).  Then the value function $V$ and an optimal control are given by
\begin{align}
V(t,x)
	&=	H(t,x) , &
u_t^*
	&=	g(t,X_t^{u^*}) .
\end{align}
\end{theorem}

\noindent
The Verification Theorem tells us that, if we solve the HJB PDE, and the solution satisfies some regularity and integrability conditions, then the solution \emph{is} the value function and the Markov control \emph{is} the optimal control.  We will not prove the verification theorem.

\section{Solving the HJB PDE}
Solving the HJB PDE essentially involves two steps
\begin{align}
\text{Step 1 (easy part)}:&&
\text{find}&&
u^*(t,x)
	&=	\operatorname*{argmax}_{u \in \Ub} ( F(t,x,u) + \Ac^u V(t,x) ) . \\
\text{Step 2 (hard part)}:&&
\text{solve}&&
0
	&=	( \d_t + \Ac^{u^*}(t) ) V + F^{u^*} , &
V(T,\cdot)
	&=	\varphi .
\end{align}
We illustrate this process with a simple, explicitly-solvable example.

\begin{example}[Merton problem]
\label{ex:merton}
The \emph{Merton problem}, due to Nobel Prize winner Robert Merton, is a classical problem in mathematical finance whereby an investor seeks to optimize his expected utility at a fixed future date $T$ by investing in a stock.  Specifically, suppose a stock $S=(S_t)_{0 \leq t \leq T}$ follows a geometric Brownian motion
\begin{align}
\dd S_t
	&=	\mu S_t \dd t + \sig_t S_t \dd W_t .
\end{align}
Let $X^u = (X_t^u)_{0 \leq t \leq T}$ be the wealth process of an investor who invests $u_t$ dollars in $S$ at time $t$ and keeps the rest of his money in a bank account (which we assume pays no interest).  The dynamics of $X^u$ are given by
\begin{align}
\dd X_t^u
	&=	\frac{u_t}{S_t} \dd S_t 
	=		\mu u_t \dd t + \sig_t u_t \dd W_t .
\end{align}  
The investor's control problem is to maximize
\begin{align}
J(u)
	&=	\max_u \Eb \, U(X_T^u) ] , &
U(x)
	&=	\frac{1}{1-\gam} x^{1-\gam} , &
\gam
	&\in (0,1) \cup (1,\infty) . \label{eq:control.problem}
\end{align}
The function $U$ is called the investor's \emph{utility function}; it is intended to map the investor's wealth to his happiness.  Since ``more money'' $=$ ``more happy'' the function $U$ is strictly increasing.  It is also concave since, an additional dollar means less to somebody with a wealth of 1 million dollars than it does to somebody with a wealth of 10 dollars.  For the investor's control problem \eqref{eq:control.problem}, we identify
\begin{align}
F^u
	&=	0 , &
\Ac^u
	&=	u \mu \d_x + \tfrac{1}{2} u^2 \sig^2 \d_x^2 , 
\end{align}
In order to maximize $(F^u + \Ac^u V)$ we simply differentiate this quantity, set it equal to zero, and solve for $u$.  We have
\begin{align}
0
	&=	\d_u \( F^u + \Ac^u V \) |_{u^*}
	=		\mu \d_x V + u \sig^2 \d_x^2 V |_{u^*} &
	&\Rightarrow&
u^*
	&=	\frac{- \mu \d_x V}{\sig^2 \d_x^2 V} .
\end{align}
Next, we try to solve
\begin{align}
0
	&=	( \d_t + \Ac^{u^*} ) V + F^{u^*} \\
	&=		\d_t V + \( - \frac{\mu \d_x V}{\sig^2 \d_x^2 V} \) \mu \d_x V 
			+ \tfrac{1}{2} \( - \frac{\mu \d_x V}{\sig^2 \d_x^2 V} \)^2 \sig^2 \d_x^2 V \\ 
	&=	\d_t V - \frac{1}{2} \lam^2 \frac{(\d_x V)^2}{\d_x^2 V}, \qquad
\lam
	=	\frac{\mu}{\sig} , &
V(T,x)
	&=	U(x) . \label{eq:HJB.investment}
\end{align}
where we have introduced the \emph{Sharpe ratio} $\lam$.  Of course, if you have never seen the horrible looking non-linear PDE before you would never know how to solve it.  However, if you have spent some time around optimal investment problems you would know that the correct thing to do at this point is to guess that $V$ is of the form
\begin{align}
V(t,x)
	&=	f(t) U(x) , &
f(T)
	&=	1 . \label{eq:guess}
\end{align}
Such a guess clearly satisfies the terminal condition $V(T,x) = U(x)$.  If we insert the guess \eqref{eq:guess} into \eqref{eq:HJB.investment} we find
\begin{align}
0
	&=	\Big( f' + \tfrac{1}{2} \lam^2 \tfrac{1-\gam}{\gam} f \Big) U ,
\end{align}
which, together with the terminal condition $f(T) = 1$ implies that
\begin{align}
f(t)
	&=	\exp \Big( \tfrac{1}{2} \lam^2 \tfrac{1-\gam}{\gam}(T-t) \Big) .
\end{align}
Thus, we have obtained the value function $V(t,x) = U(x) f(t)$.  The optimal Markov control is given by
\begin{align}
u^*(t,x)
	&=	- \frac{\mu}{\sig} \frac{\d_x V(t,x)}{\d_x^2 V(t,x)}
	=		\frac{ \mu x }{\gam \sig^2} .
\end{align}
Thus, the total proportion of wealth the investor should keep in the stock is constant: $\frac{u_t^*}{X_t^{u^*}} = \frac{\mu}{\gam \sig^2}$.
\end{example}

\begin{example}[Linear quadratic regulator]
Suppose the dynamics of a controlled process $X^u = (X_t^u)_{0 \leq t \leq T}$ are given by
\begin{align}
\dd X_t^u
	&=	( a X_t^u + b u_t ) \dd t + c \dd W_t ,
\end{align}
and consider the following cost functional, which we wish to minimize
\begin{align}
J(u)
	&=	\Eb \Big( \int_0^T ( q (X_t^u)^2 + r u_t^2 ) \dd t + h (X_T^u)^2 \Big) . \label{eq:J-LQR}
\end{align}
One could imagine, for example, that $X^u$ represents the position of a particle that we are attempting to keep near the origin.  We can adjust the drift of the  particle through the control $u$, but this has a cost, which is quadratic in $u$.  Likewise, there is a quadratic cost for allowing $X^u$ to be away from the origin.  This control problem is known as the \textit{Linear Quadratic Regulator} (LQR) because the dynamics of $X$ are linear in the state process $X^u$ and control process $u$ and the costs are quadratic in $X^u$ and $u$.  The HJB equation associated with cost functional \eqref{eq:J-LQR} is
\begin{align}
0
	&=	V_t + \sup_{u \in \Rb} ( \Ac^u V + F^u ) , &
V(t,x)
	&=	h x^2 ,
\end{align}
where the operator $\Ac^u$ and the function $F^u$ are given by
\begin{align}
\Ac^u
	&=	( a x + b u ) \d_x + \tfrac{1}{2} c^2 \d_x^2 , &
F^u
	&=	q x^2 + r u^2 .
\end{align}
Step one in solving the HJB PDE is finding the optimal control $u^*$ in feedback form.  To this end, we have
\begin{align}
0
	&=	\d_u \Big[ \Ac^u V + F^u \Big|_{u = u^*}  =	b V_x + 2 r u^* , &
	&\Rightarrow&
u^*
	&=	\frac{- b V_x}{2 r }	 .
\end{align}
The HJB PDE thus becomes
\begin{align}
0
	&=	V_t + \Ac^{u^*} V + F^{u^*} \\
	&=	V_t + \Big[ a x + b \Big( \frac{- b V_x}{2 r } \Big) \Big] V_x + \tfrac{1}{2} c^2 V_{xx} + q x^2 + r \Big( \frac{- b V_x}{2 r } \Big)^2 \\
	&=	V_t + a x V_x + \tfrac{1}{2} c^2 V_{xx} + q x^2 - \frac{ b^2 V_x^2}{4 r } . \label{eq:nonlin.pde}
\end{align}
As with all non-linear PDEs with explicit solutions, the method to solve \eqref{eq:nonlin.pde} is to guess.  In this case, the correct guess is
\begin{align}
V(t,x)
	&=	P(t) x^2 + Q(t) , &
P(T)
	&=	h , &
Q(T)
	&=	0 ,
\end{align}
where the terminal conditions for $P$ and $Q$ will ensure that $V(T,x) = h x^2$.  We have
\begin{align}
V_t
	&=	P' x^2 + Q' , &
V_x
	&=	2 x P , &
V_{xx}
	&=	2 P . \label{eq:V.derivatives}
\end{align}
Inserting \eqref{eq:V.derivatives} into \eqref{eq:nonlin.pde} we obtain
\begin{align}
0
	&=	(P' x^2 + Q') + a x ( 2 x P ) + \tfrac{1}{2} c^2 (2 P)  + q x^2 - \frac{ b^2 }{4 r } (2 x P )^2 \\
	&=	x^2 \Big( P' + 2 a P + q - (b^2/r) P^2 \Big)
			+ \Big( Q' + c^2 P \Big) .
\end{align}
Collecting terms of like order in $x$ we obtain
\begin{align}
\Oc(x^0):&&
0
	&=	Q' + c^2 P , &
Q(T)
	&=	0 , \\
\Oc(x^2):&&
0
	&=	P' + 2 a P + q - (b^2/r) P^2 , &
P(T)
	&=	h .
\end{align}
Thus, we have obtained a coupled system of ODEs for $(P,Q)$.  The $\Oc(x^2)$ equation is a \textit{Riccati equation}, which can be solved analytically (though, the solution is a tad messy and not worth writing down here).  Once one has obtained an expression for the function $P$ one can obtain $Q$ as an integral
\begin{align}
Q(t)
	&=	c^2 \int_t^T P(s) \dd s .
\end{align}
Finally, the optimal control is given by
\begin{align}
u_t^*
	&=	u^*(t,X_t^{u^*})
	=		\frac{- b }{2 r } V_x(t,X_t^{u^*})
	=		\frac{ - b }{r} P(t) X_t^{u^*} .
\end{align}
\end{example}

\section{HJB equations associated with other cost functionals}
\label{sec:other-HJB}
In this section we suppose the dynamics of a controlled diffusion $X^u$ are time-homogeneous
\begin{align}
\dd X_t^u
	&=	\mu( X_t^u, u_t) \dd t + \sig(X_t^u, u_t) \dd W_t . \label{eq:dXu}
\end{align}
In what follows, we state without proof the HJB PDEs associated with  the indefinite time and infinite time cost functionals $J(u)$, which were introduced in Section \ref{sec:formulation}.
\begin{enumerate*}
\item Consider the \emph{Indefinite time} cost functional
\begin{align}
J(u)
	&=	\Eb \Big[ \int_0^{\tau^u} F(X_t^u,u_t) \dd t +  \varphi(X_{\tau^u}^u) \Big], &
\tau^u
	&=	\inf\{ t \geq 0 : X_t^u \notin A \} ,
\end{align}
where $A \subset \Rb^d$ is some open set, $F: A \times \Ub \to \Rb$ and $\varphi : \d A \to \Rb$.
The associated HJB PDE is
\begin{align}
0
	&=	\sup_{u \in \Ub} ( \Ac^u V + F^u ) , &
x
	&\in A , \label{eq:HJB-indefinite-time} \\
V(x)
	&=	\varphi(x) , &
x
	&\in \d A .
\end{align}
where the operator $\Ac^u$ is the generator of the process $X^u$ defined in \eqref{eq:dXu} with $u_t = u$ fixed.  In the one-dimensional case, we have
\begin{align}
\Ac^u
	&=	\mu(x,u)\d_x + \tfrac{1}{2} \sig^2(x,u) \d_x^2 .
\end{align}
\item Consider the \emph{Infinite time} cost functional
\begin{align}
J(u)
	&=	\Eb \int_0^\infty \ee^{- \lam t} F(X_t^u,u_t) \dd t ,
\end{align}
where $F : \Rb^d \times \Ub \to \Rb$ and $\lam > 0$.
The associated HJB PDE is
\begin{align}
0
	&=	\sup_{u \in \Ub} ( \Ac^u V  +  F^u ) - \lam V . \label{eq:HJB-infinite-time}
\end{align}
\end{enumerate*}
We leave the formal proof of the above HJB equations as an exercise for the reader.


\section{Exercises}
\label{sec:exercises-11}

\begin{exercise}
\label{exer:11.1}
Derive the HJB equations in Section \ref{sec:other-HJB}.
\end{exercise}

\begin{exercise}
\label{exer:11.2}
[Optimal investment with consumption]
Consider the setting of the Merton problem in Example \ref{ex:merton}, but suppose now that the investor, in addition to investing $u_t$ dollars in the stock $S$ at time $t$, also consumes his wealth at a rate $c_t \geq 0$.  The dynamics of the investor's wealth $X^{u,c} = (X_t^{u,c})_{t \geq 0}$ are then
\begin{align}
\dd X_t^{u,c}
	&=	\frac{u_t}{S_t} \dd S_t - c_t \dd t
	=		( \mu u_t - c_t ) \dd t + \sig u_t \dd W_t .
\end{align}
The investor seeks to maximize his expected discounted utility of consumption over an infinite horizon
\begin{align}
J(u,c)
	&=	\Eb \int_0^\infty \ee^{- \lam t} U(c_t) \dd t , &
U(c)
	&=	\frac{1}{1-\gam} c^{1-\gam} , &
\gam
	&\in (0,1) \cup (1,\infty) ,
\end{align}
where $\lam > 0$.  Find the investor's value function
\begin{align}
V(x)
	&:=	\sup_{u,c} \Eb \Big[ \int_0^\infty \ee^{- \lam t} U(c_t) \dd t \Big| X_0^{u,c} = x \Big] ,
\end{align}
and the optimal controls $u_t^*$ and $c_t^*$.  You may assume that
\begin{align}
\nu
	&:=	\frac{1}{\gam} \Big( \lam - \frac{1-\gam}{\gam} \frac{\mu^2}{2 \sig^2} \Big)
	>	0 .
\end{align}
\textit{Hint}: use the HJB equation for the infinite time cost functional in Section \ref{sec:other-HJB} and try a solution of the form $V(x) = A \, U(x)$ for some constant $A > 0$.
\end{exercise}

%
%

\chapter{Optimal Stopping}
\label{ch:optimal.stopping}
Notes from this chapter closely follow \cite[Chapter 9]{vanhandel}.

\section{Strong Markov processes}
\label{sec:strong.markov}
We have previously defined the notion of a Markov process.  For optimal stopping, we will need to consider a subset of these processes called \emph{strong Markov processes}.

\begin{definition}
\label{def:strong.markov}
Let  $(\Om,\Fc,\Pb)$  be a probability space, let $T$ be a fixed positive number, and let $\Fb=(\Fc_t)_{t \geq 0}$ be a filtration of sub-$\sig$-algebras of $\Fc$.
Consider an $\Fb$-adapted stochastic process $X=(X_t)_{t \geq 0}$. 
We say $X$ is a \emph{strong Markov process} or simply ``$X$ is strong Markov,'' if, for any $\Fb$ stopping time $\tau$ and $t \geq 0$ we have
\begin{align}
\Eb[ f(X_{\tau+t}) | \Fc_\tau ]
	&=	\Eb[ f(X_{\tau+t}) | X_\tau ] .
\end{align}
\end{definition}

\begin{proposition}
Suppose that on a probability space $(\Om,\Fc,\Pb)$ equipped with a filtration $\Fb = (\Fc_t)_{t \geq 0}$ the process $X$ is the unique strong solution of the time-homogeneous SDE
\begin{align}
\dd X_t
	&=	\mu(X_t) \dd t + \sig( X_t) \dd W_t , &
X_0
	\in \Rb^d ,
\end{align}
where $W = (W_t)_{t \geq 0}$ is a Brownian motion in $\Rb^n$ and $\mu : \Rb^d \to \Rb^d$ and $\sig : \Rb^d \to \Rb^{d \times n}$ satisfy \eqref{eq:lin.growth}-\eqref{eq:lipschitz.continuous}.  Then $X$ is a strong Markov process.
\end{proposition}

\section{Optimal Stopping Problem Formulation}
\label{sec:optimal.stopping}
On a probability space $(\Om,\Fc,\Pb)$ equipped with a filtration $\Fb = (\Fc_t)_{t \geq 0}$ we consider a time-homogeneous diffusion $X = (X_t)_{t \geq 0}$, whose dynamics are given by
\begin{align}
\dd X_t
	&=	\mu(X_t) \dd t + \sig(X_t) \dd W_t , &
X_0
	&\in \Rb^d . \label{eq:thd}
\end{align}
As usual, $W = (W_t)_{t \geq 0}$ is an $m$-dimensional Brownian motion.  Accordingly, we have
\begin{align}
\mu
	&: \Rb^d \to \Rb^d , &
\sig
	&: \Rb^d \to \Rb^{d \times m} .
\end{align}
Although we have specified a time-homogeneous diffusion $X$, our framework is general enough to incorporate time-inhomogeneous diffusions by considering the augmented the $\Rb^{1+d}$-dimensional process $Y = (Y_s)_{s \geq 0}$ given by $Y_t = (t,X_t)$.
\\[0.5em]
For  and $\Fb$-stopping time $\tau$, we introduce the \emph{cost functional} (again, some authors prefer the phrase \emph{gain functional} if the goal is to maximize)
\begin{align}
J(\tau)
	&=	\Eb \Big[ \int_0^\tau \ee^{-\lam t} F(X_t) \dd t + \ee^{-\lam \tau} \varphi(X_\tau) \Big] , \label{eq:J.tau}
\end{align}
where $F,\varphi : \Rb^d \to \Rb$ and $\lam \geq 0$.  Expression \eqref{eq:J.tau} is not the only possible choice for a cost functional.  However, \eqref{eq:J.tau} is general enough to encompass many optimal stopping problems of practical interest.

\begin{definition}
\label{def:tau.admissible}
We call a $\Fb$-stopping time $\tau$ \emph{admissible} if \eqref{eq:J.tau} is well-defined.
\end{definition}

\noindent We will denote by $\Tc$ the set of admissible stopping times.  The goal of \emph{optimal stopping} is to find the admissible stopping time $\tau^*$ that maximizes or minimizes $J(\tau)$.  Thus, we make the following definition:

\begin{definition}
For a given cost functional $J : \Tc \to \Rb$ (not necessarily \eqref{eq:J.tau}), we define the \emph{optimal stopping time} $\tau^*$, if it exists, by 
\begin{align}
\tau^*
	&:=	\operatorname*{argmax}_{\tau \in \Tc} J(\tau) , &
	&\text{if our goal is to maximize $J(\tau)$}, \\
\tau^*
	&:=	\operatorname*{argmin}_{\tau \in \Tc} J(\tau) , &
	&\text{if our goal is to minimize $J(\tau)$}.
\end{align}
\end{definition}

\noindent
In principle, an optimal stopping time $\tau^*$ may be arbitrarily complicated functional of the path of $X$.  However, there is a class of simple stopping times that play the same role in optimal stopping as Markov controls played in optimal controls.  

\begin{definition}
\label{def:exit.time}
We call an $\Fb$-stopping to $\tau$ a \emph{first exit time} if $\tau$ can be expressed as
\begin{align}
\tau
	&=	\inf\{ t \geq 0 : X_t \notin D \} , \label{eq:exit.time}
\end{align}
where $D \subset \Rb^d$ is an open set.
\end{definition}

\noindent
The set $D$ in \eqref{eq:exit.time} is call the \emph{continuation region} for the stopping time $\tau$.  Conveniently, optimal stopping rules of the form \eqref{eq:exit.time} are Markov. 

\section{Variational inequality}
In this section we will attempt to find optimal stopping time $\tau^*$ that maximizes the cost functional $J(\tau)$ given in \eqref{eq:J.tau} (the minimization problem is completely analogous).  We will restrict our attention to first exit time strategies.  To this end, we make the following assumption:

\begin{assumption}
\label{ass:tau}
We assume there exists an optimal stopping time $\tau^*$ which is a first exit time.  We denote by $D$ the continuation region of $\tau^*$.  We also assume that $X$ is a strong Markov process.
\end{assumption}

\noindent
A priori, there is no reason that the optimal stopping time $\tau^*$ must be an exit time.  But, we will go with it for the time being.  To begin, we define the \emph{reward-to-go} functional (or, \emph{cost-to-go}, if the goal is to minimize)
\begin{align}
J^\tau(x)
	&:=	\Eb \Big[ \int_0^\tau \ee^{-\lam t} F(X_t) \dd t + \ee^{-\lam \tau} \varphi(X_\tau) \Big| X_0 = x \Big] .
\end{align}
Observe, $J^\tau(x)$ is the cost of the stopping rule $\tau$ when $X_0 = x$.  
We define the \emph{value function} as the optimal cost-to-go
\begin{align}
V(x)
	&:=	J^{\tau^*}(x) = \max_{\tau \in \Tc} J^{\tau}(x) .
\end{align}
We will attempt to find an equation for $V(x)$.  To this end, let $\tau$ be any admissible stopping time and define
\begin{align}
\tau'
	&:=	\inf\{ t \geq \tau : X_t \notin D \} ,
\end{align}
where $D$ is the continuation region of $\tau^*$.  
Observe that $\tau' \geq \tau$ by construction and $\tau' = \tau^*$ after time $\tau$.  Thus, we have
%
\begin{align}
V(X_0)
	&=	J^{\tau^*}(X_0) \\
	&\geq J^{\tau'}(X_0) \\
	&=	\Eb \Big[ \int_0^{\tau'} \ee^{-\lam t} F(X_t) \dd t + \ee^{-\lam \tau'}\varphi(X_{\tau'}) \Big| X_0 \Big] \\
	&=	\Eb \Big[ \int_0^\tau \ee^{-\lam t} F(X_t) \dd t 
			+ \Eb \Big( \int_\tau^{\tau'} \ee^{-\lam t} F(X_t) \dd t + \ee^{-\lam \tau'} \varphi(X_{\tau'})\Big| X_\tau \Big) \Big| X_0 \Big] \\
	&=	\Eb \Big[ \int_0^\tau \ee^{-\lam t} F(X_t) \dd t 
			+ \ee^{- \lam \tau }\Eb \Big( \int_\tau^{\tau'} \ee^{-\lam (t-\tau)} F(X_t) \dd t + \ee^{ -\lam (\tau' - \tau) }\varphi(X_{\tau'})\Big| X_\tau \Big) \Big| X_0 \Big] \\
	&=	\Eb \Big[ \int_0^\tau \ee^{-\lam t} F(X_t) \dd t + \ee^{- \lam \tau } V(X_\tau ) \Big| X_0 \Big] . \label{eq:series.of.equations}
\end{align}
Here, we are using the fact that $X$ is a strong Markov process.
Now, supposing that $V$ is sufficiently smooth to apply It\^o's Lemma, we have
\begin{align}
\ee^{- \lam \tau } V(X_\tau )
	&=	V(X_0) + \int_0^\tau \ee^{-\lam t} (\Ac - \lam ) V(X_t) \dd t + \text{local martingale} , \label{eq:gen.plus.mart}
\end{align}
where $\Ac$ is the infinitesimal generator of the $X$ process.  If the local martingale above is a true martingale, then inserting \eqref{eq:gen.plus.mart} into \eqref{eq:series.of.equations} we obtain
\begin{align}
0
	&\geq \Eb \Big[ \int_0^\tau \ee^{-\lam t} \Big( F(X_t) + (\Ac - \lam ) V(X_t) \Big) \dd t \Big| X_0 \Big] , &
	&\text{(with equality if $\tau' = \tau^*$).} \label{ineq}
\end{align}
We will consider the following special case:
\\[0.5em]
\underline{Case: $\tau \leq \tau^*$}.  In this case, we have by construction that $\tau' = \tau^*$ and in inequality \eqref{ineq} becomes an equality
\begin{align}
0
	&=\Eb \Big[ \int_0^\tau \ee^{-\lam t} \Big( F(X_t) + (\Ac - \lam ) V(X_t) \Big) \dd t \Big| X_0 \Big] .
\end{align}
As the above must hold for any $\tau \leq \tau^*$ we must have
\begin{align}
0
	&=	F(X_0) + (\Ac - \lam ) V(X_0) , &
	&\text{provided that}&
\tau^*
	&> 0 .
\end{align}
If $\tau^* = 0$ then we clearly have $V(x) = \varphi(x)$.  Noting that $\tau^* > 0$ if and only if $X_0 = x \in D$, we obtain
\begin{align}
0
	&=	F(x) + (\Ac - \lam ) V(x) , &
x
	&\in D , &
V(x)
	&=	\varphi(x) , &
x
	&\notin D . \label{eq:VandD}
\end{align}
Note that, at this point, \emph{we do not know $D$}.  Without knowing $D$ we cannot solve for $V$ using \eqref{eq:VandD}.
Now, we consider the general case.
\\[0.5em]
\underline{Case: $\tau$ is arbitrary}.  In this case, it is possible that $\tau' \neq \tau^*$, and the inequality \eqref{ineq} must hold for any $\tau$.  Hence, we must have
\begin{align}
0
	&\geq	F(X_0) + (\Ac - \lam ) V(X_0) , &
	&\text{provided that}&
\tau
	&> 0 .
\end{align}
In the case of $\tau = 0$ we have $J(\tau^*) \geq J(0)$ and hence $V(x) \geq \varphi(x)$.  Putting everything together we obtain
\begin{align}
0
	&=	\max\{ F(x) + (\Ac - \lam ) V(x) ,\varphi(x) - V(x) \} . \label{ineq:variational}
\end{align}
The PDE \eqref{ineq:variational} is called a \emph{variational inequality}.  Note that if $V$ is the solution to \eqref{ineq:variational} we can reconstruct the continuation region $D$ as follows
\begin{align}
D
	&=	\{ x \in \Rb^d : V(x) > \varphi(x) \} . \label{eq:D.explicit}
\end{align}
Note that the formal derivation of the variational inequality \eqref{ineq:variational} \emph{assumes} that the optimal stopping rule $\tau^*$ is a first exit time.  What we would \emph{like} is to simply find $V(x)$ by solving \eqref{ineq:variational} and then \emph{conclude} that the optimal stopping rule $\tau^*$ is a first exit time with continuation region $D$ given by \eqref{eq:D.explicit}.  This is the subject of the following theorem.

\begin{theorem}[Verification]
\label{thm:optimal.stopping}
Fix probability space $(\Om,\Fc,\Pb)$ and a filtration $\Fb = (\Fc_t)_{t \geq 0}$ and let $X = (X_t)_{t \geq 0}$ be the $d$-dimensional time-homogeneous diffusion given by \eqref{eq:thd}.
Let $J(\tau)$ be the cost functional defined in \eqref{eq:J.tau}.  Suppose there is a function $V : \Rb^d \to \Rb$ that is sufficiently smooth to apply It\^o's Lemma and satisfies the variational inequality \eqref{ineq:variational}.  Let $D$ be given by \eqref{eq:D.explicit}.  Denote by $\Tc$ be the class of admissible stopping times (see Definition \ref{def:tau.admissible}) that satisfy
\begin{align}
\Eb \sum_{i=1}^d \sum_{j=1}^m \int_0^\tau \ee^{-\lam s} \d_{x_i} V(X_t) \cdot \sig^{(ij)}(X_t) \dd W_t^{(j)} \label{eq:tau.condition}
	&=	0 .
\end{align}
Define $\tau^* := \inf\{ t \geq 0 : X_t \notin D \}$ and suppose $\tau^* \in \Tc$.  Then $\tau^*$ is optimal: $J(\tau^*) \geq J(\tau)$ for all $\tau \in \Tc$.  Moreover, the optimal expected cost-to-go satisfies $J^{\tau^*}(x) = V(x)$.
\end{theorem}

\noindent
We will not prove Theorem \ref{thm:optimal.stopping}.  Rather, we will simply comment on the theorem's usefulness.  What Theorem \ref{thm:optimal.stopping} tells us is that, if we find a function $V$ that satisfies the variational inequality \eqref{ineq:variational} and define $\tau^* = \inf \{ t \geq 0 : X_t \notin D \}$ where $D$ is given by \eqref{eq:D.explicit}, then if $V$ is sufficiently smooth and satisfies \eqref{eq:tau.condition} with $\tau = \tau^*$, then we can conclude that $\tau^*$ is the optimal stopping rule.  

\begin{remark}[Principle of smooth fit]
\label{rmk:smooth.fit}
There will often be more than one solution to the variational inequality \eqref{ineq:variational}.  Typically, however, only one of these solutions is sufficiently smooth to apply It\^o's Lemma, namely, the solution for which both $V$ and its gradient $\nabla V$ are continuous on the boundary $\d D$.  Thus, when looking for solutions of the variational inequality, we often impose that the solution $V$ of \eqref{ineq:variational} and its gradient $\nabla V$ be continuous on the boundary $\d D$.  We call this requirement the \emph{principle of smooth fit}.
\end{remark}

\section{Example: Optimal stopping of resource extraction}
Consider a company that extracts a resource (e.g., metal, oil, gas) from the earth.  Denote by $R = (R_t)_{t \geq 0}$ the total amount of resource left to be extracted.  We suppose the dynamics of $R$ are given by
\begin{align}
\dd R_t
	&=	- \gam R_t \dd t .
\end{align}
Suppose the market price $P = (P_t)_{t \geq 0}$ of the resource follows a geometric Brownian motion
\begin{align}
\dd P_t
	&=	b P_t \dd t + \sig P_t \dd W_t ,
\end{align}
where $W = (W_t)_{t \geq 0}$ is a Brownian motion.  The company incurs a fixed cost $c$ per unit time to extract $R$ from the earth.  Thus, if the company operates from time $t=0$ to time $t = \tau$ it generates a profit of
\begin{align}
\text{Profit}
	&=	\int_0^\tau ( \gam R_t P_t - c ) \dd t .
\end{align}
Let us define $X = (X_t)_{t \geq 0}$ by $X_t = R_t P_t$.  One easily deduces that the dynamics of $X$ are given by
\begin{align}
\dd X_t
	&=	(b - \gam) X_t \dd t + a X_t \dd W_t .
\end{align}
Suppose the company wishes to choose a stopping rule $\tau$ that maximizes its expected profit.  The cost functional $J(\tau)$ is given by
\begin{align}
J(\tau)
	&=	\Eb \int_0^\tau ( \gam X_t - c ) \dd t . \label{eq:J.tau.2}
\end{align}
Comparing \eqref{eq:J.tau.2} with \eqref{eq:J.tau} we identify
\begin{align}
F(x)
	&=	\gam x - c , &
\varphi(x)
	&=	0 , &
\lam
	&=	0 .
\end{align}
The infinitesimal generator $\Ac$ of the process $X$ is given by
\begin{align}
\Ac
	&=	(b - \lam) x \d_x + \tfrac{1}{2} a^2 x^2 \d_x^2 .
\end{align}
From \eqref{ineq:variational}, we see that our candidate value function $V$ is a solution to
\begin{align}
0
	&=	\max \{ \gam x - c + (b - \lam) x V'(x) + \tfrac{1}{2} a^2 x^2 V''(x)  , - V(x) \} .
\end{align}
Thus, we must have
\begin{align}
0
	&=	\gam x - c + (b - \gam) x V'(x) + \tfrac{1}{2} a^2 x^2 V''(x) , &
0
	&> - V(x) , &
x
	&\in D , \label{eq:on.D}\\
0
	&\geq	\gam x - c + (b - \gam) x V'(x) + \tfrac{1}{2} a^2 x^2 V''(x) , &
0
	&= - V(x) , &
x
	&\in D^c . \label{eq:on.Dc}
\end{align}
Let us focus first on \eqref{eq:on.Dc}.  We have
\begin{align}
x  
	&\in D^c , &
	&\Rightarrow&
V(x)
	&=	0 , \\  & &
	&\Rightarrow&
\gam x - c + (b - \gam) x V'(x) + \tfrac{1}{2} a^2 x^2 V''(x)
	&=	\gam x - c \leq 0 , \\ & &
	&\Rightarrow&
x
	&\leq  c / \gam , \\ & &
	&\Rightarrow&
D
	&\subseteq ( c/\gam , \infty) . \label{eq:D.subset}
\end{align}
Now, let us try to solve for $V$ on $D$.  We guess a solution of the form
\begin{align}
V(x)
	&=	\frac{ c \log x }{ b - \gam - a^2/2 } - \frac{\gam x}{b - \gam} - d =: V_d(x), &
x
	&\in D
\end{align}
where $d \in \Rb$ is a constant to be determined.  One can easily verify by direct substitution that
\begin{align}
0
	&=	x - c + (b - \gam) x V_d'(x) + \tfrac{1}{2} a^2 x^2 V_d''(x) , &
\forall \, x 
	&> 0 .
\end{align}
From \eqref{eq:on.D} we must have
\begin{align}
x 
	&\in D &
	&\Rightarrow&
V_d(x)
	&> 0 , &
	&\Rightarrow&
b
	&< \gam .
\end{align}
Intuitively, it makes sense that we must have $b<\gam$.  If we had $b \geq \gam$ then the growth of $X$ would be non-negative, and it would make sense to extract resources for ever (i.e., $\tau^* = \infty$).  We will assume that $b < \gam$, which will lead to an optimal stopping rule $\tau^*$ which is finite.
\\[0.5em]
The continuation region $D$ must be of the form
\begin{align}
D
	&=	(x^*, \infty) , &
x^*
	&\geq c/\gam ,
\end{align}
where the restriction $x^* \geq c/\gam$ comes from \eqref{eq:D.subset}.
In order to find $x^*$ and  $d$ we shall use the principle of smooth fit (see Remark \ref{rmk:smooth.fit}).  Namely, we shall require that $V$ and $V'$ be continuous at the boundary of $D$.  We have
\begin{align}
&&
\lim_{\eps \to 0+} V(x^*-\eps)
	&= \lim_{\eps \to 0+}V(x^*+\eps) , &
\lim_{\eps \to 0+}V'(x^*-\eps)
	&= \lim_{\eps \to 0+}V'(x^*+\eps) , \\
&\Rightarrow&
0
	&=	V_d(x^*) , &
0
	&=	V_d'(x^*)  . \label{eq:Vd.and.Vdprime}
\end{align}  
We can use the two equations in \eqref{eq:Vd.and.Vdprime} to solve for $d$ and $x^*$.  We find
\begin{align}
x^*
	&= \frac{c}{\gam} \frac{b-\gam}{b-\gam-a^2/2}, &
d
	&=	\frac{c \log x^* }{b - \gam - a^2/2} - \frac{ \gam x^*}{b -\gam} =: d^*.
\end{align}
Thus, we have obtained a formal solution for value function $V$ and optimal stopping rule $\tau^*$
\begin{align}
V(x)
	&=	\left\{ \begin{aligned}	
			&V_{d^*}(x) &x > x^* , \\
			&0			&x \leq x^* ,
			\end{aligned} \right. &
\tau^*
	&=	\inf\{ t \geq 0 : X_t \leq x^* \} .
\end{align}
We leave it as an exercise for the reader to check that the value function $V$ satisfies the conditions of Theorem \ref{thm:optimal.stopping} (and therefore, $V$ \emph{is} the value function and $\tau^*$ \emph{is} the optimal stopping rule).

\section{Exercises}
\label{sec:exercises-12}

\begin{exercise}
\label{exer:12.1}
Fix a probability space $(\Om,\Fc,\Pb)$ and a filtration $\Fb$.
Let $X$ be a strong Markov process with dynamics of the form
\begin{align}
\dd X_t
	&=	\mu X_t \dd t + \sig X_t \dd W_t ,
\end{align}
where $W$ is a $(\Pb,\Fb)$-Brownian motion.  Define
\begin{align} 
V(x)
	&:= \sup_{\tau \in \Tc} \Eb( \ee^{- \lam \tau} (K - X_\tau)^+ | X_0 = x ) .
\end{align}
Find $V(x)$ and the optimal stopping time $\tau^*$.  You may need to make some assumption about $\mu$ to ensure a non-trivial solution.  If so, state this condition.
\end{exercise}

\begin{exercise}
\label{exer:12.2}
Fix a probability space $(\Om,\Fc,\Pb)$ and a filtration $\Fb$.
Let $X$ be a strong Markov process with dynamics of the form
\begin{align}
\dd X_t
	&=	\mu(X_t) \dd t + \sig(X_t) \dd W_t ,
\end{align}
where $W$ is a $(\Pb,\Fb)$-Brownian motion.  Define
\begin{align}
u(x)
	&:=	\Eb \Big( \int_0^\infty \ee^{- \lam t} g(X_t) \dd t \Big| X_0 = x \Big) , \\
V(x)
	&:=	\sup_{\tau \in \Tc} \Eb \Big( \int_0^\tau \ee^{- \lam t} f(X_t) \dd t + \int_\tau^\infty \ee^{- \lam t} g(X_t) \dd t \Big| X_0 = x \Big) .
\end{align}
Derive an ODE satisfied by the function $u$ and a variational inequality for $V$.  Your variational inequality should involve the function $u$.
\end{exercise}

%
%

\chapter{Monte Carlo Methods}
\label{ch:monte-carlo}
Notes from this chapter are derived from \cite{glasserman2013monte}.

\section{Overview of Monte Carlo methods}
Let $X$ be a random variable defined on a probability space $(\Om,\Fc,\Pb)$.
Suppose we wish to compute the expectation $\mu := \Eb X$ but that we do \emph{not} have an explicit expression for the distribution $F_X$ or the characteristic function $\phi_X$ of $X$.  Then we cannot compute $X$ using
\begin{align}
\mu
	&=	\int_\Rb F_X(\dd x) \, x, &
	&\text{or}&
\mu
	&=	-\ii \phi_X'(0) .
\end{align}
Suppose, however that we have the ability to generate a sequence of iid random variables $(X_i)_{i \in \Nb}$ where $X_i \sim F_X$.  Let us define
\begin{align}
\muh_n
	&:=	\frac{1}{n} \sum_{i=1}^n X_i .
\end{align}
Note that $\muh_n$ is a random variable.  By the SLLN we have
\begin{align}
\muh_n
	&\to \mu , &
n
	&\to \infty .
\end{align}
Thus, it seems reasonable to take $\muh_n$ as an estimate for $\mu$.
A natural question to ask, then, is: \emph{for a given $n$, how good is our estimate $\muh_n$ of $\mu$}?  To answer this question, let us define $\sig^2 := \Vb X$.  We compute
\begin{align}
\Vb \mu_n
	&=	\sum_{i=1}^n \Vb( X_i/n)
	=		\frac{1}{n^2} \sum_{i=1}^n \Vb X_i
	=		\frac{1}{n^2} n \Vb X 
	=		\frac{1}{n} \sig^2 .
\end{align}
So the variance of $\mu_n$ scales like $1/n$ and the standard deviation $\sqrt{\Vb \mu_n}$ scales like $1/\sqrt{n}$.
From the CLT, we know that $\muh_n$ is asymptotically normal
\begin{align}
\muh_n
	&\stackrel{\Dc}{\approx} Z_n \sim \Nc( \mu, \sig^2/n) .
\end{align}
Using this, we can construct the \emph{$c$\% confidence interval} ($c$-C.I.) for $\mu$
\begin{align}
\text{$c$-C.I.}
	&=	\Big( \muh_n - z_c \frac{\sig}{\sqrt{n}},\muh_n + z_c \frac{\sig}{\sqrt{n}}\Big) , &
\int_{-z_c}^{z_c} \dd x \, \frac{1}{\sqrt{2\pi}}\ee^{-x^2/2} 
	&=	c .
\end{align}
Thus, if we want to decrease the size of the $c$-C.I. for $\mu$ by a factor of $n$, then we need to multiply the number of samples we take of $X$ by a factor of $n^2$.
\\[0.5em]
One quick side note: without an expression for $F_X$ or $\phi_X$ we cannot know $\sig^2$.  However, we can estimate $\sig^2$ with the \emph{sample variance} 
\begin{align}
\sigh_n^2
	&:=	\frac{1}{n-1} \sum_{i=1}^n (X_i - \muh_n )^2 .
\end{align}
When $\sig^2$ is unknown and $n$ is large, we can construct the $c$-C.I. by replacing $\sig$ with $\sigh_n$.
\\[0.5em]
An alternative method of estimating $\mu$ can be described as follows.  Suppose we can generate a sequence of iid random variables $(Y_i)$ and we know that $\Eb Y_i = \mu$.  Then the following random variable
\begin{align}
\nuh_n
	&:=	\frac{1}{n} \sum_{i=1}^n Y_i ,
\end{align}
is an alternative estimator of $\mu$.  In particular, $\nuh_n$ is unbiased (as $\Eb \nuh_n = \mu$) and satisfies (by the SLLN) $\nuh_n \to \mu$.  Moreover, if $\Vb Y_i < \Vb X$, then it will follow that $\Vb \nuh_n < \Vb \muh_n$.  Choosing $(Y_i)$ cleverly so that $\Vb \nuh_n < \Vb \muh_n$ is called \emph{variance reduction}.
\\[0.5em]
The methods described above for estimating a mean $\mu$ of a random variable $X$ by generating a sequence of iid random variables -- either $(X_i)$ or $(Y_i)$ -- is the \emph{Monte Carlo} method in a nutshell.  More generally, we may be interested in computing, not just the expectation of a random variable, but the expectation $\Eb F[X]$ of a functional $F$ of an entire stochastic process $X = (X_t)_{t \geq 0}$.  For example, we may wish to compute the time-weighted average of $X$, given by
\begin{align}
\Eb F[X]
	&:=	\Eb \Big( \frac{1}{T} \int_0^T X_t \dd t \Big) .
\end{align}
To do this, we will need to know how to generate iid random \emph{process} $X^{(i)} = (X_t^{(i)})_{t \geq 0}$, $i=1,2,\ldots$.  Note that a collection of processes $(X^{(i)})$ is iid if the processes are independent and have the same finite-dimensional distributions
\begin{align}
X^{(i)}
	&\ind X^{(j)} , &
i
	&\neq j , \\
(X_{t_1}^{(i)}, X_{t_2}^{(i)}, \ldots, X_{t_n}^{(i)})
	&\stackrel{\Dc}{=} (X_{t_1}^{(j)}, X_{t_2}^{(j)}, \ldots, X_{t_n}^{(j)}) , &
\forall \, 
	&0 \leq t_1 \leq t_2 \leq \ldots \leq t_n < \infty .
\end{align}
With the above in mind, the focus of this chapter will be an answering the follow two questions:
\begin{enumerate}
\item How can we generate iid random variables?  Or, more generally, how can we generate iid random processes?
\item  How can we ``intelligently sample'' (i.e., construct estimators of $\mu$) in order to reduce the variance of our estimates?
\end{enumerate}
The first question will occupy Sections \ref{sec:Gen.RVs} through \ref{sec:Gen.JD}.  Variance reduction methods will be the focus on Section \ref{sec:VarReduction}.

\section{Generating Random variables}
\label{sec:Gen.RVs}
Although programs such as Mathematica and Matlab can generate random variables with a variety of distributions, it is worth knowing how this is accomplished.
\begin{assumption}
Throughout Section \ref{sec:Gen.RVs} both $U$ and the sequence $(U_i)$ will represent iid random variables that are uniformly distributed $U_i \sim \Uc(0,1)$.
\end{assumption}
\noindent
Nearly \emph{every} computer can generate uniformly distributed random variables.

\subsection{Inverse Transform Method}
\label{sec:itm}
When the distribution function $F_X$ of a random variable $X$ is known in closed-form, we can use the \emph{Inverse Transform Method} to generate iid random variables with the same distribution as $X$.  The following theorem describes how this is accomplished.
\begin{theorem}[Inverse Transform Method]
\label{thm:itm}
Let $F_X$ be the distribution function of a random variable $X$ defined on $(\Om,\Fc,\Pb)$.  
Define
\begin{align}
F_X^{-1}(u)
	&:=	\inf\{ x : F(x) \geq u \} . \label{eq:F.inverse}
\end{align}
Then we have
\begin{align}
F_X^{-1}(U)
	&\sim F_X . \label{eq:itm}
\end{align}
\end{theorem}

\begin{proof}
For simplicity, let us Suppose $F_X$ is strictly increasing (this is always the case for continuous random variables).  
Noting that $F_X(x) \in (0,1)$ for all $x \in \Rb$ and $\Pb( U \leq y) = y$ for all $y \in (0,1)$ we have
\begin{align}
\Pb( F_X^{-1}(U) \leq x )
	&=	\Pb( U \leq F_X(x) )
	=		F_X(x) ,
\end{align}
which establishes \eqref{eq:itm}.
\end{proof}

\begin{example}
Let us apply Theorem \ref{thm:itm} to generate an exponential random variable $X \sim \Ec(\lam)$.  For $x \in (0,\infty)$ we have
\begin{align}
F_X(x)
	&=	\int_0^x \dd x \, \lam \ee^{- \lam x} 
	=		1 - \ee^{-\lam x} , &
	&\Rightarrow&
F_X^{-1}(u)
	&=	\frac{-1}{\lam} \log (1 - u) .
\end{align}
Using the fact that $U \sim 1-U$, we find that
\begin{align}
\frac{-1}{\lam} \log U 
	&\sim \Ec(\lam).
\end{align}
\end{example}

\begin{example}
Suppose $X$ is a discrete random variable with $\Pb( X = x_i) = p_i$ for $i = 1, 2, \ldots, n$.  Then
\begin{align}
F_X(x)
	&=	\sum_{i:x_i \leq x} p_i , &
	&\Rightarrow&
F_X^{-1}(u)
	&=	\sum_{i=1}^n x_i \Ib_{\{P_{i-1} < u \leq P_{i}\}} , &
P_i
	&=	\sum_{j=1}^i p_j . \label{eq:F.inverse.discrete}
\end{align}
where $P_0 := 0$.  To understand \eqref{eq:F.inverse.discrete}, it may help to draw a graph of $F_X$ and $F_X^{-1}$.
\end{example}

\begin{theorem}
Let $F_X$ be the distribution of a random variable $X$ defined on $(\Om,\Fc,\Pb)$.  Let $Y:= X | a< X \leq b$ be a random variable whose distribution $F_Y$ is the conditional distribution of $X$ given $a < X \leq b$.  Define
\begin{align}
V
	&:=	F_X(a) + \Big( F_X(b) - F_X(a) \Big) U .
\end{align}
Then, with $F_X^{-1}$ as defined in \eqref{eq:F.inverse}, we have
\begin{align}
F_X^{-1}(V) 
	&\sim F_Y = F_{X|a<X\leq b}. \label{eq:FX.V}
\end{align}
\end{theorem}

\begin{proof}
Assume $x \in (a,b)$.  Then we have
\begin{align}
\Pb( F_X^{-1}(V) \leq x )
	&=	\Pb( V \leq F_X(x) )
	=		\Pb( F_X(a) + \( F_X(b) - F_X(a) \) U \leq F_X(x) ) \\
	&=	\Pb\( U \leq \frac{F_X(x) - F_X(a)}{F_X(b) - F_X(a)} \) \\
	&=	\frac{F_X(x) - F_X(a)}{F_X(b) - F_X(a)}
	=		\frac{ \Pb( X \leq x , a < X \leq b) }{\Pb(a < X \leq b)} \\
	&=	\Pb( X \leq x | a < X \leq b ) 
	=		F_{X|a<X\leq b}(x) ,
\end{align}
which proves \eqref{eq:FX.V}.
\end{proof}

\subsection{Acceptance-Rejection method}
The inverse transform method described in Section \ref{sec:itm} works only when the function $F_X^{-1}$ is known in closed form.  Unfortunately, there are many random variables for which $F_X^{-1}$ is not known in closed-form.  An alternative method of generating random variables with a particular distribution $F_X$ is the \emph{Acceptance-Rejection Method}, which we describe in the following theorem.

\begin{theorem}[Acceptance-Rejection Method]
\label{thm:arm}
Suppose $X \in D \subset \Rb^d$ and $Y \in D$ are random vectors defined on $(\Om,\Fc,\Pb)$ with densities $f_X: D \to \Rb_+$ and $f_Y: D \to \Rb_+$, respectively.  Suppose further that there exists a constant $c>0$ such that $f_Y(x) \leq c f_X(x)$ for all $x \in D$.  
Let $U \sim \Uc(0,1)$ be independent of $X$.  Then
\begin{align}
Y
	&\stackrel{\Dc}{=} X | U \leq \frac{f_Y(X)}{c f_X(X)} . \label{eq:arm}
\end{align}
\end{theorem}

\begin{proof}
Suppose $A \subset \Bc(D)$.  Then we have
\begin{align}
\Pb \( X \in A | U \leq \tfrac{1}{c} f_Y(X)/f_X(X) \) 
	&=	\frac{ \Pb \( X \in A , U \leq \tfrac{1}{c} f_Y(X)/f_X(X) \) }{ \Pb \( U \leq \tfrac{1}{c} f_Y(X)/f_X(X) \) } . \label{eq:arm.1}
\end{align}
For the numerator in \eqref{eq:arm.1}, we have
\begin{align}
\Pb \( X \in A , U \leq \tfrac{1}{c} f_Y(X)/f_X(X) \)
	&=	\Eb \Pb \( X \in A , U \leq \tfrac{1}{c} f_Y(X)/f_X(X) | X \) \\
	&=	\Eb \Ib_{\{X \in A \}} \tfrac{1}{c} f_Y(X)/f_X(X) \\
	&=	\int_A \dd x \, f_X(x) \tfrac{1}{c} f_Y(x)/f_X(x) \\
	&=		\int_A \dd x \, \tfrac{1}{c} f_Y(x) 
	=		\tfrac{1}{c} \Pb( Y \in A ) . \label{eq:arm.3}
\end{align}
The denominator of \eqref{eq:arm.1} is obtained from \eqref{eq:arm.3} by considering $A = D$, which yields
\begin{align}
\Pb \( U \leq \tfrac{1}{c} f_Y(X)/f_X(X) \)
	=		\tfrac{1}{c} \Pb( Y \in D ) = \tfrac{1}{c} . \label{eq:arm.2}
\end{align}
Inserting \eqref{eq:arm.3} and \eqref{eq:arm.2} into \eqref{eq:arm.2} yields
\begin{align}
\Pb \( X \in A | U \leq \tfrac{1}{c} f_Y(X)/f_X(X) \) 
	&=	\Pb( Y \in A ) ,
\end{align}
which proves \eqref{eq:arm}.
\end{proof}
\noindent
Using Theorem \ref{thm:arm}, if we can generate a random variable $X$ with density $f_X$, then we can generate a random variable $Y$ with density $f_Y$ using the following algorithm
\begin{enumerate*}
\item Generate $X$. \label{itm:gen.X}
\item Generate $U$.
\item If $U \leq \tfrac{1}{c} f_Y(X)/f_X(X)$, return $X$.  Otherwise, go back to Step \ref{itm:gen.X}.
\end{enumerate*}
In order for the acceptance rejection method to work efficiently, we would like $\Pb( U \leq \tfrac{1}{c} f_Y(X)/f_X(X) )$ to be as close to one as possible.  Thus, we should choose $f_X$ to be similar to $f_Y$ and we should choose $c$ to be as small as possible.

\subsection{Generating Gaussian RVs}
In this section, we will show how to generate independent and correlated Gaussian RVs.

\subsubsection{Box-Muller Method}
The \emph{Box-Muller Method} of generating two independent Gaussian random variables is summarized in the following theorem.
\begin{theorem}[Box-Muller Method]
\label{thm:bmm}
For $i=1,2$, define
\begin{align}
X
	&=	\sqrt{- 2 \log (U_1)} \cos (2 \pi U_2), &
Y
	&=	\sqrt{- 2 \log (U_1)} \sin (2 \pi U_2)
\end{align}
Then $X \ind Y$ and $X,Y \sim \Nc(0,1)$.
\end{theorem}
\begin{proof}
Define $R := \sqrt{ - 2 \log (U_1) }$.  Then
\begin{align}
\Pb( R \leq r )
	&=	\Pb( - 2 \log(U) \leq r^2)
	=		\Pb( \log(U) \geq - r^2/2 ) \\
	&=		1 - \Pb( U < \ee^{-r^2/2} )
	=		1 - \ee^{-r^2/2} .
\end{align}
Denoting by $f_R$ the density of $R$ we have
\begin{align}
f_R(r)
	&=	\frac{\dd}{\dd r} F_R(r) = r \ee^{-r^2} , &
r
	&> 0 .
\end{align}
Next, define $\Phi := 2 \pi U_2$.  It is obvious that $\Phi \sim \Uc(0,2\pi)$ and thus
\begin{align}
f_\Phi(\phi)
	&=	\frac{1}{2\pi} , &
\phi
	&\in (0,2\pi) .
\end{align}
As $U_1 \ind U_2$ by assumption, it follows that $R \ind \Phi$ and thus
\begin{align}
f_{R,\Phi}(r,\phi)
	&=	f_R(r) f_\Phi(\phi) = \frac{r}{2\pi}\ee^{-r^2/2} .
\end{align}
Next, define $g : \Rb_+ \times (0,2\pi) \to (0,1) \times (0,1)$ by
\begin{align}
g(r, \phi)
`	&:=	( r \cos \phi, r \sin \phi ) , &
	&\Rightarrow&
g^{-1}(x,y)
	&=	( \sqrt{ x^2 + y^2 } , \tan^{-1}(y/x) ) .
\end{align}
Noting that $(X,Y) =	g(R,\Phi)$, we have
\begin{align}
f_{X,Y}(x,y)
	&=	f_{R,\Phi}( g^{-1}(x,y) ) \frac{1}{ | g'(g^{-1}(x,y) ) | } &
g'(r,\phi)
	&=	\( \begin{array}{cc}
			\d g_1 / \d r & \d g_1 / \d \phi \\
			\d g_2 / \d r & \d g_2 / \d \phi
			\end{array} \) , \\
	&=	\frac{1}{\sqrt{2\pi}} \ee^{-x^2/2} \frac{1}{\sqrt{2\pi}} \ee^{-y^2/2}
	=		f_Z(x) \cdot f_Z(y) , &
f_Z(z)
	&:=	\frac{1}{\sqrt{2\pi}} \ee^{-z^2/2} ,
\end{align}
where we have omitted a some nontrivial by straightforward algebraic manipulations.
\end{proof}
\noindent
Note that we generate a a normal random variable $X \sim \Nc (\mu, \sig^2)$ from a standard normal as follows
\begin{align}
X
	&:=	\mu + \sig Z \sim \Nc(\mu,\sig), &
Z
	&\sim \Nc(0,1) .
\end{align}

\subsubsection{Generating correlated Gaussian RVs}
We now know how to generate \emph{independent} Gaussian random variables.  In this section, we will show how to generate \emph{correlated} Gaussian random variables.  The density $f_X$ of \emph{multivariate Gaussian vector} $X \in \Rb^d$ is given by
\begin{align}
f_X(x)
	&=	\frac{1}{\sqrt{(2\pi)^d |\Sigma|}}\exp\( - \tfrac{1}{2} (x - \mu)^\text{T} \Sigma^{-1} (x - \mu) \) ,
\end{align}
where we have defined the \emph{mean vector} $\mu$ and the \emph{covariance matrix} $\Sigma$
\begin{align}
\mu
	&=	\( \begin{array}{c}
			\Eb X_1 \\ \Eb X_2 \\ \vdots, \\ \Eb X_d 
			\end{array} \) , &
\Sigma
	&=	\( \begin{array}{cccc}
			\Eb (X_1-\mu_1)(X_1-\mu_1) 	& \Eb (X_1-\mu_1)(X_2-\mu_2) & \ldots & \Eb (X_1-\mu_1)(X_d-\mu_d) \\
			\Eb (X_2-\mu_2)(X_1-\mu_1) 	& \Eb (X_2-\mu_2)(X_2-\mu_2) & \ldots & \Eb (X_2-\mu_2)(X_d-\mu_d) \\
			\vdots											&	\vdots										& \ddots & \vdots \\
			\Eb (X_d-\mu_d)(X_1-\mu_1) 	& \Eb (X_d-\mu_d)(X_2-\mu_2) & \ldots & \Eb (X_d-\mu_d)(X_d-\mu_d)
			\end{array} \) . \\[-2em]
\end{align}
If $I$ is the $d$-dimensional identity matrix then, by definition, $Z \sim \Nc( 0, I)$ is a $d$-dimensional Gaussian with \emph{independent} components.
Suppose now that $A \in \Rb^{d \times d}$.  One can show by direct computation that
\begin{align}
X
	&:=	\mu + A Z \sim \Nc( \mu, A A^\text{T} ) , &
Z
	&\sim \Nc( 0, I ) .
\end{align}
Thus, if we wish to generate a $d$-dimensional Gaussian random vector $X \sim (\mu,\Sig)$ we need only to find a matrix $A$ such that $A A^\text{T} = \Sigma$.  Such a matrix $A$ can be constructed via the \emph{Cholesky Factorization}, which we describe in the following theorem.

\begin{theorem}[Cholesky Decomposition]
\label{thm:cholesky}
Suppose $\Sig$ is a $d \times d$ symmetric positive definite matrix.  Define a $d \times d$ lower triangular (i.e. $A_{ij}=0$ if $j>i$) matrix $A$ by
\begin{align}
A_{ii}
	&=	\( \Sigma_{ii} - \sum_{k=1}^{i-1} A_{ik}^2 \)^{1/2} , &
i
	&=	1, 2, \ldots, d , \\
A_{ij}
	&=	\frac{1}{A_{jj}} \( \Sigma_{ij} - \sum_{k=1}^{j-1} A_{ik} A_{jk} \) , &
j
	&=	1, 2, \ldots, i-1 .
\end{align}
Then $A A^\text{T} = \Sigma$.
\end{theorem}
\noindent
We will not provide a proof of Theorem \ref{thm:cholesky}.  Instead, we direct the interested reader to \cite[Chapter 2]{glasserman2013monte}.

\section{Simulating Discrete-Time Markov Chains}
Let $X = (X_i)_{i \in \Nb_0}$ be a time-homogeneous discrete-time Markov chain with state space $S = \{ 1, 2, \ldots, |S|\}$ (note that $|S|$ may be infinite).  Define $p: S \times S \to [0,1]$ by
\begin{align}
p(i,j)
	&=	\Pb( X_{n+1} = j | X_n = i ) .
\end{align}
We can simulate a path of $X$ by setting
\begin{align}
X_{n+1}
	&=	Y_{n+1} , &
n
	&=	0, 1, 2, \ldots ,
\end{align}
where $Y_{n+1}$ is a discrete random variable with
\begin{align}
\Pb( Y_{n+1} = j | X_n )
	&=	p(X_n,j) , &
j
	&=	1, 2, \ldots, |S| .
\end{align}
Note that $Y_{n+1}$ can be generated using the Inverse Transform Method, described in Section \ref{sec:itm}. 

\section{Simulating Continuous-Time Markov Chains}
Let $X = (X_t)_{t \geq 0}$ be a time-homogeneous continuous-time Markov chain with state space $S = \{ 1, 2, \ldots, |S|\}$ (note that $|S|$ may be infinite).  Denote by $\Gv = (g(i,j))_{i,j}$ the generator of $X$.  Let us denote by $\tau_n$ the time of the $n$th jump of $X$ with $\tau_0 := 0$.  Recall that
\begin{align}
\tau_{n+1} - \tau_n | X_{\tau_n}
	&\sim \Ec( -g(X_{\tau_n},X_{\tau_n}) ) , &
\Pb( X_{\tau_{n+1}} = j | X_{\tau_n} )
	&=	\frac{-g(X_{\tau_n},j)}{g(X_{\tau_n},X_{\tau_n})} .
\end{align}
With this in mind, we can simulate a path of $X$ by setting
\begin{align}
\tau_{n+1}
	&=	\tau_n + T_{n+1} &
X_{\tau_{n+1}}
	&=	Y_{n+1} , &
n
	&=	0, 1, 2, \ldots ,
\end{align}
where $T_{n+1}$ is an exponentially distributed random variable and $Y_{n+1}$ is a discrete random variable with
\begin{align}
T_{n+1} | X_{\tau_n}
	&\sim \Ec(-g(X_{\tau_n},X_{\tau_n})) , &
\Pb( Y_{n+1} = j | X_{\tau_n} )
	&=	\frac{-g(X_{\tau_n},j)}{g(X_{\tau_n},X_{\tau_n})} , &
j
	&=	1, 2, \ldots, |S| .
\end{align}
Note that both $T_{n+1}$ and $Y_{n+1}$ can be generated using the Inverse Transform Method, described in Section \ref{sec:itm}. Finally, as $X$ is constant between jumps, we have
\begin{align}
X_{\tau_j + s}
	&=	X_{\tau_j} , &
s
	&=	[0, \tau_{j+1}) .
\end{align}

\section{Simulating Diffusions}
\label{sec:Gen.Diffs}
In this section, we shall assume that $X = (X_t)_{t \geq 0}$ is the solution of an SDE of the form
\begin{align}
\dd X_t
	&=	\mu(t,X_t) \dd t + \sig(t, X_t) \dd W_t , \label{eq:sde.mc}
\end{align}
or, in integral form
\begin{align}
\dd X_{t_{i+1}}
	&=	X_{t_i} + \int_{t_i}^{t_{i+1}} \mu(t,X_t) \dd t + \int_{t_i}^{t_{i+1}} \sig(t,X_t) \dd W_t .
\end{align}
In principle, $X$ and $W$ could live in $\Rb^d$ and $\Rb^m$, respectively.  However, for simplicity, we shall assume that $d = m = 1$.
\begin{assumption}
Throughout Sections \ref{sec:Gen.Diffs} and \ref{sec:Gen.JD} the random variable $Z$ and the sequence $(Z_i)$ will represent iid random variables that are Normally distributed $Z, Z_i \sim \Nc(0,1)$.
\end{assumption}

\subsection{Non-random coefficients}
Let us consider the case when the drift and diffusion coefficients $\mu$ and $\sig$ of the diffusion $X$ are deterministic functions of time.  In this case we have
\begin{align}
X_{{t_{i+1}}}
	&=	X_{t_i} + \int_{t_i}^{t_{i+1}} \mu(t) \dd t + \int_{t_i}^{t_{i+1}} \sig(t) \dd W_t . \label{eq:dX.time.dep.coeff}
\end{align}
From Proposition \ref{prop:normal.ito} we have that the It\^o integral is independent of $X_{t_i}$ and normally distributed
\begin{align}
\int_{t_i}^{t_{i+1}} \sig(t) \dd W_t
	&\sim \Nc(0,v^2(t_i,t_{i+1})) , &
v^2(t_i,t_{i+1})
	&:=	\int_{t_i}^{t_{i+1}} \sig^2(t) \dd t .
\end{align}
Thus, if we wish to simulate the value of $X$ on a grid $\Pi = \{ t_0 = 0, t_1, t_2, \ldots \}$, we can do so using the following algorithm
\begin{align}
X_{{t_{i+1}}}
	&=	X_{t_i} + \int_{t_i}^{t_{i+1}} \mu(t) \dd t + v(t_i,t_{i+1}) Z_{i+1} , \label{eq:dX.exact}
\end{align}
where we have used $v Z_{i+1} \sim \Nc(0,v^2)$.  Notice that the algorithm \eqref{eq:dX.exact} used to generate the value of $X$ on the grid $\Pi$ is \emph{exact}, meaning the distributions of $(X_{t_0},X_{t_1},X_{t_2},\ldots)$ as generated using \eqref{eq:dX.time.dep.coeff} and \eqref{eq:dX.exact} are the same.

\begin{example}
Suppose $S = (S_t)_{t \geq 0}$ is a geometric Brownian motion with time-dependent coefficients
\begin{align}
\dd S_t
	&=	b(t) S_t \dd t + a(t) S_t \dd W_t .
\end{align}
Note that the drift and diffusion coefficients of $S$ are random.  Nevertheless, we can simulate $S$ exactly as follows.  First define $X = \log S$.  Then, by It\^o's Lemma, we have
\begin{align}
\dd X_t
	&=	\( b(t) - \tfrac{1}{2} a^2(t) \) \dd t + a(t) \dd W_t .
\end{align}
As the drift and diffusion coefficients of $X$ are deterministic functions of time, we can simulate $(X_{t_0},X_{t_1},X_{t_2},\ldots)$ exactly.  Then we can set
$(S_{t_0},S_{t_1},S_{t_2},\ldots) =	(\ee^{X_{t_0}},\ee^{X_{t_1}},\ee^{X_{t_2}},\ldots)$.
\end{example}

\subsection{Ornstein-Uhlenbeck and other Gaussian Processes}
Suppose that $X=(X_t)_{t \geq 0}$ is an OU process
\begin{align}
\dd X_t
	&=	\kappa ( \theta - X_t) \dd t + a \dd W_t .
\end{align}
As we showed in Example \ref{ex:OU}, we can solve the SDE for $X$ explicitly
\begin{align}
X_{t_{i+1}}
	&=	\theta + \ee^{- \kappa (t_{i+1} - t_i)} ( X_{t_i} - \theta ) + \int_{t_i}^{t_{i+1}} \ee^{-\kappa(t_{i+1}-t)} a \dd W_t .
\end{align}
From Proposition \ref{prop:normal.ito} we have that the It\^o integral is independent of $X_{t_i}$ and normally distributed
\begin{align}
\int_{t_i}^{t_{i+1}} \ee^{-\kappa(t_{i+1}-t)} a \dd W_t
	&\sim \Nc(0,v^2(t_i,t_{i+1})) , &
v^2(t_i,t_{i+1})
	&=	\int_{t_i}^{t_{i+1}} \ee^{- 2 \kappa(t_{i+1}-t)} a^2 \dd t \\ & &
	&=	\frac{a^2}{2\kappa} \Big( 1 - \ee^{-2 \kappa (t_{i+1} - t_i)} \Big) . \label{eq:v.ou}
\end{align}
Thus, if we wish to simulate the value of $X$ on a grid $\Pi = \{ t_0 = 0, t_1, t_2, \ldots \}$, we can do so using the following algorithm
\begin{align}
X_{{t_{i+1}}}
	&=	\theta + \ee^{- \kappa (t_{i+1} - t_i)} ( X_{t_i} - \theta ) + v(t_i,t_{i+1}) Z_{i+1} ,
\end{align}
where $v$ is given by \eqref{eq:v.ou}.  Note that this simulation is exact.  More generally, one can simulate exactly, any process of the form
\begin{align}
\dd X_t
	&=	\Big( f(t) + g(t) X_t \Big) \dd t + h(t) \dd W_t ,
\end{align}
where $(f,g,h)$ are deterministic functions of time.  We leave it as an exercise for students to show that $X_{t_{i+1}} | X_{t_i}$ is normally distributed. 

\subsection{Euler Discretization of an SDE}
\label{sec:euler}
Unfortunately, most SDEs cannot be simulated exactly.  In such cases, we must resort to approximation schemes.  The work-horse of approximations schemes for SDEs is the \emph{Euler Discretization Scheme}, which we now describe.  To begin, we fix a time horizon $T>0$ and a time grid $\Pi$ as follows
\begin{align}
\Pi
	&=	\{ t_0, t_1, \ldots, t_n \} , &
t_i
	&=	i \del , &
\del
	&=	T/n . \label{eq:time.grid}
\end{align}
For $\del \ll 1$ is seems reasonable to approximate the process $X$ in \eqref{eq:sde.mc} on the grid $\Pi$ by a process $\Xh$ given by
\begin{align}
\Xh_{t_{i+1}}
	&=	\Xh_{t_i} + \int_{t_i}^{t_{i+1}} \mu(t_i,\Xh_{t_i}) \dd t + \int_{t_i}^{t_{i+1}} \sig(t_i,\Xh_{t_i}) \dd W_t \\
	&=	\Xh_{t_i} + \mu(t_i,\Xh_{t_i}) \( t_{i+1} - t_i \)  + \sig(t_i,\Xh_{t_i}) \( W_{t_{i+1}} - W_{t_i} \) . \label{eq:dXhat}
\end{align}
From the definition of the It\^o integral, as $n \to \infty$ we know that the process $\Xh$ converges (in $L^2(\Om,\Fc,\Pb)$ and therefore, in probability) to $X$.  Moreover, for a fixed $n$, we can simulate the process $\Xh$ exactly using
\begin{align}
\Xh_{t_{i+1}}
	&=	\Xh_{t_i} + \mu(t_i,\Xh_{t_i}) \del  + \sig(t_i,\Xh_{t_i}) \sqrt{\del} Z_{i+1} ,
\end{align}
where we have used $t_{i+1}-t_i=\del$ and  $\( W_{t_{i+1}} - W_{t_i} \) \sim \Nc(0,\del)$.  The question naturally arises then: \emph{how fast does $\Xh$ converge to $X$ as $\del \to 0$}?  The give a quantifiable response to this question, we need to define some reasonable measures of convergence.
\begin{definition}
Let $X$ and $\Xh$ be as given in \eqref{eq:sde.mc} and \eqref{eq:dXhat}, respectively.
We say that \emph{$\Xh$ converges strongly to $X$ with order $\beta>0$} if there exists a constant $C$ such that
\begin{align}
\Eb | \Xh_{n\del} - X_T |
	&\leq	C \del^\beta , &
	&\text{as $\del \to 0$}. \label{eq:strong.order}
\end{align}
We say that \emph{$\Xh$ converges weakly to $X$ with order $\beta>0$} if there exists a constant $C$ such that
\begin{align}
| \Eb f(\Xh_{n\del}) - \Eb f(X_T) |
	&\leq	C \del^\beta , &
	&\text{as $\del \to 0$} &
\forall \, f
	&\in	C_P^{2\beta+1} , \label{eq:weak.order}
\end{align}
where $C_P^k$ is the set of functions whose derivatives up to order $k$ are polynomially bounded.
\end{definition}
\noindent
Although we will not prove it, the Euler discretization $\Xh$ converges strongly to $X$ with order $\beta = 1/2$ if the drift and diffusion coefficients of $X$ satisfy the linear growth and Lipchitz conditions required for existence and uniqueness of a strong solution of an SDE (see Theorem \ref{thm:exist.unique}) and additionally satisfy
\begin{align}
|\mu(t,x) - \mu(s,x) | + |\sig(t,x) - \sig(s,x) |
	&\leq C (1 + |x|)\sqrt{t-s} , &
\forall \, 0
	&\leq s < t < \infty ,
\end{align}
for some $C$.  Under the additional conditions that $\mu, \sig \in C_P^{2\beta+1}$, one can show that $\Xh$ converges weakly to $X$ with order $\beta = 1$.  For detailed statements concerning weak and strong convergence of $\Xh$ to $X$ we refer the reader to \cite{kloeden-platen}.
\\[0.5em]
Statements concerning order of convergence give us information about how quickly the approximation $\Xh$ approaches $X$.  However, as the constants $C$ in \eqref{eq:strong.order} and \eqref{eq:weak.order} are typically unknown, an order of convergence statement does not provide much guidance for choosing the size of $\del$ to obtain a desired level of accuracy.  From a practical standpoint, one typically must resort to trial-and-error in order to find an appropriate $\del$.  For example, if one wishes to estimate $\Eb f(X_T)$ one could simulate $m$ sample paths $(\Xh^{(1)}, \Xh^{(2)},\ldots,\Xh^{(m)})$ with $\del$ fixed and then compute
\begin{align}
\frac{1}{m} \sum_{i=1}^m f(\Xh_{n\del}^{(i)}) . \label{eq:example.sum}
\end{align}
One could then repeat this procedure with smaller and smaller $\del$'s until the value of \eqref{eq:example.sum} does not change ``too much,'' where ``too much'' depends on the desired level of accuracy.

\section{Simulating Jump-Diffusions}
\label{sec:Gen.JD}
In this Section, we suppose that $X = (X_t)_{t \geq 0}$ is the solution of the following SDE
\begin{align}
\dd X_t
	&=	\mu(t,X_t) \dd t + \sig(t,X_t) \dd W_t + \int_\Rb \gam(t,X_{t-},z) N(\dd t, \dd z) , \label{eq:dX.jumps}
\end{align}
where $N$ is a finite-activity Poisson random measure
\begin{align}
\Eb N(\dd t, \dd z)
	&=	\nu(\dd z) \dd t = \lam F(\dd z) \dd t, &
\nu(\Rb)
	&=	\lam  .
\end{align}
We could have alternatively written the dynamics of $X$ using a compensated Poisson random measure $\Nt(\dd t, \dd z) = N(\dd t, \dd z) - \nu(\dd z) \dd t$ as follows
\begin{align}
\dd X_t
	&=	\( \mu(t,X_t) + \int_\Rb \gam(t,X_{t-},z) \nu(\dd z) \) \dd t + \sig(t,X_t) \dd W_t + \int_\Rb \gam(t,X_{t-},z) \Nt(\dd t, \dd z) .
\end{align}
However, for the purposes of simulating sample paths of $X$, it is preferable to work with the expression \eqref{eq:dX.jumps}.  We will present two Euler-like discretizations of \eqref{eq:dX.jumps}

\subsection{Euler Discretization on a fixed time grid}
As in Section \ref{sec:euler}, let us begin by fixing a fixing a time horizon $T>0$ and a time grad
\begin{align}
\Pi
	&=	\{ t_0, t_1, \ldots, t_n \} , &
t_i
	&=	i \del , &
\del
	&=	T/n . 
\end{align}
We can approximate $X$ on $\Pi$ by a process $\Xh$ given by
\begin{align}
\Xh_{t_{i+1}}
	&=	X_{t_i} + \int_{t_i}^{t_{i+1}} \mu(t_i,\Xh_{t_i}) \dd t + \int_{t_i}^{t_{i+1}} \sig(t_i,\Xh_{t_i}) \dd W_t 
			+ \int_{t_i}^{t_{i+1}} \int_\Rb \gam(t_i,\Xh_{t_i},z) N(\dd t, \dd z) \\
	&=	X_{t_i} + \mu(t_i,\Xh_{t_i}) (t_{i+1}-t_i) + \sig(t_i,\Xh_{t_i}) (W_{t_{i+1}} - W_{t_i} ) \\ &\quad
			+ \int_\Rb \gam(t_i,\Xh_{t_i},z) \big( N(t_{i+1}, \dd z) - N(t_i,\dd z) \big) . \label{eq:dXhat.2}
\end{align}
We know that $t_{i+1} - t_i = \del$ and $W_{t_{i+1}} - W_{t_i} \sim \Nc(0,\del)$. Now, for the jump-term we note that
\begin{align}
\Pb \Big( N(t_{i+1}, \Rb) - N(t_i, \Rb) = 0 \Big)
	&=	1 - \lam \del + \Oc(\del^2) , \\
\Pb \Big( N(t_{i+1}, \Rb) - N(t_i, \Rb) = 1 \Big)
	&=	\lam \del + \Oc(\del^2) , \\
\Pb \Big( N(t_{i+1}, \Rb) - N(t_i, \Rb) \geq 2 \Big)
	&=	\Oc(\del^2) .
\end{align}
Thus, we can approximate the jump-term in \eqref{eq:dXhat.2} using a Bernoulli random variable $B_{i+1}$ and a random variable $Y_{i+1}$ with distribution $F$ as follows
\begin{align}
\int_\Rb \gam(t_i,\Xh_{t_i},z) \big( N(t_{i+1}, \dd z) - N(t_i,\dd z) \big)
	&\stackrel{\Dc}{\approx} \gam(t_i,\Xh_{t_i},Y_{i+1}) B_{i+1} , &
Y_{i+1}
	&\sim F , &
B_{i+1}
	&\sim \Bc( \lam \del ) .
\end{align}
Putting everything together, we can approximately simulate $\Xh$ using the following algorithm
\begin{align}
\Xh_{t_{i+1}}
	&=	\Xh_{t_i} + \mu(t_i,\Xh_{t_i}) \del + \sig(t_i,\Xh_{t_i}) \sqrt{\del} Z_{i+1} + \gam(t_i,\Xh_{t_i},Y_{i+1}) B_{i+1} + \Oc(\del^2) 
\end{align}
Note that our simulation of $\Xh$ is not exact due to the error term of $\Oc(\del^2)$.  However, as we are already making an error of $\Oc(\del)$ be approximating $X$ with $\Xh$, the $\Oc(\del^2)$ term is not problematic.  

\subsection{Euler Discretization on a random time grid}
An alternative Euler-like Discretization of $X$ can be constructed as follows.  First, we fix a time-horizon $T < \infty$.  Noting that $N(t,\Rb)$ is a Poisson process with intensity $\lam = \nu(\Rb)$ we can simulate the jump times $\tau_i$ of $X$ as follows
\begin{align}
\tau_{i+1}
	&=	\tau_i + T_{i+1} , &
\tau_0
	&:=	0 , &
T_{i+1}
	&\sim \Ec(\lam) .
\end{align}
Next, we construct a time grid
\begin{align}
\Pi
	&=	\{ t_0 = 0, t_1, t_2, \ldots, t_i, \tau_1-,\tau_1, t_{i+1}, \ldots, t_{j}, \tau_2-,\tau_2, t_{j+1}, \ldots, t_n = T \} , &
t_k
	&=	k \del , &
\del
	&=	T / n ,
\end{align}
where the jump times $(\tau_i)$ and the fixed times $(t_i)$ are placed in chronological order.
Now, from time $\tau_i$ to $\tau_{i+1}-$ the dynamics of $X$ are given by
\begin{align}
\dd X_t
	&=	\mu(t,X_t) \dd t + \sig(t,X_t) \dd W_t , &
t
	&\in [\tau_i,\tau_{i+1}) .
\end{align}
Thus, in between jump dates, we can simulate $X$ using the Euler scheme described in Section \ref{sec:euler}.  At the jump dates we have
\begin{align}
X_{\tau_i}
	&=	X_{\tau_i-} + \int_\Rb \gam(\tau_i,X_{\tau_{i-}},z) \Del N(\tau_i,\dd z) 
	\stackrel{\Dc}{=} X_{\tau_i-} + \gam(\tau_i,X_{\tau_{i-}},Y_{i+1}) , &
Y_{i+1}
	&\sim F .
\end{align}
One advantage of simulating sample paths of $X$ on a random time grid is that we obtain for each sample path, the exact jump times for $X$.  When we simulate $X$ on a fixed time grid, we can only specify intervals $[t_i, t_{i+1})$ during which jumps occur.

\section{Variance Reduction Techniques}
\label{sec:VarReduction}
Having thoroughly discuses how generate iid sample paths for random processes, we now turn our attention to variance reduction techniques.  The basic idea of this section can be summarized as follows.  Suppose we wish to estimate $\mu = \Eb X$.  If we can generate iid random variables $(Y_i)$ such that $\Eb Y_i = \mu$ (where $Y_i$ may have a different distribution that $X$), we can construct an unbiased estimator of $\mu$ as follows
\begin{align}
\nuh_n
	&=	\frac{1}{n} \sum_{i=1}^n Y_i .
\end{align}
The variance of our estimator is given by
\begin{align}
\Vb \nuh_n
	&=	\sum_{i=1}^n \Vb (Y_i/n)
	=		\frac{1}{n^2} \sum_{i=1}^n \Vb Y_i 
	=		\frac{1}{n} \Vb Y .
\end{align}
The goal, then, is to cleverly choose $Y$ so that $\Vb Y$ is small.

\subsection{Control Variates}
Let us suppose that $(X,Y)$ are correlated random variables and that we wish to compute $\mu_Y := \Eb Y$.  Let us further suppose that we have a method of generating a sequence $(X_i,Y_i)$ of iid random vectors with the same distribution as $(X,Y)$.  Our usual method of estimating $\mu_Y$ is to compute
\begin{align}
\muh_n
	:=	\frac{1}{n} \sum_{i=1}^n Y_i .
\end{align} 
As $\Vb \muh_n = \sig_Y^2/n$ where $\sig_Y^2 := \Vb Y $, the only way to improve our estimate $\muh_n$ of $\mu_Y$ is to make $n$ larger.  Now, suppose we \emph{know} that $\Eb X = \mu_X$.  We can use this information to create a new estimator $\muh_n(b)$, which is given by
\begin{align}
\muh_n(b)
	&:=	\frac{1}{n} \sum_{i=1}^n Y_i(b) , &
Y_i(b)
	&:=	 Y_i - b ( X_i - \mu_X )  , &
b
	&\in \Rb .
\end{align}
Note that $\muh_n(b)$, like $\muh_n$ ,is an unbiased estimator of $\mu_Y$ as $\Eb \muh_n(b) = \mu$.  Moreover, we will show that, by choosing $b$ appropriately, we can obtain $\Vb \muh_n(b) \leq \Vb \mu_n$.  To see how this is so, we compute
\begin{align}
\Vb Y_i(b)
	&=	\Eb (Y_i(b) - \Eb Y_i(b) )^2 \\
	&=		\sig_Y^2 - 2 b \rho_{XY} \sig_X \sig_Y + b^2 \sig_X^2 =: \sig_Y^2(b),
\end{align}
where $\sig_X^2 = \Vb X$ and $\rho_{XY}$ is the correlation between $X$ and $Y$.  Note that $\sig_Y(b)$ is quadratic and convex as a function of $b$ and achieves a minimum at
\begin{align}
b^* 
	&= \frac{\rho_{XY} \sig_Y }{ \sig_X } 
	=		\frac{ \text{Co}\Vb(X,Y) }{ \Vb X} .
\end{align}
As $b^*$ is the minimizer of $\sig_Y(b)$, we know that $\sig_Y(b^*) \leq \sig_Y(0) = \sig_Y$.  Thus, we have
\begin{align}
\Vb \muh_n(b^*)
	&=	\frac{1}{n} \sig_Y(b^*)
	\leq	\frac{1}{n} \sig_Y
	=		\Vb \muh_n .
\end{align}
Thus, by using $\muh_n(b^*)$ in place of $\muh_n$, we can achieve a lower variance estimator of $\mu_Y$. 
The variable $X$ whose expectation $\mu_X$ is known, is called a \emph{control variate}. 
\\[0.5em]
In may practical scenarios, it may not be realistic to assume that we know $\text{Co}\Vb(X,Y)$ and $\Vb X$, which are needed to compute $b^*$.  However, we can always estimate $b^*$ using
\begin{align}
\widehat{b}_n^*
	&=	\frac{ \sum_{i=1}^n (X_i - \Eb X)( Y_i - \muh_n) }{\sum_{i=1}^n (X_i-\Eb X)^2} .
\end{align}
One cautionary note: unlike $\muh_n(b)$, which is an unbiased estimator of $\mu_Y$ for any fixed $b \in \Rb$, the estimator $\muh_n(\widehat{b}_n^*)$ is a biased estimator of $\mu_Y$ (because $\Eb \muh_n(\widehat{b}_n^*) \neq \mu_Y$).

\begin{example}
Suppose we wish to estimate $\Eb f(S_T)$ where $S = (S_t)_{t \geq 0}$ is the solution of the following SDE
\begin{align}
\dd S_t
	&=	a(S_t) \dd W_t .
\end{align}
Note that $S$ is a martingale and thus $\Eb S_T = S_0$.  Assuming we can simulate iid paths of $S$ (either exactly or approximately), then we can estimate $\Eb f(S_T)$ using the control variate method with $(X_i,Y_i) = (S_T^{(i)},f(S_T^{(i)}))$.
\end{example}

\subsection{Importance Sampling}
\label{sec:importance}
Consider a random variable $X$ defined on $(\Om,\Fc,\Pb)$.  Let $\Pbt$ be a probability measure satisfying $\Pbt(A) = 0 \Rightarrow \Pb(A) = 0$ (i.e., $\Pb \ll \Pbt$).  We can compute $\mu := \Eb \varphi(X)$ in two ways
\begin{align}
\Eb \varphi(X)
	&=	\Ebt Z \varphi(X) , &
Z
	&=	\frac{\dd \Pb}{\dd \Pbt} .
\end{align}
Let us denote by $F$ the distribution of $X$ under $\Pb$ and by $\Ft$ the distribution of $(X,Z)$ under $\Pbt$.  Suppose we can generate iid random variables $(X_i)$ with $X_i \sim F$ and iid random vectors $(X_i,Z_i)$ with $(X_i,Z_i) \sim \Ft$.  Then we can estimate $\mu $ in two ways
\begin{align}
\muh_n
	&=	\frac{1}{n} \sum_{i=1}^n \varphi(X_i) , &
X_i
	&\sim F , \\
\nuh_n
	&=	\frac{1}{n} \sum_{i=1}^n Z_i \varphi(X_i) , &
(X_i,Z_i)
	&\sim \Ft .
\end{align}
Both $\muh_n$ and $\nuh_n$ are unbiased estimators of $\mu$ as $\Eb \muh_n = \Eb \nuh_n = \mu$.  To see which estimator has a lower variance, we compute
\begin{align}
\Vb \varphi(X)
	&=	\Eb \varphi^2(X) - (\Eb \varphi(X))^2 , &
\Vbt Z \varphi(X)
	&=	\Ebt Z^2 \varphi^2(X) - (\Ebt Z \varphi(X))^2 .
\end{align}
Noting that $\Eb \varphi(X) = \Ebt Z \varphi(X)$, it follows that
\begin{align}
	&\text{if}&
\Eb \varphi^2(X) 
	&\geq \Ebt Z^2 \varphi^2(X) = \Eb Z \varphi^2(X), &
	&\text{then}&
\Vb \muh_n 
	&\geq \Vbt \nuh_n .
\end{align}
Now, let us consider the case where $Z = p(X)/\pt(X)$ where $p$ is the density of $X$ under $\Pb$ and $\pt$ is some other density.  Note that, with this choice for $Z$, we have $\Pbt( X \in \dd x) = \pt(x) \dd x$ because
\begin{align}
\Pb(X \in A)
	&= \Eb \Ib_{\{X \in A\}}
	=		\Ebt Z \Ib_{\{X \in A\}}
	=		\Ebt \frac{p(X)}{\pt(X)} \Ib_{\{X \in A\}} 
	=		\int_A \dd x \, \frac{p(x)}{\pt(x)} \pt(x)
	=		\int_A \dd x \, p(x) .
\end{align}
If $\pt$ were not the density of $X$ under $\Pbt$, then we would not have $\Pb(X \in A) = \int_A \dd x \, p(x)$.
The question arises: \emph{how should we choose $\pt$ in order to minimize $\Eb Z \varphi^2(X)$}?  
If we set
\begin{align}
\pt(x)
	&=	\frac{1}{c} p(x) |\varphi(x)| , &
c
	&=	\int \dd x \, p(x) |\varphi(x)| = \Eb |\varphi(X)|, \label{eq:p-tilde}
\end{align}
then we have
\begin{align}
\Eb Z \varphi^2(x)
	&=	\Eb \frac{p(X)}{\pt(X)} \varphi^2(X) 
	=		c \Eb  \frac{p(X)}{p(X) |\varphi(X)|} |\varphi^2(X) |
	=	( \Eb |\varphi(X)| )^2 
	\leq \Eb \varphi^2(X) .
\end{align}
Thus, if we choose $Z = p(X)/\pt(X)$ with $\pt$ given by \eqref{eq:p-tilde}, we will obtain an estimator $\nuh_n$ of $\mu$ with a lower variance than $\muh_n$.  Of course, from a practical standpoint, if we knew $p$, then we could compute $\Eb \varphi(X)$ without the need for Monte Carlo simulation.  Nevertheless, this exercise gives us insight as to how we might choose $Z$.  Namely, we should try to choose $Z$ so that $\pt(x)$ is large when $|\varphi(x)| p(x)$ is large.

\begin{example}
Let $X = (X_t)_{t \geq 0}$ be a process on $(\Om,\Fc,\Pb)$ defined by
\begin{align}
\dd X_t
	&=	a(X_t) \dd W_t ,
\end{align}
where $W$ is a Brownian motion under $\Pb$.  Suppose we wish to compute $\Pb( X_T > x ) = \Eb \Ib_{\{ X > x \} }$ where $x \gg X_0$.
We could approximate $X$ in a grid $\Pi = \{ t_0 = 0, t_1, t_2, \ldots, t_m = T \}$ with a process $\Xh$ where
\begin{align}
\Xh_{t_{i+1}}
	&=	\Xh_{t_i} + a(\Xh_{t_i}) (W_{t_{i+1}} - W_{t_i}) , &
i
	&=	0, 1, \ldots, m-1 .
\end{align}
Next, we can estimate $\mu:= \Pb( X_T > x ) = \Eb \Ib_{\{ X > x \} }$ using
\begin{align}
\muh_n
	&=	\frac{1}{n} \sum_{i=1}^n \Ib_{\{ \Xh_T^{(i)} > x \} } , &
	&\text{where $\Xh^{(i)}$ are iid realizations of $\Xh$ under $\Pb$.}
\end{align}
Alternatively, consider the change of measure
\begin{align}
Z
	&=	\frac{\dd \Pb}{\dd \Pbt}
	=		\ee^{ -\tfrac{1}{2} \gam^2 T - \gam \Wt_T } 
	=		\ee^{ -\tfrac{1}{2} \gam^2 T } \prod_{i=0}^{m-1} \ee^{ \gam (\Wt_{t_{i+1}} - \Wt_{t_i} ) } ,
\end{align}
for some $\gam > 0$.
The dynamics of $X$ can be written as
\begin{align}
\dd X_t
	&=	a(X_t) (\gam \dd t + \dd \Wt_t ) ,
\end{align}
where $\Wt$ is a Brownian motion under $\Pbt$.  We could approximate $X$ in the grid $\Pi$ with a process $\Xh$ where
\begin{align}
\Xh_{t_{i+1}}
	&=	\Xh_{t_i} + a(\Xh_{t_i}) \Big( \gam (t_{t+1} - t_i ) + (\Wt_{t_{i+1}} - \Wt_{t_i}) \Big) , &
i
	&=	0, 1, \ldots, m-1 .
\end{align}
Then, using
\begin{align}
\Pb( X_T >x)
	&=	\Ebt \frac{\dd \Pb}{ \dd \Pbt} \Ib_{\{ X_T > x \}}
	=		\Ebt Z \Ib_{\{ X_T > x \}} ,
\end{align}
we could compute
\begin{align}
\nuh_n
	&=	\frac{1}{n} \sum_{i=1}^n Z^{(i)} \Ib_{\{ \Xh_T^{(i)} > x \} } , &
	&\text{where $(\Xh^{(i)},Z{(i)})$ are iid realizations of $(\Xh,Z)$ under $\Pbt$.} , \\
Z^{(i)}
	&=	\ee^{ -\tfrac{1}{2} \gam^2 T } \prod_{j=0}^{m-1} \ee^{ \gam (\Wt_{t_{j+1}}^{(i)} - \Wt_{t_j}^{(i)} ) } .
\end{align}
As $X$ has a positive drift $\gam a(X)$ under $\Pbt$, it is likely that $\Vb \nuh_n < \Vb \muh_n$.
\end{example}

\subsection{Antithetic Variates}
Suppose a random variable $Y^+$ can be expressed as a function $g$ of a sequence of $m$ iid Standard Normal RVs $(Z_i)$
\begin{align}
Y^+ 
	&=	g(Z_1, Z_2, \ldots, Z_m).
\end{align}
Because $Z_i \sim - Z_i$ it follows that
\begin{align}
Y^-
	&:= g(-Z_1, -Z_2, \ldots, -Z_m)
	\stackrel{\Dc}{=} Y^+ .
\end{align}
We call the pair $(Y^+,Y^-)$ \emph{antithetic variables}.  Note that $Y^+$ and $Y^-$ are \emph{not} independent.
\\[0.5em]
Now, suppose we wish to compute $\mu := \Eb Y^+$.
As we know how to generate iid standard normal random variables $(Z_i)$ we can generate iid random variable $(Y_i^+)$ and $(Y_i^-)$.  
Thus, we can construct two unbiased estimators of $\mu$ as follows
\begin{align}
\muh_n
	&:=	\frac{1}{2n} \sum_{i=1}^{n} ( Y_i^+ + Y_{n+i}^+ ) , &
\nuh_n
	&:=	\frac{1}{2n} \sum_{i=1}^{n} ( Y_i^+ + Y_i^- ) ,
\end{align}
We would like to know, under which conditions we have $\Vb \nuh_n < \Vb \muh_n$.
Noting that
\begin{align}
\Vb (  Y_i^+ + Y_{n+i}^+ )
	&=	\Vb Y_i^+ + \Vb Y_{n+i} ^+
	=		2 \Vb Y^+ , \\
\Vb (  Y_i^+ + Y_i^- )
	&=	 \Vb Y_i^+ + \Vb Y_i^- + 2 \text{Co}\Vb( Y_i^+,Y_i^- ) 
	 =		2 \Vb Y^+ + 2 \text{Co}\Vb( Y^+,Y^- ) ,
\end{align}
we see that, if $\text{Co}\Vb( Y^+,Y^- ) < 0$, then we will have $\Vb \nuh_n < \Vb \muh_n$.  Note that we will have $\text{Co}\Vb( Y_i^+,Y_i^- ) < 0$ if the function $g$ is monotonic in each of its components, i.e., if 
\begin{align}
z_j
	&< y_j &
	&\Rightarrow&
g(z_1, \ldots, z_j, \ldots, z_m)
	&\leq g(z_1, \ldots, y_j, \ldots, z_m) , &
\forall \, j 
	&= 1, 2, \ldots m .
\end{align}

\begin{example}
Suppose we wish to compute $\Eb \varphi(X_T)$ where $X = (X_t)_{t \geq 0}$ satisfies the SDE $\dd X_t = \mu(t,X_t) \dd t + \sig(t,X_t) \dd W_t$.  
We can generate an approximate path of $X$ on a grid $\Pi=\{t_0=0, t_1, \ldots, t_m = T \}$ with a process $\Xh^+$ defined by
\begin{align}
\Xh_{t_{i+1}}^+
	&=	\Xh_{t_i}^+ + \mu( t_i, \Xh_{t_i}^+) (t_{i+1} - t_i) + \sig( t_i, \Xh_{t_i}^+) \sqrt{t_{i+1} - t_i} Z_{i+1} , &
i
	&=	0, 1, \ldots, m-1 .
\end{align}
Alternatively, we can generate an approximate path of $X$ on $\Pi$ with a process $\Xh^-$ defined by
\begin{align}
\Xh_{t_{i+1}}^-
	&=	\Xh_{t_i}^- + \mu( t_i, \Xh_{t_i}^-) (t_{i+1} - t_i) + \sig( t_i, \Xh_{t_i}^-) \sqrt{t_{i+1} - t_i} (- Z_{i+1} ), &
i
	&=	0, 1, \ldots, m-1 .
\end{align}
We can now construct antithetic variables $(Y^+,Y^-) = (\varphi(\Xh_T^+),\varphi(\Xh_T^-))$.
\end{example}


\section{Exercises}
\label{sec:exercises-13}
\blu{To do.}

%
%

%
%

%
%

\appendix

%
%

\chapter{Solutions to Exercises}


\section{Exercises from Section \ref{sec:exercises-1}}

\subsection*{Solution to Exercise \ref{exer:1.1}}
We would like to show that $\Gc$ contains the empty set, is closed under complementation, and is closed under countable unions. To this end, note that because $\emptyset \in \Fc$ then $\emptyset = \emptyset \cap B \in \Gc$. To see that $\Gc$ is closed under complementation, first choose a set $A \in \Fc$ and consider the set $A \cap B \subset B$. We would like to show that the complement of $A \cap B$ \textit{with respect to $B$} is a member of $\Gc$. The compliment of $A \cap B$ with respect to $B$ can be written as $\{ x \in B \mid x \not\in A \}= A^c \cap B.$ Since $A^c \in \Fc$, we have that $(A \cap B)^c \in \Gc$. Finally, let us show that $\Gc$ is closed under countable unions. Suppose that $\{ C_k \}_{k=1}^\infty \subset \Gc$ is a countable sequence of sets. For each $C_k$, there corresponds a set $A_k \in \Fc$ such that $C_k = A_k \cap B$. Let $A = \cup_{k=1}^\infty A_k$, and observe that
\begin{align}
\bigcup_{k=1}^\infty C_k 
	= \bigcup_{k = 1}^\infty \( A_k \bigcap B \) 
	= \( \bigcup_{k=1}^\infty A_k \) \bigcap B
	= A \bigcap B.
\end{align}
Since $A \in \Fc$, $\cup_{k=1}^\infty C_k \in \Gc$. Hence, $\Gc$ is a $\sig$-algebra over $B$. 

\subsection*{Solution to Exercise \ref{exer:1.2}}
\textbf{(a)}
We would like to show that $\Fc \cap \Gc$ is a $\sig$-algebra i.e. $\Fc \cap \Gc$ contains the empty set, is closed under complementation, and is closed under countable unions. First, since $\emptyset \in \Fc$ and $\emptyset \in \Gc$, then $\emptyset \in \Fc \cap \Gc$. Next, choose a set $A \in \Fc \cap \Gc$. Since $A$ is in both $\Fc$ and $\Gc$ (and each are closed under complementation), then $A^c \in \Fc$ and $A^c \in \Gc$. Therefore, $A^c \in \Fc \cap \Gc$. Finally, suppose that $\{ A_k \}_{k=1}^\infty \subset \Fc \cap \Gc$ is a countable sequence of sets and denote $A = \cup_{k=1}^\infty A_k$. Since $A \subset \Fc$ and $A \subset \Gc$, then $A \in \Fc$ and $A \in \Gc$. Therefore, $A \in \Fc \cap \Gc$. Hence, $\Fc \cap \Gc$ is a $\sig$-algebra. 
\textbf{(b)}
Recall that in general, if $A \subset \Om$, then $\{ \es, \Om, A, A^c \}$ is a $\sig$-algebra. Let $\Om = \{ a, b , c \}$ and let $\Fc$ and $\Gc$ be the $\sig$-algebras over $\Om$ given by
\begin{align}
&\Fc
	= \{ \es, \Om, \{ a \}, \{ b, c \} \},
&\Gc
	= \{ \es, \Om, \{ b \}, \{ a, c \} \}.
\end{align}
Note that the sets $\{ a \}$ and $\{ b \}$ are both members of $\Fc \cup \Gc$. However, $\{ a, b\} = \{ a \} \cup \{ b \} \not\in \Fc \cup \Gc$ which implies that $\Fc \cup \Gc$ is not closed under countable unions. Thus, $\Fc \cup \Gc$ is not a $\sig$-algebra. 

\subsection*{Solution to Exercise \ref{exer:1.3}}
\textbf{(a)}
Let $p$ be the probability that the biased coin is $H$, and $q = 1 - p$ the probability that this coin is $T$. The sample space corresponding to the experiment of flipping this coin three times is given by
\begin{align}
\Om
	= \{ 
		\om= \om_1 \om_2 \om_3  \mid  \om_k \in \{ H, T\}, k = 1,2,3
	\}.
\end{align}
Since $\Om$ is finite, $\Fc$ can be taken to be $2^\Om$ without issue. Furthermore, the finite cardinality of $\Om$ allows $\Pb$ to be completely determined once it is defined for singleton sets. To this end, since the coin tosses are independent
\begin{align}
&\Pb\( \{ \om \} \)
	= p^{ h(\om) } q^{3 - h(\om)},
& h(\om)
	= \sum_{k=1}^3 \Ib_{\{w_k = H\}} (\om).
\end{align}
\textbf{(b)}
The sample space for this event is given by
\begin{align}
\Om
	= \{ 
		\om = \om_1 \om_2 \mid \om_k \in \{ R, B\}, k = 1,2
	\}.
\end{align}
As with the last exercise, since $\Om$ is finite $\Fc$ can be taken to be $2^\Om$ without issue. Furthermore, because $\Om$ is finite once $\Pb$ is specified on singleton sets $\Pb$ is completely determined. It is known that the two color urn experiment without replacement follows the \textit{hypergeometric distribution}. Suppose that we take drawing a red ball to be a success (choosing blue as a success does not change the distribution). In our case, $N = 4$ is the population size (4 balls), $K = 2$ is the number of successes (there are 2 red balls), and there are $n =2 $ draws from the urn without replacement. For $\om \in \Om$, let $k(\om)$ be the number of red draws i.e. successes. It is known that the probability of $k$ successes for a hypergeometric random variable is given by
\begin{align}
\Pb( \{ \om \mid k(\om)= \kappa\} ) 
	= \dfrac{ \binom{K}{\kappa} \binom{N - K}{n - \kappa}}{ \binom{ N } { n} }.
\end{align}
Therefore, the probability measure is given by
\begin{align}
\Pb\( \{ \om \} \) 
	= \dfrac{ \binom{K}{k(\om)} \binom{N - K}{n - k(\om)}}{ \binom{ N } { n} \binom{n}{k(\om)}}
	= \frac{1}{6} \binom{2}{2-k(\om)}
	= \frac{1}{3 (2-k(\om))! k(\om)!} .
\end{align}
\textbf{(c)}
Let $A_k$ be the event that success happens on the $k$-th trial for $k = 1,2,3,\cdots$ and let $A_\infty$ 
denote the event that success never occurs. Then
\begin{align}
\Om
	= \( \bigcup_{k=1}^\infty A_k \) \cup \{ A_\infty \}
\end{align}
Since $\Om$ is countable, then $\Fc$ can be taken to be $2^\Om$ without issue, and once $\Pb$ is defined on each $A_k$, $k = 1,2,3,\cdots$, then $\Pb$ is completely determined. This process follows the \textit{geometric distribution}, therefore
\begin{align}
\Pb( A_k)
	= p(1 - p)^{k-1},
\end{align}
where $p$ is the probability that the coin flip is $H$. Using continuity of probability measures, one can show that 
$\Pb( A_\infty) = 0$. 

\subsection*{Solution to Exercise \ref{exer:1.4}}
\textbf{(a)}
Assume that $g$ is defined on all of $\Rb$. Since $g$ is a strictly increasing, continuous function then $g(\Rb) = I = (a,b)$, where $I$ is a potentially infinite (or half-infinite) interval. Furthermore, $g^{-1}$ exists on $I$. Define 
\begin{align}
v(x) 
	= \left\{
	\begin{array}{lr}
	 -\infty, & x \leq a \\
	 g^{-1}(x), & a < x < b \\
	 \infty, & x \geq b
	\end{array}
	\right. .
\end{align}
Then 
\begin{align}
F_Y(x)
	=\Pb( g(X) \leq x)
	= \Pb( X \leq v(x) )
	= F_X(v(x)).
\end{align}
\textbf{(b)}
Since $v$ is continuous, $F_Y$ is continuous. Therefore, 
\begin{align}
f_Y(x) 
	= F_Y'(x) 
	= F_X(v(x))'
	= f_X(v(x))\cdot v'(x).
\end{align}

\subsection*{Solution to Exercise \ref{exer:1.5}}
\textbf{(a)}
Suppose that $Y  = X^2$. Then for $y > 0$,
\begin{align}
F_Y(y)
	&= \Pb( X^2 \leq y) \\
	&= \Pb(- \sqrt y \leq X \leq \sqrt y) \\
	&= \Pb(- \sqrt y <  X \leq \sqrt y) \\
	&= F_X(\sqrt y) - \lim_{u\to y^-}F_X(- \sqrt u) \\
	&= F_X(\sqrt y) - F_X ( - \sqrt y). 
\end{align}
The last equality is justified by the continuity of $X$. Thus, 
\begin{align}
F_Y(y)
	= \left \{
	\begin{array}{lr}
	0, & y \leq 0 \\
	F_X(\sqrt y) - F_X(- \sqrt y), & y > 0
	\end{array}
	\right. .
\end{align}
\textbf{(b)}
Suppose that $Y = \sqrt{ | X |}$. Then for $y > 0$, 
\begin{align}
F_Y(y)
	= \Pb( \sqrt{ | X |} \leq y)
	= \Pb( - y^2 \leq X \leq y^2)
	= F_X(y^2) - F_X(- y^2).
\end{align}
Thus, 
\begin{align}
F_Y(y)
	= \left \{
	\begin{array}{lr}
	0, & y \leq 0 \\
	F_X(y^2) - F_X(- y^2), & y > 0
	\end{array}
	\right. .
\end{align}
\textbf{(c)}
Let $Y = \sin X$. For $y \leq -1$, $F_Y(y) = 0$. Furthermore, if $y \geq 1$, $F_Y(y) = 1$. Now suppose that $-1 \leq y 
\leq 1$ and let $\arcsin(y) \leq X < \arcsin(y) + 2\pi$. Then
\begin{align}
\{ X \in [\arcsin(y), \arcsin(y) + 2\pi) \mid \sin(X)\leq y\}
	= \{X \in [\pi - \arcsin(y), 2\pi + \arcsin(y))\}.
\end{align}
Thus, for general $X$ and $-1 \leq y \leq 1$,
\begin{align}
\label{eqn.sintrans1}
F_Y(y)
	= \sum_{n = -\infty}^\infty 
		\left( F_X ((2n +2)\pi + \arcsin(y)) - F_X((2n+1)\pi - \arcsin(y))
			\right)
\end{align}
\textbf{(d)}
Note the if $y < 0$ or $y > 1$, then $F_Y(y) = 0$. Suppose that $0 \leq y \leq 1$. Since $X$ is a continuous random 
variable, $F_X(x)$ is a continuous function. This fact, along with the monotonicity of $F_X$, implies that there exists 
a largest number $\xi(y)$ such that $F_X(\xi(y)) = y$. Thus, $F_X(X) \leq y$ if and only if $X \leq \xi(y)$. Hence, 
$F_Y(y) = F_X(\xi(y))$.
\\[0.5em]
\textbf{\underline{Supplement}}:\\
In this question, $F_X(X)$ is continuous, but it is not necessarily invertible (e.g. $F_X(X)$ is constant in some intervals). Therefore, we should define a function of $y$ instead of taking inverse of $F_X$ directly. This is one way:\\
$$\xi(y) = \text{sup} \{ a : F_X(a) \le y \}$$
And you can check this is true:
$$\Pb(F_X(X) \le y) = \Pb(X \le \xi(y))  $$
If $X$ is a continuous random variable, there is no difference to take maximum or supremum. We use supremum here because this is more general for discrete random variables. 

\subsection*{Solution to Exercise \ref{exer:1.6}}
\textbf{(a)} Since $f$ and $g$ are the densities of a random variables, we clearly have $Z = g(X)/f(X) \geq 0$.  We also have
\begin{align}
\Eb Z
	&=	\Eb g(X) = \int \frac{g(x)}{f(x)}f(x) \dd x = \int g(x) \dd x = 1.
\end{align}
Thus, $Z$ defines a Radon-Nikod\'ym derivative.
\textbf{(b)} We have
\begin{align}
\Pbt(X \leq x)
	&=	\Eb Z \Ib_{\{X \leq x\}}
	=		\int_{-\infty}^x \frac{g(x)}{f(x)} f(x) \dd x
	-		\int_{-\infty}^x g(x) \dd x .
\end{align}
Thus, under $\Pbt$ the random variable $X$ has density $g$.

\subsection*{Solution to Exercise \ref{exer:1.7}}
Assume that $g:[0,1] \to \Rb_+$ is invertible.  Then we have
\begin{align}
1 - \ee^{-x}
	&=	\Pb( Y \leq x)
	=		\Pb( g(X) \leq x )
	=		\Pb( X \leq g^{-1}(x) )
	=		g^{-1}(x) ,
\end{align}
where, in the last step we have used the fact that $X \sim \Uc([0,1[)$.  Hence, we find $g(x) = - \log (1- x)$.  Since $X$ has the same distribution as $1-X$, an alternative and equally valid solution is $g(x) = - \log x$.


\section{Exercises from Section \ref{sec:exercises-2}}

\subsection*{Solution to Exercise \ref{exer:2.1}}
(a)
Let $A = \{ a, b\}$ and $B = \{c ,d\}$. For $\omega \in \Omega$, the value of $X(\omega)$ can be determined by 
determining the membership of $\omega$ to either $A$ or $B = A^c$. Therefore, $\sigma(X) = \{ \emptyset, \Omega, 
A,B\}$. 
(b)
To compute the value of $\Eb[ Y | X](\omega)$ for each $\omega \in \Omega$, one can appeal to the formula
\begin{align}
&\Eb[Y | X] 
	= \psi(X),
&&\psi(x)
	= \Eb[ Y | X = x].
\end{align}
Let $A = \{ a,b \}$ and $B = \{ c,d \}$. Note that 
\begin{align}
\Eb[ Y | X = 1]
	= Y(a) \Pb( \{ a \} | A) + Y(b) \Pb( \{ b \} | A)
	= 1 \( \frac{1}{3} \) - 1 \( \frac{2}{3}\) 
	= -\frac{1}{3},
\end{align}
and
\begin{align}
\Eb[ Y | X = -1]
	= Y(c) \Pb( \{ c \} | B)
		+ Y(d) \Pb( \{ d \} | B)
	= 1 \( \frac{1}{2} \) - 1 \( \frac{1}{2} \)
	= 0.
\end{align}
Therefore, 
\begin{align}
&\Eb[Y | X](a) 
	= \Eb[ Y | X ] (b) 
	= - \frac{1}{3},
&&\Eb[ Y | X](c)
	= \Eb[ Y | X](d)
	= 0.
\end{align}
The partial averaging property holds trivially for the empty set. Furthermore, since $\sigma(X) \subset \Fc$, 
\begin{align}
\Eb[ \Ib_\Omega \Eb[ Y | X]]
	= \frac{1}{6} \(- \frac{1}{3} \)
	+ \frac{1}{3} \( - \frac{1}{3} \)
	+ \frac{1}{4} \( 1 \)
	+ \frac{1}{4} \( -1 \) 
	= - \frac{1}{6},	
\end{align}
and
\begin{align}
\Eb[ \Ib_\Omega Y]
	= \frac{1}{6} \(1 \)
	+ \frac{1}{3} \( - 1 \)
	+ \frac{1}{4} \( 1 \)
	+ \frac{1}{4} \( -1 \)  
	= - \frac{1}{6}.
\end{align}
It just remains to be checked for sets $A$ and $B$. To this end, observe
\begin{align}
&\Eb[ \Ib_A \Eb[ Y | X ]]
	= \frac{1}{6}\Eb[ Y | X](a)  + \frac{1}{3}\Eb[ Y|X](b) 
	= -\frac{1}{6},
&&\Eb[ \Ib_B \Eb[Y|X]]
	= 0,
\end{align}
and
\begin{align}
&\Eb[ \Ib_A Y]
	= \frac{1}{6} Y(a) + \frac{1}{3} Y(b)
	= - \frac{1}{6},
&&\Eb[ \Ib_B Y]
	= \frac{1}{4} Y(c)
		+ \frac{1}{4} Y(d)
	= 0.
\end{align}
Hence, partial averaging holds. 
(c)
First, $\Eb[ Z | X]$ can be rewritten as
\begin{align}
\Eb[ Z | X] 
	= \Eb[ Y | X] + \Eb[X | X]
	= \Eb[ Y | X] + X.
\end{align}
Thus, 
\begin{align}
&\Eb[Z |X](a)
	= \frac{2}{3},
&&\Eb[Z |X](b)
	= \frac{2}{3},
&\Eb[Z |X](c)
	= -1,
&&\Eb[Z |X](d)
	= -1.
\end{align}
To check that the partial averaging property holds,
\begin{align}
\Eb[ \Ib_\Omega \Eb[ Z | X]]
	= \frac{1}{6} \(\frac{2}{3} \)
	+ \frac{1}{3} \( \frac{2}{3} \)
	+ \frac{1}{4} \( -1 \)
	+ \frac{1}{4} \( -1 \) 
	= - \frac{1}{6},
\end{align}
and
\begin{align}
\Eb[\Ib_\Omega Z]
	= \frac{1}{6} \(2 \)
	+ \frac{1}{3} \( 0 \)
	+ \frac{1}{4} \( 0 \)
	+ \frac{1}{4} \( -2 \) 
	= - \frac{1}{6}.
\end{align}
It just remains to be checked for sets $A$ and $B$. To this end, observe
\begin{align}
&\Eb[ \Ib_A \Eb[Z|X]]
	= \frac{1}{3} \( \frac{2}{3} \)
		+ \frac{2}{3} \( \frac{2 }{3}\)
	= \frac{2}{3},
&&\Eb[ \Ib_B \Eb[Z|X]]
	= \frac{1}{2} \( -1 \) 
		+ \frac{1}{2} \( -1 \)
	= -1,
\end{align}
and
\begin{align}
&\Eb[ \Ib_A Z]
	= \frac{1}{3} \( 2 \)
		+ \frac{2}{3} \( 0\)
	= \frac{2}{3},
&&\Eb[\Ib_B Z]
	= \frac{1}{2} \( 0 \) 
		+ \frac{1}{2}\( -2 \) 
	= - 1.
\end{align}
Hence, partial averaging holds. 

\subsection*{Solution to Exercise \ref{exer:2.2}}
Define $\mu := \Eb(Y-X)$.  Then we have
\begin{align}
\Vb (Y-X)
	&=	\Eb (Y-X-\mu)^2 \\
	&=	\Eb \Big( (Y - \Eb[Y|\Gc]) + (\Eb[Y|\Gc]-X-\mu) \Big)^2 \\
	&=	\Eb (Y - \Eb[Y|\Gc])^2 + 2 \Eb (Y - \Eb[Y|\Gc])(\Eb[Y|\Gc]-X-\mu) + \Eb (\Eb[Y|\Gc]-X-\mu)^2 \\
	&=	\Vb (Y - \Eb[Y|\Gc]) + 2 \Eb (Y - \Eb[Y|\Gc])(\Eb[Y|\Gc]-X-\mu) + \Eb (\Eb[Y|\Gc]-X-\mu)^2 ,
\end{align}
where we have used $\Eb (Y - \Eb[Y|\Gc])^2 = \Vb (Y - \Eb[Y|\Gc])$, which follows from $\Eb(Y - \Eb[Y|\Gc])=0$.  The third term in the last line above is clearly non-negative.  Thus, we have
\begin{align}
\Vb (Y-X)
	&\geq	\Vb (Y - \Eb[Y|\Gc]) + 2 \Eb (Y - \Eb[Y|\Gc])(\Eb[Y|\Gc]-X-\mu) .
\end{align}
From the definition of Conditional expectation we have $\Eb Z(Y - \Eb[Y|\Gc])=0$ for any $Z \in \Gc$.  As $(\Eb[Y|\Gc]-X-\mu) \in \Gc$, the second term on the right-hand side is zero.  And thus $\Vb (Y-X) \geq	\Vb (Y - \Eb[Y|\Gc])$.

\subsection*{Solution to Exercise \ref{exer:2.3}}
Let $\Omega = \{ 0, 1, 2, 3\}$,  $X(\om) = \om$, $f(x) = \Ib_{[0,1]}(x) - \Ib_{[2,3]}(x)$, and $g(x) \equiv 1$. Note that since $X$ takes 
a distinct value for each $\omega \in \Omega$, $\sig (X) = \sig( \{ \{0 \}, \{1 \}, \{ 2 \}, \{ 3 \} \}) =  2^\Omega$. 
The function $f(X)$ can be completely resolved by determining membership to the sets $A = \{ 0,1\}$ and $A^c = \{ 
2,3\}$. Therefore, $\sig(f(X)) = \sig( \{ A, A^c \} ) = \{ \emptyset, \Omega, A, A^c \}$. Finally, since $g$ is 
constant, $g ( X) = 1$ for all $\omega \in \Omega$. Therefore, $\sig(g(X)) = \{ \emptyset, \Omega \}$. 

\subsection*{Solution to Exercise \ref{exer:2.4}}
Choose integers $m$ and $n$ such that $n \geq 0$ and $m \geq 1$. We have that $\Fc_n \subset \Fc_{n+m}$. By definition, 
$X_{n+m} = \Eb[ X |\Fc_{n+m}]$. Therefore, from iterated conditioning,
\begin{align}
\Eb[ X_{n+m} | \Fc_n]
	= \Eb[ \Eb[ X | \Fc_{n+m} ] | \Fc_n]
	= \Eb[ X | \Fc_n]
	= X_n.
\end{align}
Hence, $X_n$ is a martingale with respect to the filtration $(\Fc_n)_{n \geq 1}$.

\subsection*{Solution to Exercise \ref{exer:2.5}}
(a) Let $q := 1-p$.  Choose an integer $n \geq 0$. Then
\begin{align}
\Eb[ Z_{n+1} | \Fc_n]
	&= \Eb\[ r^{2 S_{n+1} - (n+1)} \middle\vert \Fc_{n}\]
	 = \Eb\[ r^{2S_n - n} r^{2X_{n+1} - 1}  \middle\vert \Fc_{n}\]
	 = Z_n  \Eb\[ r^{2X_{n+1} -1} \middle\vert \Fc_{n}\].
\end{align}
Since $X_{n+1}$ is independent of $\Fc_n$,
\begin{align}
\Eb[ Z_{n+1} | \Fc_n]
	&= Z_n \Eb\[ r^{2X_{n+1} -1} \]
	 = Z_n \Big( p r + \frac{q}{r} \Big) .
\end{align}
As $Z_n > 0$, we have $\Eb[ Z_{n+1} | \Fc_n] = Z_n$ (resp. $\geq$, $\leq$) if and only if $p r + q/r = 1$ (resp. $\geq 1$, $\leq 1$).  Multiplying by $r > 0$ and noting that
\begin{align}
p r^2 - r + q
	&=	(r - 1)(p r - q) ,
\end{align}
we find that the one-step condition $\Eb[ Z_{n+1} | \Fc_n] = Z_n$ holds if and only if $r = 1$ or $r = q/p = (1-p)/p$.  Likewise, the one-step submartingale condition $\Eb[ Z_{n+1} | \Fc_n] \geq Z_n$ holds if and only if $(r-1)(pr - q) \geq 0$, i.e., $r \leq \min(1,q/p)$ or $r \geq \max(1,q/p)$, and the one-step supermartingale condition $\Eb[ Z_{n+1} | \Fc_n] \leq Z_n$ holds if and only if $(r-1)(pr - q) \leq 0$, i.e., $\min(1,q/p) \leq r \leq \max(1,q/p)$.
(Note that $r = 1$ gives the trivial martingale $Z_n \equiv 1$.)
\\[0.5em]
To see that the one-step conditions imply the corresponding conditions over multiple steps, consider the martingale case and choose an integer $m \geq 1$. Since $n$ was arbitrary,
\begin{align}
\Eb [ Z_{n+m} | \Fc_n]
	&= \Eb[ \Eb[ Z_{n+m} | \Fc_{n+m-1} ] | \Fc_n] 
	 = \Eb[ Z_{n + m - 1} | \Fc_n] \\
	&= \Eb[ \Eb[ Z_{n+m-1} | \Fc_{n + m - 2}] | \Fc_n] 
	 = \Eb[ Z_{n+ m - 2} | \Fc_n] \\
	&\qquad \vdots \\
	&= \Eb[ Z_{n +2 } | \Fc_n]  
	 = \Eb[ \Eb[ Z_{n+2} | \Fc_{n+1} ] | \Fc_n ] \\
	&= \Eb[ Z_{n+1} | \Fc_n] 
	 = Z_n.
\end{align}
Thus, when $r = 1$ or $r = (1-p)/p$, the process $Z$ is a martingale with respect to the filtration $( \Fc_n)_{n\geq 0}$.  The same argument, with equalities replaced by inequalities, handles the sub- and supermartingale cases.
\\[0.5em]
(b) Yes.  Since $S_0 = 0$ and $S_k = \sum_{i=1}^k X_i$, we have $\Gc_n = \sig(X_1,\ldots,X_n) \subseteq \Fc_n$.  Moreover, $Z_n$ is a function of $S_n$ and is therefore $\Gc_n$-measurable.  Using the tower property, we obtain
\begin{align}
\Eb[ Z_{n+1} | \Gc_n ]
	&=	\Eb[ \Eb[ Z_{n+1} | \Fc_n ] | \Gc_n ]
	=		\Eb[ Z_n | \Gc_n ]
	=		Z_n .
\end{align}
Thus, $Z$ is also a $(\Pb,\Gb)$ martingale.
\\[0.5em]
(c) Repeating the computation in part (a) under $\Pbt$, with $\widetilde{q} := 1 - \pt$, we obtain
\begin{align}
\Ebt[ Z_{n+1} | \Fc_n]
	&=	Z_n \Big( \pt r + \frac{\widetilde{q}}{r} \Big) ,
\end{align}
which is equal to $Z_n$ if and only if $r = 1$ or $r = \widetilde{q}/\pt$.  If $r = 1$, then $Z_n \equiv 1$ is trivially a $(\Pbt,\Fb)$ martingale.  If $r = q/p \neq 1$, then $Z$ is a $(\Pbt,\Fb)$ martingale if and only if $q/p = \widetilde{q}/\pt$, which holds if and only if $\pt = p$.  Since $\pt \neq p$, we conclude that $Z$ is a $(\Pbt,\Fb)$ martingale only in the trivial case $r = 1$.


\section{Exercises from Section \ref{sec:exercises-3}}

\subsection*{Solution to Exercise \ref{exer:3.1}}
We have
\begin{align}
\Eb[ s^X]
	&= \Eb[ \Eb[ s^X | U] ] \\
	&= \Eb[ \( (1-U) + s U \)^n ] \\
	&= \Eb\[ \sum_{k=0}^n \binom{n}{k} s^k U^k (1-U)^{n-k} \] \\
	&= \sum_{k=0}^n \binom{n}{k} s^k \Eb[ U^k (1-U)^{n-k} ] \\
	&= \sum_{k=0}^n \binom{n}{k} s^k \frac{n! (n-k)!}{\Gamma(2 +n )} \\
	&= \sum_{k=0}^n \frac{1}{n+1} s^k.
\end{align}
This implies that $\Pb( X = k) = (n+1)^{-1}$ i.e. $X$ is uniformly distributed on $\{0,1,2,\cdots,n\}$. 

\subsection*{Solution to Exercise \ref{exer:3.2}}
Let $k = n - m$ and denote the number of descendants of agent $i$ from step $m$ to $n$ as $Y_{m,i}$. Then,
\begin{align}
Z_n
	= Y_{m,1} + Y_{m,2} + Y_{m,3} + \cdots + Y_{m,Z_m},
\end{align}
the $Y_{m,i}$ are i.i.d. with $Y_{m,i} \sim Z_k$ and each $Y_{m,i}$ is independent of $Z_m$, for $i = 1,2,3,\cdots Z_m$. Therefore, 
\begin{align}
\Eb[ Z_m Z_n]
	&= \Eb[ \Eb[ Z_m Z_n | Z_m] ] \\
	&= \Eb[ Z_m \Eb[ Z_n | Z_m] ] \\
	&= \Eb\[ Z_m \sum_{i=1}^{Z_m} \Eb[ Y_{m,i} | Z_m ] \] \\
	&= \Eb\[ Z_m \sum_{i=1}^{Z_m} \Eb[ Y_{m,i} ] \] \\
	&= \Eb[ Z_m^2 \mu^{k} ] \\
	&= \mu^{n-m} \Eb[ Z_m^2].
\end{align}
Let $\Vb Z_1 = \sigma^2 > 0$ and assume that $\mu \neq 1$. Note that
\begin{align}
\operatorname{Cov}[ Z_m Z_n]
	= \Eb[ Z_m Z_n] - \Eb[Z_m] \Eb[ Z_n]
	&= \mu^{n-m} \Eb[ Z_m^2] - \mu^{n+m}\\
	&= \mu^{n-m} ( \Vb[ Z_m] + \Eb[ Z_m]^2) - \mu^{n+m} \\
	&= \mu^{n-m} \(\frac{ \sig^2 (\mu^m - 1) \mu^{m-1}}{\mu-1} + \mu^{2m}\) - \mu^{n+m} \\
	&= \frac{\sigma ^2 \left(\mu ^m-1\right) \mu ^{n-1}}{\mu -1}
\end{align}
and
\begin{align}
\sqrt{\Vb Z_n} \sqrt{\Vb Z_m}
	= \sqrt{\frac{ \sig^2 (\mu^n - 1) \mu^{n-1}}{\mu-1}
		\frac{ \sig^2 (\mu^m - 1) \mu^{m-1}}{\mu-1}
	}
	= \frac{\sig^2}{|\mu -1|} 
	\sqrt{ \mu^{n+m-2} (\mu^n - 1)(\mu^m - 1)}
\end{align}
Thus,
\begin{align}
\rho(Z_m,Z_n)
	= \text{sign}(\mu - 1) \frac{\left(\mu ^m-1\right) \mu ^{n-1}}{\sqrt{ \mu^{n+m-2} (\mu^n - 1)(\mu^m - 1)}}
\end{align}
Now supposed that $\mu = 1$. Then 
\begin{align}
\operatorname{Cov}
	= \Eb[Z_n Z_m] - \Eb[Z_m]\Eb[Z_n]
	&= \Eb[Z_m]^2 - 1 \\
	&= \Vb[ Z_m] + \Eb[Z_m]^2 - 1 \\
	&= m\sig^2
\end{align}
and $\sqrt{\Vb Z_n} \sqrt{\Vb Z_m} = \sig^2 \sqrt{ mn}$. Thus,
\begin{align}
\rho(Z_m,Z_n)
	= \sqrt{ \frac{m}{n}}.
\end{align}

\subsection*{Solution to Exercise \ref{exer:3.3}}
Suppose we choose a node at the $n$-th level, $i_{n_1}$ who has ancestor $i_{r_k}$ at the $r$-th generation, $n > r$. 
Conditioning on $Z_r$, the probability that $i_{n_2}$ drawn at random from generation $n$ (with replacement) has 
ancestor $i_{r_k}$ is $Z_r^{-1}$ i.e. $\Pb( L \geq r | Z_r) = Z_r^{-1}$. Then
\begin{align}
\Pb( L \geq r)
	= \Eb[ \Ib_{L \geq r} ]
	= \Eb[ \Eb[ \Ib_{L \geq r} | Z_r] ]
	= \Eb[ Z_r^{-1}].
\end{align}
The event of interest is $\{ L = r\} = \{ L \geq r \} \setminus \{ L \geq r + 1 \}$. Since $\{ L \geq r+1 \} \subset \{ 
L \geq r \}$, 
\begin{align}
\Pb ( L = r )
	= \Pb( L \geq r ) - \Pb(L \geq r+1) 
	= \Eb[ Z_r^{-1}] - \Eb[ Z_{r+1}^{-1}].
\end{align}

\subsection*{Solution to Exercise \ref{exer:3.4}}
Let $Z_{n} = \widehat Z_n + Y_n$, where $\widehat Z_n$ is the number of offspring produced by generation $n-1$ i.e. $\widehat Z_n = 
\sum_{i=1}^{Z_{n-1}} X_{n-1,i}$. Following the proof Theorem 3.12 from the notes, it can be seen that $G_{\widehat Z_n}(s) 
= G_{Z_{n-1}} ( G_X(s))$. 
Now note that $Z_1 = \widehat Z_1 + Y_1$. Since $\widehat Z_1$ is simply the number of offspring produced in the first 
iteration, 
\begin{align}
&\Pb ( Z_1 = k)
	= \sum_{j=0}^k \Pb( \widehat Z_1 = j) \Pb ( Y_1 = k - j)
&&\Rightarrow
&& G_{Z_1}(s)
	= G_Y(s)G_{\widehat Z_1}(s)
	= G_Y(s)G_X(s).	
\end{align}
Furthermore, by the same reasoning as above, $G_{Z_2}(s) = G_Y(s) G_{\widehat Z_2} (s)$. So, 
\begin{align}
G_{Z_2}(s)
	= G_Y(s) G_{\widehat Z_2}(s)
	= G_Y(s) G_{Z_1} (G_X(s))
	= G_Y(s) G_Y( G_X(s)) G_X ( G_X(s)).
\end{align}
In general, $G_{Z_n} = G_Y(s) G_{\widehat Z_n}(s)$. Let the superscript $f^{(n)}(s) = f(f(f(\cdots f(s) \cdots )$ denote 
composition $n$ times with the convention that $f^{(0)}(s) = s$. Then 
\begin{align}
G_{Z_n}(s)
	&= G_Y(s) G_{\widehat Z_n}(s) \\
	&= G_Y(s) G_{Z_{n-1}}(G_X(s)) \\
	&= G_Y(s) G_Y(G_X(s)) G_{\widehat Z_{n-1}} (G_X(s)) \\
	& \hspace{2cm} \vdots \\
	&= G_X^{(n)}(s) \prod_{j=0}^{n-1} G_Y\( G_X^{(j)}(s) \)
\end{align}

\subsection*{Solution to Exercise \ref{exer:3.5}}
Note that
\begin{align}
\Eb[ e^{it X^2}]
	= \int_\Rb dx \frac{1}{\sqrt{ 2 \pi \sig^2}} e^{- (x - \mu)^2 / (2 \sig^2)} e^{i t x^2}
	= \int_\Rb dx \frac{1}{\sqrt{ 2 \pi \sig^2}} e^{ - \frac{(1 - 2 i t \sig^2)x^2 - 2 \mu x + \mu^2}{2 \sig^2}}
\end{align}
Letting $y = \alpha x$, where $\alpha = \sqrt{1 - 2i t \sig^2}$, and completing the square we get
\begin{align}
\frac{1}{\alpha} \int_{\Rb} dy\ \frac{1}{\sqrt{ 2 \pi \sig^2}}
	e^{ - \frac{ y^2 - 2 \mu \alpha^{-1} y + \mu^2}{2 \sig^2}}
	&=\frac{1}{\alpha} \int_{\Rb} dy\ \frac{1}{\sqrt{ 2 \pi \sig^2}}
	e^{ - \frac{ (y - \mu \alpha^{-1})^2}{2 \sig^2}} e^{- \frac{(\alpha^2 -1 )\mu^2}{2 \alpha^2 \sig^2}} \\
	&= \frac{1}{\alpha} e^{- \frac{(\alpha^2 -1 )\mu^2}{2 \alpha^2 \sig^2}} \\
	&= \frac{e^{-\frac{\mu ^2 t}{2 \sigma ^2 t+i}}}{\sqrt{1-2 i \sigma ^2 t}}
\end{align}

\subsection*{Solution to Exercise \ref{exer:3.6}}
Note that $F_{X_n}(x)$ is continuous. We would first like to show that $\lim_{x\to-\infty} F_{X_n}(x) = 0, 
\lim_{x\to\infty} F_{X_n}(x) = 1$ and $F_{X_n}(x)$ is monotonically increasing. The two limits are obvious. Note that 
\begin{align}
\label{eqn.der1}
F_{X_n}'(x)
	= 
\left\{
\begin{array}{lr}
0, & x \in (-\infty,0) \cup (1,\infty) \\
1 - \cos(2 \pi n x), & x \in [0,1]
\end{array}
\right. .
\end{align}
It can be seen from \eqref{eqn.der1} that $F_{X_n}'(x)$ is non-negative and continuous on $\Rb$. Therefore, 
$F_{X_n}(x)$ is monotonically increasing. Since $F_{X_n}$ is continuously differentiable, $f_{X_n} = F_{X_n}'$ is the 
density function of $X_n$. 
Since $|\sin(2 \pi n x)|$ is bounded by 1, $\lim_{n\to\infty} \sin(2 \pi n x) / (2 \pi n) = 0$. We therefore have 
pointwise convergence of 
\begin{align}
F_{X_n}(x) 
	\to F(x) := x \Ib_{0 \leq x \leq 1} + \Ib_{1 < x}.
\end{align} 
The function $F$ is the CDF of $U \sim \Uc ((0,1))$ since $F(x) = \int_{-\infty}^x dy\ \Ib_{[0,1]}$. However, $f_{X_n} 
\not\to f_U$ pointwise. To see this, choose $\alpha = 1 / (2 \pi)$, then $\lim_{n\to\infty} f_{X_n}(\alpha)$ does not 
exist.

\subsection*{Solution to Exercise \ref{exer:3.7}}
Let $X_i$ be i.i.d. with $X_i \sim X \sim \text{Geo}(p)$ for $i = 1,2,3,\cdots$. Then we can write
\begin{align}
N
	= \sum_{i=1}^k X_i.
\end{align}
Then the characteristic function of $2 N p$ is
\begin{align}
\phi_{2 N p} (t)
	&= \phi_{2X_1 p + \cdots + 2 X_k p}(t) \\
	&= \phi_{2Xp}^k (t ) \\
	&= \phi_{X}^k ( 2 t p ) \\
	&= G_{X}^k\( e^{i 2 p t} \) \\
	&= \( \frac{p}{1 - (1 - p) e^{2 i t p} }\)^k \\
	&= \( \frac{p}{( 1 - 2 i t)p + O(p^2)} \)^k \\
	&\to \( \frac{ 1}{1 -  2 i t}\)^k \label{eqn.gammachar}
\end{align}
as $p \to 0$. The expression \eqref{eqn.gammachar} is the characteristic function of $Z$ where $Z \sim \Gamma\(k, 
\frac{1}{2}\)$.  Therefore, by the continuity theorem, distribution of $2 N p $ tends to the distribution function of 
$Z$ where $Z \sim \Gamma\( k,\frac{1}{2} \)$.  

\subsection*{Solution to Exercise \ref{exer:3.8}}
Let $X_j$ be i.i.d Bernoulli trials with $p = 2^{-1}$. Then $\Eb[ X_j] = 2^{-1}$ and $\Vb [ X_j] = 2^{-2}$. Define
\begin{align}
&S_n
	:= \sum_{j=1}^n X_j,
&& U_n
	:= \frac{ S_n - n 2^{-1} }{\sqrt{ n 2^{-2}}}.
\end{align}
Given $x$, suppose we would like to compute $\Pb ( -x \leq U_n \leq x)$. First, recall that $A_n = \{ k : |k - n2^{-1}| 
\leq x \sqrt{n 2^{-2}}\}$. The $k \in A_n$ correspond to the values of $S_n$ such that $-x \leq U_n \leq x$. Since $S_n 
\sim \text{Bin}(n,p =2^{-1})$, computing $\Pb ( -x \leq U_n \leq x)$ is done by computing the sum 
\begin{align}
\sum_{ k \in A_n} \binom{n}{k} \( \frac{1}{2}\)^{k} \(1- \frac{1}{2}\)^{n - k}
	= 2^{-n} \sum_{k \in A_n} \binom{n}{k}.
\end{align}
By the CLT,  $U_n$ converges in distribution to $Z$ such that $Z \sim N(0,1)$ as $n \to 
\infty$. Therefore,  as $n\to \infty$,
\begin{align}
2^{-n} \sum_{k \in A_n} \binom{n}{k}
	= \Pb( -x \leq U_n \leq x)
	\to 
	\int_{-x}^x dy\ \frac{1}{\sqrt{2\pi}} e^{-y^2 / 2}
\end{align}

\subsection*{Solution to Exercise \ref{exer:3.9}}
Using the properties for linear combinations of characteristic functions of independent variables and using the fact that for a variable $X_j\sim N(0,1)$, $\phi_{X_j}(t)=\ee^{-\frac{1}{2}t^2}$:
\begin{align}
\phi_{Y_n}(t)
	&=	\prod_{j=1}^n\phi_{\frac{\sqrt{j}}{n}X_j}(t)=\prod_{j=1}^n\phi_{X_j}(\frac{\sqrt{j}}{n}t)
	=		\prod_{j=1}^n \ee^{-\frac{1}{2}(\frac{\sqrt{j}}{n}t)^2}
	=		\ee^{-\frac{t^2}{2n^2}\sum_{j=1}^n j} .
\end{align}
Then, noticing that:
\begin{align}
\lim_{n\rightarrow \infty}\frac{1}{n^2}\sum_{j=1}^n j
	&=	\lim_{n\rightarrow \infty}\frac{1}{n^2}\frac{n^2+n}{2}=\frac{1}{2} ,
\end{align}
it follows that 
\begin{align}
\lim_{n\rightarrow \infty}\phi_{Y_n}(t)
	&=	\ee^{-\frac{t^2}{4}}=\phi_Y(t), &
Y
	&\sim N(0,\frac{1}{2}) .
\end{align}
Then, by the Continuity Theorem \ref{thm:continuity}, we have
\begin{align}
Y_n
	&\stackrel{\Dc}{\rightarrow} Y \sim N(0,1/2) .
\end{align}

\subsection*{Solution to Exercise \ref{exer:3.10}}
(a) We have
\begin{align}
\phi_Z(t)
	&=	\Eb \ee^{\ii t Z}
	=		\sum_{k=0}^\infty (\ee^{\ii t k}) \frac{\lam^k \ee^{-\lam}}{k!}
	=		\ee^{-\lam} \sum_{k=0}^\infty \frac{(\lam \ee^{\ii t})^k }{k!}
	=		\exp \Big( \lam(\ee^{\ii t} - 1 ) \Big) , \\
\phi_X(t)
	&=	\Eb \ee^{\ii t X}
	=		p \ee^{\ii t} + (1-p) .
\end{align}
(b)
We have
\begin{align}
\phi_{S_n}(t)
	&:=	\Eb \ee^{\ii t S_n}
	=		\Eb \ee^{\ii t \sum_{j=1}^n X_j^{(n)}}
	=		\prod_{j=1}^n \Eb \ee^{\ii t X_j^{(n)}}
	=		\prod_{j=1}^n \Big( 1 - \lam/n + \lam \ee^{\ii t} / n \Big)
	=		\Big( 1 + \lam ( \ee^{\ii t} -1 ) / n \Big)^n .
\end{align}
Noting that $\phi_{S_n}(t) \to \phi_Z(t)$ as $n \to \infty$ and that $\phi_Z$ is continuous at $t=0$, it follows by the Continuity Theorem \ref{thm:continuity} that $S_n \stackrel{\Dc}{\to} Z$.


\section{Exercises from Section \ref{sec:exercises-4}}

\subsection*{Solution to Exercise \ref{exer:4.1}}
\begin{enumerate}
	\item[(a.)] $X$ takes on finitely many states, and the transition density of $X_n$ to $X_{n+1}$ depends only on 
the current value $X_n$. Thus, $X$ is a Markov chain. The one step transition matrix is given by
\begin{align}
P = 
\begin{pmatrix}
1/6 & 1/6 & 1/6 & 1/6 & 1/6 & 1/6 \\
   & 2/6 & 1/6 & 1/6 & 1/6 & 1/6 \\
   &    & 3/6 & 1/6 & 1/6 & 1/6 \\
   &    &    & 4/6 & 1/6 & 1/6 \\
   &    &    &    & 5/6 & 1/6 \\
   &    &    &    &    & 6/6 
\end{pmatrix}.
\end{align}
\item[(b.)] The process $X_n$ can be given by 
	\begin{align}
	X_{n} 
		= \left\{
		\begin{array}{lr}
		X_{n-1}& n\text{th roll is not 6} \\
		X_{n-1} +1 & n\text{th roll is 6}
		\end{array}
		\right.
	\end{align}
	Hence, $X_n$ is a Markov chain. The transition matrix is given by
	\begin{align}
	P =
	\begin{pmatrix}
	5/6 & 1/6 \\
	    & 5/6 & 1/6 \\
	    &     & 5/6 & 1/6 \\
	    &     &     &  \ddots
	\end{pmatrix}.
	\end{align}
\item[(c.)] Given that $X_{n} = i$, $X_{n+1}$ can take either the value $j = i+1$ or $j= 0$. Having additional 
information about past states of $X$ changes neither the possible states $X_{n+1}$ nor the probability of reaching each 
state. Therefore, $\Pb( X_{n+1} = j | \Fc_n) = \Pb( X_{n+1} | X_n)$. $X$ is a Markov chain. The transition matrix is 
given by
\begin{align}
P = 
\begin{pmatrix}
1/6 & 5/6 \\
1/6 &  & 5/6 \\
1/6 &  &  & 5/6 \\
1/6 &  &  &  & 5/6 \\
\vdots & & & & & \ddots
\end{pmatrix}.
\end{align}
\item[(d.)] Let $Y$ be a geometrically distributed random variable with parameter $1/6$ independent of $X_n$ 
for all $n$. Note that $X_n$ can be written as 
	\begin{align}
	X_{n}
		= \left\{
		\begin{array}{lr}
		X_{n-1} - 1 & X_{n-1} > 0 \\
		Y & X_{n-1} = 0
		\end{array}
		\right. .
	\end{align}
	Since the probability distribution of $X_n$ depends only on $X_{n-1}$, $X$ is a Markov chain. The transition 
matrix is given by
\begin{align}
P = 
\begin{pmatrix}
\tfrac{1}{6} & (\tfrac{5}{6})\tfrac{1}{6} & (\tfrac{5}{6})^2\tfrac{1}{6} & (\tfrac{5}{6})^3\tfrac{1}{6} & \hdots \\
1 & 0\\
 & 1 & 0\\
 & & 1 & 0 \\
 & & & 1 \\
 & & & & \ddots
\end{pmatrix}.
\end{align}
\end{enumerate}

\subsection*{Solution to Exercise \ref{exer:4.2}}
\begin{enumerate}
	\item[(a.)] Note that for $X_n$, the one step transition matrix is
	\begin{align}
	p(i,j)
		= \left\{
		\begin{array}{lr}
		p,& j = i+1\\
		q,& j = i-1\\
		0,& \text{else}
		\end{array}
		\right. .
	\end{align}
	In particular, the above expression holds for $n= 1$. So, computing the matrix $P^2$ will give the one 
step transition matrix $\Pi = \pi(i,j)$ for $Y_n$, since $Y_n$ is a Markov chain. Hence, the one step transition matrix 
$\pi_n$ for $Y_n = X_{2n}$ is 
\begin{align}
\pi(i,j)
	= \left\{ 
	\begin{array}{lr}
	p^2, & j = i+2 \\
	2pq, & j = i \\
	q^2, & j = i-2 \\
	0,& \text{else}
	\end{array}
	\right. .
\end{align}
	\item[(b.)] If $G$ is the generating function for an individual, and $G_n(s)$ is the function $G$ composed with 
itself $n$ times, then transition matrix for $X_n$ is given by
	\begin{align}
	&p(i,j)
		= \[\frac{1}{j!} \frac{d^j}{ds^j} G^i(s)\]_{s=0}
	\end{align}
	i.e. the $j$th coefficient of $G^{i}(s)$. Hence, the transition matrix $\pi(i,j)$ for $Y$ is 
	\begin{align}
	\pi(i,j)
		= \[\frac{1}{j!} \frac{d^j}{ds^j} G_2^i(s)\]_{s=0}.
	\end{align}
\end{enumerate}

\subsection*{Solution to Exercise \ref{exer:4.3}}
Choose $s \in S$. Since $s$ communicates with $k$, there exists a smallest integer $n_s \geq 1$ for which $p_{n_s}(s,k) 
>0$. Note that the event  $\{ X_0 =s $ and $X_{n_s} = k\}$ is a subset of $\{X_0 =s$ and then $X$ never returns to 
$s\}$. Therefore, 
\begin{align}
0 < 
	\Pb( \{X_0 =s \text{ and } X_{n_s} = k \}) 
	\leq \Pb( X_0 = s \text{ and } X \text{ never returns to } s\} ).
\end{align}
Thus, $s$ is a transient state. 

\subsection*{Solution to Exercise \ref{exer:4.4}}
Let $V$ be the number of times $X$ visits $j$ before $i$ given that $X_0 = i$, and let 
\begin{align}
p 
	:= \Pb( \tau_j < \tau_i | X_0=i) 
	= \Pb( \tau_i < \tau_j | X_0=j).
\end{align}
It is clear that $\Pb( V \geq 1 | X_0 = i) = p$. Additionally,
\begin{align}
\Pb( V \geq 2 | X_0 = i)
	&= \Pb( \tau_j < \tau_i | X_0 = i) \Pb (\tau_j < \tau_i |X_0 = j) \\
	&= \Pb( \tau_j < \tau_i | X_0 = i) (1 - \Pb( \tau_i < \tau_j | X_0 =j)) \\
	&= p (1-p).
\end{align}
Following this process, we see that 
\begin{align}
\Pb(V \geq k | X_0 = i)
	= p ( 1- p)^{k-1}. 
\end{align}
Hence, 
\begin{align}
\Eb[ V | X_0 = i]
	= \sum_{k=1}^\infty p ( 1- p)^{k-1}
	= p \sum_{k=0}^\infty (1 - p)^k 
	= 1. 
\end{align}

\subsection*{Solution to Exercise \ref{exer:4.5}}
Note that for 
\begin{align}
&V 
	= \left(
\begin{array}{ccc}
 1 & 1 & -1 \\
 1 & -1 & 0 \\
 1 & 1 & 1 \\
\end{array}
\right)
&&\text{and}
&& D 
	= 
\left(
\begin{array}{ccc}
 1 & 0 & 0 \\
 0 & 1-4 p & 0 \\
 0 & 0 & 1-2 p \\
\end{array}
\right),
\end{align}
we have that $\Pv = V D V^{-1}$. Thus, 
\begin{align}
\Pv^n
	&= (V D V^{-1})^n
	= V D^n V^{-1} \\
	&= 
\left(
\begin{array}{ccc}
 1 & 1 & -1 \\
 1 & -1 & 0 \\
 1 & 1 & 1 \\
\end{array}
\right)
\left(
\begin{array}{ccc}
 1 & 0 & 0 \\
 0 & (1-4 p)^n & 0 \\
 0 & 0 & (1-2 p)^n \\
\end{array}
\right)
\left(
\begin{array}{ccc}
 \frac{1}{4} & \frac{1}{2} & \frac{1}{4} \\
 \frac{1}{4} & -\frac{1}{2} & \frac{1}{4} \\
 -\frac{1}{2} & 0 & \frac{1}{2} \\
\end{array}
\right) \\
	&=
\left(
\begin{array}{ccc}
 \frac{1}{4} \left(2 (1-2 p)^n+(1-4 p)^n+1\right) & \frac{1}{2} \left(1-(1-4 p)^n\right) & \frac{1}{4} \left(-2 (1-2 
p)^n+(1-4 p)^n+1\right) \\
 \frac{1}{4} \left(1-(1-4 p)^n\right) & \frac{1}{2} \left((1-4 p)^n+1\right) & \frac{1}{4} \left(1-(1-4 p)^n\right) \\
 \frac{1}{4} \left(-2 (1-2 p)^n+(1-4 p)^n+1\right) & \frac{1}{2} \left(1-(1-4 p)^n\right) & \frac{1}{4} \left(2 (1-2 
p)^n+(1-4 p)^n+1\right) \\
\end{array}
\right) .
\end{align}
Furthermore, observe that the first row of $V^{-1}$ is the unique left eigenvector corresponding to the non-repeated 
eigenvalue 1. This is the unique stationary distribution i.e. $\piv = (1/4,1/2,1/4)$. By theorem 4.27,
\begin{align}
&\bar \tau_1 	
	= 4,
&&\bar \tau_2
	= 2,
&&\bar \tau_3
	= 4.
\end{align}

\subsection*{Solution to Exercise \ref{exer:4.6}}
Let $S_n \sim \text{Bin}(X_n, p)$ and $Y_n$ be i.i.d with $Y_n \sim Y \sim \text{Poisson}(\lambda)$. Furthermore, 
suppose that $S_n \perp Y_n$. Then we can write $X_{n+1} = S_n + Y_n$. Let $G_n(s)$ be the probability generating 
function of $X_n$, and let $G_Y$ be the probability generating function of $Y$. Thus,
\begin{align}
G_{n+1}(s)
	= G_{S_n + Y_n}(s) 
	= G_{Y}(s) G_{S_n}(s) 
	= G_Y(s) \Eb[ s^{S_n} ] 
	&= G_Y(s) \Eb[ \Eb[ s^{S_n} | X_{n} ] ] \\
	&= G_Y(s) \Eb[ (1 - (1-p)(1-s))^{X_n}] \\
	&= G_Y(s) G_n(1 - (1-p)(1-s)).
\end{align}
Note that
\begin{align}
&G_Y(s) G_n(1 - (1-p)(1-s))\\
	= &G_Y(s) G_Y(1 - (1-p)(1-s)) G_{n-1}(1 - (1-p)^2(1-s)) \\
	= &G_Y(s) G_Y(1 - (1-p)(1-s)) G_{Y}(1 - (1-p)^2(1-s)) G_{n-2}(1 - (1-p)^3(1-s)) \\
	& \hspace{2cm} \vdots \\
	= &\prod_{k=0}^{n+1} G_Y(1 - (1-p)^k (1-s))
\end{align}
It is known that $G_Y(s) = \exp \left( \lambda ( s -1 ) \right)$. Thus, 
\begin{align}
&\lim_{n\to\infty} \prod_{k=0}^{n+1} G_Y(1 - (1-p)^k (1-s)) \\
	= &\lim_{n\to\infty} \prod_{k=0}^{n+1} \exp\(\lambda(s-1) (1-p)^k\) \\ 
	= &\exp \( \lambda (s-1) \sum_{k=0}^\infty (1-p)^k \) \\
	= &\exp \( \frac{\lambda}{p} (s-1) \).
\end{align}
Thus, the stationary distribution of $X_n$ as $n\to\infty$ is Poisson with parameter $\lambda / p$.

\subsection*{Solution to Exercise \ref{exer:4.7}}
Note that the identity matrix satisfies the criteria of this problem, since every unit vector is a stationary 
distribution. However, the Markov chain consisting of only absorbing states is too boring. Let 
\begin{align}
\Pv
	= 
\left(
\begin{array}{cccc}
 \frac{1}{2} & \frac{1}{2} & 0 & 0 \\
 \frac{1}{2} & \frac{1}{2} & 0 & 0 \\
 0 & 0 & \frac{1}{2} & \frac{1}{2} \\
 0 & 0 & \frac{1}{2} & \frac{1}{2} \\
\end{array}
\right).
\end{align}
Note that $\piv_1 = (1/2,1/2,0,0)$ and $\piv_2 = (0,0,1/2,1/2)$ are both stationary distributions. This Markov 
chain is composed of two closed sets of states $C_1 = \{ 1,2\}$ and $C_2=\{ 3,4\}$.

\subsection*{Solution to Exercise \ref{exer:4.8}}
For the Markov chain to be reversible, we need that $\pi_{i} p(i,j) = \pi_j p(j,i)$ for $i,j = 0,1,2,\cdots,n$. Of 
interest is $\pi_i p(i,i+1) = \pi_{i+1} p(i+1,i)$. This implies that for $i = 1,2,3,\cdots,n-1$,
\begin{align}
\pi_{i+1} 
	= \pi_i \frac{p(i,i+1)}{p(i+1,i)}
	= \pi_{i}\frac{\lambda_i}{\mu_{i+1}}
	= \pi_{i-1} \frac{\lambda_{i-1}\lambda_i}{\mu_i \mu_{i+1}}
	= \cdots
	= \pi_0 \frac{\lambda_0 \lambda_1 \cdots \lambda_{i}}{\mu_1 \mu_2 \cdots \mu_{i+1}}.
\end{align}
We can choose $\pi_0$ as
\begin{align}
\pi_0
	= \frac{1}{1 + \sum_{i = 1}^{n} \frac{\lambda_0 \lambda_1 \cdots \lambda_{i-1}}{\mu_1 \mu_2 \cdots \mu_{i}}}
\end{align}
so that $\pi_0 + \pi_1 + \cdots + \pi_n = 1$. It remains to be checked that $\piv$ gives a stationary distribution. 
The matrix $\Pv$ has the form
\begin{align}
\Pv
	=
\begin{pmatrix}
1 - \lambda_0 & \lambda_0 \\
\mu_1 & 0 & \lambda_1 \\
 & \mu_2 & 0 & \lambda_2 \\
 & & \ddots & \ddots & \ddots \\
 & & & \mu_{n-1} & 0 & \lambda_{n-1} \\
 & & & & \mu_n & 1- \mu_n
\end{pmatrix}.
\end{align}
To this end, observe
\begin{align}
&(\piv \Pv)_0
	= \pi_0 ( 1- \lambda_0) + \pi_1 \mu_1 
	= \pi_0,
&&\text{and}
&& (\piv \Pv)_n
	= \pi_{n-1} \lambda_{n-1} + \pi_n (1 - \mu_n)
	= \pi_n.
\end{align}
For $i = 1,\cdots, n-1$, we have that
\begin{align}
(\piv \Pv)_i
	&= \pi_{i-1} \lambda_{i-1} + \pi_{i+1}\mu_{i+1} \\
	&= \pi_0 \frac{\lambda_0 \cdots \lambda_{i-2} \lambda_{i-1}}{\mu_1 \cdots \mu_{i}}
		+ \pi_0 \frac{\lambda_0 \cdots  \lambda_{i} \mu_{i+1}}{\mu_1 \cdots \mu_{i}\mu_{i+1}} \\
	&= \pi_0 \frac{ \lambda_0 \cdots \lambda_{i-1} ( \lambda_i + \mu_i) }{\mu_1 \cdots \mu_i }\\
	&= \pi_i.
\end{align}


\section{Exercises from Section \ref{sec:exercises-5}}

\subsection*{Solution to Exercise \ref{exer:5.1}}
The number of patients in the hospital at a given time $X = (X_t)_{t\geq 0}$ can be modeled as a birth-death process 
with ``death intensity'' $\mu$ and ``birth intensity'' $\lambda$. The generator of such a process is given in the notes 
as Example 5.14. The generator is thus
\begin{align}
\Gv
	= 
	\begin{pmatrix}
	-\lambda & \lambda & 0 & 0 & 0 & \cdots \\
	\mu & - ( \mu + \lambda) & \lambda & 0 & 0 & \cdots \\
	0 & \mu & - ( \mu + \lambda) & \lambda & 0  & \cdots\\
	0 & 0 & \mu & - ( \mu + \lambda) & \lambda & \cdots \\
	\vdots & \vdots & \vdots & \vdots & \ddots 
	\end{pmatrix}.
\end{align}
A distribution $\piv$ is invariant if and only if $\piv \Gv = 0$. Component-wise, this amounts to 
\begin{align}
0
	&= -\lambda \pi_0 + \mu \pi_1 \\
0
	&= \lambda \pi_{i-1}  - ( \mu + \lambda) \pi_i + \mu \pi_{i+1},
	& i \geq 1.
\end{align}
Using direct substitution, we find that
\begin{align}
&\pi_i
	= \frac{\lambda^i}{\mu^i} \pi_0,
&i \geq 1.
\end{align}
Additionally, 
\begin{align}
&1 
	= \sum_{i=0}^\infty \frac{ \lambda^i}{\mu^i} \pi_0
&&\Rightarrow
&\pi_0 
	= \( \sum_{i =0}^\infty  \frac{ \lambda^i}{\mu^i} \)^{-1}.
\end{align}
The above sum converges if and only if $\lambda < \mu$. In this case, 
\begin{align}
\pi_0
	= 1 - \frac{\lambda}{\mu}.
\end{align}
Intuition tells us that the expected amount number of people in the hospital at steady-state is equal to the expected 
rate of arrival times the expected rate of processing time. This is known as Little's Law. Thus, 
\begin{align}
\frac{\lambda }{\mu - \lambda}
	= \sum_{i=0}^\infty \( \frac{\lambda}{\mu}\)^i \pi_0 i
	= \Eb[ X]
	= \lambda \Eb[\tau].
\end{align}
Thus, 
\begin{align}
\Eb[\tau]
	= \frac{1}{\mu - \lambda}.
\end{align}

\subsection*{Solution to Exercise \ref{exer:5.2}}
The state space of process $Y_n$ is the same as that of $(X_t)_{t\geq 0}$. Now,
\begin{align}
\Pb( Y_{n} = j \mid Y_0 = i)
	= \Pb( X_{\tau_{n}^+} = j \mid X_{0} = i ).
\end{align}
Suppose that $\Qv_t$ is the transition matrix for $X_t$ and $\Pv^{(n)}$ is the transition matrix for $Y_n$. Then 
$\Qv_{\tau_n^+} = \Pv^{(n)}$. Therefore, 
\begin{align}
\piv \Pv^{(n)}
	= \piv \Qv_{\tau_n^+}
	= \piv. 
\end{align}
This holds for all $n \geq 1$. Thus, $Y$ has the same stationary distribution as $X$.

\subsection*{Solution to Exercise \ref{exer:5.3}}
Appealing to the Kolmogorov Forward equation, we have
\begin{align}
\frac{d}{dt} p_t (j,0)
	&= -\lambda p_t(j,0) + \mu p_t(j,1), \\
\frac{d}{dt} p_t(j,i)
	&= \lambda p_t(j,i-1) - (\lambda + i \mu) p_t(j,i) + (i+1) \mu p_t(j,i+1).
\end{align}
Multiplying both sides of the above equations by $s^i$ and summing gives
\begin{align}
&\frac{d}{dt} G_{X_t}(s)
	= \lambda (s - 1) G_{X_t} (s) + \mu( 1 - s) \frac{d}{ds} G_{X_t}(s),
&& G_{X_0}(s) = s^j	
\end{align}
This PDE has the solution
\begin{align}
\Eb[ s^{X_t}  \mid X_0 = j]
	= G_{X_t}(s) 
	= \left((s-1) e^{- \mu t }+1\right)^j e^{\frac{\lambda  (s-1) e^{- \mu t} \left(e^{\mu  t}-1\right)}{\mu }}.
\end{align}
Note that 
\begin{align}
\lim_{t\to\infty} G_{X_t}(s) 
	= e^{\lambda( s- 1) / \mu}.
\end{align}
Therefore, $X_t$ is Poisson distributed with parameter $\lambda \mu^{-1}$ as $t$ tends to infinity. 

\subsection*{Solution to Exercise \ref{exer:5.4}}
Looking at Definition 5.2 and Definition 5.9 in the notes, we see that the generator for $N$ should be
\begin{align}
\Gv(t)
	= 
	\begin{pmatrix}
	-\lambda(t) & \lambda (t) \\
	 & - \lambda(t) & \lambda(t) \\
	 & & - \lambda(t) & \lambda(t) \\
	 & & & \ddots && \ddots
	\end{pmatrix}.
\end{align}
Let $p(t,j;s,i) := \Pb(N(t) = j \mid N(s) = i)$ for $t \geq s$. Then, the Kolmogorov forward equations 
are 
\begin{align}
\d_t p(t,i; s,i)
	&= - \lambda(t) p(t,i;s,i), \\
\d_t p(t,j;s,i)
	&= \lambda(t) p(t,j-1;s,i) - \lambda(t) p(t,j;s,i), 
	&& i < j
\end{align}
and the Kolmogorov backwards equations are
\begin{align}
-\d_s p(t,i;s,i)
	&= -\lambda(s) p(t,i;s,i) \\
-\d_s p(t,j;s,i)
	&= \lambda(s) p(t,j;s,i+1) - \lambda(s) p(t,j;s,i)
	&& i<j \\
	&= \lambda(s) p(t,j-1;s,i) - \lambda(s) p(t,j;s,i),
\end{align}
where the fact that $p(t,j;s,i+1) = p(t,j-1;s,i)$ was used. Let 
\begin{align}
G(t,x;s,i)
	:= \sum_{j=0}^\infty x^j p(t,j;s,i)
	= \sum_{j=i}^\infty x^j p(t,j;s,i)
\end{align}
be the generating function for $N$. Multiplying by both sides of the Backward and Forward Kolmogorov equations by $x^i$ 
and summing gives
\begin{align}
&\text{K. Forward:}
& \d_t G(t,x;s,i)
	&= \lambda(t) ( x-1) G(t,x;s,i),
&& G(s,x;s,i) 
	= x^i, \\
&\text{K. Backward:}
& -\d_s G(t,x;s,i)
	&= \lambda(s) (x - 1) G(t,x;s,i),
&& G(t,x;t,i)
	= x^i. 
\end{align}
These equations have the solutions
\begin{align}
&\text{K. Forward:}
& G(t,x;s,i)
	&= x^i e^{(x - 1) \int_s^t du \lambda (u) } \\
&\text{K. Backward:}
& G(t,x;s,i)
	&= x^i e^{-(x - 1) \int_s^t du \lambda (u) }.
\end{align}
We have that 
\begin{align}
\Pb (\tau_1 > t)
	= p(t,0;0,0)
	= G(t,0;0,0)
	= e^{- \int_0^t du \lambda(u)}.
\end{align}
Thus,
\begin{align}
f_{\tau_1}(t) 
	= \frac{d}{dt} e^{- \int_0^t du \lambda(u)}
	= \lambda(t) e^{- \int_0^t du \lambda(u)}.
\end{align}
Therefore, for $\lambda(t) = c(1 + t)^{-1}$, we have that
\begin{align}
\Eb[ \tau_1]
	&= \int_0^\infty dv\ v \frac{c}{1+t} e^{-\int_0^v du\ c (1 + u)^{-1}} \\
	&= \int_0^\infty dv\ \frac{c v}{(1 + v)^{c+1}}.
\end{align}
The above integral is finite if and only $c > 1$. 

\subsection*{Solution to Exercise \ref{exer:5.5}}
Note that
\begin{align}
\Eb[ s^{N_t} ]
	= \Eb[ \Eb[ s^{N_t} \mid \Lambda] ]
	= \Eb[ e^{ t \Lambda( s -1 )} ]
	= p e^{\lambda_1 t (s -1 )} + (1 - p ) e^{\lambda_2 t (s - 1) }.
\end{align}
So,
\begin{align}
\Eb[ N_t]
	= \frac{d}{ds} G(1)
	= p t \lambda_1 + (1 - p) t \lambda_2,
\end{align}
and
\begin{align}
\Vb [ N_t] 
	&= \frac{d^2}{dt^2} G(1) + \frac{d}{ds} G(1) - \[ \frac{d}{ds} G(1) \]^2 \\
	&= \lambda _1^2 p t^2+\lambda _2^2 (1-p) t^2+\lambda _1 p t-\left(\lambda _1 p t+\lambda _2 (1-p) 
t\right){}^2+\lambda _2 (1-p) t.
\end{align}

\subsection*{Solution to Exercise \ref{exer:5.6}}

Proof.\\
(i) $\Rightarrow$ (ii):\\
Using (i) we have
\begin{align}
& \ln\left( \frac{g(s_0, s_1)}{g(s_1, s_0)} \right) + \ln\left( \frac{g(s_2, s_1)}{g(s_1, s_2)} \right) + \cdots + \ln\left( \frac{g(s_{n-1}, s_n)}{g(s_n,  s_{n-1})} \right) = \ln \pi(s_n) - \ln \pi(s_0)\\
&= \ln\left( \frac{\pi(s_0) g(s_0, s_1)}{\pi(s_1) g(s_1 , i_0)} \right) + \ln\left( \frac{ \pi(s_1)  g(s_1, s_2)}{\pi(s_2) g(s_2, s_1)} \right) + \cdots + \ln\left( \frac{\pi(s_{n-1}) g(s_{n-1}, s_n)}{\pi(s_{n})g(s_n, s_{n-1})} \right) +\ln \pi(s_n) - \ln \pi(s_0) \label{p-id-w}\\
&= \ln\pi(s_n) - \ln\pi(s_0)
\end{align}
This means the term in \eqref{p-id-w} is independent of the path; it is completely determined by the initial state $i \equiv s_0$ and final state $j \equiv s_n$. \\
Physical interpretation of statement (ii): it is similar to the path-independent work in a conservative force field, which implies the existence of a potential function. In this discrete case, \eqref{p-id-w}  indicates that the potential function is the minus logarithm of the stationary distribution. Equilibrium probability distribution is the exponential function of a potential energy is known as the Boltzmann's law in statistical mechanics.\\
(ii) $\Rightarrow$ (iii):\\
Equation \eqref{p-id-w} is equivalent to 
$$ \pi(i) g(i,j ) =   \pi(j) g(j,i)$$
which implies 
$$ \pi(i) p_t (i,j)  = \pi(j) p_t (j,i) $$
That is, the process $X$ is time reversible.\\
(iii) $\Rightarrow$ (i):\\
Reversible Markov process has joint distribution $\pi(i) p_t (i,j)  = \pi(j) p_t (j,i)$ for any $t$.
Therefore,
\begin{align}
\pi(i) \left(e^{t \Gv}  \right)_{ij} &= \pi(j) \left(e^{t \Gv}  \right)_{ji} \\
\pi(i) \left( \delta_{ij} + tg(i,j) + \frac{t^2}{2} ({\Gv}^2)_{ij} \cdots  \right)_{ij} &= \pi(j) \left( \delta_{ji} + t g(j,i) + \frac{t^2}{2} ({\Gv}^2)_{ji} \cdots \right)_{ji} 
\end{align}
When $i \neq j$,
\begin{align}
\pi(i) \left( t g(i,j) + \frac{t^2}{2} ({\Gv}^2)_{ij} \cdots  \right)_{ij} &= \pi(j) \left( t g(j,i) + \frac{t^2}{2} ({\Gv}^2)_{ji} \cdots \right)_{ji} 
\end{align}
Letting $t=0$ we have $\pi(i)g(i,j) = \pi(j)g(j,i)$\\
(ii) $\Rightarrow$ (iv):\\
We simply choose $i=j$, then statement (ii) becomes statement (iv).\\
(iv) $\Rightarrow$ (i):\\
This is the $(\Leftarrow)$ direction of Theorem \ref{thm:kolmogorov.cycle.cond}.\\
(i) $\Leftrightarrow$ (v):\\
Simply choose the positive diagonal matrix $\Pi = diag(\pi(1), \pi(2), \cdots)$. Conversely, if there exists a set of positive number $\nu(i)$ such that $\nu(i)g(i,j) = \nu(j)g(j,i)$, then it is an eigenvector of $\Gv$ with eigenvalue $0$. Therefore since $\Gv$ has a unique eigenvalue $0$, the normalized
$$ \frac{\nu(i)}{\sum_{k\in S} \nu(k)} = \pi(i) $$
(i) $\Leftrightarrow$ (vi):\\
The stationary entropy production rate, $e^*_p $, is defined by:
$$e^*_p = \frac{1}{2} \sum_{i, j} \left(\pi(i) g(ij)- \pi(j)g(j, i)\right) 
\ln\left(\frac{\pi(i) g(i,j)}{\pi(j)g(j, i)}\right)    $$
in which every term inside the summation $ \left( \pi(i) g(i,j) - \pi(j) g(j,i) \right) \ln \left( \frac{\pi(i) g(i, j)}{\pi(j)g(j, i)} \right)  \geq 0 $. Therefore, the entire sum is zero if and only if $\pi(i) g(i,j) - \pi(j) g(j, i) =0$. This is detailed balance. 



\section{Exercises from Section \ref{sec:exercises-6}}

\subsection*{Solution to Exercise \ref{exer:6.1}}
Since $Y_n = \frac{1}{n}X_n$, for any positive real number, we can write 
\begin{align}
 P(Y_n \leq y) &= P(X_n \leq ny)\\
                      &= 1 - (1 - \frac{\lambda}{n})^{\lfloor ny \rfloor}
\end{align}
where $\lfloor ny \rfloor$ is the largest integer less than or equal to $ny$. We then write 
\begin{align}
 \lim_{n\rightarrow \infty} F_{Y_n}(y) &= \lim_{n \rightarrow \infty}1 - (1 - \frac{\lambda}{n})^{\lfloor ny \rfloor}\\
    &= 1-  \lim_{n\rightarrow \infty} (1 - \frac{\lambda}{n})^{\lfloor ny \rfloor}\\
    &= 1 - e^{- \lambda y}
\end{align}
The last equality holds because $ny - 1 \leq  \lfloor ny \rfloor \leq ny$

\subsection*{Solution to Exercise \ref{exer:6.2}}
For any $\om \in (0,1]$, we have 
\begin{align}
\lim_{n\rightarrow \infty} X_n( \omega) &= \lim_{n \rightarrow \infty}\left[ \frac{n}{n+1}\omega + (1-\omega)^n\right]\\
     &= \omega= X(\omega)
\end{align} 
However, if $\omega=0$, then 
\begin{align}
\lim_{n\rightarrow \infty} X_n(0) &= \lim_{n \rightarrow \infty}\left[ \frac{n}{n+1}0 + 1^n\right]\\
     &= 1 \neq X(1)
\end{align} 
Thus, we conclude 
$$ \lim_{n\rightarrow \infty} X_n(\omega) = X(\omega) \quad \forall \omega\in (0.1]$$
Since $P((0,1]) = 1$, we conclude $X_n \stackrel{a.s.}{\rightarrow}X$

\subsection*{Solution to Exercise \ref{exer:6.3}}
Define the sets $A$ and $B$ as follows:
\begin{align}
 A &= \{ \omega \in \Omega : \lim_{n \rightarrow \infty}X_n(\omega) = X(\omega)  \}\\
 B &= \{ \omega\in \Omega : \lim_{n \rightarrow \infty}Y_n(\omega) = Y(\omega) \}
\end{align}
By definition of almost sure convergence, we conclude $P(A) = P(B) = 1$. Therefore, $P(A^c) = P(B^c) = 0$. We conclude
\begin{align}
 P(A \cap B) &= 1 - P(A^c \cup B^c) \\
                   &\geq 1 - P(A^c) - P(B^c)\\
                   &= 1 
\end{align}
Now, consider the sequence $\{ Z_n, n=1,2,... \}$, where $Z_n = X_n + Y_n$, and define the set $C$ as
$$C = \{  \omega \in \Omega : \lim_{n\rightarrow \infty} Z_n(\omega) = X(\omega) + Y(\omega)   \} $$
We claim $A \cap B \subset C$. Specifically, if $\omega \in A \cap B$, then we have 
$$ \lim X_n(\omega) = X  \quad  \lim Y_n(\omega) =Y $$
Therefore, 
\begin{align}
 \lim Z_n(\omega) &= \lim [ X_n(\omega) + Y_n(\omega) ] \\
                  &= \lim X_n(\omega) + \lim Y_n(\omega) \\
                  &= X(\omega) + Y(\omega)
\end{align}
Thus, $\omega\in C$. Thus
$$ P(C) \geq P(A \cap B) =1 $$
which implies $P(C) = 1$. This means that $Z_n  \stackrel{a.s.}{\rightarrow} X+Y $

\subsection*{Solution to Exercise \ref{exer:6.4}}
For $n \in \Nb$, define the following events 
\begin{align}
 A_n &=\{ |X_n - X| < \epsilon/2   \}\\
 B_n &=\{ |Y_n - Y| < \epsilon/2   \}
\end{align}
Since $X_n \rightarrow_p X$ and $Y_n \rightarrow_{p} Y$, we have for all $\epsilon > 0$
\begin{align}
 \lim P(A_n) &= 1\\
 \lim P(B_n) &= 1
\end{align}
We also can write 
\begin{align}
 P(A_n \cap B_n) &= P(A_n) + P(B_n) - P(A_n \cup B_n)\\
                          & \geq P(A_n) + P(B_n) - 1
\end{align}
Therefore,
 $$\lim P(A_n \cap B_n) = 1 $$
 Now, let us define the events $C_n$ and $D_n$ as follows:
 \begin{align}
 C_n &= \{  |X_n - X| + |Y_n - Y| < \epsilon  \}\\
 D_n &= \{  |X_n + Y_n - X  - Y| < \epsilon   \} 
  \end{align}
  Now, note that $(A_n \cap B_n) \subset C_n$, thus $P(A_n \cap B_n) \leq P(C_n)$. Also, by the triangle inequality for absolute values, we have 
  $$ 	|X_n + Y_n - X  - Y| < |X_n - X| + |Y_n - Y|   $$
Therefore, $C_n \subset D_n$, which implies
$$P(C_n) \leq P(D_n)$$
We conclude 
$$P(A_n \cap B_n) \leq P(C_n) \leq P(D_n) $$
Since $\lim P(A_n \cap B_n) =1$, we conclude $\lim P(D_n)=1$. This by definition means that $X_n + Y_n  \stackrel{\Pb}{\rightarrow} X+Y$ 

\subsection*{Solution to Exercise \ref{exer:6.5}}
Define $F_n(x) = P(|X_n| \leq x)$, the distribution function of $|X_n|$. Set $b_n = F_n^{-1}(1 - n^{-2})$ for $n=1,2, . . .,$ and let ${a_n}$ be any sequence with $a_n \rightarrow \infty$. Let $c_n = a_nb_n$. Then, for any $\epsilon > 0$ we have $\epsilon a_n \geq 1$ for $n \geq N_{\epsilon}$ and hence, using $F_n \circ F_n^{-1}(t) \geq t$ for all $0 < t< 1$, 
\begin{align}
P(|X_n| > \epsilon c_n) 
    &= P(|X_n| > \epsilon a_n b_n)\\
	&\leq P(|X_n| > b_n)\\
	&= 1 - F_n(b_n) = 1 - F_n \circ F_n^{-1}(1 - n ^{-2})\\
	&\leq n^{-2}, n \geq N_{\epsilon}
\end{align}
Hence by the first Borel - Cantelli lemma, $P(|X_n| > \epsilon c_n \ i.o.) = 0$ for every $\epsilon > 0 $, that is , $X_n/c_n \stackrel{a.s.}{\rightarrow} 0 $ 

\subsection*{Solution to Exercise \ref{exer:6.6}}
Here, for any $0 < \epsilon < 1, \ P(X_n > \epsilon) = P(x_n = 1)=p_n$.\\
Thus $X_n \stackrel{\Pb}{\rightarrow} 0$ if and only if $p_n \rightarrow 0$.\\
By the first and second Borel-Cantelli lemmas and the first part, for any $0 < \epsilon < 1$ we have\\
$$P(|X_n| > \epsilon \ i.o.) = 
   \left \{
    \begin{array}{c l}
       0\\
       1
     \end{array}
   \right \}
\text{ according as }
\sum_{n=1}^{\infty} p_n
\left \{
  \begin{array}{c l}
  < \infty\\
  = \infty
  \end{array}
\right \}
$$\\
It follows that $X_n \stackrel{a.s.}{\rightarrow}0$ if and only if $\sum_{n=1}^{\infty} p_n< \infty$

\subsection*{Solution to Exercise \ref{exer:6.7}}
Let $c > 0$. Then
\begin{align}
  P(\log X_n > c\log n) = P(X_n > n^c) = n^{-5c}
\end{align}
Hence by the first and second Borel-Cantelli lemmas\\
$$  P(\log X_n > c\log n \quad i.o.) =
\left \{
  \begin{array}{c }
  0 \quad \text{if} \ c > 1/5\\
  1 \quad \text{if} \ c\leq 1/5
  \end{array}
\right \}
$$\\
Hence limsup$_{n\rightarrow \infty} (\log X_n)/ \log n = 1/5 $ \quad almost surely


\section{Exercises from Section \ref{sec:exercises-7}}

\subsection*{Solution to Exercise \ref{exer:7.1}}
\begin{align}
\Eb[(W_t^2 - t) | \Fc_s ]
     &= \Eb[W_t^2 | \Fc_s ] - t\\
     &= \Eb[(W_t - W_s + W_s)^2 | \Fc_s ] - t\\
     &= \Eb[(W_t - W_s)^2 | \Fc_s] + 2\Eb[(W_t - W_s)W_s | \Fc_s] + \Eb[W_s^2 | \Fc_s] - t\\
     &= \Eb[(W_t - W_s)^2] + 2W_s \Eb[(W_t - W_s)| \Fc_s] + W_s^2 - t\\
     &= (t-s) + 0 - t + W_s^2\\
     &= W_s^2 - s
\end{align}

\subsection*{Solution to Exercise \ref{exer:7.2}}
\begin{align}
\phi_{W(N(t))}(u) &= \Eb{e^{i u W(N(t))}}\\
     &=  \Eb\Eb[{e^{i u W(N(t))} | N(t)}]\\
     &= \Eb[e^{-\frac{1}{2} u^2 N(t) }]\\
			&=	\exp \Big( \lam t ( \ee^{-\frac{1}{2} u^2} -1 ) \Big) ,
\end{align}
where we have used
\begin{align}
\Eb s^{N(t)}
	&=	\exp \Big( \lam t (s-1) \Big) .
\end{align}

\subsection*{Solution to Exercise \ref{exer:7.3}}
Suppose that $V(1,T)$ is finite. Then for any $\epsilon > 0$, there exists an $N \geq 1$, $\sum_{j=0}^{n-1} | 
W_{t_{j+1}} - W_{t_j}| < V(1,T) + \epsilon$ for all $n \geq N$. Thus, 
\begin{align}
V_T(2,W)
	&= \lim_{\| \Pi\| \to 0} \sum_{j=0}^{n-1} | W_{t_{j+1}} - W_{t_j}|^2 \\
	&\leq \lim_{\| \Pi\| \to 0}  
	\max_{1\leq k \leq n-1} |W_{t_{k+1}} - W_{t_{k}}| 
	\sum_{j=0}^{n-1} |W_{t_{j+1}} - W_{t_j}| \\
	&\leq \lim_{\| \Pi\| \to 0}  ( V(1,W) + \epsilon) \cdot 
	\max_{1\leq k \leq n-1} |W_{t_{k+1}} - W_{t_{k}}| 
\end{align}
Since $W_t$ is continuous in $t$, as $n\to \infty$, $\max_{1\leq k \leq n-1} |W_{t_{k+1}} - W_{t_{k}}| \to 0$. This 
implies that $V_T(2,W) = 0$, which is a contradiction. Therefore, $V(1,W)$ must be infinite. 
\\[0.5em]
Now, since $V(2,W) = T$, it must be the case that for each $\epsilon > 0$, there exists an $M \geq 1$ such that 
$\sum_{j=0}^{n-1} | W_{t_{j+1}} - W_{t_j}| < T + \epsilon$ for all $n \geq M$. So
\begin{align}
V(3,W)
	&= \lim_{\| \Pi\| \to 0} \sum_{j=0}^{n-1} | W_{t_{j+1}} - W_{t_j}|^3 \\
	&\leq \lim_{\| \Pi\| \to 0}  
	\max_{1\leq k \leq n-1} |W_{t_{k+1}} - W_{t_{k}}| 
	\sum_{j=0}^{n-1} |W_{t_{j+1}} - W_{t_j}|^2 \\
	&\leq  \lim_{\| \Pi\| \to 0} (T + \epsilon) \max_{1\leq k \leq n-1} |W_{t_{k+1}} - W_{t_{k}}|. \\
\end{align}
By the continuity of $W_t$, we see that $V(3,W) = 0$ almost surely. 

\subsection*{Solution to Exercise \ref{exer:7.4}}
First, note that
First, note that
\begin{align}
\Eb( \ee^{ \sig X_T - (\mu \sig + \tfrac{1}{2}\sig^2) T} | \Fc_t )
	&=	\Eb( \ee^{ \sig (\mu T + W_T)  - (\mu \sig + \tfrac{1}{2}\sig^2) T } | \Fc_t ) \\
	&=	\ee^{ \sig (\mu T + W_t)  - (\mu \sig + \tfrac{1}{2}\sig^2) T } \Eb( \ee^{ \sig (W_T - W_t) } | \Fc_t ) \\
	&=	\ee^{ \sig ( \mu T + W_t )  - (\mu \sig + \tfrac{1}{2}\sig^2) T + \tfrac{1}{2}\sig^2 (T-t)} \\
	&=	\ee^{ \sig W_t - \tfrac{1}{2} \sig^2 t  } \\
	&=	\ee^{ \sig X_t - ( \mu \sig + \tfrac{1}{2} \sig^2 ) t  } .
\end{align}
Thus, the process $Z$ defined by $Z_t := \exp ( \sig X_t - ( \mu \sig + \tfrac{1}{2} \sig^2 ) t )$ is a martingale.
Let $Z_t^{(m)} := Z_{t \wedge \tau_m}$ where $\tau_m := \inf\{t \geq 0 : X_t = m \}$.  The process $Z_t^{(m)}$ is also martingale. Thus,
\begin{align}
1
	= \Eb[ Z_0^{(m)}]
	&= \Eb\[ 
	      \exp ( \sigma X_{t \wedge \tau_m} - (\mu \sigma+ \sigma^2 / 2)(t\wedge \tau_m) )
	     \].
\end{align}
Since $\lim_{t\to\infty} t \wedge \tau_m = \tau_m$ and $\lim_{t\to\infty} X_{t \wedge \tau_m} = m$, we have
\begin{align}
1 
	&= \lim_{t\to\infty} \Eb\[ 
	      \exp ( \sigma X_{t \wedge \tau_m} - (\mu \sigma+ \sigma^2 / 2)(t\wedge \tau_m) )
	     \] \\
	&= \Eb\[\lim_{t\to\infty}  
	      \exp ( \sigma X_{t \wedge \tau_m} - (\mu \sigma+ \sigma^2 / 2)(t\wedge \tau_m) )
	     \] \\
	&= \Eb\[ 
	       \exp ( \sigma m - (\mu \sigma+ \sigma^2 / 2)\tau_m )
	      \]
\end{align}
Setting $\sigma = \sqrt{2 \alpha +\mu ^2}-\mu$, it is seen that
\begin{align}
\Eb[ \exp( - \alpha \tau_m)]
	= \exp\( \(\mu - \sqrt{2 \alpha + \mu^2} \)m\).
\end{align}

\subsection*{Solution to Exercise \ref{exer:7.5}}
We have
\begin{align}
| [W,f]_T |
	&=	\Big| \lim_{\norm{\Pi} \to 0} \sum_{j=0}^{n-1} (W_{t_{j+1}} - W_{t_j} ) (f(t_{j+1}) - f(t_j) ) \Big| \\
	&=	\lim_{\norm{\Pi} \to 0} \Big|  \sum_{j=0}^{n-1} (W_{t_{j+1}} - W_{t_j} ) f'(t_j^*) (t_{j+1} - t_j ) \Big| \\
	&\leq \lim_{\norm{\Pi} \to 0} \max_j |W_{t_{j+1}} - W_{t_j}| \sum_{j=0}^{n-1} |f'(t_j^*) | (t_{j+1} - t_j ) \\
	&=	0 \cdot \int_0^T |f'(t)| \dd t
	=	0 .
\end{align}

\subsection*{Solution to Exercise \ref{exer:7.6}}
(a) Fix $0 \leq s \leq t \leq T$.  We have
\begin{align}
\Eb( Z_t | \Fc_s )
	&=	\Eb( \ee^{u W_t - u^2 t / 2} | \Fc_s )
	=		\ee^{u W_s - u^2 t / 2} \Eb( \ee^{u ( W_t - W_s)} | \Fc_s )
	=		\ee^{u W_s - u^2 t / 2} \ee^{u^2 (t-s)/2}
	=		\ee^{u W_s - u^2 s / 2}
	=		Z_s .
\end{align}
(b) Using $\ee^x = 1 + x + x^2 / 2 + x^3/3! + \ldots$ we have
\begin{align}
\ee^{u W_t - u^2 t / 2}
	&=	1 + ( u W_t - u^2 t / 2 ) + (u W_t - u^2 t / 2 )^2 / 2 + (u W_t - u^2 t / 2 )^3 / 3! + \ldots \\
	&=	1 + u W_t + u^2 (W_t^2 - t) / 2 + u^3 (W_t^3 - 3 t W_t) / 3! + \Oc(u^4) \\
	&=	1 + u W_t + u^2 M_t / 2 + u^3 N_t / 3! + \Oc(u^4)
\end{align}
As $Z$ is a martingale for any $u$ it follows that terms of like powers of $u$ must also be martingales.  Thus, we have that $M$ and $N$ are martingales.
\\[0.5em]
More generally, the process $M^{(n)}=(M_t^{(n)})_{0 \leq t \leq T}$, defined by
\begin{align}
M_t^{(n)}
	&:=	\d_u^n \ee^{u W_t - u^2 t / 2} \Big|_{u = 0} ,
\end{align}
is a martingale for every $n \in \Nb_0$.

\subsection*{Solution to Exercise \ref{exer:7.7}}
First, we note that $\Eb \ee^{\ii u W_t} = \ee^{ - u^2 t / 2}$.  Thus, we have
\begin{align}
\Eb( \ee^{\ii u W_t + u^2 t / 2} | \Fc_s )
	&=	\ee^{ \ii u W_s + u^2 t / 2} \Eb( \ee^{\ii u ( W_t - W_s ) } | \Fc_s )
	=		\ee^{ \ii u W_s + u^2 t / 2} \ee^{- u^2 ( t-s )/2 }
	=		\ee^{ \ii u W_s + u^2 s / 2} .
\end{align}
Hence, we see that $Z_t := \ee^{\ii u W_t + u^2 t / 2}$ is a $\Cb$-valued martingale.  Now, note that
\begin{align}
Z_t
	&=	\ee^{ u^2 t / 2} ( \cos u W_t + \ii \sin u W_t ) .
\end{align}
Since $Z$ is a $\Cb$-valued martingale, its real part $M_t = \ee^{ u^2 t / 2} \cos u W_t$ and its imaginary part $N_t = \ee^{ u^2 t / 2} \sin u W_t$ are martingales.

\subsection*{Solution to Exercise \ref{exer:7.8}}
(a) For $m > 0$ we integrate out $w$.  Using the substitution $x = 2m - w$, we obtain
\begin{align}
f_{\Wb_t}(m)
	&=	\int_{-\infty}^m \dd w \, \frac{2(2m-w)}{t \sqrt{2\pi t}} \ee^{-(2m-w)^2/(2t)} 
	=		\int_m^\infty \dd x \, \frac{2x}{t \sqrt{2\pi t}} \ee^{-x^2/(2t)} \\
	&=		\frac{2}{\sqrt{2\pi t}} \Big[ - \ee^{-x^2/(2t)} \Big]_m^\infty
	=		\frac{2}{\sqrt{2\pi t}} \ee^{-m^2/(2t)} .
\end{align}
On the other hand, for $m > 0$ we have
\begin{align}
\Pb( |W_t| \leq m )
	&=	\Pb( -m \leq W_t \leq m )
	=		\int_{-m}^m \dd w \, f_{W_t}(w) 
	=		\int_{-\infty}^m \dd w \, f_{W_t}(w) - \int_{-\infty}^{-m} \dd w \, f_{W_t}(w) . 
\end{align}
Thus, we have
\begin{align}
f_{|W_t|}(m)
	&=	\frac{\dd}{\dd m} \Pb( |W_t| \leq m )
	=		f_{W_t}(m) + f_{W_t}(-m)
	= 	2 f_{W_t}(m)
	=		\frac{2}{\sqrt{2\pi t}} \ee^{-m^2/(2t)} .
\end{align}
As $f_{\Wb_t} = f_{|W_t|}$, the random variables $\Wb_t$ and $|W_t|$ have the same distribution.
\\[0.5em]
(b) Using the density from part (a), we have
\begin{align}
\Eb \Wb_t
	&=	\int_0^\infty \dd m \, m \frac{2}{\sqrt{2\pi t}} \ee^{-m^2/(2t)}
	=		\frac{2 t}{\sqrt{2\pi t}} \Big[ - \ee^{-m^2/(2t)} \Big]_0^\infty
	=		\frac{2t}{\sqrt{2\pi t}}
	=		\sqrt{\frac{2t}{\pi}} .
\end{align}
(c) First, we compute the conditional density of $\Wb_t$ given $W_t = w$.  For $y > 0$ and $y \geq w$ we have
\begin{align}
f_{\Wb_t|W_t}(y,w)
	&=	\frac{f_{W_t,\Wb_t}(w,y)}{f_{W_t}(w)}
	=		\frac{2(2y-w)}{t} \exp \Big( - \frac{(2y-w)^2 - w^2}{2t} \Big)
	=		\frac{2(2y-w)}{t} \ee^{-2y(y-w)/t} .
\end{align}
Noting that $\frac{\dd}{\dd y} \ee^{-2y(y-w)/t} = - \frac{2(2y-w)}{t} \ee^{-2y(y-w)/t}$, we obtain
\begin{align}
\Pb( \Wb_t \geq m | W_t = w )
	&=	\int_m^\infty \dd y \, \frac{2(2y-w)}{t} \ee^{-2y(y-w)/t}
	=		\Big[ - \ee^{-2y(y-w)/t} \Big]_m^\infty
	=		\ee^{-2m(m-w)/t} .
\end{align}
As a sanity check, note that when $w = m$ we obtain $\Pb( \Wb_t \geq m | W_t = m ) = 1$, as it must be.


\section{Exercises from Section \ref{sec:exercises-8}}

\subsection*{Solution to Exercise \ref{exer:8.1}}
Let $f(x) = x^4$ and note that by Theorem 7.6,
\begin{align}
f(W_T) 
	= W_T^4
	= \int_0^T 4 W_t^3\, \dd W_t + \int_0^T 6 W_t^2\, \dd t.
\end{align}
Since integrals with respect to $W_t$ are martingales, 
\begin{align}
\Eb[ W_T^4]
	&= \Eb \int_0^T 6 W_t^2\, \dd t
	= 6 \int_0^T \Eb[ W_t^2]\, \dd t 
	= 6 \int_0^T t\, \dd t
	= 3 T^2.
\end{align}
Let $g(x) = x^6$. Then 
\begin{align}
W_T^6 
	= \int_0^T 6 W_t^5\, \dd W_t + \int_0^T 15 W_t^4\, \dd t.
\end{align}
Thus, 
\begin{align}
\Eb[ W_T^6]
	= \Eb \int_0^T 15 W_t^4\, \dd t
	= 15 \int_0^T \Eb [W_t^4]\, \dd t 
	= 45 \int_0^T t^2\, \dd t  
	= 15 T^3.
\end{align}

\subsection*{Solution to Exercise \ref{exer:8.2}}
Note that 
\begin{align}
\dd (Z_t Y_t)
	&= Z_t\, \dd Y_t + Y_t\, \dd Z_t + \dd Y_t\, \dd Z_t \\
	&= Z_t \( r\, \dd t + \alpha Y_t\, \dd W_t\) + Y_t \( \alpha^2 Z_t\, \dd t - \alpha Z_t\, \dd W_t\) 
		- \alpha^2 Y_t Z_t\, \dd t \\
	&= r Z_t\, \dd t.
\end{align}
Therefore, 
\begin{align}
Y_t 
	&= r Z_t^{-1}\int_0^t Z_s\, \dd s + Y_0 Z_t^{-1} \\
	&= r \int_0^t e^{ \alpha( W_t - W_s) -(1/2)\alpha^2(t-s)}\, \dd s 
		+ Y_0 e^{\alpha W_t - (1/2) \alpha^2 t}.
\end{align}

\subsection*{Solution to Exercise \ref{exer:8.3}}
First, we note that 
\begin{align}
\dd \Delta_t
	= \( \d_{xt} + \frac{1}{2} \sig^2 X_t^2 \d_{xxx} \) f(t,X_t)\, \dd t + \sig X_t \d_{xx} f(t,X_t)\, \dd W_t.
\end{align}
Thus, 
\begin{align}
\dd \Pi_t
	&= \Delta_t\, \dd X_t + X_t\, \dd \Delta_t + \dd \Delta_t\, \dd X_t \\
	&= \( \sig \Delta_t X_t + \sig X_t^2 \d_{xx} f(t,X_t)\)\dd W_t 
		+ \( X_t \( \d_{xt} + \frac{1}{2} \sig^2 X_t^2 \d_{xxx} \) f(t,X_t) + \sig^2 X_t^2 \d_{xx} f(t,X_t) \)\, dt 
\label{eq:pidyn}
\end{align}
Notice that by differentiating both sides of 
\begin{align}
\( \d_t + \frac{1}{2} \sig^2 x^2 \d_{xx} \) f(t,x) 
	= 0,
\end{align}
we get that for all $(t,x)$
\begin{align}
\( \d_{xt} + \frac{1}{2} \sig^2 x^2 \d_{xxx} \) f(t,x) 
	= - \sig^2 x \d_{xx} f(t,x). \label{eq:newpde}
\end{align}
Inserting the right hand side of \eqref{eq:newpde} into \eqref{eq:pidyn} yields 
\begin{align}
\dd \Pi_t
	= \( \sig \Delta_t X_t + \sig X_t^2 \d_{xx} f(t,X_t)\)\dd W_t.
\end{align}
Since integrals with respect to Brownian motions are martingales, the result follows.

\subsection*{Solution to Exercise \ref{exer:8.4}}
First, define the differential operator
\begin{align}
\Ac(t,x) 
	:= \frac{1}{2} \sig^2(t,x) \d_{xx} + \mu(t,x) \d_x.
\end{align}
Note that by Ito's lemma,
\begin{align}
\dd f(t,X_t)
	= \[\d_t + \Ac(t,X_t)\] f(t,X_t)\, \dd t +  \sig(t,X_t) \d_x f(t,X_t) \, \dd W_t.
\end{align}
Furthermore,
\begin{align}
\dd \int_0^t \[ \d_s + \Ac(s,X_s) \] f(s,X_s)\, \dd s
	=  \[ \d_t + \Ac(t,X_t) \] f(t,X_t)\, \dd t.
\end{align}
Thus, 
\begin{align}
\dd M_t^f 
	&= \dd f(t,X_t)- \dd f(0,X_0) - \dd \int_0^t \[ \d_s + \Ac(s,X_s) \] f(s,X_s)\, \dd s \\
	&= \sig(t,X_t) \d_x f(t,X_t) \, \dd W_t. 
\end{align}
This shows that $M_t^f$ is a martingale with respect to $\Fc_t$, the filtration for $W$.

\subsection*{Solution to Exercise \ref{exer:8.5}}
Suppose 
\begin{align}
\frac{\dd \tilde{\Pb}}{\dd \Pb} =  ( -\int_0^T \frac{1}{2} \gamma_t^2 \dd t - \int_0^T \gamma_t \dd W_t)
\end{align}
And let 
\begin{align}
\dd \tilde{W_t} = \gamma_t \dd t + \dd W_t
\end{align}
Then we have
\begin{align}
\dd X_t &= K(\theta - X_t)\dd t + \sigma (\dd W_t)\\
            &= K(\theta - X_t)\dd t + \sigma (\tilde{W_t} - \gamma_t \dd t)\\
            &= K(\theta - \frac{\sigma \gamma_t}{K} - X_t)\dd t + \sigma \tilde{W_t}\\
            & = K(\theta^* - X_t)\dd t + \sigma \tilde{W_t}\\
\end{align}
The last equality is given by
\begin{align}
\gamma_t = \frac{-K(\theta^* - \theta)}{\sigma}
\end{align}
Therefore, we obtain
\begin{align}
\frac{\dd \tilde{\Pb}}{\dd \Pb} =  ( -\int_0^T \frac{1}{2}  (\frac{-K(\theta^* - \theta)}{\sigma})^2 \dd t - \int_0^T  \frac{-K(\theta^* - \theta)}{\sigma} \dd W_t)
\end{align}

\subsection*{Solution to Exercise \ref{exer:8.6}}
(1)
\begin{align}
\Eb[X_{(1-t)}] &= 0\\
Cov[X_{(1-t)}, X_{1-s}] &= Cov[W_{(1-t)} - tW_1, W_{1-s} - sW_1]\\
    &= Cov[W_{(1-t)}, W_{1-s} ] - Cov[W_{(1-t)},   sW_1] -  Cov[W_{(1-s)},   tW_1] + Cov[tW_1, sW_1]\\
    &= (1-t) \wedge (1-s) - (1-t)(1-s)\\
    &= t \wedge s - st
\end{align}
(2)
\begin{align}
\Eb[ (1-t)W_{t/(1-t)}] &= 0\\
Cov[(1-t)W_{t/(1-t)}, (1-s)W_{s/(1-s)}] &= (1-t) (1-s) ( t/(1-t) \wedge  s/(1-s)  )\\
    &= t \wedge s - st
\end{align}
\begin{align}
\Eb[ tW_{(1/t) - 1}] &= 0\\
Cov[tW_{t(1/t) - 1}, sW_{s(1/s) - 1}] &= (1-t) (1-s) ( t(1/t) - 1 \wedge  s(1/s) - 1  )\\
    &= t \wedge s - st
\end{align}
(3)
\begin{align}
&\Eb[X_{\overline{t+s}} - X_s] = 0 \\
&Cov [ X_{\overline{t_1 + s}}  -  X_s,  X_{\overline{ t_2 + s}} - X_s] = Cov[{X_{\overline{t_1 + s}},  X_{\overline{ t_2 + s}}}] - Cov[X_s, X_{\overline{ t_2 + s}} ]- Cov[X_{\overline{t_1 + s}}, X_s ] + Cov[X_s, X_s]\\
               &= (\overline{t_1 + s}  \wedge \overline{t_2 + s})- \overline{t_1 + s}*\overline{t_2 + s} - (\overline{t_2 + s} \wedge s) + s(\overline{t_2 + s}) - (\overline{t_1 + s} \wedge s) + s(\overline{t_1 + s}) + s(1-s)\\
               &= t_1 \wedge t_2 - t_1 t_2
\end{align}

\subsection*{Solution to Exercise \ref{exer:8.7}}
We seek a Radon-Nikod\'ym derivative of the form
\begin{align}
\frac{\dd \Pbt}{\dd \Pb}
	&=	Z_T , &
Z_t
	&:=	\exp \Big( - \int_0^t \theta_s \dd W_s - \tfrac{1}{2} \int_0^t \theta_s^2 \dd s \Big) ,
\end{align}
for a process $\theta = (\theta_t)_{t \geq 0}$ to be determined.  By Girsanov's Theorem, the process $\Wt = (\Wt_t)_{t \geq 0}$ defined by
\begin{align}
\Wt_t
	&:=	W_t + \int_0^t \theta_s \dd s
\end{align}
is a Brownian motion under $\Pbt$, and $\dd W_t = \dd \Wt_t - \theta_t \dd t$.  Substituting into the dynamics of $X$, we obtain
\begin{align}
\dd X_t
	&=	\mu_t \dd t + \sig_t ( \dd \Wt_t - \theta_t \dd t )
	=		( \mu_t - \sig_t \theta_t ) \dd t + \sig_t \dd \Wt_t .
\end{align}
Matching the drift with $\widetilde{\mu}_t$ requires $\mu_t - \sig_t \theta_t = \widetilde{\mu}_t$, i.e.
\begin{align}
\theta_t
	&=	\frac{\mu_t - \widetilde{\mu}_t}{\sig_t} .
\end{align}
Thus, the desired change of measure is
\begin{align}
\frac{\dd \Pbt}{\dd \Pb}
	&=	\exp \Big( - \int_0^T \frac{\mu_s - \widetilde{\mu}_s}{\sig_s} \dd W_s - \tfrac{1}{2} \int_0^T \Big( \frac{\mu_s - \widetilde{\mu}_s}{\sig_s} \Big)^2 \dd s \Big) .
\end{align}

\subsection*{Solution to Exercise \ref{exer:8.8}}

By Girsanov's Theorem, we have that

\begin{align}\mathrm{d} \widetilde{W}_{t}=\lambda \sqrt{R_{t}} \mathrm{~d} t+\mathrm{d} W_{t},
\end{align}
in which $\widetilde{W}_{t}$ is a Brownian motion under the probability $\widetilde \Pb$. 

Under $\Pbt$, we therefore have
\begin{align}\mathrm{d} R_{t}=\left(\kappa \theta-\kappa R_{t}-\delta \lambda R_{t}\right) \mathrm{d} t+\delta \sqrt{R_{t}} \mathrm{~d} \widetilde{W}_{t}.
\end{align}

Comparing this with our aim \begin{align}\mathrm{d} R_{t}=\widetilde{\kappa}\left(\widetilde{\theta}-R_{t}\right) \mathrm{d} t+\delta \sqrt{R_{t}} \mathrm{~d} \widetilde{W}_{t},
\end{align} we identify
\begin{align}\widetilde{\theta}=\frac{\kappa \theta}{\kappa+\delta \lambda}, \quad \widetilde{\kappa}=\kappa+\delta \lambda
\end{align}

\subsection*{Solution to Exercise \ref{exer:8.9}}
(a) Using the application of the two-dimensional It\^o formula for the product rule, we obtain
\begin{align}
\dd (X_t Z_t)
	&=	X_t \dd Z_t+Z_t \dd X_t + \dd [X,Z]_t	\\
	&=	X_t (-\alpha Z_t \dd W_t) + Z_t (\frac{1}{X_t} \dd t + \alpha X_t \dd W_t ) - \alpha^2 X_t Z_t \dd t \\
	&= \Big( \frac{Z_t}{X_t}  - \alpha^2X_t Z_t \Big) \dd t .
\end{align}
(b) Using It\^o's formula with $f(X_t Z_t)=(X_t Z_t)^2$, we have using part (a) that
\begin{align}
\dd ((X_t Z_t)^2)
	&=	2 X_t Z_t \dd (X_t Z_t) + \dd [XZ,XZ]_t \\
	&=	2 X_t Z_t \Big( \frac{Z_t}{X_t}  - \alpha^2X_t Z_t \Big) \dd t + 0 \\
	&=	(2Z_t^2 -2\alpha^2(X_t Z_t)^2) \dd t .
\end{align}
(c) By It\^o's formula:
\begin{align}
\dd \log Z_t 
	&=	\frac{1}{Z_t}\dd Z_t - \frac{1}{2Z_t^2} \dd [Z,Z]_t = -\frac{1}{2}\alpha^2 \dd t - \alpha \dd W_t .
\end{align}
Integrating both sides, and using $\log Z_0 = 0$,  we obtain
\begin{align}
\log Z_t 
	&=	 -\frac{1}{2}\alpha^2 t - \alpha W_t &
	&\Longrightarrow&
Z_t
	&=	\ee^{-\frac{1}{2}\alpha^2t - \alpha W_t} . \label{eq:exer.8.9.Z}
\end{align}
(d) Using \eqref{eq:exer.8.9.result} with
\begin{align}
Y_t
	&=(X_t  Z_t )^2, &
G_t 
	&= 2Z_t^2, &
D 
	&= 2 \alpha^2, &
Y_0
	&=	x^2
\end{align}
we obtain
\begin{align}
X_t
	&=	\frac{1}{Z_t}\ee^{-\alpha^2 t}\sqrt{x^2 + 2\int_0^t \ee^{2\alpha^2s} Z_s^2 \dd s} ,
\end{align}
where $Z_t$ is given by \eqref{eq:exer.8.9.Z}.

\subsection*{Solution to Exercise \ref{exer:8.10}}
By Girsanov's Theorem, the process
\begin{align}
\widetilde{W}_t = W_t + \int_0^t \sigma_s \dd s
\end{align}
is a standard Brownian motion under $\Pbt$.  Thus, we have
\begin{align}
\Eb Z_t \log Z_t
	&=	\Eb \log Z_t \Eb[ Z_T | \Fc_t] 
	=		\Eb Z_T \log Z_t 
	=		\Ebt \log Z_t \\
	&=	\Ebt \Big( -\frac{1}{2}\int_0^t\sigma_s^2 \dd s - \int_0^t \sigma_s \dd W_s \Big) \\
	&=	\Ebt \Big(+\frac{1}{2}\int_0^t\sigma_s^2 \dd s - \int_0^t \sigma_s \dd \widetilde{W}_s \Big)
	=		\Ebt \Big(\frac{1}{2}\int_0^t \sigma_s^2 \dd s \Big) ,
\end{align}
where we have used the fact that It\^o integrals are martingales.

\subsection*{Solution to Exercise \ref{exer:8.11}}
(a) We have $Z = \ee^{a W_T}/ \Eb \ee^{a W_T} \geq 0$ and $\Eb Z = \Eb \ee^{a W_T} / \Eb \ee^{a W_T} = 1$.
\\[0.5em]
(b) First, we note that
\begin{align}
Z
	&=	\ee^{a W_T - a^2 T / 2} .
\end{align}
Thus, we have
\begin{align}
\Ebt W_T
	&=	\Ebt \d_b \ee^{b W_T} \Big|_{b=0}
	=		\Eb \d_b \ee^{(a + b) W_T - a^2 T / 2} \Big|_{b=0}
	=		\d_b \Eb \ee^{(a + b) W_T - a^2 T / 2} \Big|_{b=0} \\
	&=		\d_b \ee^{( a + b )^2 T / 2 - a^2 T / 2} \Big|_{b=0} 
	=		\d_b \ee^{b^2 T / 2 + a b T} \Big|_{b=0}
	=		( b T + a T ) \ee^{b^2 T / 2 + a b T} \Big|_{b=0}
	=		a T .
\end{align}
Or, alternatively
\begin{align}
\Ebt W_T
	&=	\Eb Z W_T
	=		\Eb \ee^{a W_T - a^2 T / 2} W_T
	=		\d_a \Eb \ee^{a W_T - a^2 T / 2} + a T \Eb \ee^{a W_T - a^2 T / 2}
	=		\d_a 1 + a T
	=		a T ,
\end{align}
where we have used the fact that $\Eb \ee^{a W_T - a^2 T / 2} = 1$.
\\[0.5em]
(c) First we note that
\begin{align}
\Eb( \ee^{a W_T - a^2 T / 2} | \Fc_t ) 
	&=	\ee^{a W_t - a^2 t / 2}
\end{align}
Thus, we have
\begin{align}
\Ebt W_t
	&=	\Eb Z W_t
	=		\Eb \ee^{a W_T - a^2 T / 2} W_t
	=		\Eb \Eb( \ee^{a W_T - a^2 T / 2} | \Fc_t) W_t
	=		\Eb \ee^{a W_t - a^2 t / 2} W_t
	=		a t ,
\end{align}
where in the last step we have used the result from part (b).

\subsection*{Solution to Exercise \ref{exer:8.12}}
(a) Define $r(x,y) = \sqrt{x^2 + y^2}$.  We have
\begin{align}
\d_x r(x,y)
	&=	\frac{x}{\sqrt{x^2+y^2}} 
	=		\frac{x}{r(x,y)} , &
\d_y r(x,y)
	&=	\frac{y}{r(x,y)} , \\
\d_x^2 r(x,y)
	&=	\frac{1}{r(x,y)} + \frac{-x^2}{\sqrt{x^2+y^2}^3} 
	=		\frac{r^2(x,y) - x^2}{r^3(x,y)}
	=		\frac{y^2}{r^3(x,y)} , &
\d_y^2 r(x,y)
	&=	\frac{x^2}{r^3(x,y)} , \\
\d_x \d_y r(x,y)
	&=	\frac{- x y}{\sqrt{x^2+y^2}^3}
	=		\frac{- x y}{r^3(x,y)} .
\end{align}
Thus, by It\^o's Lemma, we have
\begin{align}
\dd R_t
	&=	\d_x r(X_t,Y_t) \dd X_t + \d_y r(X_t,Y_t) \dd Y_t 
			+ \frac{1}{2} \Big( \d_x^2 r(X_t,Y_t) \dd [X,X]_t + 2 \d_x \d_y r(X_t,Y_t) \dd [X,Y]_t + \d_y^2 r(X_t,Y_t) \dd [Y,Y]_t \Big) , \\
	&=	\frac{X_t}{R_t} \dd W_t + \frac{Y_t}{R_t} \dd B_t
			+ \frac{1}{2} \frac{X_t^2 + Y_t^2}{R_t^3} \dd t \\
	&=	\cos \Theta_t \dd W_t + \sin \Theta_t \dd B_t + \frac{1}{2 R_t} \dd t ,
\end{align}
where we have used $\dd[X,X]_t = \dd[Y,Y]_t = \dd t$ and $\dd[X,Y]_t = 0$.
\\[0.5em]
(b) Define $\theta(x,y) = \arctan(y/x)$.  We have
\begin{align}
\d_x \theta(x,y)
	&=	\frac{- y }{r^2(x,y)}, &
\d_y \theta(x,y)
	&=	\frac{x}{r^2(x,y)}, \\
\d_x^2 \theta(x,y) 
	&=	\frac{2 x y }{r^4(x,y)}, &
\d_y^2 \theta(x,y)
	&=	\frac{-2 x y }{r^4(x,y)}, \\
\d_x \d_y \theta(x,y)
	&=	\frac{y^2 - x^2 }{r^4(x,y)}.
\end{align}
Thus, by It\^o's Lemma, we have
\begin{align}
\dd \Theta_t
	&=	\d_x \theta(X_t,Y_t) \dd X_t + \d_y \theta(X_t,Y_t) \dd Y_t 
			+ \frac{1}{2} \Big( \d_x^2 \theta(X_t,Y_t) \dd [X,X]_t + 2 \d_x \d_y \theta(X_t,Y_t) \dd [X,Y]_t + \d_y^2 \theta(X_t,Y_t) \dd [Y,Y]_t \Big) , \\
	&=	\frac{-Y_t}{R_t^2} \dd W_t + \frac{X_t}{R_t^2} \dd B_t \\
	&=	\frac{-\sin \Theta_t}{R_t} \dd W_t + \frac{\cos \Theta_t}{R_t} \dd B_t .
\end{align}

\subsection*{Solution to Exercise \ref{exer:8.13}}
(a) Inserting the dynamics of $S$ and $M$ into the dynamics of $X$, we obtain
\begin{align}
\dd X_t
	&=	\Del_t S_t ( \mu - r ) \dd t + \Del_t \sig S_t \dd W_t + r X_t \dd t .
\end{align}
By the It\^o formula, we have
\begin{align}
\dd \frac{1}{S_t}
	&=	\frac{1}{S_t} \Big( ( \sig^2 - \mu ) \dd t - \sig \dd W_t \Big) .
\end{align}
Thus, using the product rule, we find
\begin{align}
\dd \frac{X_t}{S_t}
	&=	X_t \dd \frac{1}{S_t} + \frac{1}{S_t} \dd X_t + \dd \Big[ X , \frac{1}{S} \Big]_t \\
	&=	\frac{X_t}{S_t} \Big( ( \sig^2 - \mu ) \dd t - \sig \dd W_t \Big)
			+ \Del_t ( \mu - r ) \dd t + \Del_t \sig \dd W_t + r \frac{X_t}{S_t} \dd t - \Del_t \sig^2 \dd t \\
	&=	( \mu - r - \sig^2 ) \Big( \Del_t - \frac{X_t}{S_t} \Big) \dd t + \sig \Big( \Del_t - \frac{X_t}{S_t} \Big) \dd W_t .
\end{align}
Now, for a constant $\gam$, define a change of measure
\begin{align}
\frac{\dd \Pbh}{\dd \Pb}
	&:=	\exp \Big( - \gam W_T - \tfrac{1}{2} \gam^2 T \Big) .
\end{align}
By Girsanov's Theorem, the process $\Wh$ defined by $\dd \Wh_t = \gam \dd t + \dd W_t$ is a $(\Pbh,\Fb)$ Brownian motion.  Writing the dynamics of $X/S$ in terms of $\Wh$, we obtain
\begin{align}
\dd \frac{X_t}{S_t}
	&=	( \mu - r - \sig^2 - \sig \gam ) \Big( \Del_t - \frac{X_t}{S_t} \Big) \dd t + \sig \Big( \Del_t - \frac{X_t}{S_t} \Big) \dd \Wh_t .
\end{align}
The $\dd t$ term vanishes for every $\Del$ if and only if
\begin{align}
\gam
	&=	\frac{\mu - r}{\sig} - \sig ,
\end{align}
in which case $X/S$ is a $(\Pbh,\Fb)$ martingale.  Thus, the desired change of measure is
\begin{align}
\frac{\dd \Pbh}{\dd \Pb}
	&=	\exp \Big( - \Big( \frac{\mu - r}{\sig} - \sig \Big) W_T - \frac{1}{2} \Big( \frac{\mu - r}{\sig} - \sig \Big)^2 T \Big) .
\end{align}
(b) Using $\dd (1/S_t)$ from part (a) and $\dd M_t = r M_t \dd t$, we have
\begin{align}
\dd \frac{M_t}{S_t}
	&=	M_t \dd \frac{1}{S_t} + \frac{1}{S_t} \dd M_t
	=		\frac{M_t}{S_t} \Big( ( \sig^2 - \mu + r ) \dd t - \sig \dd W_t \Big) \\
	&=	\frac{M_t}{S_t} \Big( ( \sig^2 - \mu + r + \sig \gam ) \dd t - \sig \dd \Wh_t \Big)
	=		- \sig \frac{M_t}{S_t} \dd \Wh_t .
\end{align}
Thus, $\Theta_t = 0$ and $\Sig_t = - \sig M_t / S_t$.  Note that $M/S$ is a $(\Pbh,\Fb)$ martingale, which is consistent with part (a), since $M$ itself has the form of $X$ with $\Del \equiv 0$.


\section{Exercises from Section \ref{sec:exercises-9}}

\subsection*{Solution to Exercise \ref{exer:9.1}}
Suppose that $u \in C^{1,2}(\Rb)$. Then 
\begin{align}
\dd u(t,X_t)
	&= \d_t u(t,X_t)\, \dd t
	+ \d_x u(t,X_t)\, \dd X_t 
	+ \half \d_x^2 u(t,X_t^2)\, \dd X_t \dd X_t \\
	&= \( \d_t + \mu(t,X_t,Y_t,Z_t) \d_x + \half \sig(t,X_t,Y_t) \d_x^2 \) u(t,X_t)\, \dd t
	+ \sig(t,X_t,Y_t) \d_x u(t,X_t)\, \dd W_t. \label{eq:du}
\end{align}
By our assumption that $Y_t = u(t,X_t)$, we see that 
\begin{align}
\dd u(t,X_t)
	=\dd Y_t
	= - f(t,X_t,u(t,X_t),Z_t)\, \dd t + Z_t\, \dd W_t. \label{eq:dY}
\end{align}
If we assume that $\mu$ and $\sig$ are continuous, we can equate the expressions in front of the $\dd t$ and $\dd W_t$ terms in \eqref{eq:du} and \eqref{eq:dY}. Equating the $\dd W_t$ terms gives
\begin{align}
Z_t
	= \sig( t,X_t, u(t,X_t)) \d_x u(t,X_t). \label{eq:Z}
\end{align}
By definition of $Y_t$, $u(T,X_T) = \varphi(X_T)$. Furthermore, the equality corresponding to the $\dd t$ term in \eqref{eq:du} and \eqref{eq:dY} holds for all paths of $X_t$. Letting $\Ac^u(t)$ and $f^u(t)$ be defined as in the problem statement, we see that for all $(t,x) \in [0,T] \times \Rb$, 
\begin{align}
&(\d_t + \Ac^u(t) ) u(t,\cdot) + f^u (t)
	= 0
& u(T,\cdot)
	= \varphi.
\end{align}

\subsection*{Solution to Exercise \ref{exer:9.2}}
We note first that $u(t,X_t)$ is \textit{not} a martingale. Let $\Fc_t$ be the filtration generated by $W_t$. Then for $s$, $0 \leq s \leq t$,
\begin{align}
\Eb[ u(t,X_t) | \Fc_s]
	&= \Eb\[ 
		\Eb \[ 
			\exp \(- \int_t^T \dd u\, X_u\) \mid \Fc_t
		\] \mid \Fc_s
	      \] \\
	&= \Eb\[ 
			\exp \(- \int_t^T \dd u\, X_u\)  \mid \Fc_s
	      \] \\
	&\neq \Eb\[ 
			\exp \(- \int_s^T \dd u\, X_u\)  \mid \Fc_s
	      \] \\
	&= u(s,X_s).
\end{align}
However, the process
\begin{align}
&M_t
	:= e^{-\alpha(0,t)} u(t,X_t),
& \alpha(s,t) := \exp\( \int_s^t \dd u\, X_u\).
\end{align}
is a martingale since
\begin{align}
\Eb[M_t \mid \Fc_s]
	&= \Eb\[e^{-\alpha(0,t)} u(t,X_t) | \Fc_s\]\\
	&= \Eb\[ 
		e^{-\alpha(0,T)} \mid \Fc_s
	      \] \\
	&= e^{-\alpha(0,s)} \Eb\[ e^{-\alpha(s,T)}\mid \Fc_s \] \\ 
	&= M_s.
\end{align}
Taking the differential of $M_t$ yields
\begin{align}
\dd M_t	
	&= \dd \[ e^{-\alpha(0,t)} \] u(t,X_t)
		+ e^{-\alpha(0,t)}\, \dd u(t,X_t) \\
	&= -X_t e^{-\alpha(0,t)} u(t,X_t)\, \dd t
		+ e^{-\alpha(0,t)} 
		\[
\( \d_t + \kappa ( \theta- X_t) \d_x + \half\delta^2 X_t \d_x^2\)\, \dd t + \delta \sqrt X_t\, \dd W_t
		\] \\
	&= e^{-\alpha(0,t)}  \( \d_t + \kappa ( \theta- X_t) \d_x + \half \delta^2 X_t \d_x^2 - X_t\) 
		u(t,X_t)\, \dd t 
	  + e^{-\alpha(0,t)}\delta \sqrt X_t\, \dd W_t.
\end{align}
Noting that
\begin{align}
u(T,x)
	= \Eb\[ \exp\( - \int_T^T \dd u\, X_u\) \mid X_T = x\]
	= 1,
\end{align}
since $M_t$ is a martingale for any path $X_t$ takes we see that $u$ satisfies the problem
\begin{align}
&(\d_t + \Ac) u
	= 0,
&u(T,\cdot) 
	\equiv 1 \label{prob1}
\end{align}
where
\begin{align}
\Ac 
	= \kappa (\theta - x) \d_x + \half \delta^2 x \d_x^2 - x.
\end{align}
Let us make the guess $u(t,x) = e^{- x A(t) - B(t)}$ for some functions $A(t)$ and $B(t)$. Plugging this function into \eqref{prob1} that gives
\begin{align}
&x \left[\half \delta ^2 A^2(t)+\kappa  A(t)-A'(t)-1\right]-\[\theta  \kappa  A(t)+B'(t)\]
	= 0, 
& (t,x) \in [0,T] \times \Rb
	\\
& A(T) 
	= 0 
& B(T)
	= 0.
\end{align}
Since this most hold for all $x$, we have the coupled ODE's
\begin{align}
&A'(t)
	=\half \delta ^2 A^2(t)+\kappa  A(t)-1,
&A(T)
	= 0 \label{prob:odeA}\\
&B'(t)
	=-\theta  \kappa  A(t),
&B(T)
	= 0.\label{prob:odeB}
\end{align}
It can be seen that for $\beta := \sqrt{ 2 \delta^2 + \kappa^2}$, 
\begin{align}
A(t) 
	&=-\delta ^{-2}\kappa +\delta ^{-2}\beta 
		\tanh \left(
			\frac{1}{2} (T-t) \beta
			+  \tanh ^{-1}\left(\frac{\kappa }{\beta}\right)
		\right)		
\end{align}
satisfies the problem \eqref{prob:odeA}. Thus,
\begin{align}
B(t) 
	= \delta^{-2} \theta  \kappa  \left(\log \left(\frac{2 \delta ^2}{\beta^2}\right)+2 
	\log \left(\cosh \left(\frac{1}{2} \beta (T-t)+ \tanh ^{-1}\left(\frac{\kappa }{\beta}\right)
	\right)\right)+\kappa  (t-T)\right)
\end{align}

\subsection*{Solution to Exercise \ref{exer:9.3}}
We will use Levy's characterization of Brownian motion to show that $B_t$ is a Brownian motion. $B_t$ is a sum of integrals with respect to independent Brownian motions. Therefore, $B_0 = 0$, $B_t$ has continuous sample paths and $B_t$ is a martingale. What remains to be shown is that $\dd [B, B]_t = \dd t$. Note that since $ \dd W_t^{(i)}\, \dd W_t^{(j)} = \delta_{ij}\, \dd t$,
\begin{align}
\dd B_t\, \dd B_t
	&= \frac{1}{ R_t }
		\sum_{i=1}^d \sum_{j=1}^d X_t^{(i)} X_t^{(j)} \dd W_t^{(i)}\, \dd W_t^{(j)}\\
	&= \frac{1}{ R_t } \sum_{i=1}^d \( X_t^{(i)}\)^2\, \dd t \\
	&= \dd t.
\end{align}
Therefore, $B_t$ is a Brownian motion. Note that
\begin{align}
\sqrt R_t \, \dd B_t
	= \sum_{i=1}^d X_t^{(i)}\, \dd W_t^{(i)}.
\end{align}
Thus, using the definition of $X_t^{(i)}$,
\begin{align}
\dd R_t
	= \sum_{i=1}^d \dd \( X_t^{(i)}\)^2
	&= \sum_{i=1}^d 2 X_t^{(i)} \dd X_{t}^{(i)} + \dd X_{t}^{(i)} \\
	&= \sum_{i=1}^d \[ -b (X_t^{(i)})^2 + \frac{1}{4}\sig^2\]\, \dd t + \sig X_t^{(i)}\, \dd W_t^{(i)} \\
	&= \[ -b R_t + \frac{d}{4}\sig^2\]\, \dd t + \sig \sqrt R_t\, \dd B_t.
\end{align}

\subsection*{Solution to Exercise \ref{exer:9.4}}
Let $\Fc_t$ be the canonical filtration for the process $(W_t, B_t)$. Let $W_t^{(1)} = W_t$ and let $W_t^{(2)} = \rho \, \dd W_t + \sqrt{ 1 - \rho^2}\, \dd B_t$. Note that $W_0^{(2)} =0$, and $W_t^{(2)}$ is continuous. Furthermore, for $s$, $0 \leq s \leq t$, 
\begin{align}
\Eb[ W_t^{(2)} \mid \Fc_s]
	= \Eb[ \rho \, \dd W_t + \sqrt{ 1 - \rho^2}\, \dd B_t \mid \Fc_s]
	= \rho \, \dd W_s + \sqrt{ 1 - \rho^2}\, \dd B_s
	= W_s^{(2)}
\end{align}
and
\begin{align}
\dd W_t^{(2)}\, \dd W_t^{(2)}
	= \( \rho \, \dd W_t + \sqrt{ 1 - \rho^2}\, \dd B_t \) \(\rho \, \dd W_t + \sqrt{ 1 - \rho^2}\, \dd B_t\)
	= \dd t.
\end{align}
Thus, $W_t^{(2)}$ is a Brownian motion with respect to $\Fc_t$ such that $\dd W_t^{(1)}\dd W_t^{(2)} = \rho\dd t$, and 
\begin{align}
\dd Z_t
	= \kappa (\theta - Z_t)\, \dd t
	+ \sig \sqrt Z_t\, \dd W_t^{(2)}.
\end{align}
Now, let us turn our attention to $u(t,X_t,Z_t)$. For $s$ as before,
\begin{align}
\Eb[ u(t,X_t,Z_t) \mid \Fc_s]
	= \Eb[ \Eb[ \varphi(X_T) \mid \Fc_t] \mid \Fc_s]
	= \Eb[\varphi(X_T) \mid \Fc_s]
	= u(s,X_s,Z_s).
\end{align}
So, $u(t,X_t,Z_t)$ is a martingale. Taking the differential of $u(t,X_t, Z_t)$ gives
\begin{align}
\dd u(t,X_t, Z_t)
	&= \d_t u(t,X_t, Z_t)\, \dd t \\ &\quad
	  + \d_x u(t,X_t, Z_t)\, \dd X_t \\ &\quad
	  + \half \d_{x}^2 u(t,X_t, Z_t)\, \dd X_t \dd X_t \\ &\quad 
	  + \d_{xz}^2 u(t,X_t, Z_t)\, \dd X_t \dd Z_t \\ & \quad
	  + \d_z u(t,X_t, Z_t)\, \dd Z_t \\ & \quad
	  + \half \d_{z}^2 u(t,X_t, Z_t)\, \dd Z_t \dd Z_t \\
	&= \d_t u(t,X_t, Z_t)\, \dd t \\ &\quad
	  + \d_x u(t,X_t, Z_t)\( -\half Z_t \, \dd t + \sqrt Z_t\, \dd W_t^{(1)}\) \\ &\quad
	  + \half \d_{x}^2 u(t,X_t, Z_t) Z_t\, \dd t \\ &\quad 
	  + \d_{xz}^2 u(t,X_t, Z_t) \rho \sig Z_t\, \dd t \\ & \quad
	  + \d_z u(t,X_t, Z_t)\( \kappa (\theta - Z_t)\, \dd t + \sig \sqrt Z_t \, \dd W_t^{(2)}\) \\ & \quad
	  + \half \d_{z}^2 u(t,X_t, Z_t)\sig^2 Z_t\, \dd W_t^{(2)} \\
\end{align}
Recalling that $u(t,X_t,Z_t)$ is a martingale and collecting $\dd t$ terms yields
\begin{align}
\( \d_t 
+ \half Z_t \d_x^2
+ \half\sig^2 Z_t \d_z^2
+ \rho \sig Z_t \d_{xz}^2
- \half Z_t \d_x 
+ \kappa (\theta - Z_t) \d_z
\) u(t,X_t,Z_t)
	= 0.
\end{align}
The fact that the above equation holds for all paths $X_t$ and $Z_t$ along with 
\begin{align}
u(T,x,z)
	= \Eb[ \varphi(X_T)\mid X_T = x, Z_T =z]
	= \varphi(x),
\end{align}
gives that $u$ is the solution to the problem 
\begin{align}
&\( \d_t + \Ac\) u
	= 0,
& u(T,\cdot,\cdot)
	= \varphi \label{pde:heston}
\end{align}
with
\begin{align}
\Ac 
	= \half z \d_x^2
	+ \half\sig^2 z \d_z^2
	+ \rho \sig z \d_{xz}^2
	- \half z \d_x 
	+ \kappa (\theta - z) \d_z
\end{align}
Taking the Fourier transform of \eqref{pde:heston} gives
\begin{align}
&\( \d_t + \widehat \Ac\) \widehat u
	= 0,
&\widehat u(T,\cdot,\cdot)
	= \widehat \varphi, \label{pde:hestontrans}
\end{align}
with
\begin{align}
&\widehat \Ac
	= \half \sig^2 z \d_z^2
	+ ( \theta  \kappa +z (\kappa -i \xi  \rho  \sigma )) \d_z
	- \half z (  \xi^2 + i \xi).
\end{align}
Inserting our guess $\widehat u(t,\xi,z) = \widehat \varphi (\xi) e^{A(t,\xi) + z B(t,\xi)}$ into \eqref{pde:hestontrans} gives
\begin{align}
 \\
0 = &\( \d_t + \widehat \Ac\) \widehat u(t,\xi,z) \\
	= &\frac{1}{2} \widehat{\varphi }(\xi ) e^{A(t,\xi )+z B(t,\xi )} \left(2 \d_t A(t,\xi ) 
	+2 \theta  \kappa  B(t,\xi )\right)\\ & 
	+\frac{1}{2} z \widehat{\varphi }(\xi ) e^{A(t,\xi )+z B(t,\xi )} \left(2 \d_t B(t,\xi )-2 \kappa  B(t,\xi )+2 i \xi  \rho  \sigma  B(t,\xi )+\sigma ^2 B^2(t,\xi )-\xi  (\xi -i)\right).
\end{align}
Since this relation must hold for all $(t,\xi,z)$, we have the pair of ODEs 
\begin{align}
&B'(t)-\kappa  B(t)+i \xi  \rho  \sigma  B(t)+\frac{1}{2} \sigma ^2 B^2(t)-\frac{\xi ^2}{2}-\frac{i \xi }{2}
	= 0,
&B(T,z)
	= 0 \label{ode:B}\\
&\d_t A(t,\xi ) 
	+ \theta  \kappa  B(t,\xi )
	= 0,
&A(T,\xi)
	= 0 \label{ode:A}.
\end{align}
With 
\begin{align}
&\alpha(\xi)
	:= \sqrt{-\kappa ^2+2 i \kappa  \xi  \rho  \sigma +\xi  \sigma ^2 \left(\xi  \left(\rho ^2-1\right)-i\right)},
&\beta(\xi)
	:= \kappa - i \xi  \rho  \sigma, \\
&\gamma(\xi)
	:= \sqrt{\frac{\beta (\xi )^2}{\alpha (\xi )^2}+1}
\end{align}
we see that 
\begin{align}
B(t,\xi)
	= \sigma ^{-2}\beta (\xi )-\sigma ^{-2}\alpha (\xi ) 
	\tan \left(\frac{1}{2} \alpha (\xi ) (t-T)+\tan ^{-1}\left(\beta(\xi)\alpha^{-1} (\xi )\right)\right)
\end{align}
satisfies \eqref{ode:B}. Thus,
\begin{align}
A(t,\xi)
	&=
	-\theta  \kappa \sigma ^{-2} 
	\left(\log \left(\gamma^2(\xi) \right)
	+2 \log \left(
	\frac{\alpha (\xi ) \cos \left(\frac{1}{2} \alpha (\xi ) (t-T)\right)-\beta (\xi ) \sin \left(\frac{1}{2} \alpha (\xi ) (t-T)\right)
	}{\alpha (\xi ) \gamma(\xi) }\right)+\beta (\xi ) (t-T)\right)
\end{align}
solves \eqref{ode:A}.

\subsection*{Solution to Exercise \ref{exer:9.5}}
Let $Y_t = \log X_t$ then we can derive the diffusion process of $Y_t$
$$\dd Y_t = (\mu - \frac{1}{2}\sigma^2) \dd t + \sigma \dd W_t  $$
where the endpoints become $(\log l, \log r)$. And the generator of this process will be
$$ \mathcal{A}= (\mu - \frac{1}{2} \sigma^2)\partial_x + \frac{1}{2}\sigma^2 \partial^2 x  $$
In the case of killing boundaries we have
$$ \mathcal{A} \psi_n = \lambda_n \psi_n, \quad \quad \psi_n(\log l) = \psi_n(\log r) = 0   $$
Let $\psi_n =  ke^{cy} \sin \left( \frac{n \pi}{\log r - \log l} (y- \log l) \right)$  which satisfies the boundary conditions. And plug it into $ \mathcal{A} \psi_n = \lambda \psi_n$ to solve $c$ and $\lambda_n$ to get
\begin{align}
c &= -\frac{ 2\mu - \sigma^2 }{2\sigma^2}\\
\lambda_n &= -\frac{1}{2 \sigma^2}(\mu - \frac{1}{2} \sigma^2)^2 - \frac{1}{2}\sigma^2 (\frac{n\pi}{\log r - \log l})^2\\
\psi_n(y) &=  ke^{ (-\frac{2\mu - \sigma^2 }{2 \sigma^2})y} \sin \left( \frac{n \pi}{\log r - \log l} (y- \log l) \right).
\end{align}
In addition, we need to find the speed density
$$ m(y) = \frac{2}{\sigma^2}\exp (\int \dd y \frac{2(\mu - \sigma^2)}{\sigma^2}) = \frac{2}{\sigma^2} e^{(\frac{2\mu-\sigma^2}{\sigma^2} y)}.$$
Next, we normalize the eigenfunctions $(\phi_n)_{n \geq0}$ so that they satisfy $\<\phi_n,\phi_l\>_m = \delta_{nl}$. It is easy to check that setting $k = \frac{1}{\sqrt{2}}$ accomplishes this task.  By Corollary 9.5.5, we can express the $\Gamma_Y (t,x; T,y)$ by $\lambda_n, \psi_n(y)$, and $m(y)$ as
$$ \Gamma_Y (t,x; T,y) = m(y) \sum_n e^{(T-t)\lambda_n} \psi_n(y) \psi_n(x)$$
Now, to obtain the transition density of $X$ we recall $Y=\log X$. We have
\begin{align}
\Gamma_{X}(t,x; T,y) \dd y = \Gamma_{Y}(t,\log x; T,\log y) \dd \log y  = \Gamma_{Y}(t,\log x; T,\log y) \frac{1}{y} \dd y .
\end{align}
Therefore, we can express the explicit form of solution as
$$ \Gamma_{X}(t,x; T,y) = \frac{1}{y} m(\log y) \sum_n e^{(T-t)\lambda_n} \psi_n(\log y) \psi_n( \log x) $$
where $m$, $\lam_n$, and $\psi_n$ are as given above.

\subsection*{Solution to Exercise \ref{exer:9.6}}
(1)
\begin{align}
(\frac{1}{2} \partial^2_x + \frac{1}{2}\partial^2_y ) u(x,y) &=0 \\
u(x,a) &=\phi(x)
\end{align}
(2)
\begin{align}
 \widehat{\phi}(\om) &= \int dx e^{-i\om x}\phi(x)\\
 \Eb_{t,x} \phi(X_\tau) &= \frac{1}{2 \pi} \int d\om \widehat{\phi}(\om) \Eb_{t,x}e^{i\om X_\tau}\\
u(x,y) = \Eb_{t,x,y}\phi(X_\tau) &= \frac{1}{2\pi} \int d\om \widehat{\phi }(\om) \Eb_{t,x,y}[ \Eb_{t,x}e^{i\om X_\tau} | \tau]\\
     &=\frac{1}{2\pi} \int d\om \widehat{\phi} (\om)\Eb_{t,x,y}(e^{i\om x - \frac{1}{2} \om^2 (\tau -t)  })\\
     &= \frac{1}{2\pi} \int d\om \widehat{\phi }(\om) e^{i\om x - |a-y| \sqrt{\frac{2}{2} \om^2}}\\
     &=  \frac{1}{2\pi} \int d\om \widehat{\phi }(\om) e^{i\om x - |a-y| |\om|}
\end{align}
Let $\phi(x) = \delta_z(x)$ to derive 
\begin{align}
\Pb(X_\tau \in \dd z | X_t=x, Y_t= y) &=  \Eb_{t,x,y}\delta_z(X_\tau)\\
  &=     \frac{1}{2\pi} \int d\om e^{-i\om z} e^{i\om x - |a-y| |\om|} \\ 
  &= \frac{|a-y|}{(a-y)^2 + (x-z)^2}dz 
\end{align}
here we apply \em Fourier transform  \em: $\phi(x)=\delta_z(x) \Rightarrow \widehat{\phi}(\om) = e^{-i\om z}$ 
\\
(3)
\begin{align}
\frac{1}{2}(\partial^2_x + \partial^2_y) e^{i\om x - |a-y| |\om|}&= (-\om^2 + \om^2)e^{i\om x - |a-y| |\om|} = 0\\
 u(x,a) &=  \frac{1}{2\pi} \int d\om \widehat{\phi }(\om) e^{i\om x}
 =\phi(x)  
\end{align}

\subsection*{Solution to Exercise \ref{exer:9.7}}
First, we derive the dynamics of $R$ in the interior $(0,\infty)$.  For $d \geq 2$, define $r(w) := ( \sum_{i=1}^d w_i^2 )^{1/2}$ for $w \in \Rb^d \setminus \{0\}$.  We have
\begin{align}
\d_{w_i} r(w)
	&=	\frac{w_i}{r(w)} , &
\d_{w_i}^2 r(w)
	&=	\frac{1}{r(w)} - \frac{w_i^2}{r^3(w)} .
\end{align}
Thus, using the multivariate It\^o formula and $\dd [W^i,W^j]_t = \del_{ij} \dd t$, we obtain
\begin{align}
\dd R_t
	&=	\sum_{i=1}^d \frac{W_t^i}{R_t} \dd W_t^i + \frac{1}{2} \sum_{i=1}^d \Big( \frac{1}{R_t} - \frac{(W_t^i)^2}{R_t^3} \Big) \dd t
	=		\frac{d-1}{2 R_t} \dd t + \dd B_t , &
B_t
	&:=	\sum_{i=1}^d \int_0^t \frac{W_s^i}{R_s} \dd W_s^i .
\end{align}
The process $B$ is a continuous martingale with $\dd [B,B]_t = \sum_{i=1}^d (W_t^i)^2 / R_t^2 \dd t = \dd t$, and thus, by L\'evy's characterization of Brownian motion, $B$ is a Brownian motion.  When $d = 1$ we have $R = |W|$, and away from the origin $\dd R_t = \text{sgn}(W_t) \dd W_t =: \dd B_t$, which is again a Brownian motion.  Thus, in all cases, $R$ is a time-homogeneous diffusion on $(0,\infty)$ with
\begin{align}
\mu(r)
	&=	\frac{d-1}{2r} , &
\sig(r)
	&=	1 .
\end{align}
Using \eqref{eq:s.and.m}, the scale and speed densities of $R$ are
\begin{align}
s(r)
	&=	\exp \Big( - \int \dd r \, \frac{d-1}{r} \Big) = r^{1-d} , &
m(r)
	&=	\frac{2}{s(r)} = 2 r^{d-1} .
\end{align}
Now fix $x \in (0,\infty)$ and consider the endpoint $l = 0$.
\\[0.5em]
\textit{Case $d = 1$.}  We have $s(z) = 1$ and $m(z) = 2$.  Thus $S(0,z] = z$ and $S[z,x] = x - z$, so that
\begin{align}
I_0
	&=	\int_0^x \dd z \, 2 z
	=		x^2
	<		\infty , &
J_0
	&=	\int_0^x \dd z \, 2 (x - z)
	=		x^2
	<		\infty .
\end{align}
Hence, the origin is a regular endpoint.
\\[0.5em]
\textit{Case $d \geq 2$.}  We have $S(0,z] = \int_0^z \dd y \, y^{1-d} = \infty$ for every $z > 0$, since $1 - d \leq -1$.  Hence $I_0 = \infty$.  On the other hand, for $0 < z < x$ we have
\begin{align}
S[z,x]
	&=	\int_z^x \dd y \, y^{1-d}
	=		\left\{ \begin{aligned}
			&\log (x/z) , & & d = 2 , \\
			&\frac{z^{2-d} - x^{2-d}}{d-2} \leq \frac{z^{2-d}}{d-2} , & & d \geq 3 .
			\end{aligned} \right.
\end{align}
Therefore
\begin{align}
J_0
	&=	\int_0^x \dd z \, 2 z \log(x/z)
	=		\frac{x^2}{2}
	<		\infty , &
	&(d = 2) , \\
J_0
	&\leq	\int_0^x \dd z \, 2 z^{d-1} \frac{z^{2-d}}{d-2}
	=		\frac{x^2}{d-2}
	<		\infty , &
	&(d \geq 3) .
\end{align}
As $I_0 = \infty$ and $J_0 < \infty$, the origin is an entrance point when $d \geq 2$.  This is consistent with the well-known fact that a $d$-dimensional Brownian motion with $d \geq 2$ started away from the origin never hits the origin, whereas it can be started from the origin.

\subsection*{Solution to Exercise \ref{exer:9.8}}
(a) Using the It\^o formula, we have $\dd \log S_t = - \tfrac{1}{2} \sig^2 \dd t + \sig \dd \Wt_t$, and thus
\begin{align}
S_T
	&=	S_t \exp \Big( - \tfrac{1}{2} \sig^2 (T-t) + \sig ( \Wt_T - \Wt_t ) \Big) .
\end{align}
As $S_t \in \Fc_t$ and $\Wt_T - \Wt_t$ is independent of $\Fc_t$, we obtain
\begin{align}
\Ebt( p S_T^q | \Fc_t )
	&=	p S_t^q \ee^{ - \tfrac{1}{2} q \sig^2 (T-t)} \Ebt \ee^{ q \sig ( \Wt_T - \Wt_t ) }
	=		p S_t^q \ee^{ - \tfrac{1}{2} q \sig^2 (T-t)} \ee^{ \tfrac{1}{2} q^2 \sig^2 (T-t)} ,
\end{align}
where we have used the hint with $Z = \Wt_T - \Wt_t \sim \Nc(0,T-t)$ and $a = q \sig$.  Thus
\begin{align}
u(t,s)
	&=	p s^q \exp \Big( \tfrac{1}{2} \sig^2 q (q-1) (T-t) \Big) .
\end{align}
(b) By the tower property, the process $(u(t,S_t))_{0 \leq t \leq T}$ is a $(\Pbt,\Fb)$ martingale.  By the It\^o formula,
\begin{align}
\dd u(t,S_t)
	&=	\Big( \d_t + \tfrac{1}{2} \sig^2 s^2 \d_s^2 \Big) u(t,s) \Big|_{s = S_t} \dd t + \sig S_t \d_s u(t,S_t) \dd \Wt_t .
\end{align}
As $u(t,S_t)$ is a martingale, the $\dd t$ term must be zero.  Moreover, $u(T,S_T) = \Ebt( p S_T^q | \Fc_T ) = p S_T^q$.  Thus, $u$ satisfies
\begin{align}
0
	&=	\Big( \d_t + \tfrac{1}{2} \sig^2 s^2 \d_s^2 \Big) u(t,s) , &
u(T,s)
	&=	p s^q .
\end{align}
To check that the expression from part (a) satisfies this PDE, we compute
\begin{align}
\d_t u(t,s)
	&=	- \tfrac{1}{2} \sig^2 q (q-1) u(t,s) , &
\tfrac{1}{2} \sig^2 s^2 \d_s^2 u(t,s)
	&=	\tfrac{1}{2} \sig^2 q (q-1) u(t,s) ,
\end{align}
which sum to zero.  Clearly, we also have $u(T,s) = p s^q$.

\subsection*{Solution to Exercise \ref{exer:9.9}}
(a) Using
\begin{align}
M_t
	&=	\ee^{- \lam t} u(t,X_t) + \int_0^t \ee^{- \lam s} g(X_s) \dd s , 
\end{align}
we have that
\begin{align}
M_t
	&=	\ee^{- \lam t} \Eb \Big( \ee^{- \lam (T-t)} h(X_T) + \int_t^T \ee^{- \lam (s-t)} g(X_s) \dd s \Big| \Fc_t \Big) + \int_0^t \ee^{- \lam s} g(X_s) \dd s \\
	&=	\Eb \Big( \ee^{- \lam T} h(X_T) + \int_0^T \ee^{- \lam s} g(X_s) \dd s \Big| \Fc_t \Big) ,
\end{align}
which is a martingale because, for $0 \leq t_1 \leq t_2 \leq T$, we have using the short-hand $\Eb_t := \Eb( \, \cdot \, | \Fc_t )$ that
\begin{align}
\Eb_{t_1} M_{t_2}
	&=	\Eb_{t_1} \Eb_{t_2} \Big( \ee^{- \lam T} h(X_T) + \int_0^T \ee^{- \lam s} g(X_s) \dd s \Big) \\
	&=	\Eb_{t_1} \Big( \ee^{- \lam T} h(X_T) + \int_0^T \ee^{- \lam s} g(X_s) \dd s \Big) 
	=		M_{t_1} .
\end{align}
(b) We compute
\begin{align}
\dd M_t
	&=	\ee^{- \lam t} ( \d_t - \lam + \Ac) u(t,X_t) \dd t + ( \ldots ) \dd W_t + \ee^{- \lam t} g(X_t) \dd t , \\
\Ac
	&=	\mu(t,x) \d_x + \tfrac{1}{2} \sig^2(t,x) \d_x^2 .
\end{align}
As $M$ is a martingale, the sum of the $\dd t$ terms must be zero.  Thus, we have
\begin{align}
0
	&=	( \d_t - \lam + \Ac) u + g , &
u(T,x) 
	&=	h(x) ,
\end{align}
where the terminal condition comes from the definition of $u$
\begin{align}
u(T,X_T)
	&=	\Eb \Big( \ee^{- \lam (T-T)} h(X_T) + \int_T^T \ee^{- \lam (s-t)} g(X_s) \dd s \Big| \Fc_T \Big)
	=		h(X_T) .
\end{align}

\subsection*{Solution to Exercise \ref{exer:9.10}}
(a) As both $X$ and $Y$ are It\^o integrals of deterministic integrands, they are mean-zero Gaussian random variables, so $\mu = (0,0)$.  Their variances and covariance follow from the It\^o isometry:
\begin{align}
\Vb X
	&=	\int_0^T g^2(t) \dd t , &
\Vb Y
	&=	\int_0^T h^2(t) \dd t , &
\text{Co}\Vb(X,Y)
	&=	\Eb X Y
	=		\int_0^T g(t) h(t) \dd t .
\end{align}
The covariance is obtained by approximating $g$ and $h$ by simple functions and passing to the limit, exactly as in the derivation of the It\^o isometry.  Thus
\begin{align}
\Sig
	&=	\( \begin{array}{cc}
			\int_0^T g^2(t) \dd t & \int_0^T g(t) h(t) \dd t \\
			\int_0^T g(t) h(t) \dd t & \int_0^T h^2(t) \dd t
			\end{array} \) .
\end{align}
(b) We have $X_T = X_t + W_T - W_t$.  For $Y$, integrating by parts,
\begin{align}
Y_T
	&=	Y_t + \int_t^T X_s \dd g(s)
	=		Y_t + \int_t^T \dd ( X_s g(s) ) - \int_t^T g(s) \dd X_s \\
	&=	Y_t + X_T g(T) - X_t g(t) - \int_t^T g(s) \dd W_s \\
	&=	Y_t + X_t ( g(T) - g(t) ) + g(T) ( W_T - W_t ) - \int_t^T g(s) \dd W_s .
\end{align}
(c) The Kolmogorov backward equation is
\begin{align}
0
	&=	(\d_t + \Ac(t)) \Gam( t,x,y;T,\eta) , &
\Gam(T,x,y;T,\eta)
	&=	\del_{\eta}(y) , &
\Ac(t)
	&=	\tfrac{1}{2} \d_x^2 + g'(t) x \d_y .
\end{align}
Given $(X_t,Y_t) = (x,y)$, the random variable $Y_T$ is Gaussian.  Its mean and variance are
\begin{align}
\mu
	&=	y + x ( g(T) - g(t) ) , \\
\sig^2
	&=	\Eb \Big( g(T) (W_T-W_t) - \int_t^T g(s) \dd W_s \Big)^2
	=		g^2(T)(T-t) + \int_t^T g^2(s) \dd s - 2 g(T) \int_t^T g(s) \dd s .
\end{align}
Therefore $\Gam(t,x,y;T,\eta) = f_{\mu,\sig^2}(\eta)$, the Gaussian density with the above mean and variance:
\begin{align}
f_{\mu,\sig^2}(\eta)
	&:=	\frac{1}{ \sqrt{2\pi\sig^2 }} \exp \Big( -\frac{(\eta - \mu)^2}{2 \sig^2} \Big) .
\end{align}

\subsection*{Solution to Exercise \ref{exer:9.11}}
(a) Note that
\begin{align}
p(x,y;t) \dd t
	&= 	\Pb( Y_\tau \in \dd t | X_0 = x, Y_0 = y ) 
	=	\Eb[ \del_t(Y_\tau) | X_0 = x , Y_0 = y ] \dd t.
\end{align}
Letting $\Ac$ denote the generator of $(X,Y)$, the function $p$ satisfies
\begin{align}
0
	&=	\Ac \, p(x,y;t) , &
\Ac
	&=	\mu(x) \d_x + \tfrac{1}{2} \sig^2(x) \d_x^2 + \d_y , &
p(m,y;t)
	&=	\del_t(y) .
\end{align}
(b) When $\mu(x) = 0$ and $\sig(x) = 1$ we have $\dd X_t = \dd W_t$, so $\tau$ is the first hitting time of a Brownian motion to level $m$ starting from $X_0 = x$.  Using the density \eqref{eq:tau.m.density} of the first hitting time of Brownian motion and noting that $Y_\tau = y + \tau$, we obtain
\begin{align}
p(x,y;t)
	&=	\frac{|m-x|}{(t-y)\sqrt{2 \pi (t-y)}} \exp \Big( \frac{ - (m - x)^2 }{2(t - y)} \Big) , &
t
	&\geq y .
\end{align}
One can verify directly that $(\tfrac{1}{2} \d_x^2 + \d_y) p(x,y;t) = 0$ and $\lim_{x \to m} p(x,y;t) = \del_t(y)$.

\subsection*{Solution to Exercise \ref{exer:9.12}}
(a) Using It\^o's Lemma (see also Exercise \ref{exer:8.12}), we find
\begin{align}
\dd R_t
	&=	\cos (\Phi_t) \dd W_t + \sin (\Phi_t) \dd B_t + \frac{1}{2 R_t}  \dd t , \\
\dd \Phi_t
	&=	\frac{- \sin (\Phi_t)}{R_t} \dd W_t + \frac{ \cos( \Phi_t)}{R_t} \dd B_t .
\end{align}
(b) Let $u(r,\phi) := \Pb( R_\tau = a | R_0 = r, \Phi_0 = \phi )$.  The function $u$ satisfies
\begin{align}
0
	&=	\Ac u , &
u(a,\phi)
	&= 1 , &
u(b,\phi)
	&=	0 , &
\Ac
	&=	\frac{1}{2 r} \d_r + \frac{1}{2} \d_r^2 + \frac{1}{2r^2} \d_\phi^2 .
\end{align}
One checks that $u$ depends on $r$ only and is given by
\begin{align}
u(r,\phi) \equiv u(r)
	&=	\frac{\log b - \log r}{\log b - \log a} = \Pb( R_\tau = a | R_0 = r, \Phi_0 = \phi) .
\end{align}
(c) By It\^o's Lemma,
\begin{align}
\dd f(R_t,\Phi_t)
	&=	\Ac f(R_t,\Phi_t) \dd t 
			+ f_r \Big( \cos (\Phi_t) \dd W_t + \sin (\Phi_t) \dd B_t \Big) 
			+	f_\phi \Big( \frac{- \sin (\Phi_t)}{R_t} \dd W_t + \frac{ \cos( \Phi_t)}{R_t} \dd B_t \Big) .
\end{align}
Integrating from $t= 0$ to $t = \tau$, taking expectations, and using that integrals with respect to $W$ and $B$ are martingales, we obtain \eqref{eq:exer.9.12.dynkin}.
\\[0.5em]
(d) Choosing $f(r,\phi) = r^2$ gives $\Ac f = \frac{1}{2r} 2 r + \frac{1}{2} \cdot 2 = 2$, so \eqref{eq:exer.9.12.dynkin} becomes $\Eb_{r,\phi} R_\tau^2 = r^2 + 2 \Eb_{r,\phi} \tau$.  Using part (b),
\begin{align}
\Eb_{r,\phi} R_\tau^2
	&=	a^2 \Big( \frac{\log b - \log r}{\log b - \log a} \Big) + b^2 \Big( \frac{\log r - \log a}{\log b - \log a} \Big) ,
\end{align}
and therefore
\begin{align}
\Eb_{r,\phi} \tau
	&=	\frac{1}{2} \Big( a^2 \Big( \frac{\log b - \log r}{\log b - \log a} \Big) + b^2 \Big( \frac{\log r - \log a}{\log b - \log a} \Big) - r^2 \Big) .
\end{align}


\section{Exercises from Section \ref{sec:exercises-10}}

\subsection*{Solution to Exercise \ref{exer:10.1}}
(a)  A Poisson process is a special case of a compound Poisson process where all jumps are of size one.  As the intensity of jumps is $\lam$, it follows that the L\'evy measure $\nu$ corresponding to $P$ is $\nu = \lam \del_1$ where $\del_1$ is a Dirac mass at one.\\
(b)  As $X_T = x + P_T - P_t$ and $P_T - P_t \sim \text{Poi}(\lam(T-t))$ we have
\begin{align}
u(t,x)
	&=	\sum_{k=0}^\infty \varphi( x + k ) \frac{\lam^k(T-t)^k}{k!}\ee^{\lam(T-t)} .
\end{align}
The Kolmogorov backward equation is for $u$ is
\begin{align}
(\d_t + \Ac) u(t,x)
	&=	0 , &
u(T,x)
	&=	\varphi(x) , &
\Ac f(x)
	&=	\int_\Rb \Big( f(x+z) - f(x) \Big) \nu(\dd z) = \lam ( f(x+1) - f(x) ) .
\end{align}
One can verify by direct substitution that the function $u$ satisfies the KBE.  We omit the details.

\subsection*{Solution to Exercise \ref{exer:10.2}}
(a)  We check that $Z$ satisfies the definition of a L\'evy process.  First, we have $Z_0 = X_{Y_0} = X_0 = 0$.  Second, for $0 \leq t_1 \leq t_2 \leq t_3 \leq t_4 < \infty$, we have 
\begin{align}
Z_{t_4} - Z_{t_3} 
	&= X_{Y_{t_4}} - X_{Y_{t_3}} &
	&\ind&
Z_{t_2} - Z_{t_1} 
	&= X_{Y_{t_2}} - X_{Y_{t_1}} ,
\end{align}
because $Y$ is increasing and both $X$ and $Y$ are L\'evy processes (and have independent non-overlapping increments).  Third, we have
\begin{align}
Z_{t_2} - Z_{t_1} 
	&= X_{Y_{t_2}} - X_{Y_{t_1}}
	\sim  X_{Y_{t_2} - Y_{t_1}}
	\sim 	X_{Y_{t_2 - t_1}} 
	=		Z_{t_2 - t_1} .
\end{align}
Thus, assuming $Z$ is continuous in probability (which we are allowed to assume for this problem), we have that $Z$ is a L\'evy process.\\
(b) Note that
\begin{align}
\Eb \ee^{\ii \xi X_t}
	&=	\ee^{t \psi_X(\xi)} , &
\psi_X(\xi)
	&=	\ii b \xi - \tfrac{1}{2} \sig^2 \xi^2 + \int_\Rb \nu(\dd z) \Big( \ee^{i \xi z } - 1 - \ii \xi z \Big) , \\
\Eb \ee^{\ii \xi Y_t}
	&=	\ee^{t \psi_Y(\xi)} , &
\psi_Y(\xi)
	&=	\ii \gam \xi + \int_{\Rb_+} \mu(\dd z) \Big( \ee^{i \xi z } - 1 \Big) .
\end{align}
Next, we compute
\begin{align}
\Eb \ee^{\ii \xi Z_t}
	&=	\Eb \Eb ( \ee^{\ii \xi X_{Y_t}} | Y_t )
	=		\Eb \ee^{Y_t \psi_X(\xi)} \\
	&=		\Eb \ee^{\ii ( -\ii \psi_X(\xi) Y_t )}
	=		\ee^{t \psi_Y( - \ii \psi_X(\xi))} 
	=		\ee^{t \psi_Z(\xi)}
\end{align}
Thus, we have $\psi_Z(\xi) = \psi_Y( - \ii \psi_X(\xi))$.\\
(c) When $\nu \equiv 0$ we have
\begin{align}
-\ii \psi_X(\xi)
	&=	b \xi + \ii \tfrac{1}{2} \sig^2 \xi^2 .
\end{align}
Using $\psi_Z(\xi) = \psi_Y( - \ii \psi_X(\xi))$, it follows in that
\begin{align}
\psi_Z(\xi) 
	&=	\ii \gam \Big( b \xi + \ii \tfrac{1}{2} \sig^2 \xi^2 \Big) + \int_{\Rb^+} \mu(\dd z) \Big( \ee^{\ii (b \xi + \ii \tfrac{1}{2} \sig^2 \xi^2) z} - 1 \Big) \\
	&=	\ii \gam b \xi - \tfrac{1}{2} \gam \sig^2 \xi^2 + \int_{\Rb_+} \mu(\dd z) \int_\Rb \dd x \, f_{b z, \sig^2 z}(x) \Big( \ee^{\ii \xi x} - 1 \Big) \\
	&=	\ii \gam b \xi - \tfrac{1}{2} \gam \sig^2 \xi^2 + \int_\Rb \dd x \int_{\Rb_+} \, \mu(\dd z)  f_{b z, \sig^2 z}(x) \Big( \ee^{\ii \xi x} - 1 \Big) , \\
f_{m,s^2}(x)
	&=	\frac{1}{\sqrt{2 \pi s^2}} \exp \Big( \frac{-(x-m)^2}{2 s^2} \Big) .
\end{align}
Hence we have
\begin{align}
A
	&=	\gam b , &
\Sig
	&=	\sqrt{\gam} \sig , &
\Pi(\dd x)
	&=	\int_{\Rb_+} \mu(\dd z) \, f_{b z, \sig^2 z}(x) \dd x .
\end{align}
The L\'evy-It\^o decomposition of $Z$ is then
\begin{align}
Z_t
	&=	A t + \Sig B_t + \int_{\Rb} z P(t,\dd z) , &
\Eb P(\dd t, \dd z)
	&=	\Pi(\dd z ) \dd t .
\end{align}

\subsection*{Solution to Exercise \ref{exer:10.3}}
(a) Let $f(X_t,Y_t) = X_t/Y_t$.  Using subscripts for partial derivatives and omitting the arguments of $f$ (i.e., $f_{xy}:= \d_x\d_y f(X_t,Y_t)$, etc.), we have
\begin{align}
f_x
	&=	\frac{1}{Y_t} = \frac{Z_t}{X_t}, &
f_y
	&=	\frac{-X_t}{Y_t^2} = \frac{-Z_t}{Y_t} , &
f_{xx}
	&=	0 , &
f_{yy}
	&=	\frac{2 X_t}{Y_t^3} = \frac{2 Z_t}{Y_t^2}, &
f_{xy}
	&=	\frac{-1}{Y_t^2} = \frac{-Z_t}{X_t Y_t}. 
\end{align}
Additionally, we note that
\begin{align}
\Del f(X_t,Y_t)
	&:=	f(X_{t-} + ( \ee^{\gam_t(z)}-1) X_{t-}, Y_{t-} + ( \ee^{g_t(z)}-1) Y_{t-} ) - f(X_{t-} , Y_{t-} )
	=		\frac{\ee^{\gam_t(z)}X_{t-}}{\ee^{g_t(z)}Y_{t-}} - \frac{X_{t-}}{Y_{t-}} \\
	&=		(\ee^{\gam_t(z)-g_t(z)} - 1 )Z_{t-} . 
\end{align}
Using the multidimensional version of It\^o's formula, we have
\begin{align}
\dd \Big( \frac{X_t}{Y_t} \Big)
	&=	\Big( f_x \mu_t X_t + f_y b_t Y_t + \tfrac{1}{2} f_{xx} \sig_t^2 X_t^2 + \tfrac{1}{2} f_{yy} a_t^2 Y_t^2 + f_{xy} \sig_t X_t a_t Y_t \Big) \dd t \\ & \quad
			+ \Big( f_x \sig_t X_t + f_y a_t Y_t \Big) \dd W_t \\ & \quad
			+ \int_\Rb \Big( \Del f(X_t,Y_t) \Big) \Nt(\dd t, \dd z) \\ & \quad
			+ \int_\Rb \Big( \Del f(X_t,Y_t) - (\ee^{\gam_t(z)}-1) X_{t-} f_x - (\ee^{g_t(z)}-1) Y_{t-} f_y \Big) \nu(\dd z)\dd t \\ 
	&=	\Big( \mu_t Z_t - b_t Z_t + 0 + a_t^2 Z_t - \sig_t a_t Z_t \Big) \dd t \\ & \quad
			+ \Big( \sig_t Z_t - a_t Z_t \Big) \dd W_t \\ & \quad
			+ \int_\Rb \Big( (\ee^{\gam_t(z)-g_t(z)} - 1 )Z_{t-} \Big) \Nt(\dd t, \dd z) \\ & \quad
			+ \int_\Rb \Big( (\ee^{\gam_t(z)-g_t(z)} - 1 )Z_{t-} - (\ee^{\gam_t(z)}-1) Z_t + (\ee^{g_t(z)}-1) Z_t \Big) \nu(\dd z) \dd t \\
	&=  \Big( \mu_t - b_t - a_t^2  - \sig_t a_t \Big) Z_t \dd t 
			+ \Big( \sig_t + a_t \Big) Z_t \dd W_t \\ & \quad
			+ \int_\Rb ( \ee^{\gam_t(z)-g_t(z)} - 1 )Z_{t-} \Nt(\dd t, \dd z) 
			+ \int_\Rb ( \ee^{\gam_t(z)-g_t(z)} - 1 - \ee^{\gam_t(z)}+ \ee^{g_t(z)} ) Z_t \nu(\dd z) \dd t .
\end{align}
(b) In order for $Z$ to be a martingale, the $\dd t$-term must be zero.  Thus, we must have
\begin{align}
\mu_t
	&=	b_t - a_t^2  + \sig_t a_t - \int_\Rb ( \ee^{\gam_t(z)-g_t(z)} - 1 - \ee^{\gam_t(z)}+ \ee^{g_t(z)} ) \nu(\dd z) .
\end{align}

\subsection*{Solution to Exercise \ref{exer:10.4}}
When $\eta$ is a Brownian motion $W$, we recognize $X$ as an OU process.  Thus, we shall attempt to solve the SDE for $X$ using the same procedure as we did for an OU process.  First, we define
\begin{align}
Y_t
	&=	X_t - \theta , &
	& \Rightarrow&
\dd Y_t
	&=	- \kappa Y_t \dd t + \dd \eta_t .
\end{align}
Next, we define
\begin{align}
Z_t
	&=	\ee^{\kappa t} Y_t , &
	&\Rightarrow&
\dd Z_t
	&=	\kappa \ee^{\kappa t} Y_t \dd t + \ee^{\kappa t} \dd Y_t
	=		\kappa \ee^{\kappa t} Y_t \dd t + \ee^{\kappa t} \Big( - \kappa Y_t \dd t + \dd \eta_t  \Big) \\ & & & &
	&=	\ee^{\kappa t} \dd \eta_t . 
\end{align}
Thus, we have
\begin{align}
Z_t
	&=	Z_0 + \int_0^t \ee^{\kappa s} \dd \eta_s , \\
Y_t
	&=	\ee^{-\kappa t} Y_0 + \int_0^t \ee^{-\kappa (t-s)}  \dd \eta_s \\
X_t
		&=	\theta + \ee^{-\kappa t} (X_0 - \theta) + \int_0^t \ee^{-\kappa (t-s)}  \dd \eta_s .
\end{align}
(b) When $\eta_t = \sig W_t + \int_\Rb z \Nt(t, \dd z)$ we have
\begin{align}
X_t
	&=	\theta + \ee^{-\kappa t} (X_0 - \theta) + \int_0^t \ee^{-\kappa (t-s)} \Big( \sig \dd W_s + \int_\Rb z \Nt(\dd s, \dd z) \Big), \\
m(t) 
	&= \Eb X_t
	 =	\theta + \ee^{-\kappa t} (X_0 - \theta) ,
\end{align}
where we have used the fact that integrals with respect to $W$ and $\Nt$ are martingales.  Finally, assuming $t < T$ and defining $I_t := X_t - m(t)$. Let $I_t = e^{-\kappa t}J_t$ where $J_t = \int_0^t \ee^{\kappa s} \Big( \sig \dd W_s + \int_\Rb z \Nt(\dd s, \dd z) \Big)$ and $J_t$ is martingale (note: $I_t$ is NOT martingale). Then we have
\begin{align}
C(t,T) 
	&= \Eb I_t I_T 
	=e^{-\kappa t-\kappa T}\Eb J_t J_T 
	= e^{-\kappa t-\kappa T}(\Eb J_t\Eb (J_T - J_t | \Fc_t ) + \Eb J_t^2) 
	=	e^{-\kappa t-\kappa T}(\Eb J_t^2) \\
	&= e^{-\kappa t-\kappa T} \int_0^t \ee^{2 \kappa s} \Big(  \sig^2 +  \int_\Rb z^2 \nu(\dd z)  \Big) \dd s 
	 =	\frac{e^{-\kappa(T-t)}-\ee^{-\kappa (T+t)}}{2 \kappa } \Big(  \sig^2 +  \int_\Rb z^2 \nu(\dd z)  \Big) .
\end{align}
where we have used the L\'evy-It\^o Isometry.

\subsection*{Solution to Exercise \ref{exer:10.5}}
Using It\^o's Lemma we compute
\begin{align}
\dd \varphi(X_t) 
	&=	\Big( \varphi'(X_t) \mu(t,X_t) + \tfrac{1}{2}\varphi''(X_t) \sig^2(t,X_t) \Big) \dd t \\ & \quad
			+ \int_\Rb \Big( \varphi( X_{t-} + \gam(t,X_{t-},z) ) - \varphi(X_{t-}) - \gam(t,X_t,z) \varphi'(X_t) \Big) \nu(\dd z)  \dd t \\ & \quad
			+ \ldots \dd W_t + \int_\Rb \ldots \Nt( \dd t, \dd z) .
\end{align}
Next, noting that integrals with respect to $W$ and $\Nt$ are martingales, and using the shorthand notation $\Eb_{t,x} \, \cdot \, = \Eb( \, \cdot \, | X_t = x)$ we have
\begin{align}
\Ac(t) \varphi(x)
	&=	\lim_{s \searrow t} \frac{1}{s-t} \Big( \Eb[ \varphi(X_s) | X_t = x ] - \varphi(x) \Big) \\
	&=	\lim_{s \searrow t} \frac{1}{s-t} \int_t^s \Eb_{t,x} \Big( \varphi'(X_u) \mu(u,X_u) + \tfrac{1}{2}\varphi''(X_u) \sig^2(u,X_u) \Big) \dd u \\ & \quad
			+ \lim_{s \searrow t} \frac{1}{s-t} \int_t^s \Eb_{t,x} \int_\Rb \Big( \varphi( X_{u-} + \gam(u,X_{u-},z) ) - \varphi(X_{u-}) - \gam(u,X_u,z) \varphi'(X_u) \Big) \nu(\dd z)  \dd u \\
	&=	\Big( \varphi'(x) \mu(t,x) + \tfrac{1}{2}\varphi''(x) \sig^2(t,x) \Big) \\ & \quad
			+ \int_\Rb \Big( \varphi( x + \gam(t,x,z) ) - \varphi(x) - \gam(t,x,z) \varphi'(x) \Big) \nu(\dd z) .
\end{align}

\subsection*{Solution to Exercise \ref{exer:10.6}}
(a) Applying the It\^o formula for L\'evy-It\^o processes (Theorem \ref{thm:ito.formula.levy}) to $S_t = S_0 \ee^{X_t}$, we obtain
\begin{align}
\dd S_t
	&=	S_{t-} \Big( \mu + \tfrac{1}{2} \sig^2 + \int_\Rb ( \ee^z - 1 - z ) \nu(\dd z) \Big) \dd t
			+ \sig S_{t-} \dd \Wt_t + S_{t-} \int_\Rb ( \ee^z - 1 ) \Nt(\dd t, \dd z) .
\end{align}
As $1/M_t = \ee^{-rt}/M_0$ is continuous and of finite variation, we have $\dd (1/M_t) = - r (1/M_t) \dd t$ and $\dd [S,1/M]_t = 0$.  Thus
\begin{align}
\dd \frac{S_t}{M_t}
	&=	\frac{S_{t-}}{M_t} \Big( \mu + \tfrac{1}{2} \sig^2 + \int_\Rb ( \ee^z - 1 - z ) \nu(\dd z) - r \Big) \dd t
			+ \sig \frac{S_{t-}}{M_t} \dd \Wt_t + \frac{S_{t-}}{M_t} \int_\Rb ( \ee^z - 1 ) \Nt(\dd t, \dd z) .
\end{align}
For $S/M$ to be a martingale under $\Pbt$, the $\dd t$ term must vanish.  Thus
\begin{align}
\mu
	&=	r - \tfrac{1}{2} \sig^2 - \int_\Rb ( \ee^z - 1 - z ) \nu(\dd z) .
\end{align}
(b) As $V/M$ is a $\Pbt$-martingale and $V_T = S_T^p$, we have
\begin{align}
V_t
	&=	M_t \Ebt \Big( \frac{S_T^p}{M_T} \Big| \Fc_t \Big)
	=		\ee^{-r(T-t)} \Ebt ( S_T^p | \Fc_t )
	=		\ee^{-r(T-t)} S_t^p \, \Ebt \ee^{p (X_T - X_t)} ,
\end{align}
where we have used $S_T = S_t \ee^{X_T - X_t}$ and the fact that $X_T - X_t$ is independent of $\Fc_t$.  Since $X$ is a L\'evy process, we have
\begin{align}
\Ebt \ee^{p (X_T - X_t)}
	&=	\exp \Big( (T-t) \Big( p \mu + \tfrac{1}{2} \sig^2 p^2 + \int_\Rb ( \ee^{pz} - 1 - p z ) \nu(\dd z) \Big) \Big) .
\end{align}
Inserting the value of $\mu$ from part (a), we conclude that $V_t = v(t,S_t)$ where
\begin{align}
v(t,s)
	&=	s^p \exp \Big( (T-t) \Big( (p-1) r + \tfrac{1}{2} \sig^2 p (p-1) + \int_\Rb ( \ee^{pz} - p \ee^z + p - 1 ) \nu(\dd z) \Big) \Big) .
\end{align}
As a check, note that when $p = 1$ we obtain $v(t,s) = s$, as it must be, since in this case $V = S$.

\subsection*{Solution to Exercise \ref{exer:10.7}}
Under $\Pb$, we have that $P$ is a pure-jump L\'evy process with L\'evy measure $\nu = \lam \del_1$. 
As such, we can write the dynamics of $P$ as follows
\begin{align}
P_t
	&=	\int_\Rb z N( t, \dd z) 
	=		\int_\Rb z \nu(\dd z) t  + \int_\Rb z \Nt(t , \dd z) , &
\nu(\dd z)
	&=	\lam \del_1(z) \dd z .
\end{align}
If we define a change of measure
\begin{align}
\frac{\dd \Pbh}{\dd \Pb}
	&:=	\exp \Big( - \int_0^T \int_\Rb \( \ee^{\eta_s(z)} - 1 - \eta_s(z) \) \nu(\dd z) \dd s + \int_0^T \int_\Rb \eta_s(z) \Nt(\dd s , \dd z) \Big) ,
\end{align}
we have under $\Pbh$ that
\begin{align}
\Nh(\dd t, \dd z)
	&=	N(\dd t, \dd z) - \nuh_t(\dd z) \dd t , &
\nuh_t(\dd z)
	&:=	\ee^{\eta_t(z)} \nu(\dd z) 
\end{align}
is a compensated Poisson random measure.  Suppose we set $\eta_t(z) = \log ( \widetilde{\lam} / \lam )$.  Then we have
\begin{align}
\nuh_t(\dd z)
	&= \widetilde{\lam} \del_1(z) \dd z .
\end{align}
and thus, under $\Pbh$, $P$ is a Poisson process with intensity $\widetilde{\lam}$.  As such, the change of measure we are seeking is
\begin{align}
\frac{\dd \Pbh}{\dd \Pb}
	&:=	\exp \Big( - \int_0^T \int_\Rb \( \ee^{\log ( \widetilde{\lam} / \lam )} - 1 - \log ( \widetilde{\lam} / \lam )  \) \nu(\dd z) \dd s + \int_0^T \int_\Rb \log ( \widetilde{\lam} / \lam ) \Nt(\dd s , \dd z) \Big) \\
	&=	\exp \Big( - T \int_\Rb \( \widetilde{\lam} / \lam - 1 - \log ( \widetilde{\lam} / \lam )  \) \nu(\dd z)  + \int_\Rb \log ( \widetilde{\lam} / \lam ) \Big( N(T , \dd z) - \nu(\dd z) T  \Big) \Big) \\
	&=	\exp \Big( - T \int_\Rb \( \widetilde{\lam} / \lam - 1 \) \nu(\dd z)  + \int_\Rb \log ( \widetilde{\lam} / \lam ) N(T , \dd z) \Big)	\\
	&=	\exp \Big( - T  \( \widetilde{\lam} -  \lam \) + \log ( \widetilde{\lam} / \lam ) P_T \Big)	 .
\end{align}

\subsection*{Solution to Exercise \ref{exer:10.8}}
(a) Note that $X$ is a L\'evy process with $(\mu,\sig,\nu)=(0,0,\lam_1 \del_1)$ and 
$Y$ is a L\'evy process with $(\mu,\sig,\nu)=(\mu,0,\lam_2 \del_1)$.
Thus, we have
\begin{align}
\Eb \ee^{\ii u X_t}
	&=	\ee^{t \psi_1(u)}, &
\psi_1(u)
	&=	\lam_1 ( \ee^{\ii u} - 1 ) , \\
\Eb \ee^{\ii u Y_t}
	&=	\ee^{t \psi_2(u)} , &
\psi_2(u)
	&=	\ii \mu u + \lam_2 ( \ee^{\ii u} - 1 ) .
\end{align}
(b) Using the result from part (a), we have
\begin{align}
\Eb \ee^{\ii u Z_t}
	&=	\Eb \Eb[ \ee^{\ii u X_{Y_t}} | Y_t ]
	=		\Eb \ee^{\psi_1(u) Y_t}
	=		\Eb \ee^{\ii ( - \ii \psi_1(u) ) Y_t }
	=		\ee^{t \psi_2( - \ii \psi_1(u))} 
	=		\ee^{t \psi(u)} , 
\end{align}
where we have defined
\begin{align}
\psi(u)
	&:= \psi_2( - \ii \psi_1(u))
	=		\mu \lam_1 ( \ee^{\ii u} - 1 ) + \lam_2 ( \ee^{ \lam_1 ( \ee^{\ii u} - 1 )} - 1 ) .
\end{align}
(c)  Clearly, $Z_t \in \Nb_0$.  Its generating function $G_{Z_t}$ is given by
\begin{align}
G_{Z_t}(s)
	&:=	\Eb s^{Z_t}
	=		\Eb \ee^{\ii (-\ii \log s) Z_t}
	=		\ee^{\psi(- \ii \log s) t } 
	=		\exp \Big( t \big( \mu \lam_1 ( s - 1 ) + \lam_2 ( \ee^{ \lam_1 ( s - 1 )} - 1 ) \big) \Big) .
\end{align}
For $|s| \leq 1$ we have
\begin{align}
G_{Z_t}(s)
	&=	\sum_{n=0}^\infty s^n \Pb( Z_t = n )
	=		\sum_{n=0}^\infty \frac{1}{n!} G_{Z_t}^{(n)}(0) s^n .
\end{align}
where $G_{Z_t}^{(n)}(0) = \d_s^n G_{Z_t}(s) |_{s=0}$.  Thus, we find
\begin{align}
\Pb( Z_t = n )
	&=	\frac{1}{n!} G_{Z_t}^{(n)}(0) .
\end{align}
And thus
\begin{align}
\Pb(Z_t = 0) = G_{Z_t}(0)
	&=	\exp \Big(  - t \big( \mu \lam_1  + \lam_2 ( 1 - \ee^{ - \lam_1 } ) \big) \Big) .
\end{align}

\subsection*{Solution to Exercise \ref{exer:10.9}}
With $f(X_t) = X_t^p$ we have
\begin{align}
\d_x f(X_t)
	&=	pX_t^{p-1} , &
\d_x^2 f(X_t)
	&=	p(p-1)X_t^{p-2} , & 
f(X_{t-} + X_{t-} (\ee^{\gam(z)}-1) ) - f(X_{t-})
	&=	X_{t-}^p(\ee^{\gamma(z)p}-1) .
\end{align}
And thus, we have
\begin{align}
\dd Y_t = \dd f(X_t)
	&=	\Big( p \mu X_t X_t^{p-1}  + \tfrac{1}{2} p(p-1)\sig^2 X_t^2 X_t^{p-2} 
			- \int (p X_t^{p-1}X_t(\ee^{\gamma(z)}-1) - X_t^p(\ee^{\gamma(z)p}-1)) \nu(\dd z) \Big) \dd t \\ &\quad
			+ p\sig X_tX_t^{p-1}\dd W_t + \int X_{t-}^p(\ee^{\gamma(z)p}-1)\Nt(\dd t, \dd z) \\
  &=	\Big( p \mu Y_t  + \tfrac{1}{2} p(p-1) \sig^2 Y_t 
			- \int (p Y_t(\ee^{\gamma(z)}-1) - Y_t(\ee^{\gamma(z)p}-1)) \nu(\dd z) \Big) \dd t \\ & \quad
    + p\sig Y_t\dd W_t + \int Y_{t-}(e^{\gamma(z)p}-1) \Nt(\dd t, \dd z) .
\end{align}

\subsection*{Solution to Exercise \ref{exer:10.10}}
(a) Set $f(S_t) = \log S_t$, so that $f'(S_t) = 1/S_t$, $f''(S_t) = -1/S_t^2$ and $f( S_{t-} + S_{t-} (\ee^z - 1) ) - f(S_{t-}) = z$.  Using the It\^o formula for L\'evy-It\^o processes,
\begin{align}
\dd X_t
	&=	- \frac{1}{2} \sig_t^2 \dd t + \sig_t \dd W_t + \int_\Rb z \Nt(\dd t, \dd z) - \int_\Rb ( \ee^z - 1 - z) \nu( \dd z ) \dd t .
\end{align}
Computing the quadratic variation of $X$, we obtain
\begin{align}
\dd Y_t
	&=	\dd [ X, X ]_t
	=		\sig_t^2 \dd t + \int_\Rb z^2 N(\dd t, \dd z) .
\end{align}
(b) Separate $X$ into a continuous component $X^c$ and a jump component $X^j$, with
\begin{align}
\dd X_t^c
	&=	- \frac{1}{2} \sig_t^2 \dd t + \sig_t \dd W_t , &
\dd X_t^j
	&=	\int_\Rb z \Nt(\dd t, \dd z) - \int_\Rb ( \ee^z - 1 - z) \nu( \dd z ) \dd t .
\end{align}
Since $X^c \ind X^j$ and $S_t = \ee^{X_t}$, we have $Q_t = Q_t^c Q_t^j$ with $Q_t^c := \Eb( \ee^{q X_T^c} | \Fc_t )$ and $Q_t^j := \Eb( \ee^{q X_T^j} | \Fc_t )$.  As $Q^c$ does not jump, we focus on $Q^j$:
\begin{align}
\Del Q_t^j
	&=	\Eb( \ee^{q X_T^j} | \Fc_t ) - \Eb( \ee^{q X_T^j} | \Fc_{t-} ) \\
	&=	\ee^{q (X_{t-}^j + \Del X_t^j)} \Eb( \ee^{q (X_T^j -  X_{t-}^j)} | \Fc_{t-} ) - \ee^{q X_{t-}^j} \Eb( \ee^{q (X_T^j -  X_{t-}^j)} | \Fc_{t-} ) \\
	&=	\Big( \ee^{q \Del X_t^j} - 1 \Big) \Eb( \ee^{q X_T^j} | \Fc_{t-} ) 
	=		Q_{t-}^j  \int_\Rb ( \ee^{qz} - 1 ) \Del N(t, \dd z) .
\end{align}
Finally, using $\Del Q_t = Q_t^c \Del Q_t^j$ and $Q_t^c = Q_{t-}^c$ (as $Q^c$ is continuous),
\begin{align}
\Del Q_t
	&=	Q_{t-} \int_\Rb ( \ee^{qz} - 1 ) \Del N(t, \dd z) .
\end{align}
(c) When $\nu \equiv 0$, the coefficient $\sig_t$ is not given as a deterministic function of $(t,X_t)$, so $X$ is not Markov in general.  Hence there is no function $u$ with $u(t,X_t) = Q_t$ in general.
\\[0.5em]
(d) When $\sig \equiv 0$, the process $X$ is a L\'evy process, hence Markov, so such a function $u$ exists.  Using the L\'evy-Khintchine formula,
\begin{align}
Q_t
	&=	\Eb( \ee^{q X_T} | \Fc_t ) 
	=		\ee^{q X_t} \Eb( \ee^{q X_{T-t}} ) 
	=		\ee^{q X_t + (T-t) \psi(-\ii q) }
	=:	u(t,X_t) , 
\end{align}
where $\psi$ is the L\'evy exponent of $X$,
\begin{align}
\psi(\xi)
	&:= - \ii \xi \int_\Rb (\ee^z - 1 - z)\nu(\dd z) + \int_\Rb (\ee^{\ii \xi z} - 1 - \ii \xi z) \nu( \dd z) .
\end{align}

\subsection*{Solution to Exercise \ref{exer:10.11}}
(a) As $X$ is a non-decreasing L\'evy process, it is of the form $X_t = \mu t + \int_0^\infty z N(t,\dd z)$ with $\mu \geq 0$, and its characteristic function satisfies $\Eb \ee^{\ii \xi X_t} = \ee^{t \psi(\xi)}$ with
\begin{align}
\psi(\xi)
	&=	\ii \mu \xi + \int_0^\infty (\ee^{\ii \xi z } - 1 ) \nu(\dd z) .
\end{align}
Therefore
\begin{align}
\phi(\lam)
	&=	- \psi(\ii \lam)
	=		\mu \lam + \int_0^\infty (1 - \ee^{- \lam z} ) \nu(\dd z) .
\end{align}
(b) From \eqref{eq:laplace.tau.m}, the first hitting time of a Brownian motion to level $t > 0$ satisfies $\Eb \ee^{- \lam \tau_t} = \ee^{- t \sqrt{2 \lam}}$.  By the strong Markov property and the stationary, independent increments of $W$, the process $\tau = (\tau_t)_{t \geq 0}$ has stationary independent increments and is non-decreasing, i.e., it is a non-decreasing L\'evy process (a subordinator).  Writing $\Eb \ee^{- \lam \tau_t} = \ee^{- t \phi(\lam)}$, we have $\phi(\lam) = \sqrt{2\lam}$.  Using the given identity with $r = 1/2$ and $v = \lam$,
\begin{align}
\phi(\lam)
	&=	\sqrt{2 \lam}
	=		\sqrt{2} \, \frac{1/2}{\Gam(1/2)} \int_0^\infty \frac{1}{z^{3/2}}(1 - \ee^{- \lam z}) \dd z .
\end{align}
Comparing with part (a) (which has $\mu = 0$ here), we identify the L\'evy measure
\begin{align}
\nu(\dd z)
	&=	 \frac{1}{\sqrt{2} \, \Gam(1/2)} \frac{1}{z^{3/2}} \dd z .
\end{align}

\subsection*{Solution to Exercise \ref{exer:10.12}}
(a) Using the definition of $M$,
\begin{align}
M_t
	&=	\ee^{- \lam t} u(t,X_t) + \int_0^t \ee^{- \lam s} g(X_s) \dd s
	=		\Eb \Big( \ee^{- \lam T} h(X_T) + \int_0^T \ee^{- \lam s} g(X_s) \dd s \Big| \Fc_t \Big) ,
\end{align}
which is a martingale because, for $0 \leq t_1 \leq t_2 \leq T$, the tower property gives
\begin{align}
\Eb[ M_{t_2} | \Fc_{t_1} ]
	&=	\Eb \Big( \ee^{- \lam T} h(X_T) + \int_0^T \ee^{- \lam s} g(X_s) \dd s \Big| \Fc_{t_1} \Big)
	=		M_{t_1} .
\end{align}
(b) Applying the It\^o formula for L\'evy-It\^o processes,
\begin{align}
\dd M_t
	&=	\ee^{- \lam t} ( \d_t - \lam + \Ac) u(t,X_t) \dd t + ( \ldots ) \dd W_t + \int_\Rb ( \ldots ) \Nt(\dd t, \dd z)  + \ee^{- \lam t} g(X_t) \dd t , \\
\Ac
	&=	\mu(t,x) \d_x + \tfrac{1}{2} \sig^2(t,x) \d_x^2 + \int_\Rb \nu(\dd z) \Big( \theta_{\gam(t,x,z)} - 1 - \gam(t,x,z) \d_x \Big) ,
\end{align}
where $\theta_\gam$ is the shift operator $\theta_\gam \varphi(x) := \varphi(x + \gam)$.  As $M$ is a martingale, the sum of the $\dd t$ terms must vanish, giving the PIDE and terminal condition
\begin{align}
0
	&=	( \d_t - \lam + \Ac) u + g , &
u(T,x) 
	&=	h(x) .
\end{align}
The terminal condition follows from the definition of $u$, since $u(T,X_T) = h(X_T)$.


\section{Exercises from Section \ref{sec:exercises-11}}

\subsection*{Solution to Exercise \ref{exer:11.1}}
We give formal derivations, assuming that the value function $V$ is sufficiently smooth to apply the It\^o formula and that limits and expectations may be interchanged.  In both cases, we write $\Eb_x \, \cdot \, := \Eb( \, \cdot \, | X_0^u = x)$.
\\[0.5em]
\textit{Indefinite time.}  Define $V(x) := \sup_{u} J(u)$, where $J(u)$ is the indefinite time cost functional with $X_0^u = x$.  If $x \in \d A$ then $\tau^u = 0$ and thus $V(x) = \varphi(x)$.  Now, fix $x \in A$, a constant $u \in \Ub$ and a small $\del > 0$.  Consider the control that is equal to $u$ on $[0,\del \wedge \tau^u)$ and is optimal thereafter.  As the dynamics \eqref{eq:dXu} are time-homogeneous and $X^u$ is Markov, the expected cost accumulated after time $\del \wedge \tau^u$ under this control is $\Eb_x V(X_{\del \wedge \tau^u}^u)$.  As this control is not necessarily optimal, we have
\begin{align}
V(x)
	&\geq \Eb_x \Big[ \int_0^{\del \wedge \tau^u} F(X_t^u,u) \dd t + V(X_{\del \wedge \tau^u}^u) \Big] , \label{eq:DPP.indefinite}
\end{align}
with equality if $u$ is replaced by an optimal control on $[0,\del \wedge \tau^u)$.  By the It\^o formula, we have
\begin{align}
V(X_{\del \wedge \tau^u}^u)
	&=	V(x) + \int_0^{\del \wedge \tau^u} \Ac^u V(X_t^u) \dd t + \int_0^{\del \wedge \tau^u} \d_x V(X_t^u) \sig(X_t^u,u) \dd W_t .
\end{align}
Inserting this into \eqref{eq:DPP.indefinite} and using the fact that the expectation of the It\^o integral is zero, we obtain
\begin{align}
0
	&\geq \Eb_x \int_0^{\del \wedge \tau^u} \Big( F(X_t^u,u) + \Ac^u V(X_t^u) \Big) \dd t .
\end{align}
As $A$ is open and $X^u$ is continuous, we have $\Pb(\tau^u > \del) \to 1$ as $\del \to 0$.  Thus, dividing by $\del$ and taking the limit as $\del \to 0$, we find $0 \geq F(x,u) + \Ac^u V(x)$ for every $u \in \Ub$, with equality for the optimal control.  Hence
\begin{align}
0
	&=	\sup_{u \in \Ub} ( \Ac^u V + F^u ) , &
x
	&\in A , &
V(x)
	&=	\varphi(x) , &
x
	&\in \d A ,
\end{align}
which is the HJB equation \eqref{eq:HJB-indefinite-time}.
\\[0.5em]
\textit{Infinite time.}  Define $V(x) := \sup_u J(u)$, where $J(u)$ is the infinite time cost functional with $X_0^u = x$.  Fix a constant $u \in \Ub$ and a small $\del > 0$, and consider the control that is equal to $u$ on $[0,\del)$ and is optimal thereafter.  As the dynamics \eqref{eq:dXu} are time-homogeneous, the expected discounted cost accumulated after time $\del$ under this control is
\begin{align}
\Eb_x \Eb \Big[ \int_\del^\infty \ee^{-\lam t} F(X_t^u,u_t) \dd t \Big| \Fc_\del \Big]
	&=	\Eb_x \Big[ \ee^{-\lam \del} \Eb \Big[ \int_0^\infty \ee^{-\lam s} F(X_{\del + s}^u,u_{\del + s}) \dd s \Big| \Fc_\del \Big] \Big]
	=		\Eb_x \ee^{-\lam \del} V(X_\del^u) .
\end{align}
Thus, we have
\begin{align}
V(x)
	&\geq \Eb_x \Big[ \int_0^\del \ee^{-\lam t} F(X_t^u,u) \dd t + \ee^{-\lam \del} V(X_\del^u) \Big] , \label{eq:DPP.infinite}
\end{align}
with equality if $u$ is replaced by an optimal control on $[0,\del)$.  By the It\^o formula, we have
\begin{align}
\ee^{-\lam \del} V(X_\del^u)
	&=	V(x) + \int_0^\del \ee^{-\lam t} ( \Ac^u - \lam ) V(X_t^u) \dd t + \int_0^\del \ee^{-\lam t} \d_x V(X_t^u) \sig(X_t^u,u) \dd W_t .
\end{align}
Inserting this into \eqref{eq:DPP.infinite} and using the fact that the expectation of the It\^o integral is zero, we obtain
\begin{align}
0
	&\geq \Eb_x \int_0^\del \ee^{-\lam t} \Big( F(X_t^u,u) + \Ac^u V(X_t^u) - \lam V(X_t^u) \Big) \dd t .
\end{align}
Dividing by $\del$ and taking the limit as $\del \to 0$, we find $0 \geq F(x,u) + \Ac^u V(x) - \lam V(x)$ for every $u \in \Ub$, with equality for the optimal control.  Hence
\begin{align}
0
	&=	\sup_{u \in \Ub} ( \Ac^u V + F^u ) - \lam V ,
\end{align}
which is the HJB equation \eqref{eq:HJB-infinite-time}.

\subsection*{Solution to Exercise \ref{exer:11.2}}
The controlled wealth process is time-homogeneous with generator (for fixed controls $u$ and $c$)
\begin{align}
\Ac^{u,c}
	&=	( \mu u - c ) \d_x + \tfrac{1}{2} \sig^2 u^2 \d_x^2 .
\end{align}
Thus, from the HJB equation \eqref{eq:HJB-infinite-time} for the infinite time cost functional (with $F^{u,c} = U(c)$ and with $\sup$ taken over both controls), the value function $V$ formally satisfies
\begin{align}
0
	&=	\sup_{u \in \Rb , c \geq 0} \Big( ( \mu u - c ) V' + \tfrac{1}{2} \sig^2 u^2 V'' + U(c) \Big) - \lam V .
\end{align}
Assuming $V' > 0$ and $V'' < 0$, the supremum over $u$ and $c$ can be computed separately.  We have
\begin{align}
u^*
	&=	- \frac{\mu V'}{\sig^2 V''} , &
\sup_{u} \Big( \mu u V' + \tfrac{1}{2} \sig^2 u^2 V'' \Big)
	&=	- \frac{\mu^2}{2 \sig^2} \frac{(V')^2}{V''} , \\
c^*
	&=	(V')^{-1/\gam} , &
\sup_{c \geq 0} \Big( U(c) - c V' \Big)
	&=	\frac{\gam}{1-\gam} (V')^{(\gam-1)/\gam} ,
\end{align}
where we have used $U'(c) = c^{-\gam}$.  Thus, the HJB equation becomes
\begin{align}
0
	&=	- \frac{\mu^2}{2 \sig^2} \frac{(V')^2}{V''} + \frac{\gam}{1-\gam} (V')^{(\gam-1)/\gam} - \lam V .
\end{align}
Following the hint, we try $V(x) = A U(x) = A x^{1-\gam}/(1-\gam)$.  Then $V' = A x^{-\gam}$ and $V'' = - \gam A x^{-\gam-1}$, so that
\begin{align}
- \frac{(V')^2}{V''}
	&=	\frac{A}{\gam} x^{1-\gam} , &
(V')^{(\gam-1)/\gam}
	&=	A^{(\gam-1)/\gam} x^{1-\gam} .
\end{align}
Inserting these expressions into the HJB equation and dividing by $A x^{1-\gam}/(1-\gam)$, we obtain
\begin{align}
0
	&=	\frac{1-\gam}{\gam} \frac{\mu^2}{2 \sig^2} + \gam A^{-1/\gam} - \lam &
	&\Rightarrow&
A^{-1/\gam}
	&=	\frac{1}{\gam} \Big( \lam - \frac{1-\gam}{\gam} \frac{\mu^2}{2 \sig^2} \Big)
	=		\nu ,
\end{align}
which is positive by assumption.  Thus, we have obtained the value function
\begin{align}
V(x)
	&=	\nu^{-\gam} \, \frac{x^{1-\gam}}{1-\gam} ,
\end{align}
and the optimal (Markov) controls
\begin{align}
u^*(x)
	&=	- \frac{\mu V'(x)}{\sig^2 V''(x)}
	=		\frac{\mu x}{\gam \sig^2} , &
c^*(x)
	&=	(V'(x))^{-1/\gam}
	=		A^{-1/\gam} x
	=		\nu x .
\end{align}
Thus, as in the Merton problem without consumption (Example \ref{ex:merton}), the investor keeps a constant proportion $\mu/(\gam \sig^2)$ of his wealth in the stock.  Additionally, he consumes at a rate that is a constant proportion $\nu$ of his wealth.


\section{Exercises from Section \ref{sec:exercises-12}}

\subsection*{Solution to Exercise \ref{exer:12.1}}
We assume $\lam > 0$, $\sig > 0$ and $K > 0$.  Throughout, we write $\varphi(x) := (K - x)^+$ and $\Ac := \mu x \d_x + \tfrac{1}{2} \sig^2 x^2 \d_x^2$ for the generator of $X$.  Note that $X$ lives on $(0,\infty)$.  From \eqref{ineq:variational} (with $F = 0$), the value function $V$ satisfies the variational inequality
\begin{align}
0
	&=	\max \{ (\Ac - \lam) V(x) , \varphi(x) - V(x) \} .
\end{align}
As it is never optimal to stop when $\varphi(x) = 0$, and as the payoff is largest when $X$ is small, we guess that the continuation region is of the form $D = (x^*,\infty)$ for some $x^* \in (0,K)$ and that the optimal stopping time is $\tau^* = \inf \{ t \geq 0 : X_t \leq x^* \}$.
\\[0.5em]
\textit{Continuation region.}  For $x > x^*$ we must have $(\Ac - \lam) V = 0$.  This is an Euler-type ODE.  Trying $V(x) = x^\beta$, we find that $\beta$ must satisfy
\begin{align}
\tfrac{1}{2} \sig^2 \beta (\beta - 1) + \mu \beta - \lam
	&=	0 &
	&\Rightarrow&
\beta_\pm
	&=	\frac{ ( \tfrac{1}{2} \sig^2 - \mu ) \pm \sqrt{ ( \mu - \tfrac{1}{2} \sig^2 )^2 + 2 \sig^2 \lam } }{\sig^2} .
\end{align}
As $\lam > 0$, the product of the roots $\beta_+ \beta_- = - 2 \lam / \sig^2$ is negative, and thus $\beta_- < 0 < \beta_+$.  As $V(x) \leq K$ for all $x$, the coefficient of $x^{\beta_+}$ must be zero.  Thus, writing $\beta := \beta_- < 0$, we have $V(x) = C x^\beta$ for $x > x^*$.
\\[0.5em]
\textit{Value matching and smooth fit.}  To determine $C$ and $x^*$ we require that $V$ be continuous and continuously differentiable at $x^*$ (the so-called \emph{smooth-fit} condition)
\begin{align}
C (x^*)^\beta
	&=	K - x^* , &
C \beta (x^*)^{\beta - 1}
	&=	- 1 .
\end{align}
Dividing the first equation by the second, we obtain $x^*/\beta = (K - x^*)/(-1)$, and thus
\begin{align}
x^*
	&=	\frac{\beta K}{\beta - 1} \in (0,K) , &
C
	&=	(K - x^*) (x^*)^{-\beta} .
\end{align}
Thus, we have
\begin{align}
V(x)
	&=	\left\{ \begin{aligned}
			&K - x , & & 0 < x \leq x^* , \\
			&(K - x^*) \Big( \frac{x}{x^*} \Big)^\beta , & & x > x^* ,
			\end{aligned} \right. &
\tau^*
	&=	\inf \{ t \geq 0 : X_t \leq x^* \} .
\end{align}
\textit{Verification.}  It remains to check that $V$ satisfies the variational inequality.  For $x > x^*$, we have $(\Ac - \lam) V = 0$ by construction.  Moreover, as $\beta < 0$ the function $V$ is convex on $(x^*,\infty)$, and by smooth fit its tangent line at $x^*$ is $K - x$.  Thus $V(x) \geq K - x$ for $x > x^*$, and since also $V > 0$, we have $V \geq \varphi$ there.  For $0 < x \leq x^*$, we have $V = \varphi = K - x$ and
\begin{align}
(\Ac - \lam) V(x)
	&=	- \mu x - \lam ( K - x )
	=		- \lam K + (\lam - \mu) x .
\end{align}
This expression is linear in $x$ and equals $-\lam K < 0$ at $x = 0$.  At $x = x^*$ we have
\begin{align}
- \lam K + (\lam - \mu) \frac{\beta K}{\beta - 1} \leq 0
	&\quad \Leftrightarrow \quad
(\lam - \mu) \beta \geq \lam ( \beta - 1 )
	\quad \Leftrightarrow \quad
\mu \beta \leq \lam ,
\end{align}
where we have used $\beta - 1 < 0$.  From the quadratic equation for $\beta$ we have $\mu \beta - \lam = - \tfrac{1}{2} \sig^2 \beta ( \beta - 1 ) < 0$, since $\beta < 0$.  Thus $(\Ac - \lam) V \leq 0$ on $(0,x^*]$, and the variational inequality is satisfied.
\\[0.5em]
\textit{Assumptions.}  The above construction requires only $\lam > 0$; no further restriction on $\mu$ is needed.  If instead $\lam = 0$, the roots are $\beta = 0$ and $\beta = 1 - 2\mu/\sig^2$, and a non-trivial solution requires a negative root, i.e., $\mu > \tfrac{1}{2} \sig^2$.  (If $\lam = 0$ and $\mu \leq \tfrac{1}{2} \sig^2$, then $X$ gets arbitrarily close to zero with probability one, and $V \equiv K$, which is not attained by any stopping time.)

\subsection*{Solution to Exercise \ref{exer:12.2}}
Let $\Ac := \mu(x) \d_x + \tfrac{1}{2} \sig^2(x) \d_x^2$ be the generator of $X$.  The function $u$ is of the form considered in Exercise \ref{exer:9.9} (with $h = 0$ and an infinite horizon).  Since $X$ is time-homogeneous, $u$ does not depend on $t$, and the same martingale argument (applied to $M_t = \ee^{-\lam t} u(X_t) + \int_0^t \ee^{-\lam s} g(X_s) \dd s$) shows that $u$ satisfies the ODE
\begin{align}
(\Ac - \lam) u + g
	&=	0 .
\end{align}
Next, using the strong Markov property of $X$, we have
\begin{align}
V(x)
	&=	\sup_{\tau \in \Tc} \Eb \Big( \int_0^\tau \ee^{- \lam t} f(X_t) \dd t
			+ \ee^{- \lam \tau} \Eb \Big( \int_\tau^\infty \ee^{- \lam (t-\tau)} g(X_t) \dd t \Big| \Fc_\tau \Big)  \Big| X_0 = x  \Big) \\
	&=	\sup_{\tau \in \Tc} \Eb \Big( \int_0^\tau \ee^{- \lam t} f(X_t) \dd t + \ee^{- \lam \tau} u(X_\tau) \Big| X_0 = x \Big) .
\end{align}
This is an optimal stopping problem with running reward $F = f$ and stopping reward $\varphi = u$.  Thus, we obtain from \eqref{ineq:variational} that
\begin{align}
0
	&=	\max\{ f(x) + (\Ac - \lam ) V(x) , u(x) - V(x) \} .
\end{align}


\section{Exercises from Section \ref{sec:exercises-13}}
\blu{To do.}

\end{document}